%% file: HermKstableI.tex
\documentclass[10pt]{amsart}

\usepackage{xr-hyper}

\input{Packages.tex}

\input{Theorem-styles.tex}

\input{Commands.tex}

\input{Authors.tex}

\input{Toc-tweak.tex}

\newcommand{\refone}[1]{\ref{#1}}
\newcommand{\eqrefone}[1]{\eqref{#1}}
\newcommand{\refoneitem}[1]{\ref{#1}}

\newcommand{\reftwo}[1]{\cite{Part-two}.\ref{II-#1}}

\newcommand{\refthree}[1]{\cite{Part-three}.\ref{III-#1}}

\input{External.tex}

\title[Hermitian K-theory for stable $\infty$-categories I: Foundations]{Hermitian K-theory for stable $\infty$-categories I:\\
Foundations}

\input{Dedication.tex}

\date{\today}

\begin{document}

\begin{abstract}
\input{Abstract-one.tex}

\end{abstract}

\maketitle
\tableofcontents

\section*{Introduction}
\input{Intro-one.tex}

\introsubsection{Acknowledgements}
\input{Acknowledgements.tex}

\input{Acknowledgements-one.tex}
\medskip

\input{Support.tex}

\section{Poincaré categories}
\label{section:poincare-cats}%
\input{PoincareCats.tex}

\section{Poincaré objects}
\label{section:poincare-objects}%
\input{PoincareObjects.tex}

\section{Poincaré structures on module categories}
\label{section:modules}%

\input{Modules.tex}

\section{Examples}
\label{section:examples}%
\input{Examples.tex}

\section{Monoidal structures and multiplicativity}
\label{section:multiplicative}%
\input{Multiplicative.tex}

\section{The category of Poincaré categories}
\label{section:cat-of-cats}%
\input{CatofCats.tex}

\section{Hyperbolic and metabolic Poincaré categories}
\label{section:metabolic}%
\input{Metabolic.tex}

%obsolete bibtex bibliography
%\bibliographystyle{amsalpha}
%\bibliography{bib}

%using amsrefs

%\bibsection*{Bibliography}

%\end{bibdiv}

\end{document}

%% file: Packages.tex
\usepackage[utf8]{inputenc}         % for utf8 encoding to avoid tex accents, for example
\usepackage[T1]{fontenc}            % for T1 out encoding

\usepackage{a4wide}                 % for a4 paper size
\usepackage{amsmath, amssymb}       % various ams packages
\usepackage{amsthm}                 % ams theorem styles
\usepackage{appendix}               % appendices where you want
\usepackage{tikz}                   % graphics and diagrams
\usepackage{tikz-cd}                % commutative diagrams
\tikzcdset{ % for the style of arrow tips
arrow style=tikz
}

\usepackage{graphicx}               % check if used ??
\usepackage{stackrel}               % improved stacking of things

\usepackage{color}                  % colors for our comments, for example, and for links
\usepackage[normalem]{ulem}         % long underline across lines 
\usepackage{enumitem}               % configure appearance of lists and references
\setlist[itemize]{leftmargin=*}
\setlist[enumerate]{leftmargin=*,label=\roman*),ref=\roman*)}
\newlist{subenumerate}{enumerate}{2}
\setlist[subenumerate]{leftmargin=*,label=\alph*),ref=\alph*)}

\usepackage{blkarray}               % for arrays and tables ?? should go

\definecolor{darkblue}{rgb}{0,0,0.6} % uses color package
\usepackage[ocgcolorlinks,colorlinks=true, citecolor=darkblue, filecolor=darkblue, linkcolor=darkblue, urlcolor=darkblue]{hyperref}         % for hyperlinks
\usepackage[alphabetic]{amsrefs}

\usepackage{bbm}                    % blackboard style fonts 
\usepackage{stix}                   % just used for nsimeq 
\usepackage{eucal}                  % for Lurie's C, actually redefines the whole mathcal script
\DeclareFontFamily{T1}{cbgreek}{}
\DeclareFontShape{T1}{cbgreek}{m}{n}{<-6>  grmn0500 <6-7> grmn0600 <7-8> grmn0700 <8-9> grmn0800 <9-10> grmn0900 <10-12> grmn1000 <12-17> grmn1200 <17-> grmn1728}{}
\DeclareSymbolFont{quadratics}{T1}{cbgreek}{m}{n}
\DeclareMathSymbol{\qoppa}{\mathord}{quadratics}{19}
\DeclareMathSymbol{\Qoppa}{\mathord}{quadratics}{21}

\usepackage{comment}                % for comments that can be toggled in or out (mainly used by Fabian)
\usepackage{xspace}                 % to deal with text macros and the following trailing (or not) space. used for words additive, etc.  

%% file: Theorem-styles.tex
\newtheoremstyle{thms}
	{}{}{\itshape}{}{\bfseries }{.}{ }
	{\thmname{#1} \thmnumber{#2}. \thmnote{\bfseries{[#3]}}}
\newtheoremstyle{thms2}
	{}{}{\itshape}{}{\bfseries }{.}{ }
	{}
\newtheoremstyle{ithreethms}
	{}{}{\itshape}{}{\bfseries }{}{ }
	{\thmname{#1} \thmnumber{#2}. \thmnote{\bfseries{[#3]}}}

\newtheoremstyle{name}
	{}{}{\itshape}{}{\bfseries }{.}{ }
	{\thmname{#1}\thmnumber{#2}\thmnote{\bfseries{[#3]}}}

\newtheoremstyle{defs}
	{}{}{\normalfont}{}{\bfseries }{.}{ }
	{\thmname{#1} \thmnumber{#2}. \thmnote{\bfseries{(#3)}}}
\newtheoremstyle{defs2}
	{}{12pt}{\normalfont}{}{\bfseries }{.}{ }
	{\thmname{#1}\thmnumber{#2}. \thmnote{\bfseries{(#3)}}}
\newtheoremstyle{ithreedefs}
	{}{}{\normalfont}{}{\bfseries }{}{ }
	{\thmname{#1} \thmnumber{#2}. \thmnote{\bfseries{(#3)}}}

\newtheoremstyle{rmk}
	{}{}{\normalfont}{}{\itshape }{.}{ }
        {}

\newtheoremstyle{claim}
	{}{}{\normalfont}{}{\itshape}{.}{ }
        {\thmname{#1} \thmnumber{#2}. \thmnote{#3}}

\theoremstyle{thms2}

\theoremstyle{thms2}

\newtheorem*{namedthm}{\namedthmname}
\newcounter{namedthm}

\makeatletter
{\endnamedthm}
\makeatother

\theoremstyle{rmk}

\theoremstyle{ithreethms}

\theoremstyle{ithreedefs}

\theoremstyle{thms2}

\swapnumbers         % so that numbers appear befor statements title, like 3.1.3 Theorem

\newcounter{rthree}
\theoremstyle{defs}
\newtheorem{definition-r-three}{Definition}[rthree]
\newtheorem{notation-r-three}[definition-r-three]{Notation}
\newtheorem{remark-r-three}[definition-r-three]{Remark}
\newtheorem{example-r-three}[definition-r-three]{Example}

\theoremstyle{thms}
\newtheorem{proposition-r-three}[definition-r-three]{Proposition}
\newtheorem{corollary-r-three}[definition-r-three]{Corollary}
\newtheorem{theorem-r-three}[definition-r-three]{Theorem}

\newcounter{rfour}
\theoremstyle{defs}
\newtheorem{definition-r-four}{Definition}[rthree]
\newtheorem{notation-r-four}[definition-r-four]{Notation}
\newtheorem{remark-r-four}[definition-r-four]{Remark}
\newtheorem{example-r-four}[definition-r-four]{Example}

\theoremstyle{thms}
\newtheorem{proposition-r-four}[definition-r-four]{Proposition}
\newtheorem{corollary-r-four}[definition-r-four]{Corollary}
\newtheorem{theorem-r-four}[definition-r-four]{Theorem}

\theoremstyle{thms}
\newtheorem{proposition}{Proposition}[subsection]
\newtheorem{theorem}[proposition]{Theorem}
\newtheorem{lemma}[proposition]{Lemma}
\newtheorem{corollary}[proposition]{Corollary}

\newtheorem{observation}[proposition]{Observation}

\theoremstyle{defs}
\newtheorem{definition}[proposition]{Definition}%[section]
\newtheorem{notation}[proposition]{Notation}
\newtheorem{variant}[proposition]{Variant}
\newtheorem{construction}[proposition]{Construction}
\newtheorem{example}[proposition]{Example}
\newtheorem{examples}[proposition]{Examples}
\newtheorem{remark}[proposition]{Remark}
\newtheorem{remarks}[proposition]{Remarks}
\newtheorem{warning}[proposition]{Warning}

\theoremstyle{defs2}

\theoremstyle{rmk}

\theoremstyle{claim}

%% file: Commands.tex
\newcommand{\defi}[1]{\emph{#1}}                 % highlight or emphasize a word \emph{that is being defined} (not just any emphasis).

\newcommand{\spacefont}{\mathcal}                % font used for the spaces

\newcommand{\signper}[1]{\({#1}\sig\)-oriented}  % sigma-oriented

\newcommand{\alp}{\alpha}
\newcommand{\bet}{\beta}
\newcommand{\gam}{\gamma}
\newcommand{\Del}{\Delta}
\newcommand{\del}{\delta}
\newcommand{\Lam}{\Lambda}
\newcommand{\Om}{\Omega}
\newcommand{\om}{\omega}
\newcommand{\Sig}{\Sigma}
\newcommand{\sig}{\sigma}
\newcommand{\vphi}{\varphi}
\newcommand{\eps}{\epsilon}

\newcommand{\lrar}{\longrightarrow}

\newcommand{\hrar}{\hookrightarrow}

\newcommand{\wtl}{\widetilde}

\newcommand{\st}{\stackrel}
\newcommand{\ovl}{\overline}
\newcommand{\ol}{\overline} % redundant, but ..

\newcommand{\adj}{\mathbin{% for adjunctions
\begin{tikzpicture}[baseline,thick] 
\coordinate (source) at (0ex,.5ex);
\coordinate (target) at (3ex,.5ex);
\draw[->] ([yshift=1ex]source) -- ([yshift=1ex]target); 
\draw[->] ([yshift=-.5ex]target) -- ([yshift=-.5ex]source);
\node at (1.5ex,.8ex) {$\scriptscriptstyle \perp$};
\end{tikzpicture}%
}}

\newcommand{\cocolon}{\nobreak \mskip6mu plus1mu \mathpunct{}\nonscript\mkern-\thinmuskip {:}\mskip2mu \relax}
\makeatletter
\newcommand{\dovl}[1]{\overline{\dbl@overline{#1}}}
\newcommand{\dbl@overline}[1]{\mathpalette\dbl@@overline{#1}}
\newcommand{\dbl@@overline}[2]{%
  \begingroup
  \sbox\z@{$\m@th#1\overline{#2}$}%
  \ht\z@=\dimexpr\ht\z@-2\dbl@adjust{#1}\relax
  \box\z@
  \ifx#1\scriptstyle\kern-\scriptspace\else
  \ifx#1\scriptscriptstyle\kern-\scriptspace\fi\fi
  \endgroup
}
\newcommand{\dbl@adjust}[1]{%
  \fontdimen8
  \ifx#1\displaystyle\textfont\else
  \ifx#1\textstyle\textfont\else
  \ifx#1\scriptstyle\scriptfont\else
  \scriptscriptfont\fi\fi\fi 3
}
\makeatother

\newcommand{\RR}{\mathbb{R}}               % real numbers
\newcommand{\CC}{\mathbb{C}}               % complex numbers
\newcommand{\ZZ}{\mathbb{Z}}               % integers
\newcommand{\QQ}{\mathbb{Q}}               % rationnal numbers
\newcommand{\FF}{\mathbb{F}}               % finite field

\newcommand{\im}{\mathrm{im}}              % image
\newcommand{\coker}{\mathrm{coker}}        % cokernel
\newcommand{\Mat}{\mathrm{Mat}}            % ring of matrices

\newcommand{\spec}{\mathrm{spec}}          % spectrum of a ring

\newcommand{\pt}{\mathrm{pt}}                    % one point set

\newcommand{\Hom}{\operatorname{Hom}}            % Hom functor \warn{set Hom?}
\newcommand{\End}{\operatorname{End}}            % End functor
\newcommand{\Aut}{\operatorname{Aut}}            % Automorphisms

\newcommand{\fib}{\operatorname{fib}}            % homotopy fibre
\newcommand{\hofib}{\operatorname{hofib}}        % also homotopy fiber, used only once in discrete rings
\newcommand{\cof}{\operatorname{cof}}            % homotopy quotient, cofiber

\newcommand{\colim}{\mathop{\mathrm{colim}}}     % colimits
\newcommand{\Map}{\operatorname{Map}}            % Mapping space 

\newcommand{\grpcr}{\iota}                       % groupoid core in general 
\newcommand{\core}{\mathrm{Cr}}                  % groupoid core functor on $\Catx$ or $\Catp$
\newcommand{\lag}{\operatorname{lag}}            % \warn{What is it?} 

\newcommand{\Ho}{\operatorname{Ho}}              % homotopy category 

\newcommand{\Pro}{\operatorname{Pro}}            % category of pro-objects 
\newcommand{\Ind}{\operatorname{Ind}}            % category of ind-objects
\newcommand{\coev}{\mathrm{coev}}                % universal coevaluation $C \times \I \to C_\I$

\newcommand{\op}{^\mathrm{op}}                   % opposite category
\newcommand{\mop}{\mathrm{op}}                   % opp functor 

\newcommand{\Ar}{\operatorname{Ar}}              % arrow cat
\newcommand{\Twar}{\operatorname{TwAr}}          % twisted arrow cat
\newcommand{\dom}{\mathrm{dom}}                  % domain functor
\newcommand{\cod}{\mathrm{cod}}                  % codomain functor
\newcommand{\target}{\mathrm{t}}                 % target functor
\newcommand{\source}{\mathrm{s}}                 % source functor

\newcommand{\ran}{\mathrm{ran}}                  % right adjoint of evaluation at an object in diagram categories

\newcommand{\cp}{\omega}                         % compact decoration 
\newcommand{\perf}{\mathrm{p}}                   % perfect decoration (for example in paramatrized spectra) 

\newcommand{\Mon}{\operatorname{Mon}}            % Monoid objects
\newcommand{\grp}{\mathrm{grp}}                  % group object decoration

\newcommand{\Bs}{\mathrm{B}}                     % classifying space 

\newcommand{\rF}{\mathrm{F}}   % free
\newcommand{\rU}{\mathrm{U}}   % forgetful

\newcommand{\id}{\mathrm{id}}                % identity
\newcommand{\const}{\mathrm{const}}          % constant function

\newcommand{\Spa}{{\mathcal Sp}}             % spectra  \warn{To be changed}
\newcommand{\Spaf}{\Spa^{\free}}             % finite spectra
\newcommand{\Sps}{\mathcal S}                % spaces  \warn{To be changed}
\newcommand{\Sfinast}{{\mathcal S^\mathrm{fin}_*}}% finite pointed spaces
\newcommand{\map}{\operatorname{hom}}          % same thing \warn{one of them should go.}

\renewcommand{\SS}{\mathbb{S}}                 % Sphere spectrum \warn{Renewcommand, could be dangerous!}

\newcommand{\sph}{\mathrm{S}}                  % sphere (not spectrum) 

\newcommand{\h}{\mathrm{h}}                    % homotopy decoration
\newcommand{\Ct}{\mathrm{C_2}}                 % group $\Ct$
\newcommand{\BC}{\Bs \Ct}                      % classifying space $\BC$
\newcommand{\hC}{{\h\Ct}}                      % decoration for homotopy $\Ct$ fixed points/orbits
\newcommand{\tC}{{\mathrm{t}\Ct}}              % decoration for $\Ct$ Tate
\newcommand{\geofix}{{\varphi\Ct}}             % geometric fixed points
\newcommand{\gC}{{\mathrm{g}\Ct}}              % genuine $\Ct$ decoration of a category, 
\newcommand{\Spagc}{\Spa^\gC}                  % genuine $\Ct$ Spectra 
\newcommand{\gCt}{\gC}                         % genuine $\Ct$ fixed points

\newcommand{\Th}{\mathrm{Th}}		       % Thom spectrum of a spherical fibration

\newcommand{\N}{\mathrm{N}}                    % norm
\newcommand{\Aeff}{\mathrm{A^{\hspace{-.2ex}eff}}} % span of finite $\Ct$ sets \warn{This is denoted by $\Span$ elsewhere, which is probably nicer.} 
\newcommand{\uAeff}{\underline{\Aeff}}         % $\Ct$ Burnside $\infty$-cat of Denis 

\newcommand{\OCt}{\mathbf{O}_{\Ct}}                        % transitive $\Ct$-sets.
\newcommand{\uCatp}{\underline{\mathrm{Cat}}^{\mathrm{p}}}
\newcommand{\uCath}{\underline{\mathrm{Cat}}^{\mathrm{h}}}
\newcommand{\uFunq}{\underline{\mathrm{Fun}}^{\mathrm{q}}}

\newcommand{\Ab}{\mathcal{A}b}                   % category of Abelian groups

\newcommand{\CMonoids}{\mathrm{CMon}}            % symmetric monoidal category of comm. monoids
\newcommand{\CMon}{\CMonoids}                    % duplicate \warn{should go}

\newcommand{\Mod}{\operatorname{Mod}}            % category of Modules
\newcommand{\LMod}{\mathrm{LMod}}                % category of left modules \warn{Why specify? used in Examples.}

\newcommand{\fg}{{\mathrm{f}}}                   % finitely generated decoration

\newcommand{\Proj}{\mathrm{Proj}}                % category of projective modules, finitely generated 
\newcommand{\Prof}[1]{\Proj(#1)}                 % same (used in part one) \warn{should be removed.} Example: \Prof{R} 
\newcommand{\Modcp}{\Mod^{\cp}}                  % compact Modules

\newcommand{\Modp}[1]{\Modcp_{#1}}               % compact Modules Example: \Modp{A}
\newcommand{\Modf}[1]{\mathrm{Mod}^\fg_{#1}}     % finitely generated modules Example: \Modf{A}
\newcommand{\gl}{\mathrm{gl}}                    % invertible modules
\newcommand{\Unimod}{\mathrm{Unimod}}            % groupoid of unimodular forms over a ring

\newcommand{\Einf}{{\mathrm{E}_\infty}}          % E infinity 
\newcommand{\Eone}{{\mathrm{E}_1}}               % E one 

\newcommand{\Alg}{\mathrm{Alg}}                  % Algebras (E1 or Einfty) \warn{should a command be added for both?)} 

\newcommand{\Ch}{\operatorname{Ch^b}}            % Chain complexes

\newcommand{\HHom}{\operatorname{Hom}^{\mathrm{cx}}} % Hom complex \warn{To be changed}
\newcommand{\qiso}{\mathrm{qIso}}                % quasi-isomorphisms
\newcommand{\GEM}{\mathrm{H}}                    % EM-functor

\newcommand{\Der}{\mathcal D}                    % derived category of a ring
\newcommand{\Dfree}{{\mathcal D}^\fg}            % perfect derived category of a ring
\newcommand{\Dperf}{{\mathcal D}^\perf}          % perfect derived category of a ring \warn{redundant}
\newcommand{\free}{\mathrm{f}}                   % decoration restricting a functor to the free derived subcat (used in verdier.tex)

\newcommand{\Psh}{\operatorname{Psh}}            % presheaves $\Fun(-,\Sps)$ \warn{could be changed.}
\newcommand{\Fr}{\mathrm{Fr}}                    % Tate Frobenius map

\newcommand{\Com}{\mathrm{Com}}       % commutative operad
\newcommand{\Day}{\mathrm{Day}}       % don't know
\newcommand{\MCom}{\mathrm{MCom}}     % operad of pairs: commutative algebra + module over it. 
\newcommand{\Op}{\mathrm{Op}_\infty}  % infinity-category of operads 

\newcommand{\Mul}{\mathrm{Mul}}       % multi-mapping space 

\newcommand{\Set}{\mathrm{Set}}                  % sets
\newcommand{\Top}{\mathrm{Top}}                 % topological spaces

\newcommand{\Cat}{\mathrm{Cat}_\infty}                  % small infinity categories
\newcommand{\CAT}{\mathrm{CAT}_\infty}                  % large infinity categories
\newcommand{\PrL}{\mathrm{Pr}^{\mathrm{L}}}  % presentable infinity categories
\newcommand{\Catx}{\Cat^\ex}     % small stable infinity categories
\newcommand{\Catxast}{\mathrm{Cat}^\ex_{\infty,\ast}}                % pointed stable infinity categories
\newcommand{\Catrex}{\mathrm{Cat}_\infty^\mathrm{rex}}   % infinity categories with finite colimits and right exact functors 

\newcommand{\CatSpa}{\mathrm{Cat}_{\Spa}}               % infinity cats enriched in spectra
\newcommand{\CatSpaPt}{\mathrm{Cat}_{\Spa,\ast}}               % infinity cats enriched in spectra and pointed
\newcommand{\Cath}{\mathrm{Cat}^\mathrm{h}_\infty}      % hermitian cats
\newcommand{\Catsb}{\mathrm{Cat}^\mathrm{sb}_\infty}    % cats with symmetric bilinear functor
\newcommand{\Catps}{\mathrm{Cat}^\mathrm{ps}_\infty}    % cats with perfect symm bil functor
\newcommand{\Catp}{\mathrm{Cat}^{\mathrm p}_\infty}     % poincare cats
\newcommand{\Span}{\operatorname{Span}}          % Span category
\newcommand{\Seq}{\operatorname{Seq}}            % category of short exact sequences
\newcommand{\seq}{\mathrm{seq}}                  % bisheaf associated to $\Seq$

\newcommand{\Fun}{\operatorname{Fun}}            % functor cat
\newcommand{\Nat}{\operatorname{Nat}}            % space of natural trafos
\newcommand{\nat}{\operatorname{nat}}            % spectrum of natural trafos
\newcommand{\ex}{\mathrm{ex}}                    % decoration for exact functors
\newcommand{\lex}{\mathrm{lex}}                  % decoration for left exact functors
\newcommand{\rex}{\mathrm{rex}}                  % decoration for right exact functors
\newcommand{\FunR}{\operatorname{Fun^R}}         % right adjoint functors 
\newcommand{\sifted}{\mathrm{sif}}               % sifted colimit preserving functor decoration
\newcommand{\filt}{\mathrm{filt}}                % filtered (co)-limit preserving decoration
\newcommand{\fin}{\mathrm{fin}}                  % finite colimit preserving decoration
\newcommand{\Funrex}{\Fun^{\rex}}		 % right exact functors
\newcommand{\diag}{\Del}                         % decoration for restriction alng diagonal

\newcommand{\BiFun}{\operatorname{BiFun}}        % functors reduced in both variables
\newcommand{\red}{\mathrm{red}}                  % reduced decoration for associated bireduced functor
\newcommand{\Funexc}{\Fun^{2-\mathrm{exc}}}      % excisive functors
\newcommand{\Funpoly}{\Fun^{2-\mathrm{poly}}_*}  % reduced polynomial functors of degree 2
\newcommand{\App}{\mathrm{P}}                    % 2-excisive approximation
\newcommand{\rT}{\mathrm{T}}                     % used to construct 2-excisive approximation
\newcommand{\BiFib}{\mathrm{BiFib}}              % bifibrations
\newcommand{\RFib}{\operatorname{RFib}}          % right fibrations over a category
\newcommand{\LaxAr}{\mathrm{LaxAr}}             % cartesian fibrations

\newcommand{\Funx}{\operatorname{Fun^{ex}}}      % exact functors
\newcommand{\Funq}{\operatorname{Fun^q}}         % quadratic functors
\newcommand{\Funqhom}{\operatorname{Fun^{hom}}}  % homogeneous quadratic functors
\newcommand{\Funqcoh}{\operatorname{Fun^{coh}}}  % cohomogeneous quadratic functors
\newcommand{\Funnq}{\operatorname{Fun^{nq}}}     % repble quadratic functors
\newcommand{\Funs}{\operatorname{Fun^s}}         % sym, bilinear functors
\newcommand{\Funb}{\operatorname{Fun^b}}         % bilinear functors
\newcommand{\Funpb}{\operatorname{Fun^{pb}}}	 % perfect bilinear functors
\newcommand{\Funnb}{\operatorname{Fun^{nb}}}    % non-degenerate bilinear functors
\newcommand{\Funh}{\operatorname{Fun^h}}         % hermitian functors
\newcommand{\Funp}{\operatorname{Fun^p}}         % Poincare functors
\newcommand{\Funns}{\operatorname{Fun^{ns}}}     % representable sym, bilinear functors
\newcommand{\Funps}{\operatorname{Fun^{ps}}}     % perfect sym bilinear functors

\newcommand{\Pairings}{\operatorname{Pair}}      % pairings 
\newcommand{\pair}{\mathrm{pair}}                % subscript pairings

\newcommand{\Catb}{\mathrm{Cat^{b}}}             % right bi-exact bi-sheaf functors to spaces
\newcommand{\Catpb}{\mathrm{Cat^{pb}}}           % perfect bilinear categories
\newcommand{\BiFibp}{\mathrm{BiFib^{p}}}         % perfect bifibrations

\newcommand{\Motp}[1][]{{\if\relax\detokenize{#1}\mathrm{Mot^p}\relax\else\mathrm{Mot^p_{#1}}\fi}}  % poincare motives
\newcommand{\Motpun}[1][]{\mathrm{Mot^p_{un\if\relax\detokenize{#1}\relax\else{,}#1\fi}}}           % unstable poincare motives
\newcommand{\Motpbord}[1][]{{\if\relax\detokenize{#1}\mathrm{BMot^p}\relax\else\mathrm{BMot^p_{#1}}\fi}} % poincare-witt motives
\newcommand{\Mot}[1][]{{\if\relax\detokenize{#1}\mathrm{Mot}\relax\else\mathrm{Mot_{#1}}\fi}}            % noncommutative motives

\newcommand{\rO}{\mathrm{O}}                     % orthogonal group 
\newcommand{\Un}{\mathrm{U}}                     % unitary group 
\newcommand{\SO}{\mathrm{SO}}                    % special orthogonal group (needs precise definition) 
\newcommand{\Sp}{\mathrm{Sp}}                    % symplectic group
\newcommand{\Spin}{\mathrm{Spin}}                % spin group
\newcommand{\String}{\mathrm{String}}             % string group 
\newcommand{\ku}{\mathrm{ku}}
\newcommand{\KU}{\mathrm{KU}}

\newcommand{\MO}{\mathrm{MO}}                    % MO spectrum
\newcommand{\MSO}{\mathrm{MSO}}                  % MSO spectrum (Thom)
\newcommand{\MSp}{\mathrm{MSp}}                  % MSp spectrum 
\newcommand{\MU}{\mathrm{MU}}                    % MU spectrum (Cobordism)
\newcommand{\MSpin}{\mathrm{MSpin}}                % MSO spectrum (Thom)
\newcommand{\MString}{\mathrm{MString}}              % MSp spectrum 

\newcommand{\LA}{\mathrm{LA}}                    % Weiss and Williams spectrum
\newcommand{\Pic}{\operatorname{Pic}}            % Picard group 

\newcommand{\ko}{\mathrm{ko}}                    % ko spectrum \warn{used in dicrete rings}
\newcommand{\K}{\operatorname K}                 % K-Theory \warn{spectrum or group?}
\newcommand{\KO}{\mathrm{KO}}                    % How symmetric GW appears in literature
\renewcommand{\L}{\operatorname L}               % L-theory groups? \warn{Renewcommand, dangerous!}
\newcommand{\W}{\operatorname W}                 % Witt group
\newcommand{\Witt}{\W}                           % Also Witt group \warn{redundant}

\newcommand{\GWspace}{\operatorname{\spacefont{GW}}} % GW space
\newcommand{\GW}{\operatorname{GW}}              % Grothendieck-Witt spectrum 
\newcommand{\U}{\operatorname{U}}                % Karoubi's U-theory spectrum fib(K->GW)
\newcommand{\fgt}{\mathrm{fgt}}                  % forgetful functor, mostly $\GW \to \K$ \warn{but also $\GW \to \L$ sometimes} 
\newcommand{\WW}{\textrm{WW}}                    % Weiss-Williams decoration for $\L \to \K^{\tC}$

\newcommand{\Q}{\operatorname{Q}}                % Q-construction
\newcommand{\arr}{\mathrm{ar}}                   % decoration as in Q arr, used once in Cobcats.tex
\newcommand{\fpm}{\lambda}			 %form parameter (formally \Lambda and Q)
\newcommand{\dfpm}{\check{\fpm}}			 %dual form parameter
\newcommand{\fpmg}{Q}				 %value group in a form parameter
\newcommand{\gfpm}{{\g\fpm}}
\newcommand{\gdfpm}{{\g\dfpm}}

\newcommand{\B}{\mathrm{B}}                      % cross effect \warn{wrongly used sometimes}

\newcommand{\Bil}{\mathrm{B}}               % associated bilinear form or actually any bilinear form (wrong use)
\newcommand{\Biltwo}{\Bil'}                      % some other bilinear form
\newcommand{\Lin}{\Lambda}               % associated linear part

\newcommand{\Dual}{\mathrm{D}}              % associated duality
\newcommand{\ev}{\mathrm{ev}}                    % evaluation = double dual inclusion

\newcommand{\ad}{\mathrm{ad}}                    % adjoint map 

\newcommand{\inv}[1]{\overline{#1}}              % involution Example: \inv{(-)}

\newcommand{\iq}{{\operatorname{nat}_\QF^\QFtwo}}% internal quad functor C,C'  \warn{Maybe to be changed.}
\newcommand{\tf}{\mathrm{id}}                    % taut form in $\Funx((\C,\QF),(\C,\QF))$

\newcommand{\qshift}[1]{^{[#1]}}                 % shift (post-composition with suspension) of a hermitian category, Example: \qshift{1} 
\newcommand{\Hyp}{\operatorname{Hyp}}            % Hyperbolic categories
\newcommand{\gHyp}{\operatorname{gHyp}}          % genuine refinement of $\Hyp$ (on Catp) 
\newcommand{\Met}{\operatorname{Met}}            % Metabolic categories
\newcommand{\ilag}{\mathrm{can}}                 % natural map $\Hyp \to \Met$
\newcommand{\met}{\mathrm{met}}                  % natural map $\Met \to Id$ \warn{??} 
\newcommand{\triv}{\mathrm{triv}}                % trivial inclusion $(\C,\QF) \to \Met(\C,\QF)$ 
\newcommand{\dlag}{\mathrm{dlag}}                % functor $\Met(\C,\QF) \to \Hyp(\C,\QF)$
\newcommand{\dilag}{\mathrm{dcan}}                % functor $\Met(\C,\QF) \to \Hyp(\C)

\newcommand{\hyp}{\mathrm{hyp}}                  % decoration for the hyperbolic quadratic functor $\QF_\hyp$

\newcommand{\cl}{\mathrm{cl}}                    % Short for classic, used for classical GW in subscript: $\GW_\cl$

\newcommand{\sym}{\mathrm{s}}                    % decoration for symmetric quadratic functor as in $\QF^\sym$
\newcommand{\s}{\mathrm{s}}                      % same as above \warn{should go.}
\newcommand{\vis}{\mathrm{v}}                    % decoration for visible quadratic functor as in $\QF^\vis$ 
\newcommand{\qdr}{\mathrm{q}}                     % decoration for quadratic quadratic functor as in $\QF^\qdr$ 
\newcommand{\g}{\mathrm{g}}                      % decoration for genuine quadratic functor as in $\QF^\g$
\newcommand{\gs}{\mathrm{gs}}                    % decoration for genuine symmetric quadratic functor as in $\QF^\gs$
\newcommand{\gq}{\mathrm{gq}}                    % decoration for genuine quadratic quadratic functor as in $\QF^\gq$
\newcommand{\gev}{\mathrm{ge}}                   % decoration for even genuine quadratic functor as in $\QF^\gev$
\newcommand{\tate}{\mathrm{t}}                   % decoration for tate quadratic functor

\newcommand{\uni}{\mathrm{u}}                    % universal decoration for a quadratic functor 

\newcommand{\swap}{\mathrm{swap}}                % decoration for quadratic functor on pairings
\newcommand{\proj}{\mathrm{proj}}                % decoration for quadratic functor 

\newcommand{\Poinc}{\mathrm{Pn}}                 % space of poincare objects
\newcommand{\Poincdel}{\Poinc^{\partial}}        % space of poincare objects with boundary

\newcommand{\refl}{\mathrm{refl}}                % subcategory of reflexive objects
\newcommand{\catforms}{\mathrm{He}}              % category of hermitian objects 
\newcommand{\spsforms}{\mathrm{Fm}}              % space of hermitian forms 
\newcommand{\spzero}{\mathrm{sp}_0}     % WW set of Poincaré objects
\newcommand{\R}{\mathcal{R}}			% WW retractive spaces
\newcommand{\sR}{\mathrm{s}\mathcal{R}} % WW stable retractive spaces
\newcommand{\fd}{\mathrm{fd}}			% decoration finitely dominated

\newcommand{\X}{\mathcal{X}}               % a space

\newcommand{\cP}{\mathcal{P}}              % a poset
\newcommand{\T}{\mathcal{T}}               % poset of subsets of $[n]$  

\newcommand{\I}{\mathcal{I}}               % a diagram category
\newcommand{\J}{\mathcal{J}}               % a diagram category
\newcommand{\cK}{\mathcal{K}}              % a diagram category

\newcommand{\A}{\mathcal{A}}               % an infty cat
\newcommand{\cB}{\mathcal{B}}              % an infty cat
\newcommand{\E}{\mathcal{E}}               % an infty cat (often stable)

\newcommand{\M}{\mathcal{M}}               % an infty cat (often total space of a fibration)

\newcommand{\cO}{\mathcal{O}}              % an infty operad

\newcommand{\C}{\mathcal C}                % a stable infty cat
\newcommand{\Ctwo}{{{\mathcal C}'}}        % a stable infty cat
\newcommand{\Cthree}{{{\mathcal C}''}}     % a stable infty cat
\newcommand{\D}{\mathcal{D}}               % a stable infty cat
\newcommand{\F}{\mathcal{F}}               % a functor, usually from $\Catp$
\newcommand{\G}{\mathcal{G}}               % a functor, usually from $\Catp$

\newcommand{\x}{x}                         % an object in a stable cat
\newcommand{\xtwo}{x'}                     % an object in a stable cat
\newcommand{\y}{y} 			   			   % an object in a stable cat
\newcommand{\ytwo}{y'}			   		   % an object in a stable cat
\newcommand{\z}{z}                         % an object in a stable cat
\newcommand{\ztwo}{z'}					   % an object in a stable cat
\newcommand{\zthree}{z''}				   % an object in a stable cat
\newcommand{\w}{w}                         % an object in a stable cat
\newcommand{\wtwo}{w'}					   % an object in a stable cat
\newcommand{\cob}{w}                       % an object, mostly a Lagrangian
\newcommand{\cobtwo}{\cob'}				   % an object, mostly a Lagrangian

\newcommand{\RF}{\mathcal{R}}              % an arbitrary reduced functor

\newcommand{\QF}{\Qoppa}                   % a quadratic functor   
\newcommand{\QFtwo}{{\QF'}}                % another quadratic functor 
\newcommand{\QFthree}{{\QF''}}             % another quadratic functor 
\newcommand{\QFD}{\Phi}                    % another quadratic functor 
\newcommand{\QFE}{\Psi}                    % another quadratic functor

\newcommand{\qone}{q}                      % a Poincaré form 
\newcommand{\qtwo}{q'}                     % another Poincaré form
\newcommand{\pone}{p}					   % another Poincaré form
\newcommand{\ptwo}{p'}					   % another Poincaré form

\newcommand{\Lag}{\mathcal{L}}             % isotropic subcat (lagrangian...)

%% file: Authors.tex
\author[Calmès]{Baptiste Calmès}
\address{Université d'Artois, Laboratoire de Mathématiques de Lens (LML), UR 2462, Lens, France}
\email{baptiste.calmes@univ-artois.fr}

\author[Dotto]{Emanuele Dotto}
\address{University of Warwick; Mathematics Institute; Coventry, United Kingdom}
\email{emanuele.dotto@warwick.ac.uk}

\author[Harpaz]{Yonatan Harpaz}
\address{Université Paris 13; Institut Galilée; Villetaneuse, France}
\email{harpaz@math.univ-paris13.fr}

\author[Hebestreit]{Fabian Hebestreit}
\address{RFWU Bonn; Mathematisches Institut; Bonn, Germany}
\email{f.hebestreit@math.uni-bonn.de}

\author[Land]{Markus Land}
\address{Department of Mathematical Sciences, University of Copenhagen, 2100 Copenhagen, Denmark}
\email{markus.land@math.ku.dk}

\author[Moi]{Kristian Moi}
\address{KTH; Institutionen för matematik; Stockholm, Sweden}
\email{kristian.moi@gmail.com}

\author[Nardin]{Denis Nardin}
\address{Universität Regensburg; Mathematisches Institut; Regensburg, Germany}
\email{denis.nardin@ur.de}

\author[Nikolaus]{Thomas Nikolaus}
\address{WWU Münster; Mathematisches Institut; Münster, Germany}
\email{nikolaus@uni-muenster.de}

\author[Steimle]{Wolfgang Steimle}
\address{Universität Augsburg; Institut für Mathematik; Augsburg, Germany}
\email{wolfgang.steimle@math.uni-augsburg.de}

%% file: Toc-tweak.tex
\makeatletter

\renewcommand{\tocsection}[3]{%
\indentlabel{\@ifnotempty{#2}{\parbox[b]{3ex}{\bfseries\ignorespaces#1 #2}}}\bfseries#3} 

\renewcommand{\tocsubsection}[3]{%
\indentlabel{\@ifnotempty{#2}{\hspace{1.6em}\parbox[b]{5ex}{\ignorespaces#1 #2}}}#3}

\renewcommand{\tocsubsubsection}[3]{%
\indentlabel{\@ifnotempty{#2}{\hspace{3.9em}\parbox[b]{5ex}{\ignorespaces#1 #2}}}#3}

\makeatother

\DeclareRobustCommand{\SkipTocEntry}[5]{} 

\newcommand{\introsubsection}[1]{\addtocontents{toc}{\SkipTocEntry}\subsection*{#1}}
\newcommand{\bibsubsubsection}[1]{\addtocontents{toc}{\SkipTocEntry}\subsubsection*{#1}}

%% file: External.tex
\newcommand{\papertwo}{Paper \cite{Part-two}\xspace}
\newcommand{\paperthree}{Paper \cite{Part-three}\xspace}
\newcommand{\paperfour}{Paper \cite{Part-four}\xspace}

\newcommand{\refundefined}[1]{\textcolor{red}{Undefined ref}}

%% file: Dedication.tex
\dedicatory{To Andrew Ranicki.}

%% file: Abstract-one.tex
This paper is the first in a series in which we offer a new framework for hermitian \(\K\)-theory in the realm of stable \(\infty\)-categories. Our perspective yields solutions to a variety of classical problems involving Grothendieck-Witt groups of rings and clarifies the behaviour of these invariants when \(2\) is not invertible.

In the present article we lay the foundations of our approach by considering Lurie's notion of a Poincaré \(\infty\)-category, which permits an abstract counterpart of unimodular forms called Poincaré objects. We analyse the special cases of hyperbolic and metabolic Poincaré objects, and establish a version of Ranicki's algebraic Thom construction.
For derived \(\infty\)-categories of rings,
we classify all Poincaré structures and study in detail the process of deriving them from classical input, thereby locating the usual setting of forms over rings within our framework. We also develop the example of visible Poincaré structures on \(\infty\)-categories of parametrised spectra, recovering the visible signature of a Poincaré duality space.

We conduct a thorough investigation of the global structural properties of Poincaré \(\infty\)-categories, showing in particular that they form a bicomplete, closed symmetric monoidal \(\infty\)-category.
We also study the process of tensoring and cotensoring a Poincaré \(\infty\)-category over a finite simplicial complex, a construction featuring prominently in the definition of the \(\L\)- and Grothendieck-Witt spectra that we consider in the next instalment.

Finally, we define already here the 0-th Grothendieck-Witt group of a Poincaré \(\infty\)-category using generators and relations. We extract its basic properties, relating it in particular to the 0-th \(\L\)- and algebraic \(\K\)-groups,
a relation upgraded in the second instalment to a fibre sequence of spectra which plays a key role in our applications.

%% file: Intro-one.tex
Quadratic forms are among the most ubiquitous notions in mathematics. In his pioneering paper~\cite{witt}, Witt suggested a way to understand quadratic forms over a field \(k\)
in terms of an abelian group \(\W^{\qdr}(k)\), now known as the \emph{Witt group} of quadratic forms. By definition, the Witt group is generated by isomorphism classes \([V,q]\) of finite dimensional \(k\)-vector spaces
equipped with a unimodular quadratic form \(q\),
where we impose the relations \([V \oplus V',q \perp q'] = [V,q]+[V',q']\) and declare as trivial the classes of \emph{hyperbolic forms} \([V \oplus V^*,h]\) given by the canonical pairing between \(V\) and its dual \(V^*\).
In arithmetic geometry the Witt group became an important invariant of fields, related to their Milnor \(\K\)-theory and Galois cohomology via the famous Milnor conjecture.

The definition of the Witt group naturally extends from fields to commutative rings \(R\), where one replaces vector spaces by finitely generated projective \(R\)-modules. %
More generally, instead of starting with a commutative ring \(R\) and taking \(R\)-valued forms, one can study unimodular hermitian forms valued in an invertible \((R \otimes R)\)-module \(M\) equipped with an involution, a notion which makes sense also for non-commutative \(R\). This includes for example the case of a ring \(R\) with anti-involution by considering \(M=R\), and also allows to consider %
skew-quadratic forms by changing the involution on \(M\) by a sign.
Quadratic forms at this level of generality also show up naturally in the purely geometric context of \emph{surgery theory} through the
 \emph{quadratic \(\L\)-groups} of the group ring \(\ZZ[\pi_1(X)]\) for a topological space \(X\). The latter groups, whose name, coined by Wall, suggests their relation with algebraic \(\K\)-theory, are a sequence of groups \(\L^\qdr_i\)
associated to a ring with anti-involution \(R\), or more generally, a ring equipped with an invertible \((R \otimes R)\)-module with involution \(M\) as above, with \(\L^{\qdr}_0(R,M)\) being the Witt group of \(M\)-valued quadratic forms over \(R\).
They are \(4\)-periodic, or more precisely, satisfy
the skew-periodic relation \(\L^{\qdr}_{n+2}(R,M) \cong \L^{\qdr}_{n}(R,-M)\), where \(-M\) is obtained from \(M\) by twisting the involution by a sign. In particular, for a ring with anti-involution \(R\) the even quadratic \(\L\)-groups consist of the Witt groups of quadratic and skew-quadratic forms.

To obtain richer information about quadratic forms over a given \(R\), the Witt group \(\W^{\qdr}(R,M)\)
was often compared to the larger group generated by the isomorphism classes of unimodular quadratic \(M\)-valued forms \([P,q]\) over \(R\) under the relation \([P \oplus P',q \perp q'] = [P,q]+[P',q']\), but without taking the quotient by hyperbolic forms. The latter construction leads to the notion of the \emph{Grothendieck-Witt group} \(\GW^{\qdr}_0(R,M)\) of quadratic forms. The Witt and Grothendieck-Witt groups are then related by an exact sequence
\begin{equation}
\label{equation:short}%
{\K_0(R)}{_{\Ct}} \xrightarrow{\hyp} \GW^{\qdr}_0(R,M) \to \W^{\qdr}(R) \to 0,
\end{equation}
where the first term denotes the orbits for the \(\Ct\)-action on the \(\K\)-theory group \(\K_0(R)\) which sends the class of a finitely generated projective \(R\)-module \(P\) to the class of its \(M\)-dual \(\Hom_R(P,M)\). The left hand map then sends \([P]\) to the class of the associated hyperbolic form on \(P \oplus \Hom(P,M)\), and is invariant under this \(\Ct\)-action. %
The sequence \eqrefone{equation:short} can often be used to compute \(\GW^{\qdr}_0(R,M)\) from the two outer groups, and consequently obtain more complete information about quadratic forms.
For example, in the case of the integers this sequence is split short exact and we have an isomorphism \(\W^{\qdr}(\ZZ) \cong \ZZ\) given by taking the signature divided by \(8\)
and an isomorphism \(\K_0(\ZZ)_{\Ct} \cong \ZZ\) given by the dimension.

In this paper we begin a four-part investigation revisiting classical questions about Witt, Grothendieck-Witt, and \(\L\)-groups of rings from a new perspective. %
One of our main motivating applications %
is to extend the short exact sequence \eqrefone{equation:short} to a long exact sequence involving Quillen's higher \(\K\)-theory and the higher Grothendieck-Witt groups \(\GW^\qdr_i(R,M)\) introduced by  Karoubi and Villamayor~\cite{karoubi-villamayor}, see below for more details. %
In this paper we, among many other things, define abelian groups \(\L_i^{\gq}(R,M)\), called \emph{genuine quadratic \(\L\)-groups}, which are the correct higher Witt groups from this point of view: We show in \paperthree that we have \(\L_0^{\gq}(R,M) = \W^{\qdr}(R,M)\) and that the sequence \eqrefone{equation:short} can be extended to a long exact sequence involving the groups \(\L_i^{\gq}(R,M)\) which starts off as
\[
\ldots \to \GW^{\qdr}_1(R,M)  \to \L_1^{\gq}(R,M) \to {\K_0(R,M)}{_{\Ct}} \xrightarrow{\hyp} \GW^{\qdr}_0(R,M) \to \L_0^{\gq}(R,M) \to 0 \ .
\]
The groups \(\L_i^{\gq}(R,M)\) are generally different from Wall's quadratic \(\L\)-groups, and in particular are usually not \(4\)-periodic. They are, however, relatively accessible for study by means of \emph{algebraic surgery}. Combining this with the above long exact sequence allows us to obtain many new results about the Grothendieck-Witt groups \(\GW^\qdr_i(R)\) of rings in \paperthree. For example, we obtain an essentially complete calculation of these groups in the case of the integers \(R = \mathbb{Z}\). In what follows we give more background, outline our approach and its main applications, and elaborate more on the content of the present paper. \\

\introsubsection{Background}

The higher Grothendieck-Witt groups \(\GW^\qdr_i(R,M)\) mentioned above %
were first defined by Karoubi and Villamayor~\cite{karoubi-villamayor} by applying Quillen's foundational techniques from algebraic \(\K\)-theory. This is done by producing a homotopy-theoretical refinement of the 0-th Grothendieck-Witt group into
a \defi{Grothendieck-Witt space} and then defining \(\GW_i(R,M)\) as the \(i\)-th homotopy group of this space.
Given \(R\) and \(M\) as above, one organises the collection of
unimodular quadratic %
\(M\)-valued forms \((P,q)\) into a groupoid \(\Unimod^{\qdr}(R,M)\), %
which may be viewed as an \(\Einf\)-space using the symmetric monoidal structure on \(\Unimod^{\qdr}(R)\) arising from the orthogonal sum. One can then take its group completion to obtain an \(\Einf\)-group
\[
\GWspace^{\qdr}_{\cl}(R,M) := \Unimod^{\qdr}(R,M)^\grp,
\]
whose group of components is the Grothendieck-Witt group described above. Here the subscript \(\cl\) stands for classical, and is meant to avoid confusion with the constructions of the present paper series. This construction can equally well be applied for other interesting types of forms, such as symmetric bilinear, or symmetric bilinear forms which admit a quadratic refinement, also known as \defi{even} forms, and these can be taken with values in an arbitrary invertible module with involution \(M\) as above. Taking the polarisation of a quadratic form determines maps
\[
\GWspace^\qdr_{\cl}(R,M) \longrightarrow \GWspace^\ev_{\cl}(R,M) \longrightarrow \GWspace^\sym_{\cl}(R,M),
\]
which are equivalences if \(2\) is a unit in \(R\).
In this latter case Grothendieck-Witt groups are generally much more accessible. For example, if \(2\) is invertible in $R$, %
Schlichting~\cite{schlichting-derived} has produced a (generally non-connective) delooping of the Grothendieck-Witt space to a Grothendieck-Witt \defi{spectrum} \(\GW_{\cl}(R,M)\), in which case the forgetful and hyperbolic maps can be refined to spectrum level \(\Ct\)-equivariant maps
\[
\K(R) \xrightarrow{\hyp} \GW_{\cl}(R,M) \xrightarrow{\fgt} \K(R).
\]
He then showed in loc.\ cit.\ %
that the cofibre of the induced map
\begin{equation}
\label{equation:hyp}%
\K(R)_{\hC} \to \GW_{\cl}(R,M)
\end{equation}
has \(4\)-periodic homotopy groups, whose even values are given by the Witt groups \(\Witt(R,M)\) and \(\Witt(R,-M)\). More precisely, Schlichting's identification of these homotopy groups matches the \(\L\)-groups of Ranicki-Wall, which has lead to the folk theorem that, if \(2\) is a unit in \(R\), then the cofibre of~\eqrefone{equation:hyp} is naturally equivalent to Ranicki's \(\L\)-spectum \(\L^q(R,M)\) from \cite{Ranickiblue}.
This allows one to produce an extension of \eqrefone{equation:short} to a long exact sequence, whenever \(2\) is invertible, and obtain information about higher Grothendieck-Witt groups from information about higher \(\K\)- and \(\L\)-groups. A closely related connection between Grothendieck-Witt spaces with coefficients in \(\pm M\) when \(2\) is invertible was established by Karoubi in his influential paper~\cite{Karoubi-Le-theoreme-fondamental},
where he proved what is now known as
\emph{Karoubi's fundamental theorem}, forming one of the conceptual pillars of hermitian \(\K\)-theory, as well as part of its standard tool kit. It permits, for example, to inductively deduce results on higher Grothendieck-Witt groups from information about algebraic \(\K\)-theory and about the low order Grothendieck-Witt groups \(\GW_{0}(R,\pm M)\) and \(\GW_1(R,\pm M)\).

By contrast, when \(2\) is not invertible none of these assertions hold as stated. In particular, the relation between Grothendieck-Witt theory and \(\L\)-theory remained, in this generality,
completely mysterious.
Karoubi, in turn, conjectured in~\cite{karoubi-periodicity} that his fundamental theorem should have an extension to general rings, relating Grothendieck-Witt spaces for two different form parameters,
as was also suggested earlier by Giffen~\cite{williams-quadratic}.
In the context of motivic homotopy theory, crucial properties such as dévissage and
\(\mathbf{A}^1\)-invariance of Grothendieck-Witt theory were only known to hold when \(2\) is invertible by the work of Schlichting and Hornbostel~\cite{hornbostel}, \cite{hornbostel-schlichting}.
Consequently, hermitian \(\K\)-theory was available to study as a motivic spectrum exclusively over \(\ZZ[\tfrac 1 2]\), see~\cite{hornbostel-motivic}.
Finally, while all the above tools could be used to calculate Grothendieck-Witt groups of rings in which \(2\) is invertible, such as the ring \(\ZZ[\tfrac 1 2]\) whose Grothendieck-Witt groups were calculated by Berrick and Karoubi in~\cite{berrick-karoubi}, higher Grothendieck-Witt groups of general rings remain largely unknown.

\introsubsection{Hermitian \(\K\)-theory of Poincaré \(\infty\)-categories}

The goal of the present paper series is to offer new foundations for hermitian \(\K\)-theory in a framework that unites its algebraic and surgery theoretic incarnations and that is robustly adapted to handle the subtleties involved when \(2\) is not invertible.
We begin by situating hermitian \(\K\)-theory in the general framework of \emph{Poincaré \(\infty\)-categories}, a notion suggested by Lurie in his treatise of \(\L\)-theory~ \cite{Lurie-L-theory}.
A Poincaré \(\infty\)-category consists of a stable \(\infty\)-category \(\C\) together with a functor \(\QF\colon \C\op \to \Spa\) which is \defi{quadratic} in the sense of Goodwillie calculus and satisfies a suitable unimodularity condition, the latter determining in particular a duality \(\Dual_{\QF}\colon \C\op \xrightarrow{\simeq}\C\) on \(\C\). We refer to such a \(\QF\) as a \defi{Poincaré structure} on \(\C\).
Roughly speaking, the role of the Poincaré structure \(\QF\) is to encode the flavour of forms that we want to consider. For example, for a commutative ring \(R\) one may take \(\C=\Dperf(R)\) to be the perfect derived category of \(R\). One should then think of the mapping spectrum \(\map_{\Dperf(R)}(X \otimes_R X, R)\) as the spectrum of bilinear forms on the chain complex \(X\), which acquires a natural \(\Ct\)-action by flipping the components in the domain term. In this case the Poincaré structure
\[
\QF^{\sym}_R(X) = \map_{\Dperf(R)}(X \otimes_R X, R)^{\hC}
\]
encodes a homotopy coherent version of the notion of symmetric bilinear forms, while
\[
\QF^{\qdr}_R(X) = \map_{\Dperf(R)}(X \otimes_R X, R)_{\hC}
\]
encodes a homotopy coherent version of quadratic forms. Both these Poincaré structures have the same underlying duality, given by \(X \mapsto \HHom_R(X,R)\).

Alternatively, as we develop in the present paper, one may also obtain Poincaré structures on \(\Dperf(R)\) by taking a \emph{non-abelian derived functor} associated to a quadratic functor \(\Prof{R}\op \to \Ab\) from finitely generated projective modules to abelian groups. For example, taking the functors which associate to a projective module \(P\) the abelian groups of quadratic, even and symmetric forms on \(P\) one obtains Poincaré structures \(\QF^{\gq}_R,\QF^{\gev}_R\) and \(\QF^{\gs}_R\) on \(\Dperf(R)\), respectively. We call these the \emph{genuine} quadratic, even and symmetric functors, and consider them as encoding the classical, rigid notions of hermitian forms in the present setting, whereas \(\QF^{\qdr}_R\) and \(\QF^{\sym}_R\) encode their homotopy coherent counterparts. More generally, one can apply this construction to any associative ring \(R\) equipped with an invertible \((R \otimes R)\)-module with involution \(M\) as above.
The resulting Poincaré structures are then all related by a sequence of natural transformations
\[
\QF^{\qdr}_M \Rightarrow \QF^{\gq}_M \Rightarrow \QF^{\gev}_M \Rightarrow \QF^{\gs}_M \Rightarrow \QF^{\sym}_M,
\]
which encode
the polarisation map between the quadratic, even and symmetric flavours of hermitian forms and at the same time %
the comparison between homotopy coherent and rigid variants of such forms. The fact that these two types of distinctions are not entirely unrelated leads to some of the more surprising applications of our approach. When \(2\) is invertible in \(R\), all these maps are equivalences.

The fundamental invariant of a Poincaré \(\infty\)-category is its space \(\Poinc(\C,\QF)\) of \defi{Poincaré objects}, which are pairs \((\x,\qone)\) consisting of an object \(\x \in \C\) and a point \(q \in \Om^{\infty}\QF(\x)\) whose associated map \(q_{\sharp}\colon \x \to \Dual_{\QF}(\x)\) is an equivalence. These are the avatars in the present context of the notion of a unimodular hermitian form. From this raw invariant one may produce two principal spectrum valued invariants - the \defi{Grothendieck-Witt spectrum} \(\GW(\C,\QF)\) and \defi{\(\L\)-theory spectrum \(\L(\C,\QF)\)}. The \(\L\)-theory spectrum was transported by Lurie from the classical work of Wall-Ranicki to the context of Poincaré \(\infty\)-categories in~\cite{Lurie-L-theory}. In particular, the \(\L\)-theory spectra \(\L^{\qdr}(R,M) := \L(\Dperf(R),\QF^{\qdr}_M)\) and \(\L^{\sym}(R,M) := \L(\Dperf(R),\QF^{\sym}_M)\) coincide with Ranicki's \(4\)-periodic quadratic and symmetric \(\L\)-theory spectra, respectively. When applied to the genuine Poincaré structures this yields new types of \(\L\)-theory spectra \(\L^{\gq}(R,M), \L^{\gev}(R,M)\) and \(\L^{\gs}(R,M)\). It turns out that these are in fact not entirely new: We show in \paperthree that for the genuine symmetric structure %
the homotopy groups of \(\L^{\gs}(R,M)\) coincide with Ranicki's original non-periodic variant of symmetric \(\L\)-groups, as defined in~\cite{RanickiATS1}. Somewhat surprisingly, the genuine quadratic \(\L\)-theory spectrum \(\L^{\gq}(R,M)\) is a \(4\)-fold shift of \(\L^{\gs}(R)\). %

The Grothendieck-Witt spectrum \(\GW(\C,\QF)\) of a Poincaré \(\infty\)-category is defined in \papertwo, though in the present paper we already introduce its zeroth homotopy group \(\GW_0(\C,\QF)\), namely, the Grothendieck-Witt group. The underlying infinite loop space
\[
\GWspace(\C,\QF) := \Om^{\infty}\GW(\C,\QF)
\]
is then called the Grothendieck-Witt space of \((\C,\QF)\). If \(2\) is invertible in \(R\), we show in \papertwo that \(\GW(R,M) := \GW(\Dperf(R),\QF_M)\) is equivalent to the Grothendieck-Witt spectrum defined by Schlichting in~\cite{schlichting-derived} (where \(\QF_M\) is any of the Poincaré structures considered above, which coincide due to the invertibility condition on \(2\)). When \(2\) is not invertible, the fourth and ninth author show in the companion paper~\cite{comparison} that the Grothendieck-Witt spaces of \(\Dperf(R)\) with respect to the genuine Poincaré structures \(\QF^{\gq}_M,\QF^{\gev}_M\) and \(\QF^{\gs}_M\) coincide with the classical Grothendieck-Witt spaces of quadratic, even and symmetric \(M\)-valued forms, respectively. On the other hand, the Grothendieck-Witt spectra of \((\Dperf(R),\QF^{\qdr}_M)\) and \((\Dperf(R),\QF^{\sym}_M)\) are actually new invariants of rings, which are based on the homotopy coherent avatars of quadratic and symmetric forms. These sometimes have better formal properties. For example, in the upcoming work~\cite{motives},
the first, third and seventh authors show that the \(\GW\)- and \(\L\)-theory spectra associated to the symmetric Poincaré structures \(\QF^{\sym}_R\) satisfy \(\mathbf{A}^{1}\)-invariance, and can further be encoded via motivic spectra over the integers. This statement does not hold for any of the other Poincaré structures above, including the genuine symmetric one.

One of the principal results we prove in \papertwo is that the relation between Grothendieck-Witt-, \(\L\)- and algebraic \(\K\)-theory is governed by the fundamental fibre sequence
\begin{equation}
\label{equation:tate}%
\K(\C)_{\hC} \to \GW(\C,\QF) \to \L(\C,\QF) ,
\end{equation}
where the first term is the homotopy orbits of the algebraic \(\K\)-theory spectra of \(\C\) with respect to the \(\Ct\)-action induced by the duality of \(\QF\). In the case of the genuine symmetric Poincaré structure \(\QF^{\gs}_M\), this gives a relation between classical symmetric Grothendieck-Witt groups and Ranicki's non-periodic symmetric \(\L\)-groups, which to our knowledge is completely new. In the case of the genuine quadratic structure the consequence is even more surprising: The resulting long exact sequence in homotopy groups extends the classical exact sequence~\eqrefone{equation:short} to a long exact sequence involving a shifted copy of Ranicki's non-periodic \(\L\)-groups.

The main role of the present instalment is to lay down the mathematical foundations that enable the arguments of the next three papers, and eventually their fruits, to take place. In particular, we carefully develop the main concepts of Poincaré \(\infty\)-categories and Poincaré objects, discuss hyperbolic objects and Lagrangians, and prove a version of Ranicki's algebraic Thom construction in the present setting. We also define the \(\L\)-groups and zeroth Grothendieck-Witt group of a Poincaré \(\infty\)-category, and conduct a thorough investigation of the global structural properties enjoyed by the \(\infty\)-category of Poincaré \(\infty\)-categories. In addition to the general framework, we also introduce and study important constructions of Poincaré \(\infty\)-categories,
which give rise to our motivating examples of interest. %
In particular:
\begin{enumerate}
\item
We classify all Poincaré structures in the case where \(\C\) is the \(\infty\)-category of perfect modules over a ring spectrum, and show that they can be efficiently encoded by the notion of a module with genuine involution.
\item
When \(\C\) is the perfect derived category of a discrete ring, we develop the procedure of \defi{deriving} Poincaré structures used to produce the genuine Poincaré structures above. Here, we pick up on some recent ideas of Brantner, Glasman and Illusie, and show that Poincaré structures on \(\C\) are in fact uniquely determined by their values on projective modules. This allows for the connection between the present set-up and Grothendieck-Witt theory of rings in \cite{comparison}, through which the applications of \papertwo and \paperthree to classical problems can be carried out.
\item
We develop in some detail the example of visible Poincaré structures on \(\infty\)-categories of parametrised spectra, which allows us to reproduce visible \(\L\)-theory as well as \(\LA\)-theory of Weiss-Williams in the present setting. This leads to applications in surgery theory, which we will pursue in future work.
\item
Following Lurie's treatment of \(\L\)-theory we study the process of tensoring and cotensoring a Poincaré \(\infty\)-category over a finite simplicial complex. This construction is later exploited in \papertwo to define and study the Grothendieck-Witt spectrum.
\item
We show that the \(\infty\)-category of Poincaré \(\infty\)-categories has all limits and colimits. This enables one, for example, to produce new Poincaré \(\infty\)-categories by taking fibres and cofibres of Poincaré functors, and enables the notion of \emph{additivity}, which lies at the heart of Grothendieck-Witt theory, to be properly set-up in \papertwo.
\item
We show that Poincaré \(\infty\)-categories can be tensored with each other. This can be used to produce new Poincaré \(\infty\)-categories from old, but also to identify additional important structures, such as a Poincaré symmetric monoidal structure, which arises in many examples of interest and entails the refinement of their Grothendieck-Witt and \(\L\)-theory spectra to \(\Einf\)-rings. This last claim is proven in \paperfour, though we construct the consequent multiplicative structure on the Grothendieck-Witt and \(\L\)-groups already in the present paper.
\end{enumerate}

\introsubsection{Applications}

Our framework of Poincaré \(\infty\)-categories is motivated by a series of applications which we extract in the following instalments, many of which pertain to classical questions in hermitian \(\K\)-theory. To give a brief overview of what's ahead, we first mention that a key feature of the Grothendieck-Witt spectrum we construct in \papertwo is its additivity. In the setting of Poincaré \(\infty\)-categories, this can be neatly phrased by saying that the functor \((\C,\QF) \mapsto \GW(\C,\QF)\) sends split bifibre sequences
\[
(\C,\QF) \to (\Ctwo,\QFtwo) \to (\Cthree,\QFthree)
\]
of Poincaré \(\infty\)-categories to bifibre sequences of spectra, where a split bifibre sequence is one in which \(\Ctwo \to \Cthree\) admits both a left and a right adjoint. One of the main results of \papertwo is that \(\GW\) is additive, and is, furthermore, universally characterised by this property as initial among additive functors from Poincaré \(\infty\)-categories to spectra equipped with a natural transformation from \(\Sig^{\infty}\Poinc\). This is analogous to the universal property characterising algebraic \(\K\)-theory of stable \(\infty\)-categories established in~\cite{BGT}. In fact, we show in \papertwo that \(\GW\) is not only additive but also \defi{Verdier localising}, a property formulated as above but with the splitness condition removed. This will be used in~\cite{motives} by the first, third and seventh author in order to show that the \(\GW\)-spectrum satisfies \emph{Nisnevich descent} over smooth schemes. It also plays a key role in the study of Grothendieck-Witt theory of Dedekind rings in \paperthree.

One major consequence of additivity is that the hyperbolic and forgetful maps fit to form %
the \emph{Bott-Genauer} sequence
\[
\GW(\C,\QF\qshift{-1}) \xrightarrow{\fgt} \K(\C) \xrightarrow{\hyp} \GW(\C,\QF),
\]
where \(\QF\qshift{n} = \Sig^n\QF\) is the shifting operation on Poincaré functors. Such a sequence was established in the setting of rings in which \(2\) is invertible by Schlichting~\cite{schlichting-derived}, who used it to produce another proof of Karoubi's fundamental theorem. The same argument then yields a version of Karoubi's fundamental theorem in the setting of Poincaré \(\infty\)-categories. When applied to the genuine Poincaré structures we construct in the present paper, this yields an extension of Karoubi's fundamental theorem to rings in which \(2\) is not assumed invertible, establishing, in particular, a conjecture of Karoubi and Giffen. %

The fundamental fibre sequence~\eqrefone{equation:tate} is heavily exploited in \paperthree to obtain applications for classical Grothendieck-Witt groups of rings. In particular, improving a comparison bound of Ranicki we show in \paperthree that, if \(R\) is Noetherian of global dimension \(d\), the maps
\[
\L^{\gq}(R,M) \longrightarrow \L^{\gev}(R,M) \longrightarrow \L^{\gs}(R,M) \longrightarrow \L^{\sym}(R,M)
\]
are equivalences in degrees past \(d+2,d\) and \(d-2\), respectively. Thus, even though the genuine \(\L\)-theory spectra are not \(4\)-periodic, they become so in degrees sufficiently large compared to the global dimension. In addition, when combined with the fundamental fibre sequence~\eqrefone{equation:tate} this implies that the maps of classical Grothendieck-Witt spaces
\[
\GWspace^{\qdr}_{\cl}(R,M) \rightarrow \GWspace^{\ev}_{\cl}(R,M) \rightarrow \GWspace^{\sym}_{\cl}(R,M)
\]
are isomorphisms on homotopy groups in sufficiently high degrees. %
This is a new and quite unexpected result about classical Grothendieck-Witt groups and, to the best of our knowledge, it is the first time that the global dimension of a ring has been related in any way to the gap between its quadratic and symmetric \(\GW\)-groups.
Combined with our extension of Karoubi's fundamental theorem this implies that in the case of finite global dimension Karoubi's fundamental theorem holds in its classical form in sufficiently high degrees, allowing for many of the associated arguments to be picked up in this context. In a different direction, for such rings one can eventually deduce results about classical symmetric \(\GW\)-groups from results on the corresponding homotopy coherent symmetric \(\GW\)-groups, allowing one to exploit some of the useful properties of the latter,  %
such as a dévissage property we prove in \paperthree and the \(\mathbf{A}^1\)-invariance, which will be established in~\cite{motives}, for the benefit of the former. We exploit these ideas in \paperthree to solve the homotopy limit problem for number rings, show that their Grothendieck-Witt groups are finitely generated, and produce an essentially complete calculation of the quadratic and symmetric Grothendieck-Witt groups (and their skew variants) of the integers, affirming, in particular, a conjecture of Berrick and Karoubi from \cite{berrick-karoubi}.

\introsubsection{Organisation of the paper}

Let us now describe the structure and the content of the present paper in more detail. In~\S\refone{section:poincare-cats} we define \defi{Poincaré \(\infty\)-categories}. As indicated before, a Poincaré \(\infty\)-category is a stable \(\infty\)-category \(\C\) equipped with a \emph{quadratic} functor \(\QF\colon \C\op \to \Spa\) which is \emph{perfect} in a suitable sense. We give the precise definition in \S\refone{subsection:hermitian-and-poincare-cats}, after a discussion of quadratic functors in \S\refone{subsection:quadratic}. We also consider the weaker notion of a \defi{hermitian \(\infty\)-category}, obtained by removing the perfectness condition on \(\QF\), and explain how to extract from a Poincaré structure \(\QF\) a \defi{duality} \(\QF_{\QF}\colon \C\op \xrightarrow{\simeq} \C\). In \S\refone{subsection:classification} we describe how one can classify hermitian and Poincaré structures on a given stable \(\infty\)-category in terms of their \emph{linear} and \emph{bilinear} parts. Finally, in \S\refone{subsection:functoriality} we discuss the functorial dependence of hermitian structures on the underlying stable \(\infty\)-category, and relate it to the classification discussed in \S\refone{subsection:classification}.

In~\S\refone{section:poincare-objects} we define the notion of a \defi{Poincaré object} in a given Poincaré \(\infty\)-category \((\C,\QF)\). Such a Poincaré object consists of an object \(\x \in \C\) together with a map \(q\colon \SS \to \QF(\x)\), to be thought of as a form in \(X\), such that a certain induced map \(q_{\sharp}\colon \x \to \Dual_{\QF} \x\) is an equivalence.
The precise definition is given in \S\refone{subsection:hermitian-poincare-objects}. We then proceed to discuss hyperbolic Poincaré objects in \S\refone{subsection:hyp-and-sym-poincare-objects}, and in \S\refone{subsection:metabolic-and-L} the slightly more general notion of metabolic Poincaré objects, that is, Poincaré objects that admit a Lagrangian. We show how one can understand metabolic Poincaré objects via Poincaré objects in a certain Poincaré \(\infty\)-category \(\Met(\C,\QF)\) constructed from \((\C,\QF)\). The notion of metabolic Poincaré objects is the main input in the definition of the \defi{\(\L\)-groups} of a given Poincaré \(\infty\)-category, which we also given in this section. Finally, in \S\refone{subsection:GW-group} we define the Grothendieck-Witt group \(\GW_0(\C, \Qoppa)\) of a given Poincaré \(\infty\)-category and develop its basic properties.

In~\S\refone{section:modules} we study Poincaré structures on the \(\infty\)-category \(\Modp{A}\) of perfect modules over a ring spectrum \(A\). %
To this end, we introduce the notion of a module with involution in \S\refone{subsection:modules-with-involution} and show how it can be used to model bilinear functors on module \(\infty\)-categories. We then refine this notion \S\refone{subsection:genuine-modules} to a module with \defi{genuine} involution, that  allows us to encode not only bilinear functors but also hermitian and Poincaré structures. Then, in~\S\refone{subsection:restriction-induction} we discuss the basic operations of restriction and induction of modules with genuine involution along maps of ring spectra.

In~\S\refone{section:examples} we discuss several examples of interest of Poincaré \(\infty\)-categories in further detail.
We begin in \S\refone{subsection:universal} with the important example of the universal Poincaré \(\infty\)-category \((\Spaf,\QF^{\uni})\), which is characterised by the property that Poincaré functors out of it pick out Poincaré objects in the codomain.
In \S\refone{subsection:discrete-rings} we consider perfect derived \(\infty\)-categories of \emph{ordinary rings} and show how to translate the classical language of forms on projective modules into that of the present paper via the process of \emph{deriving} quadratic functors. In~\S\refone{subsection:visible} and \S\refone{subsection:parametrised-spectra} we explain how to construct Poincaré structures producing visible \(\L\)-theory as studied by Weiss~\cite{Weiss}, Ranicki~\cite{Ranickiblue}, and more recently Weiss-Williams~\cite{WWIII}.

In~\S\refone{section:multiplicative} we show that the tensor product of stable \(\infty\)-categories refines to give a symmetric monoidal structure on the \(\infty\)-category \(\Catp\) of Poincaré \(\infty\)-categories. %
The precise definition and main properties of this monoidal product are spelt out in~\S\refone{subsection:tensor-product} and
\S\refone{subsection:monoidal-structure}.
In \S\refone{subsection:symmetric-monoidal-poincare} we analyse what it means for a Poincaré \(\infty\)-category to be an algebra with respect to this structure, and use this analysis in \S\refone{subsection:examples-monoidal} in order to identify various examples of interest of symmetric monoidal Poincaré \(\infty\)-categories.

In~\S\refone{section:cat-of-cats} we study the global structural properties of the \(\infty\)-categories \(\Catp\) and \(\Cath\) of Poincaré and hermitian \(\infty\)-categories, respectively. We begin in \S\refone{subsection:limits} by showing that these two \(\infty\)-categories have all small limits and colimits, and describe how these can be computed. In \S\refone{subsection:internal} we prove that the symmetric monoidal structures on \(\Catp\) and \(\Cath\) constructed in \S\refone{subsection:monoidal-structure} are \emph{closed}, that is, admit internal mapping objects. We then show in \S\refone{subsection:cotensoring} and \S\refone{subsection:tensoring} that \(\Cath\) is tensored and cotensored over \(\Cat\). A special role is played by indexing diagrams coming from the poset of faces of a finite simplicial complex, which we study in \S\refone{subsection:finite-tensors-cotensors} and \S\refone{subsection:finite-complexes}, showing in particular that in this case this procedure preserves Poincaré \(\infty\)-categories. The cotensor construction is used in \papertwo to define the hermitian \(\Q\)-construction and eventually Grothendieck-Witt theory, while the tensor construction plays a role in proving the universal property of Grothendieck-Witt theory.

In~\S\refone{section:metabolic} we consider the relationship between \(\Catp\) and \(\Cath\), and between both of them and various coarser variants, such as bilinear and symmetric bilinear \(\infty\)-categories. By categorifying the relationship between Poincaré forms, hermitian forms and bilinear forms we construct in \S\refone{subsection:bilinear-and-pairings} and \S\refone{subsection:thom} left and right adjoints to all relevant forgetful functors. In \S\refone{subsection:thom} we also prove a generalised version of the algebraic Thom construction, which is used in \papertwo for the formation of algebraic surgery. In \S\refone{subsection:mackey-functors} we use this to study \(\Catp\) and \(\Cath\) from the perspective of \(\Ct\)-category theory as developed by Barwick and collaborators, and set-up some of the foundations leading to the \emph{genuine \(\Ct\)-refinement} of the Grothendieck-Witt spectrum we construct in \papertwo. Finally, in \S\refone{subsection:GW-L-multiplicative} we show that the Grothendieck-Witt group and the \(\L\)-groups are lax symmetric monoidal functors with respect to the tensor product of Poincaré \(\infty\)-categories.

%% file: Acknowledgements.tex
For useful discussions about our project,
we heartily thank
Tobias Barthel,
Clark Barwick,
Lukas Brantner,
Mauricio Bustamante,
Denis-Charles Cisinski,
Dustin Clausen,
Uriya First,
Rune Haugseng,
André Henriques,
Lars Hesselholt,
Gijs Heuts,
Geoffroy Horel,
Marc Hoyois,
Max Karoubi,
Daniel Kasprowski,
Ben Knudsen,
Manuel Krannich,
Achim Krause,
Henning Krause,
Sander Kupers,
Wolfgang Lück,
Ib Madsen,
Cary Malkiewich,
Mike Mandell,
Akhil Mathew,
Lennart Meier,
Irakli Patchkoria,
Nathan Perlmutter,
Andrew Ranicki,
Oscar Randal-Williams,
George Raptis,
Marco Schlichting,
Peter Scholze,
Stefan Schwede,
Graeme Segal,
Markus Spitzweck,
Jan Steinebrunner,
Georg Tamme,
Ulrike Tillmann,
Maria Yakerson,
Michael Weiss,
and
Christoph Winges.

Furthermore, we owe a tremendous intellectual debt to Jacob Lurie for creating the framework we exploit here, and to Søren Galatius for originally spotting the overlap between various separate projects of ours; his insight ultimately led to the present collaboration.

%% file: Support.tex
The authors would also like to thank the Hausdorff Center for Mathematics at the University of Bonn, the Newton Institute at the University of Cambridge, the University of Copenhagen and the Mathematical Research Institute Oberwolfach for hospitality and support while parts of this project were undertaken.
\medskip

BC was supported by the French National Centre for Scientific Research (CNRS) through a ``délégation'' at LAGA, University Paris 13. 
ED was supported by the German Research Foundation (DFG) through the priority program ``Homotopy theory and Algebraic Geometry'' (DFG grant no.\ SPP 1786) at the University of Bonn and WS by the priority program ``Geometry at Infinity'' (DFG grant no.\ SPP 2026) at the University of Augsburg. 
YH and DN were supported by the French National Research Agency (ANR) through the grant ``Chromatic Homotopy and K-theory'' (ANR grant no.\ 16-CE40-0003) at LAGA, University of Paris 13.
FH is a member of the Hausdorff Center for Mathematics at the University of Bonn (DFG grant no.\ EXC 2047 390685813) and TN of the cluster ``Mathematics Münster: Dynamics-Geometry-Structure'' at the University of Münster (DFG grant no.\ EXC 2044 390685587). 
FH, TN and WS were further supported by the Engineering and Physical Sciences Research Council (EPSRC) through the program ``Homotopy harnessing higher structures'' at the Isaac Newton Institute for Mathematical Sciences (EPSRC grants no.\ EP/K032208/1 and EP/R014604/1). 
FH was also supported by the European Research Council (ERC) through the grant ``Moduli spaces, Manifolds and Arithmetic'' (ERC grant no.\ 682922) and KM by the grant ``$\K$-theory, $\L^2$-invariants, manifolds, groups and their interactions'' (ERC grant no.\ 662400). 
ML and DN were supported by the collaborative research centre ``Higher Invariants'' (DFG grant no.\ SFB 1085) at the University of Regensburg. 
ML was further supported by the research fellowship ``New methods in algebraic K-theory'' (DFG grant no.\ 424239956) and by the Danish National Research Foundation (DNRF) through the Center for Symmetry and Deformation (DNRF grant no.\ 92) and the Copenhagen Centre for Geometry and Topology (DNRF grant no.\ 151) at the University of Copenhagen. 
KM was also supported by the K\&A Wallenberg Foundation.

%% file: PoincareCats.tex
In this section we introduce the principal notion of this paper, namely that of \defi{Poincaré \(\infty\)-categories}. These were first defined by Lurie in \cite{Lurie-L-theory}, though no name was chosen there. Succinctly stated, Poincaré \(\infty\)-categories are stable \(\infty\)-categories \(\C\) equipped with a \emph{quadratic} functor \(\QF\colon \C\op \to \Spa\) to spectra, which is \emph{perfect} in a sense we explain below. We then refer to \(\QF\) as a \defi{Poincaré structure} on \(\C\). It is also convenient to consider the more general setting where \(\QF\) is not necessarily perfect, leading to a notion that we call a \emph{hermitian \(\infty\)-category}. We present both of these in \S\refone{subsection:hermitian-and-poincare-cats}, after devoting \S\refone{subsection:quadratic} to surveying quadratic functors and their basic properties. In \S\refone{subsection:classification} we describe how one can classify hermitian and Poincaré structures on a given stable \(\infty\)-category in terms of their \emph{linear} and \emph{bilinear} parts. This is a particular case of the general structure theory of Goodwillie calculus, but we elaborate the details relevant to the case at hand, as we rely on this classification very frequently, both in explicit constructions of examples and in general arguments.
Finally, in \S\refone{subsection:functoriality} we discuss the functorial dependence of hermitian structures on the underlying stable \(\infty\)-category, and relate it to the classification discussed in \S\refone{subsection:classification}.

\subsection{Quadratic and bilinear functors}
\label{subsection:quadratic}%

In this subsection we will recall the notions of quadratic and bilinear functors, and survey their basic properties. These notions fit most naturally in the context of \emph{Goodwillie calculus}, as adapted to the \(\infty\)-categorical setting in~\cite[\S 6]{HA}. Our scope of interest here specializes that of loc.\ cit.\ in two ways: first, we will only consider the Goodwillie calculus up to degree \(2\), and second, we will focus our attention on functors from a stable \(\infty\)-category \(\C\) to the stable \(\infty\)-category \(\Spa\) of spectra. This highly simplifies the general theory, and will allow us to give direct arguments for most claims, instead of quoting~\cite[\S 6]{HA}. The reader should however keep in mind that the discussion below is simply a particular case of Goodwillie calculus, to which we make no claim of originality.

Recall that an \(\infty\)-category \(\C\) is said to be \defi{pointed} if it admits an object which is both initial and terminal. Such objects are then called \defi{zero objects}. A functor \(f\colon \C \to \D\) between two pointed \(\infty\)-categories is called \defi{reduced} if it preserves zero objects. Given two pointed \(\infty\)-categories \(\C,\D\) we will denote by \(\Fun_{\ast}(\C,\D) \subseteq \Fun(\C,\D)\) the full subcategory spanned by the reduced functors. A \defi{stable} \(\infty\)-category is by definition a pointed \(\infty\)-category which admits pushouts and pullbacks and in which a square is a pushout square if and only if it is a pullback square. To avoid breaking the symmetry one then refers to such squares as \defi{exact}. Any stable \(\infty\)-category is canonically enriched in spectra, and we will denote by \(\map_{\C}(x,y)\) the mapping \emph{spectrum} from \(x\) to \(y\). It is related to the corresponding mapping \emph{space} in \(\C\) by the formula \(\Map_{\C}(x,y) = \Om^{\infty}\map(x,y)\).

A functor \(f\colon \C \to \D\) between two stable \(\infty\)-categories is called \defi{exact} if it preserves zero objects and exact squares. We note that stable \(\infty\)-categories automatically admit all finite limits and colimits, and that a functor between stable \(\infty\)-categories is exact if and only if it preserves finite colimits, and if and only if it preserves finite limits. If \(\D\) is a stable \(\infty\)-category and \(\C \subseteq \D\) is a full subcategory which is closed under finite limits and finite colimits then \(\C\) is also stable and the inclusion \(\C \subseteq \D\) is an exact functor. In this case we will say that \(\C\) is a \defi{stable subcategory} of \(\D\).
Given two stable \(\infty\)-categories \(\C,\D\) with \(\C\) small we will denote by \(\Funx(\C,\D) \subseteq \Fun(\C,\D)\) the full subcategory spanned by the exact functors. We note that when \(\C\) and \(\D\) are stable one has that \(\Fun(\C,\D)\) is also stable and \(\Fun_\ast(\C,\D)\) and \(\Funx(\C,\D)\) are stable subcategories. We will denote by \(\Catx\) the (non-full) subcategory of \(\Cat\) spanned by the stable \(\infty\)-categories and exact functors between them.

If one considers stable \(\infty\)-categories as a categorified version of a vector space, then reduced functors correspond to zero-preserving maps, while exact functors correspond to \emph{linear maps}.
If a functor \(f\colon \C \to \D\) is only required to preserves exact squares, but is not necessarily reduced, then one says that \(f\) is \defi{\(1\)-excisive}. More generally, if \(\C\) is an \(\infty\)-category with finite colimits and \(\D\) and \(\infty\)-category with finite limits, then \(f\colon \C \to \D\) is said to be \(1\)-excisive if it sends pushout squares to pullback squares. In the above analogy with linear algebra, these correspond to affine maps, that is, maps which contain a linear part and a constant term, or said differently: polynomial maps of degree 1. In the theory of \emph{Goodwillie calculus} this point of view is generalized to higher degrees as follows:

\begin{definition}
\label{definition:2-excisive}%
A \(3\)-cube \(\rho\colon(\Del^1)^3 \to \C\) is said to be \defi{cartesian} if it exhibits \(\rho(0,0,0)\) as the limit of the restriction of \(\rho\) to the subsimplicial set \((\Del^1)^3\) spanned by the complement of \((0,0,0)\). Such a \(3\)-cube \(\rho\) is called \defi{strongly cartesian} if its restriction to each \(2\)-dimensional face of \((\Del^1)^3\) is a cartesian square. In particular, strongly cartesian \(3\)-cubes are cartesian. Dually, \(\rho\) is said to be \defi{(strongly) cocartesian} if \(\rho\op\) is a (strongly) cartesian cube in \(\C\op\). A functor \(f\colon \C \to \D\) whose domain admits finite colimits and whose target admits finite limits is called \defi{\(2\)-excisive} if it sends strongly cocartesian \(3\)-cubes to cartesian \(3\)-cubes.

If \(\C\) is stable then a \(3\)-cube is (strongly) cartesian if and only if it is (strongly) cocartesian, in which case we simply say that \(\rho\) is \defi{(strongly) exact}. A functor \(f\colon \C \to \D\) between stable \(\infty\)-categories is then \(2\)-excisive if it sends strongly exact \(3\)-cubes to exact \(3\)-cubes.
\end{definition}

\begin{remark}
Though in the present paper we will focus almost entirely on the case of stable \(\infty\)-categories, we chose to formulate the above definition in the slightly more general setting where \(f\colon \C\to \D\) is a functor from an \(\infty\)-category with finite colimits to an \(\infty\)-category with finite limits. This level of generality, in which most of Goodwillie calculus can be carried out, will be used in \S\refone{subsection:discrete-rings}, but will otherwise not be needed in the present paper.
\end{remark}

We note that every \(1\)-excisive functor is in particular \(2\)-excisive. If the former are analogous to affine maps between vector spaces, the latter are then analogous to maps between vector spaces which are polynomial of degree \(2\), that is, contain a homogeneous quadratic part, a linear part, and a constant term. If we restrict attention to \(2\)-excisive functors which are reduced, then we get the analogue of maps with terms in degrees \(1\) and \(2\), but no constant term. These are going to be the functors we consider in this paper.

In the present work it will be convenient to take a slightly different route to the definition of reduced \(2\)-excisive functors, which proceeds as follows. Given a small stable \(\infty\)-category \(\C\), let us denote by \(\BiFun(\C) \subseteq \Fun_{\ast}(\C\op \times\C\op,\Spa)\) the full subcategory spanned by those reduced functors \(\Bil\colon\C\op \times \C\op \to \Spa\) such that \(\Bil(\x,\y) \simeq 0\) if either \(\x\) or \(\y\) is a zero object. Such functors may be referred to as \defi{bi-reduced}. Then \(\BiFun(\C)\) is closed under all limits and colimits in \(\Fun_{\ast}(\C\op \times \C\op,\Spa)\), and hence the inclusion of the former in the latter admits both a left and a right adjoint. These left and right adjoints are in fact canonically equivalent, and can be described by the following explicit formula: given a reduced functor \(\Bil\colon \C\op \times \C\op \to \Spa\) we have a canonically associated retract diagram
\begin{equation}
\label{equation:reducification}%
\Bil(\x,0)\oplus \Bil(0,\y) \to \Bil(\x,\y) \to \Bil(\x,0)\oplus \Bil(0,\y),
\end{equation}
where \(0 \in \C\) is a chosen zero object, and all the maps are induced by the essentially unique maps \(0 \to \x \to 0\) and \(0 \to \y \to 0\). The composition of these two maps is the identity thanks to the assumption that \(\Bil\) is reduced, that is, \(\Bil(0,0) \simeq 0\).  The above retract diagram then induces a canonical splitting
\[
\Bil(\x,\y) \simeq \Bil^{\red}(\x,\y) \oplus \Bil(\x,0) \oplus \Bil(0,\y) ,
\]
where \(\Bil^{\red}(\x,\y)\) can be identified with both the cofibre of the left map in~\eqrefone{equation:reducification} and the fibre of the right map in~\eqrefone{equation:reducification}. We note that by construction the resulting functor
\[
\Bil^{\red}(-,-)\colon \C\op \times \C\op \to \Spa
\]
is bi-reduced. The following lemma records the fact that the association \(\Bil \mapsto \Bil^{\red}\) yields both a left and a right adjoint to the inclusion \(\BiFun(\C) \subseteq \Fun_{\ast}(\C\op\times\C\op,\Spa)\).

\begin{lemma}
\label{lemma:reducification}%
The split inclusion \(\Bil^{\red}(-,-) \Rightarrow \Bil(-,-)\) is universal among natural transformations to \(\Bil\) from a bi-reduced functor, while the projection \(\Bil(-,-) \Rightarrow \Bil^{\red}(-,-)\) is universal among natural transformations from \(\Bil\) to a bi-reduced functor. In particular, the association \(\Bil \mapsto \Bil^{\red}\) is both left and right adjoint to the full inclusion \(\BiFun(\C) \subseteq \Fun_{\ast}(\C\op\times\C\op,\Spa)\).
\end{lemma}
\begin{proof}
Given that \(\Fun_\ast(\C\op\times\C\op,\Spa)\) is stable and \(\BiFun(\C)\) is a stable full subcategory, to prove both claims it suffices to show that for \(\Bil \in \Fun_\ast(\C\op\times\C\op,\Spa)\), the associated functors
\[
(x,y) \mapsto \Bil(x,0) \qquad {and} \qquad (x,y) \mapsto \Bil(0,y)
\]
considered as functors in \(\Fun_\ast(\C\op\times\C\op,\Spa)\) have a trivial mapping spectrum to any and from any bi-reduced functor. Indeed, since \(0 \in \C\op\) is both final and initial it follows that the inclusion \(\C\op \times \{0\} \subseteq \C\op \times \C\op\) is both left and right adjoint to the projection \(\C\op \times \C\op \to \C\op \times \{0\}\), and hence restricting along this inclusion is both left and right adjoint to restricting along this projection. The same statement holds for the inclusion \(\{0\} \times \C\op \subseteq \C\op \times \C\op\) of the second factor. The mapping spectrum between any bi-reduced functor and a functor restricted along either projection is consequently trivial.
\end{proof}

\begin{definition}
\label{definition:cross-effect}%
Let \(\C\) be a stable \(\infty\)-category and \(\QF \colon \C\op \rightarrow \Spa\) a reduced functor. We will denote by \(\Bil_\QF \in \BiFun(\C)\) the bi-reduced functor
\[
\Bil_{\QF}(-,-) := \QF((-) \oplus (-))^{\red} \colon \C\op \times \C\op \to \Spa
\]
obtained by taking the universal bi-reduced replacement described above of the reduced 2-variable functor \((\x,\y) \mapsto \QF(\x \oplus \y)\).
Following the terminology of Goodwillie calculus
we will refer to \(\Bil_{\QF}(-,-)\) as the \defi{cross effect} of \(\QF\). The formation of cross effects then yields a functor
\begin{equation}
\label{equation:crs}%
\Bil_{(-)}\colon \Fun_{\ast}(\C\op,\Spa) \to \BiFun(\C)
\end{equation}
sending \(\QF\) to \(\Bil_{\QF}\).
\end{definition}

\begin{remark}
In~\cite{Lurie-L-theory} the term \defi{polarization} is used for what we called above cross effect, though in~\cite[\S 6]{HA} the term cross effect is employed.
\end{remark}

\begin{remark}
\label{remark:invariance-base-change}%
If \(f,g\colon \C \to \D\) are reduced functors then the associated restriction functor
\[
(f \times g)^*\colon \Fun_\ast(\D\op \times \D\op,\Spa) \to \Fun_{\ast}(\C\op\times\C\op,\Spa)
\]
along \((f \times g)\op\colon \C\op \times \C\op \to \D\op \times \D\op\) sends the retract diagram
\[
\Bil(\x,0)\oplus \Bil(0,\y) \to \Bil(\x,\y) \to \Bil(\x,0)\oplus \Bil(0,\y)
\]
to the retract diagram
\[
\Bil(f(\x),0)\oplus \Bil(0,g(\y)) \to \Bil(f(\x),g(\y)) \to \Bil(f(\x),0)\oplus \Bil(0,g(\y)),
\]
where we have used the symbols \(\x\) and \(\y\) to distinguish the two entries.
It then follows that the universal bi-reduction procedure described above commutes with restriction (along pairs of reduced functors).
Similarly, if \(f\colon \C \to \D\) furthermore preserves direct sums, then the formation of cross effects is compatible with restriction along \(f\), that is, the square
\[
\begin{tikzcd}
\Fun_{\ast}(\D\op,\Spa) \ar[r,"{f^*}"] \ar[d,"\Bil_{(-)}"'] & \Fun_{\ast}(\C\op,\Spa) \ar[d,"\Bil_{(-)}"] \\
\BiFun(\D) \ar[r,"{(f \times f)^*}"] & \BiFun(\C)
\end{tikzcd}
\]
naturally commutes.
\end{remark}

Given a stable \(\infty\)-category \(\C\), the diagonal functor \(\Del\colon \C\op \to \C\op \times \C\op\) induces a pullback functor
\[
\Delta^* \colon \BiFun(\C) \rightarrow \Fun_{\ast}(\C\op,\Spa).
\]
In what follows, for any \(\Bil\colon \C\op \times \C\op \to \Spa\), we will denote by \(\Bil^{\diag} := \Del^*\Bil\) the restriction of \(\Bil\) along the diagonal. Now the maps \(\QF(\x \oplus \x) \to \QF(\x)\) and \(\QF(\x) \to \QF(\x \oplus \x)\) induced by the diagonal \(\Delta_{\x}\colon \x \to \x \oplus \x\) and collapse map \(\nabla_{\x}\colon \x \oplus \x \to \x\) induce natural maps
\begin{equation}
\label{equation:cross-effect}%
\Bil_{\QF}(\x,\x) \to \QF(\x) \to \Bil_{\QF}(\x,\x),
\end{equation}
which can be considered as natural transformations
\begin{equation}
\label{equation:cross-effect-2}%
\Bil_{\QF}^{\diag} \Rightarrow \QF \Rightarrow \Bil_{\QF}^{\diag}.
\end{equation}
The formation of cross effects then enjoys the following universal property:

\begin{lemma}
\label{lemma:universal-crs}%
The two natural transformations in~\eqrefone{equation:cross-effect-2} act as a unit and counit exhibiting the cross effect functor~\eqrefone{equation:crs} as left and right adjoint respectively to the restriction functor \(\Delta^* \colon \BiFun(\C) \rightarrow \Fun_{\ast}(\C\op,\Spa)\).
\end{lemma}
\begin{proof}
The direct sum functor \(\C\op \times \C\op \to \C\op\) realizes both the product and coproduct (since \(\C\op\) is stable) and is hence both left and right adjoint to \(\Del \colon \C\op \to \C\op \times \C\op\), with units and counits given by the diagonal and collapse maps of the objects in \(\C\). It then follows that restriction along the direct sum functor is both right and left adjoint to restriction along \(\Del\), with unit and counit induced by the diagonal and collapse maps. The desired result now follows from Lemma~\refone{lemma:reducification}.
\end{proof}

\begin{remark}
\label{remark:other-unit-and-counit}%
The two sided adjunction of Lemma~\refone{lemma:universal-crs} is obtained by composing a pair of two-sided adjunctions
\[
\Fun_{\ast}(\C\op,\Spa) \leftrightarrows \Fun_{\ast}(\C\op \times \C\op,\Spa) \leftrightarrows \BiFun(\C),
\]
where the one on the left is induced by the two sided adjunction \(\C\op \stackbin[\oplus]{\Del}{\leftrightarrows} \C\op \times \C\op\) witnessing the existence of biproducts in \(\C\op\), and the one of the right exhibits the full subcategory \(\BiFun(\C)\subseteq \Fun_{\ast}(\C\op\times\C\op,\Spa)\) as reflective and coreflective (Lemma~\refone{lemma:reducification}). In particular, we may express the unit and counit of the two sided adjunction
\(\Fun_{\ast}(\C\op,\Spa) \leftrightarrows \BiFun(\C)\)
which are not specified in Lemma~\refone{lemma:universal-crs} via the unit \((\x,\y) \to (\x \oplus \y,\x \oplus \y)\) of the adjunction
\(\C\op \times \C\op \adj \C\) and counit \((\x \oplus \y,\x \oplus \y) \to (\x,\y)\) of the adjunction \(\C\op \adj \C\op \times \C\op\), which are all given by the corresponding component inclusions and projections. Unwinding the definitions, we get that the unit of the adjunction \(\BiFun(\C)\adj \Fun_{\ast}(\C\op,\Spa)\) is given by the induced map
\[
\Bil(\x,\y) \to \fib[\Bil(\x\oplus\y,\x \oplus \y) \to \Bil(\x,\x) \oplus \Bil(\y,\y)]
\]
and the counit of the adjunction \(\Fun_{\ast}(\C\op,\Spa) \adj \BiFun(\C)\) is given by the induced map
\[
\cof[\Bil(\x,\x) \oplus \Bil(\y,\y) \to \Bil(\x\oplus\y,\x \oplus \y)] \to \Bil(\x,\y).
\]
\end{remark}

\begin{lemma}
\label{lemma:bilinear-symmetric}%
Let \(\QF \colon \C\op \rightarrow \Spa\) be a reduced functor. Then the cross effect \(\Bil_\QF\) is symmetric, i.e.\ it canonically refines to an element of \(\Fun(\C\op \times \C\op,\Spa)^\hC\), where the cyclic group with two elements \(\Ct\) acts by flipping the two input variables.
\end{lemma}
\begin{proof}
By~\cite[Proposition 6.1.4.3, Remark 6.1.4.4]{HA} the bi-reduction functor
\[
(-)^{\red}\colon \Fun_{\ast}(\C\op\times\C\op,\Spa) \to \BiFun(\C)
\]
discussed above refines to a compatible functor
\[
\Fun_{\ast}(\C\op\times\C\op,\Spa)^{\hC} \to \BiFun(\C)^{\hC}
\]
on \(\Ct\)-equivariant objects. It will hence suffice to show that the functor \((\x,\y) \mapsto \QF(\x\oplus \y)\) naturally refines to a \(\Ct\)-equivariant object. For this, it suffices to note that the direct sum functor \(\C\op \times \C\op \to \C\op\) is equipped with a \(\Ct\)-equivariant structure with respect to the flip action on \(\C\op \times \C\op\) and the trivial action on \(\C\op\). Indeed, this is part of the symmetric monoidal structure afforded to the direct sum, canonically determined by its universal description as the coproduct in \(\C\op\).
\end{proof}

Keeping in mind the proofs of Lemma~\refone{lemma:universal-crs} and Lemma~\refone{lemma:bilinear-symmetric}, we now note that the diagonal functor \(\Del\colon \C\op \to \C\op \times \C\op\), which is both left and right adjoint to the \(\Ct\)-equivariant direct sum functor, is also canonically invariant under the \(\Ct\)-action on the right hand side switching the two components.
This means that the associated restriction functor
\[
\Delta^* \colon \BiFun(\C) \rightarrow \Fun_*(\C\op,\Spa)
\]
is equivariant for the trivial \(\Ct\)-action on the target, and so the restricted functor \(\Bil^{\diag}_{\QF} = \Delta^*\Bil_\QF\) becomes a \(\Ct\)-object of \(\Fun(\C\op,\Spa)\). In particular, \(\Bil_\QF(\x,\x)\) is naturally a spectrum with a \(\Ct\)-action for every \(\x \in \C\). Explicitly, this action is induced by the canonical action of \(\Ct\) on \(\x \oplus \x\) by swapping the components.

\begin{lemma}
\label{lemma:equivariance}%
The natural transformations in~\eqrefone{equation:cross-effect-2} both naturally refine to \(\Ct\)-equivariant maps with respect to the above \(\Ct\)-action on \(\Bil^{\Del}\).
In particular, the maps~\eqrefone{equation:cross-effect-2} induces natural transformations
\begin{equation}
\label{equation:orbits-fixed-points}%
[\Bil^{\Del}_{\QF}]_{\hC} \Rightarrow \QF \Rightarrow [\Bil_{\QF}^{\Del}]^{\hC}.
\end{equation}
\end{lemma}
\begin{proof}
Inspecting the construction of the natural transformations in~\eqrefone{equation:cross-effect-2} we see that it will suffice to put a \(\Ct\)-equivariant structure on the diagonal and collapse natural transformations
\[
\Del\colon \id \Rightarrow \id\oplus \id \quad\text{and}\quad \nabla \colon \id \oplus \id \Rightarrow \id
\]
of functors \(\C \to \C\). This in turn follows from the fact that the direct sum monoidal structure is both cartesian and cocartesian and every object is canonically a commutative algebra object with respect to coproducts (\cite[Proposition 2.4.3.8]{HA}).
\end{proof}

\begin{definition}
\label{definition:funb-funs}%
For \(\C,\D\) and \(\E\) stable \(\infty\)-categories, we will say that a functor \(b\colon \C \times \D \to \E\) is \defi{bilinear} if it is exact in each variable separately. For a stable \(\infty\)-category \(\C\) we will denote by \(\Funb(\C) \subseteq \Fun(\C\op \times \C\op,\Spa)\) the full subcategory spanned by the bilinear functors. We note that this full subcategory is closed under, and hence inherits, the flip action of \(\Ct\).
We will then denote by \(\Funs(\C) := [\Funb(\C)]^{\hC}\) the \(\infty\)-category of \(\Ct\)-equivariant objects in \(\Funb(\C)\) with respect to the flip action in the entries,
and refer to them as \defi{symmetric bilinear functors} on \(\C\).
\end{definition}

\begin{example}
\label{example:symmetric-monoidal}%
Suppose that \(\C\) is a stable \(\infty\)-category equipped with a monoidal structure which is exact in each variable separately. Then for every object \(a \in \C\) we have an associated bilinear functor \(\Bil_a\colon \C\op \times \C\op \to \Spa\) defined by
\[
\Bil_a(\x,\y) := \map_{\C}(\x \otimes \y,a),
\]
where \(\map_{\C}(-,-)\) refers to the canonical enrichment of \(\C\) in spectra.
If the monoidal structure refines to a symmetric one then \(\Bil_a\) refines to a symmetric bilinear functor. Natural examples of interest to keep in mind are when \(\C\) is the perfect derived category of a commutative ring (or, more generally, an \(\Einf\)-ring spectrum), or the \(\infty\)-category of perfect quasi-coherent sheaves on a scheme.
\end{example}

\begin{proposition}
\label{proposition:basic-properties-quad-functors}%
Let \(\QF \colon \C\op \rightarrow \Spa\) be a reduced functor. Then the following are equivalent:
\begin{enumerate}
\item
\label{item:goodwillie}%
\(\QF\) is \(2\)-excisive;
\item
\label{item:fibre-exact}%
the cross effect \(\Bil_\QF\) is bilinear and the fibre of the natural transformation \(\QF(\x) \to \Bil_{\QF}(\x,\x)^{\hC}\)
from~\eqrefone{equation:orbits-fixed-points}
is an exact functor in \(\x\);
\item
\label{item:cofibre-exact}%
the cross effect \(\Bil_\QF\) is bilinear and the cofibre of the natural transformation \(\Bil_{\QF}(\x,\x)_{\hC} \to \QF(\x)\)
from~\eqrefone{equation:orbits-fixed-points}
is an exact functor in \(\x\).
\end{enumerate}
\end{proposition}

\begin{proof}
Since \(\Spa\) is stable the property of being reduced and \(2\)-excisive is preserved under limits and colimits of functors \(\C\op \to \Spa\).
It then follows that both \refoneitem{item:fibre-exact} and \refoneitem{item:cofibre-exact} imply \refoneitem{item:goodwillie}, since exact functors and diagonal restrictions of bilinear functors are in particular reduced and \(2\)-excisive (see \cite[Cor.~6.1.3.5]{HA}).

In the other direction, if \(\QF\) is \(2\)-excisive then its cross effect is bilinear by \cite[Pr.~6.1.3.22]{HA}. Moreover, since taking the cross effect commutes with fibres and cofibres, the functors in the statement of \refoneitem{item:fibre-exact} and \refoneitem{item:cofibre-exact} have trivial cross effect. But they are also reduced and \(2\)-excisive by the first part of the argument, and are hence
exact by \cite[Pr.~6.1.4.10]{HA}.
\end{proof}

\begin{definition}
\label{definition:herm-structure}%
We will say that \(\QF\colon \C\op \to \Spa\) is \defi{quadratic} if it satisfies the equivalent conditions of Proposition~\refone{proposition:basic-properties-quad-functors}. For a small stable \(\infty\)-category \(\C\) we will then denote by \(\Funq(\C)\subseteq \Fun(\C\op,\Spa)\) the full subcategory spanned by the quadratic functors.
\end{definition}

\begin{remark}
\label{remark:closed}%
It follows from the first criterion in Proposition~\refone{proposition:basic-properties-quad-functors} that \(\Funq(\C)\) is closed under limits and colimits in \(\Fun(\C\op,\Spa)\). Since the latter is stable it follows that \(\Funq(\C)\) is stable as well.
\end{remark}

In light of Lemma~\refone{lemma:bilinear-symmetric} and Proposition~\refone{proposition:basic-properties-quad-functors}, the cross effect functor refines to a functor
\[
\Bil_{(-)} \colon \Funq(\C) \longrightarrow \Funs(\C).
\]
We will then refer to \(\Bil_{\QF} \in \Funs(\C)\) as the \defi{symmetric bilinear part} of \(\QF \in \Funq(\C)\), and refer to the underlying bilinear functor of \(\Bil_{\QF}\) as the \defi{bilinear part} of \(\QF\).

\begin{examples}
\label{example:first-example}%
\
\begin{enumerate}
\item
\label{item:exact-is-quadratic}%
Any exact functor \(\C\op \to \Spa\) is quadratic. These are exactly the quadratic functors whose bilinear part vanishes. In particular, we have an exact full inclusion of stable \(\infty\)-categories \(\Funx(\C\op,\Spa) \subseteq \Funq(\C)\).
\item
\label{item:diag-of-bil-is-quadratic}%
If \(\Bil\colon \C\op \times \C\op \to \Spa\) is a bilinear functor then the functor \(\Bil^{\diag}(\x) = \Bil(\x,\x)\) is a quadratic functor (\cite[Cor.~6.1.3.5]{HA}). Its symmetric bilinear part is given by the symmetrization \((\x,\y) \mapsto \Bil(\x,\y) \oplus \Bil(\y,\x)\) of \(\Bil\), equipped with its canonical symmetric structure.
\end{enumerate}
\end{examples}

\begin{example}
\label{example:quadratic-symmetric}%
If \(\Bil \in \Funs(\C)\) is a symmetric bilinear functor then the functors
\[
\QF^{\qdr}_{\Bil}(\x) := \Bil^{\diag}_{\hC}(\x) = \Bil(\x,\x)_{\hC}
\]
and
\[
\QF^{\sym}_{\Bil}(\x) := (\Bil^{\diag})^{\hC}(\x) = \Bil(\x,\x)^{\hC}
\]
are both quadratic functors. Indeed, this follows from the previous example by noting that the symmetry induces a \(\Ct\)-action on \(\Bil^{\diag}\) and invoking Remark~\refone{remark:closed}. Since taking cross-effects commutes with all limits and colimits the symmetric bilinear parts of these functors are given respectively by
\[
[\Bil(\x,\y) \oplus \Bil(\y,\x)]_{\hC} \quad\text{and}\quad [\Bil(\x,\y) \oplus \Bil(\y,\x)]^{\hC},
\]
which are both canonically equivalent to \(\Bil\) itself: indeed, when \(\Bil\) is symmetric its symmetrization canonically identifies with \(\Bil[\Ct] = \Bil \oplus \Bil\) as a \(\Ct\)-object in \(\Funs(\C)\), which, since the latter is stable, is the \(\Ct\)-object both induced and coinduced from \(\Bil\).

The superscript \((-)^{\qdr}\) and \((-)^{\sym}\) above refer to the relation between these constructions and the notions of \emph{quadratic} and \emph{symmetric} forms in algebra. To see this, consider the case where \(\C := \Dperf(R)\) is the perfect derived \(\infty\)-category of a commutative ring \(R\), that is, the \(\infty\)-categorical localisation of the category bounded complexes of finitely generated projective \(R\)-modules by quasi-isomorphisms.
We then have a natural choice of a bilinear functor \(\Bil_R\colon \C\op \times \C\op \to \Spa\) given by
\[
\Bil_R(X, Y) = \map_R(X \otimes_R Y, R)
\]
where \(\otimes_R\) denotes the (derived) tensor product over \(R\).
A point \(\beta \in \Om^{\infty}\Bil(X,Y)\) then corresponds to a map \(X \otimes_R Y \to R\), which we can consider as a bilinear form on the pair \((X,Y)\). If \(X,Y\) are ordinary projective modules then \(\pi_0\Bil_R(X,Y)\) is simply the abelian group of bilinear forms on \((X,Y)\) in the ordinary sense. For a projective \(R\)-module \(X\) we may then identify the \(\Ct\)-fixed subgroup \(\pi_0\Bil_R(X,X)^{\Ct}\) with the group of symmetric bilinear forms on \(X\), while the \(\Ct\)-quotient group \(\pi_0\Bil_R(X,X)_{\Ct}\) can be identified with the group of quadratic forms on \(X\) via the map sending the orbit of bilinear form \(b\colon X \otimes_R X \to R\) to the quadratic form \(q_b(x) = b(x,x)\). In this case the quadratic functors \(\QF^{\qdr}_R := \QF^{\qdr}_{\Bil_R}\) and \(\QF^{\sym}_R := \QF^{\sym}_{\Bil_R}\) defined as above can be considered as associating to a perfect \(R\)-complex \(X\) a suitable \emph{spectrum} of quadratic and symmetric forms on \(X\), respectively.
\end{example}

\begin{remark}
\label{remark:adj-bilinear-part-diag}%
By definition the cross effect of a quadratic functor is bilinear, and on the other hand by Example~\refone{example:first-example}\refoneitem{item:diag-of-bil-is-quadratic} the diagonal restriction of any bilinear functor is quadratic. It then follows from Lemma~\refone{lemma:universal-crs} that diagonal restriction \(\Del^*\colon \Funb(\C) \to \Funq(\C)\) determines a two-sided adjoint to the bilinear part functor \(\Bil_{(-)}\colon \Funq(\C) \to \Funb(\C)\), with unit and counit given by the natural maps
\[
\Bil_{\QF}(\x,\x) \Rightarrow \QF(\x) \Rightarrow \Bil_{\QF}(\x,\x).
\]
By Remark~\refone{remark:other-unit-and-counit} the other unit and counit are given by the component inclusion and projections
\[
\Bil(\x,\y) \Rightarrow \Bil(\x,\y) \oplus \Bil(\y,\x) \Rightarrow \Bil(\x,\y)
\]
\end{remark}

As quadratic functors are only 2-excisive, but not 1-excisive, they generally don't preserve exact squares. Their failure to preserve exact squares is however completely controlled by the associated symmetric bilinear parts. More precisely, we have the following:

\begin{lemma}
\label{lemma:goodwillie}%
Let \(\QF \colon \C\op \to \Spa\) be a quadratic functor with bilinear part \(\Bil = \Bil_{\QF}\) and let
\begin{equation}
\label{equation:goodwillie}%
\begin{tikzcd}
\x \ar[r,"{\alp'}"] \ar[d,"{\bet'}"'] & \y \ar[d,"{\bet}"] \\
\z \ar[r,"{\alp}"] & \w
\end{tikzcd}
\end{equation}
be an exact square in \(\C\). Then in the diagram
\begin{equation}
\label{equation:associated-squares}%
\begin{tikzcd}
\QF(\w) \ar[r]\ar[d] & \Bil(\z,\y) \ar[d] & \Bil(\cof(\bet'),\cof(\alp')) \ar[l] \ar[d] \\
\QF(\z) \times_{\QF(\x)} \QF(\y) \ar[r] & \Bil(\z,\x) \times_{\Bil(\x,\x)}\Bil(\x,\y) & 0 \ar[l]
\end{tikzcd}
\end{equation}
both squares are exact. In particular, there is a natural equivalence
\[
\cof[\QF(\w) \to \QF(\z) \times_{\QF(\x)} \QF(\y)] \simeq \cof[\Bil(\z,\y) \to \Bil(\z,\x) \times_{\Bil(\x,\x)}\Bil(\x,\y)] \simeq \Sig\Bil(\cof(\bet'),\cof(\alp')) \simeq \Bil(\fib(\bet'),\cof(\alp')).
\]
\end{lemma}
\begin{proof}
Consider the following pair of maps between commutative squares
\begin{equation}
\label{equation:maps-between-squares}%
\begin{tikzcd}
\QF(\w) \ar[r]\ar[d] & \QF(\y) \ar[d] \\
\QF(\z) \ar[r] & \QF(\x)
\end{tikzcd}
\Rightarrow
\begin{tikzcd}
\QF(\z \oplus \y) \ar[r]\ar[d] & \QF(\x \oplus \y) \ar[d] \\
\QF(\z \oplus \x) \ar[r] & \QF(\x \oplus \x)
\end{tikzcd}
\Rightarrow
\begin{tikzcd}
\Bil(\z,\y) \ar[r]\ar[d] & \Bil(\x, \y) \ar[d] \\
\Bil(\z,\x) \ar[r] & \Bil(\x,\x)
\end{tikzcd}
\end{equation}
where the left one is induced by the strongly cocartesian cube
\begin{equation}
\label{equation:cube}%
\begin{tikzcd}
[row sep=small, column sep=small]
& x \oplus x \arrow[dl] \arrow[rr] \arrow[dd] & & \x \oplus \y \arrow[dl] \arrow[dd] \\
\x \arrow[rr, crossing over] \arrow[dd] & & y \\
& \z \oplus \x \arrow[dl] \arrow[rr] &  & \z \oplus \y. \arrow[dl] \\
\z \arrow[rr] & & \w \arrow[from=uu, crossing over]
\end{tikzcd}
\end{equation}
Here, the map \(\x \oplus \x \to \x\) is the collapse map, the map \(\z \oplus y \to \w\) is the one whose components are \(\alp\) and \(\beta\), and the maps \(\x \oplus \y \to \y\) and \(\x \oplus \z \to \z\) have one component the identity and one component \(\alp'\) or \(\beta'\), respectively. Since \(\QF\) is quadratic it is in particular \(2\)-excisive by the first characterization in Proposition~\refone{proposition:basic-properties-quad-functors}, and so \(\QF\) maps~\eqrefone{equation:cube} to a cartesian cube of spectra. This means that the first map in~\eqrefone{equation:maps-between-squares} induces an equivalence on total fibres. On the other hand, the second map in~\eqrefone{equation:maps-between-squares} also induces an equivalence on total fibres since its cofibre is the square
\[
\begin{tikzcd}
\QF(\z) \oplus \QF(\y) \ar[r]\ar[d] &  \QF(\x) \oplus \QF(\y) \ar[d] \\
\QF(\z) \oplus \QF(\x) \ar[r] & \QF(\x) \oplus \QF(\x)
\end{tikzcd}
\]
whose total fibre is trivial. We then deduce that the composite of the two maps in~\eqrefone{equation:maps-between-squares} induces an equivalence on total fibres, and hence the left square in~\eqrefone{equation:associated-squares} is exact. Finally, the right square in~\eqrefone{equation:associated-squares} is exact because \(\Bil(-,-)\) is exact in each variable separately and hence the total fibre of the right most square in~\eqrefone{equation:maps-between-squares} identifies with \(\Bil(\cof(\bet'),\cof(\alp'))\) via the natural map \(\Bil(\cof(\bet'),\cof(\alp')) \to \Bil(\z,\y)\).
\end{proof}

\begin{remark}
\label{remark:goodwillie-dual}%
Lemma~\refone{lemma:goodwillie} admits a natural dual variant. Given a quadratic functor \(\QF \colon \C\op \to \Spa\) with bilinear part \(\Bil = \Bil_{\QF}\) and an exact square as in~\eqrefone{equation:goodwillie}, one may form instead the diagram
\begin{equation}
\label{equation:associated-squares-dual}%
\begin{tikzcd}
\QF(\z) \oplus_{\QF(\w)}\QF(\y) \ar[d] & \Bil(\z,\w) \oplus_{\Bil(\w,\w)}\Bil(\w,\y) \ar[d]\ar[l]\ar[r] & 0 \ar[d] \\
\QF(\x)  & \Bil(\z,\y) \ar[l]\ar[r] & \Bil(\fib(\alp),\fib(\bet))
\end{tikzcd}
\end{equation}
obtained using the maps on the left hand side of~\eqrefone{equation:cross-effect} instead of the right. The dual of the argument in the proof of Lemma~\refone{lemma:goodwillie} then shows that~\eqrefone{equation:associated-squares-dual} consists of two exact squares, yielding a natural equivalence
\[
\fib[\QF(\z) \oplus_{\QF(\w)}\QF(\y) \to \QF(\x)] \simeq \fib[\Bil(\z,\w) \oplus_{\Bil(\w,\w)}\Bil(\w,\y) \to \Bil(\z,\y)] \simeq \Om\Bil(\fib(\alp),\fib(\bet)) \simeq \Bil(\cof(\alp),\fib(\bet)).
\]
\end{remark}

Applying Lemma~\refone{lemma:goodwillie} in the case where \(\z=0\) we obtain:
\begin{corollary}[{cf.\ \cite[Lecture 9, Theorem 5]{Lurie-L-theory}}]
\label{example:usual-sequence}%
For an exact sequence \(\x \to \y \to \w\) in \(\C\), the natural map
\[
\begin{tikzcd}
\QF(\w) \ar[r] & \mathrm{totfib}
\end{tikzcd}
\left[
\begin{tikzcd}
\QF(\y) \ar[d] \ar[r] & \QF(\x) \ar[d] \\
\Bil_{\QF}(\x,\y) \ar[r] & \Bil_{\QF}(\x,\x)
\end{tikzcd}
\right]
\]
from \(\QF(\w)\) to the total fibre of the square on the right, is an equivalence.
\end{corollary}

\begin{definition}
\label{definition:linear-part}%
For a quadratic functor \(\QF\) we will denote by \(\Lin_{\QF}\colon \C\op \to \Spa\) the cofibre of the natural transformation \((\Bil^{\diag}_{\QF})_{\hC} \Rightarrow \QF\), which is exact by Proposition~\refone{proposition:basic-properties-quad-functors}, and refer to it as the \defi{linear part} of \(\QF\). By construction, the linear part \(\Lin_{\QF}\) sits in an exact sequence
\begin{equation}
\label{equation:homogeneous-linear}%
\Bil_{\QF}(\x,\x)_{\hC} \to \QF(\x) \to \Lin_{\QF}(\x).
\end{equation}
The formation of linear parts can be organized into a functor
\begin{equation}
\label{equation:linear}%
\Lin_{(-)} \colon \Funq(\C) \longrightarrow \Fun^{\ex}(\C\op,\Spa)
\end{equation}
whose post-composition with the inclusion \(\Fun^{\ex}(\C\op,\Spa) \subseteq \Funq(\C)\) carries a natural transformation
from the identity \(\QF \Rightarrow \Lin_{\QF}\), corresponding to the second arrow in~\eqrefone{equation:homogeneous-linear}.
\end{definition}

\begin{remark}
\label{remark:base-change-linear}%
It follows from Remark~\refone{remark:invariance-base-change} that the formation of linear parts naturally commutes with restriction along an exact functor \(f\colon\C \to \D\).
\end{remark}

\begin{lemma}
\label{lemma:linear}%
The natural map \(\QF \Rightarrow \Lin_{\QF}\) is a unit exhibiting \(\Lin_{(-)}\) as left adjoint to the inclusion \(\Fun^{\ex}(\C\op,\Spa) \subseteq \Funq(\C)\).
\end{lemma}
\begin{proof}
Since \(\Fun^{\ex}(\C\op,\Spa) \subseteq \Funq(\C)\) is a full inclusion it will suffice to show that \(\QF \Rightarrow \Lin_{\QF}\) induces an equivalence on mapping spectra to every exact functor. Since \(\Funq(\C)\) is stable this is the same as saying that the fibre of \(\QF \Rightarrow \Lin_{\QF}\) maps trivially to any exact functor. This fibre is \([\Bil_{\QF}^{\diag}]_{\hC}\) by construction,
and so it will hence suffice to show that \(\Bil_{\QF}^{\diag}\) maps trivially to any exact functor. Indeed, this follows from the adjunction of Remark~\refone{remark:adj-bilinear-part-diag} since the bilinear part of every linear functor vanishes.
\end{proof}

Let us also remark that equivalences of quadratic functors can be detected on their connective covers.
\begin{lemma}
\label{lemma:coconnective-quadratic-functors}%
Let \(\C\) be a stable \(\infty\)-category and \(\QF\colon\C\op\to \Spa\) be a quadratic functor. Suppose that for every \(x\in\C\) the spectrum \(\QF(x)\) is coconnective. Then \(\QF\) is the zero functor. In particular, if a natural transformation of quadratic functors \(f\colon\QF\to\QF'\) is an equivalence after applying \(\Omega^\infty\), then it is itself an equivalence.
\end{lemma}
\begin{proof}
First suppose that \(\QF\) is exact. Then, for every \(\x\in\C\) and \(n\in\ZZ\), \(\pi_n\QF(x)=\pi_1\QF(\Sigma^{n-1}x)=0\), and so \(\QF=0\).

Let us now prove the general case. For every \(\x,\y\in\C\), the spectrum \(\Bil_\QF(\x,\y)\) is a direct summand of \(\QF(\x\oplus \y)\). In particular it is also coconnective. Hence, if we fix an \(\x\in\C\), then \(\Bil_\QF(\x,-)\colon \C\op\to\Spa\) is an exact functor taking values in coconnective spectra, and therefore the zero functor by the previous argument. The cross-effect \(\Bil_\QF\) is therefore the zero functor. In particular, \(\QF\) is exact, and is hence the zero functor by the same argument.

The final statement follows by applying the previous argument to the fibre of \(f\).
\end{proof}

We finish this subsection with a discussion of the left and right adjoints to the inclusion of quadratic functors inside reduced functors.

\begin{construction}
\label{construction:excisive-approx}%
Let \(\E\) be a stable \(\infty\)-category.
Given a quadratic functor \(\QF\colon \E\op \to \Spa\),
Lemma~\refone{lemma:goodwillie} applied in the case where both \(z\) and \(y\) are zero objects implies that the sequence
\begin{equation}
\label{equation:canonical-sequence}%
\QF(\w) \to \Om\QF(\Om\w) \to \Om\Bil_{\QF}(\Om\w,\Om\w)
\end{equation}
is exact, and hence that the natural map
\[
\QF(\w) \xrightarrow{\simeq} \Om\fib[\QF(\Om\w) \to \Bil_{\QF}(\Om\w,\Om\w)]
\]
is an equivalence. This map itself is however defined for any reduced \(\QF\), and is natural in \(\QF\). In particular, given a stable \(\infty\)-category we may define a functor
\[
\rT^{\E}_2\colon \Fun_\ast(\E\op,\Spa) \to \Fun_{\ast}(\E\op,\Spa)
\]
which sends a reduced functor \(\RF\colon \E\op \to \Spa\) to the reduced functor
\begin{equation}
\label{equation:T-2}%
\rT^{\E}_2(\RF) := \Om\fib[\RF(\Om\w) \to \Bil_{\RF}(\Om\w,\Om\w)].
\end{equation}
The operation \(\rT^{\E}_2\) is equipped with a natural map
\[
\theta_{\RF}\colon \RF \Rightarrow \rT^{\E}_2(\RF)
\]
which is an equivalence when \(\RF\) is quadratic by Lemma~\refone{lemma:goodwillie}.
Unwinding the definitions, we see that the association \(\RF \mapsto \rT^{\E}_2(\RF)\) identifies with the one defined in~\cite[Construction 6.1.1.22]{HA} for \(\C = \E\op\) and \(\D = \Spa\). Since \(\Spa\) is stable and admits small colimits it is in particular differentiable in the sense of~\cite[Definition 6.1.1.6]{HA}. By~\cite[Theorem 6.1.1.10]{HA} we may then conclude that the association
\[
\App_2(\RF) := \colim[\RF \xrightarrow{\theta_R} \rT^{\E}_2(\RF) \xrightarrow{\theta_{\rT^{\E}_2(\RF)}} \rT^{\E}_2\rT^{\E}_2(\RF) \to \cdots]
\]
gives a left adjoint to the inclusion \(\Funq(\E) \subseteq \Fun_*(\E,\Spa)\). This procedure is often referred to as \defi{\(2\)-excisive approximation}.

In a dual manner, if we use Remark~\refone{remark:goodwillie-dual} instead of Lemma~\refone{lemma:goodwillie} then we get that for a quadratic functor \(\QF\) the sequence
\begin{equation}
\label{equation:canonical-sequence-dual}%
\Sig\Bil_{\QF}(\Sig\w,\Sig\w) \to \Sig\QF(\Sig\w) \to \QF(\w)
\end{equation}
is exact, and so the natural map
\[
\Sig\cof[\Bil_{\QF}(\Sig\w,\Sig\w) \to \QF(\Sig\w)] \xrightarrow{\simeq} \QF(\w)
\]
is an equivalence.
As above, for a general reduced functor \(\RF\) we can define the functor
\begin{equation}
\label{equation:T-to-the-two}%
\rT^2_{\E}(\RF) = \Sig\cof[\Bil_{\RF}(\Sig\w,\Sig\w) \to \RF(\Sig\w)],
\end{equation}
equipped with a natural map
\[
\tau_{\RF}\colon \rT^2_{\E}(\RF) \Rightarrow \RF ,
\]
which is an equivalence when \(\RF\) is quadratic.
We may also identify \(\rT^2_{\E}\) with the result of~\cite[Construction 6.1.1.22]{HA} applied to \(\C = \E\) and \(\D=\Spa\op\). Since \(\Spa\op\) is also differentiable by the same argument it follows from~\cite[Theorem 6.1.1.10]{HA} that the association
\[
\App^2(\RF) := \lim[\cdots \to \rT^2_{\E}\rT^2_{\E}(\RF) \xrightarrow{\tau_{\rT^2_{\E}(\RF)}} \rT^2_{\E}(\RF) \xrightarrow{\tau_{\RF}} \RF]
\]
provides a \emph{right adjoint} to the inclusion \(\Funq(\E) \subseteq \Fun_*(\E,\Spa)\).
\end{construction}

\subsection{Hermitian and Poincaré $\infty$-categories}
\label{subsection:hermitian-and-poincare-cats}%

In this subsection we introduce the key player in this paper - the notion of a \defi{Poincaré \(\infty\)-category}. For this, it will be convenient to pass first through the following weaker notion:

\begin{definition}
\label{definition:herm-cat}%
A \defi{hermitian} \(\infty\)-category is a pair \((\C,\QF)\) where \(\C\) is a small stable \(\infty\)-category and \(\QF\colon \C\op \to \Spa\) is a quadratic functor in the sense of Definition~\refone{definition:herm-structure}. We will then also refer to \(\QF\) as a \defi{hermitian structure} on \(\C\).
The collection of hermitian \(\infty\)-categories can be organized into a (large) \(\infty\)-category \(\Cath\),
obtained as the cartesian Grothendieck construction of the functor
\[
(\Catx)\op \longrightarrow \CAT, \quad \C \longmapsto \Funq(\C).
\]
(here \(\CAT\) stands for the \(\infty\)-category of possibly large \(\infty\)-categories). We shall also refer to its morphisms as \defi{hermitian functors}.
\end{definition}

Unpacking this definition, we find that a hermitian functor from \((\C,\QF)\) to \((\Ctwo,\QFtwo)\) consists of an exact functor \(f \colon \C \rightarrow \Ctwo\) and a natural transformation \(\eta \colon \QF \Rightarrow f^*\QFtwo := \QFtwo \circ f\op\). We will thus generally denote hermitian functors as pairs \((f,\eta)\) of this form.
If \((f,\eta)\colon (\C,\QF) \to (\Ctwo,\QFtwo)\) is a hermitian functor then by Remark~\refone{remark:invariance-base-change} we have a natural equivalence \((f \times f)^*\Bil_{\QFtwo} \simeq\Bil_{f^*\QFtwo}\), and consequently the natural transformation \(\eta\) determines a natural transformation
\begin{equation}
\label{equation:induced-on-forms}%
\beta_{\eta}\colon \Bil_{\QF} \Rightarrow (f \times f)^*\Bil_{\QFtwo},
\end{equation}
which we then denote by \(\beta_{\eta}\).

The notion of Poincaré \(\infty\)-category is obtained from that of a hermitian \(\infty\)-category \((\C,\QF)\) by requiring \(\QF\) to satisfy two non-degeneracy conditions. Both of these conditions depend only on the underlying symmetric bilinear part \(\Bil_{\QF}\in \Funs(\C)\).
To formulate the first one we first note that the exponential equivalence
\[
\Fun(\C\op\times\C\op,\Spa) \xrightarrow{\simeq} \Fun(\C\op,\Fun(\C\op,\Spa))
\]
restricts to an equivalence
\begin{equation}
\label{equation:exp-bilinear}%
\Funb(\C) \xrightarrow{\simeq} \Funx(\C\op,\Funx(\C\op,\Spa)).
\end{equation}
We then consider the following condition:

\begin{definition}
\label{definition:nondeg-herm-cat}%
We will say that
a  bilinear functor \(\Bil \in \Funb(\C)\) is \defi{right non-degenerate} if the associated exact functor
\begin{equation}
\label{equation:ind-duality}%
\C\op \to \Funx(\C\op,\Spa) \quad\quad \y \mapsto \Bil(-,\y)
\end{equation}
takes values in the essential image of the stable Yoneda embedding
\[
\C\hrar \Fun^{\lex}(\C\op,\Sps) \simeq \Funx(\C\op,\Spa),
\]
where \(\Fun^{\lex}\) denotes left exact (that is, finite limit preserving) functors, and the equivalence to the last term is by~\cite[Corollary 1.4.2.23]{HA}. In other words, if for each \(\y \in \C\) the presheaf of spectra \(\Bil(-,\y)\) is representable by an object in \(\C\). In this case we can factor~\eqrefone{equation:ind-duality} essentially uniquely as a functor
\[
\Dual_{\Bil}\colon \C\op \to \C
\]
followed by
\(\C \hrar \Funx(\C\op,\Spa)\), so that we obtain an equivalence
\[
\Bil(\x,\y) \simeq \map_{\C}(\x,\Dual_{\Bil} \y).
\]
Similarly, \(\Bil \in \Funb(\C)\) is called \defi{left non-degenerate} if the associated exact functor
\(\x \mapsto \Bil(\x,-)\)
takes values in the essential image of the stable Yoneda embedding.
If \(\Bil \in \Funb(\C)\) is left and right non-degenerate, then it is called non-degenerate.
In this case the two resulting dualities are, essentially by definition, adjoint to each other.

We will say that a symmetric bilinear functor is \defi{ non-degenerate} if the underlying bilinear functor is. In this case it of course suffices to check that it is right non-degenerate. The two dualities are in this case equivalent and we will refer to the representing functor \(\Dual_{\Bil}\)  as the \defi{duality} associated to the non-degenerate symmetric bilinear functor \(\Bil\) (though we point out that \(\Dual_{\Bil}\) is not in general an equivalence).
Given a hermitian structure \(\QF\) on a stable \(\infty\)-category \(\C\), we will say that \(\QF\) is \defi{non-degenerate} if its underlying bilinear part is. In this case we will also say that \((\C,\QF)\) is a non-degenerate hermitian \(\infty\)-category and will denote the associated duality by \(\Dual_{\QF}\).
\end{definition}

The full subcategories of \(\Funb(\C),\Funs(\C)\) and \(\Funq(\C)\) spanned by the non-degenerate functors will be denoted \(\Funnb(\C),\Funns(\C)\) and \(\Funnq(\C)\), respectively. The bilinear exponential equivalence~\eqrefone{equation:exp-bilinear}
then restricts to an equivalence
\[
\Funnb(\C) \xrightarrow{\simeq} \FunR(\C\op,\C),
\]
where \(\FunR\) denotes the right adjoint functors. To see this it suffices to observe that \(\Bil \in \Funb(\C)\) is non-degenerate precisely if it is right non-degenerate and the resulting duality admits a left adjoint. Under this equivalence the \(\Ct\)-action on the left corresponds to the \(\Ct\)-action on the right given by passing to adjoints and taking opposites, so that we get an equivalence
\[
\Funns(\C) \xrightarrow{\simeq} \FunR(\C\op,\C)^{\hC}.
\]
Both of these equivalences will be denote by \(\Bil \mapsto \Dual_{\Bil}\). Similarly, we will also denote the composition
\[
\Funnq(\C) \xrightarrow{\Bil_{(-)}} \Funns(\C) \xrightarrow{\Dual_{(-)}} \FunR(\C\op,\C)
\]
by \(\QF \mapsto \Dual_{\QF}\).

Let us make these adjointability statements explicit: if \(\Bil \in \Funs(\C)\) is a non-degenerate symmetric bilinear functor with associated duality \(\Dual = \Dual_{\Bil}\colon \C\op \to \C\) then the symmetric structure of \(\Bil\) determines a natural equivalence
\begin{equation}
\label{equation:sym-adj}%
\map_{\C}(\x,\Dual(\y)) \simeq \Bil(\x,\y) \simeq \Bil(\y,\x) \simeq \map_{\C}(\y,\Dual(\x)) \simeq \map_{\C\op}(\Dual\op(\x),\y)
\end{equation}
where \(\Dual\op\colon \C \to \C\op\) is the functor induced by \(\Dual\) upon taking opposites. Such a natural equivalence exhibits in particular \(\Dual\op\) as left adjoint to \(\Dual\). We will denote by
\begin{equation}
\label{equation:evaluation}%
\ev\colon \id \Rightarrow \Dual\Dual\op
\end{equation}
the unit of this adjunction, and refer it as the \defi{evaluation map} of \(\Dual\). Its individual components
\begin{equation}
\label{equation:evaluation-2}%
\ev_{\x}\colon\x \to \Dual\Dual\op(\x)
\end{equation}
are then the maps corresponding to identity \(\Dual\op(\x) \to \Dual\op(\x)\) under the equivalence~\eqrefone{equation:sym-adj}. The counit of this adjunction is given again by natural transformation~\eqrefone{equation:evaluation}, but interpreted as an arrow from \(\Dual\op\Dual\) to the identity in the \(\infty\)-category \(\Fun(\C\op,\C\op) \simeq \Fun(\C,\C)\op\).

\begin{remark}
\label{remark:reproduce}%
The process of viewing the equivalence~\eqrefone{equation:sym-adj} as an adjunction between \(\Dual\) and \(\Dual\op\) and extracting its unit as above can be reversed: knowing that \(\ev\) is a unit of an adjunction we can reproduce the equivalence \(\map_{\C}(y,\Dual(\x))  \simeq \map_{\C}(\x,\Dual(y))\) as the composite
\[
\map_{\C}(y,\Dual(\x)) \simeq \map_{\C\op}(\Dual\op(\x),y) \to \map_{\C}(\Dual\Dual\op(\x),\Dual(y)) \to \map_{\C}(\x,\Dual(\y)) ,
\]
where the last map is induced by pre-composition with the evaluation map.
\end{remark}

\begin{lemma}
\label{lemma:nat-duality}%
Let \((\C,\QF),(\Ctwo,\QFtwo)\) be two non-degenerate hermitian \(\infty\)-categories with associated dualities \(\Dual_{\QF}\) and \(\Dual_{\QFtwo}\), and let \(f,g \colon \C \to \Ctwo\) be two exact functors. Then there is a natural equivalence
\[
\nat(\Bil_{\QF},(f \times g)^*\Bil_{\QFtwo}) \simeq \nat(f\Dual_{\QF}, \Dual_{\QFtwo}g\op),
\]
where \(\nat\) stands for the the spectrum of (non-symmetric) natural transformations, on the left between two spectrum valued functors on \(\C\op \times \C\op\), and on the right between two functors \(\C\op \to \C'\).
\end{lemma}
\begin{proof}
Consider the left Kan extension functor
\[
(f \times \id)_!\colon\Fun(\C\op \times\C\op,\Spa) \to \Fun(\Ctwo\op \times \C\op,\Spa),
\]
which is left adjoint to the corresponding restriction functor. Natural transformations
\[
\Bil_{\QF} \Rightarrow (f \times g)^*\Bil_{\QFtwo} \simeq (f \times \id)^*(\id \times g)^*\Bil_{\QFtwo}
\]
then correspond under this adjunction to natural transformations
\begin{equation}
\label{equation:after-kan}%
(f\times \id)_!\Bil_{\QF} \Rightarrow (\id \times g)^*\Bil_{\QFtwo}.
\end{equation}
Now for \(\y \in \C\) we have
\[
((f \times \id)_!\Bil_{\QF})|_{\Ctwo\op \times \{y\}} \simeq (f \times \{y\})_!((\Bil_{\QF})|_{\C\op \times \{y\}}),
\]
as can be seen by the pointwise formula for left Kan extension.
Since \(\Bil_{\QF}(-,\y)\) is represented by \(\Dual_{\QF}(\y)\) and left Kan extension preserves representable functors it then follows that
\[
(f\times \id)_!\Bil_{\QF}(\x',\y) \simeq \map_{\Ctwo}(\x',f(\Dual_{\QF}(\y)))
\]
for \((\x',\y) \in \Ctwo\op \times \C\op\).
On the other hand, we have \((\id \times g)^*\Bil_{\QFtwo}(\x',\y) \simeq \map_{\Ctwo}(\x',\Dual_{\QFtwo}(g(\y)))\), and so by the fully-faithfulness of the Yoneda embedding we thus obtain
\[
\nat(\Bil_{\QF},(f \times g)^*\Bil_{\QFtwo}) \simeq \nat(f\Dual_{\QF}, \Dual_{\QFtwo}g\op),
\]
as desired.
\end{proof}

\begin{definition}
\label{definition:transformation-duality}%
Given a hermitian functor \((f,\eta) \colon (\C,\QF) \to (\Ctwo,\QFtwo)\), we will denote by
\[
\tau_{\eta}\colon f\Dual_{\QF} \Rightarrow \Dual_{\QFtwo}f\op
\]
the natural transformation corresponding to the natural transformation \(\bet_{\eta}\colon \Bil_{\QF} \Rightarrow (f \times f)^*\Bil_{\QFtwo}\) of~\eqrefone{equation:induced-on-forms}, via Lemma~\refone{lemma:nat-duality}.
\end{definition}

\begin{remark}
\label{remark:recovered}%
In the situation of Definition~\refone{definition:transformation-duality}, it follows from the triangle identities of the adjunction \((f \times \id)_! \dashv (f \times \id)^*\) that the natural transformation \(\beta_{\eta}\colon \Bil_{\QF} \Rightarrow (f \times f)^*\Bil_{\QFtwo}\) can be recovered from \(\tau_{\eta}\colon f\Dual_{\QF} \Rightarrow \Dual_{\QF'}f\op\) as the composite
\[
\Bil_{\QF}(\x,\y) \simeq \map_{\C}(\x,\Dual_{\QF}(\y)) \to \map_{\D}(f(\x),f\Dual_{\QF}(\y)) \to \map_{\D}(f(\x),\Dual_{\QFtwo}f(\y)) \simeq \Bil_{\QF}(f(\x),f(\y)) ,
\]
where the two middle maps are induced by the action of \(f\) on mapping spectra and post-composition with \(\tau_{\eta}\), respectively.
\end{remark}

\begin{definition}
\label{definition:duality}%
A hermitian functor \((f,\eta) \colon (\C,\QF) \rightarrow (\Ctwo,\QFtwo)\) between non-degenerate hermitian \(\infty\)-categories is called \defi{duality preserving} if the transformation \(\tau_{\eta}\colon f\Dual_{\QF} \Rightarrow \Dual_{\QF}f\op\) constructed above is an equivalence.
\end{definition}

\begin{definition}
\label{definition:poinc-cats}%
A symmetric bilinear functor \(\Bil\) is called \defi{perfect} if the evaluation map \(\id_\C \Rightarrow \Dual_{\Bil}\Dual_{\Bil}\op\) of~\eqrefone{equation:evaluation} is an equivalence. An hermitian structure \(\QF\) is called \defi{Poincaré} if the underlying bilinear functor of \(\QF\) is perfect. In this case we will say that \((\C,\QF)\) is a \defi{Poincaré \(\infty\)-category}. We will denote by
\[
\Catp \subseteq \Cath
\]
 the (non-full) subcategory spanned by the Poincaré \(\infty\)-categories and duality preserving functors, and will generally refer to duality-preserving hermitian functors between Poincaré \(\infty\)-categories as \emph{Poincaré functors}.
For a stable \(\infty\)-category \(\C\), we will denote by
\[
\Funp(\C) \subseteq \Funq(\C)
\]
the subcategory spanned by those hermitian structures which are Poincaré, and those natural transformations \(\eta\colon \QF \Rightarrow \QFtwo\) which are duality preserving, that is, for which the associated hermitian functor \((\id,\eta)\colon (\C,\QF) \to (\C,\QFtwo)\) is Poincaré.
\end{definition}

\begin{remark}
\label{remark:duality-equivalence}%
A symmetric bilinear functor \(\Bil\) is perfect if and only if it is non-degenerate and \(\Dual_{\Bil} \colon \C\op \rightarrow \C\) is an equivalence of categories. Indeed, an adjunction consists of a pair of inverse equivalences if and only if its unit and counit are equivalences.
\end{remark}

If \(\Bil\) is a perfect bilinear functor on \(\C\) then the duality \(\Dual_{\Bil}\colon \C\op \xrightarrow{\simeq}\C\) is not just an equivalence of \(\infty\)-categories, but carries a significant amount of extra structure. To make this precise note that there is a \(\Ct\)-action on \(\Catx\) given by sending \(\C \mapsto \C\op\). This can be seen by
using simplicial sets as a model where taking the opposite gives an action on the nose. Alternatively one can also use that the space of autoequivalences of \(\Catx\) is equivalent to the discrete group \(\Ct\) as shown in \cite{ToenVersuneaxiomatisation}, see also \cite[Theorem 4.4.1]{LurieGoodwillie}.

\begin{definition}
\label{definition:perfect-duality}%
A \defi{stable \(\infty\)-category with perfect duality} is a homotopy fixed point of \(\Catx\) with respect to the \(\Ct\)-action given by taking the opposite \(\infty\)-category. Equivalently, it is, a section \(\BC \to \wtl{\Catx}\) of the fibration \(\wtl{\Catx} \to \BC\) encoding the \(\Ct\)-action on the \(\infty\)-category \(\Catx\) given by taking opposites.
\end{definition}

We note that a stable \(\infty\)-category with perfect duality consists in particular of a stable \(\infty\)-category \(\C\) and an equivalence \(\Dual\colon \C \to \C\op\), equipped with additional coherence structure of being a homotopy \(\Ct\)-fixed point. For example, the composition \(\Dual\Dual\op\) is equipped with a natural equivalence \(\ev\colon\id \simeq \Dual\Dual\op\), which itself carries higher coherence homotopies relating it with its opposite, and so forth.
By a \defi{perfect duality} on a given stable \(\infty\)-category \(\C\) we will mean a refinement of \(\C\) to a homotopy \(\Ct\)-fixed point of \(\Catx\).
We may also identify the notion of a perfect duality with that of a \(\Ct\)-fixed equivalence \(\C\op \xrightarrow{\simeq} \C\), where the \(\Ct\)-action on \(\FunR(\C\op,\C)\) is obtained via its identification with the \(\infty\)-category \(\Funnb(\C)\) of non-degenerate bilinear functors.
We will often abuse notation and denote a perfect duality simply by its underlying equivalence \(\Dual\colon \C\op \to \C\).

In their work on \(\infty\)-categories with duality, Heine-Lopez-Avila-Spitzweck prove that the duality functor \(\Dual_{\Bil}\) associated to a perfect bilinear functor \(\Bil\) on a stable \(\infty\)-category \(\C\), naturally refines to a perfect duality on \(\C\) in the above sense. Furthermore, the association \(\Bil \mapsto \Dual_{\Bil}\) determines an equivalence between perfect bilinear functors on \(\C\) and perfect dualities on \(\C\)
see~\cite[Corollary 7.3]{HeineLopez-AvilaSpitzweck}, and~\cite[Proposition 2.1]{Spitzweck-GW} for the stable variant. Together with Lemma~\refone{lemma:nat-duality}, this association determines a forgetful functor
\begin{equation}
\label{equation:forget-Q-to-D}%
\Catp \to (\Catx)^{\hC} \quad\quad (\C,\QF) \mapsto (\C,\Dual_{\QF})
\end{equation}
from Poincaré \(\infty\)-categories to stable \(\infty\)-categories with perfect duality. This provides a key link between the present setup and the existing literature on stable \(\infty\)-categories with duality.

\begin{definition}
\label{definition:symmetric-quadratic}%
Given a stable \(\infty\)-category and a symmetric bilinear functor \(\Bil\colon \C\op \times \C\op \to \Spa\) we will refer to the hermitian structures \(\QF^{\sym}_{\Bil},\QF^{\qdr}_{\Bil} \in \Funq(\C)\) of Example~\refone{example:quadratic-symmetric} as the \defi{symmetric} and \defi{quadratic} hermitian structures associated to \(\Bil\), respectively. As the symmetric bilinear parts of both \(\QF^{\sym}_{\Bil}\) and \(\QF^{\qdr}_{\Bil}\) are canonically equivalent to \(\Bil\), these hermitian structures are Poincaré if and only if \(\Bil\) is perfect.
\end{definition}

\begin{example}
\label{example:perfect-derived}%
Let \(R\) be an ordinary commutative ring and let \(\C = \Dperf(R)\) be the perfect derived \(\infty\)-category of \(R\). Similar as in Example \refone{example:quadratic-symmetric}
we may then consider the symmetric bilinear functor \(\Bil_R \in \Funb(\C)\) given by
\[
\Bil_R(X, Y) = \map_R(X \otimes_R Y, R),
\]
with symmetric structure induced by the symmetric structure of the tensor product \(\otimes\). This bilinear functor is perfect with duality given by
\[
\Dual_{R}(Y) = \HHom_R(Y,R) ,
\]
where the right hand side stands for the internal mapping complex. An element \(\beta \in \Om^{\infty}\Bil(X,Y)\) then corresponds to a map \(X \otimes_R Y \to R\), which we can consider as a bilinear form on the pair \((X,Y)\).
To this perfect bilinear functor we can associate the corresponding symmetric and quadratic Poincaré structures
\[
\QF^{\sym}_{R}(X) := \Bil_R(X,X)^{\hC} \quad\text{and}\quad
\QF^{\qdr}_{R}(X) := \Bil_R(X,X)_{\hC},
\]
as in Definition~\refone{definition:symmetric-quadratic}. The space \(\Om^{\infty}\QF^{\sym}_{R}(X)\) is then the space of (homotopy) \(\Ct\)-fixed points of \(\Om^{\infty}\Bil(X,X)\), which should be viewed as the homotopical counterpart of the notion of a symmetric form on \(X\).
The space \(\Om^{\infty}\QF^{\qdr}_{R}(X)\), on the other hand, is the space of (homotopy) \(\Ct\)-orbits of \(\Om^{\infty}\Bil(X,X)\), which we can consider as a homotopical analogue of that a \emph{quadratic form} on \(X\), see Example~\refone{example:quadratic-symmetric}.
\end{example}

\begin{example}
\label{example:anti-symmetric}%
In the situation of Example~\refone{example:perfect-derived} we could
also consider the symmetric bilinear functor \(\Bil_{-R}\) whose underlying bilinear functor is \(\Bil_{R}\) but whose symmetric structure is twisted by the sign action of \(\Ct\). In other words, the symmetry equivalence \(\Bil_{-R}(X \otimes_R Y,R) \xrightarrow{\simeq} \Bil_{-R}(Y \otimes_R X,R)\) of \(\Bil_{-R}\) is minus the one of \(\Bil_R\). This bilinear functor is again perfect with duality which coincides with \(\Dual_R\) on the level of the underlying equivalence \(\Dperf(R)\op \to \Dperf(R)\), but which has a different double dual identification.
We may then consider the corresponding symmetric and quadratic Poincaré structures
\[
\QF^{\sym}_{-R}(X) := \Bil_{-R}(X,X)^{\hC} \quad\text{and}\quad
\QF^{\qdr}_{-R}(X) := \Bil_{-R}(X,X)_{\hC},
\]
as in Definition~\refone{definition:symmetric-quadratic}. The space \(\Om^{\infty}\QF^{\sym}_{-R}(X)\) is then a homotopical counterpart of the notion of an anti-symmetric form on \(X\), while  \(\Om^{\infty}\QF^{\qdr}_{-R}(X)\) is its quadratic counterpart.
\end{example}

\begin{example}
In the spirit of Example~\refone{example:perfect-derived}, one may also fix a scheme \(X\) and consider the stable \(\infty\)-category \(\Dperf(X)\) of perfect complexes of quasi-coherent sheaves on \(X\). Given a line bundle \(L\) on \(X\) we have an associated bilinear form \(\Bil_L\) on \(\Dperf(X)\) given by
\[
\Bil_{L}(\F, \G) = \map_X(\F \otimes_X \G, L),
\]
which is perfect with duality
\[
\Dual_{L}(\F) = \HHom_X(\F,L) .
\]
To this perfect duality we can then associate the corresponding symmetric and quadratic Poincaré structures
\[
\QF^{\sym}_{L}(\F) := \Bil_L(\F,\F)^{\hC} \quad\text{and}\quad \QF^{\qdr}_{L}(\F) := \Bil_L(\F,\F)_{\hC}.
\]
\end{example}

\begin{example}
\label{example:universal-category}%
Let \(\Spaf\) be the \(\infty\)-category of finite spectra.
We define a hermitian structure on \(\Spaf\) via the pullback square
\[
\begin{tikzcd}
	\QF^{\uni}(X) \ar[r] \ar[d] & \Dual(X) \ar[d] \\ \Dual(X\otimes X)^{\hC} \ar[r] & \Dual(X\otimes X)^{\tC}
\end{tikzcd}
\]
where \(\Dual(X)=\map(X,\mathbb{S})\) denotes the Spanier-Whitehead dual
and the right vertical map is the Tate diagonal \(\Dual X\to (\Dual X\otimes \Dual X)^{\tC}\) of \(\Dual X\). This hermitian structure is then Poincaré with duality given by Spanier-Whitehead duality. We note that this Poincaré structure is neither quadratic nor symmetric (Definition~\refone{definition:symmetric-quadratic}).
The superscript \(\uni\) is suggestive for universal, see \S\refone{subsection:universal}.
The Poincaré \(\infty\)-category also functions as the unit of the symmetric monoidal structure on Poincaré \(\infty\)-categories that we will construct in \S\refone{section:multiplicative}.
\end{example}

Defining hermitian structures as spectrum valued functors allows us to easily implement various useful manipulations. One of them, which plays a recurring role in this paper, is the procedure of \emph{shifting} hermitian structures:

\begin{definition}
\label{definition:shift}%
Let \(\C\) be a stable \(\infty\)-category and \(\QF\colon\C\op \to \Spa\) a quadratic functor. For \(n \in \ZZ\) we will denote by \(\QF\qshift{n}\colon \C\op \to \Spa\) the \(n\)-fold suspension of \(\QF\), given by
\[
\QF\qshift{n}(\x) = \Sig^n\QF(\x) .
\]
We note that \(\QF\qshift{n}\) is again a quadratic functor with bilinear part \(\Bil_{\Sig^n\QF} = \Sig^n\Bil_{\QF}\) and linear part \(\Lin_{\Sig^n\QF} = \Sig^n\Lin_{\QF}\); indeed, \(\Funq(\C)\) is a stable subcategory of \(\Fun(\C\op,\Spa)\) and \(\Bil_{(-)}\) and \(\Lin_{(-)}\) are both exact functors. In particular, if \(\QF\) non-degenerate or perfect then so is \(\QF\qshift{n}\) with duality \(\Dual_{\QF\qshift{n}}(\x) = \Sig^n\Dual_{\QF}(\x)\). We refer to \(\QF\qshift{n}\) as the \defi{\(n\)-fold shift} of \(\QF\), and to the hermitian \(\infty\)-category \((\C,\QF\qshift{n})\) as the \(n\)-fold shift of \((\C,\QF)\).
\end{definition}

\begin{remark}
The hermitian \(\infty\)-category \((\C,\QF\qshift{n})\) is Poincaré if and only if \((\C,\QF)\) is.
\end{remark}

\begin{example}
In the situation of Example~\refone{example:perfect-derived}, if we shift the Poincaré structures \(\QF^{\sym}_R\) and \(\QF^{\qdr}_R\) on \(\Dperf(R)\) by \(n \in \ZZ\) then we get Poincaré structures
\[
(\QF^{\sym}_{R})\qshift{n}(X) = \map_R(X \otimes_R X,R[n])^{\hC}
\]
and
\[
(\QF^{\qdr}_{R})\qshift{n}(X) = \map_R(X \otimes_R X,R[n])_{\hC} .
\]
respectively, which we consider as encoding \(n\)-shifted symmetric and quadratic forms. Here, \(R[n]\) denotes the \(R\)-complex which is \(R\) in degree \(n\) and zero everywhere else.
\end{example}

We note that the full subcategories of \(\Funq(\C)\) spanned by non-degenerate and perfect functors respectively are not preserved under pullback along exact functors \(f\colon \C \to \Ctwo\).
For example, the hermitian structure \(\QF^{\sym}_{\QQ}\) on \(\Dperf(\QQ)\) is perfect (see Example~\refone{example:perfect-derived}), but its pullback to \(\Dperf(\ZZ)\) is not even non-degenerate.
A notable exception to this is however the following:

\begin{observation}
\label{observation:poincare-subcat}%
If \((\C,\QF)\) is a non-degenerate hermitian or Poincaré \(\infty\)-category and \(\D \subseteq \C\) is a full stable subcategory such that the duality \(\Dual_{\QF}\) maps \(\D\) to itself then \((\D,\QF|_\D)\) is again non-degenerate with \(\Dual_{\QF|_\D} = \Dual_\QF|_{\D\op}\).
In particular, if \((\C,\QF)\) is Poincaré then \((\D,\QF|_{\D})\) is again Poincaré.
\end{observation}

\begin{example}
\label{example:reflective}%
If \((\C,\QF)\) is a non-degenerate hermitian \(\infty\)-category then the full subcategory \(\C^{\refl} \subseteq \C\) spanned by those objects \(\x \in \C\) for which the evaluation map \(\ev_{\x}\colon \x \to \Dual\op\Dual(\x)\) is an equivalence is preserved under the duality by the triangle identities, and hence the hermitian \(\infty\)-category \((\C^{\refl},\QF|_{\C^{\refl}})\) is again non-degenerate, and even Poincaré, since the evaluation map is now an equivalence by construction.
\end{example}

\subsection{Classification of hermitian structures}
\label{subsection:classification}%

In this section we discuss the classification of hermitian and Poincaré structures on a fixed stable \(\infty\)-category \(\C\), in terms of their linear and bilinear parts. For the hermitian part, this is essentially the \(n=2\) case classification of \(n\)-excisive functors in Goodwillie calculus, and is also a particular instance of the structure theory of \defi{stable recollements} (see~\cite[\S A.8]{HA},~\cite{barwick2016note},~\cite{ShahQuigley}).
For the purpose of self containment we however provide full proofs of the statements that we need in the present setting. In order to formulate these statements we first need to better understand the role played by the quadratic and symmetric hermitian structures \(\QF^{\qdr}_{\Bil},\QF^{\sym}_{\Bil}\) associated to a given symmetric bilinear form \(\Bil\).

\begin{lemma}
\label{lemma:homogeneous}%
Let \(\QF\colon \C\op \to \Spa\) be a quadratic functor on a small stable \(\infty\)-category \(\C\). Then the following are equivalent:
\begin{enumerate}
\item
\label{item:Ctwo-fixed-equals-Q}%
The map \(\Bil_{\QF}(\x,\x)_{\hC} \to \QF(\x)\) of~\eqrefone{equation:orbits-fixed-points} is an equivalence for every \(\x \in \C\).
\item
\label{item:Q-is-quadratic-quadratic}%
\(\QF\) is equivalent to a quadratic functor of the form
\(\QF^{\qdr}_{\Bil}\)
for some symmetric bilinear functor \(\Bil \in \Funs(\C)\) (see Example~\refone{example:quadratic-symmetric}).
\item
\label{item:NatQL-trivial}%
The spectrum of natural transformations \(\nat(\QF,\L)\) is trivial for any exact functor \(\L\in \Funx(\C\op,\Spa)\subseteq \Funq(\C)\).
\end{enumerate}
\end{lemma}

\begin{definition}
Following the conventions of Goodwillie calculus, we will refer to quadratic functors \(\QF\colon \C\op \to \Spa\) which satisfy the equivalent conditions of Lemma~\refone{lemma:homogeneous} as \defi{homogeneous}. We will denote by \(\Funqhom(\C) \subseteq \Funq(\C)\) the full subcategory spanned by the homogeneous functors.
\end{definition}

\begin{proof}[Proof of Lemma~\refone{lemma:homogeneous}]
Clearly \refoneitem{item:Ctwo-fixed-equals-Q} \(\Rightarrow\) \refoneitem{item:Q-is-quadratic-quadratic}. If we assume \refoneitem{item:Q-is-quadratic-quadratic} then \refoneitem{item:NatQL-trivial} follows from Lemma \refone{lemma:linear} since
the linear part of \(\QF^{\qdr}_{\Bil}\) vanishes by definition. Similarly if we assume \refoneitem{item:NatQL-trivial} then \refoneitem{item:Ctwo-fixed-equals-Q} follows by Lemma \refone{lemma:linear} since the linear part vanishes.
\end{proof}

\begin{corollary}[{cf.\ \cite[Proposition 6.1.4.14]{HA}}]
\label{corollary:homogeneous}%
The functor
\begin{equation}
\label{equation:orbits-on-B}%
\Funs(\C) \to \Funq(\C) \quad\quad\quad\quad \Bil \mapsto \QF^{\qdr}_{\Bil}
\end{equation}
is fully-faithful and its essential image is spanned by those quadratic functors which are homogeneous in the above sense.
\end{corollary}
\begin{proof}
By Lemma~\refone{lemma:homogeneous} the functor~\eqrefone{equation:orbits-on-B} takes values in homogeneous functors, and hence determines a functor \(\vphi\colon\Funs(\C) \to \Funqhom(\C) \subseteq \Funq(\C)\), where the latter denotes the full subcategory spanned by homogeneous quadratic functors. On the other hand, the formation of cross effects determines a functor in the other direction \(\psi\colon \Funqhom(\C) \to \Funs(\C)\). By Lemma~\refone{lemma:homogeneous} the composed functor \(\vphi \circ \psi \colon \Funqhom(\C) \to \Funqhom(\C)\) is naturally equivalent to the identity, and by Examples~\refone{example:quadratic-symmetric} the composite \(\psi \circ \vphi\) is naturally equivalent to the identity \(\Funs(\C) \to \Funs(\C)\) as well. It then follows that \(\vphi\) is an equivalence from \(\Funs(\C)\) to \(\Funqhom(\C) \subseteq \Funq(\C)\), as desired.
\end{proof}

The notion of a homogeneous quadratic functor has a dual counterpart, which consists of the quadratic functors which have a trivial mapping spectrum \emph{from} any exact functor. The argument of Lemma~\refone{lemma:homogeneous} then runs in a completely dual manner to show that this property is equivalent to the canonical map \(\QF(\x) \to \Bil_{\QF}(\x,\x)^{\hC}\) being an equivalence and is satisfied by quadratic functors of the form \(\QF^{\sym}_{\Bil}\) for any \(\Bil \in \Funs(\C)\). We will refer to such functors as \defi{cohomogeneous}, and denote by \(\Funqcoh(\C) \subseteq \Funq(\C)\) the full subcategory spanned by the cohomogeneous functors.
The argument of Corollary~\refone{corollary:homogeneous} then runs in a completely dual manner to show that the functor
\begin{equation}
\label{equation:fixed-points-of-B}%
\Funs(\C) \to \Funq(\C) \quad\quad\quad\quad \Bil \mapsto \QF^{\sym}_{\Bil}
\end{equation}
is fully-faithful and its essential image is spanned by the cohomogeneous quadratic functors.

\begin{remark}
Given a bilinear functor \(\Bil\), the quadratic functor \(\x \mapsto \Bil(\x,\x)\) is both homogeneous and cohomogeneous.
\end{remark}

\begin{proposition}
\label{proposition:universal-homogeneous}%
The natural transformation \(\eps\colon [\Bil_{\QF}^{\Del}]_{\hC} \Rightarrow \QF\) exhibits \([\Bil_{\QF}^{\Del}]_{\hC}\) as final among homogeneous functors equipped with a map to \(\QF\). Dually, the natural transformation \(\eta\colon \QF \Rightarrow [\Bil^{\Del}_{\QF}]^{\hC}\) exhibits \([\Bil_{\QF}^{\Del}]^{\hC}\) as initial among cohomogeneous functors equipped with a map from \(\QF\).
\end{proposition}
\begin{proof}
This first statement is equivalent to Lemma~\refone{lemma:linear} and the second is dual.
\end{proof}

\begin{corollary}
\label{corollary:hom-universal}%
For a small stable \(\infty\)-category \(\C\) the functor
\[
\Bil_{(-)}\colon \Funq(\C) \to \Funs(\C)
\]
admits left and right adjoints, both of which are fully faithful, given by sending \(\Bil\) to \(\QF^{\qdr}_{\Bil}\) and \(\QF^{\sym}_{\Bil}\) respectively.
\end{corollary}

\begin{remark}
\label{remark:stable-recollement}%
By Corollary~\refone{corollary:hom-universal} with Lemma~\refone{lemma:linear} the pair of fully-faithful inclusions
\[
\Funs(\C) \xrightarrow{\QF^{\sym}_{(-)}} \Funq(\C) \leftarrow \Funx(\C\op,\Spa)
\]
form a \defi{recollement} in the sense of~\cite[Definition A.8.1]{HA}, and more precisely a \defi{stable recollement} in the sense of~\cite{barwick2016note} and~\cite{ShahQuigley} since all \(\infty\)-categories involved are stable and all functors involved are exact.
\end{remark}

Given a symmetric bilinear functor \(\Bil \in \Funs(\C)\), Lemma~\refone{lemma:homogeneous} implies that the linear part of the quadratic functor \(\QF^{\qdr}_{\Bil}\) is trivial. This however need not be the case for the quadratic functor \(\QF^{\sym}_{\Bil}\). To identify the linear part of the latter, recall that the symmetric bilinear part of \(\QF^{\qdr}_{\Bil}\) is canonically identified with \(\Bil\) itself (see Example~\refone{example:quadratic-symmetric}), and so by Corollary~\refone{corollary:hom-universal} natural transformations \(\QF^{\qdr}_{\Bil} \Rightarrow \QF^{\sym}_{\Bil}\) correspond to natural transformations \(\Bil \Rightarrow \Bil\). In particular, there is a distinguished transformation
\begin{equation}
\label{equation:orbit-to-fixed-points}%
\QF^{\qdr}_{\Bil} \Rightarrow \QF^{\sym}_{\Bil}
\end{equation}
which corresponds to the identity \(\Bil \Rightarrow \Bil\). In terms of the adjunctions of Corollary~\refone{corollary:hom-universal}, this map can also be identified with the counit of the adjunction \((-)^{\diag}_{\hC} \dashv \Bil_{(-)}\) evaluated at \(\QF^{\sym}_{\Bil}\), or the unit of \(\Bil_{(-)} \dashv ((-)^{\diag})^{\hC}\) evaluated at \(\QF^{\qdr}_{\Bil}\).

\begin{lemma}
\label{lemma:norm}%
For \(\Bil \in \Funs(\C)\) the map~\eqrefone{equation:orbit-to-fixed-points} is canonically equivalent to the \emph{trace map} associated to the \(\Ct\)-action on the object \(\Bil^{\diag} \in \Funq(\C)\).
\end{lemma}
\begin{proof}
Given the bijective correspondence between natural transformations \(\QF^{\qdr}_{\Bil} \Rightarrow \QF^{\sym}_{\Bil}\) and natural transformations \(\Bil \Rightarrow \Bil\) it will suffice to construct an identification between the map \(\Bil \Rightarrow \Bil\) induced by the trace map of \(\Bil^{\diag}\)
and the identity on \(\Bil\). For this, note that since the functor \(\Bil_{(-)}\colon \Funq(\C) \to \Funs(\C)\) preserves all limits and colimits it also sends trace maps to trace maps. In particular, the map \(\Bil \Rightarrow \Bil\) induced on bilinear parts by the trace map \(\QF^{\qdr}_{\Bil} \Rightarrow \QF^{\sym}_{\Bil}\) is itself the trace map
\[
(\Bil \oplus \Bil)_{\hC} \Rightarrow (\Bil \oplus \Bil)^{\hC}
\]
associated to the \(\Ct\)-action on the symmetric bilinear part of \(\Bil^{\Del}\), which we identify with the induced/coinduced \(\Ct\)-object \(\Bil \oplus \Bil\) as in Example~\refone{example:quadratic-symmetric}. The desired result now follows from the following completely general property of trace maps: given a semi-additive \(\infty\)-category \(\D\) and an object \(\x \in \D\), the trace map of the induced/coinduced \(\Ct\)-object \(\x \oplus \x\) identifies with the identity \(\id\colon\x \to \x\) under the canonical identifications \((\x\oplus \x)_{\hC} \simeq \x \simeq (\x \oplus \x)^{\hC}\).
\end{proof}

\begin{remark}
\label{remark:goodwillie-terms}%
Lemma~\refone{lemma:norm} implies that \(\QF\) naturally lifts to a functor with values in \emph{genuine \(\Ct\)-spectra}. We will discuss this issue in greater detail and precision in \S\refone{subsection:mackey-functors}.
\end{remark}

\begin{corollary}
\label{corollary:tate}%
For a symmetric bilinear functor \(\Bil \in \Funs(\C)\) the linear part of \(\QF^{\sym}_{\Bil}(\x) = \Bil(\x,\x)^{\hC}\) is naturally equivalent to the Tate construction \((\Bil^{\diag})^{\tC}(\x) = \Bil(\x,\x)^{\tC}\). In particular, the latter is always an exact functor.
\end{corollary}

By virtue of Lemma~\refone{lemma:norm} and Corollary~\refone{corollary:tate}, any quadratic functor \(\QF\) on \(\C\) determines a diagram of quadratic functors
\begin{equation}
\label{equation:QF-tate}%
\begin{tikzcd}
\Bil_{\QF}(\x,\x)_{\hC} \ar[r] \ar[d,equal] & \QF(\x) \ar[r] \ar[d] & \Lin_\QF(\x) \ar[d] \\
\Bil_{\QF}(\x,\x)_{\hC} \ar[r] & \Bil_\QF(\x,\x)^\hC \ar[r] & \Bil_\QF(\x,\x)^\tC
\end{tikzcd}
\end{equation}
in which the right square is exact and the right most vertical map is obtained from the middle vertical map by taking linear parts.
Conversely, by Proposition~\refone{proposition:basic-properties-quad-functors}
a symmetric bilinear functor \(\Bil \colon \C\op \times \C\op \rightarrow \Spa\), an exact functor \(\Lin \colon \C\op \rightarrow \Spa\) and a natural transformation \(\tau \colon \Lin \Longrightarrow (\Bil^{\diag})^\tC\) together determine a quadratic functor \(\QF \colon \C\op \rightarrow \Spa\) by declaring the square
\[
\begin{tikzcd}
\QF(\x) \ar[r] \ar[d] & \Lin(\x) \ar[d,"{\tau_{\x}}"] \\
\Bil(\x,\x)^\hC \ar[r] & \Bil(\x,\x)^\tC
\end{tikzcd}
\]
cartesian. This observation leads to a well-known classification of quadratic functors, which we now explain. To formulate it, let us first note that for a quadratic functor \(\QF\), Lemma~\refone{lemma:linear} tells us that the natural transformation \(\QF \Rightarrow \Lin_{\QF}\) is universally characterized by the property that it induces an equivalence on mapping spectra to every exact functor. In particular, if \(\vphi\colon\QF \Rightarrow \Lin\) is any map from \(\QF\) to an exact functor \(\Lin\) which induces an equivalence on mapping spectra to any exact functor, then \(\vphi\) factors through an equivalence \(\QF \to \Lin_{\QF} \xrightarrow{\simeq} \Lin\) in an essentially unique manner. In this case, we will also say that \(\vphi\) \defi{exhibits} \(\Lin\) as the linear part of \(\QF\). Similarly, we will say that a map \(\psi\colon \QF \to \QF'\) exhibits \(\QF'\) as the \defi{cohomogeneous part} of \(\QF\) if \(\QF'\) is cohomogeneous and \(\psi\) induces an equivalence on mapping spectra to any cohomogeneous functor. In this case, Corollary~\refone{corollary:hom-universal} tells us that \(\psi\) factors through an essentially unique equivalence \(\QF(\x) \to \Bil_{\QF}(\x,\x)^{\hC} \xrightarrow{\simeq} \QF'(\x)\). Let us now denote by \(\E \subseteq \Fun(\Del^1 \times \Del^1,\Funq(\C))\) the full subcategory spanned by those squares of quadratic functors
\begin{equation}
\label{equation:qf-tate-generic}%
\begin{tikzcd}
\QF \ar[r]\ar[d] & \Lin \ar[d] \\
\QF' \ar[r] & \Lin'
\end{tikzcd}
\end{equation}
which are exact and for which the top horizontal map exhibits \(\Lin\) as the linear part of \(\QF\) and the left vertical map exhibits \(\QF'\) as the cohomogeneous part of \(\QF\). In particular, any square in \(\E\) is equivalent to a square as of the form appearing on the right side of~\eqrefone{equation:QF-tate} in an essentially unique way. We may consider \(\E\) as the \(\infty\)-category of quadratic functors equipped with a ``cohomogeneous-linear decomposition''.

\begin{proposition}
\label{proposition:classification}%
The evaluation at \((0,0) \in \Del^1 \times \Del^1\) map \(\E \to \Funq(\C)\) sending a square as in~\eqrefone{equation:qf-tate-generic} to \(\QF\) is an equivalence of \(\infty\)-categories. In particular, every quadratic functor can be written as a pullback of cohomogeneous and exact functors.
\end{proposition}
\begin{proof}
Given Remark~\refone{remark:stable-recollement} this can be deduced from~\cite[Lemma 9]{barwick2016note}. We however spell out the details for completeness.
Let \(\E\ulcorner \subseteq \Fun(\Lam^2_0,\Funq(\C))\) be the full subcategory spanned by those \(\Lam^2_0\)-diagrams
\begin{equation}
\label{equation:qf-tate-generic-2}%
\begin{tikzcd}
\QF \ar[r]\ar[d] & \Lin \\
\QF'  &
\end{tikzcd}
\end{equation}
for which the top horizontal map exhibits \(\Lin\) as the linear part of \(\QF\) and the left vertical map exhibits \(\QF'\) the cohomogeneous part of \(\QF\). Then the restriction of any square in \(\E\) to \(\Lam^2_0 \subseteq \Del^1 \times \Del^1\) lies in \(\E\ulcorner\) and the resulting projection \(\E \to \E\ulcorner\) is a trivial Kan fibration since every square in \(\E\) is exact and hence a left Kan extension of its restriction to \(\Lam^2_0\). It will hence suffice to show that the projection \(\E\ulcorner \to \Funq(\C)\) sending a diagram as in~\eqrefone{equation:qf-tate-generic-2} to \(\QF\) is an equivalence. Now the \(\infty\)-category \(\E\ulcorner\) can be embedded in the larger \(\infty\)-category \(\E' \subseteq \Fun(\Lam^2_0,\Funq(\C))\) consisting of those diagrams as in~\eqrefone{equation:qf-tate-generic-2} for which \(\Lin\) is exact and \(\QF'\) is cohomogeneous. Then the projection \(\E' \to \Funq(\C)\) sending~\eqrefone{equation:qf-tate-generic-2} to \(\QF\) is a cartesian fibration classified by the functor sending \(\QF\) to the product of the comma category of exact functors under \(\QF\) and the comma category of cohomogeneous functors under \(\QF\). We may then identify \(\E\ulcorner\) with the full subcategory of \(\E'\) spanned by those objects which are initial in their fibres. It then follows that the projection \(\E\ulcorner \to \Funq(\C)\) is a trivial Kan fibration, and so the proof is complete.
\end{proof}

We may now deduce the classification theorem for hermitian structures (cf.\ the general classification of recollements~\cite[Proposition A.8.11]{HA}):
\begin{corollary}[Classification of hermitian structures]
\label{corollary:classification-of-quad-functors}%
The square
\begin{equation}
\label{equation:classification}%
\begin{tikzcd}
\Funq(\C) \ar[r,"\tau"] \ar[d,"\Bil"'] & \Ar(\Funx(\C\op,\Spa)) \ar[d,"{\target}"] \\
\Funs(\C) \ar[r] & \Funx(\C\op,\Spa),
\end{tikzcd}
\end{equation}
is cartesian. Here the lower horizontal functor sends \(\Bil\) to \((\Bil^{\diag})^\tC\) and the right vertical functor sends an arrow to its target.
\end{corollary}
\begin{proof}
Let \(\E\lrcorner \subseteq \Fun(\Lam^2_2,\Funq(\C))\) be the full subcategory spanned by those \(\Lam^2_2\)-diagrams
\begin{equation}
\label{equation:qf-tate-generic-3}%
\begin{tikzcd}
 & \Lin\ar[d] \\
\QF' \ar[r] & \Lin'
\end{tikzcd}
\end{equation}
for which \(\QF'\) is cohomogeneous, \(\Lin\) is exact and the bottom horizontal map exhibits \(\Lin'\) as the linear part of \(\QF'\). Then restriction along \(\Lam^2_2 \subseteq \Del^1 \times \Del^1\) sends every square in \(\E\) to a square in \(\E\lrcorner\). On the other hand, if we complete a diagram of the form~\eqrefone{equation:qf-tate-generic-3} which belongs to \(\E\lrcorner\) to a cartesian square, then this square will belong to \(\E\): indeed, this follows from the fact that a map \(\QF'' \to \Lin''\) from a quadratic to an exact functor exhibits the latter as the linear part of the former if and only if its fibre maps trivially to any exact functor, that is, if its fibre is homogeneous. We then conclude that the projection \(\E \to \E\lrcorner\) induced by restriction along \(\Lam^2_2\) is a trivial Kan fibration. On the other hand, the \(\infty\)-category \(\E\lrcorner\) is by construction a fibre product \(\Ar(\Funx(\C)) \times_{\Funx(\C)} \E''\) where \(\E''\) is the full subcategory of \(\Fun(\Del^1,\Funq(\C))\) spanned by those arrows \(\psi\colon\QF' \to \Lin'\) such that \(\QF'\) is cohomogeneous and \(\psi\) exhibits \(\Lin'\) as the linear part of \(\QF'\). As in the proof of Proposition~\refone{proposition:classification} the projection \(\E'' \to \Funqcoh(\C)\) sending \(\QF'' \to \Lin''\) to \(\QF''\) is a trivial Kan fibration onto the full subcategory \(\Funqcoh(\C) \subseteq \Funq(\C)\) spanned by the cohomogeneous functors, and the section is given by sending a cohomogeneous functor \(\QF'\) to the arrow \(\QF'\to \Lin_{\QF'}\). We hence see that the projection \(\E'' \to \Funx(\C)\) is equivalent as an arrow to the functor \(\Funqcoh(\C) \to \Funx(\C)\) taking linear parts. Finally, by Corollary~\refone{corollary:hom-universal} and Corollary~\refone{corollary:tate} the latter arrow is also equivalent to the arrow \(\Funb(\C)^{\hC} \to \Funx(\C)\) sending \(\Bil\) to \((\Bil^{\diag})^{\tC}\). Since \(\E'' \to \Funx(\C)\) is a categorical fibration the fibre product \(\E\lrcorner\) is a model for the homotopy fibre product in the square~\eqrefone{equation:classification}. The desired result now follows from Proposition~\refone{proposition:classification} and the fact that the projection \(\E \to \E\lrcorner\) is an equivalence.
\end{proof}

Finally, let us also deduce an analogous classification for Poincaré structures. For this, let us denote by \(\Funpb(\C) \subseteq \Funb(\C)\) the non-full subcategory spanned by the perfect bilinear functors and duality preserving natural transformations, that is, the natural transformations \(\beta\colon \Bil \Rightarrow \Biltwo\) for which the associated transformation \(\tau_{\beta}\colon \Dual_{\Bil} \Rightarrow \Dual_{\Biltwo}\)
is an equivalence. We then define \(\Funps(\C)\) to be the \(\infty\)-category sitting in the pullback square
\[
\begin{tikzcd}
\Funps(\C) \ar[r]\ar[d] & \Funs(\C)\ar[d] \\
\Funpb(\C) \ar[r] & \Funb(\C) \ .
\end{tikzcd}
\]
It then follows directly from the definitions that the subcategory inclusion \(\Funp(\C)\subseteq \Funq(\C)\) (see Definition~\refone{definition:poinc-cats}) features in a commutative diagram
\[
\begin{tikzcd}
\Funp(\C) \ar[r]\ar[d] & \Funq(\C) \ar[d] \\
\Funps(\C) \ar[r]\ar[d] & \Funs(\C) \ar[d] \\
\Funpb(\C) \ar[r] & \Funb(\C) \ .
\end{tikzcd}
\]
in which both squares are pullback squares.
The following is now a direct consequence of Corollary~\refone{corollary:classification-of-quad-functors}:

\begin{corollary}[Classification of Poincaré structures]
\label{corollary:classification-of-poinc-functors}%
The square
\begin{equation}
\label{equation:classification-poincare}%
\begin{tikzcd}
\Funp(\C) \ar[r,"\tau"] \ar[d,"\Bil"'] & \Ar(\Funx(\C\op,\Spa)) \ar[d,"\target"] \\
\Funps(\C) \ar[r]    & \Funx(\C\op,\Spa),
\end{tikzcd}
\end{equation}
is cartesian.
\end{corollary}

\subsection{Functoriality of hermitian structures}
\label{subsection:functoriality}%

In this subsection, we discuss the functorial dependence of \(\Funq(\C)\) on \(\C\) from the perspective of the classification described in \S\refone{subsection:classification}, not only contravariantly via restriction along exact functors, but also covariantly via left Kan extensions.
Recall that in \S\refone{subsection:hermitian-and-poincare-cats} we defined \(\Cath\) as the total \(\infty\)-category of the cartesian fibration
\begin{equation}
\label{equation:forgetful}%
\Cath \to \Catx
\end{equation}
which classifies the functor \(\C \mapsto \Funq(\C)\). In particular, being a cartesian fibration, the projection~\eqrefone{equation:forgetful} encodes the contravariance dependence of \(\Funq(\C)\) in \(\C\). We will now show that \(\Funq(\C)\) also depends \emph{covariantly} in \(\C\) via the formations of left Kan extensions. In particular, it will follow that the projection~\eqrefone{equation:forgetful} is also a \emph{cocartesian} fibration. Since we constantly work with contravariant functors to spectra let us employ the following notation: given a functor \(g\colon \D \to \E\) between \(\infty\)-categories we denote by \(g_!\colon \Fun(\D\op,\Spa) \to \Fun(\E\op,\Spa)\) the operation of left Kan extension along \(g\op\colon \D\op \to \E\op\).

\begin{lemma}
\label{lemma:kan-extension-exact-quadratic}%
\
\begin{enumerate}
\item
\label{item:kan-preserves-exact}%
If \(f\colon \C \to \D\) is an exact functor between stable \(\infty\)-categories then the associated left Kan extension functor
\[
f_!\colon \Fun(\C\op,\Spa) \to \Fun(\D\op,\Spa)
\]
sends exact functors to exact functors.
\item
\label{item:kan-preserves-biexact}%
If \(f\colon \C \to \D\) and \(g\colon \A \to \cB\) are exact functors between stable \(\infty\)-categories then the associated left Kan extension functor
\[
(f \times g)_!\colon \Fun(\C\op \times \A\op,\Spa) \to \Fun(\D\op \times \cB\op,\Spa)
\]
sends bilinear functors to bilinear functors.
\item
\label{item:kan-preserves-quadratic}%
If \(f\colon \C \to \D\) is an exact functor between stable \(\infty\)-categories then the associated left Kan extension functor
\[
f_!\colon \Fun(\C\op,\Spa) \to \Fun(\D\op,\Spa)
\]
sends quadratic functors to quadratic functors.
\end{enumerate}
\end{lemma}

\begin{proof}
We first note that left Kan extension along any functor between pointed \(\infty\)-categories preserve reduced functors by the pointwise formula for left Kan extension. Let now \(f\colon \C \to \D\) be an exact functor and \(\RF\colon \C\op \to \Spa\) a reduced functor. Consider the following commutative diagram of stable \(\infty\)-categories
\[
\begin{tikzcd}
\C\ar[r,"f"]\ar[d,"i"] & \D\ar[d,"j"]\\
\Pro(\C) \ar[r,"\wtl{f}"] & \Pro(\D)
\end{tikzcd}
\]
Then the bottom arrow admits a left adjoint \(g \colon \Pro(\D) \lrar \Pro(\C)\), which corresponds to restriction along \(\C \to \D\) under the identification of \(\Pro(-) \simeq \Ind((-)\op)\op\) with (the opposite category of) right exact functors to spaces. Now since the Yoneda embedding \(j\colon \D \to \Pro(\D)\) is fully faithful we have that \(j^*j_!\colon \Fun(\D\op,\Spa) \to \Fun(\D\op,\Spa)\) is equivalent to the identity, and so
\[
f_!\RF \simeq j^*j_!f_!\RF\simeq j^*\wtl{f}_!i_!\RF\,.
\]
Moreover, since \(g\) is left adjoint to \(\wtl{f}\) we have that the left Kan extension functor \(\wtl{f}_!\colon \Funx(\Pro(\C)\op,\Spa) \to \Funx(\Pro(\D)\op,\Spa)\) is equivalent to restriction along \(g\op \colon \Pro(\D)\op \to \Pro(\C)\op\), and so
\[
f_!\RF\simeq j^*g^*i_!\RF \simeq (gj)^*i_!\RF\,.
\]
Applying~\cite[Proposition~6.1.5.4]{HA} we now get that \(i_!\RF\) is exact (resp.\ quadratic) if \(\RF\) is exact (resp.\ quadratic). Precomposition with the exact functor \(gj\) then preserve the properties of being exact or quadratic, and so \(f_!\RF\) is exact (resp.\ quadratic) if \(\RF\) is exact (resp.\ quadratic). This proves \refoneitem{item:kan-preserves-exact} and \refoneitem{item:kan-preserves-quadratic}. To prove \refoneitem{item:kan-preserves-biexact}, we now argue as follows. By the compatibility of left Kan extensions with composition of functors we may reduce to the case where either \(f\) or \(g\) are the identity functor. By symmetry it will suffice to assume that it is \(f\) which is the identity \(\C \to \C\). For any functor \(\RF\colon \C\op \times \A\op \to \Spa\) and every \(\x \in \C\) we then have
\[
((\id \times g)_!\RF)|_{\{x\} \times \cB\op} \simeq g_!(\RF|_{\{x\} \times \A}),
\]
as can be seen by the pointwise formula for left Kan extension. We may then conclude that under the exponential equivalences
\[
\Fun(\C\op \times \A\op,\Spa) \simeq \Fun(\C\op,\Fun(\A\op,\Spa)) \quad\text{and}\quad \Fun(\C\op \times \cB\op,\Spa) \simeq \Fun(\C\op,\Fun(\cB\op,\Spa))
\]
the left Kan extension functor \((\id \times g)_!\) corresponds to post-composing with the left Kan extension functor \(g_!\colon\Fun(\A\op,\Spa) \to \Fun(\cB\op,\Spa)\). Under the same equivalence the bilinear functors correspond to those functor \(\C\op \to \Fun(\A\op,\Spa)\) which are exact and which take values in \(\Funx(\A\op,\Spa) \subseteq \Fun(\A\op,\Spa)\). Since \(g_!\) preserve exact functors by the first part of the lemma and post-composition with \(g_!\) preserves exact functors since \(g_!\) is colimit preserving (being a left adjoint), it now follows that \((\id \times g)_!\) preserves bilinear functors, as desired.
\end{proof}

\begin{corollary}
\label{corollary:cocartesian}%
The projection
\[
\Cath\lrar \Catx
\]
is also a cocartesian fibration, with pushforward along \(f \colon \C\lrar \D\) given by \(\QF\mapsto f_!\QF\).
\end{corollary}

Let us now discuss the compatibility of restriction and left Kan extensions with the decomposition of the \(\infty\)-category of quadratic functors give by Corollary~\refone{corollary:classification-of-quad-functors}.
We first observe that, given an exact functor \(f\colon \C \to \D\), the associated restriction functor \(f^*\colon \Funq(\D) \to \Funq(\C)\) respects the square~\eqrefone{equation:classification} in its entirety: indeed, taking linear and bilinear parts is compatible with restriction by Remarks~\refone{remark:invariance-base-change} and~\refone{remark:base-change-linear}, and the bottom functor in~\eqrefone{equation:classification} is also visibly compatible with restriction. We then get that if \(\QF\) is quadratic functor on \(\D\) with bilinear part \(\Bil\), linear part \(\Lin\) and structure map
\(\alp\colon \Lin \to [\Bil^{\diag}]^{\tC}\), then \(f^*\QF\) is the quadratic functor with bilinear part \((f \times f)^*\Bil\), linear part \(f^*\Lin\) and structure map
\[
f^*\alp\colon f^*\Lin \to f^*[\Bil^{\diag}]^{\tC} \simeq [((f \times f)^*\Bil)^{\diag}]^{\tC} .
\]
We now give a similar statement for left Kan extensions:

\begin{proposition}
\label{proposition:left-kan-bilinear-linear}%
Let \(f\colon\C\lrar\D\) be an exact functor between stable \(\infty\)-categories and let \(\QF\in\Funq(\C)\) be a quadratic functor on \(\C\).
Then the natural transformations
\begin{equation}
\label{equation:left-kan-bilin}%
(f \times f)_!\Bil_{\QF} \Rightarrow \Bil_{f_!\QF} \quad\quad f_!\Lin_{\QF} \Rightarrow \Lin_{f_!\QF}
\end{equation}
and
\[
f_!(\Omega^\infty \QF)\lrar\Omega^\infty(f_!\QF)
\]
are equivalences.
\end{proposition}
\begin{proof}
Let \(i\colon \C \to \Pro(\C)\) and \(j\colon \D \to \Pro(\D)\) be the respective Yoneda embeddings. Arguing as in the proof of Lemma~\refone{lemma:kan-extension-exact-quadratic} using that \((j \times j)^*(j \times j)_!\) is equivalent to the identity and that restriction commutes with taking bilinear parts (Remark~\refone{remark:invariance-base-change})
we may identify the first map in~\eqrefone{equation:left-kan-bilin} with the restriction along \(gj \times gj\colon \D \times \D \to \Pro(\C) \times \Pro(\C)\) of
\begin{equation}
\label{equation:case-pro-C}%
(i \times i)_!\colon\Bil_{\QF} \Rightarrow \Bil_{i_!\QF} .
\end{equation}
We may hence assume without loss of generality that \(\D=\Pro(\C)\) and \(f = i\).
To prove the latter special case, we see that the component of the transformation~\eqrefone{equation:case-pro-C} at a pair of pro-objects \(\{\x_{\alp}\}_{\alp\in \I},\{\y_\beta\}_{\beta \in \J}\) in \(\C\) identifies with the natural map
\[
\displaystyle\mathop{\colim}_{(\alp,\beta) \in \I\op \times \J\op} \fib[\QF(\x_\alp,\y_{\beta}) \to \QF(x_\alp)\oplus \QF(y_\beta)] \to \fib\left[\displaystyle\mathop{\colim}_{(\alp,\beta) \in \I\op \times \J\op} \QF(\x_{\alp} \oplus \y_{\beta}) \to \displaystyle\mathop{\colim}_{\alp \in \I}\QF(\x_\alp) \oplus \displaystyle\mathop{\colim}_{\beta \in \J}\QF(\y_\beta)\right].
\]
Now since \(\I\) and \(\J\) are cofiltered the projections \(\I\op \times \J\op \to \I\op\) and \(\I\op \times \J\op \to \J\op\) are cofinal and hence we can also rewrite the above map as
\[
\displaystyle\mathop{\colim}_{(\alp,\beta) \in \I\op \times \J\op} \fib[\QF(\x_\alp,\y_{\beta}) \to \QF(x_\alp)\oplus \QF(y_\beta)]
\]
\[
\to \fib\left[\displaystyle\mathop{\colim}_{(\alp,\beta) \in \I\op \times \J\op} \QF(\x_{\alp} \oplus \y_{\beta}) \to \displaystyle\mathop{\colim}_{(\alp,\beta) \in \I\op \times \J\op}\QF(\x_\alp) \oplus \displaystyle\mathop{\colim}_{(\alp,\beta) \in \I\op \times \J\op}\QF(\y_\beta)\right]
\]
and so the desired result follows from the commutation of finite limits and filtered colimits in \(\Spa\). This also implies that the second map in~\eqrefone{equation:left-kan-bilin} is an equivalence
since the formation of linear parts is obtained by \(\Lin_{\QF}(\x) := \cof[\Bil_{\QF}(\x,\x)_{\hC} \to \QF(\x)]\) and left Kan extension commutes with colimits. Finally, the proof that the map
\[
f_!(\Omega^\infty \QF)\lrar\Omega^\infty(f_!\QF)
\]
is an equivalence is obtained via the same argument by reducing to the case of \(\D=\Pro(\C)\) and using that the formation of infinite loop spaces commutes with filtered colimits.
\end{proof}

Proposition~\refone{proposition:left-kan-bilinear-linear} tells us that the formation of linear and bilinear parts is compatible with left Kan extensions.
The situation is however slightly less simple then with restriction, since the bottom arrow in~\eqrefone{equation:classification} does not commute with left Kan extensions. This is essentially due to the fact that the formation of symmetric hermitian structures
\(\Bil \mapsto \QF_{\Bil}^{\sym} = \Bil(\x,\x)^{\hC}\) does not commute with left Kan extension. Instead, given \(\Bil\in \Funs(\C)\) we have a natural map
\[
f_!\QF^{\sym}_{\Bil} \to \QF^{\sym}_{(f\times f)_!\Bil}
\]
which is generally not an equivalence. This leads to the following description of the behavior of structure maps under left Kan extensions:

\begin{corollary}
\label{corollary:structure-map-left-kan}%
Let \(f\colon \C \to \D\) be an exact functor. If \(\QF
\colon \C\op \to \Spa\) is quadratic functor on \(\C\) with bilinear part \(\Bil\), linear part \(\Lin\) and structure map
\(\alp\colon \Lin \to [\Bil^{\diag}]^{\tC}\) then \(f_!\QF\) is the quadratic functor on \(\D\) with bilinear part \(f_!\Bil\), linear part \(f_!\Lin\), and structure map the composite
\begin{equation}
\label{equation:structure-map-left-kan}%
f_!\Lin \to f_![\Bil^{\diag}]^{\tC} \to [((f \times f)_!\Bil)^{\diag}]^{\tC},
\end{equation}
which is the map induced on linear parts by the composite
\[
f_!\QF \to f_!\QF^{\sym}_{\Bil} \to \QF^{\sym}_{(f \times f)_!\Bil} .
\]
\end{corollary}

\begin{remark}
In the situation of Corollary~\refone{corollary:structure-map-left-kan} we can also identify the second map in~\eqrefone{equation:structure-map-left-kan} with the Beck-Chevalley transformation on the lax commuting square on the right
\begin{equation}
\label{equation:tate-restriction}%
\begin{tikzcd}
[row sep=5ex,column sep=8ex]
\Funs(\D) \ar[r,"{(-)^{\tC}_{\Del}}"] \ar[d] & \Funx(\D\op,\Spa) \ar[d] \\
\Funs(\C) \ar[r,"{(-)^{\tC}_{\Del}}"] & \Funx(\C\op,\Spa)
\end{tikzcd}
\quad\quad\quad
\begin{tikzcd}
[row sep=5ex,column sep=8ex]
\Funs(\C) \ar[r,"{(-)^{\tC}_{\Del}}"] \ar[d] & \Funx(\C\op,\Spa) \ar[d] \ar[dl,Rightarrow] \\
\Funs(\D) \ar[r,"{(-)^{\tC}_{\Del}}",near end] & \Funx(\D\op,\Spa)
\end{tikzcd}
\end{equation}
which is obtained from the commuting square on the left
by replacing the vertical restriction functors \((f \times f)^*\) and \(f^*\) by their left adjoints \((f \times f)_!\) and \(f_!\), respectively.
\end{remark}

%% file: PoincareObjects.tex
In this section we introduce and study another key element of the present paper, the notion of a \defi{Poincaré object} in a given Poincaré \(\infty\)-category \((\C,\QF)\). As reflected by Examples~\refone{example:perfect-derived} and~\refone{example:anti-symmetric}, we think of a Poincaré structure on a given stable \(\infty\)-category \(\C\) as a way of encoding a particular notion of hermitian form, e.g., quadratic, symmetric, or anti-symmetric forms on modules over rings. In the context of a general Poincaré \(\infty\)-category \((\C,\QF)\) and an object \(\x \in \C\), we will consequently call points in the underlying infinite loop space \(\Om^{\infty}\QF(\x)\) \defi{hermitian forms} on \(\x\). Such a form determines in particular a map \(\x \to \Dual \x\) from \(\x\) to its dual, and we say that a form is \defi{Poincaré} if this map is an equivalence. A Poincaré object is then the abstract analogue of a module equipped with (some flavour of) a hermitian form which is unimodular.

We begin in \S\refone{subsection:hermitian-poincare-objects} by introducing the main definitions and establishing a few basic consequences. One of the simplest forms of Poincaré objects are the \defi{hyperbolic} ones, which are the abstract analogue of the notion of hyperbolic quadratic forms. We discuss these types of Poincaré objects in \S\refone{subsection:hyp-and-sym-poincare-objects} and see how their formation can be encoded as the action of a suitable Poincaré functor \(\Hyp(\C) \to \C\) from a certain Poincaré \(\infty\)-category \(\Hyp(\C)\) constructed from \(\C\). The Poincaré \(\infty\)-category \(\Hyp(\C)\) displays the interesting property that its Poincaré objects correspond to just objects in \(\C\), and we study it in further depth in \S\refone{section:metabolic}. We also exploit this construction in order to prove that Poincaré objects with respect to symmetric Poincaré structures (Definition~\refone{definition:symmetric-quadratic}) correspond to \(\Ct\)-fixed objects in \(\C\) (see Proposition~\refone{proposition:compare-symmetric} below). In \S\refone{subsection:metabolic-and-L} we study another important kind of Poincaré objects - the metabolic Poincaré objects. These correspond to metabolic forms in the classical sense, that is, forms which admit a Lagrangian. Similarly to the hyperbolic case we show how one can understand metabolic Poincaré objects via Poincaré objects in a certain Poincaré \(\infty\)-category \(\Met(\C,\QF)\) constructed from \((\C,\QF)\). The notion of a metabolic Poincaré objects is the main input in the definition of the \defi{\(\L\)-groups} of a given Poincaré \(\infty\)-category (see Definition~\refone{definition:L-groups} below). These are in fact the homotopy groups of the \defi{\(\L\)-theory spectrum} which was classically defined and studied in the seminal work of Ranicki~\cite{Ranickiblue}, and transported to the context of Poincaré \(\infty\)-categories by Lurie~\cite{Lurie-L-theory}. A key technique in studying \(\L\)-group is Ranicki's \defi{algebraic Thom construction}, which we present in \S\refone{subsection:algebraic-thom} in the setting of Poincaré \(\infty\)-categories, and revisit in greater depth in \S\refone{subsection:thom}.

A key role in the present series of papers is played by the \defi{Grothendieck-Witt spectrum} of a Poincaré \(\infty\)-category, an invariant we will construct using the framework of cobordism categories in \papertwo. The zero'th homotopy group of the Grothendieck-Witt spectrum, also known as the \defi{Grothendieck-Witt group}, was classically defined in the context of rings as the group completion of the groupoid of unimodular forms (with respect to orthogonal sum). We will see how to define the Grothendieck-Witt group in the abstract setting of Poincaré \(\infty\)-categories in \S\refone{subsection:GW-group}, and extract some of its basic properties. In particular, the Grothendieck-Witt group \(\GW_0(\C,\QF)\) of a Poincaré \(\infty\)-category \((\C,\QF)\) sits in an exact sequence
\[
\K_0(\C)_{\Ct} \to \GW_0(\C,\QF) \to \L_0(\C,\QF) \to 0
\]
between the zero'th \(\L\)-group of \((\C,\QF)\) and the \(\Ct\)-orbits of the algebraic \(\K\)-theory of \(\C\). This exact sequence is in fact the tail of a long exact sequence issued from a fibre sequence of spectra
\[
\K(\C)_{\hC} \to \GW(\C,\QF) \to \L(\C,\QF)
\]
that we will construct in \papertwo. The existence of the above fibre sequence in this generality is a principal novelty of our approach to hermitian \(\K\)-theory, and yields a variety of consequences we will exploit in \papertwo, \paperthree and \paperfour.

\subsection{Hermitian and Poincaré objects}
\label{subsection:hermitian-poincare-objects}%

In this section we will present the notions of hermitian and Poincaré objects and extract some of their basic properties.

\begin{definition}
\label{definition:herm-forms}%
Let \((\C,\QF)\) be a hermitian \(\infty\)-category and \(\x \in \C\) an object. By a \defi{hermitian form} on \(\x\) we will mean a point \(q\) in the space \(\Om^{\infty}\QF(\x)\). We will then refer to the pair \((\x,q)\) as a \defi{hermitian object} in \((\C,\QF)\).
Hermitian objects can be organized into an \(\infty\)-category given by the total \(\infty\)-category of the right fibration
\[ \catforms(\C,\QF) := \int_{\x \in \C} \Om^{\infty}\QF(\x) \to \C\]
classified by the functor
\(\Omega^\infty \QF \colon \C\op \longrightarrow \Sps\).
We will refer to \(\catforms(\C,\QF)\) as the \defi{\(\infty\)-category of hermitian objects} in \((\C,\QF)\). We will denote by \(\spsforms(\C,\QF) \subseteq \catforms(\C,\QF)\) the maximal subgroupoid of \(\catforms(\C,\QF)\), and refer to it as the \defi{space of hermitian objects}.
\end{definition}

\begin{lemma}
\label{lemma:forms-functorial}%
The assignment \((\C,\QF) \mapsto \catforms(\C,\QF)\) canonically extends to a functor \(\catforms \colon \Cath \rightarrow \Cat\), together with a natural transformation to the forgetful functor \((\C,\QF) \mapsto \C\), whose component for a given \((\C,\QF) \in \Cath\) is the defining right fibration \(\catforms(\C,\QF) \to \C\).
\end{lemma}

\begin{proof}
The functor \(\Catx \hrar \Cat\), together with the composed natural transformation
\[\Funq(-) \Rightarrow \Fun((-)\op,\Spa) \st{\Om^{\infty}_*}{\Longrightarrow} \Fun((-)\op,\Sps), \]
where \(\Om^{\infty}_*\) denotes post-composition with the infinite loop space functor \(\Om^{\infty}\colon \Spa \to \Sps\), together induces under unstraightening a functor
\[\Cath \rightarrow \int_{\C \in \Cat} \Fun(\C\op,\Sps) .\]
Invoking (the dual of)~\cite[Corollary A.31]{GHNfree},
we may identify the Grothendieck construction on the right as
\[ \int_{\C\in\Cat} \Fun(\C\op,\Sps) \simeq \int_{\C\in \Cat}\RFib(\C) \simeq \RFib\]
where \(\RFib(\C)\) denotes the \(\infty\)-category of right fibrations over \(\C\) and \(\RFib \subseteq \Ar(\Cat)\) is the full subcategory of the arrow category of \(\Cat\) consisting of right fibrations. The resulting functor \(\Cath \to \RFib \to \Ar(\Cat)\) then associates to a hermitian \(\infty\)-category \((\C,\QF)\) the right fibration \(\catforms(\C,\QF) \to \C\), yielding the desired functoriality.
\end{proof}

We will mostly be interested in hermitian forms which satisfy a unimodularity condition. To formulate it, we need to assume that \((\C,\QF)\) non-degenerate. In that case any hermitian object \((\x,q)\) determines a map \(q_\sharp \colon \x \rightarrow \Dual_\QF(\x)\) as the image of \(q\) under
\[\Omega^\infty \QF(\x) \longrightarrow \Om^{\infty}\Bil_\QF(\x,\x) = \Map_\C(\x,\Dual_\QF(\x)).\]

\begin{definition}
\label{definition:poinc-forms}%
We will say that a hermitian form \(q\) on \(\x \in \C\) is
\defi{Poincaré} if the associated map \(q_\sharp\colon \x \to \Dual_{\QF}(\x)\) is an equivalence. In this case we will also say that \((\x,q)\) is a \defi{Poincaré object}.
We will denote by \(\Poinc(\C,\QF) \subseteq \spsforms(\C,\QF)\) the full subgroupoid of \(\spsforms(\C,\QF)\) spanned by the Poincaré objects. We will refer to \(\Poinc(\C,\QF) \in \Sps\) as the \emph{space of Poincaré objects} in \((\C,\QF)\).
\end{definition}

\begin{remark}
Similarly to the \(\infty\)-category \(\catforms(\C,\QF)\) one could also form an \(\infty\)-category of Poincaré objects as a full subcategory of \(\catforms(\C,\QF)\). This construction is rather poorly behaved formally and will not play any role in this paper. Therefore we will only consider the \emph{space} \(\Poinc(\C,\QF)\) of Poincaré objects here. \end{remark}

\begin{lemma}
\label{lemma:poincare-functorial}%
If \((f,\eta)\colon (\C,\QF) \to (\Ctwo,\QFtwo)\) is a duality preserving hermitian functor between non-degenerate hermitian \(\infty\)-categories then the induced functor
\[f_{\ast}\colon \spsforms(\C,\QF) \to \spsforms(\Ctwo,\QFtwo)\]
perserves Poincaré objects, that is, it maps the full subgroupoid \(\Poinc(\C,\QF) \subseteq \spsforms(\C,\QF)\) to the full subgroupoid \(\Poinc(\Ctwo,\QFtwo) \subseteq \spsforms(\Ctwo,\QFtwo)\). In particular, the association \((\C,\QF) \mapsto \Poinc(\C,\QF)\) thus extends to a functor
\[ \Poinc\colon\Catp \to \Sps. \]
\end{lemma}

It is this functor \(\Poinc\) that plays a pivotal role in the rest of the paper.

\begin{proof}[Proof of Lemma~\refone{lemma:poincare-functorial}]
By Remark~\refone{remark:recovered} the natural transformation \(\eta\colon \QF \to f^*\QFtwo\) determines a commutative diagram
\[
\begin{tikzcd}
\Om^{\infty}\QF(\x) \ar[rr,"{\Om^{\infty}\eta}"] \ar[d] && \Om^{\infty}\QFtwo(f(\x)) \ar[d] \\
\Om^{\infty}\Bil_{\QF}(\x,\x) \ar[rr,"{\Om^{\infty}\beta_{\eta}}"] \ar[d,"\simeq"'] && \Om^{\infty}\Bil_{\QFtwo}(f(\x),f(\x))\ar[d,"\simeq"] \\
\Map_{\C}(\x,\Dual_{\QF}(\x)) \ar[r,"f"] & \Map_{\Ctwo}(f(\x),f\Dual_{\QF}(\x))\ar[r,"{(\tau_{\eta})_{\ast}}"] & \Map_{\Ctwo}(f(\x),\Dual_{\QFtwo}f(\x))
\end{tikzcd}
\]
In particular, if \((f,\eta)\) is duality preserving then \(\tau_{\eta}\) is an equivalence and hence the top horizontal arrow sends Poincaré forms on \(\x\) to Poincaré forms on \(f(\x)\).
\end{proof}

\begin{remark}
\label{remark:self-dual}%
The map \(\QF(\x) \to \Bil_{\QF}(\x,\x)\) factors as \(\QF(\x) \to \Bil_{\QF}(\x,\x)^{\hC} \to \Bil_{\QF}(\x,\x)\), see~\eqrefone{equation:orbits-fixed-points}. It then follows that for any hermitian form \(q\) on \(\x\) the corresponding map \(q_{\sharp}\colon \x \to \Dual_{\QF}(\x)\) is \emph{self-dual}, that is, it is invariant under the \(\Ct\)-action on \(\map(\x,\Dual_{\QF}\x) \simeq \Bil_{\QF}(\x,\x)\).
In particular, by Remark~\refone{remark:reproduce} there is a canonical homotopy rendering the diagram
\[
\begin{tikzcd}
\x  \ar[rd,"{q_\sharp}"'] \ar[rr,"{\ev_{\x}}"] &&
\Dual_\QF\Dual_\QF\op(\x) \ar[ld,"{\Dual_\QF(q_\sharp)}"]  \\
& \Dual_\QF(\x)  &
\end{tikzcd}
\]
commutative.
\end{remark}

\begin{remark}
\label{remark:forms-on-nondeg-herm}%
Every Poincaré form \(q\) on \(\x\) gives rise to a form \(\hat{q}\) on \(\Dual_\QF(\x)\) via the inverse of the induced map \((q_\sharp)^* \colon \Omega^\infty \QF(\Dual_\QF(\x)) \rightarrow \Omega^\infty\QF(\x)\). By construction \(q_\sharp \colon (\x ,q) \rightarrow (\Dual_{\QF}(\x), \hat{q})\) is an equivalence in \(\Poinc(\C,\QF)\) and from Remark~\refone{remark:self-dual} we find \(\hat{q}_\sharp  \simeq \Dual_\QF(q_\sharp)^{-1} \colon \Dual_\QF(\x) \rightarrow \Dual_\QF \Dual_\QF(\x)\). In particular, \(q_\sharp^{-1} \simeq \ev_{\x}^{-1} \circ \hat{q}_{\sharp}\).
\end{remark}

\begin{example}
\label{example:manifold}%
Let \(M\) be a compact oriented topological \(n\)-manifold with boundary \(\partial M \subseteq M\) and fundamental class \([M]\in H_n(M,\partial M)\). Then the fundamental class together with the cup-product
induces a \((-n)\)-shifted hermitian form
\[q^{[M]}\in \Om^{\infty}\QF^{\sym[-n]}_{\ZZ}(C^*(M,\partial M)) = \Map_{\ZZ}(C^*(M,\partial M)\otimes C^*(M,\partial M),\ZZ[-n])^{\hC}\]
sending \((\varphi,\psi)\) to \((\varphi \cup \psi)([M])\). We note that we are working with homological grading conventions, so that, for example, the complex \(C^*(M,\partial M)\) is concentrated in non-positive degrees with trivial homology outside the range \([-n,0]\).
The associated map \(q^{[M]}_{\sharp}\) from \(C^*(M,\partial M)\) to its dual can then be identified with the canonical map
\[ C^*(M,\partial M) \to C^*(M) ,\]
which is an equivalence if and only if \(C^*(\partial M) \simeq 0\), i.e., if and only if \(\partial M\) is empty. In particular, the hermitian form \(q^{[M]}\) is Poincaré if and only if \(M\) is closed.
\end{example}

\subsection{Hyperbolic and symmetric Poincaré objects}
\label{subsection:hyp-and-sym-poincare-objects}%

Given a Poincaré \(\infty\)-category \((\C,\QF)\), the space of Poincaré objects \(\Poinc(\C,\QF)\) is related to the underlying space of objects \(\grpcr\C\) in two different ways. First, one can of course take a Poincaré object and forget its Poincaré form, yielding a forgetful map \(\Poinc(\C,\QF) \to \grpcr\C\). There is however also an interesting construction in the other direction, which takes a object \(\x \in \C\) and associates to it the object \(\x \oplus \Dual \x\) endowed with its \defi{hyperbolic} Poincaré form, leading to a map \(\grpcr\C \to \Poinc(\C)\). Though these constructions seem different in nature, they are in fact closely related, and will both occupy our attention in this present section. A common feature they both share is equivariance with respect to the \(\Ct\)-action on \(\grpcr\C\) induced by the duality. In the final part of this section we will show that when the Poincaré structure is symmetric the resulting map \(\Poinc(\C,\QF) \to \grpcr\C^{\hC}\) is an equivalence. We will further study the relationship between the hyperbolic and forgetful functors in \S\refone{subsection:mackey-functors} in the setting of \(\Ct\)-categories.

\begin{definition}
\label{definition:hyp-form}%
Let \((\C,\QF)\) be a Poincaré \(\infty\)-category with duality \(\Dual\). Given an object \(\x \in \C\) we will denote by \(\hyp(\x) \in \Poinc(\C,\QF)\) the Poincaré object whose underlying object is \(\x \oplus \Dual\x\) and whose Poincaré form is given by the
the image of the identity under
\[\Map_\C(\x,\x) \xrightarrow{(\ev_{\x})_*} \Map_\C(\x,\Dual\Dual(\x)) \simeq \Omega^\infty\Bil_\QF(\x,\Dual(\x)) \longrightarrow \Omega^\infty\QF(\x \oplus \Dual(\x)).\]
Unwinding the definitions one easily checks that this indeed defines a Poincaré object. We will refer to \(\hyp(\x)\) as the \defi{hyperbolic} Poincaré object on \(\x\).
\end{definition}

To understand systematically the role played by hyperbolic Poincaré objects in \(\C\) it is most useful to describe them as Poincaré objects in another Poincaré \(\infty\)-category built from \(\C\).

\begin{definition}
\label{definition:hyperbolic-cat}%
Let \(\C\) be a stable \(\infty\)-category. We define its \defi{hyperbolic category} \(\Hyp(\C)\) to be the hermitian \(\infty\)-category whose underlying stable \(\infty\)-category is \(\C\oplus \C\op\), equipped with the hermitian structure
\(\QF_{\hyp}(\x,\y) = \map_\C(\x,\y)\).
\end{definition}

Unwinding the definitions, we see that the symmetric bilinear functor associated to the hyperbolic hermitian structure is given by
\[ \Bil_{\hyp}((\x,\y),(\xtwo,\ytwo)) = \map_\C(\x,\ytwo)\oplus \map_\C(\xtwo,\y),\]
and its linear approximation is trivial. In particular, the bilinear functor \(\Bil_{\hyp}\) is perfect with duality
\( \Dual_{\hyp}(\x,\y) = (\y,\x) \)
and consequently \(\Hyp(\C)\) is always a Poincaré \(\infty\)-category.

\begin{remark}
\label{remark:hyp-symmetric-quadratic}%
By construction, the quadratic functor \(\QF_{\hyp}\) is obtained by diagonally restricting the bilinear functor \(((\x,y),(\xtwo,\ytwo)) \mapsto \map_{\C}(\x,\ytwo)\). It then follows that the canonical maps
\[ \Bil_{\hyp}((\x,\y),(\x,\y)) \to \QF_{\hyp}(\x,\y) \to \Bil_{\hyp}((\x,\y),(\x,\y)) \]
are given by the collapse and diagonal maps
\[ \map_{\C}(\x,\y) \oplus \map_{\C}(\x,\y) \to \map_{\C}(\x,\y) \to \map_{\C}(\x,\y) \oplus \map_{\C}(\x,\y) ,\]
and \(\QF_{\hyp}\) coincides with both the quadratic and symmetric Poincaré structure associated to the symmetric bilinear functor \(\Bil_{\hyp}\).
\end{remark}

\begin{remark}
\label{remark:hyp-shift-invariant}%
The Poincaré \(\infty\)-category \(\Hyp(\C)\) is \emph{shift-invariant}: for every \(n \in \ZZ\) the functor \(\Sig^n\times \id\colon \C \times \C\op \to \C \times \C\op\) refines to an equivalence \((\C \times \C\op,\QF_{\hyp}) \simeq \big(\C\times\C\op,\QF_{\hyp}\qshift{n}\big)\), see Definition~\refone{definition:shift}.
\end{remark}

For \((\C,\QF)\) a Poincaré \(\infty\)-category with duality \(\Dual\), the associated hyperbolic category \(\Hyp(\C)\) relates to \(\C\) via Poincaré functors
\begin{equation}
\label{equation:hyp-forget}%
\Hyp(\C) \xrightarrow{\hyp} (\C,\QF) \xrightarrow{\fgt} \Hyp(\C)
\end{equation}
in both directions. Here the functor on the left in~\eqrefone{equation:hyp-forget} is given by the exact functor \((\x,\y) \mapsto \x \oplus \Dual\y\), promoted to a hermitian functor via the natural transformation
\[\map_\C(\x,\y) \xrightarrow{(\ev_{\y})_*} \map_\C(\x,\Dual\Dual\y) \simeq \Bil_\QF(\x,\Dual\y) \longrightarrow \QF(\x \oplus \Dual\y),\]
while the second functor is given by the exact functor \(\x \mapsto (\x,\Dual\x)\), promoted to a hermitian functor via the natural transformation
\[ \QF(\x) \to \Bil_\QF(\x,\x) \to \map_{\C}(\x,\Dual\x) \simeq \QF_{\hyp}(\x,\Dual\x).\]
By definition, the \(\infty\)-category \(\catforms(\Hyp(\C))\) sits in a right fibration
\[\catforms(\Hyp(\C)) \to \C \oplus \C\op \]
classified by the functor \(\C\op\oplus \C \to \Sps\) sending \((\x,\y)\) to the mapping space \(\map_{\C}(\x,\y)\). But this functor is already known to classify the right fibration
\[\Twar(\C) \to \C \oplus \C\op \]
where \(\Twar(\C)\) is the \emph{twisted arrow category} of \(\C\) (see, e.g.,~\cite[\S 5.2.1]{HA}), and we consequently obtain an equivalence
\[ \catforms(\Hyp(\C)) \simeq \Twar(\C)\]
over \(\C \oplus \C\op\). In particular, hermitian objects in \(\Hyp(\C)\) are simply given by arrows \(\alp\colon\x \to \y\) in \(\C\), while morphisms between hermitian objects correspond to diagrams in \(\C\) of the form
\[
\begin{tikzcd}
\x \ar[d,"{\alp}"'] \ar[r] & \xtwo\ar[d,"{\alp'}"] \\
\y & \ytwo. \ar[l]
\end{tikzcd}
\]

\begin{proposition}
\label{proposition:poinc-of-hyp}%
For a stable \(\infty\)-category \(\C\), the composite
\[ \Poinc(\Hyp(\C)) \to \iota\C\oplus \iota\C\op \to \iota\C\]
is an equivalence of spaces.
Here the first map is induced by the forgetful functor \(\Catp\to\Catx\) and the second is given by the projection onto the first factor. Under this equivalence, the natural map \(\Poinc(\Hyp(\C)) \to \catforms(\Hyp(\C))\) corresponds to the map
\[ \iota\C \simeq \Twar(\iota\C) \to \Twar(\C) .\]
\end{proposition}
\begin{proof}
An object \([\alp\colon \x \to \y] \in \Twar(\C)\)
viewed as hermitian object \((\x,\y,\alp)\) in \(\Hyp(\C)\), has as associated self dual map \((\x,\y) \to \Dual(\x,\y) = (\y,\x)\) the map \(\alp\) on both factors (viewed as either a map \(\x \to \y\) in \(\C\) or a map \(\y \to \x\) in \(\C\op\)). Consequently, the hermitian form \((\x,\y,\alp)\) is Poincaré if and only if \(\alp\) is an equivalence. Together with the fact that right fibrations detect equivalences
we obtain that
\[ \Poinc(\Hyp(\C)) \simeq \Twar(\iota\C) \subseteq \iota\Twar(\C) \simeq \spsforms(\Hyp(\C)).\]
We finish the proof by observing that the projection \(\Twar(\iota\C) \to \iota\C\) is an equivalence since \(\iota\C\) is an \(\infty\)-groupoid.
\end{proof}

\begin{remark}
We will show in \S\refone{subsection:mackey-functors}
that the association \(\C \mapsto \Hyp(\C)\) organizes into a functor \(\Catx \to \Catp\) which is both left and right adjoint to the forgetful functor \(\Catp \to \Catx\), with unit and counit given by~\eqrefone{equation:hyp-forget}, see Corollary~\refone{corollary:hyp-is-adjoint}. Together with the corepresentability of \(\Poinc\) (Proposition~\refone{proposition:corepresentability-of-poinc}) this will give another proof of Proposition~\refone{proposition:poinc-of-hyp}.
\end{remark}

In light of Proposition~\refone{proposition:poinc-of-hyp} the Poincaré functors~\eqrefone{equation:hyp-forget} now induce a pair of maps
\begin{equation}
\label{equation:hyp-forget-2}%
\iota\C \to \Poinc(\C,\QF) \to \iota\C
\end{equation}
Unwinding the definitions we see that the functor on the left sends an object \(\x\) to the associated hyperbolic Poincaré object \(\hyp(\x)\), while, the functor on the right sends a Poincaré object \((\x,q)\) to the underlying object \(\x\). A key feature of both these maps is that they are \emph{\(\Ct\)-equivariant} with respect to the \(\Ct\)-action on \(\iota\C\) induced by the duality of \(\QF\) and the trivial action on \(\Poinc(\C,\QF)\). To make this idea precise will first construct this action on the level of the Poincaré \(\infty\)-category \(\Hyp(\C)\).

\begin{construction}
\label{construction:equivariant-hyp}%
Let \((\C,\QF)\) be a Poincaré \(\infty\)-category with associated duality \(\Dual = \Dual_{\QF}\). We construct a \(\Ct\)-action on \(\Hyp(\C) \in \Catp\) as follows. To begin, consider the equivalence of stable \(\infty\)-categories
\[ \id \times \Dual\op \colon \C \times \C \xrightarrow{\simeq} \C \times \C\op .\]
Transporting the flip \(\Ct\)-action on \(\C \times \C\) to \(\C \times \C\op\) we thus obtain a \(\Ct\)-action on \(\C \times \C\op\), given informally by the formula \((\x,\y) \mapsto (\Dual \y,\Dual\op \x)\). We wish to promote this action to \(\Ct\)-action on the Poincaré \(\infty\)-category \(\Hyp(\C)\). Since every equivalence in \(\Cath\) is a Poincaré functor may equivalently construct a \(\Ct\)-action on \(\Hyp(\C)\) as a hermitian \(\infty\)-category. By the construction of \(\Cath\) as the unstraightening of the functor \(\C \mapsto \Funq(\C)\), lifting the above \(\Ct\)-action on \(\C \times \C\op\) to \(\Hyp(\C)\) is equivalent to giving a \(\Ct\)-fixed point structure on \(\QF_{\hyp} \in \Funq(\C \times \C\op)\) with respect to the induced \(\Ct\)-action on \(\Funq(\C \times \C\op)\). Since the relevant \(\Ct\)-action was transported from the flip action on \(\C \times \C\) via the equivalence \((\id,\Dual\op)\) we may equivalently construct a \(\Ct\)-fixed point structure on the quadratic functor \((\id \times \Dual\op)^*\QF_{\hyp}\). The latter is readily discovered to be
\[[(\id \times \Dual\op)^*\QF_{\hyp}](\x,\y) = \map_{\C}(\x,\Dual\y) = \Bil_{\QF}(\x,\y) \]
and so we need to construct a \(\Ct\)-fixed point structure on \(\Bil_{\QF}\), considered as an object of \(\Funq(\C \times \C)\). But \(\Bil_{\QF}\) lies in the full subcategory \(\Funb(\C) \subseteq \Funq(\C \times \C)\), where it is equipped with a \(\Ct\)-fixed structure by virtue of Lemma~\refone{lemma:bilinear-symmetric}.
\end{construction}

\begin{remark}
\label{remark:action-on-core}%
The \(\Ct\)-action on \(\Hyp(\C)\) constructed in~\refone{construction:equivariant-hyp} induces a \(\Ct\)-action on \(\catforms(\Hyp(\C)) \simeq \Twar(\C)\). Unwinding the definitions, this actions sends an arrow \([\x \to \y] \in \Twar(\C)\) to the dual arrow \(\Dual_{\QF}\y \to \Dual_{\QF}\x\). Similarly, this \(\Ct\)-action determines an action on \(\Poinc(\Hyp(\C))\). Under the identification \(\Hyp(\C) \simeq \iota\C\) of Proposition~\refone{proposition:poinc-of-hyp} this action can be written simply by \(\x \mapsto \Dual_{\QF}\x\). Here we point out that since \(\iota\C\) is an \(\infty\)-groupoid it is canonically equivalent to its opposite via an equivalence which sends every arrow to its inverse. Hence the contravariant equivalence \(\Dual\) becomes a self-equivalence on the level of \(\iota\C\). We may also state this as follows: the \(\Ct\)-action \((-)\op\colon \Cat \to \Cat\) admits a canonical trivialization along the full subcategory \(\Sps \subseteq \Cat\) (in fact, the space of self-equivalences of \(\Sps\) is contractible by~ its universal property~\cite[Theorem 5.1.5.6]{HTT}), yielding an identification \(\Sps^{\hC} \simeq \Fun(\BC,\Sps)\).
The duality \(\Dual_{\QF}\) then induces a duality on \(\iota\C\) and hence a \(\Ct\)-action.
\end{remark}

\begin{lemma}
\label{lemma:fgt-hyp-equivariant}%
Let \((\C,\QF)\) be a Poincaré \(\infty\)-category with duality \(\Dual = \Dual_{\QF}\). Then the functors
\[ \Hyp(\C) \xrightarrow{\hyp} (\C,\QF) \xrightarrow{\fgt} \Hyp(\C) \]
both admit a distinguished refinement to \(\Ct\)-equivariant maps with respect to the \(\Ct\)-action on \(\Hyp(\C)\) constructed in~\refone{construction:equivariant-hyp} and the trivial action on \((\C,\QF)\).
\end{lemma}

\begin{corollary}
\label{corollary:fgt-hyp-equivariant}%
The induced maps on Poincaré objects (which we denote by the same name)
\[ \iota\C \xrightarrow{\hyp} \Poinc(\C,\QF) \xrightarrow{\fgt} \iota\C \]
are \(\Ct\)-equivariant with respect to the duality induced action on \(\iota\C\) and the trivial action on \(\Poinc(\C,\QF)\).
\end{corollary}

\begin{proof}[Proof of Lemma~\refone{lemma:fgt-hyp-equivariant}]
We first construct the \(\Ct\)-equivariant structure on the underlying exact functors. For this, note that since the \(\Ct\)-action on \(\Hyp\) was constructed by transporting the flip action along the equivalence \(\id \times \Dual\op\colon \C \times \C \to \C \times \C\op\) and \(\Dual\op\Dual \simeq \id\) it will suffice to to promote the resulting exact functors
\[ \C \times \C \to \C \to \C \times \C\]
to \(\Ct\)-equivariant exact functors. Indeed, these are just the diagonal and fold map of \(\C\) as an object in the semi-additive \(\infty\)-category \(\Catx\), which are both canonically \(\Ct\)-equivariant. To lift the resulting \(\Ct\)-equivariant structure on \((\id,\Dual\op)\colon \C \to \C \times \C\op\) to a \(\Ct\)-equivariant structure on the Poincaré functor \(\fgt\) we need to promote the associated natural transformation
\[ \QF \Rightarrow (\id,\Dual\op)^*\QF_{\hyp}\]
to a \(\Ct\)-equivariant map in \(\Funq(\C)\). Transporting the problem again along the equivalence \(\id \times \Dual\op\) we need to put a \(\Ct\)-equivariant structure on the natural transformation
\[ \QF \Rightarrow \Del^*\Bil\]
where \(\Del\colon \C \times \C\) is the diagonal. Indeed, this is established in Lemma~\refone{lemma:equivariance}. By the same argument we see that in order to obtain the desired \(\Ct\)-equivariant structure on \(\hyp\) it will suffice to put a \(\Ct\)-equivariant structure on the natural transformation
\[ \Bil \Rightarrow \nabla^*\QF ,\]
where \(\nabla\colon \C \oplus \C \to \C\) is the collapse functor \((\x,\y) \mapsto \x \oplus \y\).
Using the adjunction between restriction and left Kan extension we may instead put a \(\Ct\)-equivariant structure on the adjoint transformation
\[ \nabla_!\Bil  \simeq \Del^*\Bil \Rightarrow \QF .\]
where we used that left Kan extension along \(\nabla\op\) is obtained by restriction along its right adjoint \(\Del\op\colon \C\op \to \C\op \times \C\op\). The desired \(\Ct\)-equivariant structure was again established in Lemma~\refone{lemma:equivariance}.
\end{proof}

Let us now focus on the Poincaré functor \(\fgt\colon (\C,\QF) \to \Hyp(\C)\). Upon taking hermitian and Poincaré objects (and using Proposition~\refone{proposition:poinc-of-hyp}) this Poincaré functor induces a commutative diagram
\begin{equation}
\label{equation:poinc-forms-twar}%
\begin{tikzcd}
\Poinc(\C,\QF) \ar[r]\ar[d] & \spsforms(\C,\QF) \ar[r]\ar[d] &	\catforms(\C,\QF)\ar[d] & (\x,\qone) \ar[l,phantom,"{\ni}"] \ar[d,mapsto] \\
\iota\C \ar[r]& \iota\Twar(\C)  \ar[r] & \Twar(\C) & {[q_{\sharp}\colon \x \to \Dual_{\QF}\x]} \ar[l,phantom,"{\ni}"]
\end{tikzcd}
\end{equation}
in which the vertical maps inherit from \(\fgt\) a \(\Ct\)-equivariant structure with respect to the trivial action on their domains and the \(\Ct\)-action induced by the \(\Ct\)-action on \(\Hyp(\C)\) on the target. Here the left square consists only of spaces and the horizontal maps are (up to equivalence) inclusions of components: for the upper left map these are the components of \(\catforms(\C,\QF)\) consisting of Poincaré objects and for the lower left map these are the components of \(\iota\Twar(\C)\) consisting of those arrows \([\x \to \y]\) which are equivalences. Since by definition a hermitian object \((\x,q)\) is Poincaré if and only if \(q_{\sharp}\) is an equivalence we see in particular that the left square is cartesian. Now by the \(\Ct\)-equivariance above the external rectangle in~\eqrefone{equation:poinc-forms-twar} induces a commutative square
\begin{equation}
\label{equation:poinc-forms-twar-Ct}%
\begin{tikzcd}
\Poinc(\C,\QF) \ar[r]\ar[d] & 	\catforms(\C,\QF)\ar[d] & (\x,\qone) \ar[l,phantom,"{\ni}"] \ar[d,mapsto] \\
(\iota\C)^{\hC} \ar[r]&  \Twar(\C)^{\hC} & {[q_{\sharp}\colon \x \to \Dual_{\QF}\x]} \ar[l,phantom,"{\ni}"]
\end{tikzcd}
\end{equation}

\begin{proposition}
\label{proposition:compare-symmetric}%
If \(\QF = \QF^{\sym}_{\Bil}\) is a symmetric Poincaré structure of some symmetric bilinear form \(\Bil\) then the vertical maps in~\eqrefone{equation:poinc-forms-twar-Ct} are equivalences.
\end{proposition}
\begin{proof}
Consider the extended diagram
\begin{equation}
\label{equation:poincare-hermitian-3}%
\begin{tikzcd}
\Poinc(\C,\QF) \ar[r]\ar[d] & \spsforms(\C,\QF) \ar[r]\ar[d] & \catforms(\C,\QF)\ar[d]\ar[r] & \C \ar[d]\\
\Twar(\iota\C)^{\hC} \ar[d]\ar[r]& \iota\Twar(\C)^{\hC}  \ar[r]\ar[d] & \Twar(\C)^{\hC}\ar[d] \ar[r] & (\C \times \C\op)^{\hC}\ar[d] \\
\Twar(\iota\C) \ar[r]& \iota\Twar(\C)  \ar[r] & \Twar(\C) \ar[r] & \C \times \C\op
\end{tikzcd}
\end{equation}
in which the external rectangle of the left column is cartesian as observed above. In addition, since the fibres of the map \(\Twar(\iota\C) \to \iota\Twar(\C)\) are \((-1)\)-truncated then bottom left square is cartesian as well, and hence the top left square is cartesian. Similarly, since the map from the homotopy fixed point is conservative, the bottom central square and therefore the top central square, are cartesian.

Now the map \(\C \to \C \times \C\op\) is equivalent as a \(\Ct\)-equivariant arrow to the diagonal inclusion \(\C \to \C \times \C\), which exhibits the \(\Ct\)-object \(\C \times \C\) as coinduced from \(\C\). This implies in particular that the top right vertical map in~\eqrefone{equation:poincare-hermitian-3} is an equivalence.  Thus to conclude it suffices to show that the top right square is cartesian.

Now consider the right most column in~\eqrefone{equation:poincare-hermitian-3}.
Since homotopy fixed points commute with fibre products the fibres of the middle horizontal map over a fixed object in \(\C\) is the homotopy fixed points of the corresponding fibre of the bottom horizontal map in the same column. We hence obtain that the map on horizontal fibres in the top right square can be identified with the induced map
\[
\Om^{\infty}\QF(\x) \to \Map_{\C}(\x,\Dual\x)^{\hC} ,
\]
which is an equivalence by the assumption that \(\QF = \QF^{\sym}_{\Bil}\) for some \(\Bil\).
It then follows that the second and third vertical arrows in the top of~\eqrefone{equation:poincare-hermitian-3} are equivalences. The left most vertical map in that row is consequently an equivalence as well since the top left square is cartesian.
\end{proof}

\subsection{Metabolic objects and \(\L\)-groups}
\label{subsection:metabolic-and-L}%

In this section we will introduce the notion of a \defi{metabolic} Poincaré object and use it to define the \defi{\(\L\)-groups} of a Poincaré \(\infty\)-category. In the context of modules over rings these were first defined by Wall and Ranicki in their seminal work on surgery theory \cite{wall}, and were transported to the setting of Poincaré \(\infty\)-categories in~\cite{Lurie-L-theory}.
We will then develop an analogue of Ranicki's \emph{algebraic Thom construction} \cite[Proposition~3.4]{RanickiATS1} in this context. This construction will play an important role in the framework of \emph{algebraic surgery} which we will set up in \papertwo.

\begin{definition}
\label{definition:isotropic-and-lagrangian}%
Let \((\C,\QF)\) be a Poincaré \(\infty\)-category and \((\x,\qone)\) a Poincaré object. By an \defi{isotropic object} over \(\x\) we will mean a pair \((f\colon\cob \to \x,\eta)\) where \(f\colon \cob \to \x\) is a map in \(\C\) and \(\eta\colon f^*q \sim 0 \in \Om^{\infty}\QF(\cob)\) is a null-homotopy of the restriction of \(q\) to \(\cob\). We will say that an isotropic object \((\cob \to \x,\eta)\) is a \defi{Lagrangian} if the null-homotopy of \(\cob \to \x\simeq \Dual \x \to \Dual \cob\) given by the image of \(\eta\) in \(\Om^{\infty}\Bil_{\QF}(\cob,\cob) = \Map_{\C}(\cob,\Dual\cob)\) exhibits the sequence
\[ \cob \to \x \to \Dual \cob\]
as exact. We will say that \((\x,q)\) is \defi{metabolic} if it admits a Lagrangian.
\end{definition}

\begin{example}
If \((\C,\QF)\) is a Poincaré \(\infty\)-category and \(\x \in \C\) an object then the associated hyperbolic Poincaré object \(\hyp(\x)\) is metabolic with Lagrangian given by the component inclusion \(\x \to \x \oplus \Dual\x\).
\end{example}

\begin{example}
In \(\Dperf(\FF_2)\), let \(V\) be a 2-dimensional \(\FF_2\)-vector space with basis \(v,u\) equipped with the symmetric bilinear form \(b\colon V \otimes_{\FF_2} V \to \FF_2\) given by \(b(v,v)=0\) and \(b(v,u)=b(u,v)=b(u,u)=1\). Then the Poincaré object \((V,b) \in \Poinc(\Dperf(\FF_2),\QF^{\sym}_{\FF_2})\) is metabolic with Lagrangian \(L = \left<v\right> \hrar V\) but \((V,b)\) is not isomorphic to \(\hyp(U)\) for any \(U \in \Dperf(\FF_2)\). Indeed, since any object in \(\Dperf(\FF_2)\) breaks as a direct sum of shifts of \(\FF_2\) the only possible candidate is \(U = \FF_2\), but \((V,b)\) is not isomorphic to \(\hyp(\FF_2)\). In particular, not every metabolic object is hyperbolic.
\end{example}

\begin{example}
Let \(M\) be a closed oriented \(n\)-manifold with fundamental class \([M] \in H_n(M;\ZZ)\), so that we have a symmetric Poincaré form \(q^{[M]} \in \Om^{\infty}\QF^{s[-n]}_{\ZZ}(C^*(M))\) as in Example~\refone{example:manifold}. If \(W\) is now an oriented \((n+1)\)-manifold with boundary \(M\) then the relative fundamental class \([W] \in H_{n+1}(W,M)\) can be used to promote the map
\[ C^*(W) \to C^*(M) \]
to a Lagrangian of \((C^*(M),q^{[M]})\). This can be considered as an algebraic reflection of the fact that \(W\) exhibits \(M\) as a \emph{boundary}. In particular, if \((C^*(M),q^{[M]})\) is not metabolic then \(M\) is not the boundary of any oriented \((n+1)\)-manifold, that is, \(M\) is not (oriented-ly) null-cobordant.
\end{example}

As in the case of hyperbolic Poincaré objects, it would be desirable to have a description of metabolic Poincaré objects in terms of Poincaré objects in another Poincaré \(\infty\)-category constructed from \((\C,\QF)\).

\begin{definition}
\label{definition:metabolic-cat}%
For a Poincaré \(\infty\)-category \((\C,\QF)\), we define the associated \defi{metabolic category} \(\Met(\C,\QF)\) to be the hermitian \(\infty\)-category with underlying \(\infty\)-category \(\Ar(\C) = \Fun(\Delta^1,\C)\) and hermitian structure
\[\QF_\met\colon \Ar(\C)\op=\Ar(\C\op)\xrightarrow{\Ar(\QF)}\Ar(\Spa)\xrightarrow{\fib}\Spa
\]
whose value on arrows is \(\QF_\met([\cob \to \x]) = \fib(\QF(\x) \to \QF(\cob))\).
\end{definition}

Unwinding the definitions we see that the underlying symmetric bilinear functor of \(\QF_\met\) is
\[
\Bil_{\met}(\cob\to \x, \cobtwo\to \xtwo) = \fib[\Bil_{\QF}(\x,\xtwo) \to \Bil_{\QF}(\cob,\cobtwo)].
\]
From this formula we see that if \(\Bil_{\QF}\) is perfect with duality \(\Dual\) then \(\Bil_{\met}\) is perfect with duality
\[
\Dual_\met(\cob \to \x) = \big(\fib[\Dual\x \to \Dual\cob] \to \Dual \x \big),
\]
so that \(\Met(\C,\QF)\) is Poincaré whenever \((\C,\QF)\) is so.
We note that by definition a hermitian form on \([f\colon\cob \to \x]\) with respect to \(\QF_{\met}\) consists of a form \(q \in \Om^{\infty}\QF(\x)\) together with a null-homotopy \(\eta\) of \(f^*q \in \Om^{\infty}\QF(\cob)\). Such a \(\QF_{\met}\)-form \((q,\eta)\) is Poincaré if and only if the associated self dual map encoded by the horizontal maps of the square
\begin{equation}
\label{equation:self-dual}%
\begin{tikzcd}
\cob\ar[r,"{\eta_{\sharp}}"] \ar[d] & \fib[\Dual\x \to \Dual \cob] \ar[d] \\
\x \ar[r,"{q_{\sharp}}"] & \Dual \x
\end{tikzcd}
\end{equation}
is an equivalence. Here \(q_{\sharp}\) is the self dual map determined by \(q\) and we denoted by \(\eta_{\sharp}\) the map corresponding to the null-homotopy of the composed map \(\cob \to \x \to \Dual\x \to \Dual\cob\) determined by the image of \(\eta\) in \(\Om^{\infty}\Bil(\cob,\cob) = \Map(\cob,\Dual\cob)\). We then see that~\eqrefone{equation:self-dual} constitutes an equivalence between the vertical arrows if and only if \(q_{\sharp}\colon \x \to \Dual\x\) is an equivalence
and the resulting sequence \(\cob \to \x \simeq \Dual\x \to \Dual\cob\) is exact, that is, if \((f\colon \cob \to \x,\eta)\) is a Lagrangian. We may thus conclude that Poincaré objects in \(\Met(\C,\QF)\) correspond to metabolic objects in \((\C,\QF)\), or, more precisely, to Poincaré objects in \(\C\) equipped with a specified Lagrangian.

\begin{definition}
For a Poincaré \(\infty\)-category \((\C,\QF)\) we will denote by
\[\Poincdel(\C,\QF) := \Poinc(\Met(\C,\QF))\]
the space of Poincaré objects in \(\Met(\C,\QF)\), which we consider as above as the space of Poincaré objects equipped with a specified Lagrangian.
\end{definition}

\begin{lemma}
\label{lemma:maps-with-metabolic-category}%
The maps
\begin{equation}
\label{equation:metabolic-sequence}%
(\C,\QF\qshift{-1})\xrightarrow{i} \Met(\C,\QF)\xrightarrow{\met} (\C,\QF)
\end{equation}
given respectively by \(i(\x)=[\x\to 0]\) and \(\met([\cob\to \x])=\x\)
extend to morphisms in \(\Catp\).
\end{lemma}
\begin{proof}
For the first we observe that the map \(i\) is fully faithful, that \(i^*\QF_\met \simeq \Omega \QF\) and that the image of \(i\) is closed under the duality in \(\Met(\C,\QF)\), so the result follows from Observation~\refone{observation:poincare-subcat}.
For the second map we take the hermitian structure associated to the canonical map
\[
\QF_\met([\cob\to \x]) = \fib[\QF(\x) \to \QF(\cob)] \to \QF(\x).
\]
By the explicit description of the duality above we see that the resulting hermitian functor is Poincaré.
\end{proof}

Unwinding the definitions we see that the map
\begin{equation}
\label{equation:poinc-met-to-C}%
\Poincdel(\C,\QF) \to \Poinc(\C,\QF)
\end{equation}
induced by the right hand Poincaré functor in~\eqrefone{equation:metabolic-sequence} corresponds to the forgetful map which takes a Poincaré object equipped with a Lagrangian and forgets the Lagrangian. In particular, a Poincaré object in \((\C,\QF)\) is metabolic if and only if it is in the image of~\eqrefone{equation:poinc-met-to-C}

As observed earlier, every hyperbolic form is metabolic, but not every metabolic object is equivalent to the associated hyperbolic form on its Lagrangian. This relation between metabolic and hyperbolic objects is best expressed by relating the Poincaré \(\infty\)-categories \(\Met(\C,\QF)\) and \(\Hyp(\C)\) via suitable Poincaré functors.

\begin{construction}
\label{construction:hyp-to-met}%
Let \((\C,\QF)\) be a Poincaré \(\infty\)-category. We define Poincaré functors
\[
\begin{tikzcd}[row sep=tiny]
\Hyp(\C) \ar[r,"\ilag"] &  \Met(\C,\QF) & \Met(\C,\QF) \ar[r,"\lag"] & \Hyp(\C) \\
(\x,\y) \ar[r, mapsto] & (\x \to \x\oplus \Dual\y) & (\cob\to \x) \ar[r, mapsto] & (\cob, \Dual\cof[\cob \to \x])
\end{tikzcd}
\]
and
\[
\begin{tikzcd}[row sep=tiny]
\Hyp(\C) \ar[r,"\dilag"] &  \Met(\C,\QF) & \Met(\C,\QF) \ar[r,"\dlag"] & \Hyp(\C) \\
(\x,\y) \ar[r, mapsto] & (\Dual\y \to \x\oplus \Dual\y) & (\cob\to \x) \ar[r, mapsto] & (\cof[\cob \to \x],\Dual\cob)
\end{tikzcd}
\]
via the indicated formulas on the level on the underlying exact functors, and with hermitian structures as follows.
For the two functors on the left hand side the hermitian structure is obtained via the identification
\[ \map_{\C}(\x,\y) \simeq \Bil_{\QF}(\x,\Dual \y) \simeq \fib[\QF(\x\oplus \Dual\y) \to \QF(\x) \oplus \QF(\Dual\y)], \]
whose target visibly projects to both \(\fib[\QF(\x \oplus \Dual\y) \to \QF(\x)]\) and \(\fib[\QF(\x \oplus \Dual\y) \to \QF(\Dual\y)]\). For the functors on the right hand side the hermitian structure is given by the natural transformation
\[ \fib[\QF(\x) \to \QF(\cob)] \to \fib[\Bil_{\QF}(\x,\cob) \to \Bil_{\QF}(\cob,\cob)] \simeq \Bil_{\QF}(\cof[\cob\to\x],\cob) \]
where we recognize the target as naturally equivalent to both \(\lag^*\Bil_{\hyp}\) and \(\dlag^*\Bil_{\hyp}\).
The preservation of the duality by these hermitian functors is visible by the explicit descriptions of \(\Dual_{\met}\) and \(\Dual_{\hyp}\) above. We also note that the composites
\[ \Hyp(\C) \xrightarrow{\ilag} \Met(\C,\QF) \xrightarrow{\lag} \Hyp(\C) \]
and
\[ \Hyp(\C) \xrightarrow{\dilag} \Met(\C,\QF) \xrightarrow{\dlag} \Hyp(\C) \]
are naturally equivalent to the identity, and so exhibit \(\Hyp(\C)\) as a retract of \(\Met(\C,\QF)\) in \(\Catp\).
\end{construction}

\begin{remark}
\label{remark:lag-split-dual}%
The Poincaré functors \(\lag,\dlag\colon \Met(\C) \to \Hyp(\C)\) are closely related: they differ by post-composition with the Poincaré involution \(\Hyp(\C) \xrightarrow{\simeq} \Hyp(\C)\) of Construction~\refone{construction:equivariant-hyp}. Similarly, the Poincaré functors \(\ilag,\dilag\colon \Hyp(\C) \to \Met(\C)\) differ by pre-composition with this involution. It then follows from %
Lemma~\refone{lemma:fgt-hyp-equivariant} that the composite
\[\Met(\C,\QF) \to \Hyp(\C) \xrightarrow{\hyp} (\C,\QF)\]
is independent of whether the first functor is \(\lag\) or \(\dlag\). The action of this composed functor on \(\Poinc(-)\) sends a Poincaré object \((\x,\qone)\) equipped with a Lagrangian \(\cob \to \x\) to the associated hyperbolic object \(\hyp(\cob) \simeq \hyp(\Dual\cob)\). The difference between a metabolic Poincaré object and its hyperbolic counterpart plays a key role in the definition of the Grothendieck-Witt group, see \S\refone{subsection:GW-group} below. On the other hand, we also observe that the composite Poincaré functor
\[ \Hyp(\C) \to \Met(\C,\QF) \xrightarrow{\met} (\C,\QF)\]
coincides with the functor \(\hyp\) of~\eqrefone{equation:hyp-forget}, independently of whether the first functor is \(\ilag\) or \(\dilag\). On the level of Poincaré objects we may interpret this as the observation that a hyperbolic Poincaré object \(\hyp(\cob)\) can be considered as a metabolic object in two canonical ways: one via the Lagrangian \(\cob \to \cob \oplus \Dual\cob\) and one via the Lagrangian \(\Dual\cob \to \cob \oplus \Dual\cob\).
\end{remark}

A fundamental invariant of Poincaré \(\infty\)-categories is their \(\L\)-groups. To define them, we first observe that the set \(\pi_0\Poinc(\C,\QF)\) of equivalence classes of Poincaré objects carries a natural commutative monoid structure, with sum given by
\[
[\x,\qone] + [\xtwo,\qtwo] = [\x \oplus \xtwo,\qone\operp \qtwo]
\]
for \([\x,\qone],[\xtwo,\qtwo] \in \pi_0\Poinc(\C,\QF)\), where
\[
\qone \operp \qtwo \in \Om^{\infty}\QF(\x \oplus \xtwo) \simeq \Om^{\infty}\QF(\x) \times \Om^{\infty}\QF(\xtwo) \times \Om^{\infty}\Bil_{\QF}(\x,\xtwo)
\]
corresponds to the tuple \((\qone,\qtwo,0)\). Though this commutative monoid is generally not a group, every element is invertible up to the class of a metabolic object. More precisely, we have the following:

\begin{lemma}
\label{lemma:L-is-group}%
Let \((\C,\QF)\) be Poincaré \(\infty\)-category. Then the cokernel of the map \(\pi_0\Poincdel(\C,\QF) \to \pi_0\Poinc(\C,\QF)\) in the category of commutative monoids is a group. Explicitly, the inverse to \([\x,\qone]\) is given by \([\x,-\qone]\).
\end{lemma}
\begin{proof}
This follows from the fact that
\((\x \oplus \x,\qone\operp -\qone)\) is metabolic with Lagrangian given by the diagonal inclusion \(x \to x \oplus x\) with the canonical null-homotopy \(\qone + (-\qone) \sim 0\).
\end{proof}

\begin{definition}
\label{definition:L-groups}%
Let \((\C,\QF)\) be a Poincaré \(\infty\)-category. For \(n \in \ZZ\) we define the \(n\)'th \(\L\)-group of \((\C,\QF)\) by
\[ \L_n(\C,\QF) := \coker[\pi_0\Poincdel(\C,\QF\qshift{-n}) \to \pi_0\Poinc(\C,\QF\qshift{-n})]\]
which is an abelian group by Lemma~\refone{lemma:L-is-group}.
\end{definition}

\begin{remark}\label{remark:stably-metabolic}
A standard description of cokernels in commutative monoids gives that \([\x,\qone],[\xtwo,\qtwo] \in \pi_0\Poinc(\C,\QF\qshift{-n})\) map to the same class in \(\L_n(\C,\QF)\) if and only if there exists metabolic Poincaré objects \([\y,\pone],[\ytwo,\ptwo]\) such that \([\x,\qone] + [\y,\pone] = [\xtwo,\qtwo] + [\ytwo,\ptwo]\) in \(\pi_0\Poinc(\C,\QF\qshift{-n})\). In particular, \([\x,\qone]\) maps to zero in \(\L_n(\C,\QF)\) if and only it is stably metabolic, that is \([\x,\qone] + [\y,\pone]\) is metabolic for some metabolic Poincaré object \([\y,\pone]\).
In the setting of Poincaré \(\infty\)-categories this property is actually equivalent to \([\x,\qone]\) itself being metabolic. Indeed, suppose that \(\z \to \y\) is a Lagrangian for \(\pone\) and \(\ztwo \to \x \oplus \y\) is a Lagrangian for \(\qone \operp \pone\). Setting \(\cob = [\x \oplus \z] \times_{\x \oplus \y} \ztwo\) and using the null-homotopies of \(\pone|_{z}\) and \((\qone\operp\pone)|_{\ztwo}\) we then obtain
\[ \qone|_{\cob} \sim (\qone\operp \pone)|_{\cob} \sim 0, \]
and one can verify that this null-homotopy exhibits \(\cob\) as a Lagrangian for \(\qone\), so that \([\x,\qone]\) is metabolic. It then follows from Lemma~\refone{lemma:L-is-group} that \([\x,\qone],[\xtwo,\qtwo] \in \pi_0\Poinc(\C,\QF\qshift{-n})\) map to the same class in \(\L_n(\C,\QF)\) if and only if the Poincaré object \([\x,\qone] + [\xtwo,-\qtwo] = [\x \oplus \xtwo,\qone \operp -\qtwo]\) is metabolic.
\end{remark}

\begin{remark}
\label{remark:minus-operation}%
It follows from Lemma~\refone{lemma:L-is-group} that if \((f,\eta)\colon (\C,\QF) \to (\Ctwo,\QFtwo)\) is a Poincaré functor then the two induced abelian group homomorphisms
\[ (f,\eta)_*, (f,-\eta)_*\colon \L_n(\C,\QF) \to \L_n(\Ctwo,\QFtwo)\]
differ by a sign, where \(-\eta\) is an additive inverse to \(\eta\) in the \(\Einf\)-group of natural transformations \(\QF \to f^*\QFtwo\) (well-defined up to homotopy).
\end{remark}

\begin{remark}
\label{remark:L-homology}%
We note that the Poincaré functor \((i,\eta)\colon (\C,\QF\qshift{-1}) \to \Met(\C,\QF)\) constructed in Lemma~\refone{lemma:maps-with-metabolic-category} is fully-faithful and the natural transformation \(\eta\colon \QF \Rightarrow i^*\QF_{\met}\) is an equivalence. It then follows that the induced map \(\pi_0\Poinc(\C,\QF\qshift{-1}) \to \pi_0\Poinc(\Met(\C,\QF)) = \pi_0\Poincdel(\C,\QF)\) is injective. Since the essential image of \(i\) coincides with the full subcategory spanned by those objects whose image under \(\met\colon \Met(\C,\QF) \to (\C,\QF)\) is zero it follows that the sequence of monoids
\[ 0 \to \pi_0 \Poinc(\C,\QF\qshift{-1}) \to \pi_0 \Poincdel(\C,\QF) \to \pi_0 \Poinc(\C,\QF)\]
is exact. We may consequently identify \(\L_{n}(\C,\QF)\) with the ``middle homology'' of the sequence of monoids
\[ \pi_0\Poincdel(\C,\QF\qshift{-n}) \to \pi_0\Poincdel(\C,\QF\qshift{-n+1}) \to \pi_0 \Poincdel(\C,\QF\qshift{-n+2}).\]
\end{remark}

The terminology of \(\L\)-groups goes back to Wall~\cite{wall}, who defined the quadratic \(\L\)-groups of a (not necessarily commutative) ring with anti-involution. These play a key role in the surgery theoretic classification of higher dimensional manifolds. In the case of fields the zero'th quadratic \(\L\)-group was first studied by Witt~\cite{witt}, and later became known in this context as the \emph{Witt group}. The zero'th \(\L\)/Witt group of fields plays an important role in arithmetic geometry through its relation to Milnor's K-theory and Galois cohomology, as formulated in Milnor's conjecture and eventually proven by Voevodsky.
For more on the relation between classical \(\L\)-groups of rings and the ones defined above, see Example~\refone{example:ranicki-L-groups} and Remark~\refone{remark:relation-to-classical}.

The notion of \(\L\)-groups was subsequently generalized from rings to categories by various authors: by Ranicki~\cite{Ranickiblue} for categories of bounded complexes over an additive category equipped with a \emph{chain duality}, and later by Weiss and Williams~\cite{WW-duality} for Waldhausen categories equipped with a \emph{Spanier-Whitehead product}. Higher \(\L\)-groups were also studied by Balmer~\cite{balmer-witt} under the name higher Witt groups in the context of triangulated categories. When applied to the triangulated bounded derived categories of \(\ZZ[\tfrac{1}{2}]\)-linear categories with chain duality these agree with Ranicki's \(\L\)-groups, and the latter are a particular case of the construction in~\cite{WW-duality} (see Example 1.A.1 of loc.\ cit.). Let us hence briefly explain how to compare the definition of \(\L\)-groups given~\cite{WW-duality} with the one considered here.

Recall that in~\cite{WW-duality} Weiss and Williams consider a Waldhausen category \(\E\) satisfying certain additional axioms. In particular, \(\E\) is equipped with an initial object \(\emptyset \in \E\) and two distinguished classes of morphisms, called cofibrations and weak equivalences, and the additional axioms ensure that this structure exhibits \(\E\) as a category of cofibrant objects in the sense of~\cite{cisinski}. One may then view \(\E\) as a model for the \(\infty\)-category \(\E_{\infty} = \E[W^{-1}]\) obtained by localizing at the weak equivalences, and the auxiliary structure given by the cofibrations allows one to access certain aspects of \(\E_{\infty}\) while working entirely within \(\E\). For example, the initial object of \(\E\) is also initial in \(\E_{\infty}\) (this is true in any localisation) and one can compute pushouts in \(\E_{\infty}\) by forming in \(\E\) pushouts with one leg cofibration. In fact, the axioms of a category with cofibrant objects ensure that any map in \(\E_{\infty}\) can be represented by a cofibration in \(\E\) (and this is true even if we fix the domain), and so all homotopy pushouts in \(\E_{\infty}\) can be computed in this way.

In the case considered in~\cite{WW-duality} the \(\infty\)-category \(\E_{\infty}\) obtained by localising \(\E\) at the weak equivalences is assumed to be stable, so that one can contemplate the constructions of the present paper for \(\E_{\infty}\). We then point out a second difference with~\cite{WW-duality}: while here we consider forms on objects \(\x \in \E_{\infty}\),
which can then be restricted along maps, in~\cite{WW-duality} one works with the dual notion of co-forms on a given object, which roughly correspond to forms on the dual object. In particular, the rule which associates to any object its space (or spectrum) of co-forms is encoded by a covariant, rather than a contravariant functor on the underlying category. In addition, if one is working with the 1-categorical model \(\E\) for the stable \(\infty\)-category \(\E_{\infty}\) then it is natural to also work with a 1-categorical model for spaces or spectra, which in~\cite{WW-duality} is done by using the category \(\Top\) of compactly generated Hausdorff spaces and the category \(\Spa^{\Om}\) of sequential \(\Om\)-spectra in such. When localised at the class of weak homotopy equivalences these yield the \(\infty\)-categories of spaces and spectra, respectively, and we will treat the former as models for the latter. By a \emph{stable Spanier-Whitehead (SW) product} the authors of~\cite{WW-duality} then mean a symmetric functor
\[ \Bil\colon \E \times \E \to \Spa^{\Om} \]
which is invariant under weak equivalences and exact in each variable, that is, for each \(\x \in \E\) the functor \(\Bil(\x,-)\) sends initial objects to zero \(\Om\)-spectra and pushout squares with one leg cofibration to homotopy pullback squares of \(\Om\)-spectra. In particular, by the above discussion we see that a stable SW-product always descends to symmetric bilinear functor
\[ \Bil_{\infty}\colon \E_{\infty} \times \E_{\infty} \to \Spa \]
in the sense of \S\refone{subsection:quadratic}. Weiss and Williams then require that the resulting functor
\[\pi_0\Bil(-,-)\colon \Ho(\E) \times \Ho(\E) \to \Set\]
is represented in \(\Ho(\E)\) by an equivalence \(\Dual\colon \Ho(\E) \to \Ho(\E)\op\). This condition is implied by the condition that \(\Bil_{\infty}\) is perfect in the sense of Definition~\refone{definition:poinc-cats}, and so to compare~\cite{WW-duality} with the present setting we will simply assume that \(\Bil_{\infty}\) is perfect. One may then speak of symmetric Poincaré objects in \(\E\), which are objects \(\x \in \E\) equipped with a point \(\qone \in \Om^{\infty}\Bil(\x,\x)^{\hC}\) whose corresponding map \(\Dual \x \to \x\) in \(\Ho(\E)\) is an isomorphism. In~\cite{WW-duality} these are simply considered as a set \(\spzero(\E)\), which one can endow with a commutative monoid structure via orthogonal sums.

The notion of a metabolic Poincaré object can naturally be set in the context of~\cite{WW-duality} by associating to the Waldhausen category \(\E\) the Waldhausen category \(\Met(\E)\) whose objects are cofibrations \(\x \to \cob\), the weak equivalences are defined levelwise, and the cofibrations are the Reedy cofibrations. One may then put a stable SW-product on \(\Met(\E)\) by setting \(\Bil([\x \to \cob,\xtwo \to \cobtwo]) = \hofib[\Bil(\x,\xtwo) \to \Bil(\cob,\cobtwo)]\), so that one has a functor of Waldhausen categories with stable SW-products
\[ \Met(\E) \to \E \quad\quad [\x \to \cob] \mapsto \x .\]
A Poincaré object \((\x,\qone)\) in the context of Weiss-Williams is then identified as metabolic if it is the image of a Poincaré object in \(\Met(\E)\). Though not explicitly defined in this manner, the zero'th \(\L\)-group of~\cite{WW-duality} can be identified with the quotient of \(\spzero(\E)\) by the sub-monoid of metabolic objects. These constructions naturally compare with those of the present paper, and one obtains a commutative diagram of commutative monoids
\[
\begin{tikzcd}
\spzero(\Met(\E))\ar[r]\ar[d] & \spzero(\E) \ar[d]\ar[r,twoheadrightarrow] & \L_0(\E,\Bil) \ar[d] \\
\pi_0\Poinc\big(\Met\big(\E\op_\infty,\QF^{\sym}_\infty\big)\big) \ar[r] & \pi_0\Poinc\big(\E\op_{\infty},\QF^{\sym}_\infty\big) \ar[r,twoheadrightarrow] & \L_0\big(\E_\infty\op,\QF^{\sym}_{\infty}\big) , \
\end{tikzcd}
\]
where \(\QF^{\sym}_{\infty} := \QF^{\sym}_{\Bil_\infty}\) is the symmetric Poincaré structure associated to the bilinear functor \(\Bil_\infty\). Here, the composite of the two maps in every row is zero, and the bottom row is exact by Remark~\refone{remark:stably-metabolic}. We then have the following comparison statement:

\begin{proposition}\label{proposition:compare-WW}
The left most vertical map \(\L_0(\E,\Bil) \to \L_0\big(\E_\infty\op,\QF^{\sym}_\infty\big)\) is an isomorphism.
\end{proposition}
\begin{proof}
To begin, note that \(\E \to \E_{\infty}\) is essentially surjective by virtue of being a localisation. In addition, for every \(\x \in \E\) with image \(\ovl{\x} \in \E_{\infty}\) one has \(\Bil(\x,\x)^{\hC} \simeq \Bil_{\infty}(\ovl{\x},\ovl{\x})^{\hC}\) essentially by definition, and since the functor \(\E \to \Ho(\E)\) is conservative we have that a given co-form on \(\x\) is Poincaré in \(\E\) if and only if the corresponding form on \(\Dual \ovl{\x}\) is Poincaré in \(\E_\infty\op\). Combining all this yields that the map of monoids \(\spzero(\E) \to \pi_0\Poinc\big(\E\op_{\infty},\QF^{\sym}_\infty\big)\) is surjective, and so the map \(\L_0(\E,\Bil) \to \L_0\big(\E_\infty\op,\QF^{\sym}_\infty\big)\) is surjective. To show that the last map is also injective it will now suffice to show that for every \((\x,\qone) \in \spzero(\E)\), the map
\[\spzero(\Met(\E)) \times_{\spzero(\E)} \{(\x,\qone)\} \to \pi_0\Poinc(\Met(\E\op,\QF^{\sym}_{\infty})) \times_{\pi_0\Poinc(\E\op,\QF^{\sym}_{\infty})} \{(\ovl{x},\ovl{\qone})\}\]
is surjective. Arguing as for the surjectivity on Poincaré objects, it will suffice to show that the induced functor
\[ \Met(\E) \times_{\E} \{\x\} \to \Met(\E_\infty\op) \times_{\E_\infty\op}\{\ovl{x}\} \]
is essentially surjective, which amounts in this case to showing that every map \(\ovl{\x} \to \y\) out of \(\ovl{\x}\) in \(\E_\infty\) can be represented by a cofibration \(\x \to \y\) out of \(\x\) in \(\E\).
This is a general property of categories of cofibrant objects: indeed, since these admit a left calculus of fractions~\cite[\S 7.2]{cisinski} one may represent \(\ovl{\x} \to \y\) by some map \(\x \to \y\) in \(\E\), which can then be replaced by a cofibration by using the factorization axiom for such categories (which assures that every map can be written as a composite of a cofibration followed by a weak equivalence).
\end{proof}

\begin{remark}\label{remark:also-quadratic}
In~\cite{WW-duality} Weiss and Williams also define quadratic \(\L\)-groups associated to a Waldhausen category with a Spanier-Whitehead product, by replacing homotopy fixed points by homotopy orbits in the definition of co-forms. The proof of Proposition~\refone{proposition:compare-WW} then adapts verbatim to show that the quadratic \(\L\)-groups of~\cite{WW-duality} identify with the \(\L\)-groups of the Poincaré structure \(\QF^{\qdr}_{\infty} := \QF^{\qdr}_{\Bil_\infty}\).
\end{remark}

\begin{example}\label{example:ranicki-L-groups}
Let \(R\) be a commutative ring and \(\Ch(R)\) the category of bounded complexes of finitely generated projective \(R\)-modules. Then \(\Ch(R)\) has the structure of a Waldhausen category with cofibrations being the maps which are levelwise injective with projective cokernel, and weak equivalences the collection \(\qiso\) of quasi-isomorphisms. We may endow this Waldhausen category with a Spanier-Whitehead product \(\Bil^{\otimes}(X,Y) := \GEM(X \otimes_R Y)\) which sends \(X,Y\) to the Eilenberg-Maclane spectrum of \(X \otimes_R Y\). The symmetric and quadratic \(\L\)-groups associated to this Spanier-Whitehead product by~\cite{WW-duality} then coincide with Ranicki's (4-periodic) symmetric and quadratic \(\L\)-groups, respectively, see~\cite[Example 1.A.1]{WW-duality}. Let \(\Dperf(R) = \Ch(R)[\qiso^{-1}]\) be the perfect derived \(\infty\)-category of \(R\) obtained by taking the \(\infty\)-categorical localisation of \(\Ch(R)\) by the quasi-isomorphisms. Then under the equivalence \(\Dual_R\colon\Dperf(R)\op \simeq \Dperf(R)\) the descended bilinear functor \(\Bil^{\otimes}_{\infty}\) identifies with the bilinear functor \(\Bil_R\) of Example~\refone{example:perfect-derived}. By Proposition~\refone{proposition:compare-WW} and Remark~\refone{remark:also-quadratic} we may then identify Ranicki's symmetric and quadratic \(\L\)-groups with the \(\L\)-groups of the Poincaré structures \(\QF^{\sym}_R\) and \(\QF^{\qdr}_R\) described in that example. More generally, the same holds (for the same reason) over a not-necessarily-commutative ring \(R\) with respect to a fixed invertible module with involution \(M\), see \S\refone{subsection:discrete-rings}.
\end{example}

\begin{remark}
The collection of \(\L\)-groups are in fact the homotopy groups of a \emph{spectrum valued} invariant \(\L(\C,\QF)\), known as the \(\L\)-theory spectrum. A definition in the setting of Poincaré \(\infty\)-categories was given in~\cite{Lurie-L-theory}, but was defined much earlier in the setting of rings with anti-involution by Ranicki~\cite{Ranickiblue}, and plays a key role in surgery theory. We will recall the definition of this invariant in \papertwo, prove its main properties and characterize it by a universal property. The interaction between the \(\L\)-spectrum and the closely related Grothendieck-Witt spectrum is one of the principal themes of the present series of papers.
\end{remark}

\subsection{The algebraic Thom construction}
\label{subsection:algebraic-thom}%

Given a Poincaré \(\infty\)-category \((\C,\QF)\), the map
\[\Poinc(\C,\QF\qshift{-1}) \to \Poincdel(\C,\QF)\]
induced by the the Poincaré functor
of Lemma~\refone{lemma:maps-with-metabolic-category} sends a \((-1)\)-shifted Poincaré object \((\x,q)\) to the metabolic Poincaré object \(0\) equipped with \(\x\) as its Lagrangian. In particular, the data of a Lagrangians of \(0\) is equivalent to that of a Poincaré object with respect to the shifted Poincaré structure \(\QF\qshift{-1}\). This idea fits in the more general paradigm of the \emph{algebraic Thom isomorphism} developed by Ranicki \cite[Proposition~3.4]{RanickiATS1}, under which Poincaré objects in \((\C,\QF)\) equipped with a Lagrangian can equivalently be encoded via a hermitian object with respect to \(\QF\qshift{-1}\). To discuss this equivalence, it will be useful to introduce the following construction:

\begin{definition}
\label{definition:arrow-cat}%
Let \((\C,\QF)\) be a Poincaré \(\infty\)-category. We will denote by \(\Ar(\C,\QF)\) the hermitian \(\infty\)-category whose underlying stable \(\infty\)-category is the arrow category \(\Ar(\C)\) and whose quadratic functor \(\QF_{\arr}\) sits in a pullback diagram
\begin{equation}
\label{equation:defining-square-arrow}%
\begin{tikzcd}
\QF_{\arr}([f\colon\z \to \w]) \ar[r] \ar[d] & \QF(\z) \ar[d] \\
\Bil_{\QF}(\z,\w) \ar[r,"{f^*}"] & \Bil_{\QF}(\z,\z)
\end{tikzcd}
\end{equation}
where the vertical map is the canonical one from a quadratic functor to its diagonally restricted bilinear part while the bottom horizontal map is the natural transformation whose component at \(f\) is given by restriction along \(f\).
\end{definition}

\begin{remark}
\label{remark:usual-sequence-with-Q-ar}%
For \([\z \to \cob] \in \Ar(\C)\), the commutative square
\[
\begin{tikzcd}
\QF(\w)\ar[r]\ar[d] & \QF(\z) \ar[d] \\
\Bil_{\QF}(\z,\w) \ar[r,"{f^*}"] & \Bil_{\QF}(\z,\z)
\end{tikzcd}
\]
determines a natural map \(\QF(\w) \to \QF_{\arr}(\z \to \w)\).
From Corollary~\refone{example:usual-sequence} we then get that for an exact sequence \(\z \to \cob \to \x\) the associated sequence
\[
\QF(\x) \to \QF(\cob) \to \QF_{\arr}(\z \to \cob)
\]
is exact.
\end{remark}

Unwinding the definitions we see that the underlying symmetric bilinear functor of \(\QF_{\arr}\) sits in a fibre square
\[
\begin{tikzcd}
\Bil_{\arr}([\z \to \w],[\ztwo \to \wtwo])\ar[r]\ar[d] & \Bil_{\QF}(\w,\ztwo) \ar[d] \\
\Bil_{\QF}(\z,\wtwo) \ar[r] & \Bil_{\QF}(\z,\ztwo)
\end{tikzcd}
\]
from which we we see that when \(\Bil_{\QF}\) is perfect with duality \(\Dual\) then \(\Bil_{\arr}\) is perfect with duality
\[ \Dual_{\arr}([f\colon\z\to \w]) = [\Dual f\colon \Dual\w \to \Dual\z] .\]
In this case, identifying \(\Bil_{\QF}(\z,\w) \simeq \map_{\C}(\w,\Dual\z)\), we see that a hermitian form on an arrow \([f\colon\z \to \w] \in \Ar(\C)\)
consists of a triple \((q,g,\eta)\) where \(q\) is hermitian form on \(\z\) with respect to \(\QF\), \(g\colon \w \to \Dual\z\) is a map in \(\C\), and \(\eta\) is a homotopy \(q_{\sharp} \sim g\circ f\) between the resulting two maps from \(\z\) to \(\Dual\z\). The self-dual map \([\z \to \w] \to [\Dual\w \to \Dual\z]\) associated to such a triple can then be expressed as the map between the two vertical arrows in the square
\[
\begin{tikzcd}
\z \ar[r,"{\Dual g}"] \ar[d,"f"] \ar[dr,dashed,"{q_{\sharp}}"] & \Dual\w \ar[d,"{\Dual f}"] \\
\w \ar[r,"g"] & \Dual\z
\end{tikzcd}
\]
In particular, a triple \((q,g,\eta)\) constitutes a Poincaré form if and only if \(g\) is an equivalence. We now note that the map \(\QF_{\arr}([\z \to \w]) \to \QF(\z)\) appearing in the defining square~\eqrefone{equation:defining-square-arrow} promotes the domain projection
\[
\Ar(\C) \to \C \quad\quad [\z \to \w] \mapsto \z
\]
to a hermitian functor and hence determines a map
\begin{equation}
\label{equation:domain-map-on-hermitian}%
\catforms(\Ar(\C,\QF)) \to \catforms(\C,\QF)
\end{equation}
given on the level of tuples as above by \((\z \to \w,q,g,\eta) \mapsto (\z,q)\). We then have the following:

\begin{proposition}
\label{proposition:thom}%
The map~\eqrefone{equation:domain-map-on-hermitian} restricts to an equivalence
\[
\Poinc(\Ar(\C,\QF)) \to \spsforms(\C,\QF) \quad\quad (\z \to \w,q,g,\eta) \mapsto (\z,q).
\]
An explicit inverse is given by the association
\((\z,q) \mapsto (q_{\sharp}\colon \z \to \Dual\z,q,\id,\id)\).
\end{proposition}

Concisely stated, Proposition~\refone{proposition:thom} says that for a Poincaré \(\infty\)-category \((\C,\QF)\), hermitian objects in \(\C\) can be described via Poincaré objects in its arrow category. Though a direct proof is perfectly possible at the moment, we will postpone it to \S\refone{subsection:thom}, where we will prove it this statement in a more general context in Proposition~\refone{proposition:algebraic-thom}. Meanwhile, let us connect the present conclusions to the above discussion of metabolic objects.

\begin{notation}
\label{notation:seq}%
For a stable \(\infty\)-category \(\C\) we will denote by \(\Seq(\C) \subseteq \Fun(\Del^1 \times \Del^1,\C)\) the full subcategory spanned by the exact squares of the form
\begin{equation}
\label{equation:exact-square}%
\begin{tikzcd}
\z \ar[r]\ar[d] & \cob \ar[d] \\
0 \ar[r] & \x
\end{tikzcd}
\end{equation}
In other words, \(\Seq(\C)\) is the \(\infty\)-category of exact sequences \([\z \to \cob \to \x]\) in \(\C\), where we will often omit the null-homotopy encoded by the commutative square~\eqrefone{equation:exact-square} to simplify notation.
\end{notation}

We have two projections
\begin{equation}
\label{equation:projection-seq-ar}%
\begin{tikzcd}
[column sep=5ex,row sep=1ex]
\Ar(\C) & \Seq(\C)\ar[l]\ar[r] & \Ar(\C) \\
{[\z \to \cob]} & {[\z \to \cob \to \x]} \ar[l,mapsto] \ar[r,mapsto] & {[\cob \to \x]}
\end{tikzcd}
\end{equation}

\begin{lemma}
\label{lemma:arrow-is-met}%
The projections~\eqrefone{equation:projection-seq-ar} are both equivalences. In addition, the restrictions of \(\QF_{\arr}\) and \(\QF_{\met}\qshift{1}\) to \(\Seq(\C)\) via the left and right projection respectively, are naturally equivalent.
\end{lemma}
\begin{proof}
The first claim follows from the fact that exact squares as in~\eqrefone{equation:exact-square} are both left Kan extended from their restriction to \(\Lam^2_0 \subseteq \Del^1\times \Del^1\) and right Kan extended from their restriction to \(\Lam^2_2 \subseteq \Del^1 \times \Del^1\). The natural homotopy between \(\QF_{\arr}|_{\Seq(\C)}\) and \(\QF_{\met}\qshift{1}|_{\Seq(\C)}\) is then encoded by the \(4\)-fold exact sequence of quadratic functors on \(\Seq(\C)\) given by
\[ \QF_{\met}([\cob \to \x]) \to \QF(\x) \to \QF(\cob) \to  \QF_{\arr}(\z \to \cob),\]
where the last exact sequence is by Remark~\refone{remark:usual-sequence-with-Q-ar}.
\end{proof}

Lemma~\refone{lemma:arrow-is-met} identifies the Poincaré \(\infty\)-category \(\Met(\C,\QF)\) with the Poincaré \(\infty\)-category \(\Ar(\C,\QF)\) up to a shift of the Poincaré structure. We note however that \((\Ar(\C),\QF_{\arr}\qshift{1}) \simeq \Ar(\C,\QF\qshift{1})\), that is, the formation of arrow Poincaré \(\infty\)-categories commutes with shifting the Poincaré structure. Using Lemma~\refone{lemma:arrow-is-met} and Proposition~\refone{proposition:thom} we then obtain a sequence of equivalences
\[ \Poinc(\Met(\C,\QF)) \simeq \Poinc(\Ar(\C,\QF\qshift{-1})) \simeq \spsforms(\C,\QF\qshift{-1}). \]
Unwinding the definitions, this composed map sends a tuple \((\cob \to \x,q,\eta)\), consisting of a Poincaré object \((\x,q)\) equipped with a Lagrangian \((\cob \to \x,\eta)\), to the object \(\z := \fib(\cob \to \x)\), equipped with the shifted hermitian structure encoded by the pair of null-homotopies of \(q|_{\z}\) (one restricted from \(\eta\) and one induced by the null-homotopy of the composed map \(\z \to \cob \to \x\)).

\begin{corollary}[The algebraic Thom isomorphism]
\label{corollary:algebraic-thom-iso}%
The association \([\cob \to \x] \mapsto \fib(\cob \to \x)\) underlines a natural equivalence of spaces
\[ \Poinc(\Met(\C,\QF)) \simeq \spsforms(\C,\QF\qshift{-1})\]
between Poincaré objects in \(\C\) equipped with a Lagrangian and \((-1)\)-shifted hermitian objects in \(\C\).
\end{corollary}

\begin{remark}
Combining Remark~\refone{remark:L-homology} with the algebraic Thom isomorphism of Corollary~\refone{corollary:algebraic-thom-iso} we may identify the \(\L\)-groups of \((\C,\QF)\) with the homology monoids of the chain complex of monoids of the form
\[ ... \to \pi_0\spsforms(\C,\QF\qshift{-n-1}) \to \pi_0\spsforms(\C,\QF\qshift{-n}) \to \pi_0 \spsforms(\C,\QF\qshift{-n+1}) \to ...\]
where the map \(\pi_0\spsforms(\C,\QF\qshift{-i}) \to \pi_0 \spsforms(\C,\QF\qshift{-i+1})\) sends an \((-i)\)-fold hermitian object \((\x,\qone)\) to its Ranicki boundary \(\cof[\x \to \Om^i\Dual_{\QF}(\x)]\), endowed with its associated \((-i+1)\)-fold Poincaré form.
\end{remark}

We finish this section by framing the observation that the hyperbolic and arrow constructions naturally commute with each other

\begin{lemma}
\label{lemma:arrow-hyp-commute}%
For a stable \(\infty\)-category \(\C\), the natural equivalence
\begin{equation}
\label{equation:arrows-op-is-op-arrows}%
\Ar(\C\times \C\op)\simeq \Ar(\C)\times \Ar(\C\op) \simeq \Ar(\C) \times \Ar(\C)\op
\end{equation}
in which the second equivalence sends \((f\colon \z \to \w,f'\colon \ztwo\to \wtwo)\) to \((f\colon\z \to \w,\wtwo \leftarrow \ztwo\cocolon f')\), extends to an equivalence of Poincaré \(\infty\)-categories
\[ \Ar(\Hyp(\C)) \simeq \Hyp(\Ar(\C)).\]
\end{lemma}
\begin{proof}
Transporting the Poincaré structure of \(\Ar(\Hyp(\C))\) along the equivalence~\eqrefone{equation:arrows-op-is-op-arrows} yields the quadratic functor
\[(f\colon\z \to \w,f'\colon\ztwo \to \wtwo) \mapsto \QF_{\hyp}(\z,\wtwo) \times_{\Bil_{\hyp}((\z,\wtwo),(\z,\wtwo))} \Bil_{\hyp}((\z,\wtwo),(\w,\ztwo)) \simeq \]
\[\map_{\C}(\z,\wtwo) \times_{\map_{\C}(\z,\wtwo) \times \map_{\C}(\z,\wtwo)}[\map_{\C}(\z,\ztwo) \times \map_{\C}(\w,\wtwo)] \simeq \map_{\C}(\z,\ztwo) \times_{\map_{\C}(\z,\wtwo)} \map_{\C}(\w,\wtwo) \]
\[\simeq \map_{\Ar(\C)}(f,f').\]
We thus finish the proof by recognizing the last term as the quadratic functor of \(\Hyp(\Ar(\C))\).
\end{proof}

Combining Lemma~\refone{lemma:arrow-hyp-commute} with Lemma~\refone{lemma:arrow-is-met} and Remark~\refone{remark:hyp-shift-invariant} we immediately conclude

\begin{corollary}
\label{corollary:met-hyp-splits}%
For a stable \(\infty\)-category \(\C\) there is a natural equivalence of Poincaré \(\infty\)-categories
\[ \Met(\Hyp(\C)) \simeq \Hyp(\Ar(\C)) ,\]
which, on the underlying stable \(\infty\)-categories, is given by the equivalence
\[
\Ar(\C\times\C\op)\simeq  \Ar(\C)\times\Ar(\C\op)\simeq \Ar(\C)\times\Ar(\C)\op,
\]
where the second equivalence is the product of the identity of \(\Ar(\C)\) and the equivalence \(\Ar(\C\op)\simeq \Ar(\C)\op\) which sends an arrow \(\x\leftarrow \cob\colon g\) in \(\C\op\) to the canonical arrow \(\fib(g)\to \cob\).
In particular, the Poincaré functor \(\met\colon \Met\Hyp(\C)\to \Hyp(\C)\) from Lemma~\refone{lemma:maps-with-metabolic-category} (applied to \(\Hyp(\C)\)) admits a section in \(\Catp\).
\end{corollary}

\subsection{The Grothendieck-Witt group}
\label{subsection:GW-group}%

In this section we define the Grothendieck-Witt group of a Poincaré \(\infty\)-category. In the setting of ordinary rings, the Grothendieck-Witt group was classically defined as the group completion of the monoid of isomorphism classes of pairs \((P,q)\) where \(P\) is finite dimensional projective module and \(q\) is a non-degenerate hermitian form of some flavour (symmetric, quadratic, anti-symmetric, etc.). It was later extended to more general contexts such as vectors bundles over algebraic varieties \cite{knebusch1977symmetric} and forms in abstract additive categories with duality~\cite{devissage}. In doing so it was realized that the simple definition via group completion needs to be slightly modified to take into account information coming from non-split short exact sequences. In particular, one had to quotient out the group completion by the relation \([\x,\qone] \sim [\hyp(\cob)]\) identifying the class of a metabolic object with Lagrangian \(\cob\) with that of the associated hyperbolic object. In fact, the latter relation already implies the group property for the resulting quotient (as we will see below in our context), and hence can be done on the level of monoids without explicitly group completing. On the other hand, in the case of modules over rings (or, more generally, in contexts in which all short exact sequences split) the relation \([\x,\qone] \sim [\hyp(\cob)]\) automatically holds in the group completion since every metabolic object is stably hyperbolic.

In the present section we will give a definition of the Grothendieck-Witt group in the context of Poincaré \(\infty\)-categories, and extract some of its basic properties. Using the work of the fourth and ninth authors~\cite{comparison}, this definition can be compared with the classical one in the case of modules over rings using suitable Poincaré structures on \(\Dperf(R)\), see Remark~\refone{remark:comparison} below for the precise statement.

For the following definition, recall the map \(\hyp\colon \grpcr\C \to \Poinc(\C,\QF)\) induced on Poincaré objects by the Poincaré functor \(\hyp\colon \Hyp(\C) \to (\C,\QF)\) under the equivalence \(\Poinc(\Hyp(\C)) \simeq \grpcr\C\) of Proposition~\refone{proposition:poinc-of-hyp}. Explicitly, this map sends an object \(\x \in \C\) to the Poincaré object \(\hyp(\x) = \x \oplus \Dual\x\) equipped with its canonical Poincaré form, see \S\refone{subsection:hyp-and-sym-poincare-objects}.

\begin{definition}
Let \((\C,\QF)\) be a Poincaré \(\infty\)-category. We define \(\GW_0(\C,\QF)\) to be the quotient in the category of commutative monoids of \(\pi_0\Poinc(\C,\QF)\) (with its commutative monoid structure given by direct sums) by the relations
\begin{equation}
\label{equation:GW-relations}%
[\x,\qone] \sim [\hyp(\cob)]
\end{equation}
for every Poincaré object \((\x,\qone)\) with Lagrangian \(\cob \to \x\).
\end{definition}

\begin{remark}
By definition we may identify \(\GW_0(\C,\QF)\) with the coequalizer of the pair of maps
\[ \pi_0\Poinc\Met(\C,\QF) \rightrightarrows \pi_0\Poinc(\C,\QF) ,\]
where the first map is induced by the Poincaré functor \(\met\colon \Met(\C,\QF) \to (\C,\QF)\) and the second is induced by the composed Poincaré functor
\(\Met(\C,\QF) \xrightarrow{\lag} \Hyp(\C) \xrightarrow{\hyp} (\C,\QF)\)
discussed in Remark~\refone{remark:lag-split-dual}.
\end{remark}

We quickly summarize a few properties which follow directly from the definition of \(\GW_0\).

\begin{lemma}
\label{lemma:basic-GW-group}%
\ %
\begin{enumerate}
\item
\label{item:basic-1}%
For every exact sequence \(\ztwo \to \z \to \zthree\) in \(\C\) the relation
\[[\hyp(\z)] \sim [\hyp(\ztwo)] + [\hyp(\zthree)] \]
holds in \(\GW_0(\C,\QF)\).
\item
\label{item:basic-2}%
For every Poincaré object \([\x,q] \in \Poinc(\C,\QF)\) the relation
\[ [\x,q] + [\x,-q] + [\hyp(\Omega x)] \sim [\hyp(\x)] + [\hyp(\Omega \x)] \sim 0 \]
holds in \(\GW_0(\C,\QF)\).
\end{enumerate}
\end{lemma}
\begin{proof}
For~\refoneitem{item:basic-1} note that if \(\ztwo \to \z \to \zthree\) is an exact sequence in \(\C\) then \(\cob := \ztwo\oplus \Dual_{\QF}\zthree\) is naturally a Lagrangian in \(\hyp(\z)\), and \(\hyp(\cob) \simeq \hyp(\ztwo) \oplus \hyp(\zthree)\).

To prove~\refoneitem{item:basic-2}, the first identification is given by~\eqrefone{equation:GW-relations} applied to the metabolic object \((\x \xrightarrow{\Delta} \x \oplus \x, q \oplus -q)\) and the second relation is given by~\refoneitem{item:basic-1}
applied to the exact sequence \(\Omega \x \rightarrow 0 \rightarrow \x\).
\end{proof}

\begin{corollary}
\label{corollary:GW_0}%
The commutative monoid \(\GW_0(\C,\QF)\) is always a \emph{group}. We will refer to it as the \defi{Grothendieck-Witt group} of \((\C,\QF)\).
\end{corollary}

\begin{example}
\label{example:GW-of-hyp}%
For a stable \(\infty\)-category \(\C\), the isomorphism
\begin{equation}
\label{equation:poinc-of-hyp}%
\pi_0\Poinc(\Hyp(\C)) \cong \pi_0\grpcr(\C)
\end{equation}
induced on components of the equivalence of Proposition~\refone{proposition:poinc-of-hyp}, descends to a group isomorphism
\[ \GW_0(\Hyp(\C)) \cong \K_0(\C). \]
Indeed, the isomorphism~\eqrefone{equation:poinc-of-hyp} relates the class of \(\z \in \C\) to the Poincaré class of \((\z,\z) \in \C \times \C\op\) equipped with its canonical Poincaré form \(\id_{\z} \in \map_{\C}(\z,\z) = \QF_{\hyp}(\z,\z)\).
An isotropic object in \(((\z,\z),\id_{\z})\) is then given by a pair of maps \(\ztwo \to \z,\z \to \zthree\) in \(\C\), such that the composite \(\ztwo \to \z \to \zthree\), which corresponds to the pullback of the Poincaré form \(\id_{\z}\) to \(\QF_{\hyp}(\ztwo,\zthree)=\map_{\C}(\ztwo,\zthree)\), vanishes. Such an isotropic object is a Lagrangian precisely when the resulting sequences \(\ztwo \to \z \to \zthree\) is exact. Moreover the hyperbolic object on \((\ztwo,\zthree)\) is given by the object \((\ztwo \oplus \zthree, \ztwo \oplus \zthree) \in \C \times \C\op\).  It then follows that under the isomorphism~\eqrefone{equation:poinc-of-hyp}, the defining relations of the Grothendieck-Witt group can be written as \([\z] = [\ztwo] + [\zthree]\) for every exact sequence \(\ztwo \to \z \to \zthree\), which are exactly the relations defining the quotient \(\pi_0\grpcr(\C) \twoheadrightarrow \K_0(\C)\).
\end{example}

Recall from Definition~\refone{definition:L-groups} that zero'th \(\L\)-group \(\L_0(\C,\QF)\) is defined as the cokernel of the monoid homomorphism \(\pi_0\Poincdel(\C,\QF) \to \pi_0\Poinc(\C,\QF)\). In particular, the quotient map \(\pi_0\Poinc(\C,\QF) \to \L_0(\C,\QF)\) sends the class of any metabolic object (and in particular any hyperbolic object) to zero, and hence factors through a group homomorphism
\[ \GW_0(\C,\QF) \to \L_0(\C,\QF) ,\]
which is necessarily surjective, as we can identify \(\L_0(\C,\QF)\) with the quotient group of \(\GW_0(\C,\QF)\) by the subgroup spanned by the classes of metabolic objects. To obtain more information about this kernel, we note that the map
\[\pi_0\Poincdel(\C,\QF) \to \pi_0(\grpcr\C) \quad\quad [\cob \to \x] \mapsto \cob \]
induced by the Poincaré functor \(\lag\colon \Met(\C,\QF) \to \Hyp(\C)\) of Construction~\refone{construction:hyp-to-met}, is surjective: any object \(\z \in \C\) is a Lagrangian in the Poincaré object \(\hyp(\z)\). In addition, the canonical map
\[\pi_0(\grpcr\C) \to \K_0(\C)\]
to the zero'th algebraic \(\K\)-theory group of \(\C\) is surjective as well: we may identify \(\K_0(\C)\) with the quotient of \(\pi_0(\grpcr\C)\) in the category of commutative monoids by the relations \([\z] \sim [\ztwo] + [\zthree]\) for every exact sequence \(\ztwo \to \z\to \zthree\). By Lemma~\refone{lemma:basic-GW-group}\refoneitem{item:basic-1} the homomorphism \(\hyp\colon \pi_0(\grpcr\C) \to \pi_0\Poinc(\C,\QF)\) then descends to a homomorphism of abelian groups
\begin{equation}
\label{equation:hyp-K-to-GW}%
[\hyp]\colon \K_0(\C) \to \GW_0(\C,\QF),
\end{equation}
which by Example~\refone{example:GW-of-hyp} we may also identify with the map induced on \(\GW_0\) by the Poincaré functor \(\hyp\colon \Hyp(\C) \to \C\).
Comparing the relevant universal properties
we then see that \(\GW_0(\C,\QF)\) sits in a pushout square of commutative monoids of the form
\begin{equation}
\label{equation:gw-group}%
\begin{tikzcd}
\pi_0\Poincdel(\C,\QF) \ar[r,"{[\met]}"] \ar[d,"{[\lag]}"'] & \pi_0\Poinc(\C,\QF) \ar[d] \\
\K_0(\C) \ar[r,"{[\hyp]}"] & \GW_0(\C,\QF)
\end{tikzcd}
\end{equation}
where \([\lag]\) denotes the composed map \(\pi_0\Poincdel(\C,\QF) \to \pi_0(\grpcr\C) \to \K_0(\C)\).
It then follows that \(\L_0(\C,\QF)\) can equivalently by obtained as the cokernel of the map of abelian groups \([\hyp]\colon \K_0(\C) \to \GW_0(\C,\QF)\). On the other hand, by Corollary~\refone{corollary:fgt-hyp-equivariant} the map \([\hyp]\)
is \(\Ct\)-equivariant
with respect to the induced \(\Ct\)-action on \(\K_0\) and the trivial \(\Ct\)-action on \(\GW_0(\C,\QF)\). We then obtain an induced exact sequence of abelian groups
\begin{equation}
\label{equation:mini-tate}%
\K_0(\C)_{\Ct} \to \GW_0(\C,\QF) \to \L_0(\C,\QF) \to 0 .
\end{equation}
In \papertwo we will show that this sequence comes from an exact sequence of spectra
\[ \K(\C)_{\hC} \to \GW(\C,\QF) \to \L(\C,\QF),\]
the \defi{Tate exact sequence}, encoding the fundamental relationship between these three invariants. In particular, this allows one to extend~\eqrefone{equation:mini-tate} to a long exact sequence involving the higher \(\L\)-groups and higher Grothendieck-Witt groups, yielding a powerful tool for computing the latter. In \paperthree we will exploit these ideas for computing the higher Grothendieck-Witt groups of the integers.

%% file: Modules.tex
In this section we discuss hermitian and Poincaré structures on \(\infty\)-categories of modules over ring spectra.
In particular, we fix a base \(\Einf\)-ring spectrum \(k\) and consider an \(\Eone\)-algebra \(A\) in the symmetric monoidal \(\infty\)-category \(\Mod_k\) of \(k\)-module spectra. We denote by \(\Alg_{\Eone} := \Alg_{\Eone}(\Mod_k)\) the \(\infty\)-category of \(\Eone\)-algebra objects in \(\Mod_k\), which we simply refer to as \defi{\(\Eone\)-algebras}. Given an \(\Eone\)-algebra \(A \in \Alg_{\Eone}\), we denote by \(\Mod_A\) the \(\infty\)-category of left \(A\)-module objects in \(\Mod_k\), which we refer to as \defi{\(A\)-modules}.
Our goal is to describe and study hermitian and Poincaré structures on the full subcategory \(\Modp{A} \subseteq \Mod_A\) of compact \(A\)-modules in \(\Mod_k\), as these account for many examples of interest motivating the present work
(see also \S\refone{section:examples} for some more specific examples). More generally, one often wishes to consider a dense stable subcategory of \(\Modp{A}\), that is, a stable (and in particular full) subcategory which generates all of \(\Modp{A}\) under retracts. For example, one may consider the stable subcategory \(\Modf{A} \subseteq \Modp{A}\) of finitely presented \(A\)-module spectra, which is the smallest stable subcategory containing the \(A\)-module \(A\). By the work of Thomason~\cite{thomason-classification}, the dense stable subcategories of \(\Modp{A}\) are classified by subgroups \(c \subseteq \K_0(\Modp{A})\), where to such a \(c\) corresponds the dense subcategory \(\Mod^c_A \subseteq \Modp{A}\) consisting of those compact \(A\)-modules \(M\) whose class \([M] \in \K_0(\Modp{A})\) lies in \(c\). In this context, let us point out that studying hermitian and Poincaré structures on \(\Mod^c_A\) for a general \(c \subseteq \K_0(\Modp{A})\) is not really more complicated then studying hermitian and Poincaré structures in the maximal case of \(\Modp{A}\). Indeed, by Lemma~\refone{lemma:kan-extension-exact-quadratic} left Kan extension and restriction form an adjunction
\[
\Funq(\Mod^c_A) \adj \Funq(\Modp{A}),
\]
and \(\Mod^c_A \subseteq \Modp{A}\) being a dense full subcategory implies that this adjunction is an equivalence (cf.\ Remark~\reftwo{remark:idem-poincare}). As will follow from our more detailed analysis below, this equivalence identifies the \(\infty\)-category \(\Funp(\Mod^c_A)\) of Poincaré structures on \(\Mod^c_A\) with the full subcategory of \(\Funp(\Modp{A})\) spanned by those Poincaré structures whose associated duality preserves \(\Mod^c_A\), or, equivalently, preserves the subgroup \(c \subseteq \K_0(\Modp{A})\) (see Remark~\refone{remark:classification-poincare}).

The present section is organised as follows. We begin in \S\refone{subsection:modules-with-involution} by introducing the notion of a \(k\)-module with involution, and show how it can be used to construct bilinear functors on module \(\infty\)-categories. We then extend this notion in \S\refone{subsection:genuine-modules} to a \(k\)-module with \defi{genuine} involution, that which allows to refine to associated bilinear functor to a hermitian or Poincaré structure.
The basic operations of restriction and induction of modules with genuine involution along maps of ring spectra are discussed in~\S\refone{subsection:restriction-induction}.

We point out that the \(\infty\)-categories \(\Mod^c_A\) we will consider do not depend on \(k\), that is, they only depend on the underlying \(\Eone\)-ring spectrum of \(A\). However, the notion of a \(k\)-module with (genuine) involution does depend on \(k\), and affects the type of hermitian structures on \(\Mod^c_A\) one can obtain in this way. More precisely, we will see in \S\refone{subsection:classification-in-modules} that when \(k=\SS\), modules with genuine involution precisely correspond to hermitian structure on \(\Modp{A}\), while for a general \(k\) one recovers hermitian structures equipped with a certain additional \(k\)-linear structure, the precise meaning of which is described in Remark~\refone{remark:k-linear-classification} and Example~\refone{example:k-linear-structures}. The reader who wishes to avoid this additional layer of structure is invited to assume \(k=\SS\) throughout. On the other hand, the reader who prefers to reason in terms of chain-complexes instead of spectra is invited to consider \(k=\ZZ\). The latter case (and more generally, that of a complex-oriented \(k\)) also leads to better periodicity properties, as we discuss \S\refone{subsection:herm-shifts}.

\subsection{Ring spectra and involutions}
\label{subsection:modules-with-involution}%

In this section we will define the notion of a \defi{\(k\)-module with involution} over an \(\Eone\)-algebra \(A\) and show how it can be used to construct bilinear functors on
\(\Modp{A}\). We will identify the \(k\)-modules with involution which lead to perfect bilinear functors as those which are \defi{invertible} in a suitable sense. An important class of invertible \(k\)-modules with involution arise from \(\Eone\)-algebras with anti-involutions, a case for which we will present a convenient recognition criterion, see Proposition~\refone{proposition:recognition}.

To begin, we note that since \(\otimes_k\) is a symmetric monoidal structure on \(\Mod_k\) the monoidal product
\[
\otimes_k \colon \Mod_k \times \Mod_k \to \Mod_k
\]
itself refines to a monoidal functor. This monoidal functor then sends the algebra object \((A,A)\) in \(\Mod_k \times \Mod_k\) to the algebra object \(A \otimes_k A\) in \(\Mod_k\), and thus refines to a functor
\begin{equation}
\label{equation:univeral-bilinear}%
\Mod_A \times \Mod_A \to \Mod_{A \otimes_k A} \quad\quad (X,Y) \mapsto X \otimes_k Y .
\end{equation}
In addition, since \(\otimes_k\) is a symmetric monoidal structure the functor \(\otimes_k\colon \Mod_k \times \Mod_k \to \Mod_k\) is \(\Ct\)-equivariant with respect the flip action on the domain and trivial action on the target, and consequently~\eqrefone{equation:univeral-bilinear} inherits a \(\Ct\)-equivariant structure with respect to the flip action on the domain and the \(\Ct\)-action on the target induced by the flip action on \(A \otimes_k A\). Given an \((A \otimes_k A)\)-module \(M\), let us denote by
\[
\Bil_M\colon \Modp{A} \times \Modp{A} \to \Spa \quad\quad (X,Y) \mapsto \map_{A \otimes_k A}(X \otimes_k Y,M)
\]
the resulting bilinear form.
The association \(M \mapsto \Bil_M\) then assembles to form a functor
\begin{equation}
\label{equation:bil-M}%
\Bil_{(-)}\colon \Mod_{A \otimes_k A} \to \Funb(\Modp{A})
\end{equation}
which by the above inherits a \(\Ct\)-equivariant structure with respect to the flip-induced \(\Ct\)-actions on both sides.

\begin{definition}
\label{definition:module-with-involution}%
Let \(A\) be an \(\Eone\)-algebra over \(k\). By a \defi{\(k\)-module with involution} over \(A\) we will mean an object \(M\) of the \(\infty\)-category \((\Mod_{A\otimes_k A})^{\hC}\), where as above \(\Ct\) acts on \(\Mod_{A \otimes_k A}\) via its flip action on \(A \otimes_k A\). When \(k=\SS\) we will often omit \(k\) and simply write \defi{module with involution}.
\end{definition}

Concretely, a \(k\)-module with involution over \(A\) consists of a spectrum \(M\) with a \(\Ct\)-action and an \((A\otimes_k A)\)-module structure, such that the involution is linear over the ring map \(A\otimes_k A \to A\otimes_k A\) which switches the two factors. Since~\eqrefone{equation:bil-M} is \(\Ct\)-equivariant the bilinear functor
\[
\Bil_{M}(X,Y) = \map_{A\otimes_k A}(X\otimes_k Y, M).
\]
on \(\Modp{A}\) associated to \(M\) consequently inherits the structure of a \emph{symmetric} bilinear functor.

\begin{definition}
\label{definition:functors-associated-to-module-with-involution}%
Let \(M\) be a \(k\)-module with involution over an \(\Eone\)-algebra \(A\). We will denote by
\[
\QF^{\qdr}_M, \QF^{\sym}_M\colon (\Modp{A})\op\to \Spa
\]
the quadratic and symmetric hermitian structures associated to the symmetric bilinear functor \(\Bil_{M} \in \Funs(\Modp{A})\) as in Example~\refone{example:quadratic-symmetric} and Definition~\refone{definition:symmetric-quadratic}.
\end{definition}

If \(M\) is a \(k\)-module with involution over \(A\), these are given explicitly by the formulas
\[
\QF^\qdr_M(X) = \map_{A\otimes_k A}(X\otimes_k X,M)_{\hC}\ \ \ \ \ \ \ \ \mbox{and} \ \ \ \ \  \ \ \ \QF^\sym_M(X) = \map_{A\otimes_k A}(X\otimes_k X,M)^{\hC}.
\]

\begin{remark}
Given a subgroup \(c \subseteq \K_0(\Modp{A})\) we may restrict \(\Bil_M\) to obtain a symmetric bilinear functor on \(\Mod^c_A\), and similarly for the associated hermitian structures \(\QF^{\qdr}_M\) and \(\QF^{\sym}_M\). For simplicity we will generally not distinguish in notation between these functors and their respective restrictions to \(\Mod^c_A\) for a given \(c\).
\end{remark}

We may consider an \((A \otimes_k A)\)-module \(M\) as a \(k\)-module spectrum equipped with two \emph{commuting} actions of \(A\). The first \(A\)-action then promotes \(M\) to an object of \(\Mod_A\), while the second action refines to an action of \(A\) on \(M\) via \(A\)-module maps. In particular, the second \(A\)-action can be encoded via a map
\begin{equation}
\label{equation:second-action}%
A \to \map_A(M,M).
\end{equation}
If \(M\) is a \(k\)-module with involution over \(A\) then the involution determines an equivalence between the two different \(A\)-module structures. In particular, in this case it does not matter which action is considered first and which is considered second.

\begin{definition}
\label{definition:invertible}%
We will say that a \(k\)-module with involution \(M\) over \(A\) is \defi{invertible} if it
is compact as an \(A\)-module (with respect to either the first or the second \(A\)-action) and
the map~\eqrefone{equation:second-action} is an equivalence.%
\end{definition}

\begin{notation}
\label{notation:mapping-module}%
If \(M\) is a \(k\)-module with involution over \(A\) then for every \(A\)-module \(X\) we may form the mapping spectrum \(\map_A(X,M)\), where we consider \(M\) as an \(A\)-module via its first \(A\)-action. This spectrum then carries an \(A\)-action via the residual second \(A\)-action on \(M\), and so we may view \(\map_A(X,M)\) as an \(A\)-module. In what follows we will always consider \(\map_A(X,M)\) as an \(A\)-module in this manner, without indicating it explicitly in the notation.
\end{notation}

\begin{proposition}
\label{proposition:invertible-poincare}%
Let \(A\) be an \(\Eone\)-algebra and \(M\) a \(k\)-module with involution over \(A\). Then for a subgroup \(c \subseteq \K_0(\Modp{A})\) the following statements hold:
\begin{enumerate}
\item
\label{item:non-degenerate}%
The restriction of \(\Bil_M\) to \(\Mod^c_A\) is non-degenerate if and only if the \(A\)-module \(\map_A(X,M)\) (see Notation~\refone{notation:mapping-module}) belongs to \(\Mod^c_A\) for every \(X \in \Mod^c_A\) .
\item
\label{item:duality}%
When the equivalent conditions of~\refoneitem{item:non-degenerate} hold
the associated duality on \(\Mod^c_A\) is given by
\[
\Dual_M(X) := \map_A(X,M).
\]
This duality is furthermore an equivalence if and only if \(M\) is in addition invertible.
\end{enumerate}
\end{proposition}

\begin{remark}
\label{remark:invertible-poincare}%
In the situation of Proposition~\refone{proposition:invertible-poincare}, if \(\Mod^c_A\) contains \(A\) then the condition in~\refoneitem{item:non-degenerate}
necessitates in particular that \(M\) belongs to \(\Mod^c_A\) (take \(X=A\)). If \(\Mod^c_A\) is either \(\Modp{A}\) or \(\Modf{A}\) then \(M\) belonging to \(\Mod^c_A\) is also sufficient: for \(\Modf{A}\) this is because \(A\) generates \(\Modf{A}\) under finite colimits and desuspensions and for \(\Modp{A}\) it is because \(A\) generates \(\Modp{A}\) under finite colimits, desuspensions and retracts, and \(\Modp{A}\) is furthermore idempotent complete.
\end{remark}

\begin{remark}
\label{remark:invertible-criterion}%
It follows from Proposition~\refone{proposition:invertible-poincare} and Remark~\refone{remark:invertible-poincare} that \(M\) is invertible if and only if it is compact as an \(A\)-module and the contravariant functor \(X \mapsto \map_A(X,M)\) from \(\Modp{A}\) to itself is an equivalence. Now for any \(\Eone\)-algebra \(A\) the functor \(X \mapsto \map_A(X,A)\) determines an equivalence
\[
\Modp{A} \longrightarrow (\Modp{A\op})\op.
\]
In addition, for \(M\) compact the natural map \(\map_A(X,A) \otimes_A M \to \map_A(X,M)\) is an equivalence. It then follows that for a \(k\)-module with involution \(M\) which is compact as an \(A\)-module the condition of being invertible is equivalent to the condition that
\[
(-)\otimes_A M\colon \Modp{A\op} \to \Modp{A}
\]
is an equivalence of \(\infty\)-categories.
\end{remark}

\begin{proof}[Proof of Proposition~\refone{proposition:invertible-poincare}]
We begin with the first claim.
For a fixed compact \(A\)-module \(Y\), the functor \(X \mapsto \map_{A \otimes_k A}(X \otimes_k Y,M)\) from \(\Mod_A\op\) to \(\Spa\) is represented by the \(A\)-module \(\map_A(Y,M)\), where the latter is considered as an \(A\)-module via its second \(A\)-action. Since the functor \(\Mod_A \to \Ind(\Mod^c_A) = \Funx((\Mod^c_A)\op,\Spa)\) sending a module to the presheaf it represents on \(\Mod^c_A\) is fully-faithful (in fact, it is an equivalence) it follows that \(\Bil_M(-,Y)\) is representable in \(\Mod^c_A\) if and only if \(\map_A(Y,M)\) lies in \(\Mod^c_A\), in which case \(\map_A(Y,M)\) serves as a representing object already in \(\Mod^c_A\). This proves~\refoneitem{item:non-degenerate} and the first part of~\refoneitem{item:duality}.

To prove the second part of~\refoneitem{item:duality} let now suppose that the condition in~\refoneitem{item:non-degenerate} holds, that is, \(\map_A(X,M)\) belongs to \(\Mod^c_A\) for every \(X \in \Mod^c_A\). This means in particular that \(\map_A(X,M)\) is compact for every \(X \in \Mod^c_A\). Since \(\Mod^c_A\) is dense in \(\Modp{A}\) it follows that \(A\) is a retract of an object in \(\Mod^c_A\), and hence \(M = \map_A(A,M)\) is compact as well.
We now wish to show that the evaluation map
\[
X \to \Dual_M\Dual_M(X) = \map_A(\map_A(X,M),M))
\]
is an equivalence for every \(X \in \Mod^c_A\) if and only if \(M\) is furthermore invertible. Using again that \(\Mod^c_A\) generates \(\Modp{A}\) under retracts the above evaluation map is an equivalence on \(\Mod^c_A\) if and only if it is an equivalence on all of \(\Modp{A}\). On the other hand, since \(\Modp{A}\) is generated under finite colimits and retracts by \(A\) it will suffice to check the component of the evaluation map at \(X=A\), in which case it becomes the map
\[
A \to  \map_A(\map_A(A,M),M)) = \map_A(M,M)
\]
which by definition is an equivalence if and only if \(M\) is invertible. \end{proof}

We now consider some examples.
An \defi{\(\Eone\)-algebra with anti-involution} (over \(k\)) is an object of \(\Alg_{\Eone}^{\hC} = \Alg_{\Eone}(\Mod_k)^{\hC}\), where the action of \(\Ct\) on the \(\infty\)-category of \(\Eone\)-algebras is given by sending an algebra \(A\) to its opposite \(A\op\).

\begin{example}
\label{example:module-with-involution-associated-to-ring-with-involution}%
Let \(A\) be an \(\Eone\)-algebra with anti-involution \(\tau \colon A\op \to A\). %
Then \(A\) can be naturally considered as an invertible \(k\)-module with involution over itself.
Indeed, by using that the forgetful functor \(\Alg_{\Eone} \to \Spa\) is equivariant with respect to the trivial action on the target, we obtain a functor \(\Alg_{\Eone}^{\hC} \to \Fun(\BC,\Spa)\), which allows us to view \(A\) as a spectrum with \(\Ct\)-action. In addition, by construction this action switches between the canonical left and right actions of \(A\) on itself. More precisely, viewing \(A\) as an \((A \otimes_k A\op)\)-module, this \(\Ct\)-action is linear over the \(\Ct\)-action on \(A \otimes_k A\op\) which flips the two components and applies the anti-involution \(\tau\). Equivalently, the anti-involution determines an equivalence of \(\Eone\)-algebras \(A \otimes_k A\op \simeq A \otimes_k A\) intertwining the above \(\Ct\)-action with the flip action on \(A \otimes_k A\), and we may hence view \(A\) as an object of \(\Mod_{A \otimes_k A}^{\hC}\).
Informally, the \(A \otimes_k A\)-action on \(A\) is given by \((a\otimes b)\cdot x = a\cdot x \cdot \tau(b)\). We may then recover the anti-involution \(\tau\) as the induced map of \(\Eone\)-algebras \(A \to \map_A(A,A)=A\op\). The latter is therefore an equivalence, and so \(A\) is invertible.
\end{example}

\begin{example}
\label{example:commutativetrivialaction}%
The restriction of the \(\Ct\)-action on \(\Alg_{\Eone}\) to \(\Einf\)-algebras is canonically trivialised, so that we obtain a functor
\[
\Fun(\BC,\Alg_{\Einf}) \to\Alg_{\Eone}^{\hC}.
\]
In particular, \(\Einf\)-algebras with \(\Ct\)-actions give rise to an \(\Eone\)-algebra with involution. For example, any \(\Einf\)-algebra equipped with the trivial \(\Ct\)-action, determines an \(\Eone\)-algebra with involution.

More generally, when \(A\) is an \(\Einf\)-algebra, any \(A\)-module \(M\) with \(\Ct\)-action canonically defines a \(k\)-module with involution over \(A\), with \(A \otimes_k A\) acting via the multiplication map \(A \otimes_k A \to A\). By Remark~\refone{remark:invertible-criterion}, the invertibility condition from \refone{definition:invertible} is in this case equivalent to the \(\otimes_A\)-invertibility of \(M\), or in other words to \(M \in \Pic(A)\) upon disregarding
\(\Ct\)-actions.
\end{example}

\begin{example}
An important source of \(\Eone\)-algebras with anti-involution arises from the group algebra construction. We will study this example and its relation to \emph{visible \(\L\)-theory} in further detail in~\S\refone{subsection:visible}.
\end{example}

\begin{examples}
\label{example:ordinary-ring}%
Let \(R\) be an ordinary associative ring. Using the Eilenberg-Mac Lane embedding \(\GEM\colon \Ab \hrar \Spa\) we may associate to \(R\) an \(\Eone\)-ring spectrum \(\GEM R\), and by~\cite[Theorem~7.1.2.1]{HA}
we have a natural equivalence
\[
\Modp{\GEM R} \simeq \Dperf(R)
\]
between the \(\infty\)-category of compact \(\GEM R\)-module spectra and the perfect derived \(\infty\)-category of \(R\), defined as in Example~\refone{example:perfect-derived}. An anti-involution on \(\GEM R\) is then the same as an anti-involution on \(R\), that is, an isomorphism \(\tau \colon R \rightarrow R\op\) such that \(\tau^2 = \id\). We will study this case in further detail in \S\refone{subsection:discrete-rings}. Ordinary rings which carry anti-involutions are then fairly common, see Examples~\refone{examples:rings-with-anti-involution}.
\end{examples}

\begin{example}
\label{example:wall-antistructure}%
As in Example~\refone{example:ordinary-ring}, suppose that
\(R\) is an ordinary associative ring. Recall that a \defi{Wall anti-structure}~\cite{wall1970axiomatic} on \(R\) consists of an anti-automorphism \(\tau \colon R \to R\op\) and a unit \(\eps \in R^*\) such that
\(\tau(r) = \eps^{-1} r\eps\) and \(\tau(\eps)=\eps^{-1}\).
The most common type of these are the \defi{central} Wall anti-structures, namely, those in which \(\eps\) is in the center and \(\tau\) is an anti-involution.
Given a Wall anti-structure we can consider \(R\) as an \((R \otimes_{\ZZ} R)\)-module via the action \((a \otimes b)(c) = ab\tau(c)\), and endow it with an involution given by \(x \mapsto \eps\tau(x)\). Applying the Eilenberg-Mac Lane functor this results in an invertible \(k\)-module with involution over the \(\Eone\)-ring spectrum \(\GEM R\), whose underline \(\GEM R\)-module is \(\GEM R\), but which is generally not the one associated to any anti-involution on \(\GEM R\).
\end{example}

The following lemma gives a recognition criterion for \(k\)-modules with involution over \(A\) which come from anti-involutions of \(A\). Essentially, it reflects the idea that the datum of an anti-involution on \(A\) is equivalent to that of a perfect duality \(\Dual\colon \Modp{A} \to (\Modp{A})\op\) together with a symmetric Poincaré form \(u \in \map_A(A,\Dual(A))^{\hC}\) on the \(A\)-module \(A\):

\begin{proposition}
\label{proposition:recognition}%
Let \(A\) be an \(\Eone\)-algebra and \(M\) a \(k\)-module with involution over \(A\). Then \(M\) comes from an anti-involution on \(A\) as in Example~\refone{example:module-with-involution-associated-to-ring-with-involution} if and only if there exists a \(\Ct\)-equivariant map of spectra \(u\colon \SS \to M\) (where \(\Ct\)-acts trivially on \(\SS\)) such that the induced \(A\)-module map \(A \to M\) (using, say the first \(A\)-action on \(M\)), is an equivalence. In this case, the anti-involution on \(A\) can be recovered via the map \(A \to \map_A(M,M) \simeq \map_A(A,A) = A\op\) associated to the second \(A\)-action on \(M\).
\end{proposition}

The proof of Proposition~\refone{proposition:recognition} will require some preparation. Given a stable \(\infty\)-category \(\C\), the formation of mapping spectra allows one to consider \(\C\) as an \(\infty\)-category \emph{enriched} in spectra, see~\cite[Example 7.4.14]{Gepner-Haugseng}. This enrichment is functorial in \(\C\), that is, it can be organised into a functor
\begin{equation}
\label{equation:canonical-enrichement}%
\Catx \to \CatSpa \quad\quad \C \mapsto \C_{\Spa}
\end{equation}
where the latter is the \(\infty\)-category of \(\Spa\)-enriched \(\infty\)-categories, see~\cite[Proposition 4.10]{BGT} (and~\cite[Theorem 1.1]{haugseng2015rectification} for the comparison of the model categorical and \(\infty\)-categorical approaches to spectrally enriched categories). The functor~\eqrefone{equation:canonical-enrichement} is in fact fully-faithful and exhibits \(\Catx\) as an accessible localisation of \(\CatSpa\) by the collection of \defi{triangulated equivalences}, see~\cite[Theorem 4.22]{BGT}.

\begin{remark}
\label{remark:underlying}%
Every spectrally enriched \(\infty\)-category has an
``underlying \(\infty\)-category'' obtained by applying the functor \(\Om^{\infty} = \Map(\SS,-)\colon \Spa \to \Sps\) to all mapping spectra. In particular, the underlying \(\infty\)-category of the spectrally enriched category \(\C_{\Spa}\) is just \(\C\) itself (or rather, its image in \(\Cat\)). More formally, the functor~\eqrefone{equation:canonical-enrichement} constitutes a lift of the inclusion \(\Catx \hrar \Cat\) along the underlying \(\infty\)-category functor \(\CatSpa \to \Cat\).
\end{remark}

Let \(\Catxast := \Catx \times_{\Cat} (\Cat)_{\Del^0/}\) denote the \(\infty\)-category of stable \(\infty\)-categories \(\C\) equipped with a distinguished object \(\x \in \C\), and similarly, let \(\CatSpaPt = \CatSpa \times_{\Cat} (\Cat)_{\Del^0/}\) denote the \(\infty\)-category of spectrally enriched categories equipped with a distinguished object. Then we may consider the composite
\begin{equation}
\label{equation:endomorphism-spectrum}%
\Catxast \to \CatSpaPt \to \Alg_{\Eone}
\end{equation}
where the first functor is induced from~\eqrefone{equation:canonical-enrichement} and the second is the functor of~\cite[Theorem 6.3.2(iii)]{Gepner-Haugseng}, which can be described on objects as sending a pointed spectrally enriched category \((\D,\x)\) to the endomorphism spectrum \(\Map_{\D}(\x,\x)\). As shown in loc.\ cit.\, this functor has a fully-faithful left adjoint \(\B\colon \Alg_{\Eone} \to \CatSpaPt\) which sends a ring spectrum \(A\) to the pointed spectrally enriched category \((\B A,x)\) containing a single object \(\x\) whose endomorphism ring is \(A\). The essential image of \(\B\) is then given by those pointed spectrally enriched categories \((\D,\x)\) for which \(\x\) is the only object up to equivalence.

\begin{proposition}
\label{proposition:endomorphism-poincare}%
Let \((\C,\QF)\) be a Poincaré \(\infty\)-category and \((\x,\qone)\) a Poincaré object in \(\C\). Then the image \(\map_{\C}(\x,\x) \in \Alg_{\Eone}\) of \((\C,x)\) under~\eqrefone{equation:endomorphism-spectrum} inherits a canonical anti-involution, that it, it lifts to an object of \(\Alg_{\Eone}^{\Ct}\).
\end{proposition}
\begin{proof}
We first note that the spectral enrichment functor~\eqrefone{equation:canonical-enrichement} commutes with taking opposites. Indeed, since~\eqrefone{equation:canonical-enrichement} is fully-faithful and its essential image, spanned by the pre-triangulated spectrally enriched categories, is closed under opposites, the \(\mop\) action on \(\CatSpa\) induces a \(\Ct\)-action on \(\Catx\). This action then coincides with the action induced by the inclusion of \(\Catx\) in \(\Cat\) by Remark~\refone{remark:underlying} and the fact that the underlying \(\infty\)-category functor \(\CatSpa \to \Cat\) visibly commutes with opposites since it is induced by a functor \(\Spa \to \Sps\) on the level of enriching \(\infty\)-categories.

Now since \(\Del^0 \simeq (\Del^0)\op\) via an essentially unique isomorphism the \(\mop\) action on \(\Catx\) induces a \(\Ct\)-action on \(\Catxast\), which can be described on objects by the formula \((\C,\x)\mapsto (\C\op,\x)\). Similarly, the \(\mop\) action on \(\CatSpa\) induces a \(\Ct\)-action on \(\CatSpaPt\). In addition, the endomorphism functor \(\CatSpaPt \to \Alg_{\Eone}\) in~\eqrefone{equation:endomorphism-spectrum} is \(\Ct\)-equivariant with respect to the \(\mop\) actions on both sides, as can be seen by the fact that it admits a fully-faithful right adjoint \(\B\colon \Alg_{\Eone} \to \CatSpaPt\) which is itself compatible with taking opposites essentially by construction. Now consider the diagram
\[
\begin{tikzcd}
(\C,\x) \ar[r,phantom,"{\in}"] & \Catxast\ar[d]\ar[rr] && \Alg_{\Eone} & \map_{\C}(\x,\x) \ar[l,phantom,"{\ni}"] \\
\C \ar[r,phantom,"{\in}"] & \Catx &
\end{tikzcd}
\]
in which both arrows are \(\Ct\)-equivariant with respect to the \(\mop\) action. The perfect duality \(\Dual_{\QF}\) promotes \(\C\) to a \(\Ct\)-fixed object of \(\Catx\) (see~\eqrefone{equation:forget-Q-to-D}). The fibre of \(\Catxast \to \Catx\) over \(\C\) can be identified with \(\iota\C\), with its \(\Ct\)-action induced by \(\Dual_{\QF}\). To finish the proof it will suffice to show that \(\x \in \iota\C\) refines to a \(\Ct\)-fixed point. Indeed, this now follows from Corollary~\refone{corollary:fgt-hyp-equivariant}.
\end{proof}

\begin{proof}[Proof of Proposition~\refone{proposition:recognition}]
The only if direction is clear, since any anti-involution on \(A\) preserves the unit map \(u \colon \SS \to A\). To prove the other direction, suppose that \(M\) is a \(k\)-module with involution over \(A\) and we are given a \(\Ct\)-equivariant map \(u\colon \SS \to M\) with respect to the trivial \(\Ct\)-action on \(\SS\), such that the induced \(A\)-module map \(A \to M\) is an equivalence. Since \(u\) is \(\Ct\)-invariant it follows that this condition holds for both the first and second \(A\)-actions. On the other hand, it also holds for the \(A\op\)-action on \(A\) that the analogously defined map is an equivalence. The corresponding statement also holds for the \(\map_A(M,M)\)-action on \(M\). We hence obtain a commutative diagram
\begin{equation}
\label{equation:invertible}%
\begin{tikzcd}
A \ar[d,"\simeq"'] \ar[r] & \map_A(M,M) \ar[d,"\simeq"] \ar[r,phantom,"\simeq"] & A\op \\
M \ar[r,"\simeq"] & M &
\end{tikzcd}
\end{equation}
in which both vertical maps (induced by the base point \(u\colon \SS \to M\)) are equivalences, and the bottom horizontal map is the involution on \(M\). It then follows that the top horizontal map \(A \to \map_A(M,M)\) is an equivalence as well. In particular, \(M\) is invertible and hence the induced bilinear functor on \(\Modp{A}\) is perfect by Lemma~\refone{proposition:invertible-poincare}, with associated duality \(X \mapsto \map_A(X,M)\) on \(\Modp{A}\). The \(\Ct\)-equivariant map \(u\colon \SS \to M\) then determines a form
\[
q_u \in \QF^{\sym}_M(A) = \map_A(A \otimes_k A,M)^{\hC}
\]
which is Poincaré by the condition that the induced map \(A \to M = \Dual(A)\) is an equivalence. We may then identify the top horizontal equivalence in~\eqrefone{equation:invertible} as the underlying equivalence of the anti-involution on \(\map_A(A,A) = \map_A(\Dual(A),\Dual(A))\op\) induced by the Poincaré form \(q_u\) by Proposition~\refone{proposition:endomorphism-poincare}.
\end{proof}

\subsection{Modules with genuine involution}
\label{subsection:genuine-modules}%

Our goal in this section is to refine the definition of a \(k\)-module with involution studied in \S\refone{subsection:modules-with-involution} above in order to obtain a notion capable of encoding not just bilinear but also quadratic functors on \(\Modp{A}\), and similarly on \(\Mod^c_A\) for a given subgroup \(c \subseteq \K_0(\Modp{A})\).
To begin, recall that for a spectrum \(X\) there is a canonical map \(X \to (X \otimes_{\SS} X)^{\tC}\), known as the \emph{Tate diagonal}, which enjoys a variety of favourable formal properties, see~\cite{NS}. If \(X\) is now a \(k\)-module then we can consider the composite map
\begin{equation}
\label{equation:linear-tate}%
X \to (X \otimes_{\SS} X)^{\tC} \to (X \otimes_k X)^{\tC},
\end{equation}
where the second map is induced by the lax monoidal structure of the forgetful functor \(\Mod_k \to \Spa\). For \(X=k\) the map~\eqrefone{equation:linear-tate} gives the composite map
\[
\Fr\colon k \to (k \otimes_{\SS} k)^{\tC} \to k^{\tC}
\]
which is known as the \emph{Tate Frobenius map}. For a \(k\)-module spectrum \(X\) we may then consider the spectrum \((X \otimes_k X)^{\tC}\), which is naturally a \(k^{\tC}\)-module, as a \(k\)-module spectrum, by restricting structure along the Tate Frobenius. With this \(k\)-module structure the map~\eqrefone{equation:linear-tate} becomes \(k\)-linear, and we will henceforth refer to it as the \defi{\(k\)-linear Tate diagonal}.

\begin{remark}
\label{remark:tate-diagonal-omega}%
The functor \(\id\colon \Spa \to \Spa\) is stably corepresented by \(\SS\), and so for any exact functor \(F\colon \Spa \to \Spa\), the space of natural transformations \(\id \Rightarrow F\) identifies, via evaluation at \(\SS\), with the space \(\Om^{\infty} F(\SS)\) of maps \(\SS \to F(\SS)\). Under this identification, the Tate diagonal \(\id \Rightarrow ((-) \otimes (-))^{\tC}\) corresponds
to the image in \(\Om^{\infty}[\SS^{\tC}]\) of the unit \(e \in [\Om^{\infty}\SS]^{\hC} = \Om^{\infty}[\SS^{\hC}]\) (where we consider \(\Om^{\infty}\SS = \Map(\SS,\SS)\) as a monoid with respect to composition, and the \(\Ct\)-equivariant structure, as the \(\Ct\)-action on \(\SS\), is constant).
On the other hand, the functor \(\Om^{\infty}\colon \Spa \to \Sps\) is the space valued functor corepresented by \(\SS\), and hence natural transformations \(\Om^{\infty} \Rightarrow G\) identify with maps \(\ast \to G(\SS)\) for \emph{any} functor \(G \colon \Spa \to \Sps\). The point \(e \in [\Om^{\infty}\SS]^{\hC}\) then determines
a factorization of \(\Om^{\infty}\Del\colon \Om^{\infty}(-) \Rightarrow \Om^{\infty}[((-) \otimes (-))^{\tC}]\) as
\[
\Om^{\infty}(-) \Rightarrow [\Om^{\infty}((-) \otimes (-))]^{\hC} = \Om^{\infty}[((-) \otimes (-))^{\hC}] \Rightarrow \Om^{\infty}[((-) \otimes (-))^{\tC}] .
\]
Furthermore, since \(e\) is a monoid unit it lifts along \(\Om^{\infty}\SS \times \Om^{\infty}\SS  \to \Om^{\infty}\SS\) to the point \((e,e)\), and this can done also on the level of \(\Ct\)-homotopy fixed points where \(\Ct\) acts on the domain by flipping the components. We conclude that the above factorization of \(\Om^{\infty}\Del\colon \Om^{\infty}(-) \Rightarrow \Om^{\infty}[((-) \otimes (-))^{\tC}]\) refines to a factorization
\[
\Om^{\infty}(-) \Rightarrow [\Om^{\infty}(-) \times \Om^{\infty}(-)]^{\hC} \Rightarrow [\Om^{\infty}((-) \otimes (-))]^{\hC} = \Om^{\infty}[((-) \otimes (-))^{\hC}] \Rightarrow \Om^{\infty}[((-) \otimes (-))^{\tC}] ,
\]
where the first map is simply the \(\Ct\)-equivariant diagonal.
We will use this observation in the proof of Lemma~\refone{equation:tate-M} below.
\end{remark}

\begin{warning}
\label{warning:frobenius}%
Since \(k\) is the unit of \(\Mod_k\) the flip \(\Ct\)-action on \(k \otimes_k k \simeq k\) is trivial. However, the Tate Frobenius \(k \to k^{\tC}\) is generally not equivalent to the composite \(k \to k^{\hC} \to k^{\tC}\). In particular, if we were to endow \((X \otimes_k X)^{\tC}\) with the \(k\)-module structure restricted from its \(k^{\tC}\)-module structure along \(k \to k^{\hC} \to k^{\tC}\) (which would be the \(k\)-module structure we would obtain by applying to \(X \otimes_k X\) the Tate construction internally in \(\Mod_k\)) then~\eqrefone{equation:linear-tate} would generally not be \(k\)-linear.
\end{warning}

\begin{definition}
\label{definition:module-with-genuine-involution}%
Let \(A\) be a \(\Eone\)-algebra.
A \defi{module with genuine involution over \(A\)} is a triple \((M,N,\alp)\) which consists of
\begin{enumerate}
\label{definition:modgeninv}%
\item[-] a \(k\)-module with involution \(M\) over \(A\) in the sense of Definition~\refone{definition:module-with-involution},
\item[-] an \(A\)-module \(N\), and
\item[-] an \(A\)-linear map \(f\colon N \to M^{\tC}\).
\end{enumerate}
Here we view \(M^{\tC}\), which is canonically an \((A\otimes_k A)^{\tC}\)-module, as an \(A\)-module through the \(k\)-linear Tate diagonal \(A \to (A\otimes_k A)^{\tC}\).
\end{definition}

\begin{remark}
\label{remark:NormModules}%
When \(k=\SS\) the data of a module with genuine involution over \(A\) can equivalently be described as a genuine \(\Ct\)-spectrum equipped with an action of the Hill-Hopkins-Ravenel norm \(\N A\) of \(A\), first introduced in~\cite{HHR} and later reinterpreted from the \(\infty\)-categorical perspective in~\cite{barwick2016parametrized}. The genuine \(\Ct\)-spectrum \(\N A\) has as underlying \(\Ct\)-spectrum \(A\otimes_{\SS} A\) with the flip \(\Ct\)-action, geometric fixed points \(A\), and reference map \(A \to (A\otimes_{\SS} A)^{\tC}\) the Tate-diagonal considered in~\cite{NS}. It is an algebra object with respect to the symmetric monoidal structure on genuine \(\Ct\)-spectra, and the data of a module over this algebra object consists exactly of an \((A\otimes_{\SS} A)\)-module with \(\Ct\)-action (which is \(M\) in our case), a module over the geometric fixed points \(A\) (which is \(N\) in our case), and a reference map \(\alp \colon N \to M^{\tC}\) refining \(M\) to a genuine \(\Ct\)-spectrum with geometric fixed points \(N\).

More generally, for an arbitrary base \(\Einf\)-ring spectrum \(k\) and \(\Eone\)-algebra \(A\) over \(k\), the notion of a \(k\)-module with genuine involution over \(A\) is equivalent to that of a genuine \(\Ct\)-spectrum equipped with an action of the \(k\)-linear norm \(\N_k A\) of \(A\): this is the genuine \(\Ct\)-spectrum whose underlying \(\Ct\)-spectrum is \(A\otimes_k A\) with the flip \(\Ct\)-action, whose geometric fixed point spectrum is \(A\), and whose reference map \(A \to (A\otimes_k A)^{\tC}\) is the \(k\)-linear Tate-diagonal mentioned above.
\end{remark}

\begin{lemma}
\label{lemma:tate-M}%
For \(M \in \Mod_{A \otimes_k A}^{\hC}\) and \(X \in \Modp{A}\) there is an equivalence
\begin{equation}
\label{equation:tate-M}%
\map_{A \otimes_k A}(X \otimes_k X,M)^{\tC} \simeq \map_A(X,M^{\tC})
\end{equation}
natural in \(M\) and \(X\).
\end{lemma}
\begin{proof}
Consider the functor \(F\colon (\Modp{A})^{\op} \to \Spa\) given by
\(X \mapsto \map_{A \otimes_k A}(X \otimes_k X,M)^{\tC}\).
This functor is exact, thus by Morita theory it is of the form \(\map_A(X, N)\) for some \(A\)-module \(N\). Setting \(X = A\) we find that \(N = M^{\tC}\) as a spectrum. Furthermore, the right action of \(\Omega^\infty A\) on \(A\) translates under the functor \(F\) to the diagonal action of \(\Omega^\infty A\) on \(M^{\tC}\), i.e., the action through the composite
\[
\Omega^{\infty} A \xrightarrow{\Delta} (\Omega^\infty A \times \Omega^\infty A)^{\hC} \to \Omega^\infty (A \otimes_{\SS} A)^{\hC} \to \Omega^\infty (A \otimes_{\SS} A)^{\tC} \to (A \otimes_k A)^{\tC},
\]
which at the same time also underlines the \(k\)-linear Tate diagonal by Remark~\refone{remark:tate-diagonal-omega}. Applying the same observation to shifts of \(A\) and using exactness then fully identifies the \(A\)-module structure on \(M^{\tC}\), showing the claim.
\end{proof}

\begin{construction}
\label{construction:functors-associated-to-module-with-genuine-involution}%
Let \((M,N,\alp)\) be a \(k\)-module with genuine involution over \(A\). We define a quadratic functor \(\QF^{\alp}_M\) on perfect \(A\)-modules by the pullback
\begin{equation}
\label{equation:genuine-defining-pullback-square}%
\begin{tikzcd}
	\QF^{\alp}_M(X) \ar[r] \ar[d] & \map_A(X,N) \ar[d] \\
	\QF^\sym_M(X) \ar[r] & \map_A(X,M^{\tC})
\end{tikzcd}
\end{equation}
where the lower horizontal map is given by the composite
\[
\QF^\sym_M(X)= \map_{A\otimes_k A}(X\otimes_k X,M)^{\hC} \to \map_{A\otimes_k A}(X\otimes_k X,M)^{\tC} \simeq \map_A(X,M^{\tC}),
\]
where the last equivalence is via Lemma~\refone{lemma:tate-M}.

By construction the underlying bilinear part of \(\QF^{\alp}_M\) is \(\Bil_M\), and hence the condition that \(\QF^{\alp}_M\) is Poincaré depends only on \(M\), via the criterion of Proposition~\refone{proposition:invertible-poincare}. In addition, by Lemma~\refone{lemma:equivariance} the trace map \(\QF^{\qdr}_M \to \QF^{\sym}_M\) factors canonically as
\[
\QF^{\qdr}_M\longrightarrow \QF^{\alp}_M\longrightarrow \QF^{\sym}_M.
\]
If \(M^{\tC}=0\), for example if \(A\) is an  \(\SS[\tfrac{1}{2}]\)-algebra, then \(\QF^{\alp}(X) \simeq \QF_{M}^\sym(X) \oplus \map_A(X,N)\).
\end{construction}

\begin{example}
\label{example:symquadgenuine}%
Let \(M\) be a \(k\)-module with involution over an \(\Eone\)-algebra \(A\).
The \(k\)-modules with genuine involution \((M,0,0\to M^{\tC})\) and \((M,M^{\tC}, M^{\tC}=M^{\tC})\) give rise respectively to the quadratic and symmetric functors \(\QF^{\qdr}_M\) and \(\QF^{\sym}_M\) of Definition~\refone{definition:functors-associated-to-module-with-involution}.
\end{example}

\begin{example}
\label{example:truncation}%
Let \(M\) be a \(k\)-module with involution over an \(\Eone\)-algebra \(A\) and assume that \(A\) is connective (so that the truncation of an \(A\)-module admits a canonical \(A\)-module structure). Then there is for every \(m \in \ZZ\) a \(k\)-module with genuine involution given by
\[
(M, \tau_{\geq m} M^{\tC},  \tau_{\geq m} M^{\tC} \to M^{\tC})
\]
with the reference map being the \(m\)-connective cover. Applying Construction~\refone{construction:functors-associated-to-module-with-genuine-involution} these give rise to quadratic functors \(\QF_M^{\geq m}\) that sit between the quadratic and the symmetric one, i.e., there are maps
\begin{equation}
\label{equation:tower-truncations}%
\QF_M^{\qdr} \to \cdots \to \QF_M^{\geq 1} \to \QF_M^{\geq 0} \to \QF_M^{\geq -1} \to \cdots \to  \QF_M^{\sym} \ .
\end{equation}
If \(M^{\tC}\) vanishes, e.g.\ if 2 is invertible in \(A\), then all of these maps are equivalences. The limit and colimit of these diagrams of Poincaré structures are given by
\[
\lim[\cdots \rightarrow \QF^{\geq m}_{M}\rightarrow \QF^{\geq m-1}_M\rightarrow\cdots]\simeq \QF^{\qdr}_M
\]
and
\[
\colim[\cdots \rightarrow \QF^{\geq m}_{M}\rightarrow \QF^{\geq m-1}_M\rightarrow\cdots] \simeq\QF^{\sym}_M.
\]
Indeed, inspecting the defining pullback squares~\eqrefone{equation:genuine-defining-pullback-square} for \(\QF^{\geq m}_M\) and using that pullbacks and mapping spectra out of compact \(A\)-modules respect both limits and colimits this follows
from the fact that \(\lim\tau_{\geq m} M^{\tC} = 0\) while the induced map \(\colim \tau_{\geq m} M^{\tC} \to M^{\tC}\) is an equivalence. We may hence consider the tower~\eqrefone{equation:tower-truncations}
as \defi{interpolating} between the quadratic and symmetric Poincaré structures on \(\Modp{A}\). We will study this construction in further detail in \S\refone{subsection:discrete-rings} in the case where \(A\) is (the Eilenberg-Mac Lane spectrum of) an ordinary ring.
\end{example}

\begin{example}
\label{example:anti-involution}%
Let \(A\) be a \(\Eone\)-algebra with anti-involution. In Example~\refone{example:module-with-involution-associated-to-ring-with-involution} we have seen that \(A\) can be consider as an invertible \(k\)-module with involution over itself.
In order to promote \(A\) to a \(k\)-module with genuine involution we need an \(A\)-module \(N\) and a map \(N \to A^{\tC}\) of \(A\)-modules. Such a triple \((A,N,\alp)\) is called an \emph{\(\Eone\)-algebra with genuine anti-involution}.
Any such \(\Eone\)-algebra with genuine anti-involution has an underlying genuine \(\Ct\)-spectrum, the module \(N\) taking the role of the geometric fixed points. Since \(A\) is invertible as a \(k\)-module with involution by Proposition~\refone{proposition:invertible-poincare} and Remark~\refone{remark:invertible-poincare} imply that the associated hermitian structure \(\QF^{\alp}_{A}\) on \(\Modp{A}\) and on \(\Modf{A}\) are both Poincaré.
\end{example}

\begin{example}
\label{example:orthogonal-spectra}%
Let \(A\) be an orthogonal ring spectrum with anti-involution in the sense of~\cite{DMPR}. This gives rise to a genuine \(\Ct\)-spectrum whose underlying spectrum with \(\Ct\)-action is the underlying spectrum of \(A\),
whose geometric fixed points \(A^{\geofix}\) is canonically an \(A\)-module, and where the map \(\alp\colon A^{\geofix} \to A^{\tC}\) is \(A\)-linear.
We therefore obtain a ring spectrum with genuine anti-involution \((A,A^{\geofix},\alp)\).
\end{example}

\begin{example}
\label{example:sphere}%
Consider the sphere spectrum \(\SS\) with the trivial \(\Ct\)-action. We may then view \(\SS\) as an associative ring spectrum with anti-involution, and refine it to a ring spectrum with genuine anti-involution using as reference map the composite \(\SS \to \SS^{\hC} \to \SS^{\tC}\), which also agrees in this case with the Tate diagonal. The Poincaré structure associated to this genuine anti-involution on \(\SS\) is the
universal Poincaré structure \(\QF^{\uni}\) of Example~\refone{example:universal-category}.
\end{example}

\begin{example}
\label{example:tate-structure}%
Let \(A\) be an  \(\Einf\)-\(k\)-algebra equipped with a \(\Ct\)-action. We may then view \(A\) as an associative \(k\)-algebra with anti-involution, and refine it to a \(k\)-algebra with genuine anti-involution using as reference map the composite
\[
\tate\colon A \to (A \otimes_k A)^{\tC} \to A^{\tC},
\]
where the first map is the \(k\)-linear Tate diagonal and the second map is induced by the \(\Ct\)-equivariant commutative \(k\)-algebra map \(A \otimes_k A \to A\) whose restriction to the first component is the identity \(A \to A\) and restriction to the second component is given by the generator of the \(\Ct\)-action.
This yields a Poincaré structure on \(\Modp{A}\), which we denote by \(\QF^{\tate}_{A}\), and refer to as the \defi{Tate Poincaré structure} associated to the given \(\Ct\)-action. The universal Poincaré structure \(\QF^{\uni}\) on \(\Spaf\) then corresponds to the case where \(A=k=\SS\) and the \(\Ct\)-action is trivial.
\end{example}

\begin{example}
Given an \(\Eone\)-algebra \(A\) we may form the \(\Eone\)-algebra with anti-involution \(A \oplus A\op\) where the involution flips the two factors. We may then refine this involution to a genuine one by taking the zero map \(0 \to (A \oplus A\op)^{\tC} = 0\), which at the same time is also an equivalence. The resulting Poincaré structure \(\QF^0_{A \oplus A\op}\) on \(\Modp{A \oplus A\op}\) is then both quadratic and symmetric (see Example~\refone{example:symquadgenuine}). We now note that the projections \(A \oplus A\op \to A\) and \(A \oplus A\op \to A\op\) induce an equivalence
\[
\Modp{A \oplus A\op} \to \Modp{A} \times \Modp{A\op} \simeq \Modp{A} \times (\Modp{A})\op,
\]
under which the Poincaré structure in question corresponds to the \emph{hyperbolic structure} of Definition~\refone{definition:hyperbolic-cat}, and so
\[
(\Modp{A \oplus A\op},\QF^0_{A \oplus A\op}) \simeq \Hyp(\Modp{A}) .
\]
\end{example}

\subsection{Classification of hermitian structures}
\label{subsection:classification-in-modules}

We shall henceforth focus on the case where \(k=\SS\) is the sphere spectrum.
We will show that in this case, modules with genuine involution do not only provide a convenient way of producing hermitian structures on \(\Mod_{A}^\omega\)
but that these two notions become in fact equivalent. To formulate this more precisely, we proceed to organise modules with genuine involution over \(A\) into an \(\infty\)-category \(\Mod_{\N A}\), defined as the pullback
\begin{equation}
\label{equation:genuine-modules}%
\begin{tikzcd}
\Mod_{\N A} \ar[r] \ar[d] & \Ar(\Mod_A)\ar[d,"\target"]
\\
\Mod_{A\otimes_{\SS} A}^{\hC} \ar[r,"{(-)^{\tC}}"] & \Mod_A
\end{tikzcd}
\end{equation}
of the arrow category \(\Ar(\Mod_A)\) and the \(\infty\)-category \((\Mod_{A\otimes_{\SS} A})^{\hC}\) of modules with involution. The right vertical map is the projection onto the target, and the bottom horizontal map sends a module with involution \(M\) to the Tate construction \(M^{\tC}\), considered as an \(A\)-module via the Tate diagonal map \(A\to (A\otimes_{\SS} A)^{\tC}\).

We would like to relate \(\Mod_{\N A}\) and \(\Funq(\Mod^{\omega}_{A})\) by constructing a commutative diagram
\begin{equation}
\label{equation:cube-mod-involution}%
\begin{tikzcd}
& \Mod_{\N A} \arrow[dl] \arrow[rr] \arrow[dd] & & \Ar(\Mod_A) \arrow[dl] \arrow[dd] \\
\Funq(\Modp{A}) \arrow[rr, crossing over] \arrow[dd] & & \Ar(\Funx((\Modp{A})^{\op},\Spa)) \\
& \Mod_{A\otimes_{\SS} A}^{\hC} \arrow[dl] \arrow[rr] &  & \Mod_A \arrow[dl] \\
\Funs(\Modp{A}) \arrow[rr] & & \Funx((\Modp{A})^{\op},\Spa) \arrow[from=uu, crossing over]\\
\end{tikzcd}
\end{equation}
where the front square is the pullback square of Corollary~\refone{corollary:classification-of-quad-functors} which exhibits the analogous decomposition of the \(\infty\)-category of quadratic functors \(\Funq(\Mod^{\omega}_{A})\) into linear and bilinear parts. We then define the right face of~\eqrefone{equation:cube-mod-involution} to be the square induced by the Yoneda map
\[
\Mod_A \to \Funx((\Modp{A})\op,\Spa) \quad\quad N \mapsto \map_A(-,N),
\]
and the bottom arrow in the left face is the functor
\[
\Mod_{A \otimes_{\SS} A}^{\hC} \to \Funs(\Modp{A}) \quad\quad M \mapsto \map_{A \otimes_{\SS} A}((-)\otimes_{\SS}(-),M)
\]
introduced above. The commuting homotopy in the bottom square of~\eqrefone{equation:cube-mod-involution} is then provided by Lemma~\refone{lemma:tate-M}. The cube is then uniquely determined by the fact that its front face is a cartesian square.

\begin{theorem}
\label{theorem:classification-genuine-modules}%
The cube~\eqrefone{equation:cube-mod-involution}, considered as a natural transformation from its back face to its front face, is an equivalence. In particular, the resulting arrow
\[
\Mod_{\N A} \to \Funq(\Modp{A})
\]
is an equivalence of \(\infty\)-categories, whose action on objects is given by \((M,N,\alp) \mapsto \QF^{\alp}_M\).
\end{theorem}

\begin{remark}
\label{remark:k-linear-classification}%
For a general \(k\) one can still define the pullback of \(\infty\)-categories
\[
\Mod_{\N_k A} = \Mod_{A \otimes_k A}^{\hC} \times_{\Mod_A} \Ar(\Mod_A),
\]
and construct a functor
\begin{equation}
\label{equation:mod-k-NA-to-funq}%
\Mod_{\N_k A} \to \Funq(\Modp{A}),
\end{equation}
which in general will not be an equivalence. Here we use the notation \(\Mod_{\N_k A}\) to indicate that we may think of this as the \(\infty\)-category of modules in genuine \(\Ct\)-spectra over the \(k\)-linear norm \(\N_k A\) of \(A\), see Remark~\refone{remark:NormModules}. Alternatively,
we may also identify \(\Mod_{\N_k A}\) with the \(\infty\)-category of hermitian structures on \(\Modp{A}\) equipped with an action of the Tate hermitian structure \(\QF^{\tate}_k \in \Funq(\Modp{k})\) of Example~\refone{example:tate-structure}. Here we have that \(\Modp{k}\) acts on \(\Modp{A}\) by the fact that \(A\) is a \(k\)-algebra, and this induces an action of \(\Funq(\Modp{k})\) on \(\Funq(\Modp{A})\) by means of Day convolution, see \S\refone{subsection:symmetric-monoidal-poincare} and Example~\refone{example:k-linear-structures}.
The functor~\eqrefone{equation:mod-k-NA-to-funq} then corresponds to forgetting the \(\QF^{\tate}_k\)-module structure.
\end{remark}

\begin{proof}[Proof of Theorem~\refone{theorem:classification-genuine-modules}]
Since the back and front faces are both cartesian squares it will suffice to show that its right and back faces determine equivalences from their back edge to their front edge. For this, it will suffice to show that the functors
\[
\Mod_A \to \Funx((\Modp{A})\op,\Spa) \quad\text{and}\quad \Mod_{A \otimes_{\SS} A}^{\hC} \to \Funs(\Modp{A})
\]
are equivalences. For the former, we note that since \(\Spa\) is stable we have that post-composition with \(\Om^{\infty}\colon \Spa \to \Sps\) induces an equivalence
\[
\Funx((\Modp{A})\op,\Spa) \simeq \Fun^{\rex}((\Modp{A})\op,\Sps) \simeq \Ind(\Modp{A})
\]
and consequently the claim follows from the fact that \(\Mod_A\) is generated by \(\Modp{A}\) under colimits.
For the second map, by its construction it will suffice to show that the functor
\begin{equation}
\label{equation:bilinear-modules}%
\Mod_{A \otimes_{\SS} A} \to \Funb(\Modp{A}) \quad\quad M \mapsto \map_{A \otimes_{\SS} A}((-)\otimes_{\SS}(-),M)
\end{equation}
is an equivalence of \(\infty\)-categories. For this, recall that by~\cite[Theorem 4.8.5.16]{HA} and~\cite[Remark 4.8.5.19]{HA} the association \(A \mapsto \Mod_A\) refines to a symmetric monoidal functor
\[
\Theta_{\Spa}\colon \Alg_{\Eone} = \Alg_{\Eone}(\Spa) \to \LMod_{\Spa}(\PrL)
\]
from \(\Eone\)-ring spectra to the \(\infty\)-category of \(\Spa\)-modules in presentable \(\infty\)-categories, and the latter can be identified with the full subcategory of \(\PrL\) spanned by the stable presentable \(\infty\)-categories. In particular, the bilinear functor
\[
\Mod_{A} \times \Mod_{A} \to \Mod_{A \otimes_{\SS} A} \quad\quad (X,Y) \mapsto X \otimes_{\SS} Y
\]
induces an equivalence
\begin{equation}
\label{equation:A-tensor-A}%
\Mod_{A} \otimes_{\Spa} \Mod_{A} \xrightarrow{\simeq} \Mod_{A \otimes_{\SS} A},
\end{equation}
where \(\otimes_{\Spa}\) denotes the tensor product of stable presentable \(\infty\)-categories (which can also be viewed as \(\Spa\)-modules in presentable \(\infty\)-categories). Since \(\Mod_A\) and \(\Mod_{A \otimes_{\SS} A}\) are compactly generated by \(\Modp{A}\) and \(\Modp{A \otimes_{\SS} A}\) respectively and the bilinear functor \((X,Y) \mapsto X \otimes_{\SS} Y\) maps a pair of compact \(A\)-modules to a compact \((A \otimes_{\SS} A)\)-module we see that the equivalence~\eqrefone{equation:A-tensor-A} is induced on Ind-categories by the functor
\begin{equation}
\label{equation:A-tensor-A-perf}%
\Modp{A} \otimes \Modp{A} \to \Modp{A \otimes_{\SS} A},
\end{equation}
induced by the same bilinear functor \((X,Y) \mapsto X \otimes_{\SS} Y\), where \(\otimes\) is the tensor product of stable \(\infty\)-categories, and we use the fact that \(\Ind(-)\) is symmetric monoidal. It then follows that restriction along~\eqrefone{equation:A-tensor-A-perf} induces an equivalence
\[
\Funx((\Modp{A \otimes_{\SS} A})\op,\Spa) \xrightarrow{\simeq} \Funx((\Modp{A} \otimes \Modp{A})\op,\Spa) \simeq
\]
\[
\Funx((\Modp{A})\op \otimes (\Modp{A})\op,\Spa) \simeq \Funb(\Modp{A}).
\]
Since the Yoneda map \(\Mod_{A \otimes_{\SS} A} \to \Funx((\Modp{A \otimes_{\SS} A})\op,\Spa)\) is an equivalence by the argument above we may now conclude that~\eqrefone{equation:bilinear-modules} is an equivalence, and so the proof is complete.
\end{proof}

\begin{remark}
\label{remark:finite}%
As explained in the introduction of the present section,
left Kan extension and restriction determine inverse equivalences between
\(\Funq(\Mod^c_A)\) and \(\Funq(\Modp{A})\) for every subgroup \(c \subseteq \K_0(\Modp{A})\).
Theorem~\refone{theorem:classification-genuine-modules} then implies that the association \((M,N,\alp) \mapsto \QF^{\alp}_M\) also determines an equivalence \(\Mod_{\N A} \simeq \Funq(\Mod^c_A)\) for every \(c \subseteq \K_0(\Modp{A})\).
\end{remark}

\begin{remark}
\label{remark:classification-poincare}%
It follows from Theorem~\refone{theorem:classification-genuine-modules} and Proposition~\refone{proposition:invertible-poincare} that the association \((M,N,\alp) \mapsto \QF^{\alp}_M\) determines a bijective correspondence  between equivalence classes of \emph{invertible} modules with genuine involution \((M,N,\alp) \in \Mod_{\N A}\) and equivalence classes of Poincaré structures on \(\Modp{A}\). More generally, for a given subgroup \(c \subseteq \K_0(\Modp{A})\) we obtain a classification of Poincaré structures on \(\Mod^c_A\) via invertible modules with genuine involution \((M,N,\alp)\) such that the associated duality \(\Dual_M = \map_A(-,M)\) on \(\Modp{A}\) preserves \(\Mod^c_A\) (or, equivalently, preserves \(c \subseteq \K_0(\Modp{A})\)). In the case of \(\Modf{A}\) the last condition amounts to \(M\) itself being in \(\Modf{A}\), see Remark~\refone{remark:invertible-poincare}.
We note that under these correspondences maps of Poincaré structures (see Definition~\refone{definition:poinc-cats}) correspond to those maps \((M,N,\alp) \to (M',N',\alp')\) in \(\Mod_{\N A}\) for which the map \(M \to M'\) is an equivalence.
\end{remark}

\subsection{Restriction and induction}
\label{subsection:restriction-induction}%

In the present subsection we keep the assumption that \(k = \SS\) and consider the functorial dependence of \(\Mod_{\N A}\) in \(A\), and the compatibility of this functoriality with the one for hermitian structures along the equivalence of Theorem~\refone{theorem:classification-genuine-modules}.
Let \(\phi\colon A \to B\) be a map \(\Eone\)-algebras and let
\begin{equation}
\label{equation:induction}%
p_{\phi}\colon \Modp{A} \to \Modp{B}
\end{equation}
be the induction functor sending \(X\) to \(B \otimes_A X\). Then the restriction functor
\[
p_{\phi}^*\colon \Funq(\Modp{B}) \to \Funq(\Modp{A})
\]
corresponds, under the equivalence of Theorem~\refone{theorem:classification-genuine-modules}, to a functor
\begin{equation}
\label{equation:restriction-genuine}%
\phi^*\colon \Mod_{\N B} \to \Mod_{\N A},
\end{equation}
which we consider as the restriction of structure operation for modules with genuine involution. As explained in \S\refone{subsection:classification}, restriction of quadratic functors commutes with taking linear and bilinear parts and with the formation of symmetric Poincaré structures, that is, it acts compatibly on the entire pullback square
\begin{equation}
\label{equation:classification-funq}%
\begin{tikzcd}
\Funq(\C) \ar[r,"\tau"] \ar[d,"\Bil"'] & \Ar(\Funx(\C\op,\Spa)) \ar[d,"{\target}"] \\
\Funs(\C) \ar[r]                  & \Funx(\C\op,\Spa).
\end{tikzcd}
\end{equation}
Under the equivalence of Theorem~\refone{theorem:classification-genuine-modules} we obtain the same for the restriction functor~\eqrefone{equation:restriction-genuine}, that is, it extends to the entire defining squares~\eqrefone{equation:genuine-modules} for \(A\) and \(B\). On the other hand, the Yoneda equivalences \(\Mod_A \simeq \Funx((\Modp{A})\op,\Spa)\) and \(\Mod_{B} \simeq \Funx((\Modp{B})\op,\Spa)\) fit into a commutative square
\begin{equation}
\label{equation:forget-restrict}%
\begin{tikzcd}
\Mod_{B} \ar[d,"\simeq"] \ar[r,"{\phi^*}"] & \Mod_A \ar[d,"{\simeq}"] \\
\Funx((\Modp{B})\op,\Spa) \ar[r,"{p_{\phi}^*}"] & \Funx((\Modp{A})\op,\Spa)
\end{tikzcd}
\end{equation}
in which the top horizontal arrow is the forgetful functor from \(B\)-modules to \(A\)-modules and the bottom horizontal functor is restriction along \(p_{\phi}\). Indeed, the commutativity is given by the adjunction equivalence \(\map_B(p_{\phi}X,M) \simeq \map_A(X,\phi^*M)\). Similarly, the equivalences \(\Mod_{A \otimes_{\SS} A} \simeq \Funb(\Modp{A})\) and \(\Mod_{B \otimes_{\SS} B} \simeq \Funb(\Modp{B})\) fit into a commutative square
\begin{equation}
\label{equation:forget-restrict-2}%
\begin{tikzcd}
[column sep=9ex]
\Mod_{B \otimes_{\SS} B} \ar[d,"{\simeq}"] \ar[r,"{(\phi \otimes \phi)^*}"] & \Mod_{A \otimes_{\SS} A} \ar[d,"{\simeq}"] \\
\Funb(\Modp{B}) \ar[r,"{(p_{\phi} \times p_{\phi})^*}"] & \Funb(\Modp{A})
\end{tikzcd}
\end{equation}
We may thus conclude that for \((M,N,\alp) \in \Mod_{\N B}\) the restriction functor~\eqrefone{equation:restriction-genuine} is obtained by simply restricting the \((B \otimes_{\SS} B)\)-module structure on \(M\) to \(A \otimes_{\SS} A\), the \(B\)-module structure on \(N\) to \(A\), and viewing \(\alp\) as a map of \(A\)-modules by forgetting its compatibility with the \(B\)-module structures on its domain and codomain.

We now proceed to discuss how the operation of left Kan extension is mirrored along the equivalence of Theorem~\refone{theorem:classification-genuine-modules}. Recall that by Lemma~\refone{lemma:kan-extension-exact-quadratic} the operation of left Kan extensions \(\Fun((\Modp{A})\op,\Spa) \to \Fun((\Modp{B})\op,\Spa)\) preserves quadratic functors, and the resulting functor
\begin{equation}
\label{equation:left-kan-phi}%
(p_{\phi})_!\colon\Funq(\Modp{A}) \to \Funq(\Modp{B})
\end{equation}
is compatible with taking linear and bilinear parts, that is
\[
\Bil_{(p_{\phi})_!\QF} \simeq (p_{\phi} \times p_{\phi})_!\Bil_{\QF} \quad\text{and}\quad \Lin_{(p_{\phi})_!\QF} \simeq (p_{\phi})_!\Lin_{\QF}.
\]
By Lemma~\refone{lemma:kan-extension-exact-quadratic} and Corollary~\refone{corollary:structure-map-left-kan} we may then conclude the following:
\begin{corollary}
\label{corollary:induction-genuine-modules}%
Under the equivalence of Theorem~\refone{theorem:classification-genuine-modules}, the left Kan extension functor~\eqrefone{equation:left-kan-phi} corresponds to the functor
\[
\phi_!\colon \Mod_{\N A} \to \Mod_{\N B}
\]
sending a module with genuine involution \((M,N,\alp) \in \Mod_{\N A}\) to the module with genuine involution
\[
\phi_!(M,N,\alp) = (p_{\phi \otimes \phi}M,p_{\phi}N,\phi_!\alp) = ((B \otimes_{\SS} B) \otimes_{A \otimes_{\SS} A} M, B \otimes_A N, \phi_!\alp) \in \Mod_{\N B}
\]
where \(\phi_!\alp\) is given by the composite
\[
B\otimes_A N\to B \otimes_A M^{\tC}\to ((B \otimes_{\SS} B) \otimes_{A \otimes_{\SS} A}M)^{\tC}
\]
of \(B \otimes_A \alp\) and the Beck-Chevalley transformation on the lax commuting square on the right
\begin{equation}
\label{equation:tate-restriction-modules}%
\begin{tikzcd}
(\Mod_{B \otimes_{\SS} B})^{\hC} \ar[r,"{(-)^{\tC}}"] \ar[d] & \Mod_{B} \ar[d] \\
(\Mod_{A \otimes_{\SS} A})^{\hC} \ar[r,"{(-)^{\tC}}"] & \Mod_A
\end{tikzcd}
\qquad\qquad
\begin{tikzcd}
(\Mod_{A \otimes_{\SS} A})^{\hC} \ar[r,"{(-)^{\tC}}"] \ar[d] & \Mod_A \ar[d] \ar[dl,Rightarrow] \\
(\Mod_{B \otimes_{\SS} B})^{\hC} \ar[r,"{(-)^{\tC}}"] & \Mod_{B}
\end{tikzcd}
\end{equation}
which is obtained from the commuting square on the left
after replacing the vertical forgetful functors by their left adjoints \(B \otimes_A (-)\) and \((B \otimes_{\SS} B) \otimes_{A \otimes_{\SS} A} (-)\), respectively.
\end{corollary}

We now wish to apply the above discussion in order to obtain explicit data which permits us to refine \(p_{\phi}\colon \Modp{A} \to \Modp{B}\) to a hermitian functor with respect to a pair of hermitian structures coming from modules with genuine involution
\((M_A,N_A,\alp) \in \Mod_{\N A}\) and \((M_B,N_B,\beta) \in \Mod_{\N B}\).
In terms of quadratic functors, this data is
a natural transformation
\[
\eta\colon \QF^{\alp}_{M_A} \Rightarrow p_{\phi}^*\QF^{\beta}_{M_B}.
\]
Under the above equivalences, the transformation \(\eta\) corresponds to a map \((M_A,N_A,\alp) \to \phi^*(M_B,N_B,\bet)\) in \(\Mod_{\N A}\), or equivalently by adjunction, to a map \(\phi_!(M_A,N_A,\alp) \to (M_B,N_B,\bet)\). Let us summarise the situation in explicit terms as follows:

\begin{corollary}
\label{corollary:summary-hermitian-induction}%
Keeping the notation above, the data of a hermitian functor \((\Mod_{A}^\omega,\QF^{\alp}_{M_A}) \to(\Mod_{B}^\omega,\QF^{\bet}_{M_B})\) covering the induction functor \(p_{\phi}\) can be encoded by a triple \((\delta,\gam,\sig)\) where \(\delta\colon M_A \to M_B\) and \(\gamma\colon N_A \to N_B\) in \(\Mod_{A \otimes_{\SS} A}^{\hC}\) and \(\Mod_A\) respectively, and \(\sig\) is a commutation homotopy in the square
\[
\begin{tikzcd}
N_A \ar[r,"{\gamma}"] \ar[d,"{\alp}"] & N_B \ar[d,"{\beta}"] \\
M_A^{\tC} \ar[r,"{\delta^{\tC}}"] & M_B^{\tC}
\end{tikzcd}
\]
Equivalently, we may provide the adjoints \(\ovl{\delta}\colon (B \otimes_{\SS} B) \otimes_{A \otimes_{\SS} A}M_A \to M_B\) and \(\ovl{\gamma}\colon B \otimes_A N_A \to N_B\) in \(\Mod_{B \otimes_{\SS} B}^{\hC}\) and \(\Mod_B\) respectively, together with a commutative square of the form
\[
\begin{tikzcd}
B \otimes_A N_A \ar[rr,"{\ovl{\gamma}}"] \ar[d,"{\id \otimes_A \alp}"'] && N_B \ar[d,"{\beta}"] \\
B \otimes_A M_A^{\tC} \ar[r] & \big((B \otimes_{\SS} B) \otimes_{A \otimes_{\SS} A}M_A\big)^{\tC} \ar[r,"{\ovl{\delta}^{\tC}}"] & M_B^{\tC}
\end{tikzcd}
\]
where the left lower horizontal map is the Beck-Chevalley map~\eqrefone{equation:tate-restriction-modules}. It can also be identifies with the composition of the Tate diagonal and the lax monoidal structure of \((-)^{\tC}\).
\end{corollary}

\begin{lemma}
\label{lemma:Poincresscalars}%
In the situation of Corollary~\refone{corollary:summary-hermitian-induction}, the hermitian functor
\[
(p_{\phi},\eta)\colon (\Mod_{A}^\omega,\QF^{\alp}_{M_A}) \to(\Mod_{B}^\omega,\QF^{\bet}_{M_B})
\]
associated to a map \((\delta,\gamma,\sig)\colon (M_A,N_A,\alp) \to \phi^*(M_B,N_B,\beta)\) is Poincaré if and only if the composite
\[
B \otimes_A M_A \to (B \otimes_{\SS} B)\otimes_{A \otimes_{\SS} A}M_A \xrightarrow{\ovl{\del}} M_B
\]
is an equivalence, where the first map is induced by the left unit \(B \to B \otimes_{\SS} B\).
\end{lemma}

\begin{proof}
By definition, the hermitian functor \((p_{\vphi},\eta)\) is Poincaré if and only if the induced map
\[
B \otimes_A \Dual_{M_A}(X) = B \otimes_A\map_A(X,M_A)\longrightarrow \map_{B}(B\otimes_A X,M_B) = \Dual_{M_B}(B \otimes_{A} X)
\]
is an equivalence of \(B\)-modules for every perfect \(A\)-module \(X\). Since \(\Modp{A}\) is generated under finite colimits and retracts by \(A\) this map is an equivalence for every \(X \in \Modp{A}\) if and only if it is an equivalence for \(X = A\). But this is exactly the statement that the induced map
\[
B \otimes_A M_A = B \otimes_A \map_A(A,M_A) \to \map_{B}(B,M_B) = M_B
\]
is an equivalence, as desired.
\end{proof}

\begin{definition}
\label{definition:invmap}%
In the situation of Lemma~\refone{lemma:Poincresscalars}, when the condition that the induced map \(B \otimes_A M_A \to M_B\) is an equivalence holds, we say that the morphism \((\delta,\gamma,\sig)\) is \defi{\(\phi\)-invertible}. In particular, Lemma~\refone{lemma:Poincresscalars} says that the hermitian functor induced by \((\delta,\gamma,\sig)\) is Poincaré if and only if \((\delta,\gamma,\sig)\) is \defi{\(\phi\)-invertible}.
\end{definition}

\begin{example}
\label{example:map-rings-genuine-involution}%
Suppose that \((\phi,\tau)\colon (A,N_A,\alp) \to (B,N_B,\beta)\) is a map of \(\Eone\)-algebras with genuine anti-involution (Example~\refone{example:anti-involution}), so that \(\phi\colon A \to B\) is a map of rings with anti-involution and \(\tau\colon N_A\to \phi^\ast N_B\) is a map of \(A\)-modules. Then both \((A,N_A,\alp)\) and \((B,N_B,\beta)\) are invertible as modules with genuine involution over \(A\) and \(B\) respectively and \(\tau\) is \(\phi\)-invertible. In particular, in this situation we always obtain an induced Poincaré functor
\( (p_{\phi},\tau)\colon (\Modp{A},\QF^{\alp}_A) \to (\Modp{B},\QF^{\beta}_{B}) \).
\end{example}

\begin{example}
\label{example:LANquadsym}%
Suppose that \(\phi \colon A \rightarrow B\) is a map of \(\Eone\)-algebras. A module with involution \(M_A \in (\Mod_{A \otimes_{\SS} A})^{\hC}\) then determines a symmetric bilinear functor with associated quadratic hermitian structure \(\QF^{\qdr}_{M_A}\) on \(A\) encoded by the module with genuine involution \((M_A,0,0\to M_A^{\tC})\). The left Kan extension of \(\QF^{\qdr}_{M_A}\) to \(\Modp{B}\) is then encoded by the module with genuine involution \((M_B, 0,0 \to M_B^{\tC})\) for \(M_B := M_A\otimes_{A\otimes_{\SS} A}(B\otimes_{\SS} B)\), and so
\[
(p_{\phi})_!\QF^{\qdr}_{M_A} \simeq \QF^{\qdr}_{M_B}.
\]
On the other hand, the associated symmetric hermitian structure \(\QF^{\sym}_{M_A}\) is encoded by the module with genuine involution \((M_A,M_A^{\tC},\id\colon M_A^{\tC} \to M_A^{\tC})\), and hence its left Kan extension to \(\Modp{B}\) is encoded by the module with genuine involution \((M_B, B \otimes_A M_A^{\tC}, B \otimes_A M_A^{\tC} \to M_B^{\tC})\),
which is generally not the symmetric hermitian structure \(\QF^{\sym}_{M_B}\),
unless \(B\) is perfect as an \(A\op\)-module (indeed, for a fixed \(M\) the underlying map of spectra \(B \otimes_A M^{\tC} \to (B\otimes_A M)^{\tC}\) can be considered as a natural transformation between two exact functors in the argument \(B \in \Modp{A\op}\) which is an equivalence on \(B=A\op\) and hence on any perfect \(B\)). As a counter-example consider \(\SS \rightarrow \SS[\frac 1 2]\): by Lin's theorem \cite{lin} one has that \(\SS^\tC \simeq \SS^\wedge_2\) so \(\SS[\frac 1 2] \otimes_\SS \SS^\tC \simeq \GEM\QQ_2\) whereas \((\SS[\frac 1 2] \otimes_\SS \SS[\frac 1 2])^\tC  \simeq 0\) since \(2\) is invertible on \(\SS[\frac 1 2] \otimes_\SS \SS[\frac 1 2] \simeq \SS[\frac 1 2]\).
\end{example}

\subsection{Ranicki periodicity}
\label{subsection:herm-shifts}%

In this final section we will discuss the effect of shifting Poincaré structures on module \(\infty\)-categories. In the case of ordinary rings, this is the basis for the classical \(4\)-fold periodicity in \(\L\)-theory. To unravel the essence of this phenomenon, we introduce the notion of an \emph{\(n\sig\)-orientation} on modules with involution, and explain how it combines with shifting to yield an \(2n\)-fold periodicity effect on the level of \(\L\)-groups.
In \papertwo we will extend these results to the level of both \(\L\)-spectra and Grothendieck-Witt spectra, yielding in particular a generalisation of Karoubi's fundamental theorem.
For the remainder of this subsection we fix an arbitrary base \(\Einf\)-ring \(k\).

Recall that for a quadratic functor \(\QF\) on a stable \(\infty\)-category \(\C\), and an integer \(n \in \ZZ\), we denoted by \(\QF\qshift{n}\) the quadratic functor on \(\C\) given by \(\QF\qshift{n}(\x) = \Sig^n\QF(\x)\) (see Definition~\refone{definition:shift}). Since the formation of linear and bilinear parts is exact the shifted quadratic functor \(\QF\qshift{n}\) has bilinear part \(\Sig^n\Bil_{\QF}\), linear part \(\Sig^n\Lin_{\QF}\), and structure map \(\Sig^n\Lin_{\QF}(\x) \to \Sig^n\Bil(\x,\x)^{\tC}\) induced by the structure map of \(\QF\) on \(n\)-fold suspensions. If \(\C\) is now of the form \(\Modp{A}\) for some \(\Eone\)-algebra \(A\) over \(k\) and \(\QF = \QF^{\alp}_M\) for some \(k\)-module with genuine involution \((M,N,\alp)\) over \(A\), then \(\QF\qshift{n} = \Sig^n\QF^{\alp}_{M}\) is the quadratic functor associated to the \(k\)-module with genuine involution \((\Sig^n M,\Sig^n N,\Sig^n\alp)\), and we write
\begin{equation}
\label{equation:shifting-genuine-modules}%
(\QF^{\alp}_M)\qshift{n} \simeq \QF^{\Sig^n \alp}_{\Sig^n M}.
\end{equation}
A second natural operation we can perform on the quadratic functor \(\QF^{\alp}_M\) is to \emph{pre-compose} it with \(\Sig^n\). In this case the hermitian \(\infty\)-category \((\Modp{A},\QF^{\alp}_M\circ \Sig^n)\) is canonically equivalent to \((\Modp{A},\QF^{\alp}_M)\) via the functor \(\Sig^n\colon \Modp{A} \to \Modp{A}\). However, the reparametrised quadratic functor \(\QF^{\alp}_M\circ \Sig^n\) is also directly identified with the quadratic functor associated to another \(k\)-module with genuine involution. To determine it, let us introduce some notation. Given a finite dimensional real representation \(V\) for the group \(\Ct\), let us denote by \(S^V\) the associated one-point compactification of \(V\), which is a sphere of dimension \(\dim(V)\) equipped with a based \(\Ct\)-action. Similarly, we will denote by \(\SS^{V} = \Sig^{\infty}(S^V)\) the associated suspension spectrum with \(\Ct\)-action. Given a spectrum \(X\) with a \(\Ct\)-action, we will denote by \(\Sig^V X := \SS^V \otimes_{\SS} X\) the smash product of \(\SS^{V}\) and \(X\) equipped with its diagonal \(\Ct\)-action. Similarly, the \(\infty\)-category \(\Mod_{A \otimes_k A}^{\hC}\) of \(k\)-modules with involution over \(A\) is naturally tensored over spectra with \(\Ct\)-action,  and given a \(k\)-module with involution \(M\) over \(A\) we will denote by \(\Sig^VM = \SS^V \otimes_{\SS} M\) the associated tensor of \(M\) by \(\SS^V\).
We will denote by \(\sig\) the \(1\)-dimensional sign representation of \(\Ct\) and \(\rho\) the \(2\)-dimensional regular representation of \(\Ct\). In particular, \(\rho\) decomposes of the direct sum of a trivial representation and a sign representation, which we write as \(\rho = 1+\sig\). More generally, we will denote by \(a + b\rho + c\sig\) the direct sum of \(a\) copies of the \(1\)-dimensional trivial representation, \(b\) copies of \(\rho\), and \(c\) copies of \(\sigma\).

We will require the following lemma:

\begin{lemma}
\label{lemma:tate-periodic}%
Let \(X\) be a spectrum with a \(\Ct\)-action. Then the map
\[
X^{\tC} \to (\Sig^{\sig}X)^{\tC}  ,
\]
induced by the \(\Ct\)-equivariant map \(S^0 \to S^{\sig}\), is an equivalence.
In particular, the Tate construction is invariant under tensoring with the sign representation sphere spectrum \(\SS^{\sig}\).
\end{lemma}
\begin{proof}
The cofibre of the map \(S^0 \to S^{\sig}\) is given by the pointed \(C_2\)-space \(\Sigma(C_2)_+\) and so the cofibre of \(X \to \Sig^{\sig} X\) is given by \(\Sigma^{\infty+1}_+(C_2) \otimes_{\SS} X\). The claim now follows since the latter is an induced \(C_2\)-spectrum and thus its Tate construction vanishes.
\end{proof}

\begin{remark}
\label{remark:iterating-tate}%
Instead of tensoring with the map \(S^0 \to S^{\sig}\) we may also \emph{cotensor} with it, yielding a map \(\Om^{\sig}X \to X\), which induces an equivalence on Tate constructions since its fibre is an induced \(C_2\)-spectrum by the same argument as in the proof of Lemma~\refone{lemma:tate-periodic}.
Iterating these constructions we then obtain for every non-negative \(n\) equivalences of the form \(X^{\tC} \xrightarrow{\simeq} (\Sig^{n\sig}X)^{\tC}\) and \((\Om^{n\sig}X)^{\tC} \xrightarrow{\simeq} X^{\tC}\).
\end{remark}

\begin{proposition}
\label{proposition:general-equivalence-of-poincare-infty-categories}%
Let \((M,N,\alp)\) be a \(k\)-module with genuine involution over an \(\Eone\)-algebra \(A\). Then for \(n \in \ZZ\) there are equivalences of quadratic functors
\[
(\QF^{\alp}_M)\qshift{n+m} \circ \Sig^n \simeq \QF^{\Sig^{m}\alp}_{\Sig^{m-n\sig} M}
\]
where the right hand side is the quadratic functor associated to the \(k\)-module with genuine involution
\[
(\Sig^{m-n\sig}M,\Sig^{m}N,\Sig^{m}\alp)
\]
defined using the identification \((\Sig^{m-n\sig}M)^{\tC}  \simeq \Sig^{m} M^{\tC}\) issued from Lemma~\refone{lemma:tate-periodic} and Remark~\refone{remark:iterating-tate}.
In particular, the functor \(\Omega^{n}\colon \Modp{A} \to \Modp{A}\) refines to an equivalence
\[
(\Modp{A},(\QF^{\alp}_M)\qshift{2n})\xrightarrow{\simeq} (\Modp{A},(\QF^{\alp}_M)\qshift{2n}\circ (\Sig^n)\op) \simeq(\Modp{A},\QF^{\Sig^{n}\alp}_{\Sig^{n(1-\sig)} M})
\]
of hermitian \(\infty\)-categories (or Poincaré when \(M\) is invertible). The same statement holds if we replace \(\Modp{A}\) by \(\Mod^c_A\) for some subgroup \(c \subseteq \K_0(\Modp{A})\).
\end{proposition}
\begin{proof}
In light of~\eqrefone{equation:shifting-genuine-modules} it will suffice to prove for the case \(m=-n\). In this case we need to show that \(\QF^{\alp}_M \circ \Sig^n \simeq \QF^{\Sig^n\alp}_{\Sig^{n\rho}M}\). Replacing \((M,N,\alp)\) by \((\Sig^{n\rho}M,\Sig^nN,\Sig^n\alp)\) we observe that the claim for \(n\) implies the same for \(-n\), and so we may assume that \(n\) is positive. Arguing by induction, it is then enough to prove the case \(n=1=-m\).
Inspecting the explicit formula for the linear and bilinear parts one directly sees that they commute with pre-composing with \(\Sig\op\), that is, for a quadratic functor \(\QF\) we have natural equivalences
\[
\Bil_{\QF \circ \Sig\op} \simeq \Bil_{\QF} \circ (\Sig\op,\Sig\op) \quad\text{and}\quad \Lam_{\QF \circ \Sig\op} \simeq \Lam_{\QF} \circ \Sig\op .
\]
It then follows that the canonical interchange map \(\QF \circ \Sig\op \Rightarrow \Om \circ \QF\) induces an equivalences on linear parts, and on bilinear parts gives the corresponding interchange map
\[
\Bil_{\QF}(\Sig X,\Sig Y) \to \Om\Bil_{\QF}(X,Y)
\]
associated to \(\Bil_{\QF}\), where here we view \(\Bil_{\QF}\) as a functor in a single input, the tuple \((X,Y)\). If we however break this input into two separate variables and take into account the symmetric structure of \(\Bil\), then the above interchange map factors as a composite
\[
\Bil_{\QF}(\Sig X,\Sig Y) \xrightarrow{\simeq} \Om^{\rho}\Bil_{\QF}(X,Y) \to \Om\Bil_{\QF}(X,Y) ,
\]
where the second map is induced by cotensoring along the ``smash diagonal'' \(S^1 \to S^{\rho}\), and the first arises from the interchange map of each variable separately, which is an equivalence since
\(\Bil_{\QF}\) is bilinear.
If we now evaluate at \(Y=X\)
then we obtain a composite of \(\Ct\)-equivariant maps
\[
\Bil_{\QF}(\Sig X,\Sig X) \xrightarrow{\simeq} \Om^{\rho}\Bil_{\QF}(X,X) \to  \Om\Bil_{\QF}(X,X) ,
\]
which induces an equivalence on Tate construction by Lemma~\refone{lemma:tate-periodic} and Remark~\refone{remark:iterating-tate} (indeed,
under the representation sphere isomorphism \(S^{\rho} \cong S^{1+\sig}\), the smash diagonal identifies with the suspension of the map \(S^0 \to S^{\sig}\)). We therefore conclude that the interchange map \(\QF \circ \Sig\op \Rightarrow \Om \circ \QF\) identifies the reference map \(\Lam_{\QF \circ \Sig\op} \Rightarrow (\Bil_{\QF \circ \Sig\op}^{\Del})^{\tC}\) of \(\QF \circ \Sig\op\) with the desuspension of the reference map \(\alp\colon \Lam_{\QF} \Rightarrow (\Bil_{\QF} \circ \Del)^{\tC}\) of \(\QF\).
We summarize the above arguments in the following commutative diagram
\[
\begin{tikzcd}
\QF \circ \Sigma\op \ar[d]\ar[r] &
\Lambda_\QF \circ \Sigma\op \ar[d]\ar[r,equal] &
\Lambda_\QF \circ \Sigma\op \ar[d]\ar[r,"\sim"] &
\Om \circ \Lambda_\QF  \ar[d, "\Om\alp"] \\
(\Bil_\QF \circ \Delta  \circ \Sigma\op)^\hC \ar[d,"\sim"]\ar[r] &
(\Bil_\QF  \circ \Delta \circ \Sigma\op)^\tC\ar[d,"\sim"]\ar[r,"\sim"] &
(\Om \circ \Bil_\QF \circ\Delta)^\tC \ar[r,"\sim"]\ar[d,equal] &
\Om \circ (\Bil_\QF \circ\Delta)^\tC \ar[d,equal] \\
(\Om^{\rho} \circ \Bil_\QF \circ \Delta)^\hC \ar[r] &
(\Om^{\rho} \circ \Bil_\QF \circ \Delta)^\tC \ar[r,"\sim"] &
(\Om \circ \Bil_\QF \circ\Delta)^\tC \ar[r,"\sim"] &
\Om \circ (\Bil_\QF \circ\Delta)^\tC
\end{tikzcd}
\]
which may be viewed as an equivalence between the top left square and the external rectangle. As such, it identifies the classifying square of \(\QF \circ \Sig\op\) with one having bilinear term \(\Om^{\rho}\circ \Bil_{\QF}\), linear term \(\Om \circ \Lam_{\QF}\), and reference map \(\Om\alp\). Applying this in the case of \(\QF = \QF^{\alp}_M \circ \Sig\op\) we now conclude that \(\QF^{\alp}_M \circ \Sig\op\) is represented by the module with genuine involution \((\Om^{\rho}M,\Om N,\Om\alp)\) defined using the identification \((\Om^{\rho}M)^{\tC}  \simeq \Om M^{\tC}\), as desired.
\end{proof}

\begin{remark}
\label{remark:twistsign}%
We will show in~\S\refone{subsection:mackey-functors} (see Corollary~\refone{corollary:quadratic-genuine}) that any quadratic functor \(\QF \colon \C\op \to \Spa\) canonically refines to the genuine fixed points of a functor \(\wtl{\QF}\colon \C\op \to \Spagc\) to genuine \(\Ct\)-spectra. Thus one can make sense of tensoring any quadratic functor with \(\SS^{n-m\sigma}\), equipped with its genuine \(\Ct\)-structure in which the geometric fixed points are \(\SS^n\). A version of the above calculations also holds in this generality, see Proposition~\refone{prop:shiftofqingeneral}. For example, one finds that
twisting a quadratic functor by \(\SS^{1-\sig}\) corresponds to the operation \(\QF \mapsto \QF\qshift{2} \circ \Sig\), while twisting by \(\SS^{\sig-1}\) corresponds to the operation \(\QF \mapsto \QF\qshift{-2} \circ \Om\). We shall see in the examples below that these really are two inequivalent operations.
\end{remark}

Our next goal is to use the above calculations to deduce periodicity properties for the \(\L\)-groups of certain rings. For this, we will use the notion of \emph{orientation}, suitably adapted to the present setting.

\begin{construction}
\label{construction:orientation}%
Suppose given an \(\Eone\)-group \(G\) and a \defi{spherical character} of \(G\), that is, a map of \(\Eone\)-groups \(\chi\colon G \to \gl_1(\SS) = \Aut_{\Spa}(\SS)\). We may consequently view \(G\) as acting on the sphere spectrum \(\SS\), and will write \(\SS^{\chi}\) for associated \(G\)-spectrum. Moreover, since every stable \(\infty\)-category \(\C\) is canonically tensored over \(\Spaf\), the character \(\chi\) determines an action of \(G\) on any object \(X \in \C\), by identifying \(X\) with the underlying object of \(X^{\chi} := \SS^{\chi} \otimes X\). More generally, if \(X\) is already equipped with a \(G\)-action, then we may tensor \(X\) with \(\SS^{\chi}\) to obtain a new \(G\)-object \(X^{\chi} = \SS^{\chi} \otimes X\) with the same underlying object \(X\), that which can be called \defi{twisting} \(X\) by \(\chi\).
Even more generally, if \(\C\) itself admits a \(G\)-action, then we can similarly twist a \(G\)-fixed object \(X \in \C^{\h G}\) by \(\chi\) to obtain a new \(G\)-fixed object \(X^{\chi} \in \C^{\h G}\). By a \defi{\(\chi\)-orientation} on such a \(G\)-fixed object \(X\) we will then mean an equivalence
\[
X^{\chi} \xrightarrow{\simeq} X
\]
of \(G\)-fixed objects, refining the identity \(X=X\) in \(\C\). When the \(G\)-action on both \(\C\) and \(X\) is trivial and \(\C\) admits quotients for \(G\)-objects, the data of a \(\chi\)-orientation can equivalently by encoded via a factorisation
\[
X \to (X^{\chi})_{\h G} \to X
\]
of the identity \(X = X\), where the first map is the canonical map from the underlying object to the quotient.
\end{construction}

Given an \(\Eone\)-group \(G\) and a spherical character \(\chi\colon G \to \gl_1(\SS)\), we will refer to the homotopy quotient \(\SS^{\chi}_{\h G}\) as the \defi{Thom spectrum} of \(\chi\), and will denote it also by \(\Th(\chi)\). Recall that \(\B\gl_1(\SS)\) identifies with the classifying space for stable spherical fibrations. In particular, it carries the universal stable spherical fibration $\gamma$ and the definition above coincides with the usual Thom spectrum for $(\B\chi)^*(\gamma)$. In particular, when \(\chi\)
factors over the $J$-homomorphism
\[
G \longrightarrow O \xrightarrow{\mathrm{J}} \gl_1(\SS),
\]
the spectrum \(\Th(\chi)\) is indeed the Thom spectrum of the associated stable vector bundle.
If \(\C\) is now a cocomplete, stable \(\infty\)-category,
so that the action of \(\Spaf\) extends to all of \(\Spa\) in a colimit preserving manner, then we may use the Thom spectrum to describe orientations, at least for objects with trivial \(G\)-action. In particular, a \(\chi\)-orientation on such an \(X \in \C\) is equivalent to the data of a factorisation of the identity map \(X=X\) as
\[ X \to \Th(\chi) \otimes X \to X ,\]
where the first map is the tensor with \(X\) of the map to the quotient \(\SS  \to \SS^{\chi}_{\h G} = \Th(\chi)\) induced by the initial \(\Eone\)-group map \(\ast \to G\).
In that case, an orientation also determines an equivalence \(X^{\chi}_{\h G} = \Th(\chi) \otimes X \simeq \B G \otimes X = X_{\h G}\), a generalisation of the classical Thom isomorphism, where the last term refers to the homotopy quotient of \(X\) by the trivial \(G\)-action.

\begin{notation}
\label{notation:oriented-rings}%
Let \(G\) be an \(\Eone\)-group equipped with a spherical character \(\chi\).
By a \defi{\(\chi\)-oriented \(\Einf\)-ring spectrum} we will mean an \(\Einf\)-ring spectrum \(E\), which, when considered as a module over itself acted trivially upon by \(G\), is equipped with a \(\chi\)-orientation. In particular, identifying \(\map_E(E,E) = \map_{\Spa}(\SS,E)\) such a \(\chi\)-orientation is equivalent to a refinement of the unit map \(\SS \to E\) to a \(G\)-equivariant map \(\SS \to E^{\chi}\), where \(\SS\) is considered as having trivial \(G\)-action. As above, may also encode this data as a factorisation of the unit map as a composite
\(\SS \to \Th(\chi) \to E\).
\end{notation}

\begin{remark}
\label{remark:inherited-orientation}%
Let \(G\) be an \(\Eone\)-group equipped with a spherical character \(\chi\) and \(\C\) a stable presentable \(\infty\)-category (so that \(\Spa\) acts on \(\C\)).
If \(E\) is a \(\chi\)-oriented \(\Einf\)-ring spectrum (in the sense of Notation~\refone{notation:oriented-rings})
then any \(E\)-module object in \(\C^{\h G}\) is automatically \(\chi\)-oriented. Indeed, for such an \(X\) the \(\chi\)-orientation of \(E\) induces a \(\chi\)-orientation
\( X^{\chi} = \SS^{\chi} \otimes X = (\SS^{\chi} \otimes_{\SS} E) \otimes_E X \xrightarrow{\simeq} E \otimes_E X = X \).
\end{remark}

\begin{example}
\label{example:complex-oriented}%
Consider the \(\Eone\)-group \(\Un(1)\) together with its standard spherical character \(\chi\colon \Un(1) \to \gl_1(\SS)\) coming from the tautological complex representation of \(\Un(1)\). As the Thom space of the associated line bundle on \(\B \Un(1) = \CC P^{\infty}\) is equivalent via the zero section to \(\CC P^{\infty}\) one has
\[
\Th(\chi) = \MU(1) = \Sig^{\infty-2}\CC P^{\infty}.
\]
For an \(\Einf\)-ring spectrum \(E\), the data of a \(\chi\)-orientation on \(E\) (as an \(E\)-module with constant action, see Notation~\refone{notation:oriented-rings}) is given by a factorisation in the \(\infty\)-category of spectra
\[
\SS^2 \to  \Sig^{\infty}\CC P^{\infty} \to \Sig^2E
\]
of the 2-fold suspended unit \(\SS \to E\), a structure otherwise known as a \defi{complex orientation} on \(E\).
\end{example}

\begin{remark}
\label{remark:E-inf-orientation}%
Suppose that \(G\) is not just an \(\Eone\)-group but an \(\Einf\)-group, and that \(\chi\colon G \to \gl_1(\SS)\) is an \(\Einf\)-map. Then the corresponding Thom spectrum
\( \Th(\chi) = \SS^{\chi}_{\h G} \)
inherits the structure of an \(\Einf\)-ring with unit the canonical map \(\SS \to \SS^{\chi}_{\h G} = \Th(\chi)\) to the quotient. This unit then tautologically factors through \(\Th(\chi)\), and hence
\(\Th(\chi)\) is canonically \(\chi\)-oriented. This means, in particular, that any \(\Einf\)-algebra over \(\Th(\chi)\) is canonically \(\chi\)-oriented. We note however that the structure of a \(\Th(\chi)\)-algebra (that is, the structure of an \(\Einf\)-map \(\Th(\chi) \to E\)) is generally finer then that of a \(\chi\)-orientation, which in turn corresponds to a map \(\Th(\chi) \to E\) merely on the level of plain spectra under \(\SS\). This is a reflection of the fact that \(G\) has a much richer structure than just an \(\Eone\)-group, namely, an \(\Einf\)-group structure. One may take this additional structure into account by defining a suitable notion of an \(\Einf\)-\(\chi\)-orientation on an \(\Einf\)-ring spectra \(E\), in which case one would obtain a notion equivalent to that of an \(\Einf\)-algebra over \(\Th(\chi)\).
\end{remark}

\begin{example}
\label{example:thom-spectra}%
For concrete examples of the previous remark, consider the diagram of \(\Einf\)-groups
\[
\begin{tikzcd}
& \Sp\ar[d]\ar[r] & \Un \ar[d]&&&  \\
\String\ar[r] & \Spin \ar[r] & \Spin^{\mathrm{c}} \ar[r] & \SO \ar[r] & \rO  \ar[r] & \gl_1(\SS)
\end{tikzcd}
\]
yielding by restriction spherical characters on all the appearing \(\Einf\)-groups. We then have that the corresponding Thom spectra \(\MSp,\MU,\MString,\MSpin,\MSpin^{\mathrm{c}},\MSO\) and \(\MO\) are each oriented with respect to its respective spherical character. Similarly, any \(\Einf\)-algebra over any of these Thom spectra will inherit the corresponding orientation.
\end{example}

\begin{remark}
\label{remark:very-confusing}%
The notions of a complex (or \(\U(1)\)-) orientation appearing in Example~\refone{example:complex-oriented} and the notion of a \(\U\)-orientation appearing in Example~\refone{example:thom-spectra} are closely related. On the one hand, the \(\Eone\)-group map \(\U(1) \to \U\) allows one to restrict any \(\U\)-orientation to a \(\U(1)\)-orientation. In concrete terms, this means that any map \(\MU \to E\) on the level of plain spectra under \(\SS\) restricts to a map \(\MU(1) = \Sig^{\infty-2}\CC P^{\infty} \to E\) of the same nature. On the other hand, it is a well-known consequence of the splitting principle that any complex orientation \(\MU(1) \to E\) extends to a map \(\MU \to E\), which is furthermore compatible with the ring spectrum structure on both sides up to (non-coherent) homotopy, see~\cite[Part II]{adams-stable}, or~\cite[Lecture 6]{Lurie-chromatic}.
In particular, the map from the space of \(\U\)-orientations on \(E\) to the space of complex orientations on \(E\) is surjective on components, though is generally not an equivalence. Both these structures are weaker than that of being an \(\Einf\)-algebra over \(\MU\), which would correspond, as in Remark~\refone{remark:E-inf-orientation}, to a suitable \(\Einf\)-enhancement of the notion of a \(\U\)-orientation, see~\cite{commutative-orientation}.
\end{remark}

We now return to our context of modules with involutions and hermitian structures.

\begin{definition}
For \(n \geq 0\) we will use the term \defi{\(n\sig\)-orientation} to refer to an orientation with respect to the spherical character \(\chi_{n\sig}\colon \Ct \to \gl_1(\SS)\) determined by \(\SS^{n\sig-n}\).
\end{definition}

Fixing an \(\Eone\)-algebra over \(k\), we will be interested in the notion of \(n\sig\)-orientation for \(k\)-modules with involution \(M\) over an \(\Eone\)-algebra \(A\) over \(k\), viewed as \(\Ct\)-fixed objects with respect to the \(\Ct\)-action on \(\Mod_{A \otimes_k A}\). We note that for \(A=k\), the notion of an \(n\sig\)-orientation on \(M=k\) with trivial \(\Ct\)-action reduces to that of an \(n\sig\)-orientation on the \(\Einf\)-ring \(k\) in the sense of Notation~\refone{notation:oriented-rings}.

\begin{remarks}
\label{remarks:basic}%
\
\begin{enumerate}
\item
Every object in \(\Mod_{A \otimes_k A}^{\hC}\) is canonically a \(k\)-module object, and so by Remark~\refone{remark:inherited-orientation} we have that an \(n\sig\)-orientation on \(k\) itself determines such an orientation on any \(k\)-module with involution over any \(\Eone\)-algebra \(A\) over \(k\).
\item
An \(n\sig\)-orientation determines an \(N\sig\)-orientation for any multiple \(N\) of \(n\).
\item
\label{item:tate-periodic}%
By Lemma~\refone{lemma:tate-periodic} the Tate construction swallows sign suspensions, and so if \(M\) is a \signper{n} then \(M^{\tC}\) is \(n\)-periodic, that is \(\Sig^n M^{\tC} \simeq M^{\tC}\).
\end{enumerate}
\end{remarks}

Of particular importance are the cases \(n=1,2\). In particular, \(2\sig\)-orientation determines a refinement of the identity map \(M = M\) to a \(\Ct\)-equivariant map
\[
\SS^{1-\sig} \otimes_{\SS} M \simeq \SS^{\sig-1} \otimes_{\SS} M,
\]
and we shall denote this common value by \(-M\). If the \(2\sig\)-orientation comes from a \(\sig\)-orientation, then \(-M \simeq M\).

\begin{examples}\
\label{example:signperiodicmodules}%
\begin{enumerate}
\item
\label{item:example-integers}%
If \(k\) is an \(\Einf\)-algebra over \(\GEM \ZZ\) (e.g., \(k=\GEM R\) for some discrete ring \(R\))
then \(k\) admits a \(2\sig\)-orientation. Indeed,
it suffices to check this for \(k=\GEM\ZZ\). Since \(\Map_{\GEM\ZZ}(\GEM\ZZ,\GEM\ZZ) = \ZZ\) is discrete and the action of the generator of \(\Ct\) on \(\SS^{2\sig-2}\) is homotopic to the identity, it follows that the induced action on \(\SS^{2-2\sig} \otimes_{\SS} \GEM\ZZ\) is trivial, and consequently \(1 \in \ZZ\) admits a unique lift to a fixed point. For such a \(k\) we then have by Remark~\refone{remark:inherited-orientation} that any \(k\)-module with involution \(M\) over any \(k\)-\(\Eone\)-algebra \(A\) is \signper{2}. More generally, even if \(k\) itself is not an \(\Einf\)-algebra over \(\GEM \ZZ\), it is still true that any \(\GEM \ZZ\)-module object in \(\Mod_{A \otimes_k A}^{\hC}\) is canonically \(2\sig\)-oriented, see Remark~\refone{remark:inherited-orientation}.
When \(M\) is furthermore discrete (that is, in the image of the Eilenberg-Mac Lane functor)
the meaning of \(-M\)
reduces to the naive one, that simply inserts a sign into the involution.
\item
\label{item:F-2}%
If \(k\) is an \(\Einf\)-algebra over \(\GEM\FF_2\)
then the same argument applied to \(\SS^{\sigma}\) provides a \(\sig\)-orientation, since the action of the generator of \(\Ct\) on \(\SS^\sig\) is the negative of the identity. Any \(k\)-module with involution over any \(k\)-\(\Eone\)-algebra \(A\) is then \(\sig\)-oriented.
More generally, even if \(k\) itself is not an \(\Einf\)-algebra over \(\GEM \FF_2\), it still true that any \(\GEM\FF_2\)-module object in \(\Mod_{A \otimes_k A}^{\hC}\) is canonically \(\sig\)-oriented, see Remark~\refone{remark:inherited-orientation}.
\item
As a partial converse to the above we note that if the \(\Einf\)-ring \(k\) is \(\sig\)-oriented then at least the underlying \(\mathrm{E}_2\)-algebra of \(k\) refines to a \(\mathrm{E}_2\)-\(\GEM\FF_2\)-algebra, and hence any module with involution over any \(\Eone\)-\(k\)-algebra \(A\) is canonically an \(\GEM\FF_2\)-module object in \(\Mod_{A \otimes_k A}^{\hC}\). This is because any \(\sig\)-orientation on \(k\) determines a homotopy \(2 \sim 0\) on it, and by Mahowald's theorem \(\GEM\FF_2\) is the initial \(\mathrm{E}_2\)-ring spectrum with a homotopy \(2 \sim 0\) (see~\cite[Theorem 5.1]{thom-spectra} for a recent treatment, where we recall that for \(p=2\) no completion is necessary).
\item
The module with involution \(\GEM \FF_2 \otimes_{\SS} \GEM \FF_2\) over \(\GEM \FF_2\) (equipped with the flip action) does \emph{not} admit an \(n\sig\)-orientation for any \(n\), even though the underlying spectrum of \(M\) is an \(\GEM\FF_2\)-module. Indeed, by Lin's theorem \cite{lin} (see \cite[Theorem III.1.7]{NS} for a treatment in the present language), the Tate diagonal provides an equivalence
\[
\GEM \FF_2 \longrightarrow (\GEM \FF_2 \otimes_{\SS} \GEM \FF_2)^\tC.
\]
In particular, the latter spectrum is not periodic, and so \(\GEM \FF_2 \otimes_{\SS} \GEM \FF_2\) cannot be \signper{n}, see Remark~\refone{remarks:basic}\refoneitem{item:tate-periodic}. We note that this does not contradict Remark~\refone{remark:inherited-orientation}, since the \(\GEM\FF_2\)-action on the underlying spectrum of \(M\) is not compatible with its \(\Ct\)-action.
\item
\label{item:warning-many-signs}%
In a similar vein, the sphere does not admit an \(n\sig\)-orientation for any \(n\), since we would obtain an equivalence
\[
\SS^n \otimes_{\SS} (\SS^{\tC}) \simeq (\SS^n)^{\tC} \simeq (\SS^{n\sigma})^\tC \simeq \SS^{\tC},
\]
but from Segal's conjecture (a consequence of Lin's theorem) \(\SS^\tC\) is the \(2\)-completion of \(\SS\), which is not periodic.
\item
The \(\Einf\)-ring \(\SS[1/2]\) admits a unique \(2\sig\)-orientation. To see this, one can
compute \(\pi_*(\SS^{n\sigma-n}[1/2])^\hC\) using the (collapsing) homotopy fixed point spectral sequence, to see
\[
\pi_*(\SS^{n\sigma-n}[1/2])^\hC = \pi_*(\SS[1/2])^\Ct.
\]
But if \(n\) is even, the induced \(\Ct\)-action on \(\pi_*(\SS[1/2])\) is trivial, so we in particular get that the forgetful map
\[
(\SS^{2\sigma-2}[1/2])^\hC \longrightarrow \SS[1/2]
\]
is an equivalence, so \(\SS[1/2]\) admits a unique \(2\sig\)-orientation. It follows that if \(k\) is an \(\SS[1/2]\)-algebra (that is, 2 acts invertibly on \(k\)) then any \(k\)-module with involution over any \(k\)-\(\Eone\)-algebra \(A\) is \signper{2}.
More generally, it is enough to assume that \(2\) acts invertibly on \(M\) itself, rather than on \(k\). Indeed, in this case \(M\) refines to an \(\SS[1/2]\)-module object in \(\Mod_{A \otimes_k A}^{\hC}\), and is hence \signper{2} by Remark~\refone{remark:inherited-orientation}.
\item
The representation \(\sig\) being real we find that the spherical character \(\chi_{\sig}\) factors through the \(\Einf\)-group map \(\chi_{\rO}\colon \rO \to \gl_1(\SS)\), and so any \(\chi_{\rO}\)-orientation determines a \(\sig\)-orientation. By Example~\refone{example:thom-spectra} we consequently deduce that \(\MO\), as well as any \(\Einf\)-algebra over \(\MO\), is \(\sig\)-oriented.
\item\label{item:complex-oriented}
Similarly, the representation \(2\sig\) refines to a complex representation and the spherical character \(\chi_{2\sig}\) factors through \(\U\). By Example~\refone{example:thom-spectra} we consequently deduce that \(\MU\), as well as any \(\Einf\)-algebra over \(\MU\), is \(2\sig\)-oriented. Since the map \(\Un \to \rO\) factors through \(\Spin^{\mathrm{c}}\), this applies in particular to any \(\Einf\)-algebra over \(\MSpin^{\mathrm{c}}\), such as \(\ku\) or \(\KU\).
More generally, since \(\chi_{\sig}\) actually factors through \(\U(1)\), any complex orientation on \(k\) determines a \(2\sig\)-orientation on it, see Example~\refone{example:complex-oriented} and Remark~\refone{remark:very-confusing}.
In particular, if \(k\) is a complex oriented \(\Einf\)-ring spectrum (e.g., \(k=\MU, \ku,\KU\) or any of the Lubin-Tate spectra)
then any \(k\)-module with involution \(M\) over any \(k\)-\(\Eone\)-algebra \(A\) is \signper{2}.
\item
Further up, the representation \(4\sig\) is quaternion and the spherical character \(\chi_{4\sig}\) factors through \(\Sp\), and so any
\(\Einf\)-algebra over \(\MSp\), is \(4\sig\)-oriented.
Since the map \(\Sp \to \rO\) factors through \(\Spin\), this applies in particular to any \(\Einf\)-algebra over \(\MSpin\), such as \(\ko\) or \(\KO\).
Finally, the representation \(8\sig\) refines to a string representation since its Steifel-Whitney classes vanish in degrees \(< 8\) and \(H^4(\Ct,\ZZ) = H^4(\Ct,\ZZ/2)\), and so any \(\Einf\)-algebra over \(\MString\) (such as \(\mathrm{tmf}\)) is \(8\sig\)-oriented.
\end{enumerate}
\end{examples}

We can now specialise Proposition~\refone{proposition:general-equivalence-of-poincare-infty-categories} to the case of oriented modules with involutions.

\begin{corollary}
\label{corollary:periodicity-poincare}%
Let \(A\) be an \(\Eone\)-algebra and let \((M,N,\alp\colon N \to M^{\tC})\) be an invertible \(k\)-module with genuine involution over \(A\), such that \(M\) is \signper{n}, and \(c \subseteq K_0(A)\) a subgroup closed under the induced involution. Then the \(n\)-fold loop functor \(\Omega^n\) refines to an equivalence of Poincaré \(\infty\)-categories
\[
\big(\Mod^c_A, (\QF_{M}^{\alp})\qshift{2n}\big) \xrightarrow{\simeq} \big(\Mod^c_{A}, \QF^{\Sig^n\alp}_{M}\big).
\]
If \(n=2\), e.g., if \(k\) is discrete or complex oriented, or if \(2\) acts invertibly on \(M\), then the loop functor \(\Om\) furthermore refines to an equivalence
\[
\big(\Mod^c_{A}, (\QF_{M}^{\alp})\qshift{2}\big) \xrightarrow{\simeq} \big(\Mod^c_{A}, \QF^{\Sig\alp}_{-M}\big).
\]
\end{corollary}
\begin{proof}
Apply Proposition~\refone{proposition:general-equivalence-of-poincare-infty-categories} with \(m=n\) and use that \(\Sig^{n\sig-n}M \simeq M\) when \(M\) is \signper{n}. The second claim is essentially by definition (but is useful to note nonetheless, since \(-M := \Sig^{\sig-1}M\) is explicitly given by a sign twist when \(M\) is discrete).
\end{proof}

\begin{corollary}
Let \(A\) be an \(\Eone\)-algebra over \(k\) and let \(M\) be an invertible \(k\)-module with involution over \(A\), such that \(M\) is \signper{n}, and \(c \subseteq K_0(A)\) a subgroup closed under the induced involution. Then the \(n\)-fold loop functor \(\Om^n\) refines to give equivalences of Poincaré \(\infty\)-categories
\[
\big(\Mod^c_{A}, (\QF_{M}^{\sym})\qshift{2n}\big) \xrightarrow{\simeq} \big(\Mod^c_{A}, \QF^{\sym}_{M}\big) \quad \text{and} \quad \big(\Mod^c_{A}, (\QF^{\qdr}_M)\qshift{2n}\big) \xrightarrow{\simeq} \big(\Mod^c_{A}, \QF^{\qdr}_{M}\big).
\]
For \(n=2\), e.g., if \(k\) is discrete or complex oriented, or if \(2\) acts invertibly on \(M\), then the loop functor \(\Om\) refines to give equivalences of Poincaré \(\infty\)-categories
\[
\big(\Mod^c_{A}, (\QF_{M}^{\sym})\qshift{2}\big) \xrightarrow{\simeq} \big(\Mod^c_{A}, \QF^{\sym}_{-M}\big) \quad \text{and} \quad \big(\Mod^c_{A}, (\QF^{\qdr}_M)\qshift{2}\big) \xrightarrow{\simeq} \big(\Mod^c_{A}, \QF^{\qdr}_{-M}\big).
\]
\end{corollary}

\begin{corollary}
\label{corollary:periodicity-truncated}%
Let \(A\) be an \(\Eone\)-algebra over \(k\) and let \(M\) be an invertible \(k\)-module with involution over \(A\), such that \(M\) is \signper{n} and \(c \subseteq K_0(A)\) a subgroup closed under the induced involution. Then the \(n\)-fold loop functor \(\Om^n\) refines to give equivalences of Poincaré \(\infty\)-categories
\[
\big(\Mod^c_{A}, (\QF_{M}^{\geq m})\qshift{2n}\big) \xrightarrow{\simeq} \big(\Mod^c_{A}, \QF^{\geq m+n}_{M}\big)
\]
For \(n=2\), e.g., if \(k\) is discrete or complex oriented, or if \(2\) acts invertibly on \(M\), then the loop functor \(\Om\) refines to give equivalences of Poincaré \(\infty\)-categories
\[
\big(\Mod^c_{A}, (\QF_{M}^{\geq m})\qshift{2}\big) \xrightarrow{\simeq} \big(\Mod^c_{A}, \QF^{\geq m+1}_{-M}\big).
\]
\end{corollary}
\begin{proof}
This is particular case of Corollary~\refone{corollary:periodicity-poincare} since the operation \(\alp \mapsto \Sig^n\alp\) sends the \(m\)-connective cover map \(\tau_{\geq m}M^{\tC} \to M^{\tC}\) to the \((m+n)\)-connective cover map \(\tau_{\geq m+n}(\Sig^nM^{\tC}) \to \Sig^nM^{\tC}\).
\end{proof}

The following periodicity result for \(\L\)-groups immediately follows (see also Remark~\refone{remark:relation-to-classical} for the relation with classical \(\L\)-group periodicity):

\begin{corollary}[Ranicki periodicity]
\label{corollary:periodicity-L-groups}%
Let \(A\) be an \(\Eone\)-algebra over \(k\) and let \(M\) be a \signper{n} invertible \(k\)-module with involution over \(A\),
and \(c \subseteq K_0(A)\) a subgroup closed under the induced involution.
Then for \(d \in \ZZ\) we have canonical isomorphisms
\[
\L_{d-2n}\big(\Mod^c_{A}, \QF_{M}^{\sym}\big) \cong \L_d\big(\Mod^c_{A}, \QF^{\sym}_{M}\big) \quad \text{and} \quad \L_{d-2n}\big(\Mod^c_{A}, \QF^{\qdr}_M\big) \cong \L_d\big(\Mod^c_{A}, \QF^{\qdr}_{M}\big),
\]
and when \(A\) is connective, for \(d,m \in \ZZ\) we have canonical isomorphisms
\[
\L_{d-2n}\big(\Mod^c_{A}, \QF_{M}^{\geq m}\big) \cong \L_d\big(\Mod^c_{A}, \QF^{\geq m+2}_{M}\big).
\]
If \(n=2\) (e.g, \(k\) is discrete or complex oriented, or if \(2\) acts invertibly on \(M\)) then we also have canonical isomorphisms
\[
\L_{d-2}\big(\Mod^c_{A}, \QF_{M}^{\sym}\big) \cong \L_d\big(\Mod^c_{A}, \QF^{\sym}_{-M}\big) \quad \text{and} \quad \L_{d-2}\big(\Mod^c_{A}, \QF^{\qdr}_M\big) \cong \L_d\big(\Mod^c_{A}, \QF^{\qdr}_{-M}\big),
\]
and when \(A\) is connective, for \(d,m \in \ZZ\) we have canonical isomorphisms
\[
\L_{d-2}\big(\Mod^c_{A}, \QF_{M}^{\geq m}\big) \cong \L_d\big(\Mod^c_{A}, \QF^{\geq m+1}_{-M}\big).
\]
\end{corollary}

In summary, the quadratic and symmetric \(\L\)-groups of discrete (or complex oriented) rings are \(4\)-periodic, and those of \(\GEM\FF_2\)-algebras are \(2\)-periodic. In other cases, such as \(\MSp\)- or \(\MSpin\)-algebras, one has \(8\)-periodic \(\L\)-groups, while \(\MString\)-algebras have \(16\)-periodic \(\L\)-groups. For example,  \(\L_{\ast}^\sym(\ku), \L_{\ast}^\sym(\ko)\) and \(\L_{\ast}^\sym(\mathrm{tmf})\) are 4, 8 and 16 periodic, respectively, and similarly for the non-connective variants. We will upgrade these periodicity results to the spectrum level in the second instalment, and show that it neatly combines with a spectrum-level Karoubi periodicity yielding what we call \emph{Karoubi-Ranicki periodicity}, see \S\reftwo{subsection:real-alg-K-Karoubi}.

Regarding the quadratic $\L$-groups let us also mention that for a connective ring spectrum $R$ , the algebraic $\pi$-$\pi$-theorem of Weiss and Williams from \cite{WWII} implies that generally $\L_*^\qdr(R) \cong \L_*^\qdr(\pi_0(R))$, so that the quadratic $\L$-groups in fact end up $4$-periodic; we give an account of this result in Corollary \refthree{corollary:pi-pi}.

%% file: Examples.tex
In this section we will discuss several examples of interest of Poincaré \(\infty\)-categories in further detail.
We will begin in \S\refone{subsection:universal} with the important example of the universal Poincaré \(\infty\)-category \((\Spaf,\QF^{\uni})\), which is characterized by the property that Poincaré functors out of it pick out Poincaré objects in the codomain.
In \S\refone{subsection:discrete-rings} we will consider perfect derived \(\infty\)-categories of \emph{ordinary rings} and show how to translate the classical language of forms on projective modules into that of the present paper via the process of \emph{deriving} quadratic functors. These examples form the main link between the present work and classical hermitian \(\K\)-theory, and will be the main focus of applications in \paperthree. In~\S\refone{subsection:visible} we will turn to a specific family of ring spectra with anti-involution, the \defi{group ring spectra}, which carries a special interest due to its relation with surgery theory. In particular, we will explain how to construct modules with genuine involution over such ring spectra whose associated \(\L\)-groups capture various flavours of \defi{visible \(\L\)-groups}, for example, the visible symmetric \(\L\)-groups defined by~\cite{Weiss} (see Example~\refone{example:visible-symmetric}).
Finally, in~\S\refone{subsection:parametrised-spectra} we will discuss the closely related case of \defi{parametrised spectra}, which serves as a base-point-free version of group rings, and show how to construct Poincaré structures producing the parametrised spectra variant of visible \(\L\)-theory as studied in
~\cite{WWIII} (see Corollary~\refone{corollary:compare-visible-WW}).

\subsection{The universal Poincaré category}
\label{subsection:universal}%

In this section we will discuss the Poincaré \(\infty\)-category \((\Spaf,\QF^{\uni})\) of Example~\refone{example:universal-category}, which we call the \emph{universal Poincaré \(\infty\)-category}. This term is motivated by the following mapping property which we will prove below: the Poincaré \(\infty\)-category \((\Spaf,\QF^{\uni})\) corepresents the functor \(\Poinc\colon\Catp \to \Sps\) which assigns to a Poincaré \(\infty\)-category its space of Poincaré forms.
To exhibit this, consider the map
\begin{equation}
\label{equation:universal-form}%
\SS \to \QF^{\uni}(\SS)
\end{equation}
given by the identity \(\SS \to \Dual\SS = \SS\) on the linear part and by the unit map \(\SS \to (\Dual\SS \otimes \Dual\SS)^{\hC} = \SS^{\hC}\) on the bilinear part. These two maps canonically lead to the `same' map to \((\Dual\SS \otimes \Dual\SS)^{\tC} = \SS^{\tC}\) since the Tate diagonal of \(\SS\) agrees with the composite \(\SS \to \SS^{\hC} \to \SS^{\tC}\) by construction. In particular, the composite \(\SS \to \QF^{\uni}(\SS) \to \Bil_{\QF^{\uni}}(\SS,\SS) = \Dual(\SS \otimes \SS) \simeq \SS\) is the identity. We then note that the map~\eqrefone{equation:universal-form} corresponds to a hermitian form \(\qone^{\uni} \in \Om^{\infty}\QF^{\uni}(\SS)\) such that \(\qone^{\uni}_{\sharp}\colon \SS \to \Dual\SS=\SS\) is the identity, and in particular \(\qone^{\uni}\) is Poincaré. We will refer to it as the
\defi{universal Poincaré form}.

\begin{lemma}
\label{lemma:maps-from-Qg}%
For every quadratic functor \(\QF\colon (\Spaf)\op \to \Spa\), the map
\begin{align}
\label{align:map-equivalence}%
\map(\QF^{\uni}, \QF) \to \QF(\SS),
\end{align}
induced by the universal form~\eqrefone{equation:universal-form}, is an equivalence of spectra.
\end{lemma}
\begin{proof}
Since the collection of \(\QF \in \Funq(\Spaf)\) for which~\eqrefone{align:map-equivalence} is an equivalence is closed under limits it will suffice by Proposition~\refone{proposition:classification}
to prove the claim whenever \(\QF\) is either exact or of the form \(\QF(\x) = \Bil(\x,\x)\) for some symmetric bilinear functor \(\Bil\). In the former case, the claim follows since the linear part of \(\QF^{\uni}\) is \(\Dual\) by construction and the composite \(\SS \to \QF^{\uni}(\SS) \to \Dual\SS\) exhibits \(\Dual\) as (stably) represented by \(\SS\). On the other hand, if \(\QF\) is of the form \(\Bil(\x,\x)\) for some symmetric bilinear functor \(\Bil\) then the result follows from Lemma~\refone{lemma:universal-crs} since the image of the universal form~\eqrefone{equation:universal-form} in
\[
\Bil_{\QF^{\uni}}(\SS,\SS) = \Dual(\SS \otimes \SS) = \map(\SS \otimes \SS,\SS)=\map(\SS,\SS)
\]
corresponds to the identity \(\SS\to\SS\) by construction, and thus exhibits the underlying bilinear part
\[
\Bil_{\QF^{\uni}} \in \Funb(\Spaf) \simeq \Funx((\Spaf \otimes \Spaf)\op,\Spaf) \simeq \Funx((\Spaf)\op,\Spa)
\]
as stably represented by \(\SS\) as well.
\end{proof}

\begin{lemma}
\label{lemma:finite-spectra-universal}%
The sphere spectrum \(\SS \in \Spaf\) exhibits \(\Spaf\) as corepresenting the core groupoid functor
\[
\core\colon\Catx \to \Sps \quad\quad \C \mapsto \grpcr\C .
\]
\end{lemma}
\begin{proof}
Let \(\Sfinast\) be the \(\infty\)-category of finite pointed spaces. Then the inclusion \(\Sfinast \to \Spaf\) exhibits \(\Spaf\) as the Spanier-Whitehead stabilisation of \(\Sfinast\), and in particular, for every stable \(\infty\)-category \(\D\) the restriction functor
\[
\Funx(\Spaf,\D) \to \Funrex(\Sfinast,\D)
\]
is an equivalence, where \(\Funrex(-,-) \subseteq \Fun(-,-)\) denotes the full subcategory spanned by the right exact (i.e., finite colimits preserving) functors, see, e.g.,~\cite[Proposition C.1.1.7]{SAG}. On the other hand, for \(\D\) stable we may identify right exact functors \(\Sfinast \to \D\) with those which are reduced and excisive, and hence with spectrum objects in \(\D\). We thus obtain that for \(\D\) stable the evaluation functor
\[
\Funx(\Spaf,\D) \simeq \Funrex(\Sfinast,\D) \simeq \Spa(\D) \to \D
\]
is an equivalence, and hence induces an equivalence
\[
\Map(\Spaf,\D) \simeq \core\D
\]
on the level of core groupoids.
\end{proof}

We now come to the main result of this section:
\begin{proposition}[Universality of the universal Poincaré \(\infty\)-category]
\label{proposition:corepresentability-of-poinc}%
The universal Poincaré object \((\SS,\qone^{\uni})\) exhibits \((\Spaf,\QF^{\uni})\) as corepresenting the functor \(\Poinc\).
Similarly, when considered as a hermitian \(\infty\)-category, the underlying universal hermitian object exhibits \((\Spaf,\QF^{\uni})\) as corepresenting the functor \(\spsforms\).
\end{proposition}
\begin{proof}
We begin with the second claim. We need to show that for every \((\C,\QF) \in \Cath\) the map
\[
\Map_{\Cath}((\Spaf,\QF^{\uni}),(\C,\QF)) \to \spsforms(\C,\QF)
\]
sending \((f,\eta)\colon (\Spaf,\QF^{\uni}) \to (\C,\QF)\) to \((f(\SS),\eta_{\SS}(\qone^{\uni}))\) is an equivalence.
Consider the commutative square
\begin{equation}
\label{equation:square-corepresentation}%
\begin{tikzcd}
\Map_{\Cath}((\Spaf,\QF^{\uni}),(\C,\QF)) \ar[r]\ar[d] & \spsforms(\C,\QF) \ar[d] \\
\Map_{\Catx}(\Spaf,\C) \ar[r,"{\simeq}"] & \core(\C)
\end{tikzcd}
\end{equation}
furnished by Lemma~\refone{lemma:forms-functorial}, where the bottom horizontal map is the one induced by the object \(\SS \in \iota\Spaf\), which is an equivalence by Lemma~\refone{lemma:finite-spectra-universal}.
It will hence suffice to show that for every exact functor \(f\colon \Spaf \to \C\) the top horizontal map in~\eqrefone{equation:square-corepresentation} induces an equivalence on vertical fibres. Now by the construction of \(\Cath\) as  a cartesian fibration we have that the fibre of the left vertical map over \(f\colon \Spaf \to \C\) is given by the space of natural transformations \(\Nat(\QF^{\uni}, f^*\QF)\). On the other hand, the fibres of the right vertical map over \(f(\SS)\) is the space \(\Omega^\infty\QF(f(\SS))\) of hermitian forms on \(f(\SS) \in \C\) by construction. Lemma~\refone{lemma:maps-from-Qg} then implies that the induced map
\[
\Nat(\QF^{\uni},f^*\QF) \to f^*\QF(\SS) = \QF(f(\SS))
\]
on vertical fibres is an equivalence. This shows that \((\SS,\qone^{\uni})\) exhibits \((\Spaf,\QF^{\uni})\) as corepresenting the functor \(\spsforms\) in \(\Cath\). To obtain the analogous claim for Poincaré forms we observe that the corresponding map
\[
\Nat(\Bil_{\QF^{\uni}},(f \times f)^*\Bil_{\QF}) \to (f \times f)^*\Bil_{\QF}(\SS,\SS) = \Bil_{\QF}(f(\SS),f(\SS)),
\]
induced by the image of the universal Poincaré form in \(\Bil_{\QF}(\SS,\SS)\), identifies under Lemma~\refone{lemma:nat-duality} with the map
\begin{equation}
\label{equation:induced-nat-map}%
\Nat(f\Dual,\Dual_{\QF}f\op) \to \Map(f(\SS),\Dual_{\QF}f(\SS))
\end{equation}
which sends a natural transformation \(\tau\colon f\Dual\Rightarrow \Dual_{\QF}f\op\) to the composite
\[
f(\SS) \xrightarrow{f_*q^{\uni}_{\sharp}}f(\Dual\SS) \xrightarrow{\tau_{\SS}} \Dual_{\QF}f(\SS).
\]
Since \(q^{\uni}_{\sharp}\) is an equivalence this map sends natural equivalences \(f\Dual \xrightarrow{\simeq} \Dual_{\QF}f\op\) to equivalences \(f(\SS) \xrightarrow{\simeq} \Dual_{\QF}f(\SS)\). To finish the proof it will hence suffice to show that~\eqrefone{equation:induced-nat-map} also detects equivalences. Indeed, since \(\qone^{\uni}_{\sharp}\) is an equivalence, this follows from the fact that \(\SS\) generates \(\Spaf\) under finite colimits and so a natural transformation between two exact functors on \(\Spaf\) is an equivalence if and only if it evaluates to an equivalence on \(\SS\).
\end{proof}

\subsection{Ordinary rings and derived structures}
\label{subsection:discrete-rings}%

In this section we consider the case of ordinary rings
and explain how classical inputs for Grothendieck-Witt-  and \(\L\)-theory can be encoded as Poincaré structures on the \(\infty\)-category of perfect chain complexes. Our main result (see Proposition~\refone{proposition:restriction-of-2-excisive-functors} below) is that such Poincaré structures are essentially uniquely determined by their values on projective modules. This leads to the formation of derived versions of classical notions of hermitian forms, constituting the main link through which the point of view taken in this paper series interacts with its classical counterpart.

Let \(R\) be an ordinary associative ring.
Recall the Eilenberg-MacLane inclusion
\[
\GEM \colon \Ab = \Spa^{\heartsuit} \hrar \Spa
\]
of abelian groups as the heart of the canonical \(t\)-structure on spectra. Since this \(t\)-structure is compatible with tensor products the functor \(\GEM\) is naturally lax symmetric monoidal, and consequently \(\GEM R\) is naturally an \(\Eone\)-ring spectrum, or more precisely, an \(\Eone\)-algebra over the \(\Einf\)-ring spectrum \(\GEM \ZZ\).
We then have natural equivalences \cite[Theorem~7.1.2.1]{HA}
\begin{equation}
\label{equation:HZ-modules}%
\Mod_{\GEM \ZZ} \simeq \Der(\ZZ) \quad\text{and}\quad \Modp{\GEM R} \simeq \Dperf(R)
\end{equation}
between the \(\infty\)-categories of \(\GEM \ZZ\)-module spectra and compact \(\GEM R\)-modules in \(\GEM \ZZ\)-module spectra, as considered in~\S\refone{section:modules}, and the derived and perfect derived \(\infty\)-categories of \(\ZZ\) and \(R\), respectively, see Example~\refone{example:perfect-derived} and Example~\refone{example:ordinary-ring}. Similarly, we have an equivalence \(\Modf{\GEM R} \simeq \Dfree(R)\) between the \(\infty\)-category of finitely presented \(R\)-modules in \(\Mod_{\GEM \ZZ}\) and the finitely presented derived category \(\Dfree(R)\), obtained from the category of bounded complexes of finitely generated free \(R\)-modules by inverting quasi-isomorphisms.

\begin{notation}
In what follows we will need to consider both
ordinary tensor products over \(\ZZ\) and tensor product of \(\GEM \ZZ\)-module spectra over \(\GEM \ZZ\),
that which corresponds, under the equivalence~\eqrefone{equation:HZ-modules}, to the \emph{derived} tensor product \(\otimes^{\L}_\ZZ\) of complexes over \(\ZZ\). We will consequently write
\(\otimes_{\ZZ}\) for the former and \(\otimes_{\GEM \ZZ}\) for the latter. In particular, for two ordinary \(R\)-modules \(M\) and \(N\) one has a canonical map
\[
\GEM M \otimes_{\GEM \ZZ} \GEM N \to \GEM(M \otimes_\ZZ N)
\]
which is generally not an isomorphism, though it does exhibit its target as the 0-th truncation of its domain, and in particular determines an isomorphism
\[
\pi_0(\GEM M \otimes_{\GEM \ZZ} \GEM N)  \cong M \otimes_{\ZZ} N .
\]
\end{notation}

As described in \S\refone{section:modules}, we may construct bilinear functors on \(\Modp{\GEM R}\) from \(\GEM \ZZ\)-modules with involution over \(\GEM R\) (Definition~\refone{definition:module-with-involution}), and hermitian structures on \(\Modp{\GEM R}\) from \(\GEM \ZZ\)-modules with genuine involution over \(\GEM R\) (Definition~\refone{definition:module-with-genuine-involution}). In this context, it is natural to focus attention on \(\GEM \ZZ\)-modules with involution which arise by taking Eilenberg-MacLane spectra of
ordinary \(\ZZ\)-modules with involutions, by which we mean:

\begin{definition}
\label{definition:ordinary-module-with-involution}%
By a \(\ZZ\)-module with involution over \(R\) we will mean an \((R \otimes_\ZZ R)\)-module \(M\) in the ordinary sense together with an involution of abelian groups \(\sig\colon M \xrightarrow{\simeq} M\) which is linear over the flip isomorphism \(R \otimes_\ZZ R \xrightarrow{\simeq} R \otimes_\ZZ R\), that is, which satisfies \(\sig((a \otimes b)x) = (b \otimes a)\sig(x)\).
\end{definition}

We will generally write such \(\ZZ\)-modules with involution as pairs \((M,\sig)\). We will then denote by \(\Bil_M\) the bilinear functor on \(\Dperf(R)\) corresponding to the bilinear functor \(\Bil_{\GEM M}\) on \(\Mod_{\GEM R}\) via the equivalence~\eqrefone{equation:HZ-modules}. The following lemma ensures that the value of this bilinear form on projective modules is of a classical nature:

\begin{lemma}
\label{lemma:proj-perf}%
For \(P,Q \in \Prof{R}\) the canonical map
\begin{equation}
\label{equation:proj-perf-comp}%
\GEM \Hom_{R \otimes_\ZZ R}(P \otimes_\ZZ Q,M) \to \map_{\GEM R \otimes_{\GEM \ZZ} \GEM R}(\GEM P \otimes_{\GEM \ZZ} \GEM Q,\GEM M) = \Bil_M(P,Q)
\end{equation}
is an equivalence. In particular, \(\Bil_M\) sends pairs of projective modules to (Eilenberg-MacLane spectra of) abelian groups.
\end{lemma}
\begin{proof}
For \(P \in \Prof{R}\) the condition that~\eqrefone{equation:proj-perf-comp} is an equivalence is closed under direct sums and retracts, and so it suffices to check it for \(R=P\), where it can be identified with the identity map
\[
\GEM M \simeq \GEM \Hom_{R \otimes_\ZZ R}(R \otimes_\ZZ R,M) \to \map_{\GEM R \otimes_{\GEM \ZZ} \GEM R}(\GEM R \otimes_{\GEM \ZZ} \GEM R,\GEM M) \simeq \GEM M. \qedhere
\]
\end{proof}

The bilinear functor \(\Bil_M\) is perfect if and only if \(\GEM M\) is \emph{invertible}, that is, \(\GEM M\) is perfect as an \(\GEM R\)-module and the map \(\GEM R \to \End_{\GEM R}(\GEM M)\) is an equivalence (see Definition~\refone{definition:invertible}).

\begin{definition}
\label{definition:ordinary-invertible-module}%
We will say that a \(\ZZ\)-module with involution \(M\) over \(R\) (in the sense of Definition~\refone{definition:ordinary-module-with-involution}) is \emph{invertible}
if \(\GEM M\) is invertible as a \(\GEM H\)-module with involution over  \(\GEM R\) and in addition \(M\) is projective as an \(R\)-module.
\end{definition}

\begin{remark}
\label{remark:ordinary-invertible}%
For a \(\ZZ\)-module with involution \(M\) over \(R\) which is projective as an \(R\)-module we have \(\map_{\GEM R}(\GEM M,\GEM M) = \Hom_R(M,M)\), and so for such modules invertibility is equivalent to the map \(R \to \Hom_R(M,M)\) being an isomorphism.
\end{remark}

The projectivity assumption in Definition~\refone{definition:ordinary-module-with-involution} is meant to insure that the duality \(\Dual_M\colon \Dperf(R)\op \xrightarrow{\simeq} \Dperf(R)\) associated to the perfect bilinear functor \(\Bil_M\) preserves the (ordinary) full subcategory \(\Prof{R} \subseteq \Dperf(R)\) of finitely generated projective \(R\)-modules, and determines in particular a duality of ordinary categories
\[
\Dual_M\colon\Prof{R}\op \xrightarrow{\simeq} \Prof{R} \quad\quad \Dual_M(X) = \Hom_R(X,M) ,
\]
where \(\Hom_R(X,M)\) is given an \(R\)-module structure using the second \(R\)-action.
This is also consistent with the classical terminology concerning invertible modules, see, e.g., Example~\refone{example:commutative} just below.

\begin{example}
\label{example:commutative}%
If \(R\) is commutative then any \(R\)-module \(M\) can be considered as an \((R \otimes_\ZZ R)\)-module via the multiplication homomorphism \(R \otimes_\ZZ R \to R\). In particular, the two \(R\)-actions coincide, and we may endow \(M\) with the trivial involution. For a projective \(M\) this results in an invertible \(\ZZ\)-module with involution if and only if \(M\) is invertible as an object in the symmetric monoidal category \(\Prof{R}\). From the perspective of algebraic geometry, such modules correspond to \defi{line bundles} over \(\spec(R)\).
\end{example}

\begin{examples}
\label{examples:rings-with-anti-involution}%
Suppose that \(R\) is equipped with an \emph{anti-involution}, that is, an abelian group involution \(\inv{\bullet}\colon R \xrightarrow{\simeq} R\) which satisfies \(\inv{ab} = \inv{b}\inv{a}\). In this case \(R\) can be considered as a \(\ZZ\)-module with involution over itself via the \((R \otimes_\ZZ R)\)-action \((a \otimes b)x = ax\inv{b}\) and the involution \(\sig = \inv{\bullet}\). Some examples of interest of such rings include:
\begin{enumerate}
\item
Any commutative ring with an automorphism of order \(2\) gives rise to a ring with anti-involution. For example, the field \(\CC\) of complex numbers can be considered as a ring with anti-involution via complex conjugation.
\item
The group ring \(\ZZ[G]\) associated to a discrete group \(G\) carries a natural anti-involution given on additive generators by \(g \mapsto g^{-1}\). This example is a recurring one in geometric applications of \(\L\)-theory (see also \S\refone{subsection:visible}). More generally one can consider an \defi{orientation character} \(\chi\colon G \to C_2 = \{\pm 1\}\) and define the \(\chi\)-twisted anti-involution by setting \(g \mapsto \chi(g) g^{-1}\).
\item
For a commutative ring \(k\) the \(k\)-algebra of \(n \times n\)-matrices \(\Mat_n(k)\) admits an anti-involution \(A \mapsto A^{t}\) given by sending a matrix to its transpose. This more generally works for \(k\) a ring with anti-involution.
\item
For a commutative ring \(k\) and \(a,b \in k\), the quaternion \(k\)-algebra \(\Q_k(a,b)\) is the algebra generated over \(k\) by elements \(i,j\) under the relation \(i^2=a,j^2=b,ij=-ji\). It admits an anti-involution sending \(i\) to \(-i\) and \(j\) to \(-j\).
\end{enumerate}
\end{examples}

Another common source of invertible \(\ZZ\)-modules with involution %
is the following. Recall that an \defi{anti-structure} in the sense of
Wall~\cite{wall1970axiomatic} on a ring \(R\) consists of a ring isomorphism
\[
\inv{\bullet}\colon R\op\to R,
\]
together with a unit \(\epsilon\in R^*\), such that \(\inv{\epsilon}=\epsilon^{-1}\) and \(\dovl{r}=\epsilon^{-1} r\epsilon\). In particular, if \(\epsilon\) belongs to the center of \(R\) then \(\inv{\bullet}\) is an anti-involution, and this is arguably the most common case studied in the literature. Specifically, one often considers the case where \(\epsilon=\pm 1\), which, for example, in the case of the integers \(\ZZ\), are also the only possibilities.

Given a Wall anti-structure \((\inv{\bullet},\eps)\), we may consider \(R\) as an \((R \otimes_\ZZ R)\)-module with action given by \((a\otimes b) \cdot r=ar\inv{b}\). The map \(a \mapsto \epsilon\inv{a}\) is then an involution on \(R\) which is linear over the flip action on \(R \otimes_\ZZ R\), and so we obtain the structure of a \(\ZZ\)-module with involution.
This \(\ZZ\)-module with involution is always \emph{invertible}: the induced map \(R \to \Hom_R(R,R)\cong R\op\) identifies with \(r \mapsto \ovl{r}\) (see Remark~\refone{remark:ordinary-invertible}).

\begin{remark}
A Wall anti-structure captures the most general form of an invertible \(\ZZ\)-module with involution over \(R\) whose underlying \(R\)-module is \(R\). Indeed, giving the \(R\)-module \(R\) a second commuting \(R\)-action is equivalent to providing a ring homomorphism \(\inv{(\bullet)}\colon R \to \Hom_R(R,R) = R\op\), which is furthermore an isomorphism if the desired \((R \otimes_\ZZ R)\)-module is to be invertible (see Remark~\refone{remark:ordinary-invertible}). The \((R \otimes_\ZZ R)\)-action can then be written in terms of \(\inv{(\bullet)}\) by \((a \otimes b)(c) = ac\ovl{b}\). If \(\sig\colon R \to R\) is now an abelian group isomorphism which switches the two \(R\)-action then \(\sig\) is completely determined by the value \(\eps := \sig(1) \in R\), in terms of which \(\sig\) can be written as
\[
\sig(r) = \sig(r\cdot 1) = \sig(1)\ovl{r} = \eps\ovl{r}.
\]
Since \(\sig\) and \(\ovl{(-)}\) are both isomorphisms of abelian groups so is the map \(r \mapsto \eps \cdot r\), and hence \(\eps\) must be a unit. In addition, since \(\sig\) switches the two \(R\)-actions we also have
\[
r\eps= r\sig(1) = \sig(1\cdot \ovl{r}) = \eps\dovl{r}
\]
and hence \(\dovl{r} = \eps^{-1}r\eps\).
Finally, the condition that \(\sig\) is an involution implies that
\[
1 = \sig\sig(1) = \eps\ovl{\eps}
\]
and hence \(\ovl{\eps} = \eps^{-1}\). In particular, the pair \((\inv{\bullet},\eps)\) is a Wall anti-structure on \(R\) and the \(\ZZ\)-module with involution we obtain is the one associated to that structure.
\end{remark}

\begin{remark}
If \((\inv{\bullet},\eps)\) is a Wall anti-structure on a ring \(R\) and \(u \in R^*\) is a unit then we can obtain a new Wall anti-structure by replacing \(\inv{(\bullet)}\) with \(u^{-1}\inv{(-)}u\) and \(\eps\) with \(\eps\big(\ovl{u}\big)^{-1}u\). One then says that the two Wall anti-structures \((\inv{\bullet},\eps)\) and \((u^{-1}\inv{(-)}u,\eps\big(\ovl{u}\big)^{-1}u)\) are \defi{conjugated}. In this case the associated \(\ZZ\)-modules with involution over \(R\) are isomorphic via the map \(R \to R\) sending \(x\) to \(xu\). In fact, any isomorphism between the \(\ZZ\)-modules with involution associated to two Wall anti-structures is of this form, and so two Wall anti-structures are conjugated if and only if their associated \(\ZZ\)-modules with involution are isomorphic as such.
\end{remark}

The construction below, which was shared with the authors by Uriya First, gives an example of a Wall anti-structure which is not conjugated to any central Wall anti-structure:

\begin{example}
Let \(K\) be a field which admits an automorphism \(\sigma\colon K \to K \) of order \(4\). Define \(R = K[x,x^{-1};\sigma^2]\) to be the twisted Laurent polynomial ring generated over \(K\) by an invertible generator \(x\) which satisfies the relation \(x^{-1}\alp x = \sigma^2(\alpha)\) for \(\alp \in K\). We may then extend \(\sig\) to an anti-automorphism on \(R\) defined on monomials by
\[
\inv{\alp x^i} = x^{-i}\sig(\alp) = \sig^{1-2i}(\alp)x^{-i} .
\]
Then \(\dovl{\alp x^i} = \sig^2(\alp)x^i = x^{-1}(\alp x^i) x\) and \(x\ovl{x} = 1\), so that we obtain a Wall anti-structure \((\ovl{\bullet},\eps)\) with \(\eps = x\).
This Wall anti-structure is not conjugated to any central Wall anti-structure: indeed, the units of \(R\) are exactly the monomials \(\alpha x^i\), and if we conjugate the above Wall anti-structure by \(\alpha x^i\) then the new \(\eps\) will be of the form \(\beta x^{2i+1}\) for a suitable \(\beta \in K\),
and as such cannot be in the center, since it does not commute with \(K \subseteq R\).
\end{example}

Let us now also give a non-commutative example of an invertible \(\ZZ\)-module with involution whose underlying module is not the ring itself:
\begin{example}
Let \(B=\Q_{\QQ}(-5,-13)\) be the quaternion algebra over \(\QQ\) and let \(A \subseteq B\) be the subring generated over \(\ZZ\) by \(1,i,j\) and \(\beta = \frac{1+i+j+ij}{2}\). Then \(A\) is a maximal order in the quaternion algebra \(B\), that is, it is finitely generated as a \(\ZZ\)-module and is not contained in any other subring with this property.
Invertible bimodules over maximal orders are relatively well understood, and can all be realized as invertible two-sided ideals in \(A\). In particular, the two-sided ideal \(I \subseteq A\) which is generated by \(2,i-1,j-1\) is an invertible ideal of  index \(4\) in \(A\) (the quotient \(A/I\cong \FF_2[\beta]/(\beta^2-\beta+1)\) is a finite field of order \(4\)). One can then verify that this ideal is not principal by checking that \(A\) contains no elements of norm \(2\), and so this \(A\)-bimodule is not isomorphic to \(A\) as a left \(A\)-module. At the same time, the involution on \(B\) sending \(i\) to \(-i\) and \(j\) to \(-j\) restricts to an involution on \(A\), through which we can consider \(I\) as an \(A \otimes A\)-module, and since this involution preserves \(I\) it endows it with the structure of an involution which is linear over the flip map \(A \otimes A \to A \otimes A\). Then \(I\) gives an invertible \(\ZZ\)-module with involution over \(A\) which is not isomorphic to \(A\) as a left \(A\)-module.
\end{example}

We now fix a ring \(R\) and an invertible \(\ZZ\)-module with involution \(M\) over \(R\). Let us say that a Poincaré structure \(\QF\) on \(\Dperf(R)\) is \emph{compatible with \(M\)} if its associated duality is given by \(\Dual_M\) (equivalently, if its symmetric bilinear part is \(\Bil_M\)).

\begin{example}
As in Definition~\refone{definition:functors-associated-to-module-with-involution} we have the quadratic and symmetric Poincaré structures \(\QF^{\qdr}_M\) and \(\QF^{\sym}_M\) associated to the bilinear functor \(\Bil_M\), which are compatible with \(M\) by construction.
\end{example}

\begin{example}
As in Example~\refone{example:truncation}, for an integer \(m\in\mathbb{Z}\) one can consider the associated \(\GEM\ZZ\)-module with genuine involution
\[
(\GEM M,\tau_{\geq m}\GEM M^{\tC},t_m\colon \tau_{\geq m}\GEM M^{\tC}\to \GEM M^{\tC})
\]
over \(\GEM R\),
where \(t_m\colon \tau_{\geq m}\GEM M^{\tC}\to \GEM M^{\tC}\)
is the \(m\)-connective cover of \(\GEM M^{\tC}\).
The associated Poincaré structure
\[
\QF^{\geq m}_M\colon \Dperf(R)^{\op}\to\Spa.
\]
defined as in Construction~\refone{construction:functors-associated-to-module-with-genuine-involution}, is then compatible with \(M\).
\end{example}

\begin{remark}
By Theorem~\refone{theorem:classification-genuine-modules} the data of a Poincaré structure on \(\Dperf(R)\) compatible with \(M\) is equivalent to that of an \(\GEM R\)-module spectrum \(N\) equipped with a map of \(\GEM R\)-module spectra \(N \to \GEM M^{\tC}\). The above examples then correspond to the case where this map is the identity (the symmetric structure), the map from the zero module (the quadratic structure), or the various connective covers.
\end{remark}

\begin{warning}
\label{warning:frobenius-2}%
In the construction above, the \(\GEM R\)-module \(\GEM M^{\tC}\) is such that the induced \(R\)-module structure on \(\pi_0(\GEM M)^{\tC} = \hat{\mathrm{H}}_0(\Ct,M) = \coker[\rT_{\Ct}\colon M_{\Ct} \to M^{\Ct}]\) is given by
\[
r[m] = [(r \otimes r)m]
\]
where \([m]\) is the class mod norms of a \(\Ct\)-fixed element \(m \in M\).
This can lead to some confusion
when \(R\) is commutative and the involution on \(M\) is \(R\)-linear. In this case, another natural action of \(R\) on \(M\) is available by realizing the Tate construction in \(\GEM R\)-modules. The associated action on homotopy groups is then given by \(r[m] = [rm]\), which is generally \emph{different} from the action above (see also Warning~\refone{warning:frobenius}).
\end{warning}

\begin{remark}
\label{remark:linear-colinear-parts}%
By construction the linear part of the Poincaré structure \(\QF^{\geq m}_M\) is given by the formula \(X \mapsto \map_{\GEM R}(X,\tau_{\geq m}\GEM M^{\tC})\), where in the last term we have identified \(X\) with the corresponding \(\GEM R\)-module spectrum via the equivalence~\eqrefone{equation:HZ-modules}.
In particular, \(\QF^{\geq m}_M(X)\) sits in an exact sequence of spectra
\[
\QF^{\qdr}_M(X) \to \QF^{\geq m}_M(X) \to \map_{\GEM R}(X,\tau_{\geq m}\GEM M^{\tC}) ,
\]
which can be used to obtained connectivity estimates on the gap between \(\QF^{\geq m}_M\) and \(\QF^{\qdr}_M\). On the other hand, from the fibre sequence \(\tau_{\geq m}\GEM M^{\tC} \to \GEM M^{\tC} \to \tau_{\leq m-1}\GEM M^{\tC}\) we see that \(\QF^{\geq m}_M(X)\) also sits in an exact sequence of spectra
\[
\QF^{\geq m}_M(X) \to \QF^{\sym}_M(X) \to \map_{\GEM R}(X,\tau_{\leq m-1}\GEM M^{\tC}),
\]
and this can be used to obtain coconnectivity estimates on the gap between \(\QF^{\geq m}_M\) and \(\QF^{\sym}_M\).
\end{remark}

Our principal goal in this section to construct a link between Poincaré structures on \(\Dperf(R)\) compatible with a given \(M\) and the classical framework of \(M\)-valued forms on projective \(R\)-modules. For this, let \(\Ch(\Prof{R})\) be the category of bounded chain complexes of finitely generated projective \(R\)-modules, so that the \(\infty\)-category \(\Dperf(R)\) can be identified with the \(\infty\)-categorical localisation of \(\Ch(\Prof{R})\) by the collection of quasi-isomorphisms. The inclusion \(\Prof{R} \subseteq \Ch(\Prof{R})\) as chain-complexes concentrated in degree \(0\) determines a fully-faithful functor
\begin{equation}
\label{equation:proj-to-perf}%
\Prof{R} \to \Dperf(R), \quad\quad P \mapsto P[0].
\end{equation}
We also point out that the category \(\Prof{R}\) is additive, and the inclusion~\eqrefone{equation:proj-to-perf} is additive in the sense that it preserves direct sums. One can then show that~\eqrefone{equation:proj-to-perf} exhibits \(\Dperf(R)\) as the initial stable \(\infty\)-category equipped with additive functor from \(\Prof{R}\).
It is also sometimes called the \defi{stable envelope} of \(\Prof{R}\).

The following definition is originally due to Eilenberg and MacLane \cite[Theorem~9.11]{eilenberg-maclane}.
\begin{definition}
\label{definition:quadratic-additive}%
Let \(\C,\D\) be additive \(\infty\)-categories. We will say
that a reduced (that is, zero object preserving) functor \(\QF\colon \C \to \D\) is
\emph{polynomial of degree \(2\)} if its cross-effects \(\Bil_{\QF}(X,Y)\) (defined as in Definition~\refone{definition:cross-effect} as the kernel of the split surjection \(\QF(X \oplus Y) \to \QF(X) \oplus \QF(Y)\)) preserves direct sums in each variable separately. We will denote by
\[
\Funpoly(\C,\D) \subseteq \Fun(\C,\D)
\]
the full subcategory spanned by reduced functors which are polynomial of degree \(2\).
\end{definition}

We now arrive to the main result of this section:

\begin{proposition}
\label{proposition:restriction-of-2-excisive-functors}%
Let \(R\) be an associative ring. Then restriction along the inclusion \(\Prof{R} \subseteq \Dperf(R)\) yields an equivalence of \(\infty\)-categories
\[
\Funq(\Dperf(R)) \xrightarrow{\simeq} \Funpoly(\Prof{R}\op,\Spa) .
\]
\end{proposition}

Since the Eilenberg-MacLane embedding \(\GEM\colon \Ab \hrar \Spa\) is fully-faithful and additive, Proposition~\refone{proposition:restriction-of-2-excisive-functors} implies in particular that if \(\QF_{\proj}\colon \Prof{R}\op \to \Ab\) is a reduced polynomial functor of degree \(2\) then the composed functor \(\Prof{R}\op \to \Ab \to \Spa\) extends to a quadratic functor \(\QF\colon \Dperf(R)\op \to \Spa\) in an essentially unique manner. We will refer to such an extension
\[
\QF \colon \Dperf(R)\op \to \Spa
\]
as the \defi{derived quadratic functor} of \(\QF_{\proj}\).

\begin{remark}
For perfect complexes \(X \in \Dperf(R)\) concentrated in non-negative degrees one can express the value of the derived functor \(\QF\) as follows.
One first translates \(X\) to a simplicial \(R\)-module using the Dold-Kan correspondence. Applying the functor \(\QF_{\proj}\) levelwise yields a cosimplicial abelian group
which one can then re-translate into a non-positively graded chain complex over \(\ZZ\), and consequently into a spectrum. This is the classical description of non-abelian derived functors on non-negatively graded complexes due to Dold and Puppe \cite{dold-puppe}.
\end{remark}

The proof of Proposition~\refone{proposition:restriction-of-2-excisive-functors} will be given below. Before, let us explore some of its consequences. As their higher categorical counterparts,
reduced degree \(2\) polynomial functors \(\Prof{R}\op \to \Ab\) can be used to encode various types of hermitian forms.
To make this more explicit, consider for a projective module \(P \in \Prof{R}\) the abelian group \(\Hom_{R\otimes R}(P\otimes_{\ZZ} P,M)\) of
\(M\)-valued \(R\)-bilinear forms \(\beta\colon P \otimes_\ZZ P\to M\).
This abelian group carries an involution
which sends a form \(\beta\) to the form
\[
(v, u)\mapsto \sig(\beta(u,v)),
\]
where \(\sig\colon M \to M\) is the involution of \(M\). The \(\Ct\)-orbits and \(\Ct\)-fixed points of \(\sig\) are then related via the trace map
\[
\rT_{\Ct}\colon \Hom_{R\otimes_{\ZZ} R}(P\otimes_{\ZZ} P,M)_{\Ct} \to \Hom_{R\otimes_{\ZZ} R}(P\otimes_{\ZZ} P,M)^{\Ct} \quad\quad [\beta] \mapsto \beta(v,u) + \sig\beta(u,v),
\]
which sends an orbit to the sum of its representatives.

\begin{definition}
\label{definition:ordinary-forms}%
Let \(R\) be a ring and \((M,\sig)\) an invertible \(\ZZ\)-module with involution over \(R\).
We define functors \(\Prof{R}\op \to \Ab\)
by the formulas
\[
\QF^{\qdr}_{\proj}(P) = \Hom_{R\otimes_\ZZ R}(P\otimes_{\ZZ} P,M)_{\Ct},
\qquad
\QF^{\sym}_{\proj}(P) = \Hom_{R\otimes_\ZZ R}(P\otimes_{\ZZ} P,M)^{\Ct},
\]
and
\[
\QF^{\ev}_{\proj}(P) = \im\big[\QF^{\qdr}_{\proj}(P) \xrightarrow{\rT_{\Ct}} \QF^{\sym}_{\proj}(P)\big] .
\]
For \(P \in \Prof{R}\) we will refer to these as the abelian groups of \defi{\(\sig\)-quadratic}, \defi{\(\sig\)-symmetric}, and \defi{\(\sig\)-even} forms on \(P\), respectively. These functors are visibly reduced and the cross-effect of each of them is \((P,Q) \mapsto \Hom_{R \otimes_\ZZ R}(P \otimes_\ZZ Q,M)\), which is additive in each variable separately. In particular, they are polynomial of
degree \(2\).
\end{definition}

\begin{remark}
\label{remark:genuine-pi-0}%
It follows from Lemma~\refone{lemma:proj-perf} that
\[
\begin{split}
\QF^{\qdr}_{\proj}(P) = \Hom_{R \otimes_\ZZ R}(P \otimes_\ZZ P,M)_{\Ct} & \cong [\pi_0\Bil_M(P,P)]_{\Ct} \\
 & \cong \pi_0[\Bil_M(P,P)_{\hC}] = \pi_0\QF^{\qdr}_M(P[0]),
\end{split}
\]
and similarly
\[
\begin{split}
\QF^{\sym}_{\proj}(P) = \Hom_{R \otimes_\ZZ R}(P \otimes_\ZZ P,M)^{\Ct} & \cong [\pi_0\Bil_M(P,P)]^{\Ct} \\
 & \cong \pi_0[\Bil_M(P,P)^{\hC}] = \pi_0\QF^{\sym}_M(P[0]),
\end{split}
\]
\end{remark}

By definition, a \(\sig\)-symmetric form
is an \(R\)-bilinear form \(\phi\colon P \otimes_\ZZ P \to M\) such that \(\phi(b,a) = \sig(\phi(a,b))\).
On the other hand, the data of a \(\sig\)-quadratic form, or a \(\Ct\)-orbit \([\beta] \in \Hom_{R\otimes R}(P\otimes_{\ZZ} P,M)_\Ct\), is equivalent to that of a pair \((\phi,q)\), where \(\phi \in \Hom_{R\otimes R}(P\otimes_{\ZZ} P,M)^{\Ct}\) is a \(\sig\)-symmetric form
and \(q\colon P \to M_{\Ct}\) is a set-theoretic function which satisfies
\begin{enumerate}
\item \(q(v+u)-q(v)-q(u) = [\phi(v,u)] \in M_{\Ct}\) for \(v,u \in P\);
\item \(q(r v) = (r \otimes r)q(v)\) for \(v \in P\) and \(r \in R\);
\item the image of \(q(v)\) under the trace map \(M_{\Ct} \to M^{\Ct}\) is the \(\Ct\)-fixed element \(\phi(v,v)\) for \(v \in P\).
\end{enumerate}
To obtain this description, note that the abelian group of such pairs \((\phi,q)\) forms a reduced degree 2 polynomial functor \(\QFtwo\colon \Prof{R} \to \Ab\), which receives a natural transformation \(\QF^{\qdr}_{\proj} \Rightarrow \QFtwo\) sending \([\beta] \in \QF^{\qdr}_{\proj}(P)\) to the pair \((\phi,q)\), where \(\phi = \rT_{\Ct}[\beta]\) is the trace of \(\beta\) and \(q(x) = [\beta(x,x)] \in M_{\Ct}\). One can then verify in a straightforward manner that this natural transformation induces an isomorphism on cross-effects and an isomorphism on the value on \(P=R\), and is hence an isomorphism on every \(P \in \Prof{R}\) (see also~\cite[Theorem 1]{wall1970axiomatic}, where this argument is elaborated in the case where \(M\) comes from a Wall anti-structure).

When \(R\) is commutative and \(M=R\) with trivial involution the above notion of a quadratic form identifies with the usual one. In this case even forms are symmetric forms which admit a quadratic refinement in the classical sense (which is not kept as part of the structure). For example, when \(M=R=\ZZ\) with trivial involution then a symmetric bilinear form \(b\) on \(P\) admits a quadratic refinement if and only if  \(b(x,x) \in 2\ZZ\) for all \(x\in P\), hence the terminology ``even forms''. In this case the quadratic refinement is even unique, though this is by no means the case in general. If \(R\) is commutative and \(M=R\) with involution \(\sig(x)=-x\) then \(\sig\)-symmetric forms are skew-symmetric forms, while the \(\sig\)-even forms are the alternating ones.
For non-commutative \(R\) this way of viewing quadratic forms was first devised by Tits~\cite{tits1968formes} for central simple algebras, and later generalized by Wall~\cite{wall1970axiomatic} to arbitrary rings with anti-structure as above.

By Proposition~\refone{proposition:restriction-of-2-excisive-functors} each of the functors \(\QF^{\qdr}_{\proj},\QF^{\ev}_{\proj}\) and \(\QF^{\sym}_{\proj}\) extends
in an essentially unique manner to its corresponding derived quadratic functor \(\Dperf \to \Spa\). The following proposition identifies these in terms of Poincaré structures we have already encountered:

\begin{proposition}
\label{proposition:classical-derived}%
The derived quadratic functors of \(\QF^{\qdr}_{\proj},\QF^{\ev}_{\proj}\) and \(\QF^{\sym}_{\proj}\) are canonically equivalent to the Poincaré structures \(\QF^{\geq 2}_{M}, \QF^{\geq 1}_{M}\) and \(\QF^{\geq 0}_{M}\), respectively.
\end{proposition}
\begin{proof}
By Lemma~\refone{lemma:proj-perf} we have that for \(P \in \Prof{P}\) the spectrum \(\Bil_M(P[0],P[0])\) belongs to \(\Spa^{\heartsuit}\) and so \(\QF^{\qdr}_M(P[0])\) is connective and \(\QF^{\sym}_M(P[0])\) is coconnective. By Remark~\refone{remark:linear-colinear-parts} we then get that for \(m=0,1,2\) the spectrum \(\QF^{\geq m}_M(P[0])\) is both connective and coconnective, and hence lies in the heart as well.
Now consider the pair of maps
\[
\pi_0\QF^{\qdr}_M(P[0]) \to \pi_0\QF^{\geq m}_M(P[0]) \to \pi_0\QF^{\sym}_M(P[0]).
\]
By Remark~\refone{remark:linear-colinear-parts} the first map above is an isomorphism when \(m=2\) and the second map is an isomorphism when \(m=0\).
Finally,
when \(m=1\) the same remark gives that the first map is surjective and the second is injective.
In light of Remark~\refone{remark:genuine-pi-0}
we now get an identification of \(\pi_0\QF^{\geq m}_M(P[0])\) for \(m=0,1,2\) with \(\QF^{\sym}_{\proj}(P),\QF^{\ev}_{\proj}(P)\) and \(\QF^{\qdr}_{\proj}(P)\)  respectively. By the uniqueness of Proposition~\refone{proposition:restriction-of-2-excisive-functors} we may identify \(\QF^{\geq m}_M\) for \(m=0,1,2\) with the desired derived quadratic functors.
\end{proof}

\begin{notation}
\label{notation:genuine-zero-one-two-quadratic-functors}%
In light of Proposition~\refone{proposition:classical-derived} we will denote the Poincaré structures \(\QF^{\geq 2}_M\), \(\QF^{\geq 1}_M\) and \(\QF^{\geq 0}_M\) also by \(\QF^{\gq}_M\), \(\QF^{\gev}_M\) and \(\QF^{\gs}_M\), and refer to them as the \defi{genuine quadratic}, \defi{genuine even} and \defi{genuine symmetric} Poincaré structures on \(R\) associated to \(M\).
\end{notation}

Combining Proposition~\refone{proposition:classical-derived}
with Corollary~\refone{corollary:periodicity-truncated} we conclude:

\begin{corollary}
\label{corollary:shift-genuine}%
The loop functor \(\Om\colon \Dperf(R) \to \Dperf(R)\) refines to equivalences of Poincaré \(\infty\)-categories
\[
\big(\Dperf(R), (\QF_{M}^{\gs})\qshift{2}\big) \xrightarrow{\simeq} \big(\Dperf(R), \QF^{\gev}_{-M}\big)
\]
and
\[
\big(\Dperf(R), (\QF_{M}^{\gev})\qshift{2}\big) \xrightarrow{\simeq} \big(\Dperf(R), \QF^{\gq}_{-M}\big).
\]
\end{corollary}

\begin{remark}
\label{remark:relation-to-classical}%
The \(\L\)-groups of \(\Dperf(R)\) with respect to Poincaré structures \(\QF^{\qdr}_M\) and \(\QF^{\sym}_M\) identify with Ranicki's 4-periodic quadratic and symmetric \(\L\)-groups, see Example~\refone{example:ranicki-L-groups}. The \(\L\)-groups of the genuine Poincaré structures \(\QF^{\gq}_M,\QF^{\gev}_M,\QF^{\gs}_M\) are generally not 4-periodic, though one has the relation
\[
\L_{n-2}(\Dperf(R),\QF^{\gs}_M) \simeq \L_{n}(\Dperf(R),\QF^{\gev}_{-M}) \simeq \L_{n+2}(\Dperf(R),\QF^{\gq}_{M})
\]
by Corollary~\refone{corollary:shift-genuine}.
We will show in \paperthree that the \(\L\)-groups of the genuine symmetric Poincaré structure actually rediscover an older non-periodic variant of the symmetric \(\L\)-groups, originally introduced by Ranicki~\cite{RanickiATS1} using \(n\)-dimensional Poincaré complexes to define the \(n\)'th \(\L\)-group.
\end{remark}

The notions of \(\sig\)-symmetric, \(\sig\)-even and \(\sig\)-quadratic \(M\)-valued forms were generalized by Bak~\cite{bak-form-parameter} to the setting where one is given, in addition to \(M\),
a subgroup \(\fpmg \subseteq M^{\Ct}\) containing the image of the trace map \(M_{\Ct} \to M^{\Ct}\), and closed under the quadratic action of \(R\) on \(M^{\Ct}\) given by \((a,m) \mapsto (a \otimes a)m\). Here we use the term quadratic action to designate the fact that it is encoded by a \emph{quadratic} monoid map \((R,\cdot,1,0) \mapsto \End_{\ZZ}(M^{\Ct})\), rather than a ring homomorphism. In particular, the cross effect
\[
(a\perp b)m := \big((a+b) \otimes (a+b)\big)m-(a\otimes a)m - (b \otimes b)m = (a\otimes b + b \otimes a)m
\]
of this action is not trivial, but instead given by a bilinear map \(R \otimes R \to \End_{\ZZ}(M^{\Ct})\), which in this case is induced by the symmetrization of the \((R \otimes R)\)-action on \(M\). We will refer to abelian groups equipped with this kind of action as \emph{quadratic modules}.
This notion was subsequently generalized by Schlichting~\cite{SchlichtinghigherI} by removing the condition that \(\fpmg\) injects in \(M\):

\begin{definition}[Schlichting]
\label{definition:form-parameter}%
Let \(R\) be a ring equipped with invertible \(\ZZ\)-module with involution \(M\). A \emph{form parameter} \(\fpm\) on \((R,M)\) is a quadratic \(R\)-module \(\fpmg\) lying in a sequence of quadratic modules
\[
M_{\Ct} \xrightarrow{\tau} \fpmg \xrightarrow{\rho} M^{\Ct}
\]
whose composition is the trace map, and such that the cross-effect of the \(R\)-action on \(\fpmg\) satisfies \((a\perp b)x = \tau((a\otimes b)\rho(x))\).
\end{definition}

We note that any reduced degree 2 polynomial functor \(\QF\colon \Proj(R)\op \to \Ab\) with cross-effect \(\Bil_{\QF}\)
determines a form parameter
on \((R,\Bil_{\QF}(R,R))\)
by taking \(\fpmg = \QF(R)\), with the maps
\[
\Bil_{\QF}(R,R)_{\Ct} \to \QF(R) \to \Bil_{\QF}(R,R)^{\Ct}
\]
induced by the \(\Ct\)-equivariant diagonal and collapse maps relating \(R\) and \(R \oplus R\). Schlichting then proves (see~\cite[Lemma A.16]{SchlichtinghigherI}) that the association \(\QF \mapsto (\Bil_{\QF}(R,R),\QF(R))\) determines an equivalence between the category \(\Funpoly(\Proj(R)\op,\Ab)\) and the category of form parameters as above (more precisely, Schlichting %
proves this for \(\Proj(R)\) replaced with the category of finitely generated free \(R\)-modules, but \(\Funpoly(-,\Ab)\) is invariant under idempotent completion).
Explicitly, an inverse to this map is given by sending a form parameter \(\fpm = [M_{\Ct} \xrightarrow{\tau} \fpmg \xrightarrow{\rho} M^{\Ct}]\) to the functor
\[
\QF^{\fpm}_{\proj}\colon \Proj(R)\op \to \Ab
\]
which associates to a projective \(R\)-module \(P\) its group of \emph{\(\fpm\)-hermitian forms}, which by definition consists of pairs \((\phi,q)\) where \(\phi\colon P \otimes_{\ZZ} P \to M\) is a \(\sig\)-symmetric \(M\)-valued form and \(q\colon P \to \fpmg\) is a map of sets such that
\begin{enumerate}
\item \(q(r v) = rq(v)\) for \(v \in P\) and \(r \in R\), where \(rq(v)\) stands for the (quadratic) action of \(r\) on \(q(v) \in \fpmg\).
\item \(q(v+u)-q(v)-q(u) = \tau(\phi(v,u)) \in \fpmg\) for \(v,u \in P\);
\item \(\rho q(v) = \phi(v,v)\).
\end{enumerate}
The cross effect of \(\QF^{\fpm}_{\proj}\) is then given by \((P,P') \mapsto \Hom_R(P \otimes_{\ZZ} P',M)\).
In particular, for \(\fpmg=M^{\Ct}\) we get the notion of a \(\sig\)-symmetric \(M\)-valued form, for \(\fpmg=M_{\Ct}\) we get the notion of a \(\sig\)-quadratic \(M\)-valued form, and for \(\fpmg=\im[M_{\Ct} \to M^{\Ct}]\) we get the notion of a \(\sig\)-even \(M\)-valued form. Using Proposition~\refone{proposition:restriction-of-2-excisive-functors} we may extend the composite functor \(\Proj(R)\op \to \Ab \hrar \Spa\) to a quadratic functor
\[
\QF^{\gfpm}_M\colon \Dperf(R)\op \to \Spa
\]
in an essentially unique manner.
By Lemma~\refone{lemma:proj-perf} the bilinear part of \(\QF^{\gfpm}_M\) is \(\Bil_M\) and it is hence Poincaré and compatible with \(M\). Writing \(\Lin^{\gfpm}_M := \Lin_{\QF^{\gfpm}_M}\) for the linear part of \(\QF^{\gfpm}_M\) and evaluating at \(R \in \Dperf(R)\)
we obtain a fibre sequence
\[
\Bil_{M}(R,R)_\hC \to \QF_M^{\gfpm}(R) \to \Lin^{\gfpm}_{M}(R),
\]
which we may rewrite as
\[
M_\hC \to \fpmg \to \Lin^{\gfpm}_{M}(R).
\]
This determines the linear part of \(\QF^{\gfpm}_M\) to be \(\Lin^{\gfpm}_{M} \simeq \map_R(-,N)\), where the underlying spectrum of \(N\) is given by the cofibre of the composed map \(M_{\hC} \to M_{\Ct} \to \fpmg\). In particular, in the classification of Theorem~\refone{theorem:classification-genuine-modules} the Poincaré structure \(\QF^{\gfpm}_M\) corresponds to the module with genuine involution
\[
\gfpm:= (M,\cof[M_\hC \rightarrow \fpmg], \cof[M_\hC \rightarrow \fpmg] \rightarrow M^{\tC}).
\]

\begin{remark}
\label{remark:form-parameters-and-modules-with-genuine-involution}%
To clarify the relation between form parameters and modules with genuine involution one may use the framework of genuine \(\Ct\)-spectra. Recall from Remark~\refone{remark:NormModules} that modules with genuine involution over the \(\Eone\)-ring spectrum \(\GEM R\) can be identified with module objects in genuine \(\Ct\)-spectra over the
Hill-Hopkins-Ravenel norm \(\N \GEM R\), the latter being the genuine \(\Ct\)-spectrum whose underlying \(\Ct\)-spectrum is \(\GEM R\otimes_{\SS} \GEM R\) with the flip \(\Ct\)-action, whose geometric fixed point spectrum is \(\GEM R\), and whose reference map \(\GEM R \to (\GEM R\otimes_{\SS} \GEM R)^{\tC}\) is the Tate-diagonal. We then note that the stable \(\infty\)-category of genuine \(\Ct\)-spectra admits a t-structure whose connective objects are those whose underlying spectrum and geometric fixed points are both connective. In particular, \(\N \GEM R\) is a connective algebra object in \(\Ct\)-spectra, and so its \(\infty\)-category of modules inherits a t-structure created by the forget-the-action functor. Let us now switch to viewing genuine \(\Ct\)-spectra as \(\Spa\)-valued Mackey functors, that is, product preserving functors from the span \(\infty\)-category of finite \(\Ct\)-sets to \(\Spa\), see~\cite{Barwick-MackeyI} and~\cite{barwick2016parametrized}. In this framework the above t-structure is the objectwise one: indeed, the condition that the underlying spectrum and the geometric fixed points spectrum are both connective is equivalent to the condition that the underlying spectrum and genuine fixed points being both connective. As a result, the heart of this t-structure can be identified with the category of \(\Ab\)-valued Mackey functors, and similarly the heart of the induced t-structure on \(\N\GEM R\)-modules can be identified with \(\Ab\)-valued Mackey functors which are modules over \(\pi_0\N\GEM R\). The latter can be computed by translating \(\N\GEM R\) into an \(\Spa\)-valued Mackey functor and then taking \(\pi_0\) objectwise. This computation was carried out in~\cite[Proposition 5.5]{DMPR} yielding the \(\Ab\)-valued Mackey functor encoded by the arrows
\begin{equation}
\label{equation:pi-0-norm}%
(R \otimes_{\ZZ} R)_{\Ct} \xrightarrow{\tau} W \xrightarrow{\rho} (R \otimes_{\ZZ} R)^{\Ct},
\end{equation}
where \(W\) is the abelian group whose underlying set of elements is \(R \times (R \otimes_{\ZZ} R)_{\Ct}\), and whose addition law is given by
\[
(a,v) + (b,u) = (a+b,v+u+[a \otimes b]) .
\]
Here, the map \(\tau\) is the inclusion \(v \mapsto (0,v)\) and the map \(\rho\) sends \((a,v)\) to \((a \otimes a) + \rT_{\Ct}(v)\), where \(\rT_{\Ct}\) is the trace map from \((R \otimes_{\ZZ} R)_{\Ct}\) to \((R \otimes_{\ZZ} R)^{\Ct}\). Now on the level of spectral Mackey functors the tensor product of genuine \(\Ct\)-spectra is determined by Day convolution along the symmetric monoidal structure induced on the level of spans by the cartesian product of finite \(\Ct\)-sets. The Mackey object~\eqrefone{equation:pi-0-norm} inherits from \(\N \GEM R\) an associative algebra structure with respect to this monoidal structure, and this means in particular that \(W \to (R \otimes_{\ZZ} R)^{\Ct}\) is a map of associative algebras and that \((R \otimes_{\ZZ} R)_{\Ct} \to W\) is a map of \(W\)-bimodules (where \(W\) acts on \((R \otimes_{\ZZ} R)_{\Ct}\) via the action of \((R \otimes_{\ZZ} R)^{\Ct}\) on the latter inherited from its equivariant action on \(R \otimes_{\ZZ} R\)).
The ring structure on \(W\) was also calculated in~\cite[Proposition 5.5]{DMPR} and is given by
\[
(a,v)(b,u) = (ab,(a \otimes a)u + v(b \otimes b) + u\rT_{\Ct}(v)) = (ab,(a \otimes a)u + v(b \otimes b) + \rT_{\Ct}(u)v).
\]
Unwinding the definitions, we then see that the notion of a module over \(W\) corresponds exactly to that of an abelian group equipped with a quadratic \(R\)-action, and that the data of a module over~\eqrefone{equation:pi-0-norm}
in \(\Ab\)-valued Mackey functors corresponds exactly to the notion of a form parameter as in Definition~\refone{definition:form-parameter}. In summary, we may reproduce the category of form parameters as the heart under the above t-structure of the \(\infty\)-category \(\Mod_{\N\GEM R}\) of \(\N\GEM R\)-modules in genuine \(\Ct\)-spectra, or, equivalently, of modules with genuine involution over \(\GEM R\). The latter identifies with \(\Funq(\Modp{\GEM R})\) under the equivalence of Theorem~\refone{theorem:classification-genuine-modules}, and in that setting we may identify the heart with the full subcategory spanned by the derived quadratic functors, that is, those quadratic functors \(\Modp{\GEM R} \to \Spa\) which send \(\Prof{R}\) to \(\Ab = \Spa^{\heartsuit}\) (and to which Proposition~\refone{proposition:restriction-of-2-excisive-functors} pertains). In particular, the derived quadratic functor corresponding to a given form parameter \(\fpm = [M_{\Ct} \to Q \to M^{\Ct}]\) is exactly \(\QF^{\gfpm}_M\), and the corresponding \(\GEM\ZZ\)-module with genuine involution \((\GEM M,\cof[\GEM M_\hC \rightarrow \GEM\fpmg], \cof[\GEM M_\hC \rightarrow \GEM\fpmg] \rightarrow \GEM M^{\tC})\) is then obtained by extracting the underlying spectrum and geometric fixed points of the Mackey functor \(\GEM \fpm\) associated to \(\fpm\), where the \(\GEM R\)-action on the geometric fixed points \(\cof[\GEM M_\hC \rightarrow \GEM \fpmg]\) is induced by the action of \(\N \GEM R\) (whose geometric fixed points is \(\GEM R\)) on \(\GEM \fpm\).
\end{remark}

\begin{remark}
\label{remark:comparison}%
In~\cite{comparison} the fourth and ninth authors show that the Grothendieck-Witt groups (the one defined in \S\refone{subsection:GW-group} as well as the higher ones we will define in \papertwo) of derived Poincaré structures as above can be described in terms of the group completion of the corresponding monoid of Poincaré forms on projective modules.
It then follows that for an invertible \(\ZZ\)-module with involution \((M,\sig)\) over \(R\)
the Grothendieck-Witt groups of \((\Dperf(R),\QF^{\gfpm}_M)\) coincide with the Grothendieck-Witt groups of \(\fpm\)-hermitian forms as defined and studied by Bak~\cite{bak-form-parameter}, see also~\cite{karoubi-periodicity} and \cite{SchlichtinghigherI}. In particular, the Grothendieck-Witt groups of \(\Dperf(R)\) with respect to the genuine quadratic, genuine even and genuine symmetric Poincaré structures are the classical \(\sig\)-quadratic, \(\sig\)-even and \(\sig\)-symmetric Grothendieck-Witt groups of \(R\) with coefficients in \(M\).
\end{remark}

We now show how to extend the equivalences of Corollary~\refone{corollary:shift-genuine} to more general form parameters. For this, as explained in~\cite{SchlichtinghigherI}, note that to a form parameter \(M_{\Ct} \to \fpmg \rightarrow M^{\Ct}\) may associate another form parameter
\[
\dfpm = \big[(-M)_{\Ct} \longrightarrow M/\im(\rho) \longrightarrow (-M)^{\Ct}\big].
\]
Indeed, this construction is well-defined since \(\im(\rho)\) is contained in the kernel of the map \(1-\sig\colon M \to M\) and contains the image of the map \(1+\sig\colon M \to M\). Hermitian forms with respect to the form parameter \(\dfpm\) on \(-M\) then correspond to what Bak~\cite{bak-form-parameter} called \emph{\(\fpm\)-quadratic} forms.

\begin{proposition}
\label{proposition:form-parameter-shift}%
Let \((\fpmg,\rho,\tau)\) be a form parameter on \(M\) in which \(\rho\colon \fpmg \to M^{\Ct}\) is injective. Then the loop functor \(\Om\colon \Dperf(R) \to \Dperf(R)\) refines to an equivalence of Poincaré \(\infty\)-categories
\[
\big(\Dperf(R), (\QF_{M}^{\gfpm})\qshift{2}\big) \xrightarrow{\simeq} \big(\Dperf(R), \QF^{\gdfpm}_{-M}\big)
\]
\end{proposition}
\begin{proof}
Applying Corollary~\refone{example:usual-sequence} in the case where \(y=0\) we get for \(P \in \Proj(R)\) an equivalence
\[
\Sig^2\QF^{\gfpm}_M(P[1]) \simeq \Sig\fib\big[\QF^{\gfpm}_M(P[0]) \to \Bil_M(P[0],P[0])\big] = \cof\big[\QF^{\gfpm}_M(P[0]) \to \Bil_M(P[0],P[0])\big].
\]
Now since \(\rho\colon \fpmg \to M^{\Ct}\) is injective the notion of a \(\fpm\)-hermitian form \((\phi,q)\) on \(P\) simply reduces to a symmetric \(M\)-valued form \(\phi\) such that \(\phi(v,v) \in \fpmg\) for \(v \in P\) (from which the data of \(q\) is uniquely determined). It then follows that the map \(\QF^{\gfpm}_M(P[0]) \to \Bil_M(P[0],P[0])\) is injective, and so its cofiber in \(\Spa\) is just the Eilenberg-MacLane spectrum of its cokernel in \(\Ab\). By the the equivalence of Proposition~\refone{proposition:restriction-of-2-excisive-functors} we now conclude that the Poincaré structure \((\QF^{\gfpm}_M \circ \Sig)\qshift{2}\) is the derived Poincaré structure of the degree 2 polynomial functor \(\Proj(R)\op \to \Ab\) given by
\[
\QFtwo(P) = \coker\big[\QF^{\gfpm}_{\proj}(P) \to \Hom_{R \otimes_{\ZZ} R}(P \otimes_{\ZZ} P,M)\big].
\]
We now claim that \(\QFtwo\) is naturally equivalent to \(\QF^{\gdfpm}_{\proj}\). Indeed, unwinding the definitions we see that the transfer map
\[
\Hom_{R \otimes_{\ZZ} R}(P \otimes_{\ZZ} P,M) \to \QF^{\gdfpm}_{\proj}(P)
\]
vanishes on \(\QF^{\gfpm}_{\proj}(P)\) and hence determines a natural transformation
\[
\QFtwo \Rightarrow \QF^{\gdfpm}_{\proj}.
\]
One then readily verifies that this natural transformation is an equivalence on cross effects (which on both sides identify with bilinear forms to \(-M\)). It will hence suffice to check that this natural transformation evaluates to an isomorphism on the projective \(R\)-module \(R\). Indeed, its component at \(R\) is the tautological identification
\[
\QFtwo(R) = M/\fpmg = \QF^{\gdfpm}_{\proj}(R). \qedhere
\]
\end{proof}

\begin{remark}
In \papertwo we prove a version of Karoubi's fundamental theorem in the setting of Poincaré \(\infty\)-categories. When applied to the various derived Poincaré structures on \(\Dperf(R)\) mentioned above (and using the comparison statement of Remark~\refone{remark:comparison}), this yields a complete solution to both conjectures made by Karoubi in~\cite{karoubi-periodicity} concerning potential generalizations of his fundamental theorem to the case where \(2\) is not assumed invertible. More precisely, when combined with Corollary~\refone{corollary:shift-genuine} one obtains the statement of~\cite[Conjecture 1]{karoubi-periodicity}, and when combined with Proposition~\refone{proposition:form-parameter-shift} the statement of~\cite[Conjecture 2]{karoubi-periodicity}. We point out that while these conjectures pertain to the classical theory of forms as encoded by the various polynomial functors above,
the equivalences of Poincaré \(\infty\)-categories appearing in Corollary~\refone{corollary:shift-genuine} and Proposition~\refone{proposition:form-parameter-shift} are only visible after \emph{deriving} these into Poincaré structures on \(\Dperf(R)\).
\end{remark}

We now turn to the proof of Proposition~\refone{proposition:restriction-of-2-excisive-functors}. Let \(\Delta_{\leq n} \subseteq \Delta\) be the full subcategory
spanned by those totally ordered sets with at most \(n+1\) elements.

\begin{definition}
Let \(\C\) and \(\D\) be small \(\infty\)-categories.
\begin{enumerate}
\item
A diagram \(p\colon\Delta\op \to \C\) is called \emph{finite} if it is left Kan extended from its restriction to \(\Delta\op_{\leq n}\) for some \(n\). A colimit over such a diagram is called a \emph{finite geometric realization}. Dually, a diagram \(\Delta \to \C\) is called \emph{finite} if it is right Kan extended from its restriction to \(\Delta_{\leq n}\) for some \(n\). A limit over such a diagram is called a \emph{finite totalization}.
\item
An \(\infty\)-category \(\C\) is said to admit finite geometric realizations if every finite diagram \(\Del\op \to \C\) in \(\C\) admits a colimit. Dually, \(\C\) is said to admit finite totalizations if every finite diagram \(\Del \to \C\) admits a limit.
\item
A functor \(f \colon \C \to \D\) is said to preserve finite geometric realizations if its sends colimit cones \(\ovl{p}\colon(\Del\op)^{\triangleright} \to \C\) over finite diagrams \(p\colon\Del\op \to \C\) to colimit cones in \(\D\). Here, it is not required that \(f \circ p\colon \Del\op \to \D\) remains finite in the above sense. Dually, a functor \(f\) is said to preserve finite totalizations if its sends limit cones \(\ovl{p}\colon(\Del\op)^{\triangleleft} \to \C\) over finite diagrams \(p\colon\Del \to \C\) to limit cones in \(\D\). Again, it is not required that \(f \circ p\colon \Del\op \to \D\) remains finite.
\end{enumerate}
We will denote by \(\Fun^{\Delta\op_\fin}(\C,\D) \subseteq \Fun(\C,\D)\) the full subcategory spanned by the functors which preserve finite geometric realizations.
\end{definition}

\begin{remark}
\label{remark:finite-colimit}%
If an \(\infty\)-category \(\C\) admits finite colimits then it admits finite geometric realizations, since the \(\infty\)-category \(\Del_{\leq n}\op\) is finite and colimits of left Kan extended diagrams \(X \colon \Del\op \to \C\) can be calculated on their restriction to \(\Del_{\leq n}\op\). In addition, if \(\C\) and \(\D\) are \(\infty\)-categories with finite colimits and \(f\colon \C \to \D\) preserves finite colimits then \(f\) sends finite diagrams \(X\colon \Del\op \to \C\) to finite diagrams; indeed, this follows from the pointwise formula for left Kan extensions since for each \([m] \in \Delta\op\) the comma category \((\Delta\op_{\leq n})_{/_{[m]}}\) is finite. It then follows that such an \(f\) also preserves finite geometric realizations. A similar statement holds for finite totalizations under the analogous assumptions concerning finite limits.
\end{remark}

\begin{lemma}
\label{lemma:finite-geometric}%
Let \(\C\) be an \(\infty\)-category which admits finite colimits and \(\D\) an \(\infty\)-category which admits sifted colimits. Then restriction along the Yoneda embedding \(\C \to \Ind(\C)\) induces an equivalence
\[
\Fun^{\sifted}(\Ind(\C),\D) \to \Fun^{\Delta\op_\fin}(\C,\D)
\]
between the full subcategory of functors \(\Ind(\C) \to \D\) which preserve sifted colimits on the left hand side and functors \(\C \to \D\) which preserve finite geometric realizations on the right.
\end{lemma}
\begin{proof}
We first note that by the universal property of ind-categories \cite[Proposition~5.3.5.10]{HTT} we have that \(\Ind(\C)\) admits filtered colimits and restriction along \(\iota\colon\C \hrar \Ind(\C)\) determines an equivalence
\[
\Fun^{\filt}(\Ind(\C),\D) \xrightarrow{\simeq} \Fun(\C,\D),
\]
where the left hand side stands for the full subcategory of \(\Fun(\Ind(\C),\D)\) spanned by the functors which preserves filtered colimits. It will hence suffice to show that if \(f\colon \Ind(\C) \to \D\) is a functor which preserves filtered colimits then \(f\) preserves sifted colimits if and only if \(f \circ \iota\) preserves finite geometric realizations. To begin, note that if \(f\colon \Ind(\C) \to \D\) preserves sifted colimits then it preserves in particular finite geometric realizations. Since the inclusion \(\iota\colon \C \to \Ind(\C)\) preserves finite colimits it preserves finite diagrams and finite geometric realizations by Remark~\refone{remark:finite-colimit}. It then follows that in this case \(f \circ \iota\) preserves finite geometric realizations. To prove the other direction, let us now suppose that \(f\colon \Ind(\C) \to \D\) is a functor which preserves filtered colimits such that \(f\circ \iota\) preserves finite geometric realizations. We wish to show that \(f\) preserves all sifted colimits. By~\cite[Corollary 5.5.8.17]{HTT} it will suffice to show that \(f\) preserves geometric realizations. Let \(t_n\colon \Del\op_{\leq n} \to \Del\op\) be the inclusion. Then every \(\Del\op\)-diagram \(\rho\colon \Del\op \to \Ind(\C)\) can be written as a sequential (and in particular filtered) colimit of finite diagrams of the form
\[
\rho \simeq \colim_n (t_n)_!(t_n)^*\rho .
\]
Since \(f\) preserves filtered colimits it will suffice to prove that \(f\) preserves finite geometric realizations. Now let \(\rho_n\colon \Del\op_{\leq n} \to \Ind(\C)\) be a diagram indexed on \(\Del\op_{\leq n}\). We claim that \(\rho_n\) can be written as a filtered colimit of \(\Del\op_{\leq n}\)-diagrams taking values in \(\C\). To see this, consider the smallest full subcategory \(\E \subseteq \Fun(\Del\op_{\leq n},\Ind(\C))\) which contains the left Kan extended functors \((i_k)_!\iota(\x)\), where \(\x \in \C\) is an object and \(i_k\colon \{[k]\} \hrar \Del\op_{\leq n}\) is the inclusion of the object \([k]\), and closed under finite colimits. Then \(\E\) forms a collection of compact generators for \(\Fun(\Del\op_{\leq n},\Ind(\C)))\) and hence every diagram \(\rho_n\colon \Del\op_{\leq n} \to \Ind(\C)\) is a filtered colimit of diagrams in \(\E\). Since the mapping sets in \(\Del\op_{\leq n}\) are finite the functors \((i_k)_!\iota(\x)\) takes values in the image of \(\C\), and hence factor through \(\Del\op_{\leq n}\)-diagrams in \(\C\). It will hence suffice to show that \(f\) preserves finite geometric realizations of simplicial diagrams which factor through \(\iota \colon \C \hrar \Ind(\C)\). Indeed, since \(\C\) admits finite geometric realizations which are preserved by \(\iota\) this follows from the assumption that \(f \circ \iota\) preserves finite geometric realizations.
\end{proof}

\begin{lemma}
\label{lemma:quadratic-preserves-finite-realizations}%
Any quadratic functor \(\QF\colon \C\op \to \Spa\) preserves finite geometric realizations and finite totalizations.
\end{lemma}

\begin{proof}
It suffices to show the claim for the special cases where \(\QF\) is either exact or of the form \(\x \mapsto \Bil(\x,\x)\) for a bilinear functor \(\Bil \colon \C\op\times \C\op \to \Spa\), since these generate all quadratic functors under both limits and colimits:
indeed, a general quadratic functor \(\QF\) is both the fibre of a map from \(\Bil_{\QF}(-,-)^{\hC}\) to an exact functor and the cofibre of a map from an exact functor to \(\Bil_{\QF}(-,-)_{\hC}\). The case where \(\QF\) is exact follows from Remark~\refone{remark:finite-colimit}.
We hence suppose that \(\QF(\x)= \Bil(\x,\x)\) for some bilinear \(\Bil\), and let \(X\colon \Del\op \to \C\) be a finite diagram. Since \(\Del\) is sifted we then have
\[
\colim_{[n] \in \Del}\QF(X_n) \simeq \colim_{[n] \in \Del}\Bil(X_n,X_n) \simeq \colim_{[n] \in \Del}\colim_{[m] \in \Del}\Bil(X_n,X_m) .
\]
Now for each fixed \(x \in \C\) the functor \(\Bil(\x,-)\) is exact and hence preserves finite geometric realizations by Remark~\refone{remark:finite-colimit}. Similarly, for each fixed \(\y \in \C\) the functor \(\Bil(-,\y)\) preserves finite geometric realizations. We hence get that
\[
\colim_{[n] \in \Del}\colim_{[m] \in \Del}\Bil(X_n,X_m) \simeq \colim_{[n] \in \Del}\Bil(X_n,\colim_{[m]\in\Del} X_m) \simeq \Bil(\colim_{[n]\in\Del} X_n,\colim_{[m]\in\Del} X_m).
\]
Using again that \(\Del\) is sifted we consequently get \(\Bil(\colim_n X_n,\colim_m X_m) \simeq \QF(\colim_n X_n)\), as desired.
\end{proof}

Our strategy for the proof of Proposition~\refone{proposition:restriction-of-2-excisive-functors} is to break the inclusion \(\Prof{R} \hrar \Dperf(R)\) into the composite of two inclusions
\[
\Prof{R} \hrar \Dperf(R)_{\geq 0} \hrar \Dperf(R),
\]
where \(\Dperf(R)_{\geq 0} \subseteq \Dperf(R)\) is
the image under the localisation functor \(\Ch(\Prof{R}) \to \Dperf(R)\) of the full subcategory \(\Ch_{\geq 0}(\Prof{R}) \subseteq \Ch(\Prof{R})\) spanned by the complexes concentrated in non-negative degrees. To make this strategy work we would like to have a notion of a (contravariant) quadratic functor from \(\Dperf(R)_{\geq 0}\) to spectra.

Recall from Definition~\refone{definition:2-excisive} that a functor \(\QF\colon \Dperf(R)_{\geq 0} \to \Spa\op\) is said to be \(2\)-excisive if it sends strongly cocartesian cubes in \(\Dperf(R)_{\geq 0}\) to exact cubes in \(\Spa\op\). We will denote by
\[
\Funexc_*(\Dperf(R)_{\geq 0}) \subseteq \Fun(\Dperf(R)_{\geq 0},\Spa\op)\op
\]
the full subcategory spanned by the reduced \(2\)-excisive functors \(\Dperf(R)_{\geq 0} \to \Spa\op\).

\begin{remark}
\label{remark:cross-effect-colimits}%
It follows from \cite[Proposition~6.1.3.22]{HA} that the cross effect of any reduced \(2\)-excisive functor \(\QF\colon \Dperf(R)_{\geq 0} \to \Spa\op\) is reduced and 1-excisive in each variable separately. %
\end{remark}

\begin{lemma}
\label{lemma:derived-step-1}%
Restriction along the inclusion \(\Dperf(R)_{\geq 0} \subseteq \Dperf(R)\) yields an equivalence
\begin{equation}
\label{equation:restrict-to-connective}%
\Funq(\Dperf(R)) \xrightarrow{\simeq} \Funexc_*(\Dperf(R)_{\geq 0}).
\end{equation}
\end{lemma}
\begin{proof}
Consider the composite
\begin{equation}
\label{equation:composite-kan-approx}%
\Fun_*\big(\Dperf(R)_{\geq 0}\op,\Spa\big) \to \Fun_*\big(\Dperf(R)\op,\Spa\big) \xrightarrow{\App^2} \Funq(\Dperf(R))
\end{equation}
where the first functor is given by right Kan extension along the inclusion \(\iota\colon \Dperf(R)_{\geq 0}\op \hrar \Dperf(R)\op\) (a procedure which preserves reduced functors by the pointwise formula for Kan extensions) and the second by the right adjoint \(\App^2\) to the inclusion \(\Funq(\Dperf(R)) \subseteq \Fun_*(\Dperf(R)\op,\Spa)\) described in Construction~\eqrefone{construction:excisive-approx}.
Since~\eqrefone{equation:composite-kan-approx} is a composite of right adjoints it is itself right adjoint to the composite
\begin{equation}
\label{equation:composite-inclusion-restriction}%
\Funq(\Dperf(R)) \hrar \Fun_*(\Dperf(R)\op,\Spa) \to \Fun_*\big(\Dperf(R)_{\geq 0}\op,\Spa\big).
\end{equation}
Since~\eqrefone{equation:composite-inclusion-restriction} factors through the full subcategory of quadratic functors
the unit and counit of the adjunction between~\eqrefone{equation:composite-kan-approx} and~\eqrefone{equation:composite-inclusion-restriction} also yield an adjunction
\begin{equation}
\label{equation:adjunction-kan-restriction}%
\Funexc_*(\Dperf(R)_{\geq 0}) \adj  \Funq(\Dperf(R))
\end{equation}
where the right adjoint is obtained by restricting the domain of~\eqrefone{equation:composite-kan-approx} and the left adjoint is the restriction functor~\eqrefone{equation:restrict-to-connective} under consideration. We claim that the adjunction~\eqrefone{equation:adjunction-kan-restriction} is an equivalence.
To begin, we first show that for \(\QF \in \Funexc_*(\Dperf(R)_{\geq 0})\)
the counit \(\iota^*\App^2(\iota_*\QF) \to \QF\) is an equivalence. We note that this counit is given itself by a composite
\begin{equation}
\label{equation:unit-composite}%
\iota^*\App^2(\iota_*\QF) \to \iota^*\iota_*\QF \to \QF
\end{equation}
where the second map is the counit of \(\iota^* \dashv \iota_*\), which is an equivalence since \(\iota\) is fully-faithful. The first map in turn is the one obtained by applying \(\iota^*\) to the counit \(\App^2(\iota_*\QF) \to \iota_*\QF\). We hence need to show that the component
\begin{equation}
\label{equation:counit-connective}%
\App^2(\iota_*\QF)(X) \to \iota_*\QF(X)
\end{equation}
is an equivalence for \(X\) in the image of the inclusion \( \Dperf(R)_{\geq 0} \subseteq \Dperf(R)\). For this we recall that \(\App^2\) is defined via a sequential limit
\[
\App^2(\RF) := \lim[\cdots \to \T^2\T^2(\RF) \to \T^2(\RF) \to \RF]
\]
where \(\T^2(\RF)(X) = \Sig\cof[\Bil_{\RF}(\Sig X,\Sig X) \to \RF(\Sig X)]\), see Construction~\refone{construction:excisive-approx}.
It will then suffice to show that for \(X \in \Dperf(R)_{\geq 0}\) the sequence
\[
\Sig\Bil_{\QF}(\Sig X,\Sig X) \to \Sig\QF(\Sig X) \to \QF(X)
\]
is exact. Indeed, this follows by the exact same argument as the dual version of Lemma~\refone{lemma:goodwillie} (see Remark~\refone{remark:goodwillie-dual}) in the case where \(z\) and \(w\) are zero objects, using that \(\QF\) is assumed to be reduced and \(2\)-excisive on \(\Dperf(R)_{\geq 0}\).

As the counit is an equivalence, to finish the proof it will suffice to show that \(\iota^*\) is conservative, or equivalently (since its domain is stable), detects zero objects. In particular, we need to show that if \(\QF\colon \Dperf(R)\op \to \Spa\) is a quadratic functor which vanishes on \(\Dperf(R)_{\geq 0}\) then \(\QF\) is the zero functor.
Indeed, for such a \(\QF\) we will have that \((\T^2)^{\circ n}(\QF)\) vanishes on \((\Dperf(R))_{\geq -n} := \Sig^{-n}\Dperf(R)_{\geq 0}\) and hence that
\[
\QF \simeq \App^2(\QF) = \lim[\cdots \to \T^2\T^2(\QF) \to \T^2(\QF) \to \QF] = 0,
\]
since any \(X \in \Dperf(R)\) lies in \(\Dperf(R)_{\geq -n}\) for some \(n\).
\end{proof}

\begin{lemma}
\label{lemma:derived-step-2}%
Restriction along the inclusion \(\Prof{R} \subseteq \Dperf(R)_{\geq 0}\) yields an equivalence
\[
\Funexc_*(\Dperf(R)_{\geq 0}) \xrightarrow{\simeq} \Funpoly(\Prof{R}\op,\Spa)
\]
\end{lemma}
\begin{proof}
Let \(\Der(R)_{\geq 0} \subseteq \Der(R)\) denote the full subcategory spanned by the objects represented by non-negatively graded complexes.
By~\cite[Proposition 1.3.3.14]{HA} restriction along the inclusion \(\Prof{R} \to \Der(R)_{\geq 0}\) yields an equivalence
\[
\Fun^{\sifted}(\Der(R)_{\geq 0}\op,\Spa) \xrightarrow{\simeq} \Fun(\Prof{R}\op,\Spa),
\]
where the left hand side denotes the full subcategory of \(\Fun(\Der(R)_{\geq 0}\op,\Spa)\) spanned by the functors which preserve sifted limits. It then follows from Lemma~\refone{lemma:finite-geometric} that restriction along the inclusion \(\Prof{R} \to \Dperf(R)_{\geq 0}\) induces an equivalence
\begin{equation}
\label{equation:free-finite-geometric}%
\Fun^{\Delta_\fin}(\Dperf(R)_{\geq 0}\op,\Spa) \xrightarrow{\simeq} \Fun(\Prof{R}\op,\Spa)
\end{equation}
where on the left hand side we have the \(\infty\)-category of functors \(\Dperf(R)_{\geq 0}\op \to \Spa\) which preserve finite geometric realizations. By Lemma~\refone{lemma:quadratic-preserves-finite-realizations} the latter contains \(\Funexc_*(\Dperf(R)_{\geq 0})\) as a full subcategory. We now claim that under the equivalence~\eqrefone{equation:free-finite-geometric} the reduced \(2\)-excisive functors on the left hand side correspond to the functors that are reduced and polynomial of degree \(2\) on the right hand side.
First, since the inclusion \(\Prof{R} \to \Dperf(R)_{\geq 0}\) is additive and reduced \(2\)-excisive functors have cross effects which are reduced and 1-excisive in each variable separately (see \cite[Proposition~6.1.3.22]{HA})
it follows that reduced \(2\)-excisive functors restrict to reduced functors which are polynomial of degree \(2\). Conversely, suppose that \(\RF\colon \Dperf(R)_{\geq 0}\op \to \Spa\) is a functor which preserves finite geometric realizations whose restriction to \(\Prof{R}\op\) is reduced and polynomial of degree \(2\). Since \(\Prof{R}\) contains the zero object of \(\Dperf(R)_{\geq 0}\) we have that \(\RF\) is reduced. To show that \(\RF\) is \(2\)-excisive we need to show that it sends strongly cocartesian cubes in \(\Dperf(R)_{\geq 0}\) to exact cubes of spectra. For this, let \(\rho\colon (\Del^1)^3 \to \Dperf(R)_{\geq 0}\) be a strongly cocartesian cube. We want to reduce to the case where \(\rho\) takes values in the subcategory of finite projective modules and injections. Let us identify \((\Del^1)^3\) with the nerve of the poset \(\cP([2])\) of subsets of \([2] = \{0,1,2\}\). Let \(\cP_{\leq 1}([2]) \subseteq \cP([2])\) be the full subposet spanned by \(\{\},\{0\},\{1\},\{2\}\).
Since \(\Dperf(R)_{\geq 0}\) is the localisation of the cofibration category \(\Ch_{\geq 0}(\Prof{R})\) it follows from~\cite[Theorem 7.6.17]{cisinski} that we can represent \(\rho|_{\cP_{\leq 1}([2])}\) by a diagram \(\tau\colon \cP_{\leq 1}([2]) \to \Ch_{\geq 0}(\Prof{R})\) in which each map \(c_i\colon \tau(\{\}) \to \tau(\{i\})\) is levelwise injective with projective kernel. The Dold-Kan correspondence then associates to \(\tau\) a diagram \(\tau'\colon \cP_{\leq 1}([2]) \to \Prof{R}^{\Del\op}\) of simplicial \(R\)-modules such that each of the maps \(c_i\colon \tau(\{\}) \to \tau(\{i\})\) is levelwise injective with projective cokernel (see~\cite[\S II.4.12] {quillen2006homotopical}). In addition, since \(\tau\) takes values in bounded complexes we can find an \(n \geq 0\) such that it takes values in complexes concentrated in degrees \(0\) to \(n\). Under the Dold-Kan correspondence such complexes map to simplicial \(R\)-modules which are left Kan extended from their restriction to \(\Del_{\leq n}\op\). Switching the simplicial dimension with the \(\cP_{\leq 1}([2])\)-dimension we may conclude that \(\tau|_{\cP_{\leq 1}([2])}\) can be written as a finite geometric realization of a simplicial family of diagrams \(\tau'_{n}\colon \cP_{\leq 1}([2]) \to \Prof{R}\) such that each \(\tau'_n\) has the property that \(\tau'_n(\{\}) \to \tau'_n(\{i\})\) is a split injective map of projective modules. Left Kan extending from \(\cP_{\leq 1}([2])\) we obtain a representation of \(\rho\) as a finite geometric realization of a simplicial family of strongly cocartesian cubes \(\rho_n\colon \cP([2]) \to \Prof{R}\) such that each \(\rho_n(\{\}) \to \rho_n(\{i\})\) is a split injective map of projective modules. Since \(\RF\) preserves finite geometric realizations it will suffice to show that \(\RF\) sends each \(\rho_n\) to an exact cube of spectra. Now since \(\Prof{R}\) is closed in \(\Dperf(R)_{\geq 0}\) under pushouts with one leg split injective it follows that each \(\rho_n\) is a strongly cocartesian cubes in \(\Prof{R}\). More explicitly, we may pick projective modules \(X_n,Y_n^0,Y_n^1,Y_n^2\) such that \(\tau'_n(\{\}) = X^n\) and \(\tau'_n(\{i\}) = X_n \oplus Y_n^i\) for \(i=0,1,2\), in which case \(\rho_n\colon \cP([2]) \to \Dperf(R)\) is given by
\[
[2] \supseteq S \mapsto X_n \oplus \left[\displaystyle\mathop{\oplus}_{i \in S} Y^i_n\right].
\]
To finish the proof we need to show that the resulting cube of spectra
\begin{equation}
\label{equation:cube-to-prove-exact}%
S \mapsto \RF(X_n \oplus \big[\displaystyle\mathop{\oplus}_{i \in S} Y^i_n\big])
\end{equation}
is exact.
Let \(\Bil(-,-)\) be the cross-effect of \(\RF\). The assumption that \(\RF|_{\Prof{R}}\) is polynomial of degree \(2\) means that the restriction of \(\Bil\) to \(\Prof{R} \times \Prof{R}\) preserve direct sums in each variable separately. The cube of spectra \(\RF \circ \rho_n\) thus decomposes as a direct sum of four components
\[
S \mapsto \RF(X_n) \oplus \left[\oplus_{i \in S}\RF(Y_n^i)\right] \oplus \left[\oplus_{i \in S}\Bil(X_n,Y_n^i)\right] \oplus \left[\oplus_{i < j \in S} \Bil(Y_n^i,Y_n^j)\right] .
\]
The first component is constant and is hence an exact cube, and the second and third components are visibly exact.
Finally, the fourth component \(S \mapsto \oplus_{i < j \in S} \Bil(Y_n^i,Y_n^j)\) is also an exact cube since it vanishes on \(S\) of size smaller than \(2\) and exhibits its value at \(S=\{0,1,2\}\) as the product of its values on \(\{0,1\},\{1,2\}\) and \(\{0,2\}\). We have thus showed that the cube of spectra~\eqrefone{equation:cube-to-prove-exact} is exact, and so the proof is complete.
\end{proof}

\begin{proof}[Proof of Proposition~\refone{proposition:restriction-of-2-excisive-functors}]
Combine Lemma~\refone{lemma:derived-step-1} and Lemma~\refone{lemma:derived-step-2}.
\end{proof}

\subsection{Group-rings and visible structures}
\label{subsection:visible}%

In this section we will consider group rings, and more generally, group ring spectra, which provide important examples of rings with anti-involutions, and whose various types of \(\L\)-groups play a role in surgery theory and manifold classification, see Example~\refone{example:visible-symmetric}. As in the work of Weiss-Williams~\cite{WWIII}, these geometric applications are sometimes better serviced by replacing modules over group rings by parametrised spectra, a point of view we will take up in \S\refone{subsection:parametrised-spectra}.

Recall that for an ordinary group \(G\) equipped with a homomorphism \(\chi\colon G \to \{-1,1\}\), one may endow the group ring \(\ZZ G\) with the associated \(\chi\)-twisted involution \(\tau_{\chi}\colon\ZZ G \xrightarrow{\cong} (\ZZ G)\op\), defined by sending \(g \in \ZZ G\) to \(\chi(g)g^{-1}\). In what follows we will consider a generalization of this setup where \(G\)
is replaced by an \defi{\(\Eone\)-group}, that is, a group-like \(\Eone\)-monoid in spaces, and instead of \(\ZZ\) the coefficients are taken in the sphere spectrum.
More precisely, for such a \(G\) the suspension spectrum \(\SS[G] := \Sig^{\infty}_+G\) inherits the structure of an \(\Eone\)-ring spectrum, which is called the \defi{group ring spectrum} of \(G\). This \(\Eone\)-ring spectrum is characterized by the fact that its module spectra correspond to spectra with a \(G\)-action, a notion which we will call here (naive) \defi{\(G\)-spectra}, and consequently also use the notations \(\Mod_G := \Mod_{\SS[G]}\) and \(\Modp{G} := \Modp{\SS[G]}\), as well as the shorthand \(\otimes_G\) for \(\otimes_{\SS[G]}\).
In this section we consider data giving rise to Poincaré structures on \(\Modp{\SS[G]}\).
We also consider variants for the group algebras \(A[G] := A \otimes_{\SS} \SS[G]\) with coefficients in some \(\Eone\)-ring spectrum \(A\). A closely related construction involving parametrised spectra over a space will be explored in \S\refone{subsection:parametrised-spectra}.

Recall that a spectrum \(E \in \Spa\) is called \defi{invertible} if there exists a spectrum \(E'\) such that \(E \otimes_{\SS} E' \simeq \SS\). In this case \(E\) is necessarily of the form \(\Sig^n\SS\) for a unique \(n \in \ZZ\), to which we refer as the \emph{rank} of \(E\).
By a \defi{character} of an \(\Eone\)-group \(G\) we will mean a pair \((E,\chi)\) where \(E\) is an invertible spectrum and \(\chi\colon G \to \Aut(E) \simeq \gl_1(\SS)\) is a homomorphism of \(\Eone\)-groups, encoding an action of \(G\) on \(E\).

\begin{construction}
\label{construction:visible}%
Let \(G\) be an \defi{\(\Eone\)-group} equipped with a character \((E,\chi)\). We define a module with genuine involution on the group ring \(\SS[G]\) as follows.
First, since the diagonal map \(G \to G \times G\) is naturally equivariant with respect to the flip \(\Ct\)-action on \(G \times G\) and the trivial action on \(G\), it follows that the induction functor
\[
\Mod_G \to \Mod_{G \times G}
\]
refines to a functor
\[
\Mod_G \to [\Mod_{G \times G}]^{\hC},
\]
where \(\Ct\) acts on \(\Mod_{G \times G}\) via its flip action on \(G \times G\). In particular, the \(G \times G\)-spectrum
\[
E[G^{-}] :=  \SS[G \times G] \otimes_G E = \SS[G \times G] \otimes_{\SS[G]} E
\]
induced from the \(G\)-spectrum \(E\) is naturally a \(\Ct\)-equivariant object of \(\Mod_{G \times G}\). Identifying \(\SS[G] \otimes_{\SS} \SS[G]\) with \(\SS[G \times G]\) we then see that the Tate diagonal
\[
\Del_G\colon \SS[G] \to \SS[G \times G]^{\tC}
\]
is given by the composite \(\SS[G] \to \SS[G \times G]^{\hC} \to \SS[G \times G]^{\tC}\) in which the first map is induced by the \(\Ct\)-equivariant diagonal \(G \to G \times G\). In particular, the \(G\)-action on \(E[G^{-}]^{\tC}\) is induced by the \(G\)-action on \(E[G^{-}]\) restricted from \(G \times G\) along the diagonal, and we may then consider the composite map of \(G\)-spectra
\[
\eta\colon E \simeq \SS[G] \otimes_G E \xrightarrow{\Del_G \otimes_G E} (\SS[G \times G])^{\tC} \otimes_G E \to (\SS[G \times G] \otimes_G E)^{\tC} = E[G^{-}]^{\tC},
\]
and we call the resulting module with genuine involution \((E[G^{-}],E,\eta)\) the \defi{\(\chi\)-twisted visible} module.
We will denote the corresponding hermitian structure on \(\Modp{G}\) by \(\QF^{\vis}_{\chi}\) and refer to it as the \defi{\(\chi\)-twisted visible structure} on \(\Modp{G}\). Similarly, for a given subgroup \(c \subseteq \K_0(\Modp{G})\) we may restrict \(\QF^{\vis}_{\chi}\) to obtain the corresponding visible hermitian structure on \(\Mod^c_G\).
\end{construction}

\begin{lemma}
\label{lemma:visible-invertible}%
Let \(n\) be such that \(E \simeq \Sig^n\SS\) as spectra. Then the \((\SS[G] \otimes_{\SS} \SS[G])\)-module \(E[G^{-}]\) is invertible and equivalent to \(\Sig^n\SS[G]\) as an \(\SS[G]\)-module. Furthermore, under this equivalence the involution on \(E[G^{-}]\) corresponds to the \(n\)-fold suspension of a ring spectrum anti-involution
\[
\tau_{\chi}\colon \SS[G] \to \SS[G]{\op},
\]
so that \((\Sig^{-n}E[G^{-}],\Sig^{-n}E,\Sig^{-n}\eta)\) promotes \(\SS[G]\) to a genuine ring with involution. In particular, the visible structure \(\QF^{\vis}_{\chi}\) on both \(\Modp{G}\) and \(\Modf{G}\) is Poincaré (see Example~\refone{example:anti-involution}).
\end{lemma}

We will use the notation \(\SS_{\chi}[G]\) to denote the group ring of \(G\) considered with the anti-involution \(\tau_{\chi}\colon \SS[G] \to \SS[G]\op\) of Lemma~\refone{lemma:visible-invertible}. Explicitly, \(\tau_{\chi}\) is the ring map induced by the composed map
\[
G \to \gl_1(\SS) \times G\op \to \gl_1(\SS[G]\op),
\]
in which the first map is given by \(g \mapsto (\chi(g),g^{-1})\).
In particular, on \(\pi_0\SS[G] \cong \ZZ[\pi_0 G]\) this involution can be written as \(\sum_i a_i g_i \mapsto \sum_i a_i\chi_0(g_i)g_i^{-1}\), where
\[
\chi_0\colon \pi_0 G \to \pi_0\gl_1(\SS) = \{1,-1\}
\]
is the induced homomorphism on \(\pi_0\).

\begin{proof}[Proof of Lemma~\refone{lemma:visible-invertible}]
By Proposition~\refone{proposition:recognition} it will suffice to exhibit a \(\Ct\)-equivariant map \(u\colon \SS \to \Sig^{-n}E[G^{-}]\) which freely generates \(\Sig^{-n}E[G^{-}]\) as an \(\SS[G]\)-module. Indeed, this is given by the \(\Ct\)-equivariant map
\[
\SS \simeq \Sig^{-n}E\otimes_G \SS[G] \to \Sig^{-n}E \otimes_G\SS[G \times G] = \Sig^{-n}E[G^{-}],
\]
induced by the diagonal map \(G \to G \times G\). This map exhibits \(\Sig^{-n}E[G^{-}]\) as a free \(\SS[G]\)-module on a single generator since the shear map \(G \times G \to G \times G\) (given informally by \((g,h) \mapsto (g,gh)\)) is an equivalence by the group-like property.
\end{proof}

The following variant is also of interest:
\begin{variant}[Visible structures with coefficients]
\label{variant:visible-with-coefficients}%
Let \(G\) be an \(\Eone\)-group equipped with a character \((E,\chi)\) and let \((A,N,\alp)\) be an \(\Eone\)-ring spectrum with genuine involution.
The the \(\Eone\)-ring
\[
A[G] := A \otimes_{\SS} \SS[G]
\]
can be endowed with a module with genuine involution by tensoring the visible module \((E[G^{-}], E,\eta)\) of Construction~\refone{construction:visible-param-spectra} with \((A,N,\alp)\), using the monoidal structure on genuine \(\Ct\)-spectra (see Remark~\refone{remark:NormModules}). Explicitly, this yields the module with genuine involution
\[
(A_{\chi}[G^-],N \otimes_{\SS} E,\alp \otimes_{\SS} \eta) = (A \otimes_{\SS} E[G^{-}],N \otimes_{\SS} E,\alp \otimes_{\SS} \eta)
\]
where
\[
\alp \otimes_{\SS} \eta \colon N \otimes_{\SS} E \to A^{\tC} \otimes_{\SS} E[G^{-}]^{\tC} \to (A \otimes_{\SS} E[G^{-}])^{\tC}
\]
is obtained via the lax monoidal structure map of the Tate construction. We will denote the corresponding hermitian structure on \(\Modp{A[G]}\) by \(\QF^{\alp}_{A,\chi}\), and refer to it as the \defi{\(\chi\)-twisted visible structure} on \(\Modp{A[G]}\). It then follows from Lemma~\refone{lemma:visible-invertible} that this hermitian structure is Poincaré and identifies up to a shift with the Poincaré structure associated to a genuine refinement of the anti-involution \(\tau_{A,\chi}\colon A[G] \xrightarrow{\simeq} (A[G]){\op}\) induced by the anti-involution of \(A\) and \(\chi\)-twisted anti-involution of \(\SS[G]\). We will refer to \(\tau_{A,\chi}\) as the \emph{\(\chi\)-twisted anti-involution} on \(A[G]\).
\end{variant}

\begin{example}
When \(G\) is discrete a common choice of a twisting is via a sign homomorphism \(\chi\colon G \to \{-1,1\}\). This can be made into a \(G\)-character \((E,\chi)\) once we fix what is meant by the sign action on \(\SS\). Here (at least) two equally natural options are possible: one can either take \(E = \SS^{\sig-1}\) to be the desuspension of the sign representation sphere, or take \(E = \SS^{1-\sig}\) to be its inverse. These are generally not equivalent as spectra with \(\Ct\)-action, see Example~\refone{example:signperiodicmodules} \refoneitem{item:warning-many-signs}.
\end{example}

\begin{example}
\label{example:visible-quadratic}%
Let \(G\) be an \(\Eone\)-group equipped with a character \((E,\chi)\).
If \((A,0,0 \to A^{\tC})\) is an \(\Eone\)-ring spectrum with genuine involution associated to the quadratic genuine refinement of an anti-involution on \(A\), then the associated \(\chi\)-twisted visible Poincaré structure on \(\Modp{A}\) coincides with the quadratic Poincaré structure associated to the \(\chi\)-twisted duality on \(\Modp{A[G]}\).
\end{example}

\begin{example}[The visible symmetric structure]
\label{example:visible-symmetric}%
Let \(G\) be an \(\Eone\)-group equipped with a character \((E,\chi)\) and let \(A\) be an \(\Eone\)-ring spectrum with anti-involution. We refer to the \(\chi\)-twisted Poincaré structure on \(\Modp{A[G]}\) associated to the \(\Eone\)-ring spectrum with genuine involution \((A,A^{\tC},\id\colon A^{\tC} \to A^{\tC})\) as the \defi{visible symmetric} Poincaré structure, and denote it by \(\QF^{\vis-\sym}_{A,\chi} \in \Funq(\Modp{A[G]})\).
Its reference map \(A^{\tC} \otimes_{\SS} E \to A_{\chi}[G^-]^{\tC}\) can then be factored as a composite
\begin{equation}
\begin{split}
\label{align:visible-symmetric}%
A^{\tC} \otimes_{\SS} E & = A^{\tC} \otimes_{\SS} \SS[G] \otimes_{G} E \\
 & \xrightarrow{\alp} A[G \times G]^{\tC} \otimes_G E \\
 & \xrightarrow{\bet} \big(A[G \times G]\otimes_G E\big)^{\tC} = A_{\chi}[G^-]^{\tC} \ ,
\end{split}
\end{equation}
where \(\alp\) is induced by the Tate diagonal of \(\SS[G]\) and the lax symmetric monoidal structure of \((-)^{\tC}\), and \(\bet\) is the canonical interchange map between \((-) \otimes_G E\) and \((-)^{\tC}\).
When \(G\) and \(A\) are discrete the \(\L\)-groups of this Poincaré structure recover the visible symmetric \(\L\)-groups \(A[G]\), first defined by Weiss~\cite{Weiss}. Indeed, translating the definition of~\cite{Weiss} to the present setting gives a Poincaré structure with linear part classified by the \(A[G]\)-module \(A[G \times G]^{\tC} \otimes_G E\) with reference map \(\bet\) as above. To compare the two it will hence suffice to show that \(\alp\) is an equivalence when \(G\) and \(A\) are discrete. Indeed, we can factor \(\alp\) as a composite
\[
A^{\tC} \otimes_{\SS} \SS[G] \otimes_{G} E \to A[G]^{\tC} \otimes_{G} E \to A[G \times G]^{\tC} \otimes_G E
\]
where the first map is an equivalence in this case
since the Tate construction preserves filtered colimits of \emph{discrete} spectra and the second map is an equivalence since the cofibre of the map \(\SS[G] \to \SS[G \times G]\) is an induced \(\Ct\)-module when \(G\) is discrete.
\end{example}

\begin{example}[The visible genuine structure]
Let \(G\) be an \(\Eone\)-group equipped with a character \((E,\chi)\) and let \(A\) be a connective \(\Eone\)-ring spectrum with anti-involution. We will refer to the \(\chi\)-twisted Poincaré structure on \(\Modp{A[G]}\) associated to the \(\Eone\)-ring spectrum with genuine involution \((A,\tau_{\geq m}A^{\tC},t_m)\) as the \defi{visible genuine} Poincaré structure, and denote it by \(\QF^{\vis-\geq m}_{A,\chi} \in \Funq(\Modp{A[G]})\). Its linear part is then classified by the \(A[G]\)-module \(\tau_{\geq m}A^{\tC} \otimes_{\SS} E \simeq \tau_{\geq n+m}(A^{\tC} \otimes_{\SS} E)\), where \(n\) is the rank of \(E\).
\end{example}

In some circumstances, the visible symmetric Poincaré structure on \(A[G]\) identifies with the symmetric one:

\begin{lemma}
\label{lemma:visible-is-symmetric}%
Let \(A\) be an \(\Eone\)-ring spectrum with anti-involution, \(G\) an \(\Eone\)-group and \((E,\chi)\) a \(G\)-character of rank \(n\). If either \(\Bs G\) is a finitely dominated space or \(G\) is a discrete group with no non-trivial 2-torsion then the comparison map \(\QF^{\vis-\sym}_{A,\chi} \Rightarrow \QF^{\sym}_{A_{\chi}[G^-]}\) is an equivalence.
\end{lemma}

\begin{remarks}\
\begin{enumerate}
\item
The two different conditions considered in Lemma~\refone{lemma:visible-is-symmetric} are not completely unrelated: if \(G\) is a discrete group then the condition that \(\Bs G\) is finitely dominated implies that \(G\) has no non-trivial torsion (indeed, for such a \(G\) the constant module \(\ZZ\) admits a finite projective resolution, and so all subgroups of \(G\) have finite cohomological dimension, see, e.g.,~\cite[Corollary VIII.2.5]{Brown-cohomology}).
\item
The proof of Lemma~\refone{lemma:visible-is-symmetric} in the case of \(G\) discrete with no non-trivial 2-torsion extends to more general groups: the same argument works for example for Lie groups with no non-trivial 2-torsion which are sufficiently nice in the sense \(G\) is a finite dimensional \(\Ct\)-CW-complex considered with respect to the \(\Ct\)-action \(g \mapsto g^{-1}\).
\end{enumerate}
\end{remarks}

In the circumstances of Lemma~\refone{lemma:visible-invertible}, the visible genuine Poincaré structure on \(A[G]\) also identifies with the truncated structure of Example~\refone{example:truncation}, up to a shift in the truncation point by the rank of \(E\):

\begin{corollary}
Let \(A\) be a connective \(\Eone\)-ring spectrum with anti-involution, \(G\) an \(\Eone\)-group and \((E,\chi)\) a \(G\)-character of rank \(n\). If either \(\Bs G\) is a finitely dominated space or \(G\) is a discrete group with no 2-torsion then
\[
\QF^{\vis-\geq m}_{A,\chi} \simeq \QF^{\geq m+n}_{A_{\chi}[G^-]}.
\]
\end{corollary}

\begin{proof}[Proof of Lemma~\refone{lemma:visible-is-symmetric}]
The comparison claim is equivalent to the statement that under the given assumptions the reference map \(A^{\tC} \otimes_{\SS} E \to A_{\chi}[G^{-}]^{\tC}\) of the visible symmetric structure is an equivalence.

Assume first that \(G\) is a discrete group with no \(2\)-torsion.
Then we may identify \(A_{\chi}[G^{-}] \simeq \oplus_{g \in G} [A \otimes_{\SS} E]\) with the involution acting by applying the involution of \(A\) on each factor and then switching the component associated to \(g \in G\) with the component associated to \(g^{-1}\). Under the identification \(A^{\tC} \otimes_{\SS} E \simeq (A \otimes_{\SS} E)^{\tC}\) the reference map is then induced on Tate objects by the inclusion of the component of \(1 \in G\). When \(G\) has no 2-torsion the cofibre of this inclusion is an induced \(\Ct\)-module, and so the induced map on Tate objects is an equivalence.

Now assume that \(\Bs G\) is a compact space. We then factor the reference map of the visible symmetric structure as in~\eqrefone{align:visible-symmetric}. Since \(\Bs G\) is compact it follows that \(E\), which can be written as a \(\Bs G\)-indexed colimit of free \(\SS[G]\)-modules, is itself a compact \(\SS[G]\)-module, and so the arrow \(\bet\) in~\eqrefone{align:visible-symmetric} is an equivalence.
Now consider the exact functor \(T\colon \Mod_A \to \Spa\) given by \(T(X) \mapsto (\ovl{X} \otimes_A X)^{\tC}\), where \(\ovl{X}\) denotes \(X\), but considered as a right \(A\)-module via the involution on \(A\). In particular, we have a canonical equivalence \(T(A) \simeq A^{\tC}\). Then the interchange map \(T(A) \otimes_{\SS} E \to T(A \otimes_{\SS} E)\) is an equivalence since \(E\) is finite as a spectrum, but can also be factored as a composite of interchange maps
\[
T(A) \otimes_{\SS} E = T(A) \otimes_{\SS} \SS[G] \otimes_G E \to T(A[G])\otimes_G E \to T(A[G] \otimes_G E) = T(A \otimes_{\SS} E)
\]
in which the second interchange map is an equivalence since \(E\) is perfect as an \(\SS[G]\)-module. It then follows that the first interchange map
\[
A^{\tC} \otimes_{\SS} \SS[G] \otimes_G E = T(A) \otimes_{\SS} \SS[G] \otimes_G E \to T(A[G])\otimes_G E = A[G\times G]^{\tC} \otimes_{G} E
\]
is an equivalence as well. Unwinding the definitions, this map can be identified with \(\alp\), and so the desired result is obtained.
\end{proof}

\begin{remark}
\label{remark:universal-property}%
Let \(G\) be an \(\Eone\)-group equipped with a character \((E,\chi)\).  Consider the universal Poincaré structure \(\QF^{\uni}\) on \(\Spaf\) associated to the module with genuine involution \((\SS,\SS,\SS \to \SS^{\tC})\). Since any stable \(\infty\)-category is tensored over \(\Spa\) the character \(\chi\) induces an action of \(G\) on \(E \otimes \QF^{\uni} \in \Funq(\Spaf)\), and hence on \((\Spaf,E\otimes\QF^{\uni}) \in \Catp\).
Using the point of view of parametrised spectra, we will show in \S\refone{subsection:limits}
that the visible Poincaré structure \((\Modf{G},\QF^{\vis}_{\chi})\) on finitely presented \(\SS[G]\)-modules can be universally characterized as the \emph{quotient} of \((\Spaf,E\otimes\QF^{\uni})\) by \(G\) in \(\Catp\) (see Proposition~\refone{proposition:visible-is-colimit-full} and Example~\refone{example:parametrised-is-group-ring}). The Poincaré \(\infty\)-category \((\Modp{G},\QF^{\vis}_{\chi})\) can then be identified with the corresponding quotient only computed in idempotent complete Poincaré \(\infty\)-categories.

More generally, if \((A,N,\alp)\) is an \(\Eone\)-ring spectrum with genuine anti-involution then \((\Modf{A[G]},\QF^{\alp}_{A,\chi})\) represents the quotient of the associated \(G\)-action on \((\Modf{A},E \otimes \QF^{\alp}_A)\); this can be deduced from the above claim by identifying \((\Modf{A[G]},\QF^{\alp}_{A,\chi})\) with a suitable tensor product of \((\Modf{A},\QF^{\alp}_A)\) and \((\Modf{G},\QF^{\vis}_{\chi})\), see \S\refone{subsection:monoidal-structure} and Example~\refone{example:tensor-of-rings}.
\end{remark}

\begin{variant}
\label{variant:group-with-genuine-involution}%
Construction~\refone{construction:visible} can be generalized by introducing the following additional pieces of data: a \(\Ct\)-action \(\tau\colon \Ct \to \Aut(G)\) on \(G\), an additional \(\Eone\)-group \(H\), and a \(\Ct\)-equivariant map \(H \to G\), where \(\Ct\) acts trivially on \(H\). We may consider such a structure as a \defi{genuine \(\Ct\)-action} on \(G\). The data of a character \((E,\chi)\) for \(G\) then needs to be promoted to that of \(\Ct\)-equivariant character \(\chi\colon G \to \Aut(E)\) (where \(\Ct\)-acts trivially on the target), and the induced \(\Ct\)-action on the restricted character \(\chi|_{H}\) should be equipped with a trivialization. To all this data one can associate a module with genuine involution \((E[G^{-\tau}],E[G/H],\eta)\) where \(E[G^{-\tau}] := E[G\times G]\otimes_G E\) is the \((G \times G)\)-spectrum induced from \(E\), this time along the \(\tau\)-diagonal \((\tau,\id)\colon G \to G \times G\), and \(E[G/H]\) is the \(G\)-module induced from the \(H\)-character \((E,\chi|_{H})\) along the map \(H \to G\). Since the map \((\tau,\id)\) is \(\Ct\)-equivariant with respect to flip action on \(G \times G\), and \(\chi\) is \(\Ct\)-equivariant, the \((G \times G)\)-module \(E[G^{-\tau}]\) inherits an involution compatible with the flip involution of \(G \times G\), and hence the structure of a module with involution over \(\SS[G]\). At the same time, since the usual diagonal \(G \to G \times G\) is \(\Ct\)-equivariant with respect to the trivial action on the domain, and coincides with the \(\tau\)-diagonal when restricted to \(H\), it induces a \(\Ct\)-equivariant map \(\SS[G] \otimes_H E \to  \SS[G \times G]\otimes_{G} E\) with respect to the trivial action on the domain. The composed map
\[
\eta\colon E[G/H] = \SS[G] \otimes_H E \to (\SS[G \times G] \otimes_{G} E)^{\hC} \to (\SS[G \times G] \otimes_{G} E)^{\tC} = E[G^{-\tau}]^{\tC}
\]
then constitutes a structure map exhibiting \((E[G^{\tau}],E[G/H],\eta)\) as a module with genuine involution over \(\SS[G]\), yielding a hermitian structure \(\QF^{\vis}_{\chi,\tau}\) on \(\Modp{G}\), or on \(\Mod^c_G\) for some subgroup \(c \subseteq \K_0(\Modp{G})\).
The case of Construction~\refone{construction:visible} can be recovered as corresponding to the trivial involution on \(G\) with \(H=G\) and \(H \to G\) the identity map. Arguing as in the proof of Lemma~\refone{lemma:visible-invertible} we can check that this module with genuine involution is invertible and corresponds, up to a shift, to a genuine refinement of a suitable anti-involution on \(\SS[G]\), whose underlying equivalence \(\SS[G] \to \SS[G]\op\) is induced by the map \(G \to \gl_1(\SS)\times G\) given by \(g \mapsto (\chi(g),\tau(g)^{-1})\). As in Remark~\refone{remark:universal-property} the Poincaré \(\infty\)-category \((\Modf{G},\QF^{\vis}_{\chi,\tau})\) can be characterized as a certain colimit in \(\Catp\). This can be interpreted as reflecting the structure of \(\Catp\) as a \(\Ct\)-category, see \S\refone{subsection:mackey-functors}, which admits quotients by actions of genuine \(\Ct\)-groups.
\end{variant}

\subsection{Parametrised spectra}
\label{subsection:parametrised-spectra}%

In this section we will discuss Poincaré structures on the \(\infty\)-category of compact \defi{parametrised spectra} over a space \(X\), whose \(\L\)-groups reproduces the \defi{visible \(\L\)-groups} of~\cite{WWIII} (see Corollary~\refone{corollary:compare-visible-WW} below). We will then show how to construct \defi{visible signatures} for Poincaré duality spaces in this setting.

Let us begin by establishing some terminology. We write \(\Spa_X := \Fun(X,\Spa)\) for the \(\infty\)-category of functors \(X \to \Spa\), to which we will refer to as \defi{local systems} of spectra, or as \defi{parametrised spectra} over \(X\). We then denote by \(\Spa_X^{\om} \subseteq \Spa_X\) the full subcategory spanned by the compact objects, and by \(\Spaf_X \subseteq \Spa_X^{\om}\) the smallest full stable subcategory containing the local systems \(x_!E\) left Kan extended from finite spectra \(E \in \Spaf\) along point inclusions \(\{x\} \subseteq X\). We then note that \(\Spa_X^{\om}\) is idempotent complete and the inclusion \(\Spaf_X \subseteq \Spa_X^{\om}\) is dense, and so exhibits \(\Spa_X^{\om}\) is the idempotent completion of \(\Spaf_X\).
When \(X\) is connected and pointed we may identify these subcategories with those of compact and finitely presented \(\SS[\Om X]\)-modules, respectively. The situation then renders in the context of \S\refone{subsection:visible}, and we may consider the visible Poincaré structures associated to various \(\Om X\)-characters.
The point of view of parametrised spectra has however the advantage of working without a preferred base point, and not being restricted to connected spaces. The unpointed setting is more natural, for example, when the input is a Poincaré duality space, as arising in the surgery classification of manifolds. The unpointed analogue of a \(\Om X\)-character is then given by a \defi{spherical fibration}, that is, a local system \(\xi\colon X \to \Spa\) which takes values in the full subgroupoid \(\Pic(\SS) \subseteq \iota\Spa\) spanned by the invertible spectra.

\begin{construction}
\label{construction:visible-param-spectra}%
Let \(X\) be a space and \(\xi\colon X \to \Pic(\SS)\) a spherical fibration on \(X\).
We associate to \(\xi\) a symmetric bilinear functor on \(\Spa_X^{\om} = \Fun(X,\Spa)^\omega\) by setting
\begin{equation}
\label{equation:bilin-xi}%
\Bil_{\xi}(L,L') := \map_{X \times X}(L\boxtimes L,\Delta_!(\xi)),
\end{equation}
where the mapping spectra takes place in the \(\infty\)-category \(\Fun(X \times X,\Spa)\), \(\boxtimes\) is the exterior tensor product, and \(\Delta_!\) is left Kan extension along the diagonal \(\Delta \colon X \to X\times X\).
The associated symmetric and quadratic hermitian structure on \( \Spa_X^{\om}\) are then given by
\[
\QF^{\sym}_{\xi}(L) = \map_{X\times X}(L\boxtimes L,\Delta_!(\xi))^{\hC}
\]
and
\[
\QF^{\qdr}_{\xi}(L) = \map_{X\times X}(L\boxtimes L,\Delta_!(\xi))_{\hC} ,
\]
respectively. By abuse of notation, we will also denote by \(\Bil_{\xi},\QF^{\sym}_{\xi}\) and \(\QF^{\qdr}_{\xi}\) the respective restrictions of these functors to the dense subcategory \(\Spaf_X \subseteq \Spa^{\om}_X\).

The above construction is functorial in \(X\) in the following sense. Let us call a \defi{map of spaces with spherical fibrations} a pair \((f,\tau)\) where \(f\colon X \to Y\) is a map between spaces equipped with spherical fibrations \(\xi_X\) and \(\xi_Y\) respectively, and \(\tau\colon \xi_X \to f^*\xi_Y\) is a natural transformation. We may then associate to \(f\) the corresponding left Kan extension functors
\[
f_! \colon \Spaf_X \to \Spaf_Y \quad\text{and}\quad f_! \colon \Spa^{\om}_X \to \Spa^{\om}_Y .
\]
The natural transformation \(\tau\) then induces a map
\[
\begin{split}
\Bil_{\xi_X}(L,L') = & \map_{X\times X}(L \boxtimes L',\Del_!\xi_X) \xrightarrow{\tau_*} \map_{X \times X}(L \boxtimes L',\Del_!f^*\xi_Y)  \\
\to & \map_{X \times X}(L \boxtimes L',(f \times f)^*\Del_!\xi_Y) \simeq  \map_{Y \times Y}((f \times f)_!(L \boxtimes L'),\Del_!\xi_Y) \simeq \Bil_{\xi_Y}(f_!L,f_!L')
\end{split}
\]
which is natural in \(L\) and \(L'\), so that we obtained an induced symmetric functor
\begin{equation}
\label{equation:induced-symmetric-functor}%
(\Spaf_X,\Bil_{\xi_X}) \to (\Spaf,\Bil_{\xi_Y})
 \quad\text{and}\quad (\Spa^{\om}_X,\Bil_{\xi_X}) \to (\Spa^{\om},\Bil_{\xi_Y})
\end{equation}
covering the left Kan extension functor \(f_!\).
\end{construction}

Our first goal is to show that the above construction gives a perfect bilinear functor and identify the associated duality. For this, note that for a fixed local system \(L \in \Spa_X\) the association \(L'\mapsto L \boxtimes L' \in \Spa_{X \times X}\) is colimit preserving and hence admits a right adjoint
\[
\map^{\boxtimes}(L,-)\colon \Spa_{X \times X} \to \Spa_X,
\]
which we can compute explicitly to be
\begin{align}
\label{equation:compute-map-box}%
\map^{\boxtimes}(L,T)_x \simeq &\map_{X}(x_!\SS,\map^{\boxtimes}(L,T)) \simeq \\
\nonumber&\map_{X \times X}(L \boxtimes x_!\SS ,T) \simeq \map_{X \times X}((j_x)_!L,T) \simeq \map_{X}(L,T|_{X \times \{x\}})
\end{align}
where \(x\colon \ast \to X\) is the insertion of the point \(x\) and \(j_x\colon X \times \{x\} \to X \times X\) is the corresponding insertion of the horizontal slice at height \(x\).

\begin{lemma}
\label{lemma:non-degenerate}%
For a spherical fibration \(\xi\), the bilinear functor~\eqrefone{equation:bilin-xi} is non-degenerate with duality
\[
\Dual_{\xi} L := \map^{\boxtimes}(L,\Delta_!\xi) .
\]
Furthermore, if \((f,\tau)\colon (X,\xi_X) \to (Y,\xi_Y)\) is a map of spaces with spherical fibrations such that \(\tau\colon \xi_X \to f^*\xi_Y\) is an equivalence then the induced symmetric functor~\eqrefone{equation:induced-symmetric-functor} is duality preserving.
\end{lemma}
\begin{proof}
Indeed, by adjunction we have
\[
\Bil_{\xi}(L,L') = \map_{X \times X}(L \boxtimes L',\Del_!\xi) = \map_X(L',\map^{\boxtimes}(L,\Del_!\xi)),
\]
and so \(\Bil_{\xi_X}\) is non-degenerate. To see the second claim, suppose fix a map \(f\colon X \to Y\) and an equivalence \(\tau\colon \xi_X \to \xi_Y\).
Since \(\Spa^{\om}_X\) is compactly generated by the collection of objects \(x_!\SS\), for \(x\colon \ast \to X\) a point, it suffices to check that the induced map \(f_!\Dual L \to \Dual f_!L\) is an equivalence for \(L=x_!\SS\). For this, it will suffice to prove the claim for the maps \(x\colon \ast \to X\) and \(f(x)\colon \ast \to Y\). In other words, we may as well assume that \(X=\ast\) and \(f\) is the inclusion of a point \(y \in Y\). Let \(\xi_y \in \Spa\) be the value of \(\xi_Y\) at \(y\), so that the goal becomes showing that the canonical symmetric functor
\[
(\Spa^{\om},\Bil_{\xi_y}) \to (\Spa^{\om}_Y,\Bil_{\xi_Y})
\]
is duality preserving. We now observe that if \(E \in \Spa^{\om}\) is a compact spectrum then we can factor the map \(y_!\Dual_{\xi_y}(E) \to \Dual_{\xi_Y}(y_!E)\) as the composite
\[
y_!(\Dual_{\xi_Y}(E)) = y_!(\Dual_{\SS}(E) \otimes_{\SS} \xi_{y}) \simeq\Dual_{\SS}(E) \otimes_{\SS} y_!\xi_y  \to
\Dual_{\SS}(E) \otimes_{\SS} (\Del_!\xi_Y)|_{\{y\} \times Y}\simeq \map^{\boxtimes}(y_!E,\Del_!\xi_Y) = \Dual_{\xi_Y}(y_!E)
\]
where the fourth equivalence is by the formula in~\eqrefone{equation:compute-map-box} and the third arrow
is the Beck-Chevalley transformation for the square
\[
\begin{tikzcd}
\{y\} \ar[r]\ar[d] & \{y\} \times Y \ar[d] \\
Y \ar[r,"{\Del}"] & Y \times Y
\end{tikzcd}
\]
where the top horizontal map picks the point \((y,y) \in \{y\} \times Y\). This transformation is an equivalence since the square is cartesian.
\end{proof}

\begin{corollary}
\label{corollary:coswandual}%
For a spherical fibration \(\xi\), the \(\xi\)-twisted duality
\[
\Dual_{\xi} L = \map^{\boxtimes}(L,\Delta_!\xi)
\]
is perfect. In particular, the associated symmetric and quadratic hermitian structures \(\QF^{\sym}_{\xi}\) and \(\QF^{\qdr}_{\xi}\) are Poincaré, and their restrictions to \(\Spaf_X\) are Poincaré as well.
\end{corollary}
\begin{proof}
We need to show that the evaluation map \(L \to \Dual\op\Dual L\) is an equivalence for any \(L \in \Spa_X^{\om}\). Since \(\Spa^{\om}_X\) is compactly generated by the collection of objects \(x_!\SS\), for \(x\colon \ast \to X\) a point, it suffices to check this for \(L=x_!\SS\). Invoking Lemma~\refone{lemma:non-degenerate} it will suffice to show that the duality \(\Dual_{\xi_x}\) on \(\Spa^{\om}_{\{x\}} = \Spa^{\om}\) is non-degenerate, that is, we may assume that \(X\) is a point. But then for any invertible spectrum \(E \in \Pic(\SS)\) we have that \(\Dual_E(-) \simeq E \otimes_{\SS} \Dual_{\SS}\) and is hence a perfect duality, the case of \(\Dual_{\SS}\) being the usual Spanier-Whitehead duality. Finally, \(\QF^{\sym}_{\xi}\) and \(\QF^{\qdr}_{\xi}\) also restrict to Poincaré structures on \(\Spaf_X\) by Observation~\refone{observation:poincare-subcat}; indeed, the equivalence \(\Dual_{\xi}x_!E \simeq x_!\Dual_{\xi_x}E\), implies that \(\Spaf_X\) is closed under the duality in \(\Spa^{\om}_X\).
\end{proof}

Out next goal is to construct a Poincaré structure on \(\Spa^{\om}_{X}\), which in some sense interpolates between the quadratic and symmetric structures. This structure is called the \defi{visible} Poincaré structure. To construct it, we first need to identify the linear part of the symmetric Poincaré structure \(\QF^{\sym}_{\xi}\). By definition, it is given by the formula
\begin{equation}
\label{equation:tate-param-spectra}%
L \mapsto \map_{X\times X}(L\boxtimes L, \Delta_!(\xi))^{\tC},
\end{equation}
which consequently constitutes an exact (contravariant) functor on \(\Spa^{\om}_X\). Such functors are always represented by an object in \(\Ind(\Spa^{\om}_X) \simeq \Spa_X\), and so there exists a (possibly non-compact) parametrised spectrum \(N\colon X \to \Spa\) such that
\[
\map_{X\times X}(L\boxtimes L, \Delta_!(\xi))^{\tC} \simeq \map_X(L,N)
\]
for \(L \in \Spa^{\om}_X\). To identify it, it suffices to check the values of~\eqrefone{equation:tate-param-spectra} at the generators \(x_!\SS\) of \(\Spa^{\om}_X\). In particular, \(N\) is canonically identified with the parametrised spectrum
\[
x \mapsto \map_{X\times X}(x_!\SS\boxtimes x_!\SS, \Delta_!\xi)^{\tC} = (\Del_!(\xi)_{(x,x)})^{\tC} ,
\]
or simply, \(N = (\Del^*\Del_!\xi)^{\tC}\). We thus obtain a natural transformation
\begin{equation}
\label{equation:parametrised-tate}%
\map_{X\times X}(L\boxtimes L, \Delta_!\xi)^{\tC} \to \map_X(L,(\Delta^*\Delta_! \xi)^{\tC})
\end{equation}
which is an equivalence for \(L \in \Spa^{\om}_X\).

\begin{definition}
\label{definition:xi-twisted-visible-structure}%
Let \(X\) be space and \(\xi\colon X \to \Pic(\SS)\) a spherical fibration on \(X\). We define the \defi{\(\xi\)-twisted visible Poincaré structure} \(\QF^{\vis}_{\xi}\colon \Spa_X^\omega \to \Spa\) by the top pullback square
\[
\begin{tikzcd}
\QF^{\vis}_\xi(L) \ar[rr] \ar[d] && \map_X(L,\xi) \ar[d] \\
\map_{X \times X}(L \boxtimes L,\Del_!\xi)^{\hC}\ar[dr]\ar[rr] &&  \map_X(L,(\Delta^*\Delta_!\xi)^{\tC}) \\
& \map_{X\times X}(L\boxtimes L, \Delta_!\xi)^{\tC} \ar[ur, "\simeq"] &
\end{tikzcd}
\]
where the right vertical map is induced by the canonical map \(\xi \to (\Delta^*\Delta_!\xi)^{\hC} \to (\Delta^*\Delta_!(\xi))^{\tC}\), which we can also identify with the composite \(\xi \to \xi^{\tC} \to (\Delta^*\Delta_!(\xi))^{\tC}\) (the first Tate construction being taken
with respect to the trivial action). By abuse of notation, we will also denote by \(\QF^{\vis}_{\xi}\) its restriction to the duality invariant dense subcategory \(\Spaf_X \subseteq \Spa^{\om}_X\). By construction, the linear part of \(\QF^{\vis}_{\xi}\) is given by
\[
\Lin^{\vis}_{\xi}(-) := \map_X(-,\xi).
\]
\end{definition}

\begin{example}
\label{example:parametrised-is-group-ring}%
If \(X\) is connected with base point \(x \in X\) then \(\Spaf_X\) and \(\Spa_X^{\om}\) are naturally equivalent to \(\Modf{\Om_xX}\) and \(\Modp{\Om_xX}\), respectively, and \(\xi\) can be encoded via a suitable \(\Om_xX\)-character \(\chi\). Under this equivalence, the visible Poincaré structure \(\QF^{\vis}_\xi\) corresponds to the visible Poincaré structure \(\QF^{\vis}_{\chi}\) of Construction~\refone{construction:visible}.
\end{example}

\begin{remark}
\label{remark:functoriality}%
Given a map \((f,\tau)\colon (X,\xi_X) \to (Y,\xi_Y)\) of spaces with spherical fibrations, the functoriality of Construction~\refone{construction:visible-param-spectra}
together with the naturality of the maps \(\xi \to (\Delta^*\Delta_!\xi)^{\hC} \to (\Delta^*\Delta_!(\xi))^{\tC}\) in \(\xi\) combine to yield induced hermitian functors
\[
(\Spaf_X,\QF^{\vis}_{\xi_X}) \to (\Spaf_Y,\QF^{\vis}_{\xi_Y}) \quad\text{and}\quad (\Spa^{\om}_X,\QF^{\vis}_{\xi_X}) \to (\Spa^{\om}_Y,\QF^{\vis}_{\xi_Y}).
\]
There are furthermore Poincaré when \(\tau\colon \xi_X \to f^*\xi_Y\) is an equivalence by Lemma~\refone{lemma:non-degenerate}.
\end{remark}

\begin{proposition}
\label{proposition:visible-is-colimit}%
Let \(X\) be a space equipped with a spherical fibration \(\xi \colon X \to \Pic(\SS)\), and for \(x \in X\) consider the associated maps \((x,\id)\colon (\{x\},\xi_x) \to (X,\xi)\) of spaces with spherical fibrations obtained by embedding the various \(x\) inside \(X\).
Then the natural transformation
\begin{equation}
\label{equation:colimit-visible}%
\colim_{x \in X}(\lambda_x)_!\QF^{\vis}_{\xi_x} \Rightarrow \QF^{\vis}_{\xi}
\end{equation}
induced by the associated Poincaré functors \((\lambda_x,\eta_x)\colon (\Spaf_{\{x\}},\QF^{\vis}_{\xi_x}) \to (\Spa^{\om}_X,\QF^{\vis}_{\xi})\) of Remark~\refone{remark:functoriality}, is an equivalence of quadratic functors on \(\Spa^{\om}_X\).
\end{proposition}

\begin{remark}
\label{remark:finite-visible-is-colimit}%
Since \(\Spaf_X \subseteq \Spa^{\om}_X\) is a full inclusion through which all the exact functors \(\lambda_x\colon \Spaf \to \Spa^{\om}_X\) factor, the pointwise formula for left Kan extensions shows that the statement made in Proposition~\refone{proposition:visible-is-colimit} about the Poincaré structure \(\QF^{\vis}_{\xi}\) on \(\Spa^{\om}_X\) is also inherited by its restriction to \(\Spaf_X\). In other words, that restriction is in turn given by the corresponding colimit of left Kan extensions of \(\QF^{\vis}_{\xi_x}\) along the various \(\lambda_x\colon \Spaf_{\{x\}} \to \Spaf_X\).
\end{remark}

\begin{remark}
\label{remark:visible-is-colimit}%
The universal property of the visible Poincaré structure \(\QF^{\vis}_{\xi}\) described in Proposition~\refone{proposition:visible-is-colimit} also determines universal properties of a similar nature for \((\Spaf_X,\QF^{\vis}_{\xi})\) and \((\Spa^{\om}_X,\QF^{\vis}_{\xi})\) as Poincaré \(\infty\)-categories. More precisely, we will prove in \S\refone{subsection:limits} that the collection of Poincaré functors \((\lambda_x,\eta_x)\colon (\Spaf,\QF^{\vis}_{\xi_x}) \to (\Spaf_X,\QF^{\vis}_{\xi})\) exhibits \((\Spaf_X,\QF^{\vis}_{\xi})\) as the \emph{colimit in \(\Catp\)} of the diagram \(x \mapsto (\Spaf,\QF^{\vis}_{\xi_x})\), see Proposition~\refone{proposition:visible-is-colimit-full}. Similarly, \((\Spa^{\om}_X,\QF^{\vis}_{\xi})\) is the colimit of the same diagram, but calculated inside \emph{idempotent complete} Poincaré \(\infty\)-categories. We will discuss idempotent completion of Poincaré \(\infty\)-categories in further details in~\papertwo.
\end{remark}

\begin{remark}
\label{remark:useful-formula}%
Denoting by \(\SS_X\) the constant spherical fibration on \(X\) with value \(\SS\), the proof of Proposition~\refone{proposition:visible-is-colimit} given below also establishes an identification
\[
\Bil_{\xi}(\Dual_{\SS_X}(L),\Dual_{\SS_{X}}(L')) \simeq \colim_{x \in X}[L_x \otimes_{\SS} L'_x \otimes_{\SS} \xi_x]
\]
as symmetric bilinear functors \(\Spa^{\om}_X \times \Spa^{\om}_X \to \Spa\) (in \(L,L'\)), and an equivalence
\[
\Lin^{\vis}_{\xi}(\Dual_{\SS_X}(L)) \simeq \colim_{x \in X}[L_x \otimes_{\SS}\xi]
\]
as functors \(\Spa^{\om}_X \to \Spa\). Under these equivalences the defining pullback square of \(\QF^{\vis}_{\xi}(\Dual_{\SS_X}(L))\) becomes
\[
\begin{tikzcd}
\QF^{\vis}_\xi(\Dual_{\SS_X} L) \ar[r] \ar[d] &
\colim_{x \in X} [L_x \otimes_{\SS} \xi_x] \ar[d] \\
\big(\colim_{x \in X}[L_x \otimes_{\SS} L_x \otimes_{\SS} \xi_x]\big)^{\hC}\ar[r] &
\big(\colim_{x \in X}[L_x \otimes_{\SS} L_x \otimes_{\SS} \xi_x]\big)^{\tC} \\
\end{tikzcd}
\]
where the right vertical map is given by the composite
\[
\colim_{x \in X} L_x \otimes_{\SS} \xi_x \to \colim_{x \in X}[(L_x \otimes_{\SS} L_x)^{\tC} \otimes_{\SS} \xi_x] \to \big(\colim_{x \in X}[L_x \otimes_{\SS} L_x \otimes_{\SS} \xi_x]\big)^{\tC}
\]
using the Tate diagonal of \(L_x\).
\end{remark}

\begin{proof}[Proof of Proposition~\refone{proposition:visible-is-colimit}]
Being a map of quadratic functors, it suffices by Corollary~\refone{corollary:classification-of-quad-functors} to show that~\eqrefone{equation:colimit-visible} induces an equivalence on bilinear and linear parts. Since taking linear and bilinear parts commutes with colimits, and by Proposition~\refone{proposition:left-kan-bilinear-linear} also with left Kan extensions, the two maps that we need to consider are
\begin{equation}
\label{equation:colimit-bilinear-linear}%
\colim_{x \in X}(\lambda_x \times \lambda_x)_!\Bil_{\xi_x} \Rightarrow \Bil_{\xi} \quad\text{and}\quad \colim_{x \in X}(\lambda_x)_!\Lin^{\vis}_{\xi_x} \Rightarrow \Lin^{\vis}_{\xi}.
\end{equation}
Let \(\wtl{\Bil}_{\xi_x}\) and \(\wtl{\Lin}^{\vis}_{\xi_x}\) be the left Kan extensions of \(\Bil_{\xi_x}\) and \(\Lin^{\vis}_{\xi_x}\) to \(\Pro(\Spaf)\op \times \Pro(\Spaf)\op\) and \(\Pro(\Spaf)\op\), respectively. As in the proof of Lemma~\refone{lemma:kan-extension-exact-quadratic}, we may compute the left Kan extensions of \(\Bil_{\xi_x}\) and \(\Lin^{\vis}_{\xi_x}\) along \(\lambda_x \times \lambda_x\) and \(\lambda_x\) respectively by restricting \(\wtl{\Bil}_{\xi_x}\) and \(\wtl{\Lin}^{\vis}_{\xi_x}\) along a given pro-left adjoint
\[
\Spa^{\om}_X \to \Pro(\Spaf),
\]
that is, a restriction to \(\Spa^{\om}_X\) of a left adjoint to the induced functor \(\Pro(\lambda_x)\colon \Pro(\Spaf) \to \Pro(\Spa^{\om}_X)\). To produce such a pro-left adjoint, we note that for every \(x \in X\), the duality \(\Dual_{\xi_x}\colon (\Spaf)\op \xrightarrow{\simeq} \Spaf\) induces an equivalence
\[
\wtl{\Dual}_{\xi_x}\colon \Spa\op = \Ind(\Spaf)\op = \Pro((\Spaf)\op) \xrightarrow{\Pro(\Dual_{\xi_x})} \Pro(\Spaf).
\]
Now since evaluation at \(x\) is right adjoint to left Kan extension \(\Spa \to \Spa_X\) along \(\{x\} \subseteq X\) it now follows that a pro-left adjoint for \(\lambda_x\) is given by the composite
\[
\Spa^{\om}_X \xrightarrow{\Dual_{\xi}\op}  (\Spa^{\om}_X)\op \subseteq \Spa_X\op \xrightarrow{\ev_x\op} \Spa\op \xrightarrow{\wtl{\Dual}_{\xi_x}} \Pro(\Spaf).
\]
Using this pro-left adjoint to express the left Kan extension of \(\Bil_{\xi_x}\) along \((\lambda_x \times \lambda_x)\) we now obtain
\[
[(\lambda_x \times \lambda_x)_!\Bil_{\xi_x}](L,L') \simeq \wtl{\Bil}_{\xi_x}(\wtl{\Dual}_{\xi_x}((\Dual_{\xi}L)_x),\wtl{\Dual}_{\xi_x}((\Dual_{\xi}L)_x))
\]
and
\[
[(\lambda_x)_!\Lin^{\vis}_{\xi_x}](L) \simeq \wtl{\Lin}^{\vis}_{\xi_x}(\wtl{\Dual}_{\xi_x}((\Dual_{\xi}L)_x)) .
\]
To compute these terms further, let us first unwind the pre-compositions of \(\Bil_{\xi_x}\) and \(\Lin^{\vis}_{\xi_x}\) with the duality \(\Dual_{\xi_x}\). For finite spectra \(E,E',E''\) we have
\[
\Bil_{\xi_x}(\Dual_{\xi_x}E,\Dual_{\xi_x}E') = \map(\Dual_{\xi_x}E \otimes_{\SS} \Dual_{\xi_x}E, \xi_x) = (E \otimes_{\SS} \xi_x^{-1}) \otimes_{\SS} (E' \otimes_{\SS} \xi_x^{-1}) \otimes_{\SS} \xi_x
\]
and
\[
\Lin^{\vis}_{\xi_x}(\Dual_{\xi_x}E'') =  \map(\Dual_{\xi_x}E'', \xi_x) = E'' .
\]
Now the two expressions on the right hand sides are well-defined for \(E,E',E''\) which are not necessarily finite, and are colimit preserving in each of these inputs separately. These formulas hence also describe the left Kan extensions of \(\Bil_{\xi_x}(\Dual_{\xi_x}(-),\Dual_{\xi_x}(-))\) and \(\Lin^{\vis}_{\xi_x}(\Dual_{\xi_x}(-))\) to \(\Spa \times \Spa\) and \(\Spa\), respectively. We may consequently conclude that
\[
\wtl{\Bil}_{\xi_x}(\wtl{\Dual}_{\xi_x}((\Dual_{\xi}L)_x),\wtl{\Dual}_{\xi_x}((\Dual_{\xi}L)_x)) \simeq ((\Dual_{\xi}L)_x \otimes_{\SS} \xi_x^{-1}) \otimes_{\SS} ((\Dual_{\xi}L')_x \otimes_{\SS} \xi_x^{-1}) \otimes_{\SS} \xi_x
\]
and
\[
\wtl{\Lin}^{\vis}_{\xi_x}(\wtl{\Dual}_{\xi_x}((\Dual_{\xi}L'')_x)) \simeq (\Dual_{\xi}L'')_x.
\]
Combining all the above with the equivalences
\[
(\Dual_{\xi}L)_x  \simeq  \map_{X}(L,(\Del_!\xi)|_{X \times \{x\}})\simeq \map_X(L,x_!\xi_x)
\]
supplied by~\eqrefone{equation:compute-map-box} and the Beck-Chevalley property, we now conclude that for local systems \(L,L'\in \Spa^{\om}_X\) we have
\begin{align*}
[(\lambda_x \times \lambda_x)_!\Bil_{\xi_x}](L,L') & \simeq (\map_X(L,x_!\xi_x) \otimes_{\SS} \xi_x^{-1}) \otimes_{\SS} (\map_X(L',x_!\xi_x) \otimes_{\SS} \xi_x^{-1}) \otimes_{\SS} \xi_x \\
 & \simeq \map_X(L,x_!\SS)\otimes_{\SS} \map_X(L',x_!\SS) \otimes_{\SS} \xi_x,
\end{align*}
and
\[
[(\lambda_x)_!\Lin^{\vis}_{\xi_x}](L'') \simeq \map_X(L'',x_!\xi_x).
\]
We may now finally identify the maps~\eqrefone{equation:colimit-bilinear-linear} in explicit terms as follows. The first map in~\eqrefone{equation:colimit-bilinear-linear} identifies with the composite
\begin{align*}
\colim_{x \in X}\big[\map_X(L,x_!\SS)\otimes_{\SS} \map_X(L',x_!\SS) \otimes_{\SS} \xi_x\big] & \xrightarrow{\alp_{L,L'}} \colim_{x \in X}\map_{X \times X}\big(L \boxtimes L',[(x,x)_!\SS] \otimes \xi_x\big) \\
&\xrightarrow{\beta_{L,L'}} \map_{X \times X}(L \boxtimes L',\Del_!\xi),
\end{align*}
where \(\alp_{L,L'}\) is induced on colimits by the functoriality of \(\boxtimes\) and the identification \(x_!\SS \boxtimes x_!\SS \simeq (x,x)_!\SS\), and \(\beta_{L,L'}\) is given by post-composition with the map \([(x,x)_!\SS] \otimes \xi_x \simeq (x,x)_!\xi_x \simeq \Del_!(x_!\xi_x) \to \Del_!\xi\) induced by the counit maps \(x_!\xi_x \to \xi\).
In a similar manner, the second map in~\eqrefone{equation:colimit-bilinear-linear} now becomes the map
\[
\colim_{x \in X} \map_{X}(L'',x_!\xi_x)  \xrightarrow{\gam_{L''}} \map_X(L'',\xi).
\]
induced by the counit maps \(x_!\xi_x \to \xi\).
We now note that the collection of counit maps \(x_!\xi_x \to \xi\) exhibit \(\xi\) as the colimit in \(\Spa^{\om}_X\) of the diagram \(x \mapsto x_!\xi_x\). Indeed, identifying \(x_!\xi_x\) back with \(\Del_!\xi|_{X \times \{x\}}\), this is just the statement that the map \(\Del_!\xi \to \pi_1^*\xi\) exhibits \(\xi\) as the left Kan extension of \(\Del_!\xi\) along the projection \(\pi_1\colon X \times X \to X\) on the first argument, which is clear since \(\pi_1\Del \simeq \id\). We then conclude that \(\beta_{L,L'}\) and \(\gam_{L''}\) are equivalences when \(L,L',L''\) are perfect since in this case the functors \(\map_{X \times X}(L \boxtimes L',-)\) and \(\map_X(L'',-)\) are colimit preserving. To see that \(\alp_{L,L'}\) is an equivalence, we note that since it is a natural transformations of bilinear functors in \(L,L'\), to check that it is an equivalence on perfect \(L,L'\) it suffices to verify the case where
\(L=y_!E\) and \(L'=y'_!E'\)
for points \(y,y' \in X\), since local systems of this form generate all perfect local systems under finite colimits and retracts. Then \(\alp_{y_!E,y'_!E'}\) can be identified with the map
\[
\colim_{x \in X}\big[(x_!\SS)_y\otimes (x_!\SS)_{y'} \otimes_{\SS} \xi_x\big] \to \colim_{x \in X}\big[((x, x)_!\SS)_{(y,y')} \otimes_{\SS} \xi_x\big]
\]
which is an equivalence since the maps
\[
(x_!\SS)_y\otimes (x_!\SS)_{y'} = \Sigma^{\infty}_+\Map(x,y) \otimes_{\SS} \Sigma^{\infty}_+\Map(x,y') \to \Sigma^{\infty}_+(\Map(x,y) \times \Map(x,y')) \simeq ((x, x)_!\SS)_{(y,y')}
\]
are equivalences for every \(x \in X\).
\end{proof}

Let us now pause to explain the relation between the Poincaré \(\infty\)-category \((\Spa^{\om}_X,\QF^{\vis}_{\xi})\) with the construction used in Weiss-Williams~\cite{WWIII} to define visible \(\L\)-theory. To simplify the discussion, let us consider the case where \(\xi=\SS_X\), though the argument in the general case proceeds in a similar manner. To begin, recall first that the straightening-unstraightening equivalence \(\Fun(X,\Sps) \simeq \Sps_{/X}\) induces an equivalence
\[
\Spa_X =\Fun(X,\Spa) \simeq \Spa(\Fun(X,\Sps)_{*/}) \simeq \Spa(\Sps_{X//X}),
\]
where \(\Sps_{X//X} := (\Sps_{/X})_{\id_X/}\) is the \(\infty\)-category of retractive spaces over \(X\). This restricts to an equivalence between the full subcategory \(\Spa^{\om}_X \subseteq \Spa_X\) of perfect parametrised spectra and the full subcategory of \(\Spa(\Sps_{X//X})\) spanned by the suspension spectra of \emph{finitely dominated} retractive spaces \([X \to Y \to X] \in \Sps^{\fd}_{X//X}\). Here the full subcategory \(\Sps^{\fd}_{X//X} \subseteq \Sps_{X//X}\) of finitely dominated retractive spaces is the closure under retracts of the full subcategory of \(\Spa_{X//X}\) spanned by the retractive spaces \(X \to Y \to X\) such that \(X \to Y\) can be represented by an inclusion of simplicial sets with \(Y\) having only finitely many non-degenerated simplices not in \(X\).
In particular, we may identify \(\Spa_X^{\om}\) with the Spanier-Whitehead stabilisation of \(\Sps^{\fd}_{X//X}\).

If we now represent \(X\) by a CW complex then the \(\infty\)-category \(\Sps^{\fd}_{X//X}\) can be modelled by the Waldhausen category \(\R(X)\) of finitely dominated retractive topological spaces \(X \to Y \to X\) such that \(X \to Y\) is a Serre cofibration, with weak equivalences \(W\) the weak homotopy equivalences, that is, \(\Sps^{\fd}_{X//X} \simeq \R(X)[W^{-1}]\). Indeed, this is a standard consequence of the fact that \(\Sps_{X//X}\) is presented by the Serre model structure on \(\Top_{X//X}\), and \(\R(X) \subseteq \Top_{X//X}\) is a full subcategory consisting of cofibrant objects and closed under weak equivalences between cofibrant objects. From this one may obtain a Waldhausen model \(\sR(X)\) for the Spanier-Whitehead stabilisation by considering the category whose objects are pairs \((Y,k)\) where \(X \to Y \to X\) is a finitely dominated retractive topological space over \(X\) and \(k\in \ZZ\) is an integer, with
the set of maps from \((Y,k)\) to \((Y',k')\) being
\[
\Hom_{\sR(X)}((Y,k),(Y',k')) := \colim_{n\geq k,k'}\Hom_{\R(X)}(\Sig^{n-k}Y,\Sig^{n-k'}Y') .
\]
The object \((Y,k)\) should be considered as the formal \(k\)'th desuspension of \(X \xrightarrow{i} Y \xrightarrow{f} X\). Endowing \(\sR(X)\) with the collection \(sW\) of stable weak homotopy equivalences, that is, the maps \((Y,k) \to (Y',k')\) which can be represented by a weak homotopy equivalence \(\Sig^{n-k}Y \to \Sig^{n-k'}Y'\) for some \(n\), one obtains a model for perfect parametrised spectra, in the sense that \(\sR(X)[sW^{-1}] \simeq \Spa^{\om}_X\). Indeed, this is a formal consequence of the above and the commutation of localization with colimits:
\[
\big(\colim[\R(X) \xrightarrow{\Sig} \R(X) \xrightarrow{\Sig}\hdots ]\big)[sW^{-1}] \simeq \colim[\R(X)[W^{-1}] \xrightarrow{\Sig} \R(X)[W^{-1}] \xrightarrow{\Sig} \hdots] \simeq \Spa^{\om}_X,
\]
where we note that the sequential colimit on the left hand side (before localization) is given by \(\sR(X)\) both as ordinary categories and as \(\infty\)-categories.
In this presentation, the object \((Y,k)\) corresponds to the local system
\(x \mapsto \Sig^{\infty-k}Y_x\),
where \(Y_x\) denotes the homotopy fibre of \(f\colon Y \to X\) over \(x\), considered as a pointed space using the map \(X \to Y\).
The Waldhausen category \(\sR(X)\) is equipped with a stable Spanier-Whitehead product (see discussion at the end of \S\refone{subsection:metabolic-and-L}):
\[
\odot_{\bullet} \colon \sR(X) \times \sR(X) \to \Spa^{\Om}
\]
given by
\[
(Y,k) \odot_{\bullet} (Y',k') := \Sig^{\infty-k-k'}(Y \wedge^h_X Y'),
\]
that is, by (the implicit \(\Om\)-spectrum replacement of) the shifted suspension spectrum of the pointed space
\[
Y \wedge^h_X Y' := [Y \times^h_X Y'] / \big[Y \coprod_X Y'\big] =
[Y \times_X X^I \times_X Y'] / \big[Y  \coprod_X Y'\big] ,
\]
the symbol \(\slash\) standing for collapsing a subspace to a point. Here \(X^I\) is the space of paths in \(X\), used here to form an explicit functorial model for the homotopy fibre product. This Spanier-Whitehead product preserves weak equivalences in each entry and thus descends to a bilinear functor
\[
\B^{\WW}_{\vis}\colon \Spa^{\om}_X \times \Spa^{\om}_X \to \Spa ,
\]
which we can now compare with the bilinear functor on parametrised spectra considered here:
\begin{lemma}
\label{lemme:compare-bilin-WW}%
The bilinear functor \(\Bil^{\WW}_{\vis}\) is naturally equivalent to the bilinear functor
\[
(L,L') \mapsto \Bil_{\SS_X}(\Dual_{\SS_X}L,\Dual_{\SS_X}L').
\]
\end{lemma}
\begin{proof}
Since \(\Sig^\infty\) is colimit preserving and symmetric monoidal for \(\wedge\) we have
\begin{align*}
\Sig^{\infty-k-k'}(Y \wedge^h_X Y') &\simeq \Sig^{\infty-k-k'}(\colim_x Y_x \wedge Y'_x) \\
&\simeq \colim_x \Sig^{\infty-k-k'}(Y_x \wedge Y'_x) \\
&\simeq \displaystyle\mathop{\colim}_{x \in X}[\Sig^{\infty-k}Y_x  \otimes_{\SS} \Sig^{\infty-k'}Y'_x]\\
&\simeq \Bil_{\SS_X}(\Dual_{\SS_X}[x \mapsto \Sig^{\infty-k}Y_x],\Dual_{\SS_X}[x \mapsto \Sig^{\infty-k}Y_x]) \ ,
\end{align*}
where the last equivalence is by Remark~\refone{remark:useful-formula}. Since \(\sR(X) \to \Spa^{\om}_X\) is a localisation functor this natural equivalence descends to a natural equivalence \(\Bil^{\WW}_{\vis}(L,L') \simeq \Bil_{\SS_X}(\Dual_{\SS_X}L,\Dual_{\SS_X}L')\), as desired.
\end{proof}

Combining Lemma~\refone{lemme:compare-bilin-WW} with Proposition~\refone{proposition:compare-WW} and Remark~\refone{remark:also-quadratic} we now deduce
\begin{corollary}
\label{cor:compareWWquadsymL}%
The symmetric and quadratic \(\L\)-groups associated to the Spanier-Whitehead product \(\odot_\bullet\) are naturally equivalent to the \(\L\)-groups of the Poincaré structures \(\QF^{\sym}_{\SS_X}\) and \(\QF^{\qdr}_{\SS_X}\) of Construction~\refone{construction:visible-param-spectra}.
\end{corollary}

We would like to extend this comparison to visible \(\L\)-groups. These are defined in~\cite{WWIII} by replacing the notion of symmetric co-forms on \((Y,k)\), given by the homotopy fixed points spectrum \(((Y,k) \odot_\bullet (Y,k))^{\hC}\), with the corresponding \emph{genuine} fixed points \(((Y,k) \odot_\bullet (Y,k))^{\Ct}\), where \((Y,k) \odot_\bullet (Y,k)\) is considered as a genuine \(\Ct\)-spectrum by identifying it with the shifted suspension \(\Ct\)-spectrum \(\Sig^{\infty-k\rho}(Y \wedge^h_X Y)\). Here \(\rho\) is the regular \(\Ct\)-representation and \(Y \wedge^h_X Y\) is considered as a genuine \(\Ct\)-space whose \(\Ct\)-fixed points are its point-set level fixed points
\[
(Y \wedge^h_X Y)^{\Ct} \cong \big[Y \times_X X^{[0,1/2]}\big] /X \simeq Y/X .
\]
The natural map \(((Y,k) \odot_\bullet (Y,k))^{\Ct} \to ((Y,k) \odot_\bullet (Y,k))^{\hC}\) then enables one to define unimodular visible co-forms and hence visible \(\L\)-groups by applying the machinery of~\cite{WW-duality}. We note that the functor \((Y,k) \mapsto ((Y,k) \odot_\bullet (Y,k))^{\Ct}\) preserves weak equivalences and hence descends to a functor
\[
\QF^{\WW}\colon \Spa^{\om}_X \to \Spa ,
\]
which we would like to compare with the visible Poincaré structure on parametrised spectra considered here:

\begin{proposition}
\label{proposition:compare-visible-WW}%
The functor \(\QF^{\WW}\) is naturally equivalent to the quadratic functor \(L \mapsto \QF^{\vis}_{\SS_X}(\Dual_{\SS_X}L)\).
\end{proposition}

Before we come to the proof of Proposition~\refone{proposition:compare-visible-WW}, let first recall that the notion of genuine \(\Ct\)-spectra admits a variety of point-set models, among which are orthogonal \(\Ct\)-spectra~\cite{mandell-may-equivariant, schwede-global}, symmetric \(\Ct\)-spectra~\cite{mandell-symmetric}, prespectra indexed over the poset of finite dimensional sub-representations of a complete \(\Ct\)-universe~\cite{mandell-may-equivariant}, or
alternatively over a submonoid of the representation ring of \(\Ct\)~\cite{mandell-symmetric}, this last model being the one employed in~\cite{WWII} using the submonoid spanned by the regular representation. All these models are known to be Quillen equivalent to each other, see~\cite[Theorem 4.16]{mandell-may-equivariant},~\cite[Theorem 6.2]{mandell-symmetric} and~\cite[Theorem 10.2]{mandell-symmetric}. On the other, the notion of \(\Ct\)-spectra can also be defined \(\infty\)-categorically using either Mackey functors or the stabilisation of \(\Ct\)-spaces along the regular representation sphere~\cite{Barwick-MackeyI, barwick2016parametrized, Nardin-stability, shah2018parametrized}, and this approach is equivalent to the model categorical one by~\cite[Theorem A.4]{Nardin-stability}. The last equivalence between the \(\infty\)-categorical and model categorical approaches is furthermore compatible with the formation of genuine fixed points on the one hand, and the formation of suspension \(\Ct\)-spectra on the other. In particular, the construction of Weiss and Williams can be performed using any of these models by simply implementing the formation of genuine fixed points and suspension \(\Ct\)-spectra in the relevant context. We consequently may and will freely apply ideas and results from the theory of genuine \(\Ct\)-spectra in the setting of~\cite{WWIII}.

\begin{proof}[Proof of Proposition~\refone{proposition:compare-visible-WW}]
Since \(\sR(X)[sW^{-1}] \to \Spa^{\om}_X\) is a localisation functor it will suffice to construct, naturally in a given retractive space \(X \xrightarrow{i} Y \xrightarrow{f} X \in \sR(X)\), an equivalence
\[
\QF^{\vis}_{\SS_X}(\Dual_{\SS_X}[x \mapsto \Sig^{\infty}Y_x]) \simeq [(Y,k) \odot_{\bullet} (Y,k)]^{\Ct}.
\]
Now since \((Y,k) \odot_{\bullet} (Y,k)\) is a shifted suspension \(\Ct\)-spectrum its geometric fixed points are given by
\[
[(Y,k) \odot_{\bullet} (Y,k)]^{\Phi \Ct} = \Sig^{\infty-k}[Y \wedge^h_X Y]^{\Ct} \simeq  \Sig^{\infty-k}(Y/X),
\]
and its isotropy separation square
is
\[
\begin{tikzcd}[column sep=8pt]
[(Y,k) \odot_{\bullet} (Y,k)]^{\Ct}\ar[r]\ar[d] & \Sig^{\infty-k}[Y \wedge^h_X Y]^{\Ct} \ar[r,phantom,"\simeq"] \ar[d] & \Sig^{-k}\Sig^{\infty}(Y/X)  \\
{[\Sig^{\infty-k\rho}(Y \wedge^h_X Y)]}^{\hC} \ar[r] &
{[\Sig^{\infty-k\rho}(Y \wedge^h_X Y)]}^{\tC} \ar[r,phantom,"\simeq"] & \Sig^{-k}{[\Sig^{\infty}(Y \wedge^h_X Y)]}^{\tC} \ ,
\end{tikzcd}
\]
where the right vertical map is the \(k\)-fold desuspension of the composite
\[
\Sig^{\infty}(Y/X) \to [\Sig^{\infty}(Y \wedge^h_X Y)]^{\hC} \to [\Sig^{\infty}(Y \wedge^h_X Y)]^{\tC}
\]
induced by the \(\Ct\)-equivariant diagonal \(Y/X \to Y \wedge^h_X Y= [Y \times^h_X Y]/[Y \coprod_X Y]\).
To construct the desired equivalence we would like to compare this isotropy separation square
with the square
\[
\begin{tikzcd}[column sep=7pt]
\QF^{\vis}_\xi(\Dual_{\SS_X}[x \mapsto \Sig^{\infty}Y_x]) \ar[rr] \ar[d] && \displaystyle\mathop{\colim}_{x \in X} \Sig^{\infty-k}Y_x  \ar[d]\ar[r,phantom,"\simeq"] & \Sig^{-k}\displaystyle\mathop{\colim}_{x \in X} \Sig^{\infty}Y_x \\
\big(\displaystyle\mathop{\colim}_{x \in X}[\Sig^{\infty-k}Y_x \otimes_{\SS} \Sig^{\infty-k}Y_x]\big)^{\hC}\ar[rr] &&
\big(\displaystyle\mathop{\colim}_{x \in X}[\Sig^{\infty-k}Y_x \otimes_{\SS} \Sig^{\infty-k}Y_x]\big)^{\tC} \ar[r,phantom,"\simeq"] & \Sig^{-k}\big(\displaystyle\mathop{\colim}_{x \in X}[\Sig^{\infty}Y_x \otimes_{\SS} \Sig^{\infty}Y_x]\big)^{\tC} \\
\end{tikzcd}
\]
produced by Remark~\refone{remark:useful-formula}. Indeed, since in the \(\infty\)-category of pointed spaces we have \(\colim_x Y_x \simeq Y/X\) and \(\colim_x Y_x \wedge Y_x \simeq Y \wedge^h_X Y\), and since \(\Sig^\infty\) is colimit preserving and symmetric monoidal for \(\wedge\),
we may proceed as in the proof of Lemma~\refone{lemme:compare-bilin-WW} to identify the three bottom right corners in the two squares, as well as the map from the homotopy fixed points to the corresponding Tate object, which in both cases is the component of the natural transformation \((-)^{\hC} \Rightarrow (-)^{\tC}\). To obtain an induced equivalence on the top left corner it is will now suffice to identify the right vertical map in the two squares. Peeling away the \(k\)-fold desuspension and using the fact that the Tate diagonal of \(\Sig^{\infty}Y_x\) is induced by the diagonal \(Y_x \to Y_x \wedge Y_x\) on the level of pointed spaces, the desired identification is now encoded by the following canonical commutative diagram
\[
\begin{tikzcd}
\displaystyle\mathop{\colim}_{x \in X} \Sig^{\infty}Y_x \ar[r] \ar[d,equal] &
\displaystyle\mathop{\colim}_{x \in X} {[\Sig^{\infty}Y_x \otimes_{\SS} \Sig^{\infty}Y_x]}^{\hC} \ar[d]\ar[r] &
\displaystyle\mathop{\colim}_{x \in X}{[\Sig^{\infty}Y_x \otimes_{\SS} \Sig^{\infty}Y_x]}^{\tC} \ar[d] \\
\displaystyle\mathop{\colim}_{x \in X} \Sig^{\infty}Y_x \ar[r]\ar[d,"\simeq"]&{\big[\displaystyle\mathop{\colim}_{x \in X} (\Sig^{\infty}Y_x \otimes_{\SS} \Sig^{\infty}Y_x)\big]}^{\hC} \ar[r]\ar[d,"\simeq"] & {\big[\displaystyle\mathop{\colim}_{x \in X} (\Sig^{\infty}Y_x \otimes_{\SS} \Sig^{\infty}Y_x)\big]}^{\tC} \ar[d,"\simeq"] \\
\Sig^{\infty}Y/X \ar[r] & {[\Sig^{\infty}(Y \wedge^h_X Y)]}^{\hC} \ar[r] & {[\Sig^{\infty}(Y \wedge^h_X Y)]}^{\tC}  \ .
\end{tikzcd}
\]
\end{proof}

\begin{corollary}
\label{corollary:compare-visible-WW}%
The visible \(\L\)-groups of~\cite{WWIII} are naturally equivalent to the \(\L\)-groups of the visible Poincaré structure \(\QF^{\vis}_{\SS_X}\).
\end{corollary}
\begin{proof}
The visible \(\L\)-groups are (litterally) defined in~\cite{WWIII} by invoking the construction of~\cite{WWIII} and replacing everywhere symmetric Poincaré objects by visible ones. The proof of the present corollary is then identical to that of Proposition~\refone{proposition:compare-WW}, replacing everywhere symmetric Poincaré objects by visible ones, using Proposition~\refone{proposition:compare-visible-WW}.
\end{proof}

\begin{variant}
\label{variant:coefficients}%
Let \(X\) be a space with a spherical fibration \(\xi\colon X \to \Pic(\SS)\) and let \((A,N,\alp)\) be a ring spectrum with genuine involution. We may then form a Poincaré structure \(\QF^{\vis,\alp}_{\xi}\) on the \(\infty\)-category \(\Fun(X,\Mod_A)^{\om}\) of compact local systems of \(A\)-modules by the top pullback square
\[
\begin{tikzcd}
[column sep=4ex]
\QF^{\vis,\alp}_\xi(L) \ar[rr] \ar[dd] && \map_X(L,\xi \otimes_{\SS} N) \ar[d] \\
  && \map_X(L,(\Delta^*\Delta_!\xi)^{\tC} \otimes_{\SS} A^{\tC}) \ar[d] \\
\map_{X \times X}(L \boxtimes L,\Del_!\xi \otimes_{\SS} A)^{\hC}\ar[dr]\ar[rr] &&  \map_X(L,(\Delta^*\Delta_!\xi \otimes_{\SS} A)^{\tC}) \\
& \map_{X\times X}(L\boxtimes L, \Delta_!\xi \otimes_{\SS} A)^{\tC} \ar[ur, "\simeq"] &
\end{tikzcd}
\]
Here we note that the \(\Ct\)-equivariant structure of \(\Del_!\xi \otimes_{\SS} A\) with respect to the flip action on \(X \times X\) is induced by the equivariant structure of \(\Del_!\xi\) and the involution on \(A\). Arguing as in Lemma~\refone{lemma:non-degenerate} and Corollary~\refone{corollary:coswandual} one sees again that the hermitian \(\infty\)-category \((\Fun(X,\Mod_A)^{\om},\QF^{\vis,\alp}_{\xi})\) is functorial in maps \((f,\tau) \colon (X,\xi_X) \to (Y,\xi_Y)\) of spaces with spherical fibrations and that \(\QF^{\vis,\alp}_{\xi}\) is again Poincaré with underlying perfect duality \(L \mapsto \map^{\boxtimes}_A(L,\Del_!(\xi)\otimes_{\SS} A)\), where \(\map^{\boxtimes}_A(L,-)\) denotes the right adjoint of the functor
\[
L\boxtimes (-)\colon \Fun(X,\Mod_A) \to \Fun(X\times X,\Mod_{A \otimes_{\SS} A}).
\]
\end{variant}

\begin{remark}
In addition to its functoriality in maps of spaces with spherical fibrations, the construction of Variant~\refone{variant:coefficients} is also functorial in maps in \((A,N,\alp) \to (B,K,\bet)\) of rings with genuine involution. In fact, the Poincaré \(\infty\)-category \((\Fun(X,\Mod_A)^{\om},\QF^{\vis,\alp}_{\xi})\) depends functorially on the pair of Poincaré \(\infty\)-categories \((\Spa^{\om}_X,\QF^{\vis}_{\xi})\) and \((\Mod_A,\QF^{\alp}_A)\): we may identify it with their \defi{tensor product} \((\Spa^{\om}_X,\QF^{\vis}_{\xi}) \otimes_{\SS} (\Mod_A,\QF^{\alp}_A) \in \Catp\), a construction we will study in \S\refone{section:multiplicative}. As in Remark~\refone{remark:visible-is-colimit}, one can then identify \((\Fun(X,\Mod_A)^{\om},\QF^{\vis,\alp}_{\xi})\) with the \emph{colimit} in \(\Catp\) of the diagram \(x \mapsto (\Mod_A,\QF^{\vis,\alp}_{\xi_x})\). The Poincaré structures \(\QF^{\vis,\alp}_{\xi_x}\) are then the ones associated to the module with genuine involution \((A,N,\alp)\otimes_{\SS} \xi_x = (A \otimes_{\SS} \xi_x, N \otimes_{\SS} \xi_x,\alp_x)\), where \(\alp_x\) is the composed map \(N \otimes_{\SS} \xi_x \to A^{\tC} \otimes_{\SS} \xi_x^{\tC} \to (A \otimes_{\SS} \xi_x)^{\tC}\).
\end{remark}

For the remainder of this section we will show how to construct visible signatures for Poincaré duality spaces in the present context.
We take a purely homotopy theoretical approach to Poincaré duality spaces and their Spivak normal fibration following \cite{klein-dualizing}. Let \(X\) be a finite space (that is, a space which can be realized by a simplicial set with only finitely many non-degenerate simplices). The right Kan extension functor \(r_\ast\colon \Spa_X \to \Spa\) along \(r\colon X \to \ast\) is given by taking the limit along \(X\), which is in particular a finite limit and hence preserves colimits by the stability of \(\Spa_X\). It is thus equivalent to a functor of the form \(r_!(\mathcal{F}_X \otimes (-))\) for an essentially unique \(\mathcal{F}_X\), which is then called the \defi{dualizing complex} of \(X\). Here \(\otimes\) stands for the \emph{pointwise} tensor product of \(X\)-parametrised spectra. We say that \(X\) is a \defi{Poincaré duality space} if the parametrised spectrum \(\mathcal{F}_X\colon X\to \Spa\) factors through the full subgroupoid \(\Pic(\SS) \subseteq \iota\Spa\) of invertible objects. In this case we denote the resulting spherical fibration by \(\nu := \mathcal{F}_X\), and call it the \defi{Spivak normal fibration} of \(X\). We note that the identification of \(r_!(\nu \otimes (-))\) with \(r_*\) can equivalently be encoded via a unit map
\[
c_E\colon E \to r_!(\nu \otimes r^*E) \simeq r_!(\nu) \otimes_{\SS} E
\]
exhibiting \(r_!(\nu \otimes (-))\) as right adjoint to \(r^*\). Since \(\Funx(\Spaf,\Spaf) \simeq \Spaf\) via evaluation at \(\SS\) the natural transformation \(c_{(-)}\) is canonically of the form \(c \otimes_{\SS} (-)\), where
\begin{equation}
\label{equation:collapse}%
c = c_{\SS}\colon \SS \to r_!(\nu)
\end{equation}
is its component at \(\SS\). One may then say that a map \(c \colon \SS \to r_!(\nu)\) \emph{exhibits \(\nu\) as the Spivak normal fibration of \(X\)} if the associated natural transformation
\[
c \otimes_{\SS} (-) \colon (-) \Rightarrow r_!(\nu) \otimes_{\SS} (-) \simeq r_!(\nu \otimes r^*(-))
\]
exhibits \(r_!(\nu \otimes (-))\) as right adjoint to \(r^*\). The map \(c\) is then called the \defi{Thom-Pontryagin} map of \(\nu\).

\begin{example}
If \(M\) is a closed smooth manifold then the underlying space of \(M\) is a Poincaré duality space. Furthermore, if \(\iota\colon M \hrar \RR^{N}\) is a smooth embedding with normal bundle \(E = TM^{\perp} \subseteq \iota^*T\RR^N\) and associated spherical fibration \(\xi_E\), and \(M \subseteq U \subseteq \RR^N\) is a chosen tubular neighborhood, then the Thom-Pontryagin collapse map
\[
\sph^{N} \to \RR^N/(\RR^N\setminus U) \simeq M^{E}
\]
induces a map of spectra
\[
\SS \to \Sig^{\infty-N}_+M^{E} \simeq r_!\Sig^{-N}\xi_E,
\]
that exhibits \(\Sig^{-N}\xi_E\) as the Spivak normal fibration of \(M\). Since \(E \oplus TM =\iota^*T\RR^N\) is a trivial \(N\)-dimensional vector bundle we can identify \(\Sig^{-N}\xi_E \simeq \xi_{TM}^{-1}\), where \(\xi_{TM}\) is the spherical fibration underlying the tangent bundle.
\end{example}

\begin{remark}
\label{remark:adjoint}%
Let \(X\) be a finite space with Spivak normal fibration \(\nu\) and Thom-Pontryagin map \(c\). Suppose that \(\E\) is a stable \(\infty\)-category, so that \(\E\) is canonically tensored over \(\Spaf\). We may then consider the restriction functor \(r^*\colon \E \to \Fun(X,\E)\), together with the associated left and right Kan extensions \(r_!,r_*\colon \Fun(X,\E) \to \E\) (which exist because \(X\) is finite and every stable \(\infty\)-category admits finite limits and colimits). The Thom-Pontryagin map then determines a natural transformation of functors from \(\E\) to itself whose component at \(A \in \E\) is
\[
c \otimes_{\SS} A \colon A \Rightarrow r_!(\nu) \otimes_{\SS} A \simeq r_!(\nu \otimes r^*A) .
\]
This natural transformation then exhibits \(r_!(\nu \otimes (-))\) as right adjoint to \(r^*\). Indeed, for \(A \in \E\) and \(B\colon X \to \E\) the composed map
\[
\map_X(r^*A,B) \to \map_{\E}(r_!(\nu \otimes r^*A),r_!(\nu \otimes B)) \to \map_{\E}(A,r_!(\nu \otimes B))
\]
is an equivalence, as can be seen by comparing it with the composed map
\[
\map_X(r^*\SS,B^A) \to \map(r_!\nu,r_!(\nu \otimes B^A)) \to \map(\SS,r_!(\nu \otimes B^A)),
\]
where \(B^A\) is the local system \(x \mapsto \map_{\E}(A,B_x)\), and we are using that
\[
\colim_{x \in X}\big(\nu_x \otimes \map_{\E}(A,B_x)\big) \simeq \map_{\E}\big(A,\colim_{x \in X}\nu_x \otimes B_x\big)
\]
since \(X\) is finite and each \(\nu_x\) is a finite spectrum.
\end{remark}

We now proceed to construct the visible signature of a Poincaré duality spaces \(X\).
Let \(p_1,p_2\colon X \times X \to X\) denote the two projections.
Applying Remark~\refone{remark:adjoint} to \(\E = \Spa_X\) we obtain that the natural map
\begin{equation}
\label{equation:unit-p-q}%
c \otimes_{\SS} L\colon L \to r_!(\nu) \otimes_\SS L \simeq (p_2)_!((p_1)^*\nu \otimes (p_2)^*L)
\end{equation}
for \(L \in \Spa_X\) acts as a unit exhibiting \((p_2)_!((p_1)^*\nu \otimes (-))\) as right adjoint to \((p_2)^*\). Since \((p_1)^*\nu\) is invertible tensoring with it is left and right inverse to tensoring with \((p_1)^*\nu^{-1}\), and hence the above map also acts as a unit exhibiting \((p_2)_!\colon \Spa_{X \times X} \to \Spa_X\) as right adjoint to \((p_1)^*\nu \otimes (p_2)^*(-) = \nu \boxtimes (-)\). We consequently obtain a canonical equivalence
\[
\map^{\boxtimes}(\nu,-) \simeq (p_2)_!(-)
\]
of functors \(\Spa_{X \times X} \to \Spa_X\). Evaluating at \(\Del_!\xi\) for some spherical fibration \(\xi\) we obtain a canonical equivalence
\begin{equation}
\label{equation:dual-is-nu}%
\Dual_{\xi}(\nu) \simeq \map^{\boxtimes}(\nu,\Del_!\xi) \simeq (p_2)_!\Del_!\xi \simeq \xi ,
\end{equation}
In other words, for every spherical fibration \(\xi\) on \(X\), the \(\xi\)-twisted duality switches between \(\xi\) and \(\nu\). Taking \(\xi = \nu\) we then get
\[
\Dual_{\nu}(\nu) \simeq \nu,
\]
that is, \(\nu\) is equivalent to its own \(\nu\)-twisted dual. We now claim that this equivalence is induced by a canonical Poincaré form \(q \in \Om^{\infty}\QF^{\vis}_{\nu}(\nu)\). To see this, note that by Remark~\refone{remark:useful-formula} and the equivalence~\eqrefone{equation:dual-is-nu} we have a \(\Ct\)-equivariant equivalence
\[
\Bil_{\xi}(\nu,\nu) \simeq \Bil_{\xi}(\Dual_{\SS_X}(\SS_X),\Dual_{\SS_X}(\SS_X)) \simeq \colim_{x \in X} \SS \otimes_{\SS} \SS \otimes_{\SS} \xi_x \simeq \colim_{x \in X} \xi_x \simeq r_!\xi,
\]
where the \(\Ct\)-action on the third term is by switching the two \(\SS\) factors, and the \(\Ct\)-action on the last two terms is \emph{trivial}. In particular, the \(\Ct\)-action on \(\Bil_{\xi}(\nu,\nu)\) is equivalent to a constant action.
Taking \(\xi = \nu\) we then obtain:

\begin{corollary}
\label{corollary:sym-poinc-object}%
The equivalence \(\nu \to \Dual_{\nu}\nu\) constructed above canonically refines to a self-dual equivalence, yielding in particular a symmetric Poincaré form \(q_{\sym} \in \QF^{\sym}_{\nu}(\nu) = \Bil_{\nu}(\nu,\nu)^{\hC}\).
\end{corollary}

We would like to lift the symmetric form \(q_{\sym}\) from \(\QF^{\sym}_{\nu}(\nu)\) to a visible form \(q_{\vis} \in \QF^{\vis}_{\nu}\). For this, we note that by using the second formula in Remark~\refone{remark:useful-formula} we obtain an equivalence
\[
\map_X(\nu,\xi) \simeq \Lin^{\vis}_{\xi}(\Dual_{\SS_X}(\SS_X)) \simeq r_!\xi,
\]
and so we can write the pullback square defining the visible Poincaré structure evaluated at \(\nu\) as
\[
\begin{tikzcd}
\QF^{\vis}_{\xi}(\nu) \ar[r] \ar[d] & r_!\xi \ar[d] \\
(r_!\xi)^{\hC} \ar[r] & (r_!\xi)^{\tC}
\end{tikzcd}
\]
where the \(\Ct\)-fixed points and Tate construction are performed with respect to the trivial action on \(r_!\xi\), and the right vertical map is given by the composite \(r_!\xi \to (r_!\xi)^{\hC} \to (r_!\xi)^{\tC}\).
Taking again \(\xi = \nu\) we then obtain:

\begin{corollary}
\label{corollary:visible-signature}%
Then the symmetric Poincaré form of Corollary~\refone{corollary:sym-poinc-object} canonically lifts to a visible Poincaré form \(q_{\vis} \in \QF^\vis_\nu(\nu)\).
\end{corollary}

The Poincaré object \((\nu,q_{\vis}) \in \Poinc(\Spa^{\om}_X,\QF^{\vis}_{\nu})\) is called the \defi{visible signature} of the Poincaré duality space \(X\).

\begin{remark}
\label{remark:visible-finitely-presented}%
Since \(X\) is finite any spherical fibration \(\xi\) is a finite colimit \(\xi \simeq \colim_{x \in X}x_!\xi_x\), and hence belongs to \(\Spaf_X\). It then follows that the visible signature \((\nu,q_{\vis})\) is naturally a Poincaré object in \((\Spaf,\QF^{\vis}_{\xi})\).
\end{remark}

%% file: Multiplicative.tex
In this section we will show that the tensor product of stable \(\infty\)-categories refines to give symmetric monoidal structures on the \(\infty\)-categories \(\Cath\) and \(\Catp\). We give a careful analysis of algebra objects in \(\Catp\), and use it to show that many examples of interest, such as the symmetric and genuine symmetric Poincaré structures on perfect derived categories of commutative rings, admit such an algebra structure. The significance of this fact is that for an algebra object in \(\Catp\), the \(\L\)-groups and the Grothendieck-Witt group inherit the structure of rings. We will show in \paperfour that this phenomenon extends to an \(\Einf\)-structure on the level of Grothendieck-Witt- and \(\L\)-theory \emph{spectra}, a structure which plays a key role in the study of these invariants in the commutative setting.

This section is organized as follows. In \S\refone{subsection:tensor-product} we define the tensor product of hermitian and Poincaré \(\infty\)-categories and give a formula for the linear and bilinear parts of the hermitian structure on the tensor product in term of the linear and bilinear parts of the individual terms. In \S\refone{subsection:monoidal-structure} we show that this operation organizes into a symmetric monoidal structure on the \(\infty\)-categories \(\Cath\) and \(\Catp\), such that the forgetful functors \(\Catp \to \Cath \to \Catx\) are symmetric monoidal. In \S\refone{subsection:symmetric-monoidal-poincare} we analyse what it means for a hermitian or Poincaré \(\infty\)-category to be an algebra with respect to this structure, and use this analysis in \S\refone{subsection:examples-monoidal} in order to identify various examples of interest of symmetric monoidal Poincaré \(\infty\)-categories. Finally, in \S\refone{subsection:GW-L-multiplicative} we show that the Grothendieck-Witt group and zero'th \(\L\)-group are symmetric monoidal functors on \(\Catp\), and hence their values on a given symmetric monoidal Poincaré \(\infty\)-category are rings. Similarly, we show that in this case the collection of all \(\L\)-groups acquires the structure of a graded-commutative ring.

\subsection{Tensor products of hermitian \(\infty\)-categories}
\label{subsection:tensor-product}%

In this first part we define the tensor product of hermitian \(\infty\)-categories on the level of objects and show that the tensor product of Poincaré \(\infty\)-categories is again Poincaré.

We begin by recalling the tensor product of stable \(\infty\)-categories due to Lurie.

\begin{construction}
\label{construction:monoida-catx}%
Applying the construction of~\cite[\S 4.8.1]{HA} with respect to the collection \(\cK\) of finite simplicial sets yields a symmetric monoidal structure on the \(\infty\)-category \(\Catrex\), whose objects are small \(\infty\)-categories with finite colimits and whose morphisms are functors which preserve finite colimits. For a pair of \(\infty\)-categories with finite colimits \(\C\) and \(\Ctwo\), their tensor product \(\C \otimes \Ctwo \in \Catrex\) is equipped with a functor \(\C \times \Ctwo \to \C \otimes \Ctwo\) which preserves finite colimits in each variable, and is initial with this property. Now the \(\infty\)-category \(\Catx\) embeds fully-faithful in \(\Catrex\), and is furthermore a reflective subcategory: a left adjoint to the inclusion \(\Catx \subseteq \Catrex\) is given by tensoring with \(\Spaf\).
To see this, observe that for an \(\infty\)-category \(\C\) with finite colimits, the tensor product
of \(\C\) and \(\Spaf\) yields an \(\infty\)-category with finite colimits \(\C \otimes \Spaf \in \Catrex\) which is a module over \(\Spaf\). Since \(0 \in \Spaf\) is a self-dual object it acts via a self-adjoint functor and hence \(\C \otimes \Spaf\) is pointed. In addition, since the object \(\Sig\SS \in \Spaf\) is invertible we get that suspension is invertible on \(\C \otimes \Spaf\) and so the latter is furthermore stable. To see that this operation gives a left adjoint to the inclusion \(\Catx \subseteq \Catrex\) we note that for every stable \(\infty\)-category \(\D\) the restriction functor
\[
\Funrex(\C \otimes \Spaf,\D) = \Funrex(\C,\Funrex(\Spaf,\D)) = \Funrex(\C,\Funx(\Spaf,\D)) \to \Funrex(\C,\D)
\]
is an equivalence by Lemma~\refone{lemma:finite-spectra-universal}. It then follows from~\cite[Proposition 4.1.7.4]{HA} that the symmetric monoidal structure \(\otimes\) descends to \(\Catx\). In particular, the unit of \(\Catx\) is given by \(\Spaf\), and if \(\C,\Ctwo\) are stable then \(\C \otimes \Ctwo\) is universal among stable \(\infty\)-categories receiving a bilinear functor \(\beta\colon \C \times \Ctwo \to \C \otimes \Ctwo\): given a stable \(\infty\)-category \(\D\), restriction along \(\beta\) induces an equivalence between the \(\infty\)-category of exact functors from \(\C \otimes \Ctwo \to \D\) and the \(\infty\)-category of bilinear functors \(\C \times \Ctwo \to \D\). We will refer to \(\Einf\)-algebra objects in \((\Catx)^{\otimes}\) as \defi{stably symmetric monoidal} \(\infty\)-categories. Concretely, this means a symmetric monoidal \(\infty\)-category whose underlying \(\infty\)-category is stable and the tensor product is exact in each variable.
\end{construction}

We now refine this construction to the level of hermitian \(\infty\)-categories.

\begin{construction}
\label{construction:tensor-product}%
For a pair of hermitian \(\infty\)-categories \((\C,\QF)\) and \((\Ctwo,\QFtwo)\), we define their \defi{tensor product}
\[
(\C,\QF) \otimes (\Ctwo,\QFtwo) := (\C \otimes \Ctwo,\QF\otimes \QFtwo)
\]
to be the hermitian \(\infty\)-category whose underlying stable \(\infty\)-category is the tensor product of the underlying stable \(\infty\)-categories, and whose hermitian structure
\[
\QF \otimes \QFtwo := \App_2\beta_!(\QF \boxtimes \QFtwo)\colon  \C\op \otimes \Ctwo{\op} \to \Spa
\]
is obtained by applying the 2-excisive approximation functor of Construction~\refone{construction:excisive-approx} to the left Kan extension along \(\beta \colon \C \times \Ctwo \to \C \otimes \Ctwo\) of the `external' tensor product
\[
\QF \boxtimes \QFtwo \colon \C{\op} \times {\Ctwo}{\op} \xrightarrow{\QF \times \QFtwo} \Spa \times \Spa \xrightarrow{\otimes} \Spa.
\]
Here we note that \(\beta_!(\QF \boxtimes \QFtwo)\) is already reduced since \(\QF \boxtimes \QFtwo\) is reduced and \(\beta\) preserve zero objects, and so Construction~\refone{construction:excisive-approx} is applicable to it.
\end{construction}

Note that the tensor product \((\C,\QF) \otimes (\Ctwo,\QFtwo) = (\C \otimes \Ctwo, \QF \otimes \QFtwo)\) carries by design a similar universal property to the tensor product \(\C \otimes \Ctwo\) of stable \(\infty\)-categories: for any hermitian \(\infty\)-category \((\Cthree, \QFthree)\), maps from the tensor product \((\C,\QF) \otimes (\Ctwo,\QFtwo)\) to \((\Cthree, \QFthree)\) correspond to bilinear maps \(b\colon \C \times \Ctwo \to \Cthree\)
together with a natural transformation \(\QF \boxtimes \QFtwo \to b^*\QFthree\), i.e.\ a  natural transformation in the square
\[
\begin{tikzcd}
\C\op \times \Ctwo{\op}  \ar[d,"b"'] \ar[r,"{\QF \times \QFtwo}"] &  \Spa \times \Spa \ar[d,"\otimes"]  \ar[ld,Rightarrow] \\
\Cthree{\op} \ar[r,"{\QFthree}"] & \Spa \ .
\end{tikzcd}
\]

Our next goal is to identify the linear and bilinear parts of the hermitian structure on \((\C,\QF) \otimes (\Ctwo,\QFtwo)\) in more explicit terms.
To this end, let \(\Lin_\QF,\Lin_{\QFtwo}\) denote the linear parts and \(\Bil_\QF,\Bil_{\QFtwo}\) the bilinear parts of the hermitian structures \(\QF\) and \(\QFtwo\), respectively.
The functor
\[
\Lin_\QF\boxtimes\Lin_\QFtwo \colon  \C\op\times \Ctwo{\op}\xrightarrow{\Lin_\QF\times\Lin_{\QFtwo}}\Spa\times\Spa\xrightarrow{\otimes}\Spa
\]
is then bilinear, and therefore extends along \(\beta\) to a linear functor \(\Lin_\QF\otimes\Lin_{\QFtwo} \colon  \C{\op}\otimes \Ctwo{\op}\rightarrow \Spa\) in an essentially unique manner. Similarly, the multilinear functor
\[
\Bil_\QF\boxtimes\Bil_{\QFtwo} \colon \C{\op}\times\C{\op}\times \Ctwo{\op}\times \Ctwo{\op}\xrightarrow{\Bil_\QF\times\Bil_\QFtwo}\Spa\times\Spa\xrightarrow{\otimes}\Spa
\]
extends to a
bilinear functor
\[\Bil_\QF\otimes\Bil_{\QFtwo}\colon (\C{\op}\otimes  \Ctwo{\op})\times(\C{\op}\otimes  \Ctwo{\op})\rightarrow \Spa .\]
The symmetric structures of \(\Bil_{\QF}\) and \(\Bil_{\QFtwo}\) then determine a \(\Ct\)-fixed structure on \(\Bil_\QF\boxtimes\Bil_{\QFtwo}\) with respect to the \(\Ct\)-action which permutes the two \(\C\op\)-coordinates and the two \(\Ctwo{\op}\)-coordinates. This structure then descends to a symmetric structure on \(\Bil_\QF\otimes\Bil_{\QFtwo}\) by its universal characterization.

\begin{proposition}
\label{proposition:tensor}%
For hermitian \(\infty\)-categories \((\C,\QF),(\Ctwo,\QFtwo)\) there is a canonical pullback square
\begin{equation}
\label{equation:pullback-square}%
\begin{tikzcd}
\QF \otimes \QFtwo \ar[r]\ar[d] & \Lin_\QF\otimes\Lin_\QFtwo\ar[d] \\
{\left[ (\Bil_\QF\otimes\Bil_\QFtwo)^{\Delta}\right]^{\hC}} \ar[r] & {\left[(\Bil_\QF\otimes\Bil_{\QFtwo})^\Delta\right]^{\tC}}
\end{tikzcd}
\end{equation}
of functors \(\C{\op}\otimes \Ctwo{\op}\rightarrow \Spa\).
In particular the linear part of \(\QF \otimes \QFtwo\) is
given by \(\Lin_\QF\otimes\Lin_{\QFtwo}\) and its symmetric bilinear part by \(\Bil_\QF\otimes\Bil_{\QFtwo}\).
\end{proposition}
\begin{proof}
By definition the functor \(\QF \otimes \QFtwo = \App_2\beta_!(\QF \boxtimes \QFtwo)\) is characterized by the fact that for every quadratic functor \(\QFthree\colon  \C\op\otimes \Ctwo\op\to\Spa\) the natural map
\begin{equation*}
\Nat(\QF \otimes \QFtwo ,\QFthree) \to \Nat(\QF\boxtimes \QFtwo,\beta^*\QFthree)
\end{equation*}
is an equivalence. Using this universal mapping property the commutative square~\eqrefone{equation:pullback-square} is then obtained from the external square in the diagram
\[
\begin{tikzcd}
\QF \boxtimes \QFtwo \ar[r]\ar[d] & \Lin_\QF\boxtimes\Lin_\QFtwo \ar[d] \\
{\left[(\Bil_\QF)^{\Delta}\right]^{\hC}\boxtimes\left[(\Bil_\QFtwo)^{\Delta}\right]^{\hC}} \ar[r] \ar[d] & {\left[(\Bil_\QF)^{\Delta}\right]^{\tC}\boxtimes\left[(\Bil_\QFtwo)^{\Delta}\right]^{\tC}} \ar[d] \\
{\left[ (\Bil_\QF\boxtimes\Bil_\QFtwo)^{\Delta}\right]^{\hC}} \ar[r] & {\left[(\Bil_\QF\boxtimes\Bil_{\QFtwo})^\Delta\right]^{\tC}} \ ,
\end{tikzcd}
\]
where the top square is obtained by taking the external product of the classifying squares of \(\QF\) and \(\QFtwo\), and the bottom square witnesses the lax symmetric monoidal structure of the homotopy fixed points functor and the projection to the Tate construction. To show that the resulting square~\eqrefone{equation:pullback-square} is a pullback square we need to show that the induced map
\begin{equation}
\label{equation:tensorproduct-map}%
\Nat\big((\Lin_\QF\otimes\Lin_{\QFtwo})\times_{[(\Bil_\QF\otimes\Bil_\QFtwo)^{\Del}]^{\tC}} [(\Bil_\QF\otimes\Bil_{\QFtwo})^{\Del}]^{\hC} ,\QFthree\big) \to \Nat(\QF\boxtimes \QFtwo,\beta^*\QFthree)
\end{equation}
is an equivalence for any quadratic functor \(\QFthree\). Let us analyse both sides of the map \eqrefone{equation:tensorproduct-map}. We start with the following claim (see \S\refone{subsection:classification} for the terminology of homogeneous and cohomogeneous and their basic properties):
\begin{quote}
\textbf{Claim 1:} \(\Nat(\QF\boxtimes \QFtwo,\beta^*\QFthree) = 0\) if either \(\QF\) is exact and \(\QFthree\) is cohomogeneous or \(\QF\) is homogeneous and \(\QFthree\) is exact.
\end{quote}
To see this claim we note that for a fixed object \(c' \in \Ctwo\) we have that the space of natural transformations
\[
\Nat(\QF(-) \otimes \QFtwo(c') ,\QFthree(\beta(-,c')))
\]
(natural in \((-)\)) vanishes under these assumptions since left hand functor is exact (resp.\ homogeneous) and the right hand functor is cohomogeneous (resp.\ exact). But the space
\(\Nat(\QF\boxtimes \QFtwo,\beta^*\QFthree)\) can be written as a limit of these spaces over the twisted arrow category of \(\Ctwo\) so that the claim follows.
We also have the following claim
\begin{quote}
\textbf{Claim 2:} \(\Nat\big((\Lin_\QF\otimes\Lin_{\QFtwo})\times_{[(\Bil_\QF\otimes\Bil_\QFtwo)^{\Del}]^{\tC}} [(\Bil_\QF\otimes\Bil_{\QFtwo})^{\Del}]^{\hC} ,\QFthree\big) = 0\) if either \(\QF\) is exact and \(\QFthree\) is cohomogeneous or \(\QF\) is homogeneous and \(\QFthree\) is exact.
\end{quote}
which follows by the same argument as above since under the assumptions the pullback is either \(\Lin_\QF \otimes \Lin_\QFtwo\) or \([(\Bil_\QF \otimes \Bil_\QFtwo)^{\Del}]_{\hC}\).

Together the last two claims show that the map  \eqrefone{equation:pullback-square} is an equivalence under the assumption that either \(\QF\) is exact and \(\QFthree\) is cohomogeneous or \(\QF\) is homogeneous and \(\QFthree\) is exact. Now \(\QF\) can be fit in an exact sequence between a homogeneous functor and an exact functor, and
\(\QFthree\) can be fit in an exact sequence between an exact and cohomogeneous functors.
Since both sides of \eqrefone{equation:pullback-square} are natural and exact in \(\QF,\QFtwo\) and \(\QFthree\) we can thus assume without loss of generality that \(\QF\) and \(\QFthree\) are both exact or \(\QF\) is homogeneous and \(\QFthree\) is cohomogeneous. Since everything is symmetric in \(\QF\) and \(\QFtwo\) we can also make the same reduction in this variable so that we only need to show the fact that \eqrefone{equation:pullback-square} is a pullback under the following assumptions:
 \begin{quote}
\textbf{Claim 3:} Either all three functors \(\QF, \QFtwo, \QFthree\) are exact or \(\QF\) and \(\QFtwo\) are homogeneous and \(\QFthree\) is cohomogeneous.
\end{quote}
In the first case the statement unwinds to the universal property of \(\C \otimes \Ctwo\) and in the second case it unwinds (also using the universal property) to the statement that maps from a homogeneous functor to a cohomogeneous functor are equivalent to maps between the associated symmetric bilinear functors.
\end{proof}

\begin{corollary}
\label{corollary:poinc}%
If \((\C,\QF)\) and \((\Ctwo,\QFtwo)\) are Poincaré with duality functors \(\Dual_{\QF},\Dual_{\QFtwo}\) then \((\C,\QF) \otimes (\Ctwo,\QFtwo)\) is Poincaré with duality functor \(\Dual_{\QF}\otimes \Dual_{\QFtwo}\colon (\C \otimes \Ctwo)\op = \C\op \otimes \Ctwo{\op} \to \C \otimes \Ctwo\).
\end{corollary}
\begin{proof}
We get from Proposition~\refone{proposition:tensor} that the cross effect of the quadratic functor on \(\C \otimes \Ctwo\) is given by \(\Bil_\QF \otimes \Bil_\QFtwo\), which coincides with the left Kan extension of \(\Bil_{\QF} \boxtimes \Bil_{\QFtwo}\) along the map
\[
\beta \times \beta \colon (\C\op \times \Ctwo{\op}) \times (\C\op \times \Ctwo{\op}) \to (\C\op \otimes \Ctwo{\op}) \times (\C\op \otimes \Ctwo{\op}),
\]
where
\[
[\Bil_{\QF} \boxtimes \Bil_{\QFtwo}](\x,\x',\y,\y') = \Bil_{\QF}(\x,y)\otimes \Bil_{\QFtwo}(x',y') = \map_{\C}(\x,\Dual_{\QF}\y) \otimes \map_{\Ctwo}(\x',\Dual_{\QFtwo}\y').
\]
We want to show that this is represented by the functor \(\Dual_{\QF} \otimes \Dual_{\QFtwo}\).
Now the left Kan extension along \(\beta \times \beta\) can be computed by composing left Kan extensions along \(\beta \times \id\) and \(\id \times \beta\). Then for \(\y \in \C,\y'\in \Ctwo\) we then have
\[
[(\beta \times \id)_!\Bil_{\QF} \boxtimes \Bil_{\QFtwo}]|_{\C \otimes \Ctwo \times \{y\} \times \{y'\}} = \beta_![\Bil_{\QF}(-,y)\otimes \Bil_{\QFtwo}(-,y')] =
\]
\[
\beta_![\map_{\C}(-,\Dual_{\QF}\y) \otimes \map_{\Ctwo}(-,\Dual_{\QFtwo}\y')] = \map_{\C \otimes \Ctwo}(-,\Dual_{\QF}\y \otimes \Dual_{\QFtwo}\y'),
\]
and the left Kan extension along \(\beta\) of the functor \((\y,\y') \mapsto \Dual_{\QF}\y \otimes \Dual_{\QFtwo}\y'\) is \(\Dual_{\QF} \otimes \Dual_{\QFtwo}\).
\end{proof}

\begin{proposition}
Recall the universal Poincaré \(\infty\)-category \((\Spaf,\QF^{\uni})\) of \S\refone{subsection:universal}. Then we have a natural equivalence
\[
(\C, \QF) \otimes (\Spaf, \QF^{\uni}) \simeq (\C, \QF)
\]
for \((\C,\QF) \in \Cath\).
\end{proposition}
\begin{proof}
We first note that \(\Spaf\) is the unit with respect to the tensor product of stable \(\infty\)-categories (see Construction~\refone{construction:monoida-catx}), and that \((\Spaf)\op \simeq \Spaf\) through Spanier-Whitehead duality. By definition, the linear part of \(\QF^{\uni}\colon (\Spaf)\op \to \Spa\) is Spanier-Whitehead duality \(\Dual\) and the bilinear part is the composite
\[
(\Spaf){\op} \times (\Spaf){\op} \xrightarrow{\otimes_\SS} (\Spaf){\op} \xrightarrow{\Dual} \Spaf,
\]
which also corresponds to Spanier Whitehead duality \((\Spaf){\op} \to \Spaf\) under the equivalence \((\Spaf){\op} \otimes (\Spaf){\op} \simeq (\Spaf)\op\). Now using Proposition~\refone{proposition:tensor} the claim is reduced to the statement that for any exact functor \(\Lin\colon \C\op \to \Spa\) the functor
\[
\C\op = (\C \otimes \Spaf)\op \to \C\op \otimes (\Spaf){\op} \xrightarrow{\Lin\otimes \Dual}\Spa\otimes\Spa\xrightarrow{\otimes}\Spa
\]
is equivalent to \(\Lin\). Indeed, the composite \(\C\op = (\C \otimes \Spaf)\op \to \C\op \otimes (\Spaf){\op}\) sends \(c \in \C\op\) to \(c \otimes \SS\).
\end{proof}

\begin{example}
\label{example:tensor-of-rings}%
Let \(A,B\) be \(\Eone\)-ring spectra equipped with modules with genuine involutions \((M_A,N_A,\alp)\) and \((M_B,N_B,\beta)\) respectively (see \S\refone{subsection:genuine-modules}).
As discussed in the proof of Theorem~\refone{theorem:classification-genuine-modules}, it follows from~\cite[Theorem 4.8.5.16 and Remark 4.8.5.19]{HA} that the exact functor
\begin{equation}
\label{equation:tensor-A-B}%
\Modf{A} \otimes \Modf{B} \to \Modf{A \otimes_{\SS} B},
\end{equation}
induced by the bilinear functor \((X,Y) \mapsto X \otimes Y\), yields an equivalence \(\Mod_{A} \otimes_{\Spa} \Mod_{B} \xrightarrow{\simeq} \Mod_{A \otimes_{\SS} B}\) upon passing to Ind-categories, where \(\otimes_{\Spa}\) denotes the tensor product of stable presentable \(\infty\)-categories. In particular, the functor~\eqrefone{equation:tensor-A-B} is necessarily a dense full inclusion, which is therefore an equivalence since its image contains \(A \otimes B\) and its target contains no proper stable subcategories with that property. More generally, given subgroups \(c \subseteq \K_0(\Modp{A})\) and \(d \subseteq \K_0(\Modp{B})\) we obtain an equivalence
\[ \Mod^c_A \otimes \Mod^d_{B} \xrightarrow{\simeq} \Mod^e_{A \otimes_{\SS} B} \]
where \(e \subseteq \K_0(\Modp{A \otimes_{\SS} B})\) is the image of \(c \otimes d\) under the induced map \(\K_0(\Modp{A}) \otimes \K_0(\Modp{B}) \to \K_0(\Modp{A \otimes_{\SS} B})\).
We note that in general \(e\) might fail to be all of \(\K_0(\Modp{A \otimes_{\SS} B})\), even if \(c\) and \(d\) are all of \(\K_0(\Modp{A})\) and \(\K_0(\Modp{B})\), respectively. It then follows from Proposition~\refone{proposition:tensor} below that
\[
(\Mod^c_{A},\QF^{\alp}_{M_A}) \otimes (\Mod^d_{B},\QF^{\beta}_{M_B}) \simeq (\Mod^e_{A \otimes_{\SS} B},\QF^{\alp \otimes_{\SS} \beta}_{M_A \otimes_{\SS} M_B})
\]
where the reference map on the right hand side is the composite
\[
\alp \otimes_{\SS} \beta\colon N_A \otimes_{\SS} N_B \to M_A^{\tC} \otimes_{\SS} M_B^{\tC} \to (M_A \otimes_{\SS} M_B)^{\tC},
\]
obtained using the lax monoidal structure of the Tate construction. \end{example}

\subsection{Construction of the symmetric monoidal structure}
\label{subsection:monoidal-structure}%

In this section we will show that the notion of tensor product constructed in \S\refone{subsection:tensor-product} above can be enhanced to symmetric monoidal structures on \(\Cath\) and on \(\Catp\).  The construction is somewhat technical and can be skipped on a first read.

\begin{construction}
For an \(\infty\)-category \(\D\) we will denote by
\[
(\Cat)_{//\D} \to \Cat
\]
the cartesian fibration classified by the functor
\[
\Cat\op \to \Cat \qquad \C \mapsto \Fun(\C, \D)\ .
\]
We will refer to \((\Cat)_{//\D}\) as the \emph{lax slice over \(\D\)}. The objects of \((\Cat)_{//\D}\) are given by functors \(\C \to \D\) and the morphisms by diagrams
\[
\begin{tikzcd}
[row sep=1.3ex]
\C \ar[dd,swap, "f"] \ar[rd,swap, "p"',""{name=p,right}] & \\
& \D \\
\C' \ar[ru,swap, "q"]\arrow[Rightarrow, from=p, shorten <= 1.5ex, shorten >= 1.5ex]  & \\
\end{tikzcd}
\]
filled by a non-invertible 2-cell \(p \Rightarrow qf\).
The actual slice \((\Cat)_{/\D}\) is a non-full subcategory of
\((\Cat)_{//\D}\) which contains all objects but only those 1-morphisms for which the natural transformation \(p  \Rightarrow qf\) is an equivalence.
\end{construction}

\begin{remark}
\label{remark:mapping-propoerty}%
The \(\infty\)-category \((\Cat)_{//\D}\) can be characterized by the following universal mapping property:
the data of a functor from \(\E\) to \((\Cat)_{//\D}\) is equivalent to the data of a diagram in \(\Cat\) of the form
\begin{equation}
\label{equation:ECD}%
\begin{tikzcd}
\C \ar[r] \ar[d,"p"'] & \D \\
\E
\end{tikzcd}
\end{equation}
with \(p\) a cocartesian fibration.
In this description functors from \(\E\) to the actual slice \((\Cat)_{/\D}\) correspond to diagrams as above where the functor \(\C \to \D\) send \(p\)-cocartesian lifts to equivalences in \(\D\). This description also uniquely determines \((\Cat)_{//\D}\) since it describes the represented functor \(\Ho(\Cat)\op \to \Set\).
\end{remark}

\begin{lemma}
\label{lemma:symm}%
Let \(\D\) be a symmetric monoidal \(\infty\)-category. Then  \((\Cat)_{//\D}\) admits a symmetric monoidal refinement \((\Cat)_{//\D}^\otimes\) with the following properties:
\begin{enumerate}
\item
\label{item:product}%
The tensor product of \(f\colon \C \to \D\) and \(g\colon \C' \to \D\) in \((\Cat)_{//\D}\) is given by the composite
\[
\C \times \C' \xrightarrow{f \times g} \D \times \D \xrightarrow{\otimes} \D,
\]
and the tensor unit by the functor \(\pt \to \D\) corresponding to the tensor unit of \(\D\).
\item
\label{item:forgetful}%
The forgetful functor \((\Cat)_{//\D} \to \Cat\) admits a symmetric monoidal refinement \((\Cat)_{//\D}^\otimes \to \Cat^\times\) with respect to the cartesian symmetric monoidal structure on \(\Cat\).
\item
\label{item:algebras}%
For any \(\infty\)-operad \(\cO^\otimes\), the space of
\(\cO\)-algebras in \((\Cat)_{//\D}^\otimes\) is naturally equivalent to the space of diagrams of \(\infty\)-operads
\[
\begin{tikzcd}
\C^\otimes \ar[d,"p"'] \ar[r] & \D^\otimes \\
\cO^\otimes
\end{tikzcd}
\]
where \(p\) is a cocartesian fibration of \(\infty\)-operads (see~\cite[Definition 2.1.2.13]{HA}). In other words, a pair consisting of an \(\cO\)-monoidal \(\infty\)-category \(\C\) together with a lax \(\cO\)-monoidal functor from \(\C\) to \(\D\), where the latter is consider as an \(\cO\)-monoidal
\(\infty\)-category by pullback along the terminal map of \(\infty\)-operads \(\cO \to \Einf\).
\end{enumerate}
\end{lemma}

\begin{remark}
\label{remark:property-determines}%
Property~\refoneitem{item:algebras} of Lemma~\refone{lemma:symm}
determines the symmetric monoidal \(\infty\)-category \((\Cat)_{//\D}^\otimes\) uniquely. Indeed, a symmetric monoidal \(\infty\)-category is uniquely determined by its underlying \(\infty\)-operad, and Property~\refoneitem{item:algebras} determines the functor
\[
\Ho(\Op)\op \to \Set
\]
represented by the underlying \(\infty\)-operad \((\Cat)_{//\D}^\otimes\).
One can in fact show that the first two properties of the lemma are direct consequences of the third (see the arguments in the proof of Lemma~\refone{lemma:symm} below).
\end{remark}

\begin{proof}[Proof of Lemma~\refone{lemma:symm}]
For a fixed symmetric monoidal \(\infty\)-category \(\D\) with underlying \(\infty\)-operad \(\D^\otimes\), we consider the \(\infty\)-category \(\mathcal{X}\) whose objects are given by diagrams of \(\infty\)-operads of the form
\begin{equation}
\label{equation:OCD}%
\begin{tikzcd}
\C^\otimes \ar[d,"p"'] \ar[r] & \D^\otimes\\
\cO^\otimes
\end{tikzcd}
\end{equation}
where \(p\) is a cocartesian fibration of \(\infty\)-operads. The morphisms in \(\mathcal{X}\) are those maps of such diagrams which are the identity on \(\D^\otimes\) and preserve \(p\)-cocartesian arrows in the \(\C^\otimes\) component. Projecting to the \(\cO^\otimes\)-component defines a functor \(\mathcal{X} \to \Op\)
which is a cartesian fibration classified by some functor \(\Op\op \to \Cat\). Postcomposing the latter with the groupoid core functor we obtain a functor
\[
q\colon \Op\op \to \Sps
\]
sending an \(\infty\)-operad \(\cO\) to the space consisting of pairs of a cocartesian fibration \(p\colon \C^\otimes \to \cO^\otimes\) of \(\infty\)-operads together with a map \(\C^\otimes\to \D^\otimes\) of \(\infty\)-operads. We now claim that the functor \(q\) is representable by an \(\infty\)-operad. To see this we use that \(\Op\) is a presentable \(\infty\)-category by \cite[Section 2.1.4]{HA} so that we have to check that the functor \(q\) preserves limits by \cite[Proposition 5.5.2.2]{HTT}. This can be seen as follows: a cocartesian fibration \(\C^\otimes \to \cO^\otimes = \colim_i \cO_i^\otimes\) over a colimit of \(\infty\)-operads  is classified by a map of \(\infty\)-operads \(\chi\colon \colim_i \cO_i^\otimes \to \Cat^\times\) where \(\Cat^\times\) is equipped with the structure of an \(\infty\)-operad induced by the cartesian symmetric monoidal structure on \(\Cat\) (see \cite[Remark 2.4.2.6]{HA}). The space of such functors is thus a limit of the spaces of maps of \(\infty\)-operads \(\chi_i\colon \cO_i^\otimes \to \Cat^\times\). In particular the space of cocartesian fibrations over \(\colim_i \cO_i^\otimes\) is the limit of the spaces of cocartesian fibrations over \(\cO_i^\otimes\).
In addition, for every cocartesian fibration \(\C^{\otimes} \to \cO^{\otimes}\) of \(\infty\)-operads the natural map
\begin{equation}
\label{equation:colim-operads}%
\colim_i \C_i^{\otimes} \to \C^\otimes
\end{equation}
is an equivalence of \(\infty\)-operads,
where \(\C_i^\otimes := \C^\otimes\times_{\cO^\otimes} {\cO_i^\otimes}\) is the corresponding fibre product (computed in \(\Op\)).
This follows from the fact that the functor
\[
\C^\otimes \times_{\cO^\otimes} - \colon (\Op)_{/\cO^\otimes} \to (\Op)_{/\cO^\otimes}
\]
commutes with colimits of \(\infty\)-operads since it has a right adjoint given by the relative Day convolution \(\Fun^{\cO}(\C, -)^\otimes\), see~\cite[Construction 2.2.6.7 and Remark 2.2.6.8]{HA}.
The colimit description of \eqrefone{equation:colim-operads} then implies that the space of maps \(\C^\otimes \to \D^\otimes\) of \(\infty\)-operads is given the limit of the space of maps \(\C_i^\otimes \to \D^\otimes\). Together this shows that the functor \(q\) preserves limits and is thus representable.
We denote the representing object by \((\Cat^\otimes)_{//\D}\).
We note that by Remark~\refone{remark:mapping-propoerty} the underlying \(\infty\)-category \(\big((\Cat^{\otimes})_{//\D}\big)_{\left<1\right>}\) identifies with \((\Cat)_{//\D}\); indeed, when \(\cO^{\otimes}\) is the image of an \(\infty\)-category \(\E\) under the full inclusion \(\Cat \subseteq \Op\), the data of a diagram as in~\eqrefone{equation:OCD} with \(p\) a cocartesian fibration of \(\infty\)-operads reduces to that of a diagram of the form~\eqrefone{equation:ECD}, with \(p\) a cocartesian fibration of \(\infty\)-categories. To show \((\Cat^{\otimes})_{//\D}\) is a symmetric monoidal structure on \((\Cat)_{//\D}\)
we will need a more explicit description of the multi-mapping spaces in \((\Cat^\otimes)_{//\D}\).

Let \(C_n\) (sometimes called the \(n\)-corolla) be the \(\infty\)-operad freely generated by a single \(n\)-ary operation \(x_1,\ldots, x_n \to x\) with
colours \(x_1,...,x_n\) and \(x\)%
\footnote{This is in fact the nerve of an ordinary operad which can for example be seen using the theory of dendroidal sets, but we shall not need this fact here.}.%
The space of maps \(C_n^\otimes \to  (\Cat^\otimes)_{//\D}\) is then by the defining property of \((\Cat^\otimes)_{//\D}\) the classifying space of pairs of a \(C_n\)-monoidal \(\infty\)-category \(\C\) and a lax \(C_n\)-monoidal functor \(\C \to \D\). As the \(\infty\)-operad \(C_n\) is free a \(C_n\)-monoidal \(\infty\)-category is simply given by a sequence \(\{\C_1,...,\C_n;\C\}\) of \(\infty\)-categories together with a functor \(\alp\colon \C_1 \times \ldots\times \C_n \to \C\), and a lax \(C_n\)-monoidal functor from this to \(\D\) corresponds to a collection of functors \(\{f_1\colon \C_1 \to \D, \ldots, f_n\colon \C_n \to \D; f\colon \C \to \D\}\) together with a transformation in the square
\[
\begin{tikzcd}
\C_1 \times \ldots \times \C_n \ar[d,"\alp"] \ar[rr,"f_1 \times \ldots \times f_n"] && \D \times \ldots \times \D \ar[d,"\otimes"] \arrow[lld,Rightarrow, shorten <= 8ex, shorten >= 8ex] \\
\C \ar[rr,"f"] && \D \ .
\end{tikzcd}
\]
As a result we find that the corresponding multi-mapping space, which can be identified with the pullback
\[
\begin{tikzcd}
\Mul_{(\Cat^\otimes)_{//D}}(f_1,...,f_n; f) \ar[d]\ar[r] & \Map_{\Op}\left(C_n, (\Cat^\otimes)_{//D}\right) \ar[d]^{\prod_i \ev_{x_i} \times \ev_x}\\
\pt \ar[r,"{(f_1,...,f_n, f)}"] & \prod_{n+1} \Map_{\Op}\left(\mathcal{T}\mathrm{riv}^\otimes,(\Cat^\otimes)_{//\D}\right) \ ,
\end{tikzcd}
\]
is given by the spaces of maps in \((\Cat)_{//\D}\) from the object
\begin{equation}
\label{equation:corep-multi}%
\C_1 \times \ldots \times \C_n \xrightarrow{f_1 \times \ldots \times f_n} \D \times \ldots \times \D \xrightarrow{\otimes} \D
\end{equation}
to the object \(f\colon \C \to \D\).
In particular, this multi-mapping space is corepresented by the object \eqrefone{equation:corep-multi}
and so the \(\infty\)-operad \((\Cat)_{//\D}^\otimes\) is
corepresentable in the sense of \cite[Definition 6.2.4.3]{HA},  that is, the functor \((\Cat)_{//\D}^\otimes \to \mathrm{Fin}_*\) is a locally cocartesian fibration. To see that it is actually cocartesian, i.e.\ \((\Cat)_{//D}^\otimes\) is symmetric monoidal, we have to additionally verify that the induced maps from \cite[Example 6.2.4.9]{HA} are equivalences. This is however clear in the case at hand.

Finally, let us verify that \((\Cat)_{//\D}^\otimes\) satisfies the required Properties \refoneitem{item:product}-\refoneitem{item:algebras}. Indeed, Property~\refoneitem{item:algebras} is satisfied by construction and Property~\refoneitem{item:product} follows from the explicit description of multi-mapping spaces above.
To prove Property~\refoneitem{item:forgetful} we note that the \(\infty\)-operad \(\Cat^\times\) represents the functor \(\Op\op \to \Sps\) which sends an \(\infty\)-operad \(\cO^{\otimes}\) to the space of cocartesian fibrations \(\E^{\otimes} \to \cO^{\otimes}\) (this follows from \cite[Remark 2.4.2.6]{HA}). The functor \((\Cat)_{//D} \to \Cat\) which forgets the map refines to a transformation of represented functors \(\Op\op \to \Sps\) (again given by forgetting the map to \(\D^\otimes\)). Thus we get a lax symmetric monoidal structure on the functor \((\Cat)_{//\D} \to \Cat\) and by the description of the tensor product given above it follows that this functor is actually symmetric monoidal as opposed to merely lax symmetric monoidal.
\end{proof}

\begin{remark}
It is also possible to give a direct construction of the symmetric monoidal \(\infty\)-category \((\Cat)_{//\D}^\otimes\) as follows.
Let \(\LaxAr \subseteq (\Cat)_{/\Del^1}\) be the full subcategory spanned by the cartesian fibrations \(\M \to \Del^1\). We note that such a cartesian fibration encodes the data of an functor \((\Del^1)\op \to \Cat\), corresponding to an arrow \(\M_1 \to \M_0\) in \(\Cat\), where \(\M_i := \M \times_{\Del^1} \Del^{\{i\}}\) is the fibre of \(\M\) over \(i \in \{0,1\}\). By definition the morphisms in \(\LaxAr(\Cat)\) simply correspond to functors \(\M \to \M'\) over \(\Del^1\), and these are not required to preserve cartesian edges. As a result, morphisms in \(\LaxAr\) correspond to lax natural transformations of arrows, that is, to squares
which commute up to a specified transformation.
We then endow \(\LaxAr\) with the cartesian monoidal structure \(\LaxAr^{\times}\), which is simply given by fibre product over \(\Del^1\) (since cartesian fibrations are closed under fibre products) and we define \((\Cat)_{//\D}\) to be the fibre of the functor
\[
f_0\colon\mathrm{LaxAr}(\Cat) \to  \Cat \quad\quad [M \to \Delta^1] \mapsto M_0
\]
over \(\D\). Since \(f_0\) is product preserving and \(\D\) is an \(\Einf\)-monoid object in \(\Cat\) the fibre \((\Cat)_{//\D}\) inherits a symmetric monoidal structure, which we denote by \((\Cat)_{//\D}^{\otimes}\).
We now claim that the underlying \(\infty\)-operad of \((\Cat)_{//\D}^{\otimes}\) represents the same functor described in Property~\refoneitem{item:algebras}, and hence identifies with the construction given above. To see this, let us first identify the functor represented by \(\LaxAr^{\times}\).
For \(\cO\) an \(\infty\)-operad, \(\cO\)-algebra objects in \(\LaxAr^{\times}\) correspond to \(\cO\)-monoids,
which are simply functors \(\cO^{\otimes} \to \LaxAr\) in which certain diagrams are cartesian.
But the \(\infty\)-category of functors \(\cO^{\otimes} \to \LaxAr\) embeds in the \(\infty\)-category of functors \(\cO^{\otimes} \to (\Cat)_{/\Delta^1}\), and the latter correspond via unstraightening to cocartesian fibrations \(p\colon\E \to \cO^{\otimes}\) equipped with a map \(\E \to \Delta^1\) which sends \(p\)-cocartesian edges to equivalences. The condition that the associated functor
\(\cO^{\otimes} \to  (\Cat)_{/\Delta^1}\) lands in \(\LaxAr\) corresponds in these terms to the condition that for every \(x\) in \(\cO^{\otimes}\) the restricted map \(\E_x \to \Delta^1\) is a cartesian fibration. By the dual of \cite{HTT} Corollary 4.3.1.15
this is equivalent to saying that \(\E \to \Delta^1\) is a cartesian fibration whose cartesian edges all map to equivalences in \(\cO^{\otimes}\). Straightening over \(\Delta^1\), this data is equivalent to that of a map \(\E_1 \to \E_0\) of \(\infty\)-categories over \(\cO^{\otimes}\) such that the maps \(\E_0 \to \cO^{\otimes}\) and \(\E_1 \to \cO^{\otimes}\) are cocartesian fibrations. The monoid condition is then equivalent to the condition that for \(i=0,1\) the cocartesian fibration \(\E_i \to \cO^{\otimes}\) exhibits \(\E_i\) as an \(\cO\)-monoidal \(\infty\)-category and that the functor \(\E_1 \to \E_0\) preserves inert maps. We hence get that the data of an \(\cO\)-algebra object in \(\LaxAr^{\times}\) is equivalent to that of a pair of
\(\cO\)-monoidal \(\infty\)-categories \(\E_0,\E_1\) equipped with a lax \(\cO\)-monoidal functor \(\E_1 \to \E_0\). We may then conclude that the data of an \(\cO\)-algebra object in \((\Cat)_{//D}\) is equivalent to that of an \(\cO\)-monoidal \(\infty\)-category \(\E_1\) equipped with a lax \(\cO\)-monoidal functor \(\E_1 \to \D\).
\end{remark}

We will now apply the construction of Lemma~\refone{lemma:symm} to the category \(\D = \Spa\) of spectra, equipped with its symmetric monoidal structure given by the tensor product of spectra, and form the pullback along the autoequivalence \((-)\op\colon \Cat \to \Cat\). More precisely we define a symmetric monoidal \(\infty\)-category \(\left((\Cat)_{\mop//\Spa}\right)^\otimes\) as the pullback
\[
\begin{tikzcd}
\left((\Cat)_{\mop//\Spa}\right)^\otimes \ar[r,"{\simeq}"] \ar[d] & \left((\Cat)_{// \Spa}\right)^\otimes \ar[d] \\
\Cat^\times \ar[r,"{(-)\op}","{\simeq}"'] & \Cat^\times
\end{tikzcd}
\]
Objects of this symmetric monoidal \(\infty\)-category are given by pairs
\((\C, \QF)\) consisting of an \(\infty\)-category \(\C\) and a functor \(\QF\colon \C\op \to \Spa\). Morphisms
\[
(\C_1,\QF_1) \otimes \ldots \otimes (\C_n, \QF_n) \to (\C', \QFtwo)
\]
in \(\left((\Cat)_{\mop//\Spa}\right)^\otimes\) are given by
pairs \((f,\eta)\) of a functor \(f\colon \C_1 \times \ldots \times \C_n  \to \C'\)
and a natural transformation \(\eta\colon \QF_1 \boxtimes \ldots \boxtimes \QF_n \Rightarrow \QFtwo \circ f\op\). This description also holds for \(n = 0\).

\begin{construction}
\label{construction:symmetric-structure}%
We define the \(\infty\)-operad \((\Cath)^\otimes\) as the suboperad of \(\left((\Cat)_{\mop//\Spa}\right)^\otimes\) spanned by those objects \(((\C_1,\QF_1),...,(\C_n,\QF_n)) \in \left((\Cat)_{\mop//\Spa}\right)^\otimes\) such that each \(\C_i\) is stable and each \(\QF_i\) is quadratic, and those maps \((f,\eta)\) as above such that \(f\colon \C_1 \times \ldots \times \C_n \to \Ctwo\) is exact in each variable.

We note that the underlying \(\infty\)-category of \((\Cath)^\otimes\) is indeed given by \(\Cath\) since it is, essentially by definition, the Grothendieck construction of the functor \(\C \mapsto \Funq(\C)\) which is a subcategory of the lax slice \((\Cat)_{\mop//\Spa}\).
\end{construction}

By definition we have that the composed lax symmetric monoidal functor
\[
(\Cath)^{\otimes} \to \left((\Cat)_{\mop//\Spa}\right)^\otimes \to \Cat^\times
\]
factors through the suboperad
\[
(\Cath)^{\otimes} \to (\Catx)^{\otimes} \subseteq \Cat^\times
\]
spanned by the tuples of stable \(\infty\)-categories and tuples of functors which are exact in each variable.

\begin{theorem}
\label{theorem:tensorpoincare}%
\
\begin{enumerate}
\item
\label{item:thm-one}%
The \(\infty\)-operad \((\Cath)^\otimes\) is symmetric monoidal with tensor product given by the tensor product of Construction~\refone{construction:tensor-product} and tensor unit given by \((\Spaf, \QF^{\uni})\).
\item
\label{item:thm-two}%
The resulting map \(p^{\otimes}\colon (\Cath)^\otimes \to (\Catx)^\otimes\) is a cocartesian fibration of \(\infty\)-operads and in
particular a symmetric monoidal functor.
\item
\label{item:thm-three}%
This symmetric monoidal structure on \(\Cath\) restricts to a symmetric monoidal structure on the subcategory \(\Catp\).
\end{enumerate}
\end{theorem}
\begin{proof}
By design, the hermitian \(\infty\)-category
\((\C_1,\QF_1) \otimes (\C_2,\QF_2)\) of Construction~\refone{construction:tensor-product} corepresents the binary multi-mapping space functor
\[
\Mul_{\Cath}\left ((\C_1,\QF_1), (\C_2,\QF_2); - \right) \colon \Cath \to \Sps \ .
\]
Similarly, the object \((\Spaf, \QF^{\uni})\) corepresents the nullary operations
\[
\Mul_{\Cath}\left (\emptyset; -\right) \colon \Cath \to \Sps
\]
since the morphisms \((\pt, \mathbb{S}) \to (\Spaf, \QF^{\uni})\) exhibits the target as the initial hermitian \(\infty\)-category under \((\pt, \SS)\) by Proposition~\refone{proposition:corepresentability-of-poinc}. These constructions thus produce locally cocartesian lifts for the active arrows \(\langle 0 \rangle \to \langle 1 \rangle\) and \(\langle 2 \rangle \to \langle 1 \rangle\), and more generally for every active arrow \(\alp\colon \langle n\rangle \to \langle m\rangle\) whose fibres are of size \(\leq 2\). Since the latter generate all maps in \(\Com\), to prove~\refoneitem{item:thm-one} it will now suffice to show that these locally cocartesian lifts are cocartesian. Taking into account the decomposition of mapping spaces in \(\infty\)-operads it will be enough to verify that the induced maps
\begin{align*}
\Mul_{\Cath}\left ((\C_1,\QF_1) \otimes (\C_2, \QF_2),(\C_3, \QF_3),\ldots, (\C_n, \QF_n); - \right)
\to \Mul_{\Cath}\left ((\C_1,\QF_1),\ldots , (\C_n, \QF_n); - \right)
\end{align*}
and
\begin{align*}
\Mul_{\Cath}\left ((\Spaf,\QF^\uni),(\C_1, \QF_1),\ldots, (\C_n, \QF_n); - \right)
\to \Mul_{\Cath}\left ((\C_1,\QF_1),\ldots , (\C_n, \QF_n); - \right)
\end{align*}
are equivalences of functors \(\Cath \to \Sps\). We will give the argument for the first assertion, the second works similar. The first assertion unwinds to the statement that natural transformations from the functor
\[
\left(\App_2 \beta_!(\QF_1 \boxtimes \QF_2)\right) \boxtimes \QF_3 \boxtimes \ldots \boxtimes \QF_n \colon (\C_1 \otimes \C_2) \times \C_3 \times \ldots \times \C_n \to \Spa
\]
to any functor \( (\C_1 \otimes \C_2) \times \C_3 \times \ldots \times \C_n \to \Spa\) pulled back from a quadratic functor
 \(\C_1 \otimes \ldots \otimes \C_n \to \Spa\)
are equivalent to natural transformations from
\[
\QF_1 \boxtimes \QF_2  \boxtimes \QF_3 \boxtimes \ldots \boxtimes \QF_n \colon \C_1 \times \C_2 \times \C_3 \times \ldots \times \C_n \to \Spa
\]
to the restriction of the same functor. After fixing objects \(\x_3 \in \C_3, \x_4 \in \C_4,...\) it is certainly true that the space of natural transformations between the restricted functors along
\begin{align*}
& (\C_1 \otimes \C_2) \to (\C_1 \otimes \C_2) \times \C_3 \times \ldots \times \C_n &\qquad &(\x_1 \otimes \x_2) \mapsto ((\x_1 \otimes \x_2),\x_3,...,\x_n) \\
& \C_1 \times \C_2 \to \C_1 \times \C_2 \times \C_3 \times \ldots \times \C_n & \qquad & (\x_1, \x_2) \mapsto (\x_1,\x_2,\x_3,...,\x_n)
\end{align*}
agree by the universal properties of left Kan extension \(\beta_!\) and 2-excisive approximation \(\App_2\).
The claim then follows since the space of transformations is a limit over these restricted spaces.

To see \refoneitem{item:thm-two} first observe that the operad map \(p^{\otimes}\colon (\Cath)^\otimes \to (\Catx)^\otimes\) preserves cocartesian edges, as is visible by the explicit formula for the tensor product above. In particular, it is a symmetric monoidal functor. Since the functor on underlying \(\infty\)-categories \(p\colon \Cath \to \Catx\) is a cocartesian fibration by Corollary~\refone{corollary:cocartesian} it now follows from~\cite[Proposition 2.4.2.11]{HTT} that \(p^{\otimes}\) is a locally cocartesian fibration.
To show that \(p^{\otimes}\) is a cocartesian fibration one needs to additionally verify that for every arrow in \(\alp\colon \langle n \rangle \to \langle m \rangle\) in \(\Com^{\otimes}\), the associated transition functor \(\alp_!\colon (\Cath)^{\otimes}_{\langle n \rangle} \to (\Cath)^{\otimes}_{\langle m \rangle}\) sends locally \(p^{\otimes}_{\langle n \rangle}\)-cocartesian edges to locally \(p^{\otimes}_{\langle m \rangle}\)-cocartesian edges (indeed, by the explicit description of locally \(p^{\otimes}\)-cocartesian edges provided in ~\cite[Proposition 2.4.2.11]{HTT}, this would imply that these are closed under composition, and are hence all \(p^{\otimes}\)-cocartesian by~\cite[Proposition 2.4.2.8]{HTT}). Unwinding the definitions and using Corollary~\refone{corollary:cocartesian} we observe that this statement is straightforward when \(\alp\) is inert, and for \(\alp\) active amounts to verifying that for every commutative square of the form
\[\begin{tikzcd}
\C_1 \times ... \times \C_n \ar[r,"f_1 \times ... \times f_n"]\ar[d, "\beta"] & \D_1 \times ... \times \D_n \ar[d, "\beta"] \\
\C_1 \otimes ... \otimes \C_n \ar[r,"f_1 \otimes ... \otimes f_n"] & \D_1 \otimes ... \otimes \D_n \\
\end{tikzcd}
\]
and every collection of quadratic functors \(\QF_i \in \Funq(\C_i)\), the natural map
\[(f_1 \otimes ... \otimes f_n)_!\App_2\beta_![\QF_1 \boxtimes ... \boxtimes \QF_n] \to \App_2\beta_![{f_1}_!\QF_1 \boxtimes ... \boxtimes {f_n}_!\QF_n]\]
is an equivalence. %
Since \(f_1 \otimes ... \otimes f_n\) is exact, restriction along it preserves quadratic functors and hence left Kan extension along it commutes with \(\App_2\). We may consequently identify the above map with the image under \(\App_2\beta_!\) of the map
\[(f_1 \boxtimes ... \boxtimes f_n)_![\QF_1 \boxtimes ... \boxtimes \QF_n] \to {f_1}_!\QF_1 \boxtimes ... \boxtimes {f_n}_!\QF_n.\]
The latter is then easily seen to be an equivalence by the pointwise formula for left Kan extension and the fact that tensor products of spectra commute with colimits in each variable.

For Assertion~\refoneitem{item:thm-three} about \(\Catp\), since every equivalence between Poincaré \(\infty\)-categories belongs to \(\Catp\) it suffice to check that the tensor product of \((\Cath)^{\otimes}\) preserves Poincaré \(\infty\)-categories, which is Corollary~\refone{corollary:poinc}, and that the tensor unit \((\Spaf, \QF^{\uni})\) is Poincaré, which was already observed in Example~\refone{example:universal-category}.
\end{proof}

\begin{corollary}
\label{corollary:poinclax}%
The functors \(\Poinc\colon \Catp \to \Sps\) and \(\catforms\colon \Cath \to \Sps\) admit canonical lax symmetric monoidal structures.
\end{corollary}
\begin{proof}
By Theorem~\refone{theorem:tensorpoincare} and Proposition~\refone{proposition:corepresentability-of-poinc}
both of these functors are corepresented by the respective tensor units, so that the result immediately follows from \cite[Corollary 3.10]{Nik_Yoneda}.
\end{proof}

We now point out that both
\[
(\Catp)^{\otimes} \hrar (\Cath)^{\otimes} \hrar \left((\Cat)_{\mop//\Spa}\right)^\otimes
\]
are subcategory inclusions, and hence induce subcategory inclusions
\begin{equation}
\label{equation:map-on-algebras}%
\Alg_{\cO}(\Catp) \hrar \Alg_{\cO}(\Cath) \hrar \Alg_{\cO}((\Cat)_{\mop//\Spa})
\end{equation}
on the level of algebras for every \(\infty\)-operad \(\cO\). By Lemma~\refone{lemma:symm} an \(\cO\)-algebra in \(\left((\Cat)_{\mop//\Spa}\right)^\otimes\) consists of an \(\cO\)-monoidal \(\infty\)-category \(\C\) equipped with a lax \(\cO\)-monoidal functor \(\QF\colon \C\op \to \Spa\). By construction, such an \(\cO\)-algebra belongs to the essential image of \(\Alg_{\cO}(\Cath)\) if and only if the following conditions hold:
\begin{enumerate}
\item
for every colour \(t \in \cO\) the corresponding \(\infty\)-category \(\C_{t}\) is stable and the functor \(\QF_t\colon \C_t\op \to \Spa\) is quadratic;
\item
for every multi-map \(\alp\colon\{t_1,\dots, t_n\} \to  t'\) in \(\cO\) the induced functor
\[
\alp_*\colon \C_{t_1} \times \ldots \times \C_{t_n} \to \C_{t'}
\]
is exact in each variable.
\end{enumerate}
In addition, such an \(\cO\)-algebra
is further in the essential image of \(\Alg_{\cO}(\Catp)\)
if and only if for every colour \(t \in \cO\), the corresponding hermitian \(\infty\)-category \((\C_{t},\QF_t)\) is Poincaré, and for every multi-map \(\alp\colon\{t_1,\dots, t_n\} \to  t'\) in \(\cO\)
and each tuple \(\x_1 \in \C_{t_1}, \ldots, \x_n \in \C_{t_n}\) of objects the corresponding hermitian functor
\[
\alp_*\colon (\C_{t_1},\QF_{t_1}) \otimes \ldots \otimes (\C_{t_n},\QF_{t_n}) \to (\C_{t'},\QF_{t'})
\]
is Poincaré.

\begin{notation}
\label{notation:sym-monoidal-poincare}%
For an \(\infty\)-operad \(\cO\), we will refer to \(\cO\)-algebra objects \((\C,\QF)\) in \(\Cath\) with respect to the symmetric monoidal structure of Theorem~\refone{theorem:tensorpoincare} as \defi{\(\cO\)-monoidal hermitian \(\infty\)-categories}, and similarly to \(\cO\)-algebra objects in \(\Catp\) as \defi{\(\cO\)- monoidal Poincaré \(\infty\)-categories}. We will then refer to the hermitian (resp.\ Poincaré) structure \(\QF\) as an \(\cO\)-monoidal hermitian (resp.\ Poincaré) structure. When \(\cO=\Com\) we will replace as customary the term \(\cO\)-monoidal by \defi{symmetric monoidal}.
\end{notation}

\subsection{Day convolution of hermitian structures}
\label{subsection:symmetric-monoidal-poincare}%

In this section we will analyse in more explicit terms symmetric monoidal hermitian and Poincaré structures over a fixed stably symmetric monoidal \(\infty\)-category \(\C\), and show that they can be encoded in terms of their linear and bilinear parts. To this end, recall that for two symmetric monoidal \(\infty\)-categories \(\E,\D\), there is an associated \(\infty\)-operad \(\Fun(\E,\D)^{\otimes}\) with underlying \(\infty\)-category \(\Fun(\E,\D)\), called the \emph{Day convolution} \(\infty\)-operad, see~\cite{glasman2016day},~\cite[\S 2.2.6]{HA}, and~\cite{day1970} for the classical counterpart. It is characterized by the following universal property: there is an evaluation map of \(\infty\)-operads
\[
\ev\colon \E^{\otimes} \times_{\Com^{\otimes}} \Fun(\E,\D)^{\otimes} \to \D^{\otimes},
\]
refining the usual evaluation map, such that for every \(\infty\)-operad \(\cO^{\otimes}\), the composed map
\[
\Alg_{\cO}(\Fun(\E,\D)) \to \Alg_{\E \times_{\Com} \cO}(\E \times_{\Com} \Fun(\E,\D)) \xrightarrow{\ev_*} \Alg_{\E \times_{\Com} \cO}(\D)
\]
is an equivalence of \(\infty\)-categories. Here, we may identify \(\Alg_{\E \times_{\Com} \cO}(\D) \simeq \Alg_{\E \times_{\Com} \cO/\cO}(\D \times_{\Com}\cO)\) with the \(\infty\)-category of lax \(\cO\)-monoidal functors from \(\E\) to \(\D\), both considered as \(\cO\)-monoidal \(\infty\)-categories by pulling back along \(\cO \to \Com\).
In particular, for \(\cO^{\otimes}=\Com^{\otimes}\) we get an equivalence between commutative algebra objects in \(\Fun(\E,\D)\) and lax symmetric monoidal functors \(\E \to \D\). On the other hand, taking \(\cO^{\otimes}\) to be the underlying \(\infty\)-operad of a symmetric monoidal \(\infty\)-category \(\C\), we get that lax symmetric monoidal functors \(\C \to \Fun(\E,\D)\) correspond to lax symmetric monoidal functors \(\C \times \E \to \D\).

By~\cite[Corollary 2.2.6.12]{HA} the multi-mapping space in \(\Fun(\E,\D)^{\otimes}\) from a collection \(\{\vphi_i\colon \E \to \D\}_{i=1,...,n}\) to \(\psi\colon \E \to \D\) is given by the space of natural transformations in the square
\[
\begin{tikzcd}
\E \times \ldots \times \E \ar[d,"{\otimes}"'] \ar[rr,"{\vphi_1 \times \ldots \times \vphi_n}"] & &  \D \times \ldots \times \D \ar[d,"{\otimes}"]  \ar[lld,shorten <= 2ex,shorten >= 2ex,Rightarrow] \\
\E \ar[rr,"{\psi}"']  & & \D \ .
\end{tikzcd}
\]
If \(\E\) is small, \(\D\) admits small colimits, and the tensor product in \(\D\) preserves small colimits in each variable, then \(\Fun(\E,\D)^{\otimes}\) is a symmetric monoidal \(\infty\)-category, with tensor product \(\vphi_1 \otimes \vphi_2\colon \E \to \D\) given by the left Kan extension of \(\vphi_1 \boxtimes \vphi_2\colon \E \times \E \to \D \times \D \to \D\) along \(\E \times \E \to \E\) (see~\cite[Proposition 2.2.6.16]{HA}). Furthermore, in this case \(\Fun(\E,\D)\) admits small colimits and the Day convolution product preserves small colimits in each variable~\cite[Lemma 2.13]{glasman2016day}.

\begin{remark}
\label{remark:twice-nested}%
Comparing universal properties we see that if \(\C,\E\) are small symmetric monoidal \(\infty\)-categories and \(\D\) is an \(\infty\)-category with small colimits which is endowed with a symmetric monoidal structure which preserves colimits in each variable, then there is a natural equivalence of symmetric monoidal \(\infty\)-categories
\[
\Fun(\C \times \E,\D) \simeq \Fun(\C,\Fun(\E,\D)) ,
\]
where the left hand side is equipped with the Day convolution structure, and the right hand side with the twice nested Day convolution structure.
\end{remark}

We now relate this concept to the constructions we made in the previous sections.
Given a small symmetric monoidal \(\infty\)-category \(\C\), consider the pullback square of \(\infty\)-operads
\[
\begin{tikzcd}
\Fun(\C\op,\Spa)^\otimes \ar[d] \ar[r] & \Einf\ar[d,"\C"] \\
((\Cat)_{\mop//\Spa})^\otimes \ar[r] & \Cat^\times
\end{tikzcd}
\]
refining the pullback square of \(\infty\)-categories which identifies \(\Fun(\C\op,\Spa)\) as the fibre of the cartesian fibration \((\Cat)_{\mop//\D} \to \Cat\) over \(\C \in \Cat\).
We then have the following:

\begin{lemma}
\label{lemma:monoidal-fibre}%
Let \(\C\) be a small stably symmetric monoidal \(\infty\)-category. Then the \(\infty\)-operad
\(\Fun(\C\op,\Spa)\) is a symmetric monoidal \(\infty\)-category.
Furthermore, we may identify the symmetric monoidal structure on \(\Fun(\C\op,\Spa)\) with \emph{Day convolution}. In particular, for an \(\infty\)-operad \(\cO\), the data of an \(\cO\)-algebra structure on \(\QF \in \Fun(\C\op,\Spa)\) corresponds, naturally in \(\cO\), to that of a lax \(\cO\)-monoidal refinement of \(\QF\), where \(\C\) and \(\Spa\) are considered as \(\cO\)-monoidal \(\infty\)-categories by pulling back along \(\cO \to \Com\).
\end{lemma}

\begin{remark}
It follows from Lemma~\refone{lemma:monoidal-fibre} that the tensor product of \(\QF, \QFtwo \in \Fun(\C\op,\Spa)\) is given by the left Kan extension of the functor \(\QF \boxtimes \QFtwo\colon \C\op \times \C\op \to \Spa\) along \(\C\op \times \C\op \to \C\op\).
\end{remark}

\begin{proof}
Since \(\left((\Cat)_{\mop//\D}\right)^{\otimes} \to \Cat^\times\) is
a cocartesian fibration it follows that \(\Fun(\C\op,\Spa)^\otimes \to \Einf\) is also a cocartesian fibration, so that \(\Fun(\C\op,\Spa)^\otimes\) is a symmetric monoidal \(\infty\)-category.
By the description of algebra objects in \(\left((\Cat)_{\mop//\D}\right)^{\otimes}\) given in Lemma~\refone{lemma:symm}\refoneitem{item:algebras}
and the fact that taking algebra objects is compatible with limits of \(\infty\)-operads we get that
\(\cO\)-algebras in \(\Fun(\C\op,\Spa)^{\otimes}\)
are given by lax \(\cO\)-monoidal functors \(\C\op \to \Spa\), and this description is natural in \(\cO\), which allows us to identify the monoidal structure on \(\Fun(\C\op,\Spa)\) with Day convolution.
\end{proof}

We now wish to understand in similar terms the hermitian context, in which one considers \emph{quadratic} functors \(\C\op \to \Spa\). For this, let us consider the setting of Day convolution when the target \(\D\) is a presentably symmetric monoidal \(\infty\)-category, that is, it is presentable and the tensor product preserves small colimits in each variable. In this case, for a small \(\infty\)-category \(\E\), the functor category \(\Fun(\E,\D)\) is again presentably symmetric monoidal with respect to Day convolution.
Suppose now that \(\I = \{\ovl{p}_{\alp}\colon K_{\alp}^{\triangleright} \to \E\}\) is a small collection of diagrams in \(\E\).  Let \(\Fun_{\I}(\E,\D) \subseteq \Fun(\E,\D)\) be the full subcategory spanned by those functors \(\E \to \D\) which send every diagram in \(\I\) to a limit diagram. We claim that \(\Fun_{\I}(\E,\D)\) is an accessible localisation of \(\Fun(\E,\D)\). Indeed, choose a small set of generators \(T\) for \(\D\), and for \(\alp \in \I\) let us denote by \(p_{\alp} := (\ovl{p}_{\alp})|_{K_{\alp}}\) the corresponding restriction. Then for every \(a \in T\) and \(\alp \in \I\) the diagram \(\ovl{p}_{\alp}\) induces a map
\[
s_{\alp,a}\colon \colim_{K_{\alp}}jp_{\alp}\op \otimes a \to jp_{\alp}(*) \otimes a,
\]
where \(j\colon \E\op \to \Fun(\E,\Sps)\) is the Yoneda embedding and \(\otimes\) denotes the canonical bifunctor \(\Fun(\E,\Sps) \times \D \to \Fun(\E,\D)\) induced levelwise by the tensoring of the presentable \(\infty\)-category \(\D\) over spaces. Let \(S = \{s_{\alp,a}\}_{\alp \in \I, a \in T}\).
We may then identify \(\Fun_{\I}(\E,\D) \subseteq \Fun(\E,\D)\) with the full subcategory spanned by the \(S\)-local objects. Since \(S\) is a set it follows from~\cite[Proposition 5.5.4.15]{HTT} that \(\Fun_{\I}(\E,\D)\) is also presentable and its inclusion admits a left adjoint
\[
L\colon \Fun(\E,\D) \to \Fun_{\I}(\E,\D) .
\]
We then have that \(L\) exhibits \(\Fun_{\I}(\E,\D)\) as the localisation of \(\Fun(\E,\D)\) by the set of maps \(S\). Since \(L\) is a left adjoint functor we also refer to it as a \defi{left Bousfield localisation}.

We now consider the above setup in the context of Day convolution. Recall that in general, a left Bousfield localisation functor \(L\colon \A \to \cB\), with fully-faithful right adjoint \(R\colon \cB \hrar \A\), is said to be \defi{compatible} with respect to a given symmetric monoidal structure \(\A^{\otimes}\) on \(\A\), if for every \(f\colon \x \to \y\) in \(\A\) such that \(L(f)\) is an equivalence, and every \(\z \in \A\), the map \(L(\z \otimes f)\) is again an equivalence. By~\cite[Proposition 2.2.1.9]{HA} the \(\infty\)-category \(\cB\) then inherits a symmetric monoidal structure \(\cB^{\otimes}\) such that \(L\) refines to a symmetric monoidal functor \(L^{\otimes}\colon \A^{\otimes} \to \cB^{\otimes}\) and \(R\) refines to a fully-faithful inclusion of \(\infty\)-operads \(R^{\otimes}\colon \cB^{\otimes} \to \A^{\otimes}\). In addition, the symmetric monoidal functor \(L^{\otimes}\) exhibits \(\cB^{\otimes}\) as universal among symmetric monoidal \(\infty\)-categories receiving a symmetric monoidal functor from \(\A\) which inverts the maps inverted by \(L\). Indeed, the symmetric monoidal functor \(L^{\otimes}\) must factor through such a symmetric monoidal localisation
by~\cite[Proposition 4.1.7.4]{HA}, and the resulting comparison between the two symmetric monoidal \(\infty\)-categories under \(\A^{\otimes}\) is an equivalence since it is an equivalence on the level of underlying \(\infty\)-categories. We will consequently refer to \(L^{\otimes}\) as a \defi{symmetric monoidal Bousfield localisation}.

\begin{lemma}
\label{lemma:day-localisation}%
Given \(\E,\D\) and \(\I\) be as above, let us denote by \(\Fun_{\I}(\E,\D)^{\otimes} \subseteq \Fun(\E,\D)^{\otimes}\) the full suboperad spanned by \(\Fun_{\I}(\E,\D)\). Suppose that \(\I\) is closed under post-composition with \(\x \otimes (-)\colon \E \to \E\) for every \(\x \in \E\), that is, if \(p_\alp \in \I\) then \((\x \otimes (-)) \circ p_\alp\) is also in \(\I\).
Then the left Bousfield localisation functor \(L\) is compatible with Day convolution, and hence extends to a symmetric monoidal localisation Bousfield localisation functor
\[
L^{\otimes}\colon \Fun(\E,\D)^{\otimes} \to \Fun_{\I}(\E,\D)^{\otimes}.
\]
In particular, \(\Fun_{\I}(\E,\D)^{\otimes}\) inherits a symmetric monoidal \(\infty\)-category, universally obtained from \(\Fun(\E,\D)\) by inverting the set of maps \(S\).
\end{lemma}

\begin{remark}
\label{remark:localized-product}%
In the situation of Lemma~\refone{lemma:day-localisation}, the tensor product \(\otimes_{\I}\) in \(\Fun_{\I}(\E,\D)^{\otimes}\) can be expressed in terms of the tensor product \(\otimes_{\Day}\) in \(\Fun(\E,\D)^{\otimes}\) and the localisation functor. Explicitly, the tensor product of \(\vphi,\psi \in \Fun_{\I}(\E,\D)\) is given by \(L(\vphi \otimes_{\Day} \psi)\).
\end{remark}

\begin{remark}
\label{remark:localized-colimits}%
In the situation of Lemma~\refone{lemma:day-localisation}, since Day convolution preserves small colimits in each variable and \(L\) preserves colimits it follows from Remark~\refone{remark:localized-product} that the localized tensor product on \(\Fun_{\I}(\E,\D)\) also preserves small colimits in each variable.
\end{remark}

\begin{proof}[Proof of Lemma~\refone{lemma:day-localisation}]
Let \(W\) be the collection of all maps in \(\Fun(\E,\D)\) whose image under \(L\) is an equivalence. We need to show that
\(W\) is closed
under Day convolution against objects, that is, to show that for \(\tau \in W\) and \(\vphi \in \Fun(\E,\D)\), we have that \(\tau \otimes_{\Day} \vphi\) is again in \(W\).
By~\cite[Proposition 5.5.4.15]{HTT} we have that \(W\) is generated as a strongly saturated class by the set \(S = \{s_{\alp,a}\}\) above, and so it will suffice to show that \(S\) is closed under Day convolution against \(\vphi \in \Fun(\E,\D)\). Since Day convolution preserves colimits in each variable we may as well check this for a generating set of \(\Fun(\E,\D)\). Such a generating set is given, for example, by the functors of the form \(j(\x) \otimes a\), for \(\x \in \E\) and \(a \in T\), where \(j\) is the Yoneda embedding as above. Since \(j\) is symmetric monoidal
(\cite[\S 3]{glasman2016day} or \cite[Corollary 4.8.1.12 and Remark 4.8.1.13]{HA}) and using again that Day convolution preserves colimits in each variable, the closure of \(S\) under Day convolution with these generators now follows from our condition that \(\I\) is closed under post-composition with \(\x \otimes (-)\) for every \(\x \in \E\).
\end{proof}

We now apply the above ideas in the context of quadratic functors.
For a given stable \(\infty\)-category \(\C\), the \(\infty\)-category \(\Funq(\C)\) sits in a diagram
\begin{equation}
\begin{tikzcd}
\Funq(\C) \ar[r] \ar[d] & \Fun(\C\op,\Spa)  \ar[d]\ar[r] & \Delta^0 \ar[d,"\C"] \\
\Cath \ar[r]& \left((\Cat)_{\mop//\D}\right)^{\otimes} \ar[r] & \Cat \ .
\end{tikzcd}
\end{equation}
in which all squares are pullbacks.
We then refine this diagram
to a diagram of \(\infty\)-operads
\begin{equation}
\label{equation:square-dayconvolutions}%
\begin{tikzcd}
\Funq(\C)^\otimes \ar[r] \ar[d] & \ar[r] \ar[d] \Fun(\C\op,\Spa)^\otimes & \Einf\ar[d,"\C"] \\
(\Cath)^\otimes \ar[r] & \left((\Cat)_{\mop//\Spa}\right)^\otimes \ar[r] & \Cat^\otimes \ .
\end{tikzcd}
\end{equation}
by extending the lower and right part and defining the
the \(\infty\)-operads \(\Fun(\C\op,\Spa)^{\otimes}\), \(\Funq(\C)^{\otimes}\)
as the respective pullbacks. Since taking algebra objects is compatible with limits in the target it follows that for an \(\infty\)-operad \(\cO\) we have a pullback square
\[
\begin{tikzcd}
\Alg_{\cO}(\Funq(\C)) \ar[r] \ar[d] & \{\C\} \ar[d] \\
\Alg_{\cO}(\Cath) \ar[r] & \Alg_{\cO}(\Cat) \ ,
\end{tikzcd}
\]
and so \(\cO\)-algebras in \(\Funq(\C)\) correspond to \(\cO\)-monoidal hermitian \(\infty\)-categories refining the underlying \(\cO\)-monoidal \(\infty\)-category of \(\C\). On the other hand, since \(\Funq(\C) \subseteq \Fun(\C\op,\Spa)\) is a full inclusion and the monoidal structure on the latter identifies with Day convolution by Lemma~\refone{lemma:monoidal-fibre},
we may identify the data of an \(\cO\)-algebra structure on a given quadratic functor \(\QF\colon \C\op \to \Spa\) with that of a lax \(\cO\)-monoidal structure. We will refer to such \(\QF\) as \defi{\(\cO\)-monoidal hermitian refinements} of \(\C\).

\begin{corollary}
\label{corollary:monoidal-funq}%
Let \(\C\) be a small stably symmetric monoidal \(\infty\)-category.
Then the full inclusion of \(\infty\)-operads
\[
\Funq(\C)^\otimes \subseteq \Fun(\C\op,\Spa)^{\otimes}
\]
admits a symmetric monoidal left adjoint exhibiting \(\Funq(\C)^\otimes\) as a symmetric monoidal localisation of \(\Fun(\C\op,\Spa)^{\otimes}\). In particular, the \(\infty\)-operad \(\Funq(\C)^\otimes\) is a symmetric monoidal \(\infty\)-category.
\end{corollary}

\begin{remark}
\label{remark:monoidal-funq}%
It follows from Corollary~\refone{corollary:monoidal-funq} that
the tensor product of \(\QF,\QFtwo \in \Funq(\C)\) is given by their Day convolution followed by an application of the left adjoint to the inclusion \(\Funq(\C) \subseteq \Fun(\C\op,\Spa)\). Since \(\QF,\QFtwo\) are in particular reduced the result of this left Kan extension is reduced and hence the left adjoint in question can be implemented via the \(2\)-excisive approximation of Construction~\refone{construction:excisive-approx}. Comparing this description with Construction~\refone{construction:tensor-product} we may equivalently describe the tensor product of \(\QF,\QFtwo \in \Funq(\C)\) as the left Kan extension of the quadratic functor \(\QF \otimes \QFtwo \in \Funq(\C \otimes \C)\) along \(\C\op \otimes \C\op \to \C\op\).
\end{remark}

\begin{proof}[Proof of Corollary~\refone{corollary:monoidal-funq}]
This is a particular case of Lemma~\refone{lemma:day-localisation} since the condition of being a quadratic functor is equivalent to that of being reduced and 2-excisive (Proposition~\refone{proposition:basic-properties-quad-functors}) which in turn can be formulated as sending a suitable set of diagrams in \(\C\op\) (consisting of the constant diagram on \(0\) and all strongly exact \(3\)-cubes) to limit diagrams.
\end{proof}

Our next goal is to understand the monoidal structure on \(\Funq(\C)^{\otimes}\) in terms of decomposition into linear and bilinear components, as described in Corollary~\refone{corollary:classification-monoidal} below. Towards this end, we begin with the following direct application of Lemma~\refone{lemma:day-localisation}:

\begin{corollary}
\label{corollary:monoidal-exact-bilinear}%
Let \(\C\) be a small stably symmetric monoidal \(\infty\)-category \(\C\).
Then the following holds:
\begin{enumerate}
\item
\label{item:exact}%
The full suboperad \(\Funx(\C\op,\Spa)^{\otimes} \subseteq \Fun(\C\op,\Spa)^{\otimes}\) is a symmetric monoidal localisation of the Day convolution product on \(\Fun(\C\op,\Spa)\).
\item
\label{item:bilinear}%
The full suboperad \(\Funb(\C)^{\otimes} \subseteq \Fun(\C\op \times \C\op,\Spa)^{\otimes}\) is a symmetric monoidal localisation of the Day convolution product on \(\Fun(\C\op \times \C\op,\Spa)\).
\end{enumerate}
\end{corollary}

\begin{proof}
In case~\refoneitem{item:exact} we apply Lemma~\refone{lemma:day-localisation} with respect to the collection of diagrams in \(\C\op\) consisting of the constant diagram on \(0\) and all exact squares.
For~\refoneitem{item:bilinear} we apply Lemma~\refone{lemma:day-localisation} with respect to the collection of diagrams in \(\C\op \times \C\op\) consisting of the constant diagrams on \((\x,0)\) and \((0,\x)\) for all \(\x \in \C\op\) and all squares of the form \(\{\x\} \times \sig\) and \(\sig \times \{\x\}\) where \(\sig\colon \Del^1 \times \Del^1 \to \C\op\) is an exact square.
\end{proof}

\begin{remark}
As in Remark~\refone{remark:localized-product},
the tensor product of \(\Lin,\Lin'\in \Funx(\C\op,\Spa)\) is obtained by taking their Day convolution \(\Lin \otimes_{\Day} \Lin'\) and applying to it the left adjoint to the inclusion \(\Funx(\C\op,\Spa) \subseteq \Fun(\C\op,\Spa)\). Similarly, in~\refoneitem{item:bilinear} the tensor product of \(\Bil,\Bil'\in \Funb(\C)\) is obtained by taking their Day convolution \(\Bil \otimes_{\Day} \Bil'\) as functors \(\C\op \times \C\op \to \Spa\), and then applying to it the left adjoint to the inclusion \(\Funb(\C) \subseteq \Fun(\C\op \times \C\op,\Spa)\).
\end{remark}

\begin{remark}
\label{remark:monoidal-funq-colimits}%
It follows from Remark~\refone{remark:localized-colimits} that the tensor products on \(\Funq(\C),\Funx(\C\op,\Spa)\) and \(\Funb(\C)\) of Corollaries~\refone{corollary:monoidal-funq} and~\refone{corollary:monoidal-exact-bilinear} all preserve small colimits in each variable.
\end{remark}

Combining Corollary~\refone{corollary:monoidal-funq} and Corollary~\refone{corollary:monoidal-exact-bilinear} we obtain
\begin{corollary}
\label{corollary:linear-monoidal}%
The linear part functor \(\Lin_{(-)}\) refines to a symmetric monoidal localisation functor
\[
\Lin_{(-)}^{\otimes}\colon \Funq(\C)^{\otimes} \to \Funx(\C\op,\Spa)^{\otimes}.
\]
\end{corollary}
We would like to establish a similar property for the bilinear part functor. To this end,
consider the commutative diagram
\begin{equation}
\label{equation:bilinear-diagonal}%
\begin{tikzcd}
\Funb(\C) \ar[r] \ar[d] & \BiFun(\C) \ar[r] \ar[d] & \Fun_{\ast}(\C\op \times \C\op,\Spa) \ar[r] \ar[d] & \Fun(\C\op \times \C\op,\Spa) \ar[d] \\
\Funq(\C) \ar[r] & \Fun_{\ast}(\C) \ar[r,equal] & \Fun_{\ast}(\C) \ar[r] & \Fun(\C)
\end{tikzcd}
\end{equation}
in which the horizontal arrows are the relevant inclusions and the vertical arrows are all induced by restriction along the diagonal \(\C\op \to \C\op \times \C\op\). An application of Lemma~\refone{lemma:day-localisation} shows that all \(\infty\)-categories in this diagram inherit symmetric monoidal structures from the Day convolution product on the functor \(\infty\)-categories in the right most column, and such that all horizontal inclusions admit symmetric monoidal left adjoints, which are also localisation functors. In addition, the right most vertical functor, given by restriction along the diagonal \(\C\op \to \C\op \times \C\op\), admits a left adjoint via the corresponding left Kan extension. Since the diagonal is symmetric monoidal, so is the corresponding left Kan extension. In addition, since the diagonal admits itself a two sided adjoint via the direct sum functor \(\oplus\colon\C\op \times \C\op \to \C\op\) this left Kan extension is just given by restriction along \(\oplus\). By the universal property of all the appearing localisations it then follows that the vertical arrows in~\eqrefone{equation:bilinear-diagonal} admit symmetric monoidal left adjoints and we consequently obtain a diagram of symmetric monoidal \(\infty\)-categories and symmetric monoidal functors
\begin{equation}
\label{equation:bilinear-diagonal-adj}%
\begin{tikzcd}
\Funb(\C)^{\otimes} & \BiFun(\C)^{\otimes} \ar[l] & \Fun_{\ast}(\C\op \times \C\op,\Spa)^{\otimes} \ar[l] & \Fun(\C\op \times \C\op,\Spa)^{\otimes} \ar[l] \\
\Funq(\C)^{\otimes} \ar[u,"{\Bil_{(-)}}"] & \Fun_{\ast}(\C)^{\otimes} \ar[u] \ar[l] & \Fun_{\ast}(\C)^{\otimes} \ar[l,equal] \ar[u] & \Fun(\C)^{\otimes} \ar[l] \ar[u]
\end{tikzcd}
\end{equation}
in which all horizontal functors are symmetric monoidal localisations and all vertical functors are given by restriction along \(\oplus\colon\C\op \times \C\op \to \C\op\) followed by the projection to the relevant full subcategory of \(\Fun_{\ast}(\C\op\times\C\op,\Spa)\). In particular, the top middle horizontal arrow in~\eqrefone{equation:bilinear-diagonal-adj} is the bi-reduction functor
of Lemma~\refone{lemma:reducification} and the second from the left vertical functor in~\eqrefone{equation:bilinear-diagonal-adj} is the cross-effect functor of Definition~\refone{definition:cross-effect}. Since the cross-effect of any quadratic functor is bilinear this formula also holds for the left most vertical arrow, that is, we may identify it with the bilinear part functor \(\Bil_{(-)}\). Arguing as in the proof of Lemma~\refone{lemma:bilinear-symmetric} we see that this symmetric monoidal refinement of \(\Bil_{(-)}\) is \(\Ct\)-equivariant with respect to the flip action on \(\Funb(\C)^{\otimes}\) and the trivial action on \(\Funq(\C)\), and consequently refines to a symmetric monoidal functor
\begin{equation}
\label{equation:symmetric-monoidal-bil}%
\Bil_{(-)}^{\otimes}\colon \Funq(\C)^{\otimes} \to \Funs(\C)^{\otimes} ,
\end{equation}
where the target is endowed with the symmetric monoidal structure obtained by taking the \(\Ct\)-fixed points of \(\Funb(\C)^{\otimes}\) in the \(\infty\)-category of symmetric monoidal \(\infty\)-categories.

Now recall from that the symmetric bilinear part functor
\[
\Bil_{(-)}\colon \Funq(\C) \to \Funs(\C)
\]
is also a left Bousfield localisation functor, since it admits a fully-faithful right adjoint (in fact, it admits fully-faithful adjoints from both sides). We now verify that the same holds for its symmetric monoidal refinement just constructed:

\begin{lemma}
\label{lemma:smashing}%
The symmetric monoidal functor~\eqrefone{equation:symmetric-monoidal-bil} is a symmetric monoidal localisation functor.
\end{lemma}

\begin{proof}
It will suffice to show that the localisation functor \(\Bil_{(-)}\) is compatible with the symmetric monoidal structure (the resulting comparison between the corresponding localisation and \(\Bil_{(-)}^{\otimes}\) is then necessarily an equivalence since it is an equivalence on the level of underlying \(\infty\)-categories). Now the maps that are sent to equivalences by \(\Bil_{(-)}\) are exactly the maps whose cofibre is exact. Since left Kan extensions along exact functors preserve exact functors (Lemma~\refone{lemma:kan-extension-exact-quadratic}) it will suffice to show that if \(\QF,\QFtwo \in \Funq(\C)\) are such that \(\QFtwo\) is exact then \(\QF \otimes \QFtwo \in \Funq(\C \otimes \C)\) (defined as in Construction~\refone{construction:tensor-product}) is exact. Indeed, its bilinear part is \(\Bil_{\QF} \otimes \Bil_{\QFtwo} = 0\) by Proposition~\refone{proposition:tensor}.

\end{proof}

We now combine the symmetric monoidal structures of Corollaries~\refone{corollary:monoidal-funq} and~\refone{corollary:monoidal-exact-bilinear} to enhance the classification of hermitian structures of Corollary~\refone{corollary:classification-of-quad-functors} to a symmetric monoidal one:

\begin{corollary}[Monoidal classification of hermitian structures]
\label{corollary:classification-monoidal}%
The pullback square of Corollary~\refone{corollary:classification-of-quad-functors} refines to a pullback square in \(\Op\)
\begin{equation}
\label{equation:fundamental-square}%
\begin{tikzcd}
\Funq(\C)^{\otimes} \ar[r,"{\tau^{\otimes}}"] \ar[d,"{\Bil_{(-)}^{\otimes}}"'] & \Ar(\Funx(\C\op,\Spa))^{\otimes} \ar[d] \\
\Funs(\C)^{\otimes} \ar[r] & \Funx(\C\op,\Spa)^{\otimes},
\end{tikzcd}
\end{equation}
in which all corners are symmetric monoidal \(\infty\)-categories and the vertical functors are symmetric monoidal.
Here the left vertical arrow is~\eqrefone{equation:symmetric-monoidal-bil}, the bottom right corner is endowed with the symmetric monoidal structure of corollary~\refone{corollary:monoidal-exact-bilinear}, and the arrow category carries the corresponding pointwise symmetric monoidal structure. \end{corollary}

Concretely, this means that for a given \(\infty\)-operad \(\cO\), the data of an \(\cO\)-monoidal hermitian refinement of a stably symmetric monoidal \(\C\) can be identified with that of a triple \((\Bil,\Lin,\alp\colon \Lin \Rightarrow [\Bil^{\Del}]^{\tC})\) where \(\Bil\colon \C\op \times \C\op \to \Spa\) is a lax \(\cO\)-monoidal symmetric bilinear functor, \(\Lin\colon \C\op \to \Spa\) is a lax \(\cO\)-monoidal exact functor, and \(\Lin \Rightarrow [\Bil^{\Del}]^{\tC}\) is a lax \(\cO\)-monoidal transformation.

\begin{proof}[Proof of Corollary~\refone{corollary:classification-monoidal}]
Recall from Remark~\refone{remark:stable-recollement} that the pair of fully-faithful exact functors
\[
\Funs(\C) \xrightarrow{\QF^{\sym}_{(-)}} \Funq(\C) \leftarrow \Funx(\C\op,\Spa)
\]
form a \emph{stable recollement}. We may then invoke the theory of monoidal recollements as developed in~\cite{ShahQuigley} to promote our classification of hermitian structures to a monoidal one. More precisely, combining Corollary~\refone{corollary:linear-monoidal} and Lemma~\refone{lemma:smashing} it follows that the above stable recollement is symmetric monoidal in the sense of
~\cite[Definition 1.19]{ShahQuigley}.
The symmetric monoidal refinement of the classification square is then a consequence of~\cite[Proposition 1.26]{ShahQuigley}.
\end{proof}

Recall from Definition~\refone{definition:poinc-cats} that for a given stable \(\infty\)-category \(\C\), the \(\infty\)-category \(\Funp(\C)\) is the non-full subcategory of \(\Funq(\C)\) consisting of the Poincaré structures on \(\C\) and the duality preserving natural transformations between them. We can refine this inclusion to a (non-full) inclusion of \(\infty\)-operads by extending~\eqrefone{equation:square-dayconvolutions} to a commutative diagram of \(\infty\)-operads
\begin{equation}
\begin{tikzcd}
\Funp(\C)^{\otimes} \ar[r] \ar[d] & \Funq(\C)^\otimes \ar[r] \ar[d] & \Fun(\C\op,\Spa)^\otimes \ar[r] \ar[d] & \Einf \ar[d,"\C"] \\
(\Catp)^{\otimes} \ar[r] & (\Cath)^\otimes \ar[r]& \left((\Cat)_{\mop//\Spa}\right)^\otimes \ar[r] & \Cat^\times \ .
\end{tikzcd}
\end{equation}
in which all squares are pullbacks. For an \(\infty\)-operad \(\cO\), the data of an \(\cO\)-algebra structure on a given Poincaré structure \(\QF \in \Funp(\C)\) then corresponds to that of an \(\cO\)-monoidal structure on \((\C,\QF) \in \Catp\) which refines the underlying \(\cO\)-monoidal structure of \(\C\).

\begin{construction}
\label{construction:bilinear-to-duality}%
Applying Remark~\refone{remark:twice-nested} and the uniqueness of localized symmetric monoidal structures we deduce that the equivalence
\[
\Funb(\C) \simeq \Funx(\C\op,\Funx(\C\op,\Spa)) = \Funx(\C\op,\Ind(\C)),
\]
refines to a symmetric monoidal equivalence
\[
\Funb(\C)^{\otimes} \simeq  \Funx(\C\op,\Ind(\C))^{\otimes},
\]
where the the domain is endowed with the localized Day convolution product and the target with the corresponding nested localized Day convolution. Passing to non-degenerate bilinear functors we obtain an equivalence
\begin{equation}
\label{equation:bilinear-duality-operads}%
\Funnb(\C)^{\otimes} \simeq \Funx(\C\op,\C)^{\otimes}
\end{equation}
between the full suboperad of \(\Funb(\C)^{\otimes}\) spanned by the non-degenerate bilinear functors and the full suboperad of the Day convolution \(\infty\)-operad \(\C\op \to \C\) spanned by the exact functors. We note that the latter may fail to be a symmetric monoidal \(\infty\)-category since \(\C\) generally does not admit small colimits. To avoid confusion, let us try to make the equivalence~\eqrefone{equation:bilinear-duality-operads} more explicit. Given non-degenerate bilinear functors \(\Bil_1,...,\Bil_n\) and \(\Biltwo\), a multi-map \(\{\Bil_1,...,\Bil_n\} \to \Biltwo\) in \(\Funnb(\C)^{\otimes}\) is given by a natural transformation
\begin{equation}
\label{equation:multi-map-bilinear}%
(m\op \times m\op)_![\Bil_1 \boxtimes \ldots \boxtimes \Bil_n] = (\wtl{m}\op \times \wtl{m}\op)_![\Bil_1 \otimes \ldots \otimes \Bil_n]\to \Biltwo,
\end{equation}
where \(m\colon \C \times \ldots \times \C \to \C\) is the symmetric monoidal product of \(\C\), and \(\wtl{m}\colon \C \otimes \ldots \otimes \C \to \C\) is the exact functor induced by it. Here we point out that the domain in~\eqrefone{equation:multi-map-bilinear} is generally not non-degenerate, which is why \(\Funnb(\C)^{\otimes}\) is usually not a symmetric monoidal \(\infty\)-category. By adjunction we may equivalently describe the data of~\eqrefone{equation:multi-map-bilinear} via a natural transformation in the square
\begin{equation}
\label{equation:multi-beta}%
\begin{tikzcd}
[column sep=9ex]
{[\C\op \otimes \ldots \otimes \C\op] \times [\C\op \otimes \ldots \otimes \C\op]} \ar[d,"{\wtl{m}\op \times \wtl{m}\op}"'] \ar[r,"{\Bil_1 \otimes \ldots \otimes \Bil_n}"] & \Spa \otimes \ldots \otimes \Spa \ar[d,"{\otimes}"] \ar[ld,shorten <= 2ex,shorten >= 2ex,Rightarrow,"{\beta}"'] \\
\C\op \times \C\op \ar[r,"{\Biltwo}"] & \Spa \ .
\end{tikzcd}
\end{equation}
By Corollary~\refone{corollary:poinc} (and its proof), the bilinear functor \(\Bil_1 \otimes \ldots \otimes \Bil_n\) is non-degenerate and represented by \(\Dual_{\Bil_1} \otimes \ldots \Dual_{\Bil_n}\colon \C\op \otimes \ldots \otimes \C\op \to \C \otimes \ldots \otimes \C\), where \(\Dual_{\Bil_i}\colon \C\op \to \C\) is the exact functor representing the non-degenerate bilinear functor \(\Bil_i\). The corresponding multi-map on the right hand side of~\eqrefone{equation:bilinear-duality-operads} is then given by the map
\[
m_![\Dual_{\Bil_1} \boxtimes \ldots \boxtimes \Dual_{\Bil_n}] = \wtl{m}_![\Dual_{\Bil_1} \otimes \ldots \otimes \Dual_{\Bil_n}] \to \Dual_{\Biltwo}
\]
induced by the natural transformation in the square
\begin{equation}
\label{equation:multi-tau}%
\begin{tikzcd}
[column sep=9ex]
\C\op \otimes \ldots \otimes \C\op \ar[d,"{\wtl{m}\op}"'] \ar[r,"{\Dual_{\Bil_1} \otimes \ldots \otimes \Dual_{\Bil_n}}"] & \C \otimes \ldots \otimes \C \ar[d,"\wtl{m}"] \ar[ld,shorten <= 2ex,shorten >= 2ex,"\tau",Rightarrow] \\
\C\op \ar[r,"{\Dual'}"']  & \C \ .
\end{tikzcd}
\end{equation}
associated to~\eqrefone{equation:multi-beta} by Lemma~\eqrefone{lemma:nat-duality}.

We now define \(\Funpb(\C)^{\otimes} \hrar \Funnb(\C)^{\otimes}\) to be the non-full suboperad spanned by the perfect bilinear functors and those multi-maps whose corresponding natural transformations~\eqrefone{equation:multi-tau} is an equivalence.
Similarly, we define \(\Funps(\C)^{\otimes}\) to be the suboperad of the symmetric monoidal \(\infty\)-category \(\Funs(\C)^{\otimes} = [\Funb(\C)^{\otimes}]^{\hC}\) sitting in the pullback square
\[
\begin{tikzcd}
\Funps(\C)^{\otimes} \ar[r] \ar[d] & \Funs(\C)^{\otimes} \ar[d] \\
\Funpb(\C)^{\otimes} \ar[r] & \Funb(\C)^{\otimes} \ .
\end{tikzcd}
\]
We note that by construction, the underlying \(\infty\)-categories of \(\Funpb(\C)^{\otimes}\) and \(\Funps(\C)^{\otimes}\) are the subcategories \(\Funpb(\C)\) and \(\Funps(\C)\) of \(\Funb(\C)\) and \(\Funs(\C)\) respectively, spanned by the perfect (symmetric) bilinear functors respectively and duality preserving (symmetric) transformations between them.
\end{construction}

\begin{lemma}
The pullback square
\[
\begin{tikzcd}
\Funp(\C) \ar[r] \ar[d] & \Funq(\C) \ar[d] \\
\Funpb(\C) \ar[r] & \Funb(\C) \ ,
\end{tikzcd}
\]
refines to a pullback square of \(\infty\)-operads
\[
\begin{tikzcd}
\Funp(\C)^{\otimes} \ar[r] \ar[d] & \Funq(\C)^{\otimes} \ar[d] \\
\Funpb(\C)^{\otimes} \ar[r] & \Funb(\C)^{\otimes} \ ,
\end{tikzcd}
\]
\end{lemma}
\begin{proof}
Since \(\Funpb(\C)^{\otimes} \hrar \Funb(\C)^{\otimes}\) is a suboperad this follows from the fact that a hermitian functor
\[
(\wtl{m},\eta)\colon (\C,\QF_1) \otimes \ldots \otimes (\C,\QF_n) \to (\C,\QFtwo)
\]
lying over the monoidal product \(\wtl{m}\colon \C \otimes \ldots \otimes \C \to \C\) of \(\C\) is Poincaré exactly when the corresponding natural transformation
\[
\wtl{m}(\Dual_{\QF_1} \otimes \ldots \otimes \Dual_{\QF_n}) \rightarrow \Dual_{\QFtwo}\wtl{m}\op
\]
is an equivalence, where we have used the identification of the duality on a Poincaré tensor product of Corollary~\refone{corollary:poinc}.
\end{proof}

\begin{corollary}[Monoidal classification of Poincaré structures]
\label{corollary:classification-monoidal-poincare}%
The pullback square of Corollary~\refone{corollary:classification-of-poinc-functors} refines to a pullback square in \(\Op\)
\begin{equation}
\label{equation:fundamental-square-poincare}%
\begin{tikzcd}
\Funp(\C)^{\otimes} \ar[r,"{\tau^{\otimes}}"] \ar[d,"{\Bil_{(-)}^{\otimes}}"'] & \Ar(\Funx(\C\op,\Spa))^{\otimes} \ar[d] \\
\Funps(\C)^{\otimes} \ar[r] & \Funx(\C\op,\Spa)^{\otimes},
\end{tikzcd}
\end{equation}
in which all corners are symmetric monoidal \(\infty\)-categories and the vertical functors are symmetric monoidal.
\end{corollary}

\begin{corollary}
\label{corollary:poincare-sym-monoidal}%
Let \(\C\) be a stably symmetric monoidal \(\infty\)-category, \(\cO\) an \(\infty\)-operad and \(\QF\colon \C\op \to \Spa\) a lax \(\cO\)-monoidal quadratic functor, so that \(\Bil_{\QF}\) and \(\Lin_{\QF}\) inherit lax \(\cO\)-monoidal structures by Corollary~\refone{corollary:classification-monoidal}.
Then the corresponding \(\cO\)-monoidal hermitian \(\infty\)-category \((\C,\QF)\) is \(\cO\)-monoidal Poincaré if and only if the following holds:
\begin{enumerate}
\item
\label{item:objects}%
The underlying hermitian \(\infty\)-category \((\C,\QF)\) is Poincaré.
\item
\label{item:arrows}%
The lax \(\cO\)-monoidal structure on \(\Dual_{\QF}\colon \C\op \to \C\) induced from that of \(\Bil_{\QF}\) via Construction~\refone{construction:bilinear-to-duality} is (strongly) \(\cO\)-monoidal. In particular, \(\Dual_{\QF}\) is an equivalence of \(\cO\)-monoidal \(\infty\)-categories.
\end{enumerate}
\end{corollary}

Concretely, Corollary~\refone{corollary:poincare-sym-monoidal} implies that for a stably symmetric monoidal \(\infty\)-category \(\C\), providing a compatible \(\cO\)-monoidal Poincaré structure on \(\C\) is equivalent to providing a self-dual equivalence of \(\cO\)-monoidal \(\infty\)-categories \(\Dual\colon \C\op \xrightarrow{\simeq}\C\), a lax \(\cO\)-monoidal exact functor \(\Lin \colon \C\op \to \Spa\), and a lax \(\cO\)-monoidal transformation \(\Lin(-) \Rightarrow [\map(-,\Dual(-))]^{\tC}\).

\begin{proof}[Proof of Corollary~\refone{corollary:poincare-sym-monoidal}]
In light of Corollary~\refone{corollary:classification-monoidal-poincare} it will suffice to show that for a lax \(\cO\)-monoidal bilinear form \(\Bil\), the corresponding \(\cO\)-algebra object in \(\Funb(\C)^{\otimes}\) lies in the non-full subcategory \(\Alg_{\cO}(\Funpb(\C)) \subseteq \Alg_{\cO}(\Funb(\C))\) if and only if \(\Bil\) is perfect and the associated lax \(\cO\)-monoidal structure on \(\Dual_{\Bil}\) is strongly \(\cO\)-monoidal. To see this, note that under the equivalence~\eqrefone{equation:bilinear-duality-operads} of Construction~\refone{construction:bilinear-to-duality}, the suboperad \(\Funpb(\C)^{\otimes}\) corresponds to the suboperad
\(\Fun^{\circ}(\C\op,\C)^{\otimes} \subseteq \Funx(\C\op,\C)^{\otimes}\) spanned by the colours corresponding to equivalences \(\Dual\colon \C\op \xrightarrow{\simeq}\C\) and the multi-maps \(\{\Dual_1,\ldots \Dual_n\} \to \Dual'\) corresponding to natural transformations
\[
\wtl{m}(\Dual_{1} \otimes \ldots \otimes \Dual_{n}) \Rightarrow \Dual'\wtl{m}\op .
\]
which are equivalences. We now observe that by the universal property of the tensor product of stable \(\infty\)-categories, such a natural transformation of exact functors is an equivalence if and only if it is sent to an equivalence of multi-linear maps after pre-composing both sides with the functor
\[
\C \times \ldots \times \C \to \C \otimes \ldots \otimes \C .
\]
We may hence equivalently define this suboperad \(\Funx(\C\op,\C)^{\otimes}\) to be spanned by the colours corresponding to equivalences and the multi-maps \(\{\Dual_1,\ldots \Dual_n\} \to \Dual'\) corresponding to natural equivalences
\[
m(\Dual_{1} \boxtimes \ldots \boxtimes \Dual_{n}) \Rightarrow \Dual'm\op ,
\]
where \(m\) is the composite \(\C \times \ldots \times \C \to \C \otimes \ldots \otimes \C \xrightarrow{\simeq} \C\).
It then follows directly from the definitions that \(\cO\)-algebra objects in this suboperad correspond to equivalences \(\Dual\colon \C\op \to \C\) equipped with strong \(\cO\)-monoidal structures, as desired.
\end{proof}

\subsection{Examples}
\label{subsection:examples-monoidal}%

In this section we want to give some examples of algebras and modules in \(\Catp\).

\begin{example}
\label{example:universal-monoidal}%
Consider the Poincaré \(\infty\)-category \((\Spaf, \QF^{\uni})\). It is the tensor unit of \(\Catp\) and thus admits a structure of commutative algebra.
This algebra structure is given by the usual smash product on \(\Spaf\) and the canonical commutative algebra structure on the universal quadratic functor \(\QF^{\uni}\) induced by its being the image under the symmetric monoidal left adjoint \(\Fun((\Spaf)\op,\Sps) \to \Funq((\Spaf)\op,\Sps)\) of the Yoneda image of \(\SS \in \Spaf\).
\end{example}

\begin{example}
\label{example:hyp-monoidal}%
Since the forgetful functor \(U\colon \Catp \to \Catx\) is symmetric monoidal by Theorem~\refone{theorem:tensorpoincare} every symmetric monoidal Poincaré \(\infty\)-category \((\C,\QF)\) yields a stably symmetric monoidal \(\infty\)-category \(\C\) upon forgetting \(\QF\). On the other direction, we will show in \S\refone{subsection:bilinear-and-pairings} that the formation of hyperbolic categories \(\C \mapsto \Hyp(\C)\) (see \S\refone{subsection:hyp-and-sym-poincare-objects}) gives both a left and a right adjoint to \(U\), and hence \(\Hyp\) is both lax and oplax monoidal, see Remark~\refone{remark:hyp-monoidal} below. It then follows in particular that if \(\C\) is a stably symmetric monoidal \(\infty\)-category then \(\Hyp(\C)\) inherits a symmetric monoidal structure. Unwinding the definitions, this structure is given explicitly by the operation \((\x,\y) \otimes (\x',\y') = (\x \otimes \x',\y \otimes \y')\) with the canonical lax monoidal structure on \(\QF_{\hyp}(\x,\y) = \map_{\C}(\x,\y)\). Similarly, if \((\C,\QF)\) is a symmetric monoidal Poincaré \(\infty\)-category then the Poincaré functor
\[
\fgt\colon (\C,\QF) \to \Hyp(\C)
\]
sending \(\x\) to \((\x,\Dual\x)\) (see \S\refone{subsection:hyp-and-sym-poincare-objects}) is canonically a symmetric monoidal Poincaré functor, since \(\fgt\) acts as the unit of the symmetric monoidal adjunction \(U \dashv \Hyp\), see Remark~\refone{remark:units-and-counits}. On the other hand, the Poincaré functor \(\hyp\colon \Hyp(\C) \to (\C,\QF)\) sending \((\x,\y)\) to \(\x\oplus \Dual\y\) (see \S\refone{subsection:hyp-and-sym-poincare-objects}) is not symmetric monoidal in general, though we will show below that it is a morphism of \((\C,\QF)\)-module objects in \(\Catp\).
\end{example}

\begin{example}
We will show below (see Remark~\refone{remark:pairings-monoidal}) that the association \((\C,\QF) \mapsto \Ar(\C,\QF)\) from \(\Catp\) to \(\Catp\) (see Definition~\refone{definition:arrow-cat}) %
refines to a lax symmetric monoidal functor. It then follows that if \((\C,\QF)\) is a Poincaré \(\infty\)-category then \(\Ar(\C,\QF)\) inherits a symmetric monoidal structure. On the level of underlying objects this is given simply by the levelwise product of arrows: \([\x \to \y] \otimes [\x' \to \y'] = [\x \otimes \x' \to \y \otimes \y']\). The Poincaré functor
\[
\id\colon (\C,\QF) \to \Ar(\C,\QF) \quad\quad \x \mapsto [x=x]
\]
then inherits the structure of a symmetric monoidal Poincaré functor by its role as a unit of a symmetric monoidal adjunction, see Remark~\refone{remark:unit-fully-faithful}.
\end{example}

We now consider the case of modules over ring spectra studied in \S\refone{section:modules}. We will specialize to the case where \(A\) is an \(\Einf\)-ring spectrum, that is, an algebra object in \(\Spa\) over the commutative \(\infty\)-operad, and the base \(\Einf\)-ring \(k\) is just the sphere spectrum (though see Example~\refone{example:k-linear-structures} for a perspective incorporating a general \(k\)).
The \(\infty\)-category \(\Modp{A}\) then carries a symmetric monoidal structure given by tensoring over \(A\). We then consider the data needed in order to promote \((\Modp{A})^{\otimes}\) to a symmetric monoidal hermitian or Poincaré \(\infty\)-category.

To begin, by Corollary~
\refone{corollary:monoidal-funq} and Remark~\refone{remark:monoidal-funq}, the \(\infty\)-category \(\Funq(\Modp{A})\) inherits a symmetric monoidal structure, given by Day convolution followed by \(2\)-excisive approximation, such that symmetric monoidal hermitian refinements of \(\Modp{A}\) correspond to algebra objects \(\QF \in \Alg_{\Einf}(\Funq(\Modp{A}))\). By Theorem~\refone{theorem:classification-genuine-modules} we have a natural equivalence \(\Funq(\Modp{A}) \simeq \Mod_{\N A}\) between the \(\infty\)-category of hermitian structures on \(\Modp{A}\) and that of modules with genuine involution over \(A\), and so by transport of structure the symmetric monoidal structure on \(\Funq(\Modp{A})\) induces one on \(\Mod_{\N A}\).

\begin{definition}
\label{definition:genuine-O-algebra}
For an \(\Einf\)-ring spectrum \(A\) and an \(\infty\)-operad \(\cO\), we refer to an \(\cO\)-algebra object in \(\Mod_{\N A}\) with respect to the above symmetric monoidal structure as an \(\cO\)-algebra with genuine involution over \(A\).
\end{definition}

\begin{warning}
Notation~\refone{definition:genuine-O-algebra} for \(\cO = \Eone\) should not be confused with that of an \(\Eone\)-algebra with genuine anti-involution as in Example~\refone{example:anti-involution}.
\end{warning}

As described in the proof of Corollary~\refone{corollary:classification-monoidal}, the recollement \(\left(\Funs(\Modp{A}),\Funx((\Modp{A})\op,\Spa)\right)\) on \(\Funq(\Modp{A})\) is compatible with the symmetric monoidal structure and hence the same holds for the recollement \((\Mod_{A \otimes_{\SS} A}^{\hC},\Mod_A)\) on \(\Mod_{\N A}\). It then follows that the square
\[
\begin{tikzcd}
\Mod_{\N A} \ar[r] \ar[d] & \Ar(\Mod_A) \ar[d,"\target"] \\
\Mod_{A\otimes_{\SS} A}^{\hC} \ar[r,"{(-)^{\tC}}"] & \Mod_A
\end{tikzcd}
\]
refines to a square of symmetric monoidal \(\infty\)-categories and lax symmetric monoidal functors, such that the equivalence~\eqrefone{equation:cube-mod-involution} becomes a symmetric monoidal one.
To identify the resulting symmetric monoidal structure on the individual components \(\Mod_{A \otimes_{\SS} A}^{\hC},\Mod_A\) we have the following:

\begin{lemma}
Let \(A\) be an \(\Einf\)-ring spectrum.
\begin{enumerate}
\item
\label{item:exact-monoidal-modules}%
Under the equivalence \(\Funx((\Modp{A})\op,\Spa) \simeq \Mod_{A}\) the symmetric monoidal structure on the left hand side induced by Day convolution corresponds to the symmetric monoidal structure \(\otimes_A\) on \(\Mod{A}\).
\item
\label{item:bilinear-monoidal-modules}%
Under the equivalence \(\Funb(\Modp{A}) \simeq \Funx((\Modp{A \otimes_{\SS} A})\op,\Spa) \simeq \Mod_{A \otimes_{\SS} A}\) the symmetric monoidal structure on the left hand side induced by Day convolution corresponds to the symmetric monoidal structure \(\otimes_{A \otimes_{\SS} A}\) on \(\Mod_{A \otimes_{\SS} A}\).
\end{enumerate}
\end{lemma}
\begin{proof}
We begin with~\refoneitem{item:exact-monoidal-modules}.
By Corollary~\cite[Corollary 4.8.1.14]{HA} there is a unique symmetric monoidal structure on \(\Ind(\Modp{A})\) which preserves filtered colimits in each variable and such that the inclusion \(\Modp{A} \hrar \Ind(\Modp{A})\) is symmetric monoidal. These two properties are satisfied by the tensor product on \(\Mod{A}\), and by Remark~\refone{remark:monoidal-funq-colimits} the first property holds for the localized Day convolution product.
To finish the proof it is hence left to verify that for a stable \(\infty\)-category \(\E\) the stable Yoneda embedding \(\E \to \Funx(\E\op,\Spa)\) admits a symmetric monoidal structure. Indeed, it follows from~\cite[Proposition 6.3]{Nik_Yoneda} that the stable Yoneda embedding coincides with the composite
\[
\E \to \Fun(\E\op,\Sps) \xrightarrow{\Sig^{\infty}_+ \circ (-)} \Fun(\E\op,\Spa) \to \Funx(\E\op,\Spa) ,
\]
where the first arrow is the ordinary Yoneda embedding, which is naturally symmetric monoidal~\cite[\S 3]{glasman2016day}, the second is post-composition with the symmetric monoidal suspension infinity functor \(\Sig^{\infty}_+\), and the last one is the left adjoint of the inclusion \(\Funx(\E\op,\Spa) \subseteq \Fun(\E\op,\Sps)\), which is symmetric monoidal by Corollary~\refone{corollary:monoidal-exact-bilinear}.

We now prove~\refoneitem{item:bilinear-monoidal-modules}. In light of~\refoneitem{item:exact-monoidal-modules}, it will suffice to show the localized Day convolution structures on \(\Funb(\Modp{A})\) and \(\Funx(\Modp{A \otimes_{\SS} A})\) coincide. For this we consider the commutative square
\[
\begin{tikzcd}
\Funx(\Modp{A \otimes_{\SS} A}) \ar[d] \ar[r,"{\simeq}"] & \Funb(\Modp{A}) \ar[d] \\
\Fun(\Modp{A \otimes_{\SS} A}) \ar[r] & \Fun((\Modp{A})\op \times (\Modp{A})\op,\Spa)
\end{tikzcd}
\]
where the horizontal maps are induced by restriction along the bifunctor \((X,Y) \mapsto X \otimes_{\SS} Y\) from \(\Modp{A} \times \Modp{A}\) to \(\Modp{A}\), and the vertical maps are the relevant full inclusions. Passing to left adjoints and using the universal property of localisation we obtain a commutative square
\[
\begin{tikzcd}
\Funx(\Modp{A \otimes_{\SS} A}) & \Funb(\Modp{A}) \ar[l,"{\simeq}"'] \\
\Fun(\Modp{A \otimes_{\SS} A}) \ar[u] & \Fun((\Modp{A})\op \times (\Modp{A})\op,\Spa)\ar[l] \ar[u]
\end{tikzcd}
\]
in which the horizontal maps are induced by left Kan extension along \((X,Y) \mapsto X \otimes_{\SS} Y\), followed by projection to \(\Funx(\Modp{A \otimes_{\SS} A})\) in case of the top arrow. In particular, all arrows carry a natural symmetric monoidal structure, and so the top horizontal equivalence identifies the localized Day convolution structures on both sides.
\end{proof}

\begin{corollary}
\label{corollary:monoidal-hermitian}%
Let \(A\) be an \(\Einf\)-ring spectrum and \(\cO\) an \(\infty\)-operad.
Then the data of an \(\cO\)-monoidal hermitian refinement of \((\Modp{A})^{\otimes}\)
corresponds to a triple \((B,C,\alp\colon C \to B^{\tC})\) where \(B \in \Alg_{A \otimes_{\SS} A}^{\hC}\) is an \(\cO\)-\((A \otimes_{\SS} A)\)-algebra equipped with a symmetry with respect to the flip action on \(A \otimes_{\SS} A\), \(C\) is an \(\cO\)-\(A\)-algebra, and \(\alp\) is a map of \(\cO\)-\(A\)-algebras.
\end{corollary}

\begin{corollary}
\label{corollary:monoidal-poincare}%
Let \(A\) be an \(\Einf\)-ring spectrum and \(\cO\) an \(\infty\)-operad. Suppose that \(\cO\) is unital. Then \(B\), an \(\cO\)-algebra with genuine involution over \(A\), determines an \(\cO\)-monoidal Poincaré structure on \((\Modp{A})^{\otimes}\) if and only if the underlying \(\cO\)-\(A\)-algebra of \(B\) (with respect to either of the two component inclusions \(A \to A \otimes_{\SS} A\), this makes no difference due to the symmetry of \(B\)) is initial.
\end{corollary}

\begin{remark}
In the situation of Corollary~\refone{corollary:monoidal-poincare}, the initiality condition on \(B\) can be more explicitly formulated by saying that for every colour \(t \in \cO\) the composed map \(A \to A \otimes_{\SS} A \to B_t\) associated to the essentially unique null-operation \(\{\} \to t\) and either of the two component inclusions \(A \to A \otimes_{\SS} A\), is an equivalence.
\end{remark}

\begin{proof}[Proof of Corollary~\refone{corollary:monoidal-poincare}]
Let \(\QF^{\alp}_{B}\) be the \(\cO\)-monoidal hermitian structure on \(\Modp{A}\) associated to the \(\cO\)-algebra with genuine involution \((B,C,\alp\colon C \to B^{\tC})\).
By Corollary~\refone{corollary:poincare-sym-monoidal} this hermitian structure is \(\cO\)-monoidal Poincaré if and only if \(\QF^{\alp}_{B}\) is Poincaré and the associated lax \(\cO\)-monoidal functor \(\Dual_{B} = \map_A(-,B)\) is \(\cO\)-monoidal. This condition requires in particular that the map
\[
\mathbf{1}_{\Modp{A}} = A \to \map_A(A,B) = B
\]
in question is an equivalence. On the other hand, when this condition holds we get from Proposition~\refone{proposition:recognition} that \((B,C,\alp)\) is induced by an anti-involution on \(A\), and hence in particular \(\QF^{\alp}_B\) is Poincaré by Example~\refone{example:module-with-involution-associated-to-ring-with-involution}. In addition, in this case we may write the duality as a composition of the standard duality \(\Dual_A\) and the \(\Ct\)-action on \(\Modp{A}\) induced by the involution on \(B=A\), both of which are \(\cO\)-monoidal.
\end{proof}

Taking \(\cO=\Einf\) in Corollary~\refone{corollary:monoidal-poincare}, we note that the data of a symmetric \(\Einf-(A \otimes_{\SS} A)\)-algebra whose underlying \(A\)-algebra is initial is the same as the data of a \(\Ct\)-action on \(A\) as an \(\Einf\)-ring spectrum. We hence get that symmetric monoidal Poincaré refinements of \((\Modp{A})^{\otimes}\) correspond to the data of a \(\Ct\)-action on \(A\) as an \(\Einf\)-ring spectrum, together with a map of \(\Einf\)-\(A\)-algebras \(C \to A^{\tC}\). Here \(A^{\tC}\) is considered as a \(\Einf\)-\(A\)-ring spectrum via the Tate Frobenius \(A \to A^{\tC}\). We will refer to such a triple \((A,C,C \to A^{\tC})\) as an \defi{\(\Einf\)-ring spectrum with genuine involution}.

\begin{examples}
\label{examples:monoidal}%
Let \(A \in \Alg_{\Einf}^{\hC}\) be an \(\Einf\)-ring spectrum equipped with an involution. We then have the following examples of genuine refinements of interest:
\begin{enumerate}
\item
\label{item:symmetric}%
The map \(\id\colon A^{\tC} \to A^{\tC}\) determines an \(\Einf\)-ring spectrum with genuine involution \((A,A^{\tC},\id)\), and hence a symmetric monoidal refinement of the associated
symmetric Poincaré structure \(\QF^{\s}_A\) on \(\Modp{A}\).
\item
\label{item:truncated}%
If \(A\) is connective then the commutative \(A\)-algebra map \(\tau_{\geq 0}A^{\tC} \to A^{\tC}\) determines an \(\Einf\)-ring spectrum with genuine involution \((A,A^{\tC},t_0)\), and hence a symmetric monoidal refinement of the Poincaré \(\infty\)-category \((\Modp{A},\QF^{\geq 0}_A)\) of Example~\refone{example:truncation}. When \(A\) is discrete, this is the \emph{genuine symmetric} Poincaré structure, see~\refone{subsection:discrete-rings}.
\item
\label{item:tate}%
The twisted Tate Frobenius map \(\tate\colon A \to (A \otimes_{\SS} A)^{\tC} \to A^{\tC}\), which in this context is also the unit map of the \(\Einf\)-\(A\)-algebra \(A^{\tC}\), determines an \(\Einf\)-ring spectrum with genuine involution \((A,A^{\tC},\tate)\), and hence a symmetric monoidal refinement of the Poincaré structure \(\QF^{\tate}_A\) of Example~\refone{example:tate-structure}.
\end{enumerate}
\end{examples}

We now consider some examples where \(\cO\) is not the commutative \(\infty\)-operad, but instead the operad \(\MCom\) whose algebras are pairs \((A,M)\) of a commutative algebra and a module over it \cite[Definition~3.3.3.8]{HA}. In particular, the data of an \(\cO\)-monoidal hermitian (resp.\ Poincaré) \(\infty\)-category consists of a pair \(((\C,\QF),(\E,\QFE))\) where \((\C,\QF)\) is a symmetric monoidal hermitian (resp.\ Poincaré) \(\infty\)-category and \((\E,\QFE)\) is a module over \((\C,\QF)\) in \(\Cath\) (resp.\ \(\Catp\)). Consider the case where \(\E\) is \(\C\) considered as a module over itself, so that
\(\QFE\) is a module over \(\QF\) in \(\Funq(\C)^{\otimes}\).
By Corollary~\refone{corollary:classification-monoidal} the data of such a \(\QF\)-module corresponds to a choice of \(\Bil_{\QF}\)-module \(\Bil\in\Funs(\C)\), an \(\Lin_{\QF}\)-module \(\Lin \in \Funx(\C\op,\Spa)\), and an \(\Lin_{\QF}\)-module map \(\alp\colon\Lin \to [\Bil^{\Del}]^{\tC}\). For example, we may always take \(\Bil=\Bil_{\QF}\) considered as a module over itself. In this case, if \((\C,\QF)\) is a symmetric monoidal Poincaré \(\infty\)-category then all the structure maps in the \(\MCom\)-algebra structure on \(((\C,\QF),(\C,\QF^{\alp}_{\Bil}))\) are Poincaré, since on the level of the underlying bilinear forms these maps come from the structure of \((\C,\QF)\) as a module over itself. Similarly, any commutative triangle
\[
\begin{tikzcd}
\Lin\ar[dr,"{\alp}"'] \ar[rr] && \Lin' \ar[dl,"{\beta}"] \\
& {[\Bil_{\QF}^{\Del}]^{\tC}} &
\end{tikzcd}
\]
of \(\Lin_{\QF}\)-modules in \(\Funx(\C\op,\Spa)\) induces a morphism \((\C,\QF^{\alp}_{\Bil}) \to (\C,\QF^{\beta}_{\Bil})\) of \((\C,\QF)\)-module objects in \(\Catp\).

\begin{example}
\label{example:quadratic-module}%
Let \((\C,\QF)\) be a symmetric monoidal hermitian \(\infty\)-category with underlying bilinear part \(\Bil = \Bil_{\QF}\). Then the associated quadratic hermitian \(\infty\)-category \((\C,\QF^{\qdr}_{\Bil})\) is canonically a module over \((\C,\QF)\), and the natural hermitian functor
\begin{equation}
\label{equation:map-from-quadratic}%
(\C,\QF^{\qdr}_{\Bil}) \to (\C,\QF)
\end{equation}
is a map of \((\C,\QF)\)-module objects in \(\Cath\). In addition, if \((\C,\QF)\) is Poincaré then \((\C,\QF^{\qdr}_{\Bil})\) is a module over \((\C,\QF)\) in \(\Catp\) and the corresponding Poincaré functor~\eqrefone{equation:map-from-quadratic} is a morphism of \((\C,\QF)\)-module objects in \(\Catp\).
\end{example}

When \(\C = \Modp{A}\) for some \(\Einf\)-algebra \(A\) and \(\QF = \QF^{\alp}_A\) for some \(\Ct\)-action on \(A\) and some map of \(\Einf\)-\(A\)-algebras \(\alp\colon C \to A^{\tC}\) then Corollary~\refone{corollary:monoidal-poincare} tells us that the notion of a \(\Lin_{\QF}\)-module in \(\Funx(\C,\Spa)\) is equivalent to that of a \(C\)-module \(D\) equipped with a \(C\)-module map to \(A^{\tC}\).

\begin{example}
\label{example:module-truncated}%
Let \(A\) be a connective \(\Einf\)-ring spectrum equipped with an involution. Then for every \(m \in \ZZ\) the truncated Poincaré structure \((\Modp{A},\QF^{\geq m}_A) \in \Catp\) of Example~\refone{example:truncation} is canonically a module over the symmetric monoidal Poincaré \(\infty\)-category \((\Modp{A},\QF^{\geq 0}_A)\) of Example~\refone{examples:monoidal}\refoneitem{item:truncated}. Indeed, this follows from the above since \(\tau_{\geq m}A^{\tC}\) is canonically a module over the algebra \(\tau_{\geq 0}A^{\tC}\) for every \(m\). Similarly, the canonical Poincaré functors
\[
(\Modp{A},\QF^{\geq m}_A) \to (\Modp{A},\QF^{\geq m-1}_A)
\]
refines to a map of \((\Modp{A},\QF^{\geq 0}_A)\)-modules.
\end{example}

\begin{example}
\label{example:k-linear-structures}%
Let \(k\) be an \(\Einf\)-ring spectrum and \(A\) an \(\Eone\)-algebra over \(k\). Then \(\Modp{k}\) is a symmetric monoidal \(\infty\)-category which acts on \(\Modp{A}\), giving an \(\MCom^{\otimes}\)-algebra object in \(\Catx\). The fibre product \((\Cath)^{\otimes} \times_{(\Catx)^{\otimes}} \MCom^{\otimes}\) then encodes the Day convolution symmetric monoidal structure on \(\Funq(\Modp{k})\) and an action of this structure on \(\Funq(\Modp{A})\). On the level of modules with genuine involution this is an action of \(\Mod_{\N k}\) on \(\Mod_{\N A}\), which can also be encoded via the action of \(\N k\) on \(\N A\) in \(\Ct\)-spectra. By Example~\refone{examples:monoidal}\refoneitem{item:tate} the Tate Poincaré structure \(\QF^{\tate}_k\) is a commutative algebra object in \(\Funq(\Modp{k}) \simeq \Mod_{\N k}\) corresponding to the Tate Frobenius triple \(k^{\Fr} := (k,k,k \to k^{\tC})\) encoding a commutative algebra object in modules with genuine involution over \(k\), or, equivalently, a commutative \(\N k\)-algebra object in genuine \(\Ct\)-spectra. We may then identify the \(k\)-linear norm \(\N_k A\) of \(A\) (see Remark~\refone{remark:NormModules} and Remark~\refone{remark:k-linear-classification}) with the tensor product \(\N A \otimes_{\N k} k^{\Fr}\), and deduce that \(\N_k A\)-modules in genuine \(\Ct\)-spectra correspond to \(k^{\Fr}\)-modules in \(\N A\)-modules, and hence to \(\QF^{\tate}_k\)-module objects in \(\Funq(\Modp{A})\). In summary, the hermitian structures on \(\Modp{A}\) associated to \(k\)-modules with genuine involution as in~\S\refone{subsection:genuine-modules} (with base \(\Einf\)-ring \(k\)) are exactly hermitian structures on \(\Modp{A}\) which are \(\QF^{\tate}_k\)-modules, and so such \(k\)-modules with genuine involution classify refinements of \(\Modp{A}\) to \((\Modp{k},\QF^{\tate}_k)\)-module objects in \(\Cath\).
\end{example}

%% file: CatofCats.tex
In the present section we study the global structural properties of the \(\infty\)-categories \(\Catp\) and \(\Cath\). We will begin in \S\refone{subsection:limits} by showing that \(\Catp\) and \(\Cath\) have all limits and colimits, and that those are preserved by the inclusion \(\Catp \hrar \Cath\) and by the forgetful functor \(\Cath \to \Catx\).
In \S\refone{subsection:internal} we will prove that the symmetric monoidal structures on \(\Catp\) and \(\Cath\) constructed in \S\refone{subsection:monoidal-structure} are closed, that is, they are equipped with a compatible notion of internal mapping objects, which we refer to as \defi{internal functor categories}. This implies, in particular, that these monoidal structures preserve colimits in each variable separately.
We will then show in \S\refone{subsection:cotensoring} and \S\refone{subsection:tensoring} that \(\Cath\) is tensored and cotensored over \(\Cat\) in a manner compatible with its closed symmetric monoidal structure. Taken together, these properties imply that \(\Cath\) (just like \(\Cat\)) can be viewed as an \((\infty,2)\)-category, that is furthermore enriched over itself. Though this point of view does not directly extend to \(\Catp\), in some special cases, such as poset of faces of finite simplicial complexes, the tensor and cotensor constructions do preserve Poincaré \(\infty\)-categories. We will prove this in the final \S\refone{subsection:finite-complexes}, after dedicating \S\refone{subsection:finite-tensors-cotensors} to the role played by imposing certain finiteness conditions on the (co)tensoring \(\infty\)-category.

In addition to its conceptual aspect, we will also make concrete use of the material of this section in a variety of contexts, including \S\refone{subsection:thom} in the present paper, and later in \papertwo, notably via the \(\Q\)-construction which is used in the definition of the Grothendieck-Witt spectrum, and the dual \(\Q\)-construction which is used in establishing its universal property.

\subsection{Limits and colimits}
\label{subsection:limits}%

In this section we will prove that \(\Cath\) and \(\Catp\) have all small limits and colimits and that the functors \(\Catp \to \Cath\) and \(\Cath \to \Catx\) preserve these limits and colimits. Since the former is conservative it follows automatically that it also detects limits and colimits. We will then verify a few (co)limit-related results which will be useful later on, including the semi-additivity of \(\Cath\) and \(\Catp\), and the fact that the functors \(\Poinc\) and \(\spsforms\) preserve filtered colimits.

We begin by recording the following statement, to which we could not find a reference in this form:

\begin{proposition}
\label{proposition:catx-has-limits-and-colimits}%
The \(\infty\)-category \(\Catx\) admits all small limits and colimits. Limits and filtered colimits are preserved by the inclusion \(\Catx \to \Cat\).
\end{proposition}
\begin{proof}
The statement for limits is~\cite[Theorem 1.1.4.4]{HA} and for filtered colimits is given by~\cite[Proposition 1.1.4.6]{HA}. For general colimits we consider the fully-faithful embedding of \(\Catx\) inside the \(\infty\)-category \(\Catrex\) of \(\infty\)-categories with finite colimits and right exact functors between them. The latter admits small colimits by~\cite[Lemma 6.3.4.4]{HTT}. On the other hand, the embedding \(\Catx \subseteq \Catrex\) admits a left adjoint given by tensoring with \(\Spaf\), see Construction~\refone{construction:monoida-catx}.
It then follows that \(\Catx\) has small colimits obtained by forming colimits in \(\Catrex\) and then tensoring with \(\Spaf\).
\end{proof}

\begin{proposition}
\label{proposition:cath-has-limits-and-colimits}%
The \(\infty\)-category \(\Cath\) admits all small limits and colimits, and these are preserved by the forgetful functor \(\pi\colon \Cath \to \Catx\).
\end{proposition}
\begin{proof}
By construction the forgetful functor \(\pi\colon\Cath \to \Catx\) is a cartesian fibration classified by the functor \(\C \mapsto \Funq(\C)\). By Corollary~\refone{corollary:cocartesian} we have that \(\pi\) is also a cocartesian fibration and by Remark~\refone{remark:closed} the fibres \(\Funq(\C)\) of \(\pi\) admits small limits and colimits. Given an exact functor \(f\colon \C \to \D\) the associated cartesian transition functor \(f^*\colon \Funq(\D) \to \Funq(\C)\) preserves limits and the cocartesian transition functor \(f_!\colon\Funq(\C) \to \Funq(\D)\) preserves colimits (indeed \(f_!\) is left adjoint to \(f^*\)). Now the base \(\Catx\) of this bicartesian fibration admits small limits and colimits by Proposition~\refone{proposition:catx-has-limits-and-colimits}. It then follows from \cite[4.3.1.11]{HTT} and \cite[4.3.1.5.(2)]{HTT} that \(\Cath\) admits all small limits and colimits and that these are preserved by \(\pi\).
\end{proof}

\begin{remark}
\label{remark:explicit-limits-cath}%
To make the content (and proof) of Proposition~\refone{proposition:cath-has-limits-and-colimits} more explicit, let \(K\) be a simplicial set and \(p\colon K \to \Cath\) a diagram. Then the limit of \(p\) is computed as follows: one first extends the diagram \(q := \pi p\colon K \to \Catx\) of stable \(\infty\)-categories to a limit diagram \(\ovl{q}\colon K^{\triangleleft} \to \Catx\). Let \(\C_{\infty} = \ovl{q}(\infty)\) be the image of the cone point (which we denote by the symbol \(\infty\)), so that \(\ovl{q}\) exhibits \(\C_{\infty}\) as the limit of \(q\) in \(\Catx\). Interpreting \(\ovl{q}\) as a natural transformation with target \(q\colon K \to \Catx\) and domain the constant diagram \(K \to \Catx\) with value \(\C_{\infty}\), we may lift it to a pointwise \(\pi\)-cartesian natural transformation with target \(p\) and domain some diagram \(p_{\infty}\colon K \to \Cath\) which is concentrated in the fibre of \(\C_{\infty}\). In other words, \(p_{\infty}\) encodes a diagram \(p_{\infty}\colon K \to \Funq(\C_{\infty})\). The hermitian \(\infty\)-category \((\C_{\infty},\lim_{K} p_{\infty})\) is then the limit of the original diagram \(p\). Somewhat informally, though more explicitly, we may describe the diagram \(p_{\infty}\) via the formula \(k \mapsto r_k^*\QF_k\), where \((\C_k,\QF_k)\) is the hermitian \(\infty\)-category associated to \(k\) by the diagram \(p\) and \(r_k\colon \C_{\infty} \to \C_k\) is the exact functor associated to the map \(\infty \to k\) in \(K^{\triangleleft}\) by the limit diagram \(\ovl{p}\). Since \(\Funq(\C_{\infty})\) is closed in \(\Fun(\C_{\infty}\op,\Spa)\) under limits the limit of \(p_{\infty}\) can also be computed in \(\Fun(\C_{\infty}\op,\Spa)\), that is, object-wise. Similarly, in order to compute the colimit of \(p\) we first extend \(q=\pi p\) to a colimit diagram \(\ovl{q}\colon K^{\triangleright} \to \Catx\). Let \(\C_{\infty} = \ovl{q}(\infty)\) be the image of the cone point, and \(i_k\colon \C_k \to \C_{\infty}\) the exact functor associated to the arrow \(k \to \infty\) in \(K^{\triangleright}\).
Then as above we can ``push'' the diagram \(p\) to the fibre over \(\C_{\infty}\), yielding a diagram \(p_{\infty}\colon K \to \Funq(\C_{\infty})\) given by the formula \(k \mapsto (i_k)_!\QF_k\). The colimit of \(p\) in \(\Cath\) is then given by
\[
(\C_{\infty},\colim_K p_{\infty}) = (\C_{\infty},\colim_{k \in K} (i_k)_! \QF_k),
\]
where the colimit of \(p_{\infty}\) is computed in \(\Funq(\C_{\infty})\). This last colimit can also be computed in \(\Fun(\C_{\infty}\op,\Spa)\) since \(\Funq(\C_{\infty})\) is closed in \(\Fun(\C_{\infty}\op,\Spa)\) under colimits (see Remark~\refone{remark:closed}).
\end{remark}

\begin{proposition}
\label{proposition:Catp-cocomplete}%
The \(\infty\)-category \(\Catp\) has small limits and colimits and the inclusion \(\Catp\to \Cath\) preserves small limits and colimits.
\end{proposition}

The proof of Proposition~\refone{proposition:Catp-cocomplete} will require the following lemma.

\begin{lemma}
\label{lemma:colim-stable-cats}%
Let \(\ovl{p}:K^\triangleright\to \Catx\) be a colimit diagram of stable \(\infty\)-categories and let \(f\colon \C_\infty := \ovl{p}(\infty)\to \D\) be an exact functor to a cocomplete stable \(\infty\)-category. Then the canonical map
\[
\colim_{k\in K} (i_k)_!i_k^*f\to f
\]
is an equivalence, where \(i_k\colon \C_k := \ovl{p}(k)\to \C_\infty\) is the map associated to \(k \to \infty\) by \(\ovl{p}\).
\end{lemma}
\begin{proof}
For any diagram \(\ovl{p}:K^\triangleright\to \Catx\), the functor
\[
i^*:\Funx(\C_\infty,\D)\to \lim_{k\in K}\Fun^{\ex}(\C_k,\D)
\]
sending \(f\) to \(\{f\circ i_k\}_{k\in K}\) has a left adjoint given by
\[
i_!:\{f_k\}_{k\in K}\mapsto \colim_{k\in K}(i_k)_!f_k.
\]
Indeed, this follows for example by the general formula of~\cite[Theorem B]{horev-yanovski}.
Our claim is then equivalent to the statement that the counit
\[
i_!i^*f\to f
\]
is an equivalence if \(\ovl{p}\) is a colimit diagram. But \(\Funx(-,\D):(\Catx)\op\to \Catx\) is a right adjoint (in fact, its own right adjoint) and so it preserves all limits. Hence \(i^*\) is an equivalence of \(\infty\)-categories and in particular the counit is an equivalence of exact functors.
\end{proof}

\begin{proof}[Proof of \refone{proposition:Catp-cocomplete}]

We begin with the case of colimits.
Let \(p: K \to \Catp\) be a diagram. By Proposition~\refone{proposition:cath-has-limits-and-colimits} we may find a colimit diagram \(\ovl{p}\colon K^\triangleright\to \Cath\) in \(\Cath\) extending (the image of) \(p\). We will show that the image of \(\ovl{p}\) is contained in \(\Catp\) and forms a colimit diagram there. Let \((\C_{\infty},\QF_{\infty}) = \ovl{p}(\infty)\) be the image of the cone point and for \(k \in K\) let \((i_k,\eta_k)\colon (\C_k,\QF_k) \to (\C_{\infty},\QF_{\infty})\) be the hermitian functor associated to the map \(k \to \infty\) in \(K^{\triangleright}\) by \(\ovl{p}\). In particular, the collection of natural transformations \(\eta^{\ad}_k\colon (i_k)_!\QF_k \to \QF_{\infty}\) exhibits \(\QF_{\infty}\) as the colimit of the diagram \(k \mapsto (i_k)_!\QF_k\).
Our argument proceeds in two steps:
\smallskip

\textbf{Step 1.} \emph{The hermitian \(\infty\)-category \((\C_{\infty},\QF_{\infty})\) is Poincaré, and all the hermitian functors \((i_k,\eta_k)\colon (\C_k,\QF_k) \to (\C_{\infty},\QF_{\infty})\) are Poincaré.}\\
Let \(\Dual_k\colon \C_k\op\to \C_k\) be the duality associated to \(\QF_k\). Let \(q := \pi p\colon K \to \Catx\) be the underlying diagram of stable \(\infty\)-categories and let us denote by
\(q_{\mathrm{op}}\colon K \to \Catx\) the composite of \(q\) and the equivalence \((-)\op \colon \Catx \xrightarrow{\simeq} \Catx\). Since the diagram \(p\) takes values in Poincaré \(\infty\)-categories and duality preserving functors the functors \(\Dual_k\) form the components of a natural transformation \(q_{\mathrm{op}} \Rightarrow q\), which consequently induce a functor
\[
\Dual_{\infty}\colon \C_{\infty}\op \to \C_{\infty} ,
\]
where we note that \((-)\op\colon \Catx \to \Catx\) preserves colimits since it is an equivalence.
We claim that \(\Dual_{\infty}\) represents the bilinear part of \(\QF\).
Indeed, since bilinear parts commute with left Kan extensions (Lemma~\refone{proposition:left-kan-bilinear-linear}) and colimits (by Lemma~\refone{lemma:universal-crs}) the bilinear part of \(\QF_\infty\) is given by
\begin{multline*}
\displaystyle\mathop{\colim}_{k\in K} (i_k\op\times i_k\op)_!\Bil_{\QF_k}(-,-) \simeq \displaystyle\mathop{\colim}_{k\in K} (i_k\op\times i_k\op)_!\map_{\C_k}(-,\Dual_k(-))\simeq \displaystyle\mathop{\colim}_{k\in K} (\id\times i_k\op)_!\map_{\C_{\infty}}(-,i_k\Dual_k(-))\simeq \\
\simeq \displaystyle\mathop{\colim}_{k\in K} (\id\times i_k\op)_!\map_{\C_\infty}(-,\Dual_{\infty} i_k\op(-))\simeq
\displaystyle\mathop{\colim}_{k\in K} (\id\times i_k\op)_!(\id\times i_k\op)^*\map_{\C_{\infty}}(-,\Dual_{\infty}-) \simeq \map_{\C_\infty}(-,\Dual_{\infty}(-)),\\
\end{multline*}
where the last equivalence is by
Lemma~\refone{lemma:colim-stable-cats}. We have thus established that the hermitian \(\infty\)-category \((\C_{\infty},\QF_{\infty})\) is non-degenerate and all the hermitian functors \((i_k,\eta_k)\colon (\C_k,\QF_k) \to (\C_{\infty},\QF_{\infty})\) are duality preserving. To finish the proof it will hence suffice to show that \(\QF_{\infty}\) is perfect. But this follows from Remark~\refone{remark:duality-equivalence} since \(\Dual_{\infty}\) is an equivalence, being the functor induced on limits by a natural equivalence of diagrams \(q_{\mathrm{op}} \Rightarrow q\).

\smallskip

\textbf{Step 2.} \emph{If \((\D,\QFD)\) is a Poincaré \(\infty\)-category then a hermitian functor \((f,\eta)\colon (\C_{\infty},\QF_{\infty}) \to (\D,\QFD)\) is Poincaré if and only if \((f\circ i_k,i_k^*\eta)\) is Poincaré for every \(k\in K\).}\\
The hermitian functor \((f,\eta)\) is a Poincaré functor if and only if the associated natural transformation
\begin{equation}
\label{equation:colim-nat}%
f\Dual_\infty \to \Dual_\infty f\op
\end{equation}
is an equivalence. Since \(\C_{\infty}\) is the colimit of \(p\) in \(\Catx\) the natural transformation~\eqrefone{equation:colim-nat} is an equivalence if and only if the natural transformation
\[
f\Dual_\infty i_k\op \to \Dual_\infty f\op i_k\op
\]
is an equivalence for all \(k\in K\). But, since the \(i_k\) are all Poincaré functors, there is a commutative diagram
\[
\begin{tikzcd}
f\circ \Dual_\infty \circ i_k\op \ar[r] & \Dual_\infty \circ f\op \circ i_k\op \\
f\circ i_k\circ \Dual_k \ar[u,"\simeq"] \ar[ur]
\end{tikzcd}
\]
Thus, \eqrefone{equation:colim-nat} is an equivalence if and only if \((f\circ i_k,i_k^*\eta)\) is a Poincaré functor for all \(k\in K\).

The case of limits is similar and slightly easier. Indeed, as above, the natural equivalence of diagrams \(q_{\mathrm{op}} \Rightarrow q\) induced by the collection of dualities \(\Dual_k\) induces an equivalence
\[
\Dual_\infty\colon \C_{\infty}\op = \lim_K q_{\mathrm{op}} \to \lim_K q = \C_{\infty}
\]
This time, showing that \(\Dual_\infty\) represents the bilinear form of \(\QF_{\infty} = \lim_{k \in K} r_k^*\QF_k\) (where \(r_k\colon \C_{\infty} \to \C_k\) are the canonical maps) is even simpler. Indeed, since taking bilinear forms commutes with restriction (Remark~\refone{remark:invariance-base-change}), and limits (by Lemma~\refone{lemma:universal-crs}), the bilinear part of \(\QF_\infty\) is given by
\[
\lim_{k \in K} (r_k \times r_k)^*\Bil_{\QF_k} \simeq \lim_{k \in K} \map_{\C_k}(r_k(-),\Dual_k(r_k(-)) = \lim_{k \in K} \map_{\C_k}(r_k(-),r_k\circ \Dual_\infty) \simeq \map_{\C}(-,\Dual_\infty(-)),
\]
and this concludes the proof of Step 1 in the case of limits. The proof of Step 2 is completely dual to that of colimits.
\end{proof}

\begin{remark}
\label{remark:computed}%
By Proposition~\refone{proposition:Catp-cocomplete} the \(\infty\)-category \(\Catp\) has small limits and colimits, and those are preserved by the inclusion in \(\Cath\), which itself also has small limits and colimits. Since the forgetful functor \(\Catp \to \Cath\) is conservative, it also detects limits and colimits. One then says that limits and colimits in \(\Catp\) are \defi{computed} in \(\Cath\).
\end{remark}

\begin{proposition}
\label{proposition:catp-pre-add}%
The \(\infty\)-categories \(\Catx\), \(\Cath\) and \(\Catp\) are all semi-additive, i.e.\ finite products and coproducts agree.
\end{proposition}
\begin{proof}
For \(\Catx\), as in the proof of Proposition~\refone{proposition:catx-has-limits-and-colimits} we can embed  \(\Catx\) as a reflective full subcategory of \(\Catrex\). It will then suffice to verify that \(\Catrex\) is semi-additive, which in turn follows from~\cite[Lemma 7.3.3.4]{HTT}. In particular, the coproduct of \(\C\) and \(\Ctwo\) in \(\Catx\) is given by the product \(\C\times\Ctwo\), and exhibited by the two inclusions
\[
i\colon \C \times \{0\} \rightarrow \C \times \Ctwo \leftarrow \{0\} \times \Ctwo \cocolon i'.
\]
Turning to the other two cases, we first observe that since limits and colimits in \(\Catp\) are computed in \(\Cath\) (Remark~\refone{remark:computed}) it is in fact sufficient to show that \(\Cath\) is semi-additive. We first observe that \(\Cath\) is pointed. Indeed, the hermitian \(\infty\)-category \((\{0\},0)\), whose underlying stable \(\infty\)-category consists of a single object \(0\), equipped with the hermitian structure which sends \(0\) to the zero spectrum, is both initial and final in \(\Cath\).

Now let \((\C,\QF)\) and \((\Ctwo,\QFtwo)\) be two hermitian \(\infty\)-categories.
By the explicit description of Remark~\refone{remark:explicit-limits-cath} we have that the coproduct of \((\C,\QF)\) and \((\Ctwo,\QFtwo)\) in \(\Cath\) is given by \((\C \times \Ctwo,i_!\QF\oplus i'_!\QFtwo)\), while the product is given by \((\C\times \Ctwo, p^*\QF\oplus (p')^*\QFtwo)\), where \(p,p'\) are the projections on \(\C\) and \(\Ctwo\), respectively. It will hence suffice to show that the canonical natural transformation
\begin{equation}
\label{equation:semi-additive}%
i_!\QF\oplus i'_!\QFtwo \Rightarrow p^*\QF\oplus (p')^*\QFtwo
\end{equation}
is an equivalence.
But, since \(i\) is also right adjoint to \(p\) and \(i'\) is right adjoint to \(p'\) we have natural equivalences
\[
i_!\QF\simeq p^*\QF\textrm{ and }i'_!\QFtwo\cong (p')^*\QFtwo,
\]
under which the canonical map~\eqrefone{equation:semi-additive} becomes the identity, since \(p \circ i\) and \(p' \circ i'\) are the respective identities while \(p \circ i'\) and \(p' \circ i\) are the zero functors.
\end{proof}

The following result will be useful for us in \paperfour for proving that \(\Catp\) and \(\Cath\) are \emph{compactly generated} presentable \(\infty\)-categories:
\begin{proposition}
\label{proposition:forms-poinc-compact}%
The functors \(\spsforms \colon \Cath \rightarrow \Sps\) and \(\Poinc \colon \Catp \rightarrow \Sps\) commute with filtered colimits. In particular, the object \((\Spaf,\QF^{\uni})\) which corepresents these functors (Proposition~\refone{proposition:corepresentability-of-poinc}) is compact in both \(\Cath\) and \(\Catp\).
\end{proposition}
\begin{proof}
    Let us first prove the statement for \(\spsforms\). Let \(\E\to \Cat\) be the presentable fibration classified by the functor \(\Cat\to \CAT\) sending \(\C\) to \(\Psh(\C)\), so the objects of \(\mathcal{E}\) are given by a pair \((\C,S)\) where \(S\colon\C\op\to\Sps\) is a presheaf, and morphisms are pairs \((f,\eta)\) where \(f\colon\C\to\C'\) is a functor and \(\eta\colon S\to f^*S'\) is a natural transformation. There is a functor \(\Phi\colon \Cath\to \E\) lying above the inclusion \(\Catx\to\Cat\) sending \((\C,\QF)\) to \((\C,\Omega^\infty \QF)\). Then we can factor \(\spsforms(-)\) as
\[
\Cath\xrightarrow{\Phi} \E \to \Sps
\]
where the second functor is the functor corepresented by the final object \(\ast:=(\Delta^0,\ast\colon\Delta^0\to\Sps)\) in \(\E\). Note that \(\Phi\) is a morphism of cocartesian fibrations by Proposition~\refone{proposition:left-kan-bilinear-linear} and that it preserves filtered colimits fibrewise. Moreover it lies above the inclusion \(\Catx\to\Cat\), which preserves filtered colimits by~\cite[Proposition~1.1.4.6]{HA}. Hence the functor \(\Phi\) preserves filtered colimits. It will hence suffice to show that \(\ast\) is compact in \(\E\). Now by the naturality of the straightening-unstraightening equivalence as recorded in~\cite[Corollary A.31]{GHNfree}, the unstraightening procedure determines an equivalence \(\E \simeq \RFib\), where \(\RFib \subseteq \Ar(\Cat)\) is the full subcategory spanned by the right fibrations. It will hence suffice to show that the final object in \(\RFib\) is compact. This object corresponds to the identity right fibration \(\Del^0 \to \Del^0\), which is compact when considered as an object of \(\Ar(\Cat)\) since \(\Del^0\) is compact in \(\Cat\). It will hence suffice to show that the inclusion \(\RFib \subseteq \Ar(\Cat)\) is closed under filtered colimits. Indeed, the condition of being a right fibration can be phrased as being local with respect to the maps in \(\Ar(\Cat)\) encoded by the squares of the form
\[
\begin{tikzcd}
\Lam^n_i \ar[r]\ar[d] & \Del^n \ar[d] \\
\Del^n \ar[r] & \Del^n
\end{tikzcd}
\]
with \(0 < i \leq n\). This is an arrow between two compact objects of \(\Ar(\Cat)\) and hence the locality condition it defines is closed under filtered colimits.

Now we want to prove the statement for \(\Poinc\). We know that the functors \(\Catx\to\Sps\) sending \(\C\) to \(\iota\C\) and \(\C\mapsto \iota\Twar(\C) \simeq \iota\Fun(\Delta^1,\C)\) commute with filtered colimits since \(\Del^0\) and \(\Del^1\) are compact in \(\Cat\) and the inclusion \(\Catx \subseteq \Cat\) preserves filtered colimits by \cite[Proposition~1.1.4.6]{HA}.
We then consider the cartesian square (see~\eqrefone{equation:poinc-forms-twar})
\[
\begin{tikzcd}
\Poinc(\C,\QF) \ar[r] \ar[d] & \spsforms(\C,\QF) \ar[d] \\
\Twar(\iota\C) \ar[r] & \iota\Twar(\C)
\end{tikzcd}
\]
where we identify \(\Twar(\iota\C)\) with the subspace of \(\iota\Twar(\C)\) spanned by those arrows that are equivalences in \(\C\). Since all the corners of the square except \(\Poinc\) preserve filtered colimits in \(\Cath\) (and so in \(\Catp\)), so does \(\Poinc\).
\end{proof}

We finish this subsection by describing an example of a colimit in \(\Catp\) arising from visible Poincaré structures on parametrised spectra. Recall from \S\refone{subsection:parametrised-spectra} that given a space \(X\) and a spherical fibration \(\xi\colon X \to \Pic(\SS)\), one may endow the \(\infty\)-category \(\Spaf_X\) of finitely presented parametrised spectra over \(X\) with the \(\xi\)-twisted \emph{visible Poincaré structure} \(\QF^{\vis}_{\xi}\), whose \(\L\)-groups are the associated visible \(\L\)-groups studied in~\cite{WWIII} (see Corollary~\refone{corollary:compare-visible-WW}). We now show that this construction actually produces the colimit of a suitable diagram of Poincaré \(\infty\)-categories. Note first that, as inspecting the definitions reveals, for every \(x \in X\) the Poincaré structure \(\QF^{\vis}_{\xi_x}\) on \(\Spaf_{\{x\}} = \Spaf\) identifies with the tensor \(\xi_x \otimes \QF^{\uni}\), defined using the action of \(\Spaf\) on \(\Funq(\Spaf)\) canonically determined by the stability of the latter. We may then consider the composite diagram
\begin{equation}
\label{equation:spherical-diagram}%
\begin{tikzcd}
[row sep=1ex]
X \ar[r,"\xi"] & \Pic(\SS) \ar[r,hookrightarrow] & \Spaf \ar[r,"(-) \otimes \QF^{\uni}"] & \Funq(\Spaf) \ar[r] & \Cath \\
x \ar[rr,mapsto] & & \xi_x \ar[r,mapsto] & \QF^{\vis}_{\xi_x} \ar[r,mapsto] & (\Spaf,\QF^{\vis}_{\xi_x}) \\
\end{tikzcd}
\end{equation}
where the last arrow places the fibre of the cartesian fibration \(\Cath \to \Catx\) over \(\Spaf \in \Catx\). This diagram can then be viewed as a diagram in \(\Catp\): indeed, each \((\Spaf,\QF^{\vis}_{\xi_x})\) is Poincaré by Corollary~\refone{corollary:coswandual} and all the hermitian functors involved are equivalences in \(\Cath\) and hence in particular duality preserving.

\begin{proposition}
\label{proposition:visible-is-colimit-full}%
Let \(X\) be a space equipped with a spherical fibration \(\xi \colon X \to \Pic(\SS)\). Then the collection of Poincaré functors
\[
(\lambda_x,\eta_x)\colon (\Spaf,\QF^{\vis}_{\xi_x}) \to (\Spaf_X,\QF^{\vis}_{\xi}),
\]
associated by Remark~\refone{remark:functoriality} to the various point inclusions \(\{x\} \hrar X\), exhibits \((\Spaf_X,\QF^{\vis}_{\xi})\) as the colimit in \(\Catp\) of the diagram
\[
X \to \Catp \quad\quad x \mapsto (\Spaf,\QF^{\vis}_{\xi_x})
\]
constructed in~\eqrefone{equation:spherical-diagram}.
\end{proposition}

\begin{remark}
In the situation of Proposition~\refone{proposition:visible-is-colimit-full},
if \(\xi = \SS_X\) is the constant spherical fibration with value \(\SS\) then the composed diagram~\eqrefone{equation:spherical-diagram} is constant with value \((\Spaf,\QF^{\uni})\), and we may view the statement of the proposition as identifying \((\Spaf_X,\QF^{\vis}_{\SS_X})\) with the \emph{tensor} of the Poincaré \(\infty\)-category \((\Spaf,\QF^{\uni})\) and the space \(X\). We will discuss such tensors from a more general perspective in \S\refone{subsection:tensoring} below.
\end{remark}

\begin{proof}[Proof of Proposition~\refone{proposition:visible-is-colimit-full}]
By Proposition~\refone{proposition:Catp-cocomplete} it will suffice to check that the underlying hermitian functors exhibit \((\Spaf_X,\QF^{\vis}_{\xi})\) as the colimit in \(\Cath\) of the diagram in question.
By the explicit description of colimits in \(\Cath\) given in Remark~\refone{remark:explicit-limits-cath} and the fact that we have already established this colimit property on the level of quadratic functors in Proposition~\refone{proposition:visible-is-colimit} (see Remark~\refone{remark:finite-visible-is-colimit}),
all that is left is to check that the collection of exact functors \(\lambda_x\colon \Spaf \to \Spaf_X\) exhibit \(\Spaf_X\) as the colimit in \(\Catx\) of the constant \(X\)-indexed diagram with value \(\Spaf\). We now note that the equivalence \((-)\op \colon \Cat \xrightarrow{\simeq} \Cat\) preserves stable \(\infty\)-categories and exact functors, and hence restricts to an equivalence \((-)\op \colon \Catx \xrightarrow{\simeq} \Catx\). We may hence show instead that the diagram in question becomes a colimit diagram after applying \((-)\op\). We now observe that
\[
\Spa \simeq \Ind(\Spaf) \simeq \Pro((\Spaf)\op)\op,
\]
and consequently
\[
\Fun(X,\Spa)\op \simeq \Fun(X\op,\Pro((\Spaf)\op)).
\]
Under the last equivalence, the full subcategory \((\Spaf_X)\op\) of finitely presented parametrised spectra corresponds to the full subcategory \(((\Spaf)\op)_X \subseteq \Fun(X\op,\Pro((\Spaf)\op))\) generated under finite limits by the functors \(X\op \to \Pro((\Spaf)\op)\) right Kan extended along point inclusions \(\{x\} \subseteq X\op\) from functors \(\{x\} \to \Pro((\Spaf)\op)\) with value in
\((\Spaf)\op\).
The desired result now follows from the fact that the functor
\[
(\Spaf)\op \times X \to ((\Spaf)\op)_X \quad\quad (E,x) \mapsto x_*E
\]
exhibits \(((\Spaf)\op)_X\) as the \emph{tensor} of \((\Spaf)\op\) over \(X\) in \(\Catx\), see Remark~\refone{remark:universal-underlying} below.
\end{proof}

\subsection{Internal functor categories}
\label{subsection:internal}%
In the present section we will show that the symmetric monoidal structures on \(\Cath\) and \(\Catp\) constructed in~\S\refone{subsection:monoidal-structure} are \emph{closed}. In particular, for two hermitian \(\infty\)-categories \((\C,\QF)\) and \((\Ctwo,\QFtwo)\), we will promote the stable \(\infty\)-category \(\Funx(\C,\Ctwo)\) of exact functors from \(\C\) to \(\Ctwo\) to a hermitian \(\infty\)-category \(\Funx((\C,\QF),(\Ctwo,\QFtwo))\),
characterized by a natural equivalence
\[
\Map_{\Cath}((\Cthree,\QFthree),\Funx((\C,\QF),(\Ctwo,\QFtwo))) \simeq \Map_{\Cath}((\Cthree,\QFthree)\otimes (\C,\QF),(\Ctwo,\QFtwo))
\]
for \((\Cthree,\QFthree) \in \Cath\). Furthermore, if \((\C,\QF)\) and \((\Ctwo,\QFtwo)\) are Poincaré \(\infty\)-categories, then the hermitian \(\infty\)-category \(\Funx((\C,\QF),(\Ctwo,\QFtwo))\) is Poincaré as well and satisfies
\[
\Map_{\Catp}((\Cthree,\QFthree),\Funx((\C,\QF),(\Ctwo,\QFtwo))) \simeq \Map_{\Catp}((\Cthree,\QFthree)\otimes (\C,\QF),(\Ctwo,\QFtwo)).
\]
The internal functor category \(\Funx((\C,\QF),(\Ctwo,\QFtwo))\) enjoys the following useful property: its hermitian objects correspond exactly to the hermitian functors \((\C,\QF) \to (\Ctwo,\QFtwo)\), and when \((\C,\QF)\) and \((\Ctwo,\QFtwo)\) are Poincaré the Poincaré objects in
\(\Funx((\C,\QF),(\Ctwo,\QFtwo))\) correspond to Poincaré functors \((\C,\QF) \to (\Ctwo,\QFtwo)\). In particular, this allows one to view the notions of hermitian functors and hermitian objects in a unified setting, and describe the relation between hermitian and Poincaré functors in terms of that which holds between hermitian and Poincaré objects.

\begin{definition}
Let \((\C,\QF)\) and \((\Ctwo,\QFtwo)\) be two hermitian \(\infty\)-categories. We set
\[
\iq \colon \Funx(\C,\Ctwo)\op \to \Spa, \quad f \mapsto \nat(\QF,f^*\QFtwo),
\]
where \(\nat\) denotes the \emph{spectrum} of natural transformations between two spectrum valued functors.
\end{definition}

\begin{proposition}
\label{proposition:basic-properties-hermitian-functor-cats}%
The functor \(\iq\colon \Funx(\C,\Ctwo)\op \to \Spa\) is quadratic. Its bilinear part is given by
\[
\Bil_\iq(f,g) := \nat(\Bil_\QF, (f \times g)^*\Bil_\QFtwo)
\]
and its linear part makes the diagram
\[
\begin{tikzcd}
\Lin_\iq(f) \ar[rrr] \ar[d] & & & \nat(\Lin_\QF, f^*\Lin_\QFtwo) \ar[d] \\
 \Bil_\iq(f,f)^\tC \ar[r,phantom,"{\simeq}"] & \nat(\QF,f^*\Bil_{\QFtwo}^{\Delta})^\tC \ar[rr] & & \nat\left(\QF,f^*(\Bil_{\QFtwo}^\Delta)^\tC\right),
\end{tikzcd}
\]
cartesian, where the left bottom equivalence is by the adjunction of Lemma~\refone{lemma:universal-crs}. Here, as in \S\refone{subsection:quadratic}, for a bilinear form \(\Bil\) we denote by \(\Bil^{\Delta} = \Del^*\Bil\) its pre-composition with the diagonal.
If both \(\QF\) and \(\QFtwo\) are non-degenerate, so is \(\iq\) with duality given by
\[
\Dual_\iq(f) = \Dual_\QFtwo f\op\Dual_\QF\op.
\]
Finally, if both \(\QF\) and \(\QFtwo\) are perfect, then so is \(\iq\).
\end{proposition}

\begin{proof}
We begin by computing the cross effect of \(\iq\). It is given by
\begin{align*}
\Bil_{\iq}(f,g) = &\fib\Big[\nat(\QF,(f\oplus g)^*\QFtwo) \longrightarrow \nat(\QF,f^*\QFtwo) \oplus \nat(\QF,g^*\QFtwo)\Big] \\
\simeq & \nat\left(\QF,\fib\big[(f\oplus g)^*\QFtwo \longrightarrow f^*\QFtwo \oplus g^*\QFtwo \big]\right) \\
\simeq & \nat\left(\QF,((f \times g)^*\Bil_{\QFtwo})^{\Del}\right) \simeq \nat\left(\Bil_{\QF},(f \times g)^*\Bil_{\QFtwo}\right)
\end{align*}
where the last equivalence is by the adjunction of Lemma~\refone{lemma:universal-crs}.
In particular, the cross effect of \(\iq\) is bilinear. Taking \(f=g\) we similarly get that
\begin{equation}
\label{equation:bil-iq-del}%
\Bil_{\iq}^{\Del}(f) \simeq \nat\left(\Bil_{\QF},(f \times f)^*\Bil_{\QFtwo}\right) \simeq \nat(\QF,f^*\Bil^{\Del}_{\QF'}),
\end{equation}
and hence
\begin{equation}
\label{equation:nat-colinear-part}%
\fib\big[\iq(f) \to (\Bil_{\iq}^{\Del}(f))^{\hC}\big] \simeq \nat(\QF,f^*\fib[\QF' \to (\Bil^{\Del}_{\QF'})^{\hC}]).
\end{equation}
Since \(\fib[\QFtwo \to \Bil_{\QF}^{\Del}]\) is exact by Proposition~\refone{proposition:basic-properties-quad-functors} it then follows that~\eqrefone{equation:nat-colinear-part} is exact in \(f\), and hence by the same proposition we have that \(\iq\) is quadratic.

We now compute the linear part of \(\iq\). Applying
\(\nat(\QF,f^*(-))\) to the classifying square of \(\QFtwo\) and using the equivalence~\eqrefone{equation:bil-iq-del} and the equivalence \(\nat(\Lin_\QF,f^*\Lin_\QFtwo) \xrightarrow{\simeq} \nat(\QF,f^*\Lin_\QFtwo )\) given by Lemma~\refone{lemma:linear} we obtain an exact square
\[
\begin{tikzcd}
\iq(f) \ar[r] \ar[d] & \nat(\Lin_{\QF},f^*\Lin_{\QFtwo}) \ar[d] \\
\Bil_{\iq}(f,f)^{\hC} \ar[r] & \nat(\QF,f^*(\Bil_{\QFtwo}^{\Del})^{\tC})
\end{tikzcd}
\]
in which the terms on the left hand side are quadratic in \(f\) and the terms on the right hand side are exact in \(f\). The desired formula for \(\Lin_{\iq}\) is now obtained by taking linear parts.

Now suppose that \(\QF\) and \(\QFtwo\) are non-degenerate. Applying Lemma~\refone{lemma:nat-duality}  and the adjunction between \(\Dual_{\QF}\) and \(\Dual_{\QF}\op\) we obtain natural equivalences
\[
\Bil_{\iq}(f,g) \simeq \nat(\Bil_{\QF},(f \times g)^*\Bil_{\QFtwo}) \simeq \nat(f\Dual_{\QF}, \Dual_{\QFtwo}g\op) \simeq \nat(f,\Dual_{\QFtwo}g\op\Dual_{\QF}\op),
\]
which shows that \(\iq\) is non-degenerate with duality \(\Dual_{\iq}(f) \simeq \Dual_{\QF}f\op\Dual_{\QF}\op\).
If \(\QF\) and \(\QFtwo\) are in addition perfect then the dualities \(\Dual_\QF\) and \(\Dual_\QFtwo\) are equivalences and hence so is \(\Dual_{\nat_{\QF}^\QFtwo}\) by the above formula. The hermitian structure \(\nat_\QF^\QFtwo\) is then also perfect, as desired.
\end{proof}

\begin{definition}
For hermitian \(\infty\)-categories \((\C,\QF),(\Ctwo,\QFtwo)\) we will denote by
\[
\Funx((\C,\QF),(\Ctwo,\QFtwo)) := (\Funx(\C,\Ctwo),\iq)
\]
the hermitian \(\infty\)-category given by Proposition~\refone{proposition:basic-properties-hermitian-functor-cats}. We will refer to \(\Funx((\C,\QF),(\Ctwo,\QFtwo))\) as the \defi{internal functor category} from \((\C,\QF)\) to \((\Ctwo,\QFtwo)\).
\end{definition}

\begin{remark}
\label{remark:internal-poincare}%
If \((\C,\QF)\) and \((\Ctwo,\QFtwo)\) are Poincaré then \(\Funx((\C,\QF),(\Ctwo,\QFtwo))\) is Poincaré. This follows from the last part of Proposition~\refone{proposition:basic-properties-hermitian-functor-cats}.
\end{remark}

\begin{example}
\label{example:internal-universal}%
Let \((\C,\QF)\) be a hermitian \(\infty\)-category. Then under the natural equivalence \(\Funx(\Spaf,\C)\simeq\C\) given by evaluation at the sphere spectrum, the functor \(\nat_{\QF^{\uni}}^\QF\) corresponds to \(\QF\) by virtue of Lemma~\refone{lemma:maps-from-Qg}. We consequently obtain a natural equivalence
\[
\Funx((\Spaf,\QF^{\uni}),(\C,\QF)) \simeq (\C,\QF)
\]
of hermitian \(\infty\)-categories.
\end{example}

\begin{construction}
\label{construction:evaluation}%
For two hermitian \(\infty\)-categories \((\C,\QF),(\Ctwo,\QFtwo)\) we construct a hermitian functor
\[
\ev\colon (\C,\QF) \otimes \Funx((\C,\QF),(\Ctwo,\QFtwo))\to (\Ctwo,\QFtwo)
\]
as follows. We first observe that by definition of the monoidal structure on \(\Cath\), specifying such a hermitian functor
is
equivalent to specifying a functor \(f\colon\C \times \Funx(\C,\Ctwo) \to \Ctwo\) which is exact in each variable separately, together with a natural transformation
\[
\eta\colon \QF \otimes \iq \Rightarrow f^*\QFtwo
\]
of functors \(\C\op \times \Funx(\C,\Ctwo)\op \to \Spa\).
In the case at hand we then take as \(f\) the usual evaluation functor \(\ev \colon\C \times \Funx(\C,\Ctwo) \to \Ctwo\). To construct \(\eta\) we use the curry-equivalence
\[
\Fun(\C\op \times \Funx(\C,\Ctwo)\op,\Spa) \simeq \Fun(\Funx(\C,\Ctwo)\op, \Fun(\C\op,\Spa)),
\]
and pull the evaluation transformation
\[
\QF \otimes \nat(\QF,-) \Rightarrow \id_{\Fun(\C\op,\Spa)}
\]
back along the postcomposition functor
\[
\QFtwo \circ (-)\colon \Funx(\C,\Ctwo)\op \to \Fun(\C\op,\Spa) ,
\]
which is the curry of \(\ev^*\QFtwo\).
\end{construction}

\begin{proposition}
\label{proposition:internal}%
Let \((\C,\QF)\) and \((\Ctwo,\QFtwo)\) be two hermitian \(\infty\)-categories. Then the evaluation functor \((\ev,\eta)\) of Construction~\refone{construction:evaluation} exhibits \(\Funx((\C,\QF),(\Ctwo,\QFtwo))\) as the internal mapping object in the symmetric monoidal \(\infty\)-category \(\Cath\). That is, for any \((\Cthree,\QFthree)\) hermitian \(\infty\)-category the composite
\begin{multline*}
\Map_{\Cath}((\Cthree,\QFthree),\Funx((\C,\QF),(\Ctwo,\QFtwo))) \\
\longrightarrow \Map_{\Cath}((\C,\QF) \otimes (\Cthree,\QFthree),(\C,\QF) \otimes \Funx((\C,\QF),(\Ctwo,\QFtwo))) \\
\xrightarrow{(\ev,\eta)_\ast} \Map_{\Cath}((\C,\QF) \otimes (\Cthree,\QFthree),(\Ctwo,\QFtwo))
\end{multline*}
is an equivalence of spaces.
\end{proposition}

By the Yoneda lemma we immediately find:

\begin{corollary}
\label{corollary:herm-fun-cat-functorial}%
The association \((\C,\QF),(\Ctwo,\QFtwo) \mapsto \Funx((\C,\QF),(\Ctwo,\QFtwo))\) canonically extends to a functor \((\Cath)\op \times \Cath \rightarrow \Cath\) in a way that renders the evaluation map of Construction~\refone{construction:evaluation} a natural transformation.
\end{corollary}

\begin{corollary}
The symmetric monoidal structures on \(\Cath\) constructed in \S\refone{subsection:monoidal-structure} is closed.
\end{corollary}

\begin{corollary}
\label{corollary:tensor-preserve-colimit}%
The symmetric monoidal product on \(\Cath\) preserves small colimits in each variable.
\end{corollary}

\begin{remark}
\label{remark:forgetful-closed}%
Since the hermitian evaluation functor of Construction~\refone{construction:evaluation} refines by definition the usual evaluation functor of \(\Catx\), it follows that the forgetful functor \(\rU\colon \Cath \to \Catx\) is not only symmetric monoidal (Theorem~\refone{theorem:tensorpoincare}\refoneitem{item:thm-two}), but also \defi{closed} symmetric monoidal, that is, for \((\C,\QF),(\Ctwo,\QFtwo) \in \Cath\) the composed map
\[
\rU(\C,\QF) \otimes \rU\Funx((\C,\QF),(\Ctwo,\QFtwo)) \xrightarrow{\simeq} \rU((\C,\QF) \otimes  \Funx((\C,\QF),(\Ctwo,\QFtwo)) \xrightarrow{\rU(\ev)} \rU(\Ctwo,\QFtwo)
\]
exhibits \(\rU\Funx((\C,\QF),(\Ctwo,\QFtwo))\) as the internal mapping object in \(\Catx\) from \(\rU(\C,\QF)=\C\) to \(\rU(\Ctwo,\QFtwo)=\Ctwo\). In other words, the transposed map
\[
\rU\Funx((\C,\QF),(\Ctwo,\QFtwo)) \to \Funx(\rU(\C,\QF),\rU(\Ctwo,\QFtwo))
\]
is an equivalence.
\end{remark}

\begin{proof}[Proof of Proposition~\refone{proposition:internal}]
The forgetful functor \(\Cath \to \Catx\) induces a commutative diagram
\begin{equation}
\label{equation:forget}%
\begin{tikzcd}
\Map_{\Cath}((\Cthree,\QFthree),\Funx((\C,\QF),(\Ctwo,\QFtwo))) \ar[r] \ar[d] &  \Map_{\Cath}((\C,\QF) \otimes (\Cthree,\QFthree),(\Ctwo,\QFtwo)) \ar[d] \\
\Map_{\Catx}(\Cthree,\Funx(\C,\Ctwo)) \ar[r,"{\simeq}"] & \Map_{\Catx}(\C \otimes \Cthree,\Ctwo)
\end{tikzcd}
\end{equation}
in which the bottom horizontal arrow is an equivalence since on the level of stable \(\infty\)-categories, the bilinear functor \(\ev\colon \C \times \Funx(\C,\Ctwo) \to \Ctwo\) already exhibits \(\Funx(\C,\Ctwo) \in \Catx\) as the internal mapping object from \(\C\) to \(\Ctwo\) in \(\Catx\). It will hence suffice to show that the map induced by~\eqrefone{equation:forget} on vertical homotopy fibres is an equivalence. Let us thus fix a linear functor \(g\colon \Cthree \to \Funx(\C,\Ctwo)\), and write \(\ovl{g}\colon\C \otimes \Cthree \to \Ctwo\) for its image in the bottom right corner of~\eqrefone{equation:forget}. Since \(\Cath \to \Catx\) is a cartesian fibration the induced map on vertical fibres in~\eqrefone{equation:forget} can be identified with the composed map of spaces of natural transformations
\begin{equation}
\label{equation:fibers}%
\Nat(\QFthree,g^*\nat^{\QFtwo}_{\QF}) \to \Nat(\QF \otimes \QFthree,\QF \otimes g^*\nat^{\QFtwo}_{\QF}) \to \Nat(\QF \otimes \QFthree,\ovl{g}^*\QFtwo) ,
\end{equation}
where the second map is induced by the natural transformation \(\eta\colon\QF \otimes \nat^{\QFtwo}_{\QF} \Rightarrow \ev^*\QFtwo\) of Construction~\refone{construction:evaluation},
restricted along \((\id,g)\colon\C \times \Cthree \to \C \times \Funx(\C,\Ctwo)\). Let us now identify functors \((\C \times \Cthree)\op \to \Spa\) with functors \(\Cthree\op \to \Fun(\C\op,\Spa)\), and write \(\ovl{\QFtwo}\colon \Funx(\C,\Ctwo)\op \to \Fun(\C\op,\Spa)\) for the curried functor determined by \(\ev^*\QFtwo\colon(\C \times \Funx(\C,\Ctwo))\op \to \Spa\). Unwinding the definitions
we may rewrite~\eqrefone{equation:fibers} as
\begin{multline*}
\Map_{\Fun(\Cthree\op,\Spa)}(\QFthree,\nat(\QF,g^*\ovl{\QFtwo})) \\
\to \Map_{\Fun(\Cthree\op,\Fun(\C\op,\Spa))}(\QF \otimes \QFthree,\QF \otimes\nat(\QF,g^*\ovl{\QFtwo}))) \\\to \Map_{\Fun(\Cthree\op,\Fun(\C\op,\Spa))}(\QF \otimes \QFthree,g^*\ovl{\QFtwo}),
\end{multline*}
where we understand \(\nat(\QF,-)\colon\Fun(\Cthree\op,\Fun(\C\op,\Spa)) \to \Fun(\Cthree\op,\Spa)\) as the functor obtained by applying \(\nat(\QF,-)\colon\Fun(\C\op,\Spa) \to \Spa\) levelwise.
To finish the proof it will hence suffice to show that the evaluation natural transformation
\[
\QF \otimes \nat(\QF,-) \Rightarrow \id_{\Fun(\C\op,\Spa)}
\]
exhibits
\[
\QF \otimes -\colon\Fun(\Cthree\op,\Spa) \to \Fun(\Cthree\op,\Fun(\C\op,\Spa))
\]
as left adjoint to
\[
\nat(\QF,-) \colon \Fun(\Cthree,\Fun(\C\op,\Spa)) \to \Fun(\Cthree\op,\Spa).
\]
Indeed, this is simply the adjunction induced by the canonical adjunction
\[
\QF \otimes (-) \colon \Spa \adj \Fun(\C\op,\Spa) \colon \nat(\QF,-)
\]
encoding the structure of \(\Fun(\C\op,\Spa)\) as tensored over \(\Spa\).
\end{proof}

Specializing Proposition~\refone{proposition:internal} to the unit Poincaré \(\infty\)-category \((\Cthree,\QFthree) = (\Spaf,\QF^{\uni})\) and using Proposition~\refone{proposition:corepresentability-of-poinc} we recover a natural equivalence
\begin{equation}
\label{equation:hermitian-objects-internal}%
\spsforms(\Funx((\C,\QF),(\Ctwo,\QFtwo))) \xrightarrow{\simeq} \Hom_{\Cath}((\C,\QF),(\Ctwo,\QFtwo)),
\end{equation}
which by comparing with the analogous claim for \(\Catx\) we can place in a commutative square
\[
\begin{tikzcd}
\spsforms(\Funx((\C,\QF),(\Ctwo,\QFtwo))) \ar[r,"{\simeq}"] \ar[d] & \Map_{\Cath}((\C,\QF),(\Ctwo,\QFtwo)) \ar[d] \\
\iota\Funx(\C,\Ctwo) \ar[r,"{\simeq}"] & \Map_{\Catx}(\C,\Ctwo)
\end{tikzcd}
\]
In particular, the equivalence~\eqrefone{equation:hermitian-objects-internal} is of a somewhat tautological nature: a hermitian object in \(\Funx((\C,\QF),(\Ctwo,\QFtwo))\) consists of a pair \((f,\eta)\), where \(f\colon \C \to \Ctwo\) is an exact functor and \(\eta \in \Om^{\infty}\iq(f)\) is a form for \(\iq\), that which by definition means a natural transformation \(\eta \colon \QF \Rightarrow f^*\QFtwo\). The equivalence~\eqrefone{equation:hermitian-objects-internal} then associates to this hermitian object the same pair \((f,\eta)\), now considered as a hermitian functor from \((\C,\QF)\) to \((\Ctwo,\QFtwo)\). Consulting the proof of Proposition~\refone{proposition:basic-properties-hermitian-functor-cats} we observe that when \((\C,\QF)\) and \((\Ctwo,\QFtwo)\) are non-degenerate
the map
\[
\eta_\sharp\colon f\to \Dual_\iq f = \Dual_\QFtwo f\op \Dual_\QF\op
\]
associated to a hermitian object \((f,\eta)\), corresponds to the map
\[
\tau_{\eta}\colon f \Dual_{\QF} \to \Dual_{\QFtwo}f\op
\]
(see Definition~\refone{definition:transformation-duality}) via the adjunction between pre-composition with \(\Dual_{\QF}\) and pre-composition with \(\Dual_{\QF}\op\). When \(\Dual_{\QF}\) is an equivalence we thus have that \(\eta_{\sharp}\) is an equivalence if and only if \(\tau_{\eta}\) is an equivalence. In particular, for Poincaré \(\infty\)-categories \((\C,\QF),(\Ctwo,\QFtwo)\) we have that a hermitian object \((f,\eta)\) in \(\Funx((\C,\QF),(\Ctwo,\QFtwo))\) is Poincaré if and only if the corresponding hermitian functor \((f,\eta)\colon (\C,\QF) \to (\Ctwo,\QFtwo)\) is Poincaré. In particular, in this case the equivalence~\eqrefone{equation:hermitian-objects-internal} restricts to an equivalence
\begin{equation}
\label{equation:poincare-objects-internal}%
\Poinc(\Funx((\C,\QF),(\Ctwo,\QFtwo))) \xrightarrow{\simeq} \Hom_{\Catp}((\C,\QF),(\Ctwo,\QFtwo)).
\end{equation}
We may summarize the situation as follows:

\begin{corollary}
\label{corollary:forms-in-functor-cats}%
For a hermitian \(\infty\)-category \((\C,\QF)\), the hermitian object \(\tf_{(\C,\QF)} \) in \(\Funx((\C,\QF),(\C,\QF))\) corresponding to the identity under the equivalence~\eqrefone{equation:hermitian-objects-internal}, exhibits \((\C,\QF)\) as corepresenting the functor
\[
(\Ctwo,\QFtwo) \mapsto \spsforms(\Funx((\C,\QF),(\Ctwo,\QFtwo)))
\]
in \(\Cath\). In addition, if \((\C,\QF)\) is Poincaré then \(\tf_{(\C,\QF)}\) is Poincaré and exhibits \((\C,\QF)\) as corepresenting the functor
\[
(\Ctwo,\QFtwo) \mapsto \Poinc(\Funx((\C,\QF),(\Ctwo,\QFtwo)))
\]
in \(\Catp\).
\end{corollary}

\begin{corollary}
\label{corollary:internal-enhanced}%
Let \((\C,\QF),(\Ctwo,\QFtwo)\) and \((\Cthree,\QFthree)\) be hermitian \(\infty\)-categories. Then there is an equivalence of hermitian \(\infty\)-categories
\begin{equation}
\label{equation:internal-enhanced}%
\Funx((\Cthree,\QFthree),\Funx((\C,\QF),(\Ctwo,\QFtwo))) \simeq \Funx((\C,\QF) \otimes (\Cthree,\QFthree),(\Ctwo,\QFtwo)).
\end{equation}
which is natural in \((\C,\QF),(\Ctwo,\QFtwo),(\Cthree,\QFthree)\).
In addition, for \((\Cthree,\QFthree) = \Funx((\C,\QF),(\Ctwo,\QFtwo))\) this equivalence sends the hermitian object on the left hand side corresponding to the identity to the hermitian object on the right side corresponding to the evaluation functor of Construction~\refone{construction:evaluation}.
\end{corollary}
\begin{proof}
This is a formal consequence of the Yoneda lemma that holds true in any closed symmetric monoidal \(\infty\)-category. Indeed, embedding both sides in presheaves to spaces we may use Proposition~\refone{proposition:internal} and Corollary~\refone{corollary:herm-fun-cat-functorial} to construct natural equivalences between the resulting presheaves
\[
\Map_{\Cath}((-),\Funx((\Cthree,\QFthree),\Funx((\C,\QF),(\Ctwo,\QFtwo))) \simeq \Map_{\Cath}((\Cthree,\QFthree) \otimes (-),\Funx((\C,\QF),(\Ctwo,\QFtwo))) \simeq
\]
\[
\Map_{\Cath}((\C,\QF) \otimes (\Cthree,\QFthree) \otimes (-),(\Ctwo,\QFtwo)) \simeq \Map_{\Cath}(-,\Funx((\C,\QF) \otimes (\Cthree,\QFthree),(\Ctwo,\QFtwo)))
\]
The additional claim in the case \((\Cthree,\QFthree) = \Funx((\C,\QF),(\Ctwo,\QFtwo))\) can be obtained by taking \((\D,\QFD) = (\Spaf,\QF^{\uni})\) and tracing through the equivalences on both sides.
\end{proof}

Taking Poincaré objects in~\eqrefone{equation:internal-enhanced} and using the equivalence~\eqrefone{equation:poincare-objects-internal}
we thus conclude:

\begin{corollary}
\label{corollary:internal-poincare}%
Let \((\C,\QF)\) and \((\Ctwo,\QFtwo)\) be two Poincaré \(\infty\)-categories. Then the evaluation functor \((\ev,\eta)\) of Construction~\refone{construction:evaluation} is Poincaré and exhibits \((\Funx(\C,\Ctwo),\nat_\QF^\QFtwo)\) as the internal mapping objects in the symmetric monoidal category \(\Catp\). In particular, it determines an equivalence of spaces
\[
\Map_{\Catp}((\Cthree,\QFthree),\Funx((\C,\QF),(\Ctwo,\QFtwo))) \simeq \Map_{\Catp}((\C,\QF) \otimes (\Cthree,\QFthree),(\Ctwo,\QFtwo))
\]
for \((\Cthree,\QFthree) \in \Catp\).
\end{corollary}

\begin{corollary}
\label{corollary:catp-closed}%
The symmetric monoidal structure on \(\Catp\) constructed in \S\refone{subsection:monoidal-structure} is closed.
\end{corollary}

\begin{corollary}
\label{corollary:poinc-fun-cat-functorial}%
The association \((\C,\QF),(\Ctwo,\QFtwo) \mapsto \Funx((\C,\QF),(\Ctwo,\QFtwo))\) canonically extends to a functor \((\Catp)\op \times \Catp \rightarrow \Catp\) in a way that renders the evaluation map of Construction~\refone{construction:evaluation} a natural transformation.
\end{corollary}

\begin{corollary}
\label{corollary:poinc-tensor-preserve-colimit}%
The symmetric monoidal product on \(\Catp\) preserves small colimits in each variable.
\end{corollary}

\begin{remark}
\label{remark:inclusion-closed}%
Since the evaluation functor used in Corollary~\refone{corollary:internal-poincare} to exhibit internal functor categories in \(\Catp\) is the same as the one that was used to exhibit internal functor categories in \(\Cath\) it follows formally that the inclusion \(\Catp \hrar \Cath\) is not only symmetric monoidal (Theorem~\refone{theorem:tensorpoincare}\refoneitem{item:thm-three}), but also \defi{closed} symmetric monoidal, that is, preserves internal functor categories. It then follows from Remark~\refone{remark:forgetful-closed} that the composed functor
\[
\Catp \to \Cath \to \Catx
\]
is closed symmetric monoidal as well.
\end{remark}

\begin{definition}
For hermitian \(\infty\)-categories \((\C,\QF)\) and \((\Ctwo,\QFtwo)\) we define the \emph{category of hermitian functors} from \((\C,\QF)\) to \((\Ctwo,\QFtwo)\) to be
\[
\Funh((\C,\QF),(\Ctwo,\QFtwo)) := \catforms(\Funx((\C,\QF),(\Ctwo,\QFtwo))).
\]
\end{definition}

\begin{remark}
Combining Corollary~\refone{corollary:internal-enhanced} and Corollary~\refone{corollary:forms-in-functor-cats} we deduce that for hermitian \(\infty\)-categories \((\C,\QF),(\Ctwo,\QFtwo)\) and \((\Cthree,\QFthree)\), the evaluation functor of Construction~\refone{construction:evaluation} determines a natural equivalence
\[
\Funh((\Cthree,\QFthree),\Funx((\C,\QF),(\Ctwo,\QFtwo))) \simeq \Funh((\C,\QF) \otimes (\Cthree,\QFthree),(\Ctwo,\QFtwo)).
\]
\end{remark}

\begin{remark}
\label{remark:composition-functors}%
Proposition~\refone{proposition:internal} tells us that \(\Cath\) is a closed monoidal category, so we can turn it into an \(\infty\)-category enriched over itself via \cite[Cor.~7.4.10]{Gepner-Haugseng}. In particular, for three hermitian \(\infty\)-categories \((\C,\QF), (\Ctwo,\QFtwo)\) and \((\Cthree,\QFthree)\)
we have natural composition hermitian functors
\[
\Funx((\C,\QF),(\Ctwo,\QFtwo)) \otimes \Funx((\Ctwo,\QFtwo), (\Cthree,\QFthree)) \to  \Funx((\C,\QF),(\Cthree,\QFthree)).
\]
which one can of course also easily write down without using enriched technology.
Applying \cite[Cor.~5.7.6]{Gepner-Haugseng} to the lax monoidal functor \(\catforms\colon\Cath\to \Cat\) we see that \(\Cath\) is canonically endowed with an enrichment over \(\Cat\), with composition functors
\[
\Funh((\C,\QF),(\Ctwo,\QFtwo)) \times \Funh((\Ctwo,\QFtwo), (\Cthree,\QFthree)) \to \Funh((\C,\QF),(\Cthree,\QFthree)),
\]
and identities given by the tautological hermitian objects \(\tf_{(\C,\QF)}\) of Corollary~\refone{corollary:forms-in-functor-cats}.
In particular, one should consider \(\Cath\) as an \((\infty,2)\)-category, with mapping categories given by \(\Funh(-,-)\).
Though we will not make explicit use of this point of view, these hermitian functor categories will play a role in \S\refone{subsection:cotensoring}-\refone{subsection:tensoring} when we study the tensor-cotensor constructions, a structure most naturally viewed with the \((\infty,2)\)-categorical perspective in mind.
\end{remark}

\subsection{Cotensoring of hermitian categories}
\label{subsection:cotensoring}%

Given a hermitian \(\infty\)-category \((\C,\QF)\) and an \(\infty\)-category \(\I\), our goal in this section is to promote the diagram category \(\C^{\I} := \Fun(\I,\C)\) to a hermitian \(\infty\)-category \((\C,\QF)^{\I}\), which we call the \defi{cotensor of \((\C,\QF)\) by \(\I\)}.
We will characterize the resulting hermitian \(\infty\)-category by a universal property, see Proposition~\refone{proposition:cotensor} below, which can be considered as witnessing it being the cotensor structure of \(\Cath\) over \(\Cat\) with respect to the enrichment of the former in the latter described in \S\refone{subsection:internal}.
This construction will feature prominently in \papertwo via the hermitian \(\Q\)-construction.

\begin{construction}
Let \(\I\) be a small \(\infty\)-category and \((\C,\QF)\) a hermitian \(\infty\)-category. We will denote by \(\C^{\I} := \Fun(\I,\C)\) the stable \(\infty\)-category of functors \(\I \to \C\). Let \(\ev\colon \I \times \C^{\I} \to \C\) be the evaluation functor, which, under the exponential equivalence
\[
\Fun(\I \times \C^{\I},\C) \simeq \Fun(\I,\Fun(\C^{\I},\C))
\]
corresponds to the functor which associates to \(i \in \I\) the evaluation-at-\(i\) functor \(\ev_i\colon \C^{\I} \to \C\).
Define a functor \(\QF^{\I}\colon (\C^{\I})\op \to \Spa\) by
\[
\QF^{\I} := \lim_{i \in \I\op} \ev_i^*\QF.
\]
On a given diagram \(\vphi \in \C^{\I}\) the functor \(\QF^{\I}\) is given by the formula
\[
\QF^{\I}(\vphi) = \lim_{i\in\I\op}\QF(\vphi(i)).
\]
\end{construction}

\begin{proposition}
\label{proposition:basic-properties-herm-diagrams}%
The functor \(\QF^{\I}\colon \C^{\I} \to \Spa\) is quadratic. Its bilinear part is given by
\[
\Bil^{\I}(\vphi,\psi) := \lim_{i \in \I\op}\Bil_\QF(\vphi(i),\psi(i)),
\]
and the linear part \(\Lin^{\I}\) makes the square
\begin{equation}
\label{equation:linear-part-diagrams}%
\begin{tikzcd}
[column sep=3ex]
\Lin^{\I}(\vphi) \ar[rrr] \ar[d] & & & \displaystyle\mathop{\lim}_{i \in \I\op} \Lin_\QF(\vphi(i)) \ar[d] \\
 \Bil^{\I}(\vphi,\vphi)^\tC \ar[r,equal] & \left[\lim_{i \in \I\op} \Bil_{\QF}(\vphi(i),\vphi(i))\right]^\tC \ar[rr] & & \displaystyle\mathop{\lim}_{i \in \I\op} \Bil_{\QF}(\vphi(i),\vphi(i))^\tC
\end{tikzcd}
\end{equation}
cartesian. If \(\QF\) is non-degenerate and \(\C\) admits \((\I_{i/})\op\)-shaped limits for all \(i \in \I\), then \(\QF^{\I}\) is also non-degenerate with duality given by
\[
\big[\Dual^{\I}(\vphi)\big](i) = \lim_{[i \to j] \in (\I_{i/})\op} \Dual_{\QF}(\vphi(j)) .
\]
\end{proposition}

\begin{remark}
If \(\I\) is a finite \(\infty\)-category then \(\Lin^{\I}(\vphi) = \lim_{\I\op} \Lin_\QF\vphi\op\), since the bottom horizontal map in~\eqrefone{equation:linear-part-diagrams} is then an equivalence.
\end{remark}

\begin{proof}[Proof of Proposition~\refone{proposition:basic-properties-herm-diagrams}]
The functor \(\QF^{\I}\) is defined as a limit of the functors \(\ev_i^*\QF\), each of which is quadratic since \(\QF\) itself is quadratic and each \(\ev_i\) is exact. Since the collection of quadratic functors is closed under limits (Remark~\refone{remark:closed}) it follows that \(\QF^{\I}\) is quadratic. Since the formation of bilinear parts is compatible with restriction (Remark~\refone{remark:invariance-base-change}) and commutes with limits (e.g, by Lemma~\refone{lemma:universal-crs}), it follows that
\[
\Bil^{\I} \simeq \lim_{i \in \I\op}(\ev_i^* \times \ev_i^*)\Bil_{\QF},
\]
and in particular
\[
\Bil^{\I}(\vphi,\vphi) \simeq \lim_{i \in \I\op}\Bil_{\QF}(\vphi(i),\vphi(i)) .
\]
The formation of linear parts however does not commute with limits. To compute it, we apply \(\lim_{\I\op}\) to the square classifying \(\QF\) via Corollary~\refone{corollary:classification-of-quad-functors}, yielding the square
\[
\begin{tikzcd}
\QF^{\I}(\vphi) \ar[r] \ar[d] & \displaystyle\mathop{\lim}_{i \in \I\op}\Lin_{\QF}(\vphi(i)) \ar[d] \\
\Bil^{\I}(\vphi,\vphi)^{\hC} \ar[r] & \displaystyle\mathop{\lim}_{i \in \I\op}\Bil_{\QF}(\vphi(i),\vphi(i))^{\tC}
\end{tikzcd}
\]
of quadratic functors in \(\vphi\) whose left hand side consists of exact functors. Taking linear parts we then get the desired description of \(\Lin^{\I}\).

We now prove the desired formula for the duality. For this, we henceforth assume that \((\C,\QF)\) is non-degenerate and that \(\C\) admits \((\I_{/i})\op\)-indexed limits for every \(i \in \I\). Let \(s\colon\Twar(\I) \to \I\) and \(t\colon \Twar(\I) \to \I\op\) be the source and target functors, respectively.
Define \(\Dual^{\I}\colon (\C^{\I})\op \to \C^{\I}\) by the composite formula
\[
(\C^{\I})\op \xrightarrow{\Dual_{\QF}\circ(-)} \C^{\I\op} \xrightarrow{t^*} \C^{\Twar(\I)} \xrightarrow{s_*} \C^{\I}
\]
where \(s_*\) stands for right Kan extension. This right Kan extension indeed exists: since \(s\colon \Twar(\I) \to \I\) is a cartesian fibration classified by the functor \(i \mapsto (\I_{i/})\op\) this right Kan extension is given by the explicit formula
\[
[\Dual^{\I}(\vphi)](i) = \lim_{[i \to j] \in (\I_{i/})\op}\Dual_{\QF}(\vphi(j))
\]
where the required limits exist in \(\C\) by assumption. We now claim that \(\Dual^{\I}\) represents the bilinear functor \(\Bil^{\I}\). To prove this, note first that since the right Kan extension functor \(s_*\) is right adjoint to the corresponding restriction functor \(s^*\) we get that
\begin{multline*}
\nat(\vphi,\Dual^{\I}(\psi)) = \nat(\vphi,s_*t^*\Dual_{\QF}\psi) \simeq \nat(s^*\vphi,t^*\Dual_{\QF}\psi) \\
\simeq \lim_{\tiny \begin{matrix}[\sig\colon \alp \Rightarrow \beta] \in \\
\Twar(\Twar(\I))\op \end{matrix}} \map_{\C}(\vphi(s\alp),\Dual_{\QF}\psi(t\beta)) \simeq \lim_{\tiny \begin{matrix} [\sig\colon \alp \Rightarrow \beta] \in \\ \Twar(\Twar(\I))\op \end{matrix}} \Bil_{\QF}(\vphi(s\alp),\psi(t\beta)),
\end{multline*}
where we have used the standard formula for the spectrum of natural transformations as a limit over the twisted arrow category.
We wish to show that the last limit above is equivalent to \(\lim_{i \in \I\op} \Bil_{\QF}(\vphi(i),\psi(i))\). For this we will make use of several cofinality arguments. To facilitate readability in what follows, we invite the reader to visualize an object \([\sig\colon \alp \Rightarrow \beta] \in \Twar(\Twar(\I))\op\) as a diagram of the form
\[
\begin{tikzcd}
i \ar[r,"{\alp}"] \ar[d] & l  \\
j \ar[r,"{\bet}"] & k \ar[u]
\end{tikzcd}
\]
To begin, consider the commutative diagram
\[
\begin{tikzcd}
[column sep=1ex]
{[i \to j \to k \to l]} \ar[d,mapsto] \ar[r,phantom,"{\in}"] & \Twar(\Twar(\I)) \ar[r,"{\Twar(\source \times \target)}"] \ar[d,"\target"'] &[5ex] \Twar(\I \times \I\op) \ar[d] \ar[r,equal] &[0ex] \Twar(\I) \times \Twar(\I\op) \ar[d,"{\target \times \target}"] & {([i \to j],[l \leftarrow k])} \ar[l,phantom,"{\ni}"] \ar[d,mapsto] \\
{[j \to k]} \ar[r,phantom,"{\in}"] & \Twar(\I)\op \ar[r,"{(\source \times \target)\op}"] & (\I \times \I\op)\op \ar[r,equal] & \I\op \times \I & (j,k) \ar[l,phantom,"{\ni}"]
\end{tikzcd}
\]
Since \((s \times t)\colon\Twar(\C) \to \C \times \C\op\) is a right fibration it induces an equivalence on over categories. It then follows that the commutative square on the left is cartesian, and hence the induced map
\[
\Twar(\Twar(\I)) \xrightarrow{\simeq} \Twar(\I)\op \times_{[\I\op \times \I]}[\Twar(\I) \times \Twar(\I\op)] \simeq
\Twar(\I) \times_{\I\op} \Twar(\I)\op \times_{\I} \Twar(\I\op)
\]
is an equivalence. In particular, the projection
\begin{equation}
\label{equation:projection-twar}%
\Twar(\Twar(\I)) \to \Twar(\I) \times_{\I\op} \Twar(\I)\op
\end{equation}
is a cartesian fibration whose fibres have terminal objects, being pulled back from the (target) cartesian fibration \(\Twar(\I\op) \to \I\) given by \([l \leftarrow k] \mapsto k\) which has this property. By~\cite[Lemma 4.1.3.2]{HTT} it then follows that the functor~\eqrefone{equation:projection-twar} is cofinal. To avoid confusion, we note that in the fibre product in~\eqrefone{equation:projection-twar}, the map \(\Twar(\I) \to \I\op\) is the target projection \([i \to j] \mapsto j\) and the map \(\Twar(\I)\op \to \I\op\) is the opposite of the source projection \([j \to k] \mapsto j\).
Now the composed functor
\begin{equation}
\label{equation:projection-twar-2}%
\Twar(\I) \times_{\I\op} \Twar(\I)\op \to \Twar(\I) \times_{\I\op} [\I\op \times \I] = \Twar(\I) \times \I \to \I
\end{equation}
is cocartesian fibration, being a composition of a left fibration and a constant cocartesian fibration. By compatibility with base change we see that this cocartesian fibration is classified by the functor \(k \mapsto \Twar(\I) \times_{\I\op} (\I_{/k})\op\). Since the projection \(\I_{/k} \to \I\) is a right fibration it follows by the same argument as above that the map \(\Twar(\I_{/k}) \to \Twar(\I) \times_{\I\op} (\I_{/k})\op\) (induced by the target projection) is an equivalence. The cocartesian~\eqrefone{equation:projection-twar-2} is hence also classified by the equivalent functor \(k \mapsto \Twar(\I_{/k})\). Now consider the map of cocartesian fibrations (over \(\I\))
\begin{equation}
\label{equation:big-fibration}%
\Twar(\I) \times_{\I\op} \Twar(\I)\op \simeq \int_{\I}\Twar(\I_{/k}) \to \int_{\I} \I_{/k} \simeq \Ar(\I)
\end{equation}
induced by the source projections \(\Twar(\I_{/k}) \to \I_{/k}\). Then~\eqrefone{equation:big-fibration} is a map of cocartesian fibrations which is fibrewise cofinal by Lemma~\cite[Lemma 4.1.3.2]{HTT} and is hence itself cofinal. Since cofinal maps are closed under composition we may now conclude that the composed projection
\[
\Twar(\Twar(\I)) \to \Ar(\I)
\]
is cofinal. On the other hand the canonical inclusion \(\I \to \Ar(\I)\) sending \(\x\) to \(\id_x\) is also cofinal since it has a left adjoint (the target functor). We may therefore conclude that
\[
\lim_{\tiny \begin{matrix} [\sig\colon \alp \Rightarrow \beta] \in \\ \Twar(\Twar(\I))\op \end{matrix}} \Bil_{\QF}(\vphi(s\alp),\psi(t\beta)) \simeq \lim_{[i \to k] \in \Ar(\I)\op} \Bil_{\QF}(\vphi(i),\psi(k)) \simeq \lim_{i \in \I\op} \Bil_{\QF}(\vphi(i),\psi(i)) ,
\]
and so \(\Dual_{\I}\) represents \(\Bil_{\QF}\), as desired.
\end{proof}

\begin{definition}
For a hermitian \(\infty\)-category \((\C,\QF)\) and an \(\infty\)-category \(\I\) we will denote by \((\C,\QF)^{\I} := (\C^{\I},\QF^{\I})\) the \(\infty\)-category given by Proposition~\refone{proposition:basic-properties-herm-diagrams}. We will refer to it as the \defi{cotensor of \((\C,\QF)\) by \(\I\)}.
\end{definition}

\begin{remark}
If \(\I\) is a finite poset then the comma \(\infty\)-categories \(\I_{/i}\) are finite for every \(i \in \I\), and hence every stable \(\infty\)-category admits \((\I_{/i})\op\) indexed limits. In particular, in this case \((\C,\QF)^{\I}\) is non-degenerate as soon as \((\C,\QF)\) is non-degenerate.
\end{remark}

\begin{example}
\label{example:diagram-twar-del-1}%
For \(\I = \Twar(\Delta^1)\) we may identify \(\C^{\I}\) with the \(\infty\)-category of spans \(\x \leftarrow \cob \rightarrow \y\) in \(\C\), with \(\QF^{\I}\) given by
\[
\QF^{\I}([\x \leftarrow \cob \rightarrow \y]) = \QF(\x) \times_{\QF(\cob)} \QF(\y) ,
\]
and the duality (when \((\C,\QF)\) is non-degenerate) given by
\[
\Dual^{\I}([\x \leftarrow \cob \rightarrow \y]) = [\Dual_{\QF}\x \leftarrow\Dual_\QF \x \times_{\Dual_\QF \cob} \Dual_\QF \y \rightarrow \Dual_{\QF}\y].
\]
It is then straightforward to verify that this duality is perfect whenever \(\Dual_{\QF}\) is perfect, in which case \((\C,\QF)^{\I}\) is Poincaré. This example will feature prominently in subsequent parts of the present paper in the context of the cobordism category of a Poincaré \(\infty\)-category, where we will view the above duality as an algebraic incarnation of Lefschetz duality for manifolds.
\end{example}

\begin{warning}
\label{warning:not-poincare}%
For a Poincaré \((\C,\QF)\) and \(\I\) arbitrary, the hermitian \(\infty\)-category \((\C,\QF)^{\I}\) might fail to be Poincaré, even if \(\I\) is a finite poset. This happens for example if \(\I\) has a final object but is not itself equivalent to a point; indeed, in this case the image of \(\Dual^{\I}\) is the full subcategory of \(\C^\I\) spanned by the constant diagrams. On the other hand,
we will see in \S\refone{subsection:tensoring} that cotensor by \(\I\) does preserve Poincaré \(\infty\)-categories when \(\I\) is the poset of faces of a finite simplicial complex.
\end{warning}

We go on to establish the universal property of these hermitian diagram categories, from which we will also deduce their functoriality. We will require the following lemma:

\begin{lemma}
\label{lemma:mapping-into-Funh}%
Let \(\I\) be a small \(\infty\)-category and \(\C,\Ctwo\) be two hermitian \(\infty\)-categories. Then the fibre of the map
\[
\Fun(\I,\Funh((\C,\QF),(\Ctwo,\QFtwo)))\to \Fun(\I,\Funx(\C,\Ctwo))
\]
over a functor \(f\colon\I\times\C\to\Ctwo\) is naturally equivalent to the space \(\Nat(p_{\C}^*\QF, f^*\QFtwo)\)
of natural transformations \(p^*_{\C}\QF\Rightarrow f^*\QFtwo\), where \(p_{\C}\colon \I \times \C \to \C\) denotes the projection to \(\C\).
\end{lemma}
\begin{proof}
We want to describe the space of dotted lifts
\[
\begin{tikzcd}
 & \Funh(\C,\Ctwo) \ar[d] \\
\I \ar[r,pos=.7,"f"] \ar[ur,dashed] & \Funx(\C,\Ctwo)
\end{tikzcd}
\]
Recall that the vertical map above
is the right fibration classified by
functor \(\Omega^\infty\iq\colon\Funx(\C,\Ctwo)\op\to \Sps\), and so by \cite[Corollary 3.3.3.2]{HTT}, the space of sections of this right fibration coincides with the limit
\[
\lim_{i\in \I\op} \Omega^\infty\iq(f_i) = \lim_{i\in \I\op} \Nat(\QF,f_i^*\QFtwo) \simeq \Nat\big(\QF,\lim_{i\in \I\op}f^*_i\QFtwo\big).
\]
On the other hand, we have
\[
\Nat(p^*_{\C}\QF,f^*\QFtwo)\simeq \Nat(\QF,p_\ast f^*\QFtwo)\simeq \Nat(\QF,\lim_{i\in \I\op}f_i^*\QFtwo),
\]
since right Kan extensions along \(p_{\C}\) are computed by taking the limit fibrewise, as can be seen from the pointwise formula for right Kan extensions. Hence the two constructions are naturally equivalent.
\end{proof}

\begin{construction}
\label{construction:evaulation-cotensor}%
For a hermitian \(\infty\)-category \((\C,\QF)\) and an \(\infty\)-category \(\I\) we define a functor
\[
\ovl{\ev}\colon \I \to \Funh((\C,\QF)^{\I},(\C,\QF))
\]
as follows.
Let \(p \colon \I\op \times (\C^{\I})\op \to (\C^{\I})\op\) be the projection on the second factor and
\[
\ev: \I\op \times (\C^{\I})\op \to \C\op
\]
the evaluation. By Lemma~\refone{lemma:mapping-into-Funh} the additional data needed in order to define \(\ovl{\ev}\) is a natural transformation
\begin{equation}
\label{equation:tau-cotensor}%
\tau\colon p^*\QF^{\I} \Rightarrow \ev^*\QF.
\end{equation}
We then define \(\tau\) by taking the counit transformation
\[
\mathop{\const}_{\I\op} \lim_{\I\op} \Rightarrow \id_{\Fun(\I\op,\Spa)},
\]
currying it into morphism in \(\Fun(\I\op \times \Fun(\I\op,\Spa),\Spa)\) and finally pre-composing with the functor
\[
\I\op \times (\C^{\I})\op \to \I\op \times \Fun(\I\op,\Spa)
\]
induced by \(\QF\).
\end{construction}

\begin{proposition}
\label{proposition:cotensor}%
Let \((\C,\QF)\) and \((\Ctwo,\QFtwo)\) be two hermitian \(\infty\)-categories and \(\I\) a small \(\infty\)-category. Then the composite map
\[
\Funh((\Ctwo,\QFtwo),(\C,\QF)^\I) \times \I
\xrightarrow{\id \times \ovl{\ev}} \Funh((\Ctwo,\QFtwo),(\C,\QF)^\I) \times \Funh((\C,\QF)^\I,(\C,\QF))
\longrightarrow \Funh((\Ctwo,\QFtwo),(\C,\QF))
\]
defined using the functor \(\ovl{\ev}\) of Construction~\refone{construction:evaulation-cotensor} and the composition functor of Remark~\refone{remark:composition-functors}, determines an equivalence of \(\infty\)-categories
\begin{equation}
\label{equation:cotensor}%
\Funh((\Ctwo,\QFtwo),(\C,\QF)^\I) \xrightarrow{\simeq} \Fun(\I,\Funh((\Ctwo,\QFtwo),(\C,\QF)))\,,
\end{equation}
and in particular an equivalence
\begin{equation}
\label{equation:cotensor-mapping-spaces}%
\Hom_{\Cath}((\Ctwo,\QFtwo),(\C,\QF)^\I) \simeq \iota\Fun(\I,\Funh((\Ctwo,\QFtwo),(\C,\QF))).
\end{equation}
\end{proposition}

We will give the proof of Proposition~\refone{proposition:cotensor} at the end of this subsection. Before let us explore some of its consequences. First, as in the case of internal functor categories, the Yoneda lemma immediately implies:

\begin{corollary}
\label{corollary:functorial-cotensor}%
The association \((\I,(\C,\QF)) \mapsto (\C,\QF)^\I\) extends canonically to a functor \(\Cat\op \times \Cath \rightarrow \Cath\) that rendered the equivalence from proposition~\refone{proposition:cotensor} natural.
\end{corollary}

\begin{remark}
\label{remark:restriction}%
Unwinding the definitions, if \(\alp\colon \I \to \J\) is a functor of small \(\infty\)-categories and \((\C,\QF)\) is a hermitian \(\infty\)-category then the hermitian functor
\[
(\alp^*,\eta^{\alp}) \colon (\C,\QF)^{\J} \to (\C,\QF)^{\I}
\]
issued via the functoriality of Corollary~\refone{corollary:functorial-cotensor} is given by the usual restriction functor \(\alp^*\colon \C^{\J} \to \C^{\I}\) on the underlying stable \(\infty\)-categories accompanied by the usual restriction-induced map
\[
\eta^{\alp}_{\vphi}\colon \QF^{\J}(\vphi) = \lim_{j \in \J\op} \QF(\vphi(j)) \to \lim_{i \in \I\op} \QF(\vphi(\alp(i))) = \QF^{\I}(\alp^*\vphi)
\]
on limits.
\end{remark}

\begin{remark}
As pointed out in Warning~\refone{warning:not-poincare}, the cotensor construction does \emph{not} restrict to a functor \(\Cat\op \times \Catp \rightarrow \Catp\). In particular, while this construction is best understood by considering \(\Cath\) as an \((\infty,2)\)-category, the \((\infty,2)\)-categorical perspective does not seem to extend to \(\Catp\) in a meaningful manner.
\end{remark}

\begin{remark}
\label{remark:groupoid-poincare}%
It follows from Proposition~\refone{proposition:cotensor} that when \(\I\) is an \(\infty\)-groupoid the cotensor \((\C,\QF)^{\I}\) coincides with the limit in \(\Cath\) of the constant \(\I\)-diagram with value \((\C,\QF)\). In particular, it follows form Proposition~\refone{proposition:Catp-cocomplete} that for such an \(\I\) the functor \((\C,\QF) \mapsto (\C,\QF)^{\I}\) does preserve Poincaré \(\infty\)-categories.
\end{remark}

Taking \((\Ctwo,\QFtwo) = (\Spa^\omega,\QF^{\uni})\) in Proposition~\refone{proposition:cotensor} yields:

\begin{corollary}
\label{corollary:forms-in-diagram-cats}%
There is a natural equivalence \(\catforms((\C,\QF)^{\I}) \simeq \Fun(\I,\catforms(\C,\QF))\).
\end{corollary}

\begin{remark}
We know of no analogous formula for the Poincaré objects of \((\C,\QF)^\I\) when the latter happens to be Poincaré. It is certainly \emph{not} true, for example, that the individual objects of a Poincaré diagram are Poincaré objects themselves, as demonstrated by the case \(\I = \Twar(\Delta^1)\), see Example~\refone{example:diagram-twar-del-1}.
\end{remark}

The following is again a formal consequence:

\begin{corollary}
For any hermitian \(\infty\)-categories \((\C,\QF), (\Ctwo,\QFtwo)\) and any \(\infty\)-category \(\I\) there is a canonical equivalence
\[
\Funx((\C,\QF),(\Ctwo, \QFtwo)^\I) \simeq \Funx((\C,\QF),(\Ctwo,\QFtwo))^{\I}
\]
of hermitian \(\infty\)-categories.
\end{corollary}

\begin{proposition}
\label{proposition:induced-maps-duality-pres}%
Let \((\C,\QF)\) be non-degenerate and \(\alpha \colon \I \rightarrow \J\) a functor between small categories, such that \((\C,\QF)\) admits both \((\I_{i/})\op\)- and \((\J_{j/})\op\)-shaped limits for all \(i\in \I\) and \(j \in \J\). If the induced maps \(\I_{i/} \rightarrow \J_{\alpha(i)/}\) are cofinal for every \(i \in \I\) then the hermitian functor
\[
(\alpha^*,\eta^{\alp}) \colon (\C,\QF)^{\J} \rightarrow (\C,\QF)^\I
\]
is duality preserving. In particular, if \((\C,\QF)^{\J}\) and \((\C,\QF)^{\I}\) are Poincaré then \((\alpha^*,\eta^{\alp})\) is a Poincaré functor.
\end{proposition}
\begin{proof}
This follows directly from the explicit description of the duality in Proposition~\refone{proposition:basic-properties-herm-diagrams}.
\end{proof}

\begin{remark}
In the situation of Proposition~\refone{proposition:induced-maps-duality-pres}, if \(\alpha\) is a map of posets, then the given criterion for preservation of duality can be rephrased more explicitly as by saying that for every \(i \in \I\) and \(j \in \J\) with \(j \geq \alpha(i)\), the realization of the poset \(\{k \in \I \mid i \leq k, j \leq \alpha(k)\}\) is contractible.
\end{remark}

\begin{proof}[Proof of Proposition~\refone{proposition:cotensor}]
We argue similarly to the proof of Proposition~\refone{proposition:internal}. The forgetful functor determines a commutative diagram of \(\infty\)-categories
\begin{equation}
\label{equation:square-cotensor}%
\begin{tikzcd}
\Funh((\Ctwo,\QFtwo),(\C,\QF)^\I) \ar[r] \ar[d] & \Fun(\I,\Funh((\Ctwo,\QFtwo),(\C,\QF))) \ar[d] \\
\Funx(\Ctwo,\C^\I) \ar[r,"{\simeq}"] & \Fun(\I,\Funx(\Ctwo,\C))
\end{tikzcd}
\end{equation}
in which the vertical maps are right fibrations and where the bottom arrow is an equivalence since the functor \(\ev':\I \to \Funx(\C^\I,\Ctwo)\) underlying \(\ovl{\ev}\) already exhibits \(\C^\I\) as the cotensor of \(\C\) over \(\I\) in \(\Catx\). It will hence suffice to show the map induced by~\eqrefone{equation:square-cotensor} on vertical fibres is an equivalence. Let us hence fix an exact functor \(g\colon\Ctwo \to \C^\I\) and let \(g_\I := \{g_i\}: \I \to \Funx(\Ctwo,\C)\) be its image in the bottom right corner of~\eqrefone{equation:square-cotensor}. Now the fibre of the right vertical map in~\eqrefone{equation:square-cotensor} over \(g_\I\) is  the space of sections of the base change
\begin{equation}
\label{equation:base-change}%
\Funh((\Ctwo,\QFtwo),(\C,\QF)) \times_{\Funx(\Ctwo,\C)}\I \to \I,
\end{equation}
where the fibre product is taken with respect to the map \(g_\I\). By the compatibility of base change and straightening we see that~\eqrefone{equation:base-change} is the right fibration classified by the functor \(i \mapsto \Nat(\QFtwo,g_i^*\QF)\). By~\cite[Corollary 3.3.3.2]{HTT} evaluation at the various \(i \in \I\) exhibits the space of sections of~\eqrefone{equation:base-change} as the limit \(\lim_i \Nat(\QFtwo,g_i^*\QF)\). We may then identify the map induced by~\eqrefone{equation:square-cotensor} from the fibre over \(g\) to the fibre of \(g_\I\) with the map
\[
\Nat(\QFtwo,g^*\QF^\I) \to \lim_{i \in \I}\Nat(\QFtwo,g_i^*\QF)
\]
whose components \(\Nat(\QFtwo,g^*\QF^\I) \to \Nat(\QFtwo,g^*_i\QF)\) are induced by the components \(\tau_i\colon\QF^\I \rightarrow \ev_i^*\QF\) of~\eqrefone{equation:tau-cotensor}. Since pulling back functors preserve limits the desired result now follows from the fact that the collection of maps \(\QF^\I \to \ev_i^*\QF\) exhibit \(\QF^\I\) as the limit of the diagram \(\{\ev_i^*\QF\}\), by definition.
\end{proof}

\subsection{Tensoring of hermitian categories}
\label{subsection:tensoring}%

In this section we will consider the dual of the cotensor construction studied in \S\refone{subsection:cotensoring}, which we will refer to as \emph{tensoring} a Poincaré \(\infty\)-category \((\C,\QF)\) by an \(\infty\)-category \(\I\). In general this construction is somewhat less accessible then the cotensor construction, but we will be able to say more about it when \(\I\) satisfies certain finiteness conditions, see \S\refone{subsection:finite-tensors-cotensors} below. We will exploit the tensor construction in \papertwo in order to form the \emph{dual} \(\Q\)-construction, which is needed in the proof of the universal property of the Grothendieck-Witt spectrum.

\begin{construction}
\label{construction:tensor}%
Let \((\C,\QF)\) be a hermitian \(\infty\)-category and \(\I\) a small \(\infty\)-category. For \(i \in \I\) and \(\x \in \C\), let us denote by \(R_{i,\x} \colon \I\op \to \Pro(\C)\) the functor \(R_{i,\x} = (\iota_i)_*(\x)\) right Kan extended along the inclusion \(\iota_i\colon \{i\} \hrar \I\op\) of \(i\) from the functor \(\{i\} \to \Pro(\C)\) with value \(\x \in \C \subseteq \Pro(\C)\). We then let \(\C_\I \subseteq \Pro(\C)^{\I\op}\) be the smallest full subcategory containing \(R_{i,\x}\) for \(i \in \I\) and \(\x \in \C\) and closed under finite limits. Then \(\C_\I\) is also closed under suspensions (since the collection \(R_{i,\x}\) is, as suspension in \(\C\) commutes with finite limits) and is hence stable. It is also equipped by construction with a functor
\begin{equation}
\label{equation:tensor}%
\iota\colon \C \times \I \to \C_\I \quad\quad (\x,i) \mapsto R_{i,\x} .
\end{equation}
We then promote \(\C_\I\) to a hermitian \(\infty\)-category by endowing it with the quadratic functor
\[
\QF_\I\colon\C_\I\op \to \Spa
\]
obtained by taking the left Kan extension of \(p_{\C}^*\QF\colon \C\op \times \I\op \to \Spa\) along \(\iota\op\colon \C\op \times \I\op \to \C_{\I}\op\),
which results in a reduced functor, and then applying to it the \(2\)-excisive approximation of Construction~\refone{construction:excisive-approx}, left adjoint to the inclusion \(\Funq(\C_{\I}) \subseteq \Fun_{\ast}(\C_{\I}\op,\Spa)\) (in fact, we will see in the proof of Proposition~\refone{proposition:tensor-explicit} that this 2-excisive approximation is not needed, that is, the result of the left Kan extension is already 2-excisive). Here we denote by \(p_{\C}\colon \C \times \I \to \C\) projection to \(\C\). We then set
\[
(\C,\QF)_\I := (\C_\I,\QF_\I) ,
\]
and refer to it as the \defi{tensor of \((\C,\QF)\) by \(\I\)}. By construction the functor \(\QF_{\I}\) supports a natural transformation
\[
p_{\C}^*\QF\Rightarrow \iota^*\QF_\I,
\]
where \(p_{\C}\colon \C\times \I \to \C\) is the projection to \(\C\), and by Lemma~\refone{lemma:mapping-into-Funh} this transformations determines a functor
\begin{equation}
\label{equation:tensor-2}%
\coev\colon \I \to \Funh((\C,\QF),(\C_\I,\QF_{\I})),
\end{equation}
which to \(i \in \I\) associates the exact functor \(\x \mapsto R_{i,\x},\) equipped with the natural transformation \(\QF(\x) \Rightarrow \QF_\I(R_{i,\x})\) given by the construction of \(\QF_{\I}\).
\end{construction}

\begin{remark}
\label{remark:universal-underlying}%
The functor \(\coev\colon \C \times \I \to \C_{\I}\) exhibits \(\C_{\I}\) as universal among stable \(\infty\)-categories equipped with a functor from \(\C \times \I\) which is exact in the first entry. To see this, we may replace the term ``exact'' by ``finite limit preserving''. In other words, it will suffice to show that~\eqrefone{equation:tensor} is universal among maps from \(\C \times \I\) to a finitely complete \(\infty\)-category which preserve finite limits in the first variable. Such universal constructions are explicitly described in \cite[\S 5.3.6]{HTT}. In particular, it will suffice to show that \(\C_\I\) coincides with the construction appearing in the proof of \cite[Proposition~5.3.6.2]{HTT}. To see this, note that \(\Fun(\I\op,\Pro(\C))\) identifies with the full subcategory of \(\Fun(\C \times \I,\Sps)\op\) spanned by those functors which preserve finite limits in the \(\C\)-variable. This full inclusion then admits a right adjoint
\[
R\colon \Fun(\C \times \I,\Sps)\op \to \Fun(\I\op,\Pro(\C)),
\]
and what we need to check is that \(\C_\I \subseteq \Fun(\I\op,\Pro(\C))\) identifies with the full subcategory generated under finite limits by the images under \(R\) of the corepresentable functors \(\Map_{\C \times \I}((x,i),-)\colon \C \times \I \to \Sps\). In light of the definition of \(\C_\I\) this amounts to showing that \(R_{i,x} \simeq R(\Map_{\C \times \I}((x,i),-))\). Indeed, both these objects represent the functor
\[
\Fun(\I\op,\Pro(\C)) \to \Sps \quad\quad \vphi \mapsto \Map_{\Pro(\C)}(\vphi(i),x).
\]
\end{remark}

To describe \(\QF_{\I}\) more explicitly, let \(\wtl{\QF}\colon \Pro(\C)\op \to \Spa\) be the left Kan extension of \(\QF\) along the Yoneda embedding \(\C\op \hrar \Pro(\C)\op\). Then \(\wtl{\QF}\) is quadratic by~\cite[Proposition 6.1.5.4]{HA}. Its bilinear part then coincides with the essentially unique bilinear functor \(\wtl{\Bil}_{\QF}\colon \Pro(\C)\op \times \Pro(\C)\op \to \Spa\) which extends \(\Bil\) and preserve colimits in each variable separately, and its linear part is the essentially unique colimit preserving functor \(\wtl{\Lin}_{\QF}\colon\Pro(\C)\op \to \Spa\) extending \(\Lin_{\QF}\).

\begin{proposition}
\label{proposition:tensor-explicit}%
Let \((\C,\QF)\) be a hermitian \(\infty\)-category and \(\I\) a small \(\infty\)-category. Then the quadratic functor \(\QF_\I\) of Construction~\refone{construction:tensor} is given explicitly by the formula
\[
\QF_{\I}(\vphi) = \mathop{\colim}_{i \in \I} \wtl{\QF}(\vphi(i)) .
\]
Its bilinear and linear parts are given by
\[
\Bil_{\I}(\vphi,\psi) := \mathop{\colim}_{i \in \I} \wtl{\Bil}_{\QF}(\vphi(i),\psi(i)) \quad\text{and}\quad \Lin_{\I}(\vphi) := \mathop{\colim}_{i \in \I} \wtl{\Lin}_{\QF}(\vphi),
\]
respectively.
\end{proposition}
\begin{proof}
To establish the formula for \(\QF_\I\) we first note that
\(\colim_{i \in \I}\wtl{\QF}(\vphi(i))\) is quadratic, being a colimit of the quadratic functors \(\ev_i^*\wtl{\QF}\) for \(i \in \I\op\), see Remark~\refone{remark:closed}. Its linear and bilinear parts are then given by the indicated formulas since taking linear and bilinear parts commutes with colimits.
It will hence suffice to identify \(\colim_{i \in \I}\wtl{\QF}(\vphi(i))\) with the left Kan extension of \(p_{\C}^*\QF\)
along \(\iota\op\).
For this, consider the commutative square
\[
\begin{tikzcd}
(\C \times \I)\op \ar[r,"{\iota\op}"] \ar[d,"j"'] & \C_{\I}\op \ar[d] \\
\Fun(\C \times \I,\Sps) \ar[r] & \Fun(\I\op,\Pro(\C))\op
\end{tikzcd}
\]
where the left vertical map \(j\) is the Yoneda embedding and the bottom horizontal map is the left adjoint to the inclusion \(\Fun(\I\op,\Pro(\C))\op \hrar \Fun(\C \times \I,\Sps)\) induced by the inclusion \(\Pro(\C) \hrar \Fun(\C,\Sps)\op\) as the full subcategory spanned by left exact functors. Since the right vertical map is fully-faithful we may compute \(\iota_!p_{\C}^*\QF\) by further Kan extending to \(\Fun(\I\op,\Pro(\C))\op\) and then restricting back to \(\C_{\I}\op\). By the commutativity of the above square the left Kan extension to \(\Fun(\I\op,\Pro(\C))\op\) can be performed by first left Kan extending to \(\Fun(\C \times \I,\Sps)\) and then left Kan extending to \(\Fun(\I\op,\Pro(\C))\op\), the latter given by restriction along the right adjoint \(\Fun(\I\op,\Pro(\C))\op \hrar \Fun(\C \times \I,\Sps)\). Now the left Kan extension of \(p_{\C}^*\QF\) along the Yoneda embedding results in the coend construction
\[
\Fun(\C \times \I,\Sps) \ni \rho  \mapsto  \int_{\C \times \I} \rho \otimes p_{\C}^*\QF \simeq \int_{\C} (p_{\C})_!\rho \otimes \QF,
\]
where \(\otimes\) denotes the tensor of spectra over spaces.
Now, in the case where our functor \(\C \times \I \to \Sps\) is of the form \(\rho_{\vphi}(\x,i) = \Map_{\C}(\vphi(i),\x)\) for some \(\Pro(\C)\)-valued presheaf \(\vphi\colon \I\op \to \Pro(\C)\),
then its left Kan extension along \(p_{\C}\colon \C \times \I \to \C\) is given by
\[
(p_\C)_!\rho_{\vphi}(\x) = \colim_{i \in \I}\Map_{\Pro(\C)}(\vphi(i),\x) .
\]
We may then conclude that
\[
[\iota_!p_{\C}^*\QF](\vphi) = \int_{\C}(p_\C)_!\rho_{\vphi} \otimes \QF \simeq \colim_{i \in \I}\int_{\C}\Map_{\Pro(\C)}(\vphi(i),\x) \otimes \QF \simeq \colim_{i \in \I}\wtl{\QF}(\vphi(i)).
\]
As this expression already 2-excisive in \(\vphi\) we deduce that
\[
\QF_{\I}(\vphi) = \colim_{i \in \I}\wtl{\QF}(\vphi(i)),
\]
as desired.
\end{proof}

We shall now address the universal property of the tensor construction.

\begin{proposition}
\label{proposition:tensor-mapping-property}%
Let \((\C,\QF)\) and \((\Ctwo,\QFtwo)\) be hermitian \(\infty\)-categories and \(\I\) a small \(\infty\)-category. Then the composed map
\[
\I \times \Funh((\C,\QF)_\I,(\Ctwo,\QFtwo)) \xrightarrow{\coev \times \id} \Funh((\C,\QF),(\C,\QF)_\I) \times \Funh((\C,\QF)_\I,(\Ctwo,\QFtwo)) \to \Funh((\C,\QF),(\Ctwo,\QFtwo))
\]
induces an equivalence
\begin{equation}
\label{equation:tensor-mapping-categories}%
\Funh((\C,\QF)_\I,(\Ctwo,\QFtwo)) \simeq \Fun(\I,\Funh((\C,\QF),(\Ctwo,\QFtwo)) .
\end{equation}
and in particular an equivalence
\begin{equation}
\label{equation:tensor-mapping-spaces}%
\Hom_{\Cath}((\C,\QF)_{\I},(\Ctwo,\QFtwo)) \simeq \iota\Fun(\I,\Funh((\C,\QF),(\Ctwo,\QFtwo))).
\end{equation}

\end{proposition}

As in the case of the cotensor construction the universal characterization implies functoriality:
\begin{corollary}
\label{corollary:functorial-tensor}%
The association \((\I,(\C,\QF)) \mapsto (\C,\QF)_\I\) extends canonically to a functor \(\Cat \times \Cath \rightarrow \Cath\) that rendered the equivalence from proposition~\refone{proposition:tensor-mapping-property} natural.
\end{corollary}

\begin{remark}
\label{remark:comparing-1}%
Comparing universal properties, we see that there are canonical equivalences of hermitian \(\infty\)-categories
\[
\Funx((\C,\QF)_\I,(\Ctwo,\QFtwo)) \simeq \Funx((\C,\QF),(\Ctwo,\QFtwo))^\I \simeq \Funx((\C,\QF),(\Ctwo,\QFtwo)^\I)\,.
\]
\end{remark}

\begin{remark}
\label{remark:comparing-2}%
Comparing universal properties we see that there are canonical equivalences of hermitian \(\infty\)-categories
\[
(\C,\QF)_\I \otimes (\Ctwo,\QFtwo) \simeq (\C,\QF) \otimes (\Ctwo,\QFtwo)_{\I} \simeq ((\C,\QF) \otimes (\Ctwo,\QFtwo))_{\I} \,.
\]
\end{remark}

\begin{remark}
\label{remark:I-preserve-poincare}%
It follows from Remarks~\refone{remark:comparing-1},~\refone{remark:comparing-2} and~\refone{remark:internal-poincare} that for a given small \(\infty\)-category \(\I\) the conditions
\begin{enumerate}
\item
the functor \((\C,\QF) \mapsto (\C,\QF)_{\I}\) preserves Poincaré \(\infty\)-categories;
\item
the hermitian \(\infty\)-category \((\Spaf,\QF^{\uni})_{\I}\) is Poincaré;
\end{enumerate}
are equivalent, and that when these equivalent conditions hold the functor \((\C,\QF) \mapsto (\C,\QF)^{\I}\) preserves Poincaré \(\infty\)-categories as well.
\end{remark}

\begin{remark}
\label{remark:right-kan-extension}%
If \(\alp\colon \I \to \J\) is a map between small \(\infty\)-categories then the hermitian functor
\[
\C_\I \to \C_\J
\]
resulting from the functoriality of Corollary~\refone{corollary:functorial-tensor} must induce the associated restriction functor
\[
\alp^*\colon\Fun(\J,\Funh((\C,\QF),(\Ctwo,\QFtwo))) \to \Fun(\I,\Funh((\C,\QF),(\Ctwo,\QFtwo)))
\]
under the equivalence of Proposition~\refone{proposition:tensor-mapping-property}, upon mapping into any \((\Ctwo,\QFtwo)\). The underlying exact functor is consequently the essentially unique one (see Remark~\refone{remark:universal-underlying}) making the diagram
\[
\begin{tikzcd}
\C\times\I \ar[d,"\coev"]\ar[r, "\id \times \alp"] & \C\times\J\ar[d,"\coev"] \\
\C_\I \ar[r,] & \C_\J
\end{tikzcd}
\]
commute, and must therefore coincide with the restriction to \(\C_\I\) of the \emph{right Kan extension} functor
\[
\alp_*\colon \Fun(\I\op,\Pro(\C)) \to \Fun(\J\op,\Pro(\C)).
\]
By a slight abuse of notation we will denote this restriction by
\(\alp_*\colon \C_{\I} \to \C_{\J}\) as well. Using the formula of Proposition~\refone{proposition:tensor-explicit} the hermitian structure on \(\alp_*\) is then given by the natural map
\[
\QF_{\I}(\vphi) = \mathop{\colim}_{i \in \I}\wtl{\QF}(\vphi(i)) \to \mathop{\colim}_{i \in \I}\wtl{\QF}(\alp^*\alp_*\vphi(i)) \to \mathop{\colim}_{j \in \J}\wtl{\QF}(\alp_*\vphi(j)) = \QF_{\I}(\alp_*\vphi)
\]
for \(\vphi \in \C_{\I}\).
\end{remark}

\begin{remark}
\label{remark:tensor-cotensor}%
The natural equivalences~\eqrefone{equation:cotensor-mapping-spaces} and~\eqrefone{equation:tensor-mapping-spaces} exhibit \((-)_\I\) as left adjoint to \((-)^\I\). Furthermore, if \(\alp\colon\I \to \J\) is a map of finite posets then this adjunction intertwines the restriction functor \(\alp^*\colon(\C,\QF)^\J \to (\C,\QF)^\I\) with the (restricted) right Kan extension functor \(\alp_*\colon(\C,\QF)_\I \to (\C,\QF)_\J\).
\end{remark}

\begin{proof}[Proof of Proposition~\refone{proposition:tensor-mapping-property}]
The forgetful functor \(\Cath \to \Catx\) determines a commutative square of \(\infty\)-categories
\begin{equation}
\label{equation:square-tensor}%
\begin{tikzcd}
\Funh((\C,\QF)_\I,(\Ctwo,\QFtwo)) \ar[r] \ar[d] & \Fun(\I,\Funh((\C,\QF),(\Ctwo,\QFtwo))) \ar[d] \\
\Funx(\C_\I,\Ctwo) \ar[r] & \Fun(\I,\Funx(\C,\Ctwo))
\end{tikzcd}
\end{equation}
in which the vertical maps are right fibrations. By Remark~\refone{remark:universal-underlying} the bottom horizontal map is an equivalence.
It will hence suffice to show that~\eqrefone{equation:square-tensor} induces an equivalence on vertical fibres. Let \(g\colon \C_\I \to \Ctwo\) be an exact functor and let \(g_\I = \{g_i\}\colon \I \to \Funx(\C,\Ctwo)\) be its image in the bottom right corner of~\eqrefone{equation:square-tensor}. To identify the fibre on the left side, let \(g_0\colon \C \times \I \to \C_\I \st{g}{\to} \Ctwo\) be the composed functor, so that \(g_0\) corresponds to \(g_\I\) under the identification of functors \(\I \to \Fun(\C,\Ctwo)\) and functors \(\C \times \I \to \Ctwo\). In light of the definition of \(\QF_\I\) via left Kan extensions and \(2\)-excisive approximations we may identify the fibre of the left vertical arrow in~\eqrefone{equation:square-tensor} over \(g\) with \(\Nat_{\C \times \I}(p_{\C}^*\QF,g_0^*\QFtwo)\), where \(p_{\C}\colon\C \times \I \to \C\) is the projection. The map between the vertical fibres in~\eqrefone{equation:square-tensor} can then be identified with the map of spaces
\begin{equation}
\label{equation:map-fibers}%
\Nat_{\C \times \I}(p_{\C}^*\QF,g_0^*\QFtwo) \to \lim_{i \in \I\op}\Nat_{\C}(\QF,g_i^*\QFtwo) ,
\end{equation}
whose \(i\)'th component \(\Nat_{\C \times \I}(p_{\C}^*\QF,g_0^*\QFtwo) \to \Nat_{\C}(\QF,g_i^*\QFtwo)\) is given by restricting to \(\C \times \{i\}\subseteq \C \times \I\). This map is an equivalence by Lemma~\refone{lemma:mapping-into-Funh}, and so the proof is complete.
\end{proof}

\subsection{Finite tensors and cotensors}
\label{subsection:finite-tensors-cotensors}%

In this section we will consider the tensor and cotensor constructions in the case where the \(\infty\)-category \(\I\) satisfies strong finiteness conditions, e.g., when \(\I\) is a finite poset. In this case the tensor construction admits a more accessible description, and sends non-degenerate Poincaré \(\infty\)-categories to non-degenerate ones, with explicit induced duality, see Proposition~\refone{proposition:tensor-strongly-finite}. In addition, we will show that under these conditions the functor \((\C,\QF) \mapsto (\C,\QF)^{\I}\) is not only right adjoint to \((\C,\QF) \mapsto (\C,\QF)_{\I}\), but also \emph{left} adjoint to it, and extract some useful consequences. Finally,
in this case both the tensor and cotensor constructions are functorial not only in maps \(\alp\colon \I \to \J\), but also in cofinal maps \(\beta\colon \J \to \I\) going in the other direction, a phenomenon we refer to as \emph{exceptional functoriality}, see Construction~\refone{construction:exceptional}.

To begin, recall that an \(\infty\)-category \(\I\) is said to be \defi{finite} if it is categorically equivalent to a simplicial set with only finitely many non-degenerate simplices. If \(\I\) is a space then the condition that \(\I\) is finite as an \(\infty\)-category is equivalent to the condition that \(\I\) is finite as space, that is, that it is weakly equivalent to a simplicial set with finitely many non-degenerate simplices. We will use the term \defi{finite (co)limits} to refer to (co)limits indexed by finite \(\infty\)-categories. We recall that any stable \(\infty\)-category admits finite limits and colimits, and that these are preserved by any exact functor. In particular, in any stable \(\infty\)-category which admits small (co)limits, the latter automatically commute with finite (co)limits.

\begin{definition}
We will say that an \(\infty\)-category \(\I\) is \defi{strongly finite} if it is finite, and in addition for every \(i,j \in \I\) the mapping space \(\Map_{\I}(i,j)\) is finite.
\end{definition}

\begin{example}
Any finite poset is strongly finite.
\end{example}

\begin{example}
Any Reedy category with finitely many objects and finitely many morphisms is strongly finite. This follows by induction from~\cite[Proposition A.2.9.14]{HTT}. For example, the full subcategory \(\Del_{\leq n} \subseteq \Del\) spanned by the ordinals \([k]\) for \(k \leq n\) is strongly finite.
\end{example}

\begin{remark}
\label{remark:strongly-finite}%
If \(\I \to \J\) is a cartesian or cocartesian fibration such that \(\J\) is finite and the fibres \(\I_j\) are finite for every \(j \in \J\) then \(\I\) is finite. This follows from the explicit description of cartesian fibrations over the \(n\)-simplex via generalized mapping cones, see~\cite[\S 3.2.3]{HTT}. It then follows that for a strongly finite \(\infty\)-category \(\I\) the twisted arrow category \(\Twar(\I)\) is finite. Similarly, if \(\alp\colon \I \to \J\) is a functor between strongly finite \(\infty\)-categories then all the comma \(\infty\)-categories of \(\alp\) are finite.
\end{remark}

\begin{remark}
\label{remark:strongly-finite-2}%
Any localisation of a finite \(\infty\)-category by a finite set of arrows is finite, since it can be written as a pushout of finite \(\infty\)-categories. In particular, if \(\I\) is an \(\infty\)-category such that \(\Twar(\I)\) is finite then \(\I\) is finite, since \(\I\) can be written as a localisation of \(\Twar(\I)\) by a collection of arrows of the form \([f\colon \x \to \y] \to [\id\colon \x \to \x]\) where \(f\) runs over a set of representatives of equivalence types in \(\Twar(\I)\). Combining this with Remark~\refone{remark:strongly-finite} it follows that the condition that \(\I\) is strongly finite is equivalent to the condition that \(\Twar(\I)\) is finite and all mapping space in \(\I\) are finite.
\end{remark}

\begin{lemma}
\label{lemma:finite}%
Let \(\I\) be a small \(\infty\)-category and \(\C\) a stable \(\infty\)-category. Then the following holds:
\begin{enumerate}
\item
\label{item:map-finite}%
If the mapping spaces of \(\I\) are finite then \(\C_\I\) is contained in \(\Fun(\I\op,\C) \subseteq \Fun(\I\op,\Pro(\C))\).
\item
\label{item:twisted-arrow-finite}%
If the twisted arrow category \(\Twar(\I\op)\) is finite
then \(\Fun(\I\op,\C)\) is contained in \(\C_\I\).
\end{enumerate}
In particular, if \(\I\) is a strongly finite \(\infty\)-category then \(\C_{\I} = \Fun(\I\op,\C)\) as full subcategories of \(\Fun(\I\op,\Pro(\C))\).
\end{lemma}

\begin{remark}
The objects \(R_{\x,i}\) are cocompact in \(\Fun(\I\op,\Pro(\C))\) and generate it under limits. Since \(\C_{\I}\) is by definition the closure of \(R_{\x,i}\) under finite limits it follows that the inclusion \(\C_{\I} \subseteq \Fun(\I\op,\Pro(\C))\) induces an equivalence \(\Pro(\C_{\I}) \simeq \Fun(\I\op,\Pro(\C))\). When \(\I\) is strongly finite Lemma~\refone{lemma:finite} then gives an equivalence
\[
\Pro(\Fun(\I\op,\C)) \simeq \Fun(\I\op,\Pro(\C)).
\]
This generalizes~\cite[Proposition 5.3.5.15]{HTT} (in the case of \(\kappa=\omega\)) from finite posets to all strongly finite \(\infty\)-categories.
\end{remark}

\begin{proof}[Proof of Lemma~\refone{lemma:finite}]
To prove \refoneitem{item:map-finite} it will suffice to show that \(\C^{\I\op}\), which is closed under finite limits in \(\Pro(\C)^{\I\op}\), contains the objects \(R_{i,\x}\) for every \((i,x) \in \I \times \C\). Indeed \(R_{i,\x}(j) = x^{\Map_\I(i,j)}\) is contained in \(\C \subseteq \Pro(\C)\) since \(\x\) is in \(\C\), \(\Map_\I(i,j)\) is a finite space, and \(\C\) is closed inside \(\Pro(\C)\) under finite limits.

Let us now prove \refoneitem{item:twisted-arrow-finite}. We need to show that if \(\Twar(\I)\op\) is finite then any \(\C\)-valued presheaf \(\vphi\colon \I\op \to \C\) is a finite limit of cofree presheaves of the form \(R_{i,\x}\). But it is a standard fact that any presheaf \(\vphi\) is canonically the limit of the composed \(\Twar(\I)\op\)-indexed diagram
\[
\begin{tikzcd}
[row sep=1ex]
\Twar(\I)\op \ar[r] & \I\op \times \I \ar[r,"{\vphi \times \id}"] & \C \times \I \ar[r] & \Pro(\C)^{\I\op} \\
{[\alp\colon i \to j]} \ar[r,mapsto] & (i,j) \ar[r,mapsto] & (\vphi(i),j) \ar[r,mapsto] & R_{j,\vphi(i)}
\end{tikzcd}
\]
which takes values in cofree presheaves. To see this note that the \(\Twar(\I)\op\)-indexed family of maps \(c_{[i \to j]}\colon \vphi(j)\to \vphi(i)\) determines a \(\Twar(\I)\op\)-indexed family of maps \(\vphi \Rightarrow R_{j,\vphi(i)}\), and hence a map
\[
\vphi \Rightarrow \lim_{[i \to j] \in \Twar(\I)\op}R_{j,\vphi(i)}.
\]
Evaluating at \(k \in \I\op\), the resulting map
\[
\vphi(k) \Rightarrow \lim_{[i \to j] \in \Twar(\I)\op}\vphi(i)^{\Map_\I(j,k)} \simeq \lim_{[i \to j \to k] \in \Twar(\I_{/k})\op}\vphi(i)
\]
is then seen to be an equivalence by the cofinality of the functors \(\Twar(\I_{/k})\xrightarrow{\dom} \I_{/k} \leftarrow \{\id_k\}\).
\end{proof}

\begin{proposition}
\label{proposition:tensor-strongly-finite}%
Let \(\I\) be a strongly finite \(\infty\)-category (e.g., any finite poset). Then under the identification
\[
\C_\I = \Fun(\I\op,\C) \subseteq \Fun(\I\op,\Pro(\C))
\]
of Lemma~\refone{lemma:finite}, the quadratic functor \(\QF_\I\) corresponds to the functor
\[
\QF_\I(\vphi) = \mathop{\colim}_{i \in \I}\QF(\vphi(i)).
\]
Its bilinear and linear parts are then given by
\[
\Bil_{\I}(\vphi,\psi) = \mathop{\colim}_{i \in \I}\Bil_{\QF}(\vphi(i),\psi(i))  \quad\text{and}\quad  \Lin_{\I}(\vphi) =\mathop{\colim}_{i \in \I} \Lin_{\QF}(\vphi(i)),
\]
respectively. In addition, if \((\C,\QF)\) is non-degenerate
then \((\C,\QF)_\I\) is non-degenerate with duality
\begin{equation}
\label{equation:tensor-duality}%
[\Dual_{\I}\vphi](j) = \mathop{\colim}_{i\in \I}\Dual_{\QF}(\vphi(i))^{\Map_\I(i,j)}\,.
\end{equation}
\end{proposition}

\begin{proof}[Proof of Proposition~\refone{proposition:tensor-strongly-finite}]
The identification of \(\QF_{\I}\) together with its linear and bilinear parts follows directly from Proposition~\refone{proposition:tensor-explicit} and Lemma~\refone{lemma:finite}.
Now assume that \((\C,\QF)\) is non-degenerate.
To prove the formula for the duality, we need to show that for diagrams \(\vphi,\psi\colon \I\op \to \C\) there is an equivalence
\[
\Bil_\I(\vphi,\psi)\simeq \nat(\vphi,\Dual_\I\psi)
\]
natural in \(\vphi,\psi\), where \(\Dual_{\I}\) is given by ~\eqrefone{equation:tensor-duality}. Expanding the right hand side and using the standard formula for natural transformations we obtain
\begin{align*}
 \nat(\vphi,\Dual_\I\psi) &\simeq \lim_{[i\to j]\in\Twar(\I)\op} \map_\C(\vphi(j),(\Dual_\I\psi)(i))                                        &\\
                                &\simeq \lim_{[i\to j]\in \Twar(\I)\op} \map_\C\left(\vphi(j),\colim_{k\in \I}(\Dual_{\QF}\psi(k))^{\Map_\I(k,i)}\right)  &\\
                                &\simeq \colim_{k\in \I} \lim_{[i\to j]\in\Twar(\I)\op}\map_\C(\vphi(j),\Dual_{\QF}\psi(k))^{\Map_\I(k,i)}               &\\
                                &\simeq \colim_{k\in \I} \lim_{[i\to j]\in\Twar(\I_{k/})\op}\map_\C(\vphi(j),\Dual_{\QF}\psi(k))               &\\
                                &\simeq \colim_{k\in \I}\map_\C(\vphi(k),\Dual_{\QF}\psi(k))                                                          &\simeq \Bil_\I(\vphi,\psi),
\end{align*}
where we have used the finiteness of \(\Twar(\I)\) and \(\Map_\I(-,-)\) to commute limits and colimits and the cofinality of the maps \(\Twar(\I_{k/}) \xrightarrow{\cod} \I_{k/}\op \leftarrow \{\id_k\}\).
\end{proof}

We now turn our attention to some structural properties of the tensor and cotensor constructions which are special to the strongly finite case.
Recall from Remark~\refone{remark:tensor-cotensor} that for a fixed \(\infty\)-category \(\I\), the functor \((\C,\QF) \mapsto (\C,\QF)^{\I}\) is right adjoint to the functor \((\C,\QF) \mapsto (\C,\QF)_{\I}\). Our next goal is to show that when \(\I\) is strongly finite the functor \((\C,\QF) \mapsto (\C,\QF)^{\I}\) is also left adjoint to the functor \((\C,\QF) \mapsto (\C,\QF)_{\I}\). To exhibit this, consider for hermitian \(\infty\)-categories \((\C,\QF),(\Ctwo,\QFtwo)\) the evaluation hermitian functor
\begin{equation}
\label{equation:evaluation-internal}%
(\C,\QF) \otimes \Funx((\C,\QF),(\Ctwo,\QFtwo)) \to (\Ctwo,\QFtwo)
\end{equation}
from Construction~\refone{construction:evaluation}. By the universal property of internal functor categories this transposes to a hermitian functor
\[
(\C,\QF) \to \Funx(\Funx((\C,\QF),(\Ctwo,\QFtwo)),(\Ctwo,\QFtwo)),
\]
and consequently induces for \(\I \in \Cat\) a hermitian functor
\[
(\C,\QF)^{\I} \to \Funx(\Funx((\C,\QF),(\Ctwo,\QFtwo)),(\Ctwo,\QFtwo))^{\I} \simeq \Funx(\Funx((\C,\QF),(\C,\QFtwo))_{\I},(\Ctwo,\QFtwo)),
\]
where we have used the equivalence of Remark~\refone{remark:comparing-1}. The resulting functor then transposes twice to give a hermitian functor
\[
\Funx((\C,\QF),(\C,\QFtwo))_{\I} \to \Funx((\C,\QF)^{\I},(\Ctwo,\QFtwo)).
\]
On the other hand, using the equivalence of Remark~\refone{remark:comparing-2} the evaluation functor~\eqrefone{equation:evaluation-internal} induces a hermitian functor
\begin{multline*} (\C,\QF) \otimes \Funx((\C,\QF),(\Ctwo,\QFtwo))_{\I} \simeq\\
                \simeq(\C,\QF) \otimes \Funx((\C,\QF),(\Ctwo,\QFtwo)) \otimes (\Spaf,\QF^{\uni})_{\I} \to (\Ctwo,\QFtwo) \otimes (\Spaf,\QF^{\uni})_{\I}
                \simeq (\Ctwo,\QFtwo)_{\I}\end{multline*}
which transposes to give a hermitian functor
\[
\Funx((\C,\QF),(\Ctwo,\QFtwo))_{\I} \to \Funx((\C,\QF),(\Ctwo,\QFtwo)_{\I}).
\]
Combining the above constructions we hence obtain a pair hermitian functors
\begin{equation}
\label{equation:exceptional-adj}%
\Funx((\C,\QF),(\Ctwo,\QFtwo)_{\I}) \longleftarrow \Funx((\C,\QF),(\Ctwo,\QFtwo))_{\I} \longrightarrow \Funx((\C,\QF)^{\I},(\Ctwo,\QFtwo)).
\end{equation}
natural in \((\C,\QF),(\Ctwo,\QFtwo)\) and \(\I\) (indeed, all the operations used above have already been proven natural in \S\refone{subsection:internal}, \S\refone{subsection:cotensoring} and \S\refone{subsection:tensoring} through the various universal properties they encode).

\begin{proposition}
\label{proposition:exceptional-adj}%
If \(\I\) is strongly finite then the hermitian functors in~\eqrefone{equation:exceptional-adj} are equivalences of hermitian \(\infty\)-categories. Passing to hermitian objects (see~\eqrefone{equation:hermitian-objects-internal}) they then determine a natural equivalence
\[
\Map_{\Cath}((\C,\QF),(\Ctwo,\QFtwo)_{\I}) \simeq \Map_{\Cath}((\C,\QF)^{\I},(\Ctwo,\QFtwo))
\]
exhibiting \((-)^{\I}\) as left adjoint to \((-)_{\I}\).
\end{proposition}
\begin{proof}
To begin, we note that on the level of underlying stable \(\infty\)-categories both functors in~\eqrefone{equation:exceptional-adj} are equivalences by Lemma~\refone{lemma:finite}. Indeed, replacing \(\I\) with \(\I\op\) these identify with the equivalences of stable \(\infty\)-categories
\begin{equation}
\label{equation:adj-underlying}%
\Funx(\C,\Ctwo^{\I\op}) \xleftarrow{\simeq} \Funx(\C,\Ctwo)^{\I\op} \xrightarrow{\simeq} \Funx(\C_{\I\op},\Ctwo)
\end{equation}
underlying those of Remark~\refone{remark:comparing-1}.
Explicitly, the equivalence on the left hand side of~\eqrefone{equation:adj-underlying} associates to a diagram \(\vphi\colon \I\op \to \Funx(\C,\Ctwo)\) the exact functor \(g_{\vphi}\colon \C \to \Ctwo^{\I\op} = \Ctwo_{\I}\)
given by \([g_{\vphi}(\x)](i) = \vphi_i(\x)\).
Unwinding the definitions, the hermitian structure of the left hand side functor in~\eqrefone{equation:exceptional-adj} is given by the map
\[
\colim_{i \in \I}\nat(\QF,\vphi^*_i\QFtwo) \to \nat\big(\QF,\colim_{i \in \I}\vphi^*_i\QFtwo\big) = \nat\big(\QF,g_{\vphi}^*\QF_{\I}\big),
\]
which is indeed an equivalence since \(\I\) is finite and \(\nat(\QF,-)\) is an exact functor.

Similarly, the equivalence on the right hand side of~\eqrefone{equation:adj-underlying} associates to a diagram \(\vphi\colon \I\op \to \Funx(\C,\Ctwo)\) an exact functor \(h_{\vphi}\colon \C_{\I\op}=\C^{\I} \to \Ctwo\)
such that \(\vphi\) can be recovered from \(h_{\vphi}\) as \(\vphi_i(\x) = h_{\vphi}(R_{i,\x})\).
Let us denote by \(\ev_i\colon \C^{\I} \to \C\) the evaluation at \(i \in \I\) functor and by \(\ran_i\colon \C \to \C^{\I}\) its right adjoint, given by right Kan extension. In particular, we have \(\ran_i(\x) = R_{i,\x}\) by definition. Unwinding the definitions, the hermitian structure of the right hand side functor in~\eqrefone{equation:exceptional-adj} is then given by the map
\[
\colim_{i \in \I}\nat(\QF,\vphi^*_i\QFtwo)  = \colim_{i\in \I}\nat(\QF,\ran_i^*h_{\vphi}^*\QFtwo) = \colim_{i \in \I}\nat(\ev_i^*\QF,h_{\vphi}^*\QFtwo) \to \nat(\lim_i\ev_i^*\QF,h_{\vphi}^*\QFtwo) = \nat(\QF^{\I},h_{\vphi}^*\QFtwo),
\]
which is indeed an equivalence since \(\I\) is finite and \(\nat(-,h^*_{\vphi}\QFtwo)\) is an exact functor.
\end{proof}

\begin{corollary}
\label{corollary:limit-preservation}%
Let \(\I\) be a strongly finite \(\infty\)-category. Then the functor \((\C,\QF) \mapsto (\C,\QF)_{\I}\) from \(\Cath\) to itself preserves all limits and the functor \((\C,\QF) \mapsto (\C,\QF)^{\I}\) preserves all colimits.
\end{corollary}

\begin{corollary}
\label{corollary:internally-corepresented}%
For a fixed strongly finite \(\infty\)-category \(\I\), the functor \((\C,\QF) \mapsto (\C,\QF)_{\I}\) is internally corepresented by \((\Spaf,\QF^{\uni})^{\I}\). More precisely, there is an equivalence of hermitian \(\infty\)-categories
\[
(\C,\QF)_{\I} \simeq \Funx((\Spaf,\QF^{\uni})^{\I},(\C,\QF))
\]
natural in \((\C,\QF)\) and \(\I\).
\end{corollary}

\begin{corollary}
\label{corollary:I-preserve-poincare}%
For a strongly finite \(\infty\)-category \(\I\) the following conditions are equivalent:
\begin{enumerate}
\item
the operation \((\C,\QF) \mapsto (\C,\QF)_{\I}\) preserves Poincaré \(\infty\)-categories and Poincaré functors.
\item
the hermitian \(\infty\)-category \((\Spaf,\QF^{\uni})_{\I}\) is Poincaré;
\item
the operation \((\C,\QF) \mapsto (\C,\QF)^{\I}\) preserves Poincaré \(\infty\)-categories and Poincaré functors.
\item
the hermitian \(\infty\)-category \((\Spaf,\QF^{\uni})^{\I}\) is Poincaré;
\end{enumerate}
\end{corollary}

In a similar spirit, we may deduce that the criterion for duality  preservation of Proposition~\refone{proposition:induced-maps-duality-pres} holds for tensors as well in the strongly finite case:

\begin{corollary}
\label{corollary:tensor-duality-preserving}%
Let \(\alp\colon \I \to \J\) be a map of strongly finite \(\infty\)-categories and \((\C,\QF)\) a non-degenerate hermitian \(\infty\)-category. If the induced map \(\I_{i/} \to \J_{\alp(i)/}\) is cofinal for every \(i \in \I\) then the induced hermitian functor \((\alp_*,\eta_{\alp})\colon(\C,\QF)_\I \to (\C,\QF)_\J\)
is duality preserving. In particular, if \((\C,\QF)_{\I}\) and \((\C,\QF)_{\J}\) are Poincaré then \((\alp_*,\eta_{\alp})\) is a Poincaré functor.
\end{corollary}
\begin{proof}
Identify \((\C,\QF)_{\I}\) with \(\Funx((\Spaf,\QF^{\uni})_{\I},(\C,\QF))\) as a functor of \(\I\) using Corollary~\refone{corollary:internally-corepresented}
and apply Proposition~\refone{proposition:induced-maps-duality-pres}.
\end{proof}

We now describe some additional functoriality exhibited by the tensor and cotensor constructions in the strongly finite case.

\begin{construction}[Exceptional functoriality]
\label{construction:exceptional}%
Let \(\beta\colon \J \to \I\) be a functor between strongly finite \(\infty\)-categories. Then by Lemma~\refone{lemma:finite} the full subcategory \(\C_{\I} \subseteq \Fun(\I\op,\Pro(\C))\) is sent into \(\C_{\J} \subseteq \Fun(\J\op,\Pro(\C))\) by restriction along any functor \(\beta\colon \J \to \I\). On the cotensor side, the comma categories of \(\beta\) are all finite by Remark~\refone{remark:strongly-finite} and hence the operation of right Kan extension \(\beta_*\colon \C^{\J} \to \C^{\I}\) exists for any stable \(\C\).
Now suppose that \(\beta\) is \emph{cofinal}. Then we can refine \(\beta^*\colon \C_{\I} \to \C_{\J}\) and \(\beta_*\colon \C^{\J} \to \C^{\I}\) to hermitian functors as follows. In the tensor case we simply note that the cofinality of \(\beta\) yields a natural equivalence
\[
\QF_{\J}(\beta^*\vphi) = \colim_{j \in \J}\QF(\vphi(\beta(j)) \simeq \colim_{i \in \I}\QF(\vphi(i)) = \QF_{\I}(\vphi)
\]
and so we obtain a hermitian functor \((\beta^*,\vartheta^{\beta})\colon (\C,\QF)_{\I} \to (\C,\QF)_{\J}\) in which \(\vartheta^{\beta}\colon \QF_{\I} \Rightarrow \QF_{\J} \circ \beta^*\) is an equivalence. For cotensors we consider the counit \(\beta^*\beta_*\psi \Rightarrow \psi\) and use the coinitiality of \(\beta\op\) to obtain a map
\[
\QF^{\J}(\psi) \Rightarrow \QF^{\J}(\beta^*\beta_*\psi) = \lim_{j \in \J\op} \QF(\beta_*\psi(\beta(j)) \simeq \lim_{i \in \I\op}\QF(\beta_*\psi(i)) = \QF^{\I}(\beta_*\psi)
\]
giving a hermitian refinement \((\beta_*,\vartheta_{\beta})\colon (\C,\QF)^{\J} \to (\C,\QF)^{\I}\).
\end{construction}

\begin{example}
\label{example:exceptional-is-direct}%
A common source of cofinal maps are maps \(\beta\colon \J \to \I\) which admit a left adjoint \(\alp\colon \I \to \J\). Unwinding the definitions and using Remarks~\refone{remark:restriction} and~\refone{remark:right-kan-extension}, we see that in this case the exceptional hermitian functors \((\beta^*,\vartheta^{\beta})\) and \((\beta_*,\vartheta_{\beta})\) of Construction~\refone{construction:exceptional} coincide with the direct hermitian functors \((\alp_*,\eta_{\alp})\) and \((\alp^*,\eta^{\alp})\), respectively.
\end{example}

\begin{remark}
\label{remark:switch-sides}%
Let \(\I,\J\) be two strongly finite \(\infty\)-categories and \((\C,\QF)\) a Poincaré \(\infty\)-category. The identification
\begin{equation}
\label{equation:first-identification}%
(\C,\QF)^{\I} \simeq \Funx((\Spaf,\QF^{\uni})_{\I},(\C,\QF))
\end{equation}
of Remark~\refone{remark:comparing-1}, being natural in \(\I\), identifies, for every map \(\alp\colon \I \to \J\), the associated hermitian functor \((\alp^*,\eta^{\alp})\colon (\C,\QF)^{\J} \to (\C,\QF)^{\I}\) with the one induced from \((\alp_*,\eta_{\alp})\colon (\Spaf,\QF^{\uni})_{\I} \to (\Spaf,\QF^{\uni})_{\J}\)
upon taking upon taking internal functor categories to \((\C,\QF)\). Similarly, the identification
\begin{equation}
\label{equation:second-identification}%
(\C,\QF)_{\I} \simeq \Funx((\Spaf,\QF^{\uni})^{\I},(\C,\QF))
\end{equation}
of Corollary~\refone{corollary:internally-corepresented} identifies the hermitian functor \((\alp_*,\eta_{\alp})\colon (\C,\QF)_{\I} \to (\C,\QF)_{\J}\) with the one induced from \((\alp^*,\eta^{\alp})\colon (\Spaf,\QF^{\uni})^{\I} \to (\Spaf,\QF^{\uni})^{\J}\).
Unravelling all definitions, and observing the similarity between the formulas for the direct and exceptional functorialities, we see that the equivalence~\eqrefone{equation:first-identification} also identifies for every \(\beta\colon \J \to \I\) the exceptional hermitian functor \((\beta_*,\vartheta_{\beta})\colon (\C,\QF)^{\J} \to (\C,\QF)^{\I}\) with the one induced from the exceptional functor \((\beta^*,\vartheta^{\beta})\colon (\Spaf,\QF^{\uni})_{\I} \to (\Spaf,\QF^{\uni})_{\J}\)
upon taking internal functor categories to \((\C,\QF)\), and similarly the equivalence~\eqrefone{equation:second-identification} also identifies the exceptional hermitian functor \((\beta^*,\vartheta^{\beta})\colon (\C,\QF)_{\I} \to (\C,\QF)_{\J}\) with the one induced from \((\beta_*,\vartheta_{\beta})\colon (\Spaf,\QF^{\uni})^{\J} \to (\Spaf,\QF^{\uni})^{\I}\).
\end{remark}

\begin{remark}
\label{remark:compatible-with-composition}%
Like the ordinary functoriality of the tensor and cotensor constructions, the exceptional functorialities are compatible with composition. More precisely, if \(\cK \xrightarrow{\beta'} \J \xrightarrow{\beta} \I\) are a pair of composable cofinal maps then the hermitian functor \(((\beta\circ \beta')^*,\vartheta^{\beta\circ \beta'})\colon (\C,\QF)_{\I} \to (\C,\QF)_{\cK}\) is naturally equivalent to the composite of \((\beta^*,\vartheta^{\beta})\) and \(((\beta')^*,\vartheta^{\beta'})\) and the hermitian functor \(((\beta\circ \beta')_*,\vartheta_{\beta\circ \beta'})\colon (\C,\QF)^{\cK} \to (\C,\QF)^{\I}\) is naturally equivalent to the composite of \((\beta_*,\vartheta_{\beta})\) and \((\beta'_*,\vartheta_{\beta'})\). Indeed, by Remark~\refone{remark:switch-sides} it will suffice to check this for the tensor construction, where it amounts to the fact the Beck-Chevalley transformation relating restriction and colimits is compatible with composition of restriction maps.
\end{remark}

\begin{remark}
\label{remark:one-sided-inverse}%
Comparing the explicit formulas of the direct and exceptional functorialities we see that if \(\beta\colon \J \to \I\) is a cofinal map then the resulting hermitian structure on the exceptional-direct composite \((\beta^*,\eta^{\beta})\circ (\beta_*,\vartheta_{\beta})\colon (\C,\QF)^{\J} \to (\C,\QF)^{\J}\) is given by the map
\[
\QF^{\J}(\psi) \Rightarrow \QF^{\J}(\beta^*\beta_*\psi)
\]
induced by the counit \(\beta^*\beta_*\psi \Rightarrow \psi\).
Similarly, the resulting hermitian structure on the direct-exceptional composite
\((\beta^*,\vartheta^{\beta})\circ (\beta_*,\eta_{\beta})\colon (\C,\QF)_{\J} \to (\C,\QF)_{\J}\) is
induced by the counit of \(\C_{\I} \adj \C_{\J}\) in a similar manner. In particular, if \(\beta\) is fully-faithful then these counits are equivalences, in which case we get that the exceptional functoriality of \(\beta\) gives a one sided inverse to its direct functoriality. More generally, it can be checked that for any cofinal \(\beta\) the exceptional functoriality \((\beta_*,\vartheta_{\beta})\) is right adjoint to the direct functoriality \((\beta^*,\eta^{\beta})\), and the exceptional functoriality \((\beta^*,\vartheta^{\beta})\) is left adjoint to the direct functoriality \((\beta_*,\eta_{\beta})\), when these terms are understood with respect to the \((\infty,2)\)-categorical structure of \(\Cath\) determined by the internal functor categories of \S\refone{subsection:internal}.
\end{remark}

\subsection{Finite complexes and Verdier duality}
\label{subsection:finite-complexes}%

In this final section we explore the particular case when \(\I\) is the poset of simplices of a finite simplicial complex. In this context the tensor construction plays a central role in~\cite{Lurie-L-theory}, where its duality functor is identified as a form of \defi{Verdier duality}. Filling in the details in some of the arguments of loc.\ cit., we will show that tensoring and cotensoring by the poset of faces of a finite simplicial complex preserves Poincaré \(\infty\)-categories, and that the hermitian functors associated to maps of simplicial complexes (direct functoriality) and refinements of triangulations (exceptional functoriality) are Poincaré.

Recall that a finite simplicial complex \(K\) consists of a finite set of vertices \(K_0\) and a collection \(\I_{K}\) of non-empty subsets \(S \subseteq K_0\), called faces, which contain all singletons and are downwards closed in the sense that if \(S\) is a face of \(K\) and \(S'\subseteq S\) then \(S'\) is a face of \(K\) as well. The \defi{dimension} of a face \(S\) is by definition \(\dim(S) := |S|-1\). We may realize a simplicial complex \(K\) geometrically as the subspace \(|K|\) of the full simplex on \(K_0\) obtained as the union of the given faces. For a finite simplicial complex \(K\) we will consider \(\I_K\) as a poset, and consequently a category, by \emph{inverse} inclusion, that is, there is a unique morphism \(S \to S'\) if \(S'\subseteq S\). A map of simplicial complexes \(K \to L\) is by definition a map of sets \(K_0 \to L_0\) which sends every face of \(K\) to a face of \(L\). In particular, such a map determines a map of posets \(\alp\colon\I_{K} \to \I_{L}\). If \(K\) is a simplicial complex then a \defi{refinement} of \(K\) consists of a simplicial complex \(L\) together with a homeomorphism \(|L| \xrightarrow{\cong} |K|\) which carries the realization of every face of \(L\) into the realization of a face of \(K\), and such that the realization of every face of \(K\) in \(|K|\) is a union of faces of \(L\). We note that such a homeomorphism determines in particular a map of posets \(\beta\colon \I_L \to \I_K\) which sends every face of \(L\) to the face of \(K\) containing it under the homeomorphism \(|L| \cong |K|\).

\begin{proposition}[{\cite[Lecture 19, Proposition 3]{Lurie-L-theory}}]
\label{proposition:posets-of-faces}%
Let \((\C,\QF)\) be a Poincaré \(\infty\)-category and \(K\) a finite simplicial complex with poset of faces \(\I_K\).
Then the hermitian \(\infty\)-categories \((\C,\QF)_{\I_K}\) and \((\C,\QF)^{\I_K}\) are Poincaré.
\end{proposition}
\begin{proof}
By Corollary~\refone{corollary:I-preserve-poincare} it will suffice to prove the claim for \((\C,\QF)_{\I_K}\).
We need to show that for every \(\vphi\in \C_{\I_K}\) the map
\begin{equation}
\label{equation:double-dual-cotensor}%
    \vphi\to \Dual_{\I_K}\Dual_{\I_K}\vphi
\end{equation}
is an equivalence. To do so it suffices to show it for a system of objects that generate under colimits. We choose the set \(\vphi=R_{\x,S}\) for \(S\in \I_K\) and \(x\in\C\). Now the face \(S\) corresponds to an injective map of simplicial complexes \(\Del^n \to K\)
which in turn determines an inclusion of posets \(\alp\colon \I_{\Del^n} \subseteq \I_K\) such that \(\alp([n]) =S\) and
\(R_{\x,S}\) is the right Kan extension along \(\alp\op\) of the constant diagram \(\vphi_{\x}\colon \I_{\Del^n}\op \to \C\) with value \(\x\).
Now the map \(\alp\) satisfies the hypothesis of Corollary~\refone{corollary:tensor-duality-preserving} (in fact, the map that needs to be cofinal is an isomorphism of posets, see also Proposition~\refone{proposition:poset-of-faces-functoriality} below), and so the hermitian functor \((\alp_*,\eta_{\alp})\colon \C_{\I_{\Del^n}} \to \C_{\I_K}\) is duality preserving. We can hence reduce to the case of \(K=\Del^n\) and \(\vphi=\vphi_x\). Using Proposition~\refone{proposition:tensor-strongly-finite} we now have
\begin{equation}
\label{equation:dual-on-simplex}%
\Dual_{\I_{\Del^n}}\vphi_{\x}(S) =
\colim_{\emptyset \neq T \subseteq [n]}
\left\{\begin{matrix} \Dual\x & S \subseteq T \\ 0 & \text{otherwise} \\ \end{matrix}\right..
\end{equation}
To calculate this colimit let us denote by \(\iota\colon \I^S_{\Del^n} \subseteq \I_{\Del^n}\) the subposet spanned by those \(T \subseteq [n]\) such that \(S \not{\subseteq} T\). Then the functor whose colimit is calculated in~\eqrefone{equation:dual-on-simplex} can be identified with the cofibre of the map \(\iota^*\iota_!\vphi_x \to \vphi_x\), and so we get that
\[
[\Dual_{\I_{\Del^n}}\vphi_{\x}](S) = \cof\big[\colim_{\I^S_{\Del^n}}\Dual\x\to \colim_{\I_{\Del^n}}\Dual\x\big] .
\]
Now since \(\I^S_{\Del^n}\) is closed under subfaces in \(\I_{\Del^n}\) it corresponds to some subcomplex \(\Del^n\),
which we readily identify as the join \(\partial\Del^S \ast \Del^{S'}\), where \(S'= [n]- S\) is the complement of \(S'\), and we have used the notation \(\Del^S\) and \(\Del^{S'}\) to denote the corresponding faces, considered as subcomplexes. The poset \(\I^S_{\Del^n}\) is hence weakly contractible if \(S' \neq \emptyset\), that is, if \(S \neq [n]\), and is weakly equivalent to \(\partial \Del^n\) of \(S=[n]\). We thus conclude that
\[
[\Dual_{\I_{\Del^n}}\vphi_{\x}](S) =
\left\{\begin{matrix} \Sig^n\Dual\x & S = [n] \\ 0 & \text{otherwise} \\ \end{matrix}\right..
\]
By the same argument we then have
\[
[\Dual_{\I_{\Del^n}}\Dual_{\I_{\Del^n}}\vphi_{\x}](S) = \colim_{\emptyset \neq T \subseteq [n]}
\left\{\begin{matrix} \Dual\Sig^n\Dual\x & T=[n] \\ 0 & \text{otherwise} \\ \end{matrix}\right.,
\]
so that \(\Dual_{\I_{\Del^n}}\Dual_{\I_{\Del^n}}\vphi_{\x}\) is constant with value \(\Sig^n\Dual\Sig^n\Dual\x \simeq \x\), and is in particular equivalent to \(\vphi_{\x}\). We need however to make sure that specifically the evaluation map is an equivalence. Since both \(\Dual_{\I_{\Del^n}}\Dual_{\I_{\Del^n}}\vphi_{\x}\) and \(\vphi_{\x}\) are constant it will suffice to show that the component of the evaluation map at \(S=[n]\) is an equivalence. Unwinding the definitions, this is the composed map
\[
\begin{split}
\x \xrightarrow{\simeq} \Dual\Dual\x & \to \colim_{T \in \I_{\Del^n}}\lim_{T' \in \I_{\Del^n}\op}
\begin{cases}
\Dual\Dual\x & \text{if \(T=T'=[n]\)} \\
0 & \text{otherwise.}
\end{cases} \\
 & = \lim_{T' \in \I_{\Del^n}\op} \colim_{T \in \I_{\Del^n}}
\begin{cases}
\Dual\Dual\x & \text{if \(T=T'=[n]\)} \\
0 & \text{otherwise.}
\end{cases}
\end{split}
\]
whose component at level \(T'\)
is \(0\) for \(T'\neq  [n]\) and is the inclusion of the \(T=[n]\) component in the colimit otherwise. The invertibility of this map
then reduces to the fact that in a stable \(\infty\)-category an \(n\)-cube is cartesian if and only if it is cocartesian.
\end{proof}

The perfectness of hermitian structures asserted in Proposition~\refone{proposition:posets-of-faces} is accompanied by the following duality preservation statement:

\begin{proposition}{\cite[Lecture 19]{Lurie-L-theory}}
\label{proposition:poset-of-faces-functoriality}%
Let \((\C,\QF)\) be a Poincaré \(\infty\)-category.
\begin{enumerate}
\item
If \(K \to L\) is a map of finite simplicial complexes with \(\alp\colon \I_K \to \I_L\) the induced map of posets of faces then the induced hermitian functors \((\alp_*,\eta_{\alp})\colon (\C,\QF)_{\I_K} \to (\C,\QF)_{\I_L}\) and \((\alp^*,\eta^{\alp})\colon (\C,\QF)^{\I_L} \to (\C,\QF)^{\I_K}\) are Poincaré.
\item
If \(L\) is a refinement of a simplicial complex \(K\) and \(\beta\colon \I_{L} \to \I_{K}\) is the associated map of posets of faces then \(\beta\) is cofinal and the exceptional hermitian functors
\((\beta^*,\vartheta^{\beta})\colon (\C,\QF)_{\I_K} \to (\C,\QF)_{\I_L}\)
and
\((\beta_*,\vartheta_{\beta})\colon (\C,\QF)^{\I_L} \to (\C,\QF)^{\I_K}\) are Poincaré.
\end{enumerate}
\end{proposition}
\begin{proof}
By Remark~\refone{remark:switch-sides} it will suffice to prove the tensor case. For the first statement we observe that the functor \(\alp\colon\I_{K} \to \I_{L}\) satisfies the criterion of Corollary~\refone{corollary:tensor-duality-preserving}, since for every face \(S \in \I_{K}\) the functor \((\I_{K})_{S/} \to (\I_{L})_{\alp(S)/}\) admits a left adjoint sending \(T \subseteq \alp(S)\) to its inverse image in \(S\).

To prove the second statement, we begin by showing that \(\bet\) is cofinal. Indeed for any simplex \(S\in {\I_K}\) the poset \({\I_L}\times_{\I_K} (\I_K)_{S/}\) has geometric realization homeomorphic to a simplex and so it is weakly contractible.
To prove that \((\beta^\ast,\vartheta^{\beta})\) is Poincaré, it will suffice to show that for every generator \(R_{\x,S} \in \C_{\I_K}\) the associated map
\[
\beta^*\Dual_{\I_K}R_{\x,S} \to \Dual_{\I_L}\beta^*R_{\x,S}
\]
is an equivalence.
Now the face \(S\) corresponds to an injective map of simplicial complexes \(\Del^S \to K\) (where \(\Del^S\) denotes the full simplex with vertex set \(S)\), which in turn determines an inclusion of posets \(\alp\colon \I_{\Del^S} \subseteq \I_K\) such that
\(R_{\x,S}\) is the right Kan extension along \(\alp\) of the constant diagram \(\vphi_{\x}\colon \I_{\Del^S} \to \C\) with value \(\x\). The inverse image of \(\I_{\Del^S}\) in \(\I_L\) then determines a subcomplex \(L_S \subseteq L\) which is a refinement of \(\Del^S\).

Let us denote by \(\beta_S\colon \I_{L_S} \to \I_{\Del^S}\) the induced refinement map
and by \(\wtl{\alp}\colon \I_{L_S} \hrar \I_L\) the inclusion. Since \(\I_{\Del^S}\) and \(\J\) are downward closed the pointwise formula for right Kan extensions implies that the square
\[
\begin{tikzcd}
(\C,\QF)_{\I_{\Del^S}} \ar[r,"{(\beta_S^*,\vartheta^{\beta}_S)}"] \ar[d,"{(\alp_*,\eta_{\alp})}"'] & (\C,\QF)_{\I_{L_S}} \ar[d,"{(\wtl{\alp}_*,\wtl{\eta}_{\alp})}"] \\
(\C,\QF)_{\I_K} \ar[r,"{(\beta^*,\vartheta^{\beta})}"] & (\C,\QF)_{\I_L}
\end{tikzcd}
\]
commutes. Since the vertical hermitian functors are Poincaré by the first part we may reduce to the case where \(K\) is the \(n\)-simplex \(\Del^n\), \(L\) is some refinement of \(\Del^n\), and \(\vphi = \vphi_{\x}\).
Now for \(T \in \I_L\) let us denote by \(\I^T_L \subseteq \I_L\) the subposet spanned by those faces which do not contain \(T\), and by \(\J^T_L \subseteq \I^T_L\) the subposet spanned by those faces whose image in \(\I_K\) does not contain \(\beta(T)\). We note that both of these subposets are downward closed and correspond to subcomplexes of \(L\). In particular, \(\I^T_L\) corresponds to the subcomplex \(L^T \subseteq L\) obtained by removing all faces which contain \(T\), and \(\J^T_L\) corresponds to the subcomplex \(L^T_0 \subseteq L^T\) obtained by removing all faces whose image in \(K\) contains \(T\). We also note that \(L^T_0\) is a refinement of the subcomplex \(\partial\Del^{\beta(T)} \ast \Del^{[n]-\beta(T)} \subseteq \Del^n\) obtained from \(\Del^n\) by removing all the faces which contain \(\beta(T)\). Calculating as in the proof of Proposition~\refone{proposition:posets-of-faces} and using that refinement maps are cofinal as established above we may identify the cofibre of the map
\[
[\beta^*\Dual_{\I_{\Del^n}}\vphi_{\x}](T) \to [\Dual_{\I_L}\beta^*\vphi_{\x}](T)
\]
for \(T \in \I_L\) with the total cofibre of the square
\[
\begin{tikzcd}
\colim_{\J^{T}_{L}}\Dual\x \ar[r] \ar[d] & \colim_{\I_{L}}\Dual\x \ar[d,equal] \\
\colim_{\I^T_{L}}\Dual\x \ar[r] & \colim_{\I_{L}}\Dual\x
\end{tikzcd}.
\]
To finish the proof it will hence suffice to show that \(|L^0_T| \to |L_T|\) is a weak homotopy equivalence. Let \(p \in |\Del^n|\) be a point in the interior of the face \(T\) (and hence also in the interior of the face \(\beta(T)\)). Then we have a sequence of inclusions
\( |L^0_T| \subseteq |L_T| \subseteq U \)
where \(U \subseteq |\Del^n|\) is the complement of \(p\). The desired result now follows from the fact that both \(|L^0_T| \subseteq U\) and \(|L_T| \subseteq U\) are deformation retracts; this is a general property of simplicial complexes: if one removes a point from the realization of a simplicial complex then the result deformation retracts to the subcomplex spanned by all the simplices which do not contain that point.
\end{proof}

%% file: Metabolic.tex
From a conceptual view point, it is arguably tempting to regard hermitian structures on stable \(\infty\)-categories as categorified versions of hermitian forms on modules. Similarly,
Poincaré \(\infty\)-categories correspond to modules equipped with a unimodular hermitian form. Inspired by this informal perspective, in this section we will identify a surprisingly comprehensive variety of such categorified counterparts, including the categorified analogues of bilinear forms, perfect bilinear forms, hyperbolic objects, metabolic objects and the algebraic Thom construction. Beyond its conceptually pleasing effect, it turns out that many of the constructions encountered via this perspective give explicit left and right adjoints to various natural functors, a feature which we will repeatedly exploit in subsequent instalments of this project. In particular, the main practical outcomes of our categorified stroll will include the following:
\begin{enumerate}
\item
After exploring the categorified counterparts of bilinear forms in \S\refone{subsection:bilinear-and-pairings}, we will deduce that the association \(\C \mapsto \Hyp(\C)\) described in \S\refone{subsection:hyp-and-sym-poincare-objects} is both left and right adjoint to the forgetful functor \(\Catp \to \Catx\), with unit and counit given on the side of \(\Catp\) by the Poincaré functors \(\hyp\colon \Hyp(\C) \to (\C,\QF)\) and \(\fgt\colon (\C,\QF) \to \Hyp(\C)\) described in \S\refone{subsection:hyp-and-sym-poincare-objects}. In addition, we will show that \(\Hyp\) is \(\Ct\)-equivariant with respect to the \(\mop\) action on \(\Catx\) and the trivial action on \(\Catp\), and that \(\hyp\) and \(\fgt\) are equivariant natural transformations. All this information is best organised in the setting of \(\Ct\)-categories and Mackey functors, which we will explore in~\S\refone{subsection:mackey-functors}. This will also be the basis for our organization in \papertwo of algebraic \(\K\)-theory, Grothendieck-Witt theory and \(\L\)-theory into a single functor taking values in genuine \(\Ct\)-spectra, which we call the \emph{real \(\K\)-theory spectrum}.
\item
While forming the categorified analogue of metabolic objects, Lagrangians, and the algebraic Thom construction in \S\refone{subsection:thom}, we will deduce an explicit formula for the left and right adjoints to the inclusion \(\Catp \to \Cath\). This will be exploited in \papertwo when setting up the framework of \emph{algebraic surgery} in the context of Poincaré \(\infty\)-categories, and in analysing the effect of additive and bordism invariant functors applied to the \(\Q\)-construction. We will also use it \paperfour for constructing the localising analogue of the Grothendieck-Witt spectrum and for proving that \(\Catp\) is compactly generated.
\end{enumerate}

The present section is organised as follows. We begin in \S\refone{subsection:bifibrations} with some preliminary material on \emph{bifibrations}, a notion used for encoding space valued bifunctors which are covariant in one entry and contravariant in the other. By translating results from~\cite[\S 5.2.1]{HA} to the context of bifibrations we deduce in particular that the \(\infty\)-category of perfect symmetric bifibrations is equivalent to \(\Cat^{\hC}\).
In \S\refone{subsection:bilinear-and-pairings} we specialise to the setting of stable \(\infty\)-categories and replace space valued bifunctors by spectrum valued ones. This leads to the notion of bilinear \(\infty\)-categories as analogous of pairs of modules equipped with a bilinear form. We also consider the variants of requiring the form to be perfect and/or symmetric, and identify, using \S\refone{subsection:bifibrations}, the notion of a perfect symmetric \(\infty\)-category with that of a stable \(\infty\)-category equipped with a perfect duality.
All these notions accept natural forgetful functors from either \(\Catp\) (in the perfect case) or \(\Cath\). Studying all of them on equal footing allows one to efficiently identify left and right adjoints to these forgetful functors, which constitutes the main content of \S\refone{subsection:bilinear-and-pairings}.
In particular, we will see that the association \(\C \mapsto \Hyp(\C)\) described in \S\refone{subsection:hyp-and-sym-poincare-objects} gives a two-sided adjoint to the forgetful functor \(\Catp \to \Catx\).
In \S\refone{subsection:thom} we will discuss the categorified analogous of metabolic objects and Lagrangians, and will show that the categorified analogue of the Thom construction enables one to produce both a left and a right adjoint to the inclusion \(\Catp \to \Cath\). We will also prove a generalization of the non-categorified Thom construction, thus providing in particular a proof for Proposition~\refone{proposition:thom} which was stated in \S\refone{subsection:algebraic-thom}.
Finally, in \S\refone{subsection:mackey-functors} we will discuss \(\Ct\)-categories and Mackey functors, and show how to view various players in the present paper in that context. In particular, we will show that the relations between quadratic and bilinear functors, between hermitian and bilinear \(\infty\)-categories, and between Poincaré and stable \(\infty\)-categories, can all be understood on the same footing, giving another take on the categorified perspective. The results of \S\refone{subsection:mackey-functors} will be primarily used in \papertwo in order to define the real \(\K\)-theory spectrum, but the \(\Ct\)-equivariant properties of the hyperbolic construction resulting from it will be useful for a variety of other purposes as well.

\subsection{Preliminaries: pairings and bifibrations}
\label{subsection:bifibrations}%

Let \(\A,\cB\) be two \(\infty\)-categories. By a \defi{correspondence} on the pair \((\A,\cB)\) we will simply mean a functor \(b\colon\A\op \times \cB \to \Sps\). In particular, we will think of a correspondence as a space valued functor on pairs \((\x,\y)\) with \(\x \in \A\) and \(\y \in \cB\), which is contravariant in \(\x\) and covariant in \(y\). The prototypical example to have in mind is taking \(\A = \cB = \C\) for some \(\infty\)-category \(\C\), with \(b(\x,\y) = \Map_{\C}(\x,\y)\). Given a correspondence \(b\colon \A\op \times \cB \to \Sps\) and \(\y \in \cB\), evaluation at \(\y\) yields a presheaf of spaces \(b(-,\y)\) on \(\A\), which we can unstraighten to obtain a right fibration
\begin{equation}
\label{equation:first-unstraightening}%
\int^{\x \in \A} b(\x,\y) \to \A.
\end{equation}
Since the presheaf \(b(-,\y) \in \Fun(\A\op,\Sps)\) depends functorially on \(\y\), so does the domain of~\eqrefone{equation:first-unstraightening}. In fact, we can identify the arrow~\eqrefone{equation:first-unstraightening} with the map \(\int^{\x \in \A} b(\x,\y) \to \int^{\x \in \A}\ast\) associated to the terminal map of correspondences \(b \to \ast\), and so the entire arrow~\eqrefone{equation:first-unstraightening} depends functorially on \(\y\). Equivalently, we may view it as a natural transformation between two \(\Cat\)-valued functors on \(\cB\), the second of which is constant with value \(\A\). We then define
\begin{equation}
\label{equation:second-unstraightening}%
\Pairings(\A,\cB,b) := \int_{\y \in \cB} \int^{\x \in \A}b(\x,\y) \to \A \times \cB ,
\end{equation}
to be the \(\infty\)-category over \(\A \times \cB\) obtained by unstraightening~\eqrefone{equation:first-unstraightening} over \(\cB\).
The \(\infty\)-category \(\Pairings(\A,\cB,b)\) can informally be described as having objects triples \((\x,\y,\beta)\) where \(\x \in \A\) and \(\y \in \cB\) are objects and \(\bet \in b(\x,\y)\) is a \(b\)-valued pairing on \(\x\) and \(\y\).
A map from \((\x,\y,\bet)\) to \((\x',\y',\beta')\) is then given by maps \(f\colon \x \to \x',g\colon \y \to \y'\) and a homotopy \(\eta\colon g_*\bet \sim f^*\bet' \in b(\x,\y')\),
where \(f^*\) and \(g_*\) encode the contravariant and covariant dependence of \(b\) on \(\x\) and \(\y\), respectively.

We point out that since the above construction involves both the cartesian unstraightening \(\int^{\x \in \A}\) and the cocartesian straightening \(\int_{\y \in \cB}\), the resulting arrow in~\eqrefone{equation:second-unstraightening} is neither a cartesian nor a cocartesian fibration. We can nonetheless describe it as follows: recall that a \defi{bifibration} (see~\cite[Definition 2.4.7.2]{HTT}) is a pair of maps
\[
\A \xleftarrow{q} \X \xrightarrow{p} \cB
\]
consisting of a cartesian fibration \(q\colon \X \to \A\) and a cocartesian fibration \(p\colon\X \to \cB\), such that the \(q\)-cartesian edges are exactly those projecting to equivalences in \(\cB\) and the \(p\)-cocartesian edges are exactly those projecting to equivalences in \(\A\). Equivalently, the pair of maps \(p,q\) forms a bifibration if and only if \((q,p)\colon \X \to \A\times \cB\) is a map of cocartesian fibrations over \(\cB\) whose fibres are right fibrations, and if and only if \((q,p)\colon \X \to \A\times \cB\) is map of cartesian fibrations over \(\A\) whose fibres are left fibrations. In particular, one readily verifies that for a correspondence \(b\colon \A\op \times \cB \to \Sps\), the pair of projections
\[
\A \leftarrow \Pairings(\A,\cB,b) \rightarrow \cB
\]
constitutes a bifibration. In fact, this association determines an equivalence between correspondences \(b\colon \A\op \times \cB \to \Sps\) and bifibrations \(\A \leftarrow \X \rightarrow \cB\) (see ~\cite{stevenson2018model} and~\cite[\S 4]{ayala2020fibrations}), and can be considered as a bivariant form of the space-valued straightening-unstraightening equivalence. In particular, using the coherent compatibility of the straightening-unstraightening equivalence with base change as established in~\cite[Corollary A.31]{GHNfree}, this equivalence integrates to an equivalence
\begin{equation}
\label{equation:bivariant-straightening}%
(-) \leftarrow \Pairings(-,-,-) \rightarrow (-)\colon \int^{(\A,\cB) \in \Cat \times \Cat}\Fun(\A\op \times \cB,\Sps) \xrightarrow{\simeq} \BiFib,
\end{equation}
where \(\BiFib \subseteq \Fun(\Lam^2_0,\Cat)\) is the full subcategory spanned by the bifibrations.

\begin{example}
\label{example:canonical}%
In the case of \(\cB=\A=\C\) and the canonical correspondence \(m_{\C} := \Map_{\C}(-,-)\colon \C\op \times \C\to \Sps\), the \(\infty\)-category \(\Pairings(\C,\C,m_{\C})\) canonically identifies with the arrow category \(\Ar(\C) := \Fun(\Del^1,\C)\).
\end{example}

\begin{remark}
\label{remark:same-as-pairings}%
The notion of a correspondence \(\A\op \times \cB \to \Sps\) can equivalently be encoded by a right fibration \(\M \to \A \times \cB\op\). Right fibrations of this form were studied in~\cite[\S 5.2.1]{HA} under the name \defi{pairings}. In particular, the \(\infty\)-category \(\mathrm{CPair}\) of~\cite[Construction 5.2.1.14]{HA} of pairings is naturally equivalent to \(\BiFib\), since both are equivalent to \(\int^{\A,\cB \in \Cat}\Fun(\A\op \times \cB,\Sps)\). Under this equivalence, the canonical bifibration \(\C \leftarrow \Ar(\C) \rightarrow \C\) of Example~\refone{example:canonical} encoding the mapping space correspondence corresponds to the right fibration \(\Twar(\C) \to \C \times \C\op\).
\end{remark}

\begin{remark}\label{remark:ortho}
The notion of a bifibration was generalised in~\cite{Ortho} to that of an \emph{orthofibration}, which the authors show encodes the data of a functor \(\A\op \times \cB \to \Cat\). In particular, an explicit dualisation procedure transforming an orthofibration over \(\A,\cB\) to that of a cartesian fibration over \(\A \times \cB\op\) is constructed in loc.\ cit.; when applied in the case of bifibrations, this enables one to identify them directly with pairings as in Remark~\refone{remark:same-as-pairings}
without passing through straightening and unstraightening. For example, this dualisation procedure shows directly that the bifibration \(\C \leftarrow \Ar(\C) \rightarrow \C\) corresponds to the right fibration \(\Twar(\C) \to \C \times \C\op\).
The machinery of~\cite{Ortho} also shows that the identification of functors \(\A\op \times \cB \to \Cat\) with orthofibrations (and hence of functors \(\A\op \times \cB \to \Sps\) with bifibrations), is independent of whether we first unstraighten the contravariant part, and then covariant one (as we did above), or the other way around.
\end{remark}

\begin{definition}
\label{definition:right-representable}%
We will say that a
correspondence \(b\colon \A\op \times \cB \to \Sps\) is \defi{right-representable} if the associated map \(\cB \to \Psh(\A)\) factors through the image of the Yoneda embedding \(\iota\colon \A \hrar \Psh(\A)\). In this case the resulting functor \(d\colon \cB \to \A\) is characterised by a natural equivalence \(\Map_{\A}(\x,d(y)) \simeq b(\x,\y)\).
\end{definition}

\begin{remark}
\label{remark:representable}%
If a correspondence \(b\colon \A\op \times \cB \to \Sps\) is right-representable
then \(\Pairings(\A,\cB,b)\) naturally identifies with the fibre product
\(\Ar(\A) \times_{\A} \cB\) along the target projection \(\Ar(\A) \to \A\) and \(d\colon \cB \to \A\). More generally, if we Yoneda embed \(\A\) in \(\Psh(\A) = \Fun(\A\op,\Sps)\), then \(b\) becomes tautologically right representable by the functor \(\wtl{d}\colon \cB \to \Psh(\A)\) sending \(\y \in \cB\) to \(b(-,\y) \in \Psh(\A)\). We may then write
\[
\Pairings(\A,\cB,b) \simeq \A \times_{\Psh(\A)} \Ar(\Psh(\A)) \times_{\Psh(\A)} \cB
\]
where the fibre product is taken along the Yoneda embedding \(\A \to \Psh(\A)\) and \(\wtl{d}\colon \cB \to \Psh(\cB)\), and \(\Ar(\Psh(\A))\) projects to the domain on the left hand side and to the target on the right hand side.
\end{remark}

\begin{definition}
\label{definition:perfect-bifibration}%
We will say that a correspondence \(b\colon \A\op \times \cB \to \Sps\) is \defi{perfect} if \(b\) is right-representable and the associated representing functor \(d\colon \cB \to \A\) is an equivalence.
In this case, we will say that a pairing \((\x,\y,\beta) \in \Pairings(\A,\cB,b)\) is \defi{perfect} if the map \(\x \to d(\y)\) determined by \(\beta\) is an equivalence. Similarly, we will say that a bifibration \(\A \leftarrow \X \rightarrow \cB\) is perfect if its classifying correspondence \(\A\op \times \cB \to \Sps\) is. In this case we will say that an object in \(\X\) is perfect if it corresponds to a perfect pairing under the identification of \(\X\) with \(\Pairings(\A,\cB,b)\). We will then say that a map of bifibration
\[
[\A \leftarrow \X \rightarrow \cB] \to [\A'\leftarrow \X' \rightarrow \cB']
\]
is \defi{perfect} if it sends perfect objects of \(\X\) to perfect objects of \(\X'\). We will denote by \(\BiFibp \subseteq \BiFib\) the (non-full) subcategory spanned by the perfect bifibrations and the perfect maps between them.
\end{definition}

While \(\BiFibp \subseteq \BiFib\) is not a full subcategory, it does satisfy the following weaker property \cite[Definition~20.1.1.2]{SAG}:

\begin{definition}
\label{definition:sound}%
Let \(\C\) be an \(\infty\)-category and \(\C' \subseteq \C\) a subcategory. We will say that \(\C'\) is \defi{replete} if for every \(\x\in \C'\) and \(\y\in\C\) such that there is an equivalence \(\alp\colon \x \xrightarrow{\simeq}\y\) in \(\C\), then \(\y\in\C'\) and there exists an equivalence \(\beta\colon \x \xrightarrow{\simeq} \y'\) in \(\C'\) whose image in \(\C\) is homotopic to \(\alp\).
\end{definition}

\begin{remark}
The condition of being a replete subcategory \(\C'\subseteq \C\) is detected on the level of the homotopy categories: it is equivalent to saying that \(\Ho\C'\) is closed under isomorphisms and for every \(\x,\y \in \C'\) the subset \(\Hom_{\Ho\C'}(\x,\y) \subseteq \Hom_{\Ho\C}(\x,\y)\) contains all isomorphisms from \(\x\) to \(\y\) in \(\C\).
\end{remark}

\begin{example}
\label{example:perfect-sound}%
The subcategory \(\BiFibp \subseteq \BiFib\) is replete. This follows from the observation that any equivalence between perfect bifibrations is necessarily a perfect map.
\end{example}

\begin{example}
\label{example:stable}%
The subcategory \(\Catx \subseteq \Cat\) is replete.
\end{example}

\begin{remark}
\label{remark:sound-subcat-truncated}%
For any subcategory \(\C'\subseteq \C\) the induced map \(\Map_{\C'}(\x,\y) \to \Map_{\C}(\x,\y)\) is a \((-1)\)-truncated map of spaces for every \(\x,\y \in \C'\), that is, its fibres are either empty or contractible. If the subcategory \(\C'\) is replete then the induced map \(\grpcr\C'\to \grpcr\C\) on core groupoids is also \((-1)\)-truncated. This implies that every replete subcategory inclusion \(\C' \hrar \C\) is \((-1)\)-truncated as a map in \(\Cat\), that is, for every test \(\infty\)-category \(\D\) the induced map
\[
\Map_{\Cat}(\D,\C') \to \Map_{\Cat}(\D,\C)
\]
is \((-1)\)-truncated. In many contexts, this property is what makes replete subcategories behave more like ``subobjects'' than general subcategories.
\end{remark}

The following proposition records the content of~\cite[Remark 5.2.1.20]{HA} in the context of bifibrations:

\begin{proposition}
\label{proposition:projection-perfect}
The composite
\begin{equation}
\label{equation:projection-perfect}%
\begin{tikzcd}
[row sep=1ex]
\BiFibp \ar[r,hook] & \BiFib \ar[r] & \Cat \times \Cat \ar[r] & \Cat \\ (\A,\cB,b) \ar[rr,mapsto] && (\A,\cB) \ar[r,mapsto] & \A
\end{tikzcd}
\end{equation}
is an equivalence of \(\infty\)-categories. An inverse is given by \(\C \mapsto [\C \leftarrow \Ar(\C) \rightarrow \C]\).
\end{proposition}
\begin{proof}
Since the association \(\C \mapsto [\C \leftarrow \Ar(\C) \rightarrow \C]\) is visibly a one-sided inverse to~\eqrefone{equation:projection-perfect} we see that the latter is in particular essentially surjective,
and it will hence suffice to show that it is also fully-faithful,
On the other hand, if \(\A \leftarrow \X \rightarrow \cB\) is a perfect bifibration then it is equivalent in \(\BiFib\) to \(\A \leftarrow \Ar\A \rightarrow \A\) via the associated functor \(d\colon \cB \to \A\) and the natural equivalence \(b(\x,\y) \simeq \map_{\A}(\x,d(\y))\). It will hence suffice to show that for every \(\C \in \Cat\) and every perfect bifibration \(\A \leftarrow \X \rightarrow \cB\) the induced map
\[
\Map_{\BiFibp}([\A \leftarrow \X \rightarrow \cB],[\C \leftarrow \Ar(\C) \rightarrow \C]) \to \Map_{\Cat}(\A,\C)
\]
is an equivalence.
But this now follows from~\cite[Proposition 5.2.1.18]{HA} under the equivalence between bifibrations as above and pairings in the sense of~\cite[Definition 5.2.1.5]{HA}, see Remark~\refone{remark:same-as-pairings}.
\end{proof}

The \(\infty\)-category \(\BiFib \subseteq \Fun(\Lam^2_0,\Cat)\) of bifibrations carries a natural action of \(\Ct\) induced by the action of \(\Ct\) on \(\Lam^2_0\) switching the vertices \(\Del^{\{1\}}\) and \(\Del^{\{2\}}\) and post-composing with the \(\mop\) action on \(\Cat\). Explicitly, this action sends a bifibration \(\A \leftarrow \X \rightarrow \cB\) to the bifibration \(\cB\op \leftarrow \X\op \rightarrow \A\op\) (indeed, the latter is again a bifibration since taking opposites switches cartesian and cocartesian fibrations). A \(\Ct\)-fixed structure on a given bifibration \(\A \leftarrow \X \rightarrow \cB\) can then be described as a duality \(\Dual\colon \X \xrightarrow{\simeq} \X\op\) on \(\X\), an equivalence \(\A \simeq \cB\op\), and a duality-preserving refinement of \(\X \to \A \times \cB \simeq \A \times \A\op\). We will refer to such a structure as a \defi{\(\Lam^2_0\)-duality} on \(\A \leftarrow \X \rightarrow \cB\), and will call a bifibration equipped with a \(\Lam^2_0\)-duality a \defi{symmetric bifibration}.

\begin{proposition}
\label{proposition:identify-symmetric}%
The \(\Ct\)-fixed \(\infty\)-category \(\BiFib^{\hC}\) participates in a cartesian fibration
\[
\BiFib^{\hC} \to \Cat
\]
classified by the functor \(\C \mapsto \Fun(\C\op \times \C\op,\Sps)^{\hC}\).
\end{proposition}
\begin{proof}
Equipping \(\Del^{\{1\}} \coprod \Del^{\{2\}}\) with the swap action and restricting along the \(\Ct\)-equivariant inclusion \(\Del^{\{1\}} \coprod \Del^{\{2\}} \subseteq \Lam^2_0\) we get that the cocartesian fibration
\begin{equation}
\label{equation:bifib-forget}%
\BiFib \to \Cat \times \Cat,
\end{equation}
naturally refines to a \(\Ct\)-equivariant functor, where \(\Ct\) acts on the target by flipping the factors and taking opposites.
Since cartesian fibration are preserved under limits, taking \(\Ct\)-fixed points results in a cartesian fibration
\begin{equation}
\label{equation:fixed-cocar-fibration}%
\BiFib^{\hC} \to (\Cat \times \Cat)^{\hC}.
\end{equation}
We now observe that the equivalence
\[
(\id,(-)\op)\colon \Cat \times \Cat \xrightarrow{\simeq} \Cat \times \Cat
\]
intertwines the flip-\textrm{op} action on the left hand side with the flip action on the right hand side, and so the target of~\eqrefone{equation:bifib-forget} is a coinduced \(\Ct\)-object. We may consequently identify the target of~\eqrefone{equation:fixed-cocar-fibration} with \(\Cat\) and write it as a cocartesian fibration
\begin{equation}
\label{equation:fixed-cocar-fibration-2}%
\BiFib^{\hC} \to \Cat .
\end{equation}
Since taking fibres commutes with fixed points we may identify the fibres of~\eqrefone{equation:fixed-cocar-fibration-2} with the \(\Ct\)-fixed points of the corresponding fibres
of~\eqrefone{equation:bifib-forget}. Now \(\C \in \Cat\) corresponds to the fixed object \((\C,\C\op) \in \Cat \times \Cat\), and so the fibre of~\eqrefone{equation:fixed-cocar-fibration-2} over \(\C\) is the \(\Ct\)-fixed points of \(\Fun(\C\op \times \C\op,\Sps)\), as claimed.
\end{proof}

\begin{construction}
\label{construction:duality-pairing}%
By Proposition~\refone{proposition:identify-symmetric} we may identify the notion of a symmetric bifibration with that of a pair \((\C,b)\) where \(\C\) is an \(\infty\)-category and \(b\in\Fun(\C\op\times\C\op,\Sps)^\hC\) is a symmetric functor. Given \((\C,b)\), the associated symmetric bifibration is
\[
\C\op \leftarrow \Pairings(\C,\C\op,b) \rightarrow \C
\]
equipped with its \(\Lam^2_0\)-duality which we will denote by \(\Dual_{\pair}\).
It is given explicitly by the duality
\[
\Dual_{\pair}(\x,\y,\beta) = \big(\y,\x,\sig_{\x,\y}(\beta)\big) ,
\]
on \(\Pairings(\C,\C\op,b)\), where \(\sig_{\x,\y}\colon b(\x,\y) \xrightarrow{\simeq} b(\y,\x)\) is given by the symmetric structure of \(b\), equipped with the tautological duality-preserving structure of the projection \(\Pairings(\C,\C\op,b) \to \C \times \C\op\).
\end{construction}

\begin{proposition}
\label{proposition:action-perfect}%
The \(\Ct\)-action on \(\BiFib\) restricts to a \(\Ct\)-action on \(\BiFibp\). Under the equivalence \(\BiFibp \simeq \Cat\) of Proposition~\refone{proposition:projection-perfect}, this action corresponds to the \(\mop\)-action on \(\Cat\). In particular, we may identify the notion of a perfect symmetric bifibration with that of an \(\infty\)-category equipped with a perfect duality.
\end{proposition}
\begin{proof}
We claim that the functor
\[
\Cat \to \BiFib \quad\quad \C \mapsto [\C \leftarrow \Ar(\C) \rightarrow \C]
\]
admits a natural \(\Ct\)-equivariant structure, where \(\Ct\) acts on \(\Cat\) via the \(\mop\)-action. To see this, it will suffice to construct a \(\Ct\)-equivariant structure for the composite map \(\Cat \to \BiFib \hrar \Fun(\Lam^2_0,\Cat)\). This composite is by definition given by mapping out of the diagram \(e := [\Del^{\{0\}} \rightarrow \Del^1 \leftarrow \Del^{\{1\}}]\), and so it will suffice to put a \(\Ct\)-equivariant structure on \(e\colon (\Lam^2_0)\op \to \Cat\).
Such an equivariant structure is then given by the canonical duality on \(\Dual_{\Del^1}\colon \Del^1 \to (\Del^1)\op\) which switches \(\{0\}\) and \(\{1\}\) (since \(\Del^1\) is the nerve of an ordinary category not much coherence is needed in order to verify this fact).

By Proposition~\refone{proposition:projection-perfect} it now follows in particular that the \(\Ct\)-action on \(\BiFib\) preserves the subcategory \(\BiFibp\). This subcategory is \emph{replete} by Example~\refone{example:perfect-sound},
and so by Remark~\refone{remark:sound-subcat-truncated}
the \(\Ct\)-action on \(\BiFib\) restricts to an essentially unique \(\Ct\)-action on \(\BiFibp\) making the inclusion \(\BiFibp \hrar \BiFib\) equivariant. By the above this action must then coincide with the \(\mop\)-action via the equivalence \(\BiFibp \simeq \Cat\), as desired.
\end{proof}

\begin{example}
\label{example:swap-duality-arrows}%
In the situation of Construction~\refone{construction:duality-pairing}, if the symmetric correspondence \(b\colon \C\op \times \C\op \to \Sps\) is perfect with duality \(\Dual\colon \C\op \xrightarrow{\simeq} \C\)
then the associated symmetric bifibration identifies by Proposition~\refone{proposition:projection-perfect} with
\[
\C \leftarrow \Ar(\C) \rightarrow \C
\]
and by Proposition~\refone{proposition:action-perfect} (and its proof) the associated \(\Lam^2_0\)-duality \(\Dual_{\pair}\) corresponds to the arrow duality \([\x \to \y] \mapsto [\Dual\y \to \Dual\x]\) induced on the functor category \(\Ar(\C) = \Fun(\Del^1,\C)\) from the dualities of \(\Del^1\) and \(\C\).
\end{example}

The remainder of this section is devoted to producing an explicit formula expressing the mapping spaces in \(\Pairings(\A,\cB,b)\) in terms of \(b\) and the mapping spaces in \(\A \times \cB\). This will be useful for us in \S\refone{subsection:thom} when we will need to upgrade the pairings construction to the hermitian setting.
To begin, consider the following diagram in \(\Fun(\Lam^2_0,\Cat)\):
\[
\begin{tikzcd}
\left[\partial\Del^{\{0,1\}} \leftarrow \partial\Del^{\{0,1\}} \amalg \partial\Del^{\{1,2\}} \rightarrow \partial\Del^{\{1,2\}}\right]  \ar[d]\ar[r] &
\left[\partial\Del^{\{0,1\}} \leftarrow \Del^{\{0\}} \amalg \Del^{\{1\}}\amalg \Del^{\{2\}} \rightarrow \partial\Del^{\{1,2\}}\right] \ar[d] \\
\left[\Del^{\{0,1\}} \leftarrow \Del^{\{0,1\}} \coprod \Del^{\{1,2\}} \rightarrow \Del^{\{1,2\}}\right] \ar[r] &
\left[\Del^{\{0,1\}} \leftarrow \Del^2 \rightarrow \Del^{\{1,2\}}\right] \\
\left[\Del^{\{0,1\}} \leftarrow \Del^{\{0\}} \coprod \Del^{\{2\}} \rightarrow \Del^{\{1,2\}}\right] \ar[r]\ar[u] &
\left[\Del^{\{0,1\}} \leftarrow \Del^{\{0,2\}} \rightarrow \Del^{\{1,2\}}\right]\ar[u] \ .
\end{tikzcd}
\]
Here, left going internal arrows are all given on vertices by \([0 \mapsto 0],[1 \mapsto 1],[2 \mapsto 1]\) and all the right going internal arrows are given by \([0 \mapsto 1],[1 \mapsto 1], [2\mapsto 2]\), while the external arrows always preserve the vertex labels.
The upper square is cocartesian, as can be tested levelwise using the standard categorical equivalence \(\Del^{\{0,1\}} \coprod_{\Del^{\{1\}}} \Del^{\{1,2\}} \xrightarrow{\simeq} \Del^2\). In addition, all entries in this square carry compatible \(\Lam^2_0\)-dualities, induced by the canonical duality \(\Del^2 \xrightarrow{\simeq} (\Del^2)\op\) which switches between \(0\) and \(2\) and between \(\Del^{\{0,1\}}\) and \((\Del^{\{1,2\}})\op\).
Mapping out of the above diagram now yields a \(\Ct\)-equivariant functor \(\Fun(\Lam^2_0,\Cat) \to \Fun(\Lam^2_0 \times \Del^1,\Cat)\) sending \([\A \xleftarrow{q} \X \xrightarrow{p} \cB]\) to the diagram
\begin{equation}
\label{equation:span-decomposition}%
\begin{tikzcd}
\X^{\Del^{\{0\}}} \times_{\cB^{\Del^{\{1\}}}} \X^{\Del^{\{1\}}} \times_{\A^{\Del^{\{1\}}}} \X^{\Del^{\{2\}}} \ar[r] &
\X^{\Del^{\{0\}}} \times_{\cB^{\Del^{\{1\}}}} \X^{\Del^{\{1\}}} \times_{(\A \times \cB)^{\Del^{\{1\}}}}\X^{\Del^{\{1\}}} \times_{\A^{\Del^{\{1\}}}} \X^{\Del^{\{2\}}} \\
\X^{\Del^2}\displaystyle\mathop{\times}_{(\A \times \cB)^{\Del^2}}\big[\A^{\Del^{\{0,1\}}} \times \cB^{\Del^{\{1,2\}}}\big] \ar[u]\ar[d]\ar[r] &
\X^{\Del^{\{0,1\}} \amalg \Del^{\{1,2\}}} \displaystyle\mathop{\times}_{(\A \times \cB)^{\Del^{\{0,1\}}\amalg\Del^{\{1,2\}}}} \big[\A^{\Del^{\{0,1\}}} \times \cB^{\Del^{\{1,2\}}}\big] \ar[u]\ar[d] \\
\X^{\Del^{\{0,2\}}}\displaystyle\mathop{\times}_{(\A \times \cB)^{\Del^{\{0,2\}}}}\big[\A^{\Del^{\{0,1\}}} \times \cB^{\Del^{\{1,2\}}}\big] \ar[r] &
\X^{\Del^{\{0\}} \amalg \Del^{\{2\}}} \displaystyle\mathop{\times}_{(\A \times \cB)^{\Del^{\{0\}} \amalg \Del^{\{2\}}}} \big[\A^{\Del^{\{0,1\}}} \times \cB^{\Del^{\{1,2\}}}\big]
\end{tikzcd}
\end{equation}
in which the top square is cartesian. In addition, using that all entries in this diagram compatibly project to \(\X^{\Del^{\{0\}}}\) and \(\X^{\Del^{\{2\}}}\) we will view this as a diagram in \((\Cat)_{/\X \times \X}\).

\begin{proposition}
\label{proposition:pairings-mapping-spaces}%
When \([\A \leftarrow \X \rightarrow \cB]\) is a bifibration the bottom vertical arrows in~\eqrefone{equation:span-decomposition} are equivalences. In particular, inverting these and taking the external rectangle yields a \(\Ct\)-equivariant functor \(\BiFib \to \Fun(\Del^1 \times \Del^1,\Fun(\Lam^2_0,\Cat))\) which sends \([\A \leftarrow \X \rightarrow \cB]\) to a cartesian square~\eqrefone{equation:bifib-decomposition} of the form
\begin{equation}
\label{equation:bifib-decomposition}%
\begin{tikzcd}
\Ar(\X) \ar[r]\ar[d] & \X \times_{\A} \X \times_{\cB} \X \ar[d] \\
\X  \times_{\A \times \cB} [\Ar(\A \times \cB)] \times_{\A \times \cB} \X \ar[r] & \X \times_{\A} \X \times_{\A \times \cB} \X \times_{\cB} \X.
\end{tikzcd}
\end{equation}
in \((\Cat)_{/\X \times \X}\). Here, the two projections to \(\X\) are given by the domain and codomain projections in the case of \(\Ar(\X)\) and by the projection to the two extremal factors in the three other cases. In addition, this functor takes values in \(\Fun(\Del^1 \times \Del^1,\BiFib)\) after post-composing with the inclusion \((\Cat)_{/\X \times \X} \to \Fun(\Lam^2_0,\Cat)\).
\end{proposition}

Writing \(\X = \Pairings(\A,\cB,b)\) for some correspondence \(b\colon \A\op \times \cB \to \Sps\) the square of bifibrations ~\eqrefone{equation:bifib-decomposition} corresponds to a cartesian square of correspondences \(\X\op \times \X \to \Sps\) of the form
\begin{equation}
\label{equation:mapping-space-formula}%
\begin{tikzcd}
\Map_{\X}((\x,\y,\beta),(\x',\y',\beta')) \ar[d] \ar[r] & b(\x,\y') \ar[d] \\
\Map_{\A}(\x,\x') \times \Map_{\cB}(\y,\y') \ar[r] & b(\x,y') \times b(\x,\y') ,
\end{tikzcd}
\end{equation}
giving, in particular, an explicit pullback formula for the mapping spaces in \(\Pairings(\A,\cB,b)\).

\begin{remark}
\label{remark:mapping-formula-duality-preserving}%
In the situation of Proposition~\refone{proposition:pairings-mapping-spaces}, the \(\Ct\)-equivariance of the functor in question means in particular that if a bifibration \(\A \leftarrow \X \rightarrow \cB\) carries a \(\Lam^2_0\)-duality then all the entries in the square~\eqrefone{equation:bifib-decomposition} inherit such a duality and all arrows in the square are duality preserving. %
\end{remark}

\begin{proof}[Proof of Proposition~\refone{proposition:pairings-mapping-spaces}]
By definition the arrows in \(\X\) which map to equivalences in \(\A\) are exactly the \(p\)-cocartesian arrows, and the arrows which map to equivalences in \(\cB\) are exactly the \(q\)-cartesian arrows.
It then follows that the projections
\[
\X^{\Del^{\{0,1\}}} \times_{\cB^{\Del^{\{0,1\}}}}\cB \to \X^{\Del^{\{0\}}} \times_{\A^{\Del^{\{0\}}}} \A^{\Del^{\{0,1\}}}
\]
and
\[
\X^{\Del^{\{1,2\}}} \times_{\A^{\Del^{\{1,2\}}}}\A \to \X^{\Del^{\{2\}}} \times_{\cB^{\Del^{\{2\}}}} \cB^{\Del^{\{1,2\}}}
\]
are equivalences, and so the right bottom vertical arrow in~\eqrefone{equation:span-decomposition} (which is the fibre product of these two maps over the identity on \(\A \times \cB\)) is an equivalence. Similarly, the projections
\[
\X^{\Del^2} \times_{\cB^{\Del^{\{0,1\}}}} \cB \to \A^{\Del^2} \times_{\A^{\Lam^2_0}} \X^{\Lam^2_0}\times_{\cB^{\Del^{\{0,1\}}}} \cB \to \A^{\Del^2} \times_{\A^{\Del^{\{0,2\}}}} \X^{\Del^{\{0,2\}}}
\]
are both equivalences, from which it follows that the left bottom vertical arrow in~\eqrefone{equation:span-decomposition}, which is a base change of the above composite, is an equivalence.
\end{proof}

\subsection{Bilinear and symmetric \(\infty\)-categories}
\label{subsection:bilinear-and-pairings}%

In the present section we define and study the notion of \defi{bilinear} and \defi{perfect bilinear} \(\infty\)-categories. We then show that the notion of an \(\infty\)-category equipped with a symmetric bilinear form, which we call a \defi{symmetric} \(\infty\)-category, can be identified with a \(\Ct\)-fixed bilinear \(\infty\)-category, and similarly in the perfect case. Using the results of \S\refone{subsection:bifibrations} we then identify
the \(\infty\)-category of perfect bilinear \(\infty\)-categories with \(\Catx\) itself, through which we also deduce an equivalence between \(\infty\)-categories with perfect dualities and \(\Ct\)-fixed objects of \(\Catx\) with respect to the \(\mop\)-action. We then construct left and right adjoints to the forgetful functors from hermitian to symmetric to bilinear \(\infty\)-categories, and similarly in the Poincaré/perfect case, where we finally recover the hyperbolic construction \(\C \mapsto \Hyp(\C)\) acting as both left and right adjoint to the forgetful functor \(\Catp \to \Catx\).

Let \(\A,\cB\) be stable \(\infty\)-categories. We will say that a correspondence \(b\colon \A\op \times \cB \to \Sps\) is \defi{left biexact} if it preserves finite limits in each variable separately. Such a correspondence lifts in an essentially unique manner to a bilinear functor \(\A\op \times \cB \to \Spa\). More precisely, post-composition with the infinite loop space functor induces an equivalence between bilinear functors \(\A\op \times \cB \to \Spa\) and left biexact correspondences \(\A\op \times \cB \to \Sps\).

\begin{definition}
We will denote by \(\Funb(\A,\cB) \subseteq \Fun(\A\op \times \cB ,\Spa)\) the full subcategory spanned by the bilinear functors and write
\[
\Catb := \int^{(\A,\cB) \in \Catx \times \Catx}\Funb(\A,\cB)
\]
for the \(\infty\)-category of triples \((\A,\cB,\Bil)\) where \(\A,\cB\) are stable \(\infty\)-categories and \(\Bil\colon \A\op \times \cB \to \Spa\) is a bilinear functor. We will refer to the object \((\A,\cB,\Bil)\) of \(\Catb\) as \defi{bilinear categories}.
\end{definition}

\begin{example}
For a stable \(\infty\)-category \(\C\) the mapping correspondence \(m_{\C}\colon \C\op \times \C \to \Sps\) of Example~\refone{example:canonical} is left biexact. We may hence consider the triple \((\C,\C,m_{\C})\) as a bilinear \(\infty\)-category.
\end{example}

\begin{example}
\label{example:rep-biexact}%
Let \(\A,\cB\) be stable \(\infty\)-categories. If a correspondence \(b\colon \A\op \times \cB \to \Sps\) is right representable by an exact functor \(d\colon \A \to \cB\) then \(b\) is left biexact.
\end{example}

\begin{example}
\label{example:canonical-shift}%
Let \(\C\) be a stable \(\infty\)-category. Consider the \(\infty\)-category \(\Seq(\C)\) of exact sequences in \(\C\) (see Notation~\refone{notation:seq}).
The pair of projections
\[
\begin{tikzcd}
[row sep=1ex]
\C & \Seq(\C) \ar[l] \ar[r] & \C \\
\x & \left[\y \to \z \to \x\right] \ar[l,mapsto] \ar[r,mapsto] & \y
\end{tikzcd}
\]
then constitute a bifibration. Indeed, a map of exact sequences
\[
\begin{tikzcd}
\y \ar[r] \ar[d] & \z \ar[d] \ar[r] & \x \ar[d] \\
\y' \ar[r] & \z'\ar[r] & \x'
\end{tikzcd}
\]
is a cocartesian lift of \(\y \to \y'\) if and only if the left square is exact, or, equivalently, if the map \(\x \to \x'\) is an equivalence, and similarly forms a cartesian lift of \(\x \to \x'\) if and only if its component \(\y \to \y'\) is an equivalence. To identify the correspondence \(\seq_{\C}\colon \C\op \times \C \to \Sps\) associated to this bifibration we use the fact that every exact square as in~\eqrefone{equation:exact-square} extends in an essentially unique manner to a diagram
\[
\begin{tikzcd}
\x' \ar[d] \ar[r] & 0 \ar[r] \ar[d] & 0 \ar[d] \\
\y \ar[r] \ar[d] & \z \ar[r] \ar[d] & 0 \ar[d] \\
0 \ar[r] & \x \ar[r] & y'
\end{tikzcd}
\]
in which all squares except the top right one are exact. We note that such a diagram determines in particular equivalences \(\x' \simeq \Om \x\) and \(\y' \simeq \Sig \y\). At the same time, the forgetful functor sending a diagram as a above to its external square
\begin{equation}
\label{equation:shifted-map}%
\begin{tikzcd}
\x' \ar[d] \ar[r] & 0 \ar[d] \\
0 \ar[r] & \y'
\end{tikzcd}
\end{equation}
is an equivalence as well. In particular, if we denote by \(\Ar^{\Om}(\C) \subseteq \Fun(\Del^1 \times \Del^1,\C)\) the full subcategory spanned by the squares of the form~\eqrefone{equation:shifted-map} then we obtain an equivalence of bifibrations
\[
\begin{tikzcd}
\C \ar[d,"{\Om}"'] & \Seq(\C) \ar[l] \ar[r] \ar[d,"{\simeq}"] & \C \ar[d,"\Sig"] \ar[d] \\
\C & \Ar^{\Om}(\C) \ar[l] \ar[r] & \C
\end{tikzcd}
\]
The correspondence associated to \(\Ar^{\Om}(\C)\) can then be identified with \((\x',\y') \mapsto \Om\Map_{\C}(\x',\y')\), and so the correspondence associated to \(\Seq\) can be written as \(\seq_{\C}(\x,\y) = \Om\Map(\Om \x,\Sig\y)\). In particular, it is left biexact.
\end{example}

\begin{remark}
\label{remark:catb-has-all-limits-and-colimits}%
It follows from Lemma~\refone{lemma:kan-extension-exact-quadratic} that the defining cartesian fibration \(\Catb \to \Catx \times \Catx\) is also a cocartesian fibration. Applying the precise same argument as in the proof of Proposition~\refone{proposition:cath-has-limits-and-colimits} we may consequently conclude that \(\Catb\) has all small limits and colimits, and that those are preserved by the projection to \(\Catx \times \Catx\).
\end{remark}

Given a bilinear category \((\A,\cB,\Bil)\) we may consider the pairings \(\infty\)-category associated to the underlying left biexact correspondence \(\Om^{\infty}\Bil\). To simplify notation we will denote
\[
\Pairings(\A,\cB,\Bil) := \Pairings(\A,\cB,\Om^{\infty}\Bil).
\]
We note that since \(\Om^{\infty}\Bil\) is left biexact
the functor \(\y \mapsto \wtl{d}(\y) = \Om^{\infty}\Bil(-,\y)\) from \(\cB\) to \(\Psh(\A)\) takes values in the full subcategory \(\Ind(\A) \subseteq \Psh(\A)\) spanned by the left exact presheaves. In addition, in this case \(\Ind(\A)\) is also stable, the Yoneda embedding \(\A \to \Ind(\A)\) is exact, and the functor \(\wtl{d}\colon \cB \to \Ind(\A)\) given by \(\Om^{\infty}\Bil\) is exact as well. As in Remark~\refone{remark:representable} we may then identify
\begin{equation}
\label{equation:representable}%
\Pairings(\A,\cB,\Bil) \simeq \A \times_{\Ind(\A)} \Ar(\Ind(\A)) \times_{\Ind(\A)} \cB ,
\end{equation}
a description from which we see that
\(\Pairings(\A,\cB,\Bil)\) is stable and that exact squares in \(\Pairings(\A,\cB,\Bil)\) are detected in \(\A \times \cB\). This means in particular that any map of bifibrations
\[
[\A \leftarrow \Pairings(\A,\cB,\Bil) \rightarrow \cB] \to [\A' \leftarrow \Pairings(\A',\cB',\Bil') \rightarrow \cB']
\]
whose component \(\A \to \A',\cB \to \cB'\) are exact, is also exact on \(\Pairings(-,-,-)\). Invoking again the coherent compatibility of the straightening-unstraightening equivalence with base change as established in~\cite[Corollary A.31]{GHNfree}, we may consequently assemble the association \((\A,\cB,\Bil) \mapsto \Pairings(\A,\cB,\Bil)\) to a functor
\[
\Pairings(-,-,-)\colon \Catb \to \Catx ,
\]
taking values in stable \(\infty\)-categories and exact functors between them. It then follows that the bivariant straightening equivalence~\eqrefone{equation:bivariant-straightening} restricts to a (non-full) subcategory inclusion
\begin{equation}
\label{equation:catb-to-bifib}%
\Catb \to \BiFib \quad\quad (\A,\cB,\Bil) \mapsto [\A \leftarrow \Pairings(\A,\cB,\Bil) \rightarrow \cB]
\end{equation}
whose image is spanned by those bifibrations \(\A \leftarrow \X \rightarrow \cB\) for which \(\A,\cB\) are stable and the associated correspondence is left biexact and by those maps of bifibrations whose components are exact functors.

\begin{remark}
The subcategory inclusion~\eqrefone{equation:catb-to-bifib} is replete (see Definition~\refone{definition:sound}).
\end{remark}

\begin{remark}
\label{remark:catb-pre-add}%
The association in~\eqrefone{equation:catb-to-bifib} can also be viewed as a functor from \(\Catb\) to \(\Fun(\Lam^2_0,\Catx)\). As such, it is fully-faithful and its image is spanned by those diagrams \(\A \leftarrow \X \rightarrow \cB\) in \(\Catx\) which are bifibrations with associated correspondence left biexact. From this description we see that \(\Catb\) is closed under finite products in \(\Fun(\Lam^2_0,\Catx)\), and hence inherits from its the property of being semiadditive, see Proposition~\refone{proposition:catp-pre-add}.
\end{remark}

\begin{definition}
\label{definition:perfect-bilinear}%
We will say that a bilinear category \((\A,\cB,\Bil)\) is \defi{perfect} if its corresponding bifibration is perfect in the sense of Definition~\refone{definition:perfect-bifibration}, and say that a map of perfect bilinear categories is perfect if the corresponding map of bifibrations is so.
We will denote by \(\Catpb \subseteq \Catb\) the (non-full) subcategory spanned by the perfect bilinear categories and perfect maps between them.
\end{definition}

\begin{remark}
As in Example~\refone{example:perfect-sound}, the subcategory inclusion \(\Catpb \subseteq \Catb\) is replete.
\end{remark}

By construction we have a commutative diagram
\[
\begin{tikzcd}
\Catpb \ar[r]\ar[d] & \Catb \ar[d]\ar[r] & \Catx \times \Catx \ar[d] \\
\BiFibp \ar[r] & \BiFib \ar[r] & \Cat \times \Cat
\end{tikzcd}
\]
in which both squares are cartesian (the right one because (bi)exact functors to spectra are determined by their underlying space-valued functors) and all vertical arrows, as well as the horizontal arrows in the left square, are replete subcategory inclusions.

\begin{proposition}
\label{proposition:projection-perfect-stable}%
The composite functor
\begin{equation}
\label{equation:projection-perfect-stable}%
\begin{tikzcd}
[row sep=1ex]
\Catpb \ar[r,hook] & \Catb \ar[r] & \Catx \times \Catx \ar[r] & \Catx \\ (\A,\cB,b) \ar[rr,mapsto] && (\A,\cB) \ar[r,mapsto] & \A\\
\end{tikzcd}
\end{equation}
is an equivalence of \(\infty\)-categories. An inverse is given by \(\A \mapsto (\A,\A,m_{\A})\).
\end{proposition}
\begin{proof}
This follows directly from Proposition~\refone{proposition:projection-perfect} since a perfect bifibration \(\A \leftarrow \X \rightarrow \cB\) belongs to the image of \(\Catpb\) if and only if \(\A \simeq \cB\) is stable (in which case the associated correspondence is automatically left biexact by Example~\refone{example:rep-biexact}).
\end{proof}

Recall from \S\refone{subsection:bifibrations} that the \(\infty\)-category \(\BiFib \subseteq \Fun(\Lam^2_0,\Cat)\) of bifibrations
carries a natural action of \(\Ct\)
induced by the action of \(\Ct\) on \(\Lam^2_0\) switching the vertices \(\Del^{\{1\}}\) and \(\Del^{\{2\}}\) and post-composing with the \(\mop\) action on \(\Cat\). Explicitly, this action
sends a bifibration \(\A \leftarrow \X \rightarrow \cB\) to the bifibration \(\cB\op \leftarrow \X\op \rightarrow \A\op\).
Since taking opposites also preserves stable \(\infty\)-categories and exact functors, the above action induces a \(\Ct\)-action on the replete subcategory \(\Catb \hrar \BiFib\).
Explicitly, the resulting \(\Ct\)-action sends a triple \((\A,\cB,\Bil)\) to the triple \((\cB\op,\A,\Bil_{\swap})\), where \(\Bil_{\swap}\colon \cB \times \A\op \to \Sps\) is \(\Bil\) pre-composed with the swap equivalence \(\cB \times \A\op \simeq \A\op \times \cB\).

\begin{definition}
\label{definition:symmetric-cat}%
We will denote by \(\Catsb := (\Catb)^{\hC}\) the \(\infty\)-category of homotopy fixed points in \(\Catb\) with respect to the above \(\Ct\)-action.
\end{definition}

\begin{proposition}
\label{proposition:identify-symmetric-stable}%
The
\(\infty\)-category \(\Catsb\) participates in a cartesian fibration
\[
\Catsb \to \Catx
\]
classified by the functor \(\C \mapsto \Funs(\C)\).
\end{proposition}
\begin{proof}
The claim follows from its unstable counterpart Proposition~\refone{proposition:identify-symmetric} since the restriction along exact functors preserves left biexact correspondences.
\end{proof}

By Proposition~\refone{proposition:identify-symmetric-stable} we may identify the objects of \(\Catsb\) with pairs \((\C,\Bil)\) where \(\C\) is a small stable \(\infty\)-category and \(\Bil\colon \C\op \times \C\op \to \Spa\) is a symmetric bilinear functor, that is an object of the \(\infty\)-category \(\Funs(\C) = \Funb(\C)^{\hC}\). We will refer to such a pair as a \defi{symmetric} \(\infty\)-category.
We will refer to morphisms \((\C,\Bil) \to (\Ctwo,\Biltwo)\) in \(\Catsb\) as \defi{symmetric functors}. Using again Proposition~\refone{proposition:identify-symmetric-stable} we may identify these with pairs \((f,\beta)\) where \(f\colon \C \to \Ctwo\) is an exact functor and \(\beta\colon \Bil \to (f \times f)^*\Biltwo\) is a natural transformation.

\begin{definition}
\label{definition:perfect-symmetric}%
We will say that a symmetric bilinear \(\infty\)-category \((\C,\Bil)\) is \defi{non-degenerate} if \(\Bil\) is non-degenerate in the sense of Definition~\refone{definition:nondeg-herm-cat}. In this case \(\Bil\) is induced by a (possibly imperfect) duality \(\Dual_{\Bil}\colon \C \to\C\op\), and every symmetric functor \((f,\beta)\colon (\C,\Bil) \to (\Ctwo,\Biltwo)\) induces a natural transformation \(\tau_{\beta}\colon f\Dual_{\Bil} \Rightarrow \Dual_{\Biltwo}f\op\) via Lemma~\refone{lemma:nat-duality}. We will say that such a symmetric functor \((f,\beta)\) is \defi{duality preserving} if \(\tau_{\beta}\) is an equivalence. We will say that \((\C,\Bil)\) is \defi{perfect} if \(\Dual\) is an equivalence.
We then denote by \(\Catps \subseteq \Catsb\) the (non-full) subcategory spanned by the perfect symmetric \(\infty\)-categories and duality preserving functors between them.
\end{definition}

\begin{lemma}
\label{lemma:perfect-is-perfect}%
The \(\Ct\)-action on \(\Catb\) preserves the replete subcategory \(\Catpb\) of perfect bilinear \(\infty\)-categories. In addition, on the level of homotopy fixed points the resulting replete subcategory
\[(\Catpb)^{\hC} \hrar (\Catb)^{\hC} = \Catsb\]
coincides with the replete subcategory \(\Catps\) of perfect symmetric bilinear \(\infty\)-categories.
\end{lemma}
\begin{proof}
We need to verify two things:
\begin{enumerate}
\item
\label{item:first-thing}
If \((\C,\Bil)\) is a symmetric \(\infty\)-category then \(\Bil\) is perfect in the sense of Definition~\refone{definition:perfect-symmetric} if and only if the left biexact correspondence \(\Om^{\infty}\Bil\) is perfect in the sense of Definition~\refone{definition:perfect-bifibration}.
\item
\label{item:second-thing}%
A symmetric functor \((f,\beta)\colon (\C,\Bil) \to (\Ctwo,\Biltwo)\) between perfect symmetric \(\infty\)-categories is duality preserving if and only if the induced functor \(\Pairings(\C,\C\op,\Bil) \to \Pairings(\Ctwo,\Ctwo\op,\Biltwo)\) preserves perfect pairings.
\end{enumerate}
To prove \refoneitem{item:first-thing} we need to show that \(\Bil\colon \C\op \times \C\op \to \Spa\) can be written as \(\Bil(\x,\y) \simeq \map_{\C}(\x,\Dual\y)\) for some equivalence \(\Dual\colon \C\op \to \C\) if and only if \(\Om^{\infty}\Bil\) can be written as \(\Map_{\C}(\x,\Dual\y)\) for some equivalence \(\Dual\colon \C\op \to \C\). Clearly the former implies the latter, but the latter also implies the former since post-composition with \(\Om^{\infty}\) induces an equivalence between bilinear functors \(\C\op \times \C\op \to \Spa\) and left biexact correspondences \(\C\op \times \C\op \to \Sps\).
To verify \refoneitem{item:second-thing}, consider the commutative diagram
\[
\begin{tikzcd}
\Om^{\infty}\Bil(\x,\y) \ar[rr,"{\Om^{\infty}\beta}"] \ar[d,"{\simeq}"'] & & \Om^{\infty}\Biltwo(f(\x),f(\y)) \\
\Map_{\C}(\x,\Dual_{\Bil}(\y)) \ar[r] & \Map_{\Ctwo}(f(\x),f\Dual_{\Bil}(\y)) \ar[r,"{(\tau_{\beta})_{\ast}}"] & \Map_{\Ctwo}(f(\x),\Dual_{\Biltwo}f(\y)) \ar[u,"{\simeq}"]
\end{tikzcd}
\]
furnished by Remark~\refone{remark:recovered}. We then see that the map \(f\Dual_{\Bil}(y) \to \Dual_{\Biltwo}f(y)\) is an equivalence if and only if \(f\) sends the tautological perfect pairing \((\Dual_{\Bil}(\y),\y,\iota) \in \Pairings(\C,\C\op,\Bil)\) to a perfect pairing in \(\Ctwo\). Since every perfect pairing is equivalent to a tautological perfect pairing \((\Dual_{\Bil}(\y),\y,\iota)\) for some \(y\) we may thus conclude that \(\tau_{\bet}\colon f\Dual_{\Bil} \Rightarrow \Dual_{\Biltwo}f\op\) is an equivalence if and only if the induced functor \(\Pairings(\C,\C\op,\Bil) \to \Pairings(\Ctwo,\Ctwo\op,\Biltwo)\) preserves perfect pairings, as desired.
\end{proof}

\begin{remark}\label{remark:action-op-stable}
Under the equivalence \(\Catpb \simeq \Catx\) of Proposition~\refone{proposition:projection-perfect-stable}, the restricted \(\Ct\)-action on \(\Catpb\) recovers the \(\mop\)-action on \(\Catx\). Indeed, this follows from the corresponding unstable statement on the level of \(\BiFibp\) and \(\Cat\) (see Proposition~\refone{proposition:action-perfect})
since in both cases the \(\Ct\)-action is the restricted one (uniquely determined since the inclusions \(\Catpb \subseteq \BiFibp\) and \(\Catx \subseteq \Cat\) are both replete).
\end{remark}

Combining Proposition~\refone{proposition:projection-perfect-stable}, Lemma~\refone{lemma:perfect-is-perfect} and Remark~\refone{remark:action-op-stable} and we now obtain the following statement, versions of which were previously established in \cite{HeineLopez-AvilaSpitzweck} and \cite{Spitzweckreal}:

\begin{corollary}
\label{corollary:action-perfect-bilinear}%
The forgetful functor \(\Catps \to \Catx\) lifts to an equivalence
\[
\Catps \simeq (\Catx)^{\hC}
\]
where \(\Ct\) acts on \(\Catx\) by the \(\mop\)-action.
\end{corollary}

The constructions we have made so far can now be summarised by the following commutative diagram
\begin{equation}
\label{equation:cross-effect-catp}%
\begin{tikzcd}
[column sep=2ex]
\Catp \ar[d,hook] \ar[rrr] &&& \Catps \ar[d,hook] \ar[r,phantom,"{\simeq}"] & (\Catx)^{\hC} \ar[d,hook] \ar[rr] && \Catx\ar[d,hook]  \\
\Cath \ar[rrr] \ar[d] &&& \Catsb \ar[d] \ar[r,phantom,"{=}"] & (\Catb)^{\hC} \ar[d] \ar[rr] && \Catb \ar[d] \\
\Catx \ar[rrr,equal] &&& \Catx \ar[r,phantom,"{\simeq}"] & (\Catx \times \Catx)^{\hC} \ar[rr] && \Catx \times \Catx
\end{tikzcd}
\end{equation}
in which the vertical arrows in the top row are replete subcategory inclusions, the top squares are cartesian, the vertical arrows in the bottom row are cartesian fibrations, and the horizontal arrows in the middle row preserve cartesian edges. In fact, all the vertical maps in the bottom row are also cocartesian fibrations and the horizontal arrows in the middle row also preserve cocartesian edges, see Corollary~\refone{corollary:cocartesian} and Proposition~\refone{proposition:left-kan-bilinear-linear}. We also recall that the \(\Ct\)-action on \(\Catx \times \Catx\) is given by \((\A,\cB) \mapsto (\cB\op,\A\op)\) and the composite \(\Catx \to \Catx \times \Catx\) along the bottom and right sides both equivalent to the diagonal map.

A useful feature of the diagram~\eqrefone{equation:cross-effect-catp} is that all arrows in it admit both left and right adjoints, and in particular all arrows preserve all limits and colimits. For the vertical functors on the bottom row, these are all cartesian and cocartesian fibrations and their fibres admit zero objects, thus the zero section gives a two sided adjoint in all three cases. The left and right adjoints of the horizontal functors in~\eqrefone{equation:cross-effect-catp} will be studied in this section, with special interest given to the resulting adjoints of the composite arrow \(\Catp \to \Catx\) on the top row, which are both given by the construction \(\C \mapsto \Hyp(\C)\) of \S\refone{subsection:hyp-and-sym-poincare-objects}. The left and right adjoint to the top left vertical inclusion \(\Catp \hrar \Cath\) will be produced in \S\refone{subsection:thom} below using the pairings construction. Left and right adjoints to the top middle vertical inclusion then follow by formal considerations, while left and right adjoints to the top right inclusion can be constructed in a similar manner, see Remarks~\refone{remark:first-adj-bilinear} and~\refone{remark:second-adj-bilinear}.

The remainder of this section is dedicated to the construction of left and right adjoints to the horizontal functors in~\eqrefone{equation:cross-effect-catp}. We begin with the top left square:
\begin{proposition}
\label{proposition:adjunction-for-hermitian-categories}%
The functor
\[
\Cath \to \Catsb
\]
admits fully-faithful left and right adjoints, given by sending \((\C,\Bil)\) to \((\C,\QF^{\qdr}_{\Bil})\) and \((\C,\QF^{\sym}_{\Bil})\), respectively. Furthermore, the units and counits of these adjunctions project to equivalences in \(\Catx\).
\end{proposition}
\begin{proof}
We have a diagram
\[
\begin{tikzcd}
\Cath \ar[rr,"f"] \ar[rd,"p"] && \Catsb \ar[ld,"q"'] \\
& \Catx &
\end{tikzcd}
\]
where \(p\) and \(q\) are cartesian fibrations and the forgetful functor \(f\) preserves cartesian edges. The functor \(f\) has fibrewise left and right adjoints by Corollary~\refone{corollary:hom-universal}, and so
by \cite[Proposition 7.3.2.6]{HA} \(f\) admits a left adjoint whose unit transformation is sent to an equivalence in \(\Catx\).
The counit of the adjunction is given by the fibrewise counit and thus is an equivalence by Corollary~\refone{corollary:hom-universal}, according to which the fibrewise left adjoint is fully faithful. On the other hand, by Corollary~\refone{corollary:cocartesian} and Remark~\refone{remark:catb-has-all-limits-and-colimits}
the maps \(p,q\) are also cocartesian fibrations, and by Proposition~\refone{proposition:left-kan-bilinear-linear} the map \(f\) preserves cocartesian edges.
Applying the dual of~\cite[Proposition 7.3.2.6]{HA} we now conclude that \(f\) admits a right adjoint whose associated unit is mapped to an equivalence in \(\Catx\) and whose counit is an equivalence in \(\Catsb\).
\end{proof}

\begin{proposition}
\label{proposition:adjunction-for-poincare-categories}%
The fully-faithful adjoints of Proposition~\refone{proposition:adjunction-for-hermitian-categories} map \(\Catps\) to \(\Catp\), and yield fully-faithful left and right adjoints to the forgetful functor \(\Catp \to \Catps\).
\end{proposition}
\begin{proof}
Since the left and right adjoints constructed in Proposition~\refone{proposition:adjunction-for-hermitian-categories} are fully-faithful and since the condition for a hermitian \(\infty\)-category of being non-degenerate is defined via the non-degeneracy of the underlying symmetric bilinear form we see that these adjoints send non-degenerate symmetric \(\infty\)-categories to non-degenerate hermitian \(\infty\)-categories and similarly perfect symmetric \(\infty\)-categories to Poincaré \(\infty\)-categories. For the same reason these adjoints also send duality-preserving symmetric functors to duality preserving functors. To see that these now yield left and right adjoints to \(\Catp \to \Catps\) it is enough to check that all the units and counits are contained in the respective subcategories. Now on the side of \(\Catps\) these units and counits are equivalences, and so by triangle identities these units and counits in \(\Catp\) map to equivalences in \(\Catps\). This implies that they are duality preserving, as desired.
\end{proof}

\begin{definition}
\label{definition:symmetric-quadratic-categories}%
Given a perfect bilinear functor \(\Bil \in \Funs(\C)\), the images of \((\C,\Bil)\) under the fully-faithful left adjoint and right adjoint of Proposition~\refone{proposition:adjunction-for-poincare-categories} are, as in the case of Proposition~\refone{proposition:adjunction-for-hermitian-categories}, given by \((\C,\QF^{\qdr}_{\Bil})\) and \((\C,\QF^{\sym}_{\Bil})\), respectively. We will refer to Poincaré \(\infty\)-categories of this form as \defi{quadratic} and \defi{symmetric} Poincaré \(\infty\)-categories, respectively.
\end{definition}

Taking now the top external rectangle in~\eqrefone{equation:cross-effect-catp} we obtain a square
\begin{equation}
\label{equation:forgetful-square}%
\begin{tikzcd}
\Catp \ar[r]\ar[d] & \Catx \ar[d] \\
\Cath \ar[r] & \Catb
\end{tikzcd}
\end{equation}
in which the right vertical arrow is given by
\begin{equation}
\label{equation:catx-catb}%
\begin{tikzcd}
[row sep=1ex]
\Catx \ar[r,"{\simeq}"] & \Catpb \ar[r,hook] & \Catb \\
\A \ar[rr,mapsto] \ar[u,phantom,"{\rotatebox{90}{$\in$}}"] && (\A,\A,m_{\A}) \ar[u,phantom,"{\rotatebox{90}{$\in$}}"].
\end{tikzcd}
\end{equation}

\begin{proposition}
\label{proposition:generalized-hyp}%
In the square~\eqrefone{equation:forgetful-square}
both horizontal arrows admit left and right adjoints, both compatible with the vertical subcategory inclusions. In addition, both adjoints of the bottom horizontal arrow are equivalent and given by the formula
\((\A,\cB,\Bil) \mapsto (\A \times \cB\op,\Bil)\).
\end{proposition}

Pre-composing the formula of Proposition~\refone{proposition:generalized-hyp} with the functor~\eqrefone{equation:catx-catb} we conclude

\begin{corollary}
\label{corollary:hyp-is-adjoint}%
The forgetful functor \(\rU\colon \Catp \to \Catx\) admits both a left and a right adjoint. The two are equivalent and given by the association \(\A \mapsto (\A \times \A\op,m_{\A}) = \Hyp(\A)\).
\end{corollary}

\begin{remark}
\label{remark:units-and-counits}%
The unit exhibiting \(\Hyp\) as right adjoint to \(\rU\) and the counit exhibiting \(\Hyp\) as left adjoint to \(\rU\) are given respectively by the Poincaré functors
\[
\Hyp(\C) \xrightarrow{\hyp} (\C,\QF) \xrightarrow{\fgt} \Hyp(\C)
\]
of~\eqrefone{equation:hyp-forget}.
The counit \(\rU \Hyp(\C) = \C \oplus \C\op \to \C\) exhibiting \(\Hyp\) as right adjoint to \(\rU\) is then given by the projection on the first fact and the unit \(\C \to \rU\Hyp(\C) = \C \oplus \C\op\) exhibiting \(\Hyp\) as left adjoint to \(\rU\) is given by the inclusion into the first summand. Here we are using the direct sum notation keeping in mind that \(\Catx\) is semi-additive, see Proposition~\refone{proposition:catp-pre-add}.
\end{remark}

\begin{remark}
\label{remark:hyp-monoidal}%
By~\cite[Corollary 7.3.2.7]{HA} the functor \(\Hyp\) inherits a lax symmetric monoidal structure by virtue of being right adjoint to the symmetric monoidal functor \(\rU\colon \Catp \to \Catx\) (see Theorem~\refone{theorem:tensorpoincare}). In particular, \(\rU \dashv \Hyp\) is a symmetric monoidal adjunction. Applying the same argument to the symmetric monoidal functor \(\rU\op\colon (\Catp)\op \to (\Catx)\op\) we get that the adjunction \(\rU\op \dashv \Hyp\op\) (opposite to \(\Hyp \dashv \rU\)) is symmetric monoidal and \(\Hyp\) also carries an oplax symmetric monoidal structure.
\end{remark}

\begin{remark}
\label{remark:functor-category-to-hyp}%
The symmetric monoidal structure on the forgetful functor \(\rU\colon \Catp \to \Catx\) yields for every \((\C,\QF) \in \Catp\) a commutative square
\[
\begin{tikzcd}
\Catp \ar[r,"{\rU}"] \ar[d,"{(\C,\QF) \otimes (-)}"'] & \Catx \ar[d,"{\C \otimes (-)}"] \\
\Catp \ar[r,"{\rU}"] & \Catx \ .
\end{tikzcd}
\]
Passing to right adjoints using Corollary~\refone{corollary:hyp-is-adjoint} and Corollary~\refone{corollary:catp-closed} we obtain a commutative square
\[
\begin{tikzcd}
\Catp\ar[d,"{\Funx((\C,\QF),-)}"'] & \Catx \ar[d,"{\Funx(\C,-)}"] \ar[l,"{\Hyp}"] \\
\Catp & \Catx \ar[l,"{\Hyp}"].
\end{tikzcd}
\]
In particular, we have a natural equivalence of Poincaré \(\infty\)-categories
\[
\Funx((\C,\QF),\Hyp(\Ctwo)) \simeq \Hyp\Funx(\C,\Ctwo).
\]
\end{remark}

\begin{proof}[Proof of Proposition~\refone{proposition:generalized-hyp}]
We first observe that the square~\eqrefone{equation:forgetful-square} is also the external top rectangle in~\eqrefone{equation:cross-effect-catp}. Propositions~\refone{proposition:adjunction-for-hermitian-categories} and~\refone{proposition:adjunction-for-poincare-categories} provide left and right adjoints to the horizontal arrows in the top left square of~\eqrefone{equation:cross-effect-catp}, which are furthermore compatible with the vertical subcategory inclusions.
In the same diagram, the horizontal arrows in the top middle square are equivalences, and the horizontal arrows in the top right square have left and right adjoints as well. Indeed, in both cases these are the natural maps from \(\Ct\)-fixed points to underlying objects, which admit left and right adjoint since \(\Catx\) and \(\Catb\) admits products and coproducts (see Proposition~\refone{proposition:catx-has-limits-and-colimits} and Remark~\refone{remark:catb-has-all-limits-and-colimits}).
In addition, since \(\Catx\) is semiadditive (Proposition~\refone{proposition:catp-pre-add}) these products and coproducts coincide,
and we get that the left and right adjoints of \((\Catx)^{\hC} \to \Catx\) are both give by the formula \(\C \mapsto \C \times \C\op\). This implies that the left and right adjoints for the horizontal maps in the top right square in~\eqrefone{equation:cross-effect-catp} are compatible with the vertical replete subcategory inclusions, and hence we may conclude that the horizontal arrows in~\eqrefone{equation:forgetful-square}
admit left and right adjoints, both compatible with the vertical subcategory inclusions.

It is left to obtain the desired explicit formula.
First by semi-additivity, the left and right adjoints of \(\Catsb \to \Catb\) are
both given by the formula
\[
(\A,\cB,\Bil) \mapsto (\A \times \cB\op,\Bil_{\mathrm{sym}})
\]
where \(\Bil_{\mathrm{sym}}\in\Funs(\A \times \cB\op)\) is the symmetrisation of \(\Bil\), given by
\[
\Bil_{\mathrm{sym}}((a,b),(a',b')) = \Bil(a,b') \oplus \Bil(a',b) .
\]
We note that this symmetric bilinear form is both induced and coinduced from \(\Bil\). In particular, \(\QF^{\sym}_{\Bil_{\mathrm{sym}}}\) and \(\QF^{\qdr}_{\Bil_{\mathrm{sym}}}\)
both canonically identify with \(\Bil\) itself, when the latter is viewed as a single entry functor on \(\A\op \times \cB\). By Propositions~\refone{proposition:adjunction-for-hermitian-categories} the left and right adjoints of the bottom horizontal map in~\eqrefone{equation:forgetful-square} are consequently both given by the formula
\((\A,\cB,\Bil) \mapsto (\A \times \cB\op,\Bil)\).
\end{proof}

\subsection{The categorical Thom isomorphism}
\label{subsection:thom}%

In~\S\refone{subsection:algebraic-thom} we described the \emph{algebraic Thom construction}, an operation introduced by Ranicki which allows one to identify the notion of a metabolic Poincaré object \((\x,\qone)\) in \((\C,\QF)\) equipped with a prescribed Lagrangian \((\cob \to \x,\eta)\), with the data of a hermitian object \((\z,r)\) in \(\C\) with respect to the shifted Poincaré structure \(\QF\qshift{-1}\). In this section we will see that this procedure naturally fits in a more general perspective. We will begin by refining the construction of pairing \(\infty\)-categories described in \S\refone{subsection:bilinear-and-pairings} above to the context of hermitian \(\infty\)-categories. This will result in a construction which takes a hermitian \(\infty\)-category \((\C,\QF)\) and produces a Poincaré \(\infty\)-category \(\Pairings(\C,\QF)\). When \((\C,\QF)\) is Poincaré this construction reproduces that of the \emph{arrow category} \(\Ar(\C,\QF)\) described in~\S\refone{subsection:algebraic-thom}. We will then prove that Poincaré objects in \(\Pairings(\C,\QF)\) correspond to hermitian objects in \((\C,\QF)\), yielding in particular a proof of Proposition~\refone{proposition:thom} upon taking \((\C,\QF)\) to be Poincaré.

The Poincaré \(\infty\)-categories of the form \(\Pairings(\C,\QF)\) can be considered as a categorified form of the notion of a metabolic Poincaré object: they contain a stable full subcategory on which the Poincaré structure vanishes and which is equivalent to its own orthogonal complement, a property which can be considered as a categorical analogue of the notion of a Lagrangian. From this point of view, we may consider the association \((\C,\QF) \mapsto \Pairings(\C,\QF)\) as a categorified form of the Thom construction, taking a hermitian \(\infty\)-category and producing a Poincaré \(\infty\)-category with a canonical choice of Lagrangian. We will show that this process is reversible: given a Poincaré \(\infty\)-category \((\D,\QFD)\) with a Lagrangian, one can reconstruct a hermitian \(\infty\)-category \((\C,\QF)\) such that \(\Pairings(\C,\QF) \simeq (\D,\QFD)\). We consider this as a categorical form of the Thom isomorphism. Relying on these results we will then use the pairing construction in order to produce both a left and a right adjoint to the forgetful functor \(\Catp \to \Cath\). This adjunction ties together all the above results in a conceptual manner, and at the same time is quite useful in practice. In particular, the results of this section we will be used in subsequent instalments of this project, for example in setting up the theory of \emph{algebraic surgery} in \papertwo, and in proving that \(\Catp\) is compactly generated in \paperfour.

We now proceed to introduce the main construction of the current section:

\begin{construction}
\label{construction:hermitian-pairings}%
Given a stable \(\infty\)-category \(\C\), any bilinear functor \(\Bil\colon \C\op \times \C\op \to \Spa\) determines a left biexact correspondence on the pair \((\C,\C\op)\), given by \((\x,\y) \mapsto \Om^{\infty}\Bil(\x,\y)\). We will then denote by
\[
\Pairings(\C,\Bil) := \Pairings(\C,\C\op,\Bil) \in \Catx
\]
the associated \(\infty\)-category of pairings. Given a hermitian structure \(\QF\colon \C\op \to \Spa\) on \(\C\), we define an associated hermitian structure \(\QF_{\pair}\) on \(\Pairings(\C,\Bil_{\QF})\) via the pullback square
\begin{equation}
\label{equation:defining-square-swap}%
\begin{tikzcd}
\QF_{\pair}(\x,\y,\beta) \ar[r] \ar[d] & \QF(\x) \ar[d] \\
\map_{\C}(\x,\y) \ar[r] & \Bil(\x,\x)
\end{tikzcd}
\end{equation}
where the bottom horizontal map is given by the association \([f\colon \x \to \y] \mapsto f^*\beta \in \Bil(\x,\x)\), canonically extended from mapping spaces to mapping spectra. We denote the resulting hermitian \(\infty\)-category by \(\Pairings(\C,\QF) := (\Pairings(\C,\Bil_{\QF}),\QF_{\pair})\).
\end{construction}

\begin{lemma}
\label{lemma:pair-structure-perfect}%
The hermitian structure \(\QF_{\pair}\) is Poincaré. Its bilinear part sits in the pullback square of symmetric bilinear forms
\begin{equation}
\label{equation:square-bil-swap}%
\begin{tikzcd}
\Bil_{\pair}((\x,\y,\beta),(\x',\y',\beta')) \ar[r] \ar[d] & \Bil(\x,\x') \ar[d,"{(\tau_{\x,\x'},\id)}"] \\
\map_{\C}(\x,\y')\oplus \map_{\C}(\x',\y) \ar[r] & \Bil(\x',\x)\oplus \Bil(\x,\x') .
\end{tikzcd}
\end{equation}
and its associated duality coincides with the duality
\[
\Dual_{\pair}(\x,\y,\beta) = \big(\y,\x,\sig_{\x,\y}(\beta)\big) ,
\]
of Construction~\refone{construction:duality-pairing}, which encodes the symmetric structure of \(\Bil\).
\end{lemma}
\begin{proof}
Substituting in~\eqrefone{equation:defining-square-swap} the direct sum of \((\x,\y,\beta),(\x',\y',\beta') \in \Pairings(\C,\Bil_{\QF})\)
we obtain the square
\[
\begin{tikzcd}
\QF_{\pair}(\x\oplus \x',\y \oplus \y',\beta\oplus \beta') \ar[r] \ar[d] & \QF(\x \oplus \x') \ar[d] \\
\map_{\C}(\x\oplus\x',\y\oplus\y') \ar[r] & \Bil(\x\oplus\x',\x\oplus\x')
\end{tikzcd},
\]
which yields the square~\eqrefone{equation:square-bil-swap} upon passing to bireduced replacements.
On the other hand, working backwards from the required duality, we note that \(\Dual_{\pair}\) is part of a \(\Lam^2_0\)-duality of the associated bifibration
\[
\C \leftarrow \Pairings(\C,\Bil_{\QF}) \rightarrow \C\op.
\]
Applying Proposition~\refone{proposition:pairings-mapping-spaces} and Remark~\refone{remark:mapping-formula-duality-preserving} to this bifibration we obtain a cartesian square of symmetric bilinear forms
\[
\begin{tikzcd}
\map_{\Pairings(\C,\Bil_{\QF})}((\x,\y,\beta),\Dual_{\pair}(\x',\y',\beta')) \ar[r] \ar[d] & \Bil(\x,\x') \ar[d,"{(\tau_{\x,\x'},\id)}"] \\
\map_{\C}(\x,\y')\oplus \map_{\C}(\x',\y) \ar[r] & \Bil(\x',\x)\oplus \Bil(\x,\x') ,
\end{tikzcd}
\]
where we are again silently identifying exact functors valued in spectra with left exact functors valued in spaces. Comparing this square with~\eqrefone{equation:square-bil-swap}
we thus get that \(\Bil_{\pair}\) is perfect with duality
\(\Dual_{\pair}\).
\end{proof}

\begin{examples}
\label{examples:typical-pairings}%
\ %
\begin{enumerate}
\item
\label{item:zero-pairing-hyp}%
For a stable \(\infty\)-category \(\C\) equipped with the zero hermitian structure, the Poincaré \(\infty\)-category \(\Pairings(\C,0)\) is naturally equivalent to \(\Hyp(\C)\).
\item
\label{item:correspondence-pairing}%
If the hermitian \(\infty\)-category \((\C,\QF)\) is Poincaré with associated duality \(\Dual\), then the correspondence \(\Om^{\infty}\Bil_{\QF}\) on the pair \((\C,\C\op)\) is equivalent to the correspondence \(\Map_{\C}(-,-)\) on the pair \((\C,\C)\) via the
\[
(\id,\Dual)\colon \C\op \times \C\op \xrightarrow{\simeq} \C\op \times \C,
\]
and so \(\Pairings(\C,\Bil_{\QF})\) is naturally equivalent to the arrow category \(\Ar(\C)\) of Definition~\refone{definition:arrow-cat}. Under this equivalence, the Poincaré structure \(\QF_{\pair}\) directly translates to the Poincaré structure
\[
\QF_{\arr}(\x \to \y) = \QF(\x) \times_{\Bil_{\QF}(\x,\x)}\Bil_{\QF}(\x,\y),
\]
and so we obtain a natural equivalence \(\Pairings(\C,\QF) \simeq \Ar(\C,\QF)\).
\item
\label{item:pairing-arr-met}%
Combining the previous example with Lemma~\refone{lemma:arrow-is-met} we obtain a natural identification of Poincaré \(\infty\)-categories
\(\Pairings(\C,\QF\qshift{-1}) \simeq \Ar(\C,\QF\qshift{-1}) \simeq \Met(\C,\QF)\) whenever \((\C,\QF)\) is Poincaré.
\end{enumerate}
\end{examples}

\begin{remark}
\label{remark:duality-arrows}%
Combining Example~\refone{examples:typical-pairings}(ii), Example~\refone{example:swap-duality-arrows} and Lemma~\refone{lemma:pair-structure-perfect}, we obtain that for \((\C,\QF)\) Poincaré, the underlying duality \(\Dual_{\arr}\) of \(\Ar(\C,\QF)\) is equivalent to the one induced on \(\Ar(\C) \simeq \Fun(\Del^1,\C)\) by the duality \(\Dual_{\QF}\) on \(\C\) and the canonical duality of \(\Del^1\).
\end{remark}

By construction, the underlying stable \(\infty\)-category of \(\Pairings(\C,\QF)\) sits in a bifibration
\begin{equation}
\label{equation:bifibration-pairing}%
\begin{tikzcd}
[row sep=1ex]
\x\ar[d,phantom,"{\rotatebox{270}{$\in$}}"]  & (\x,\y,\beta) \ar[l,mapsto] \ar[r,mapsto] \ar[d,phantom,"{\rotatebox{270}{$\in$}}"]  & \y \ar[d,phantom,"{\rotatebox{270}{$\in$}}"]  \\
\C & \Pairings(\C,\Bil_{\QF}) \ar[l,"q"'] \ar[r,"p"] & \C\op
\end{tikzcd}
\end{equation}
so that \(q\) is a cartesian fibration and \(p\) is a cocartesian fibration. The fact that \(\C\) is pointed and \(\Bil_{\QF}\) is bireduced implies that these fibrations have fully-faithful adjoints
\[
\begin{tikzcd}
[row sep=1ex]
\x\ar[d,phantom,"{\rotatebox{270}{$\in$}}"] \ar[r,mapsto] & (\x,0,0)\ar[d,phantom,"{\rotatebox{270}{$\in$}}"]  &  \\
\C \ar[r,"j"] & \Pairings(\C,\Bil_{\QF}) & \C\op\ar[l,"i"'] \\
 & (0,\y,0) \ar[u,phantom,"{\rotatebox{90}{$\in$}}"]  & \y \ar[u,phantom,"{\rotatebox{90}{$\in$}}"] \ar[l,mapsto]
\end{tikzcd}
\]
More precisely, \(q\) has a fully-faithful right adjoint \(j\colon \C \to \Pairings(\C,\QF)\) sending \(\x\) to \((\x,0,0)\). Indeed, the canonical arrows \((\x,y,\beta) \to (\x,0,0)\) in \(\Pairings(\C,\QF)\) induce an equivalence on mapping spaces into any triple of the form \((\x',0,0)\), and thus assemble to form a unit exhibiting \(j\) as right adjoint to \(q\). Similarly, the collection of arrows \((0,\y,0) \to (\x,\y,\beta)\) assemble to form a counit exhibiting the functor \(i\colon \C\op \to \Pairings(\C,\QF)\) sending \(\y\) to \((0,\y,0)\) as left adjoint to \(p\), and we observe that the image of \(i\) coincides with the kernel of \(q\) and the image of \(j\) with the kernel of \(p\).

We now observe that \(q\colon \Pairings(\C,\QF) \to \C\) naturally extends to a hermitian functor
\begin{equation}
\label{equation:counit-hermitian}%
(q,\eta)\colon \Pairings(\C,\QF) \to (\C,\QF) \quad\quad\quad
\quad q(\x,\y,\beta) = \x \quad,\quad \eta\colon \QF_{\pair}(\x,\y,\beta) \to \QF(\x)
\end{equation}
with \(\eta\colon \QF_{\pair} \Rightarrow q^*\QF\) given by the natural projection furnished directly from the definition of \(\QF_{\pair}\). In generalisation of the algebraic Thom isomorphism we then have:
\begin{proposition}%
\label{proposition:algebraic-thom}%
For every hermitian \(\infty\)-category \((\C,\QF)\) the composite map
\[
\Poinc(\Pairings(\C,\QF)) \to \spsforms(\Pairings(\C,\QF)) \xrightarrow{(q,\eta)_*} \spsforms(\C,\QF)
\]
is an equivalence. In particular, Poincaré objects in \(\Pairings(\C,\QF)\) classify hermitian objects in \((\C,\QF)\).
\end{proposition}

\begin{remark}
Applied in the case where \((\C,\QF)\) is Poincaré, Proposition~\refone{proposition:algebraic-thom} reduces to the statement of Proposition~\refone{proposition:thom} via the identification \(\Pairings(\C,\QF) \simeq \Ar(\C,\QF)\) of Examples~\refone{examples:typical-pairings}.
\end{remark}

Proposition~\refone{proposition:algebraic-thom} is a direct consequence of the following lemma:
\begin{lemma}
\label{lemma:final-in-fiber}%
Let \((\C,\QF)\) be a hermitian \(\infty\)-category. Then the functor
\[
(q,\eta)_*\colon \catforms(\Pairings(\C,\QF)) \to \catforms(\C,\QF)
\]
induced by~\eqrefone{equation:counit-hermitian} is a cartesian fibration whose fibres admit final objects. In addition, a hermitian object \(((\x,\y,\beta),q)\) in \(\Pairings(\C,\QF)\) is final in its fibre over \(\catforms(\C,\QF)\) if and only if it is Poincaré.
\end{lemma}

Given Lemma~\refone{lemma:final-in-fiber}, the proof of Proposition~\refone{proposition:algebraic-thom} is immediate:
\begin{proof}[Proof of Proposition~\refone{proposition:algebraic-thom}]
By Lemma~\refone{lemma:final-in-fiber} the homotopy fibres of
\[
\Poinc(\Pairings(\C,\QF)) \to \spsforms(\C,\QF)
\]
are contractible, and so the desired result follows.
\end{proof}

\begin{proof}[Proof of Lemma~\refone{lemma:final-in-fiber}]
Consider the commutative square
\[
\begin{tikzcd}
\catforms(\Pairings(\C,\QF)) \ar[r] \ar[d] \ar[dr,dashed] & \catforms(\C,\QF) \ar[d] \\
\Pairings(\C,\QF) \ar[r] & \C
\end{tikzcd}
\]
where the vertical arrows are the defining right fibrations of \(\catforms(-)\). Since the bottom horizontal map is a cartesian fibration (by the construction of \(\Pairings(\C,\QF)\) as a bifibration) we obtain that the composite dotted map is a cartesian fibration, and hence we can view the top horizontal map is a map of cartesian fibrations over \(\C\) whose target is a right fibration. Such a map is automatically a cartesian fibration (up to equivalence), which gives us the first claim of the lemma. We now verify that its fibres contain final objects. Let
\[
\X := \int^{(\x,\y,\alp) \in \Pairings(\C,\Bil_{\QF})}\Map_{\C}(\x,\y) \to \Pairings(\C,\Bil_{\QF})
\]
be the right fibration classifying the contravariant functor \((\x,\y,\alp) \mapsto \Map_{\C}(\x,\y)\). The defining fibre square~\eqrefone{equation:defining-square-swap} then determines a commutative square
\begin{equation}
\label{equation:defining-rectangle}%
\begin{tikzcd}
\catforms(\Pairings(\C,\QF)) \ar[r] \ar[d] & \catforms(\C,\QF) \ar[d] \\
\X \ar[r] \ar[d] & \int^{\x \in \C}\Om^{\infty}\Bil_{\QF}(\x,\x) \ar[d] \\
\Pairings(\C,\Bil_{\QF}) \ar[r] & \C
\end{tikzcd}
\end{equation}
of \(\infty\)-categories in which (using the same argument as above) the vertical maps are right fibrations and the horizontal maps are cartesian fibrations. In addition, the top square in~\eqrefone{equation:defining-rectangle} is cartesian; indeed, base changing the right fibrations on the right from \(\C\) to \(\Pairings(\C,\Bil_{\QF})\) yields a fibre square of right fibrations which is the straightening of the defining square~\eqrefone{equation:defining-square-swap} (after taking \(\Om^{\infty}\)). It will hence suffice to show that the fibres of the middle horizontal cartesian fibration have final objects. Fix an object \((\x,\beta) \in \int^{\x \in \C}\Om^{\infty}\Bil_{\QF}(\x,\x)\). Then the fibre \(\X_{(\x,\beta)}\) of the middle horizontal map sits in a right fibration
\begin{equation}
\label{equation:fiber-x-beta}%
\X_{(\x,\beta)} \to \Pairings(\C,\Bil_{\QF})_{\x} = \int_{\y \in \C\op}\Om^{\infty}\Bil_{\QF}(\x,\y) ,
\end{equation}
where the middle term stands for the fibre of the bottom horizontal map over \(\x \in \C\). Using the pullback formula~\eqrefone{equation:mapping-space-formula} for the mapping spaces in \(\Pairings(\C,\Bil_{\QF})\) we now calculate
\[
\begin{split}
\X_{(\x,\beta)}  & = \int_{(\y,\alp) \in \Pairings(\C,\Bil_{\QF})_{\x}}\Map_{\C}(\x,\y)\times_{\Bil_{\QF}(\x,\x)}\{\beta\} \simeq
\int_{(\y,\alp) \in \Pairings(\C,\Bil_{\QF})_{\x}}\Map_{\C\op}(\y,\x)\times_{\Bil_{\QF}(\x,\x)}\{\beta\} \\
 & \simeq \int_{(\y,\alp) \in \Pairings(\C,\Bil_{\QF})_{\x}}\Map_{\Pairings(\C,\Bil_{\QF})_{\x}}((\x,\y,\alp),(\x,\x,\beta)) \simeq (\Pairings(\C,\Bil_{\QF})_{\x})_{/(\x,\x,\beta)},
\end{split}
\]
from which we see that the right fibration~\eqrefone{equation:fiber-x-beta} is visibly represented by \((\x,\x,\beta)\), so that
\(\X_{(\x,\beta)}\) has a final object.

To finish the proof we now need to verify that a hermitian object \(((\x,\y,\beta),q)\) in \(\Pairings(\C,\QF)\) is Poincaré if and only if it is final in the fibre. The top square in~\eqrefone{equation:defining-rectangle} being cartesian it will suffice to show that \(((\x,\y,\beta),q)\) is Poincaré if and only if its image in \(\X\) is final in its fibre over \(\int^{\x \in \C}\Om^{\infty}\Bil_{\QF}(\x,\x)\). Indeed, a hermitian form \(q\) on \((\x,\y,\beta)\) determines a self dual map \(q_{\sharp}\colon(\x,\y,\beta) \to \Dual_{\pair}(\x,\y,\beta) = (\y,\x,\beta)\) with components \(f\colon \x \to \y\) and \(\y \leftarrow \x\colon g\) (the latter considered as a map from \(\y\) to \(\x\) in \(\C\op\)). The form \(q\) is then Poincaré if and only if \(f\) and \(g\) are equivalences. But since \(q_{\sharp}\) is self-dual the components \(f\) and \(g\) are homotopic to each other. We then get that \(q\) is Poincaré if and only if the map \(g\) is an equivalence. But \(g\) (or \(f\)) is exactly the image of \(q\) in \(\Map_{\C}(\x,\y)\) via the vertical map in the defining square~\eqrefone{equation:defining-square-swap}, and so the image of \(((\x,\y,\beta),q)\) in \(\X\) is \(((\x,\y,\beta),g)\). The latter object lies over \((\x,g^*\beta) \in \int^{\x \in \C}\Om^{\infty}\Bil_{\QF}(\x,\x)\) and corresponds to the object
\begin{equation}
\label{equation:image-of-q}%
[(\x,\y,\beta) \xrightarrow{(\id,g)} (\x,\x,g^*\beta)] \in \X_{(\x,g^*\beta)} \simeq (\Pairings(\C,\Bil_{\QF})_{\x})_{/(\x,\x,g^*\beta)}
\end{equation}
We may then conclude that \(((\x,\y,\beta),q)\) is Poincaré if and only if \(g\) is an equivalence, and so if and only if~\eqrefone{equation:image-of-q} is final, as desired.
\end{proof}

The following almost immediate corollary of Lemma~\refone{lemma:final-in-fiber} will also be useful for us later:
\begin{corollary}
\label{corolary:left-kan}%
The natural transformation
\[
\eta\colon \QF_{\pair} \Rightarrow q^*\QF
\]
exhibits \(\QF\colon \C\op \to \Spa\) as the left Kan extension of \(\QF_{\pair}\) along \(q\op\colon \Pairings(\C,\Bil_{\QF})\op \to \C\op\).
\end{corollary}
\begin{proof}
By Proposition~\refone{proposition:left-kan-bilinear-linear} and Lemma~\refone{lemma:coconnective-quadratic-functors}
it will suffice to show that \(\Om^{\infty}\eta\) exhibits \(\Om^{\infty}\QF\) as the left Kan extension of \(\Om^{\infty}\QF_{\pair}\). Since \(q\op\) is a cocartesian fibration left Kan extensions along \(q\op\) are calculated by colimits along the fibres. In particular, we need to show that for every \(\x \in \C\) the map
\[
\displaystyle\mathop{\colim}_{(\x,\y,\alp) \in \Pairings(\C,\Bil_{\QF})\op_x}\Om^{\infty}\QF_{\pair}(\x,\y,\alp) \to \Om^{\infty}\QF(\x),
\]
induced by \(\eta\), is an equivalence of spaces. Since colimits in spaces are universal it will suffice to show that for every point \(\beta \in \Om^{\infty}\QF(\x)\) the space
\[
\displaystyle\mathop{\colim}_{(\x,\y,\alp) \in \Pairings(\C,\Bil_{\QF})\op_x}\Om^{\infty}\QF_{\pair}(\x,\y,\alp) \times_{\QF(\x)} \{\beta\}
\]
is contractible. Indeed, this space can in turn be identified with the geometric realization of the fibre of \(\catforms(\Pairings(\C,\QF))\op \to \catforms(\C,\QF)\op\) over \((\x,\beta)\), and the latter has an initial object by Lemma~\refone{lemma:final-in-fiber}, so its realization is contractible.
\end{proof}

We now take a closer look at the Poincaré \(\infty\)-categories of the form \(\Pairings(\C,\QF)\). We wish to make the argument that they constitute a categorical analogue of the notion of a metabolic Poincaré object, making the passage from \((\C,\QF)\) to \(\Pairings(\C,\QF)\) an analogue of the algebraic Thom construction.
To identify further key properties we introduce the following piece of notation:

\begin{definition}
Let \((\D,\QFD)\) be a Poincaré \(\infty\)-category and \(\Lag \subseteq \D\) a full subcategory. We will denote by \(\Lag^{\perp} \subseteq \D\) the full subcategory spanned by the objects \(\y \in \D\) such that \(\Bil_{\QFD}(\x,\y) = 0\) for every \(\x \in \Lag\). We will refer to \(\Lag^{\perp}\) as the \defi{orthogonal complement} of \(\Lag\).
\end{definition}

Using the notion of orthogonal complements, we may identify the following additional properties held by the full subcategory inclusion \(i\colon \C\op \hrar \Pairings(\C,\Bil_{\QF})\):
\begin{enumerate}
\item
\label{item:lag-1}%
The restriction of the quadratic functor \(\QF_{\pair}\) to \(\C\op\) vanishes.
\item
\label{item:lag-2}%
The inclusion \(i\C\op \subseteq (i\C\op)^{\perp}\), furnished by~\refoneitem{item:lag-1} above, is an equivalence.
\item
\label{item:lag-3}%
\(i\) admits a right adjoint \(p\colon \Pairings(\C,\Bil_{\QF}) \to \C\op\) given by \((\x,\y,\beta) \mapsto \y\) (see~\eqrefone{equation:bifibration-pairing} and the discussion below it).
\end{enumerate}
The validity of~\refoneitem{item:lag-1} above is evident from the fibre square~\eqrefone{equation:defining-square-swap} defining \(\QF_{\pair}\). To see that~\refoneitem{item:lag-2} holds, note that
\[
\Bil_{\pair}(i(\z),(\x,\y,\beta)) = \map((0,\z,0),(\y,\x,\sig_{x,y}(\beta))) = \map_{\C\op}(\z,\x) = \map_{\C}(\x,\z),
\]
and hence \((\x,\y,\beta) \in (i\C\op)^{\perp}\) if and only if \(\x=0\), i.e., if and only if \((\x,\y,\beta) \in i\C\op\). This motivates the following definition:

\begin{definition}
\label{definition:lagrangian}%
Let \((\D,\QFD)\) be a Poincaré \(\infty\)-category and \(\Lag \subseteq \D\) a full subcategory.
We will say that \(\Lag\) is a \defi{Lagrangian} in \(\D\) if it satisfies Properties~\refoneitem{item:lag-1},~\refoneitem{item:lag-2} and~\refoneitem{item:lag-3} above. In other words, if \(\QF\) vanishes when restricted to \(\Lag\), the orthogonal complement \(\Lag^{\perp}\subseteq \D\) coincides with \(\Lag\) itself and the inclusion \(\Lag \subseteq \D\) admits a right adjoint.
We will say that a Poincaré \(\infty\)-category is \defi{metabolic} if it admits a Lagrangian.
\end{definition}

In particular, the Poincaré \(\infty\)-category \(\Pairings(\C,\QF)\) contains \(\C\op\) as a Lagrangian. We now claim that this property completely characterises Poincaré \(\infty\)-categories of the form \(\Pairings(\C,\QF)\). To see this, let \((\D,\QFD)\) be a Poincaré \(\infty\)-category with underlying duality \(\Dual\), and let \(i\colon \Lag \hrar \D\) be a Lagrangian with right adjoint \(p\colon \D \to \Lag\). Let \(j\colon \Lag^{\op} \subseteq \D\) be the inclusion sending \(\z\) to \(\Dual(i(\z))\). Then
\[
\im(j) =\Dual(i\Lag) = \Dual((i\Lag)^{\perp}) = \ker(p),
\]
and \(j\) admits a left adjoint \(q\colon \D \to \Lag\op\) given by the formula \(q(\x) = p(\Dual\x)\).

We then have the following:

\begin{proposition}[Recognition principle for pairing Poincaré categories]
\label{proposition:recognize-poincare}%
Let \((\D,\QFD)\) be a Poincaré \(\infty\)-category admitting a Lagrangian \(i \colon \Lag \hrar \D\) with right adjoint \(p\colon \D \to \Lag\), and let the adjunction \(q\colon \D \adj \Lag\op \colon j\)
be as above. Let \(\Pi = q_!\QFD \in \Funq(\Lag\op)\) be the left Kan extension of \(\QFD\) along \(q\op\colon \D\op \to \Lag\) and \(\eta\colon \QFD \Rightarrow q^*\Pi\) the unit natural transformation.
Then there exists a canonical diagram of hermitian \(\infty\)-categories
\[
\begin{tikzcd}
(\Lag\op,\Pi) \ar[d,equal] & \Pairings(\Lag\op,\Pi) \ar[l] \ar[r] \ar[d,"{\simeq}"] & (\Lag,0)\ar[d,equal] \\
(\Lag\op,\Pi) & (\D,\QFD) \ar[l,"{(q,\eta)}"'] \ar[r,"{(p,0)}"]  & (\Lag,0)
\end{tikzcd}
\]
in which the middle vertical arrow is an equivalence of Poincaré \(\infty\)-categories.
Here the top row consists of the
underlying bifibration~\eqrefone{equation:bifibration-pairing} of \(\Pairings(\Lag\op,\Bil_{\Pi})\) promoted to the level of hermitian \(\infty\)-categories trivially on the right and as in~\eqrefone{equation:counit-hermitian} on the left.
\end{proposition}

\begin{remark}
\label{remark:recognition-bilinear}%
The recognition principle of Proposition~\refone{proposition:recognize-poincare} could also be formulated in the bilinear setting of \S\refone{subsection:bilinear-and-pairings}. In particular, given a stable \(\infty\)-category \(\D\), equivalences of the form \(\D \simeq \Pairings(\A,\cB,b)\) for \((\A,\cB,b) \in \Catb\) correspond to fully-faithful embeddings \(i\colon \cB \hrar \D\) which admit a right adjoint \(p\colon \D \to \cB\), in which case \(\A\) is recovered as the kernel of \(p\), and \(b\) is recovered as the restriction of the correspondence \(\seq_{\D}\colon \D\op \times \D \to \Sps\) of Example~\refone{example:canonical-shift} to \(\A\op \times \cB\). The proof of this claim essentially amounts to the first half of the proof of Proposition~\refone{proposition:recognize-poincare} below.
\end{remark}

\begin{proof}[Proof of Proposition~\refone{proposition:recognize-poincare}]
Let \(\Seq(\D,\Lag) \subseteq \Seq(\D)\) denote the full subcategory spanned by those exact sequences \([\y \to \z \to \x]\) such that \(\y \in i(\Lag)\) and \(\x \in j(\Lag\op)\). We claim that the projection
\[
\Seq(\D,\Lag) \to \D \quad\quad [\y \to \z \to \x] \mapsto \z
\]
is an equivalence. To see this, consider first the map \(\Seq(\D,\Lag) \to \Ar(\D) \times_{\D} j(\Lag\op)\) sending \([\y \to \z \to \x]\) to \([\z \to \x]\). By~\cite[Proposition 4.3.2.15]{HTT} this map is fully-faithful, with essential image the full subcategory of \(\Ar(\D) \times_{\D} j(\Lag\op)\) spanned by those arrows \(\z \to \x\) with \(\x \in j(\Lag\op)\) whose fibre lies in \(i(\Lag)\). Now since \(i(\Lag) = \ker(q)\) the condition that the fibre of \(\z \to \x\) lies in \(i(\Lag)\) is equivalent to the condition that \(q(\z) \to q(\x) \simeq \x\) is an equivalence. Now the projection \(\Ar(\D) \times_{\D} j(\Lag\op) \to \D\) sending \(\z \to \x\) to \(z\) is a cartesian fibration whose fibre over \(\z \in \D\) is equivalent to the comma category \(\D_{z/} \times_{\D} j(\Lag\op)\). This comma category is equivalent by adjunction to \(\Lag\op_{q(z)/}\), and the above condition shows that under this equivalence the full subcategory \(\Seq(\D,\Lag) \subseteq \Ar(\D) \times_{\D} j(\Lag\op)\) consists of exactly those objects which are initial in their fibres. The projection \(\Seq(\D,\Lag) \to \D\) is consequently an equivalence.

We note that under the equivalence between bifibrations and correspondences, base changes on the cartesian side correspond the restriction along the first entry, while base changes on the cocartesian side correspond to restriction in the second entry. In particular, we may identify \(\Seq(\D,\Lag)\) with \(\Pairings(\Lag\op,\Lag,(\seq_{\D})|_{\Lag \times \Lag})\) as full subcategories of \(\Pairings(\D,\D,\seq_{\D})\), and so the projections
\[
\begin{tikzcd}
[row sep=1ex]
\Lag\op & \Seq(\D,\Lag)\ar[l] \ar[r] & \Lag \\
\x & {[\y \to \z \to \x]} \ar[l,mapsto] \ar[r,mapsto] & \y
\end{tikzcd}
\]
form a bifibration classified by the restricted correspondence \((\seq_{\D})|_{\Lag \times \Lag}\colon \Lag \times \Lag \to \Sps\). Under this equivalence \(\Seq(\D,\Lag) \simeq \D\) these projections correspond to adjoints \(p\colon \D \to \Lag\) and \(q \colon \D \to \Lag\op\) to the inclusions \(i\colon \Lag\hrar \D\) and \(j\colon \Lag\op \hrar \D\).
We now observe that for \(\x \in \Lag\op\) and \(\y \in \Lag\) we have
\[
\seq_{\D}(\x,\y) \simeq \Om\Map_{\D}(\Om j(\x),\Sig i(\y)) \simeq \Om\Map_{\D}(\Om j(\x),\Sig \Dual j(\y))
\simeq \Om^{\infty+1}\Bil_{\QFD}(\Om j(\x),\Om j(\y))
\]
and so if we let \(\Pi\colon \Lag \to \Spa\) be the quadratic functor given by the formula \(\Pi(\z) = \Om\QFD(\Om j(\z))\) then we obtain an equivalence
\[
\D \simeq \Pairings(\Lag\op,\Bil_{\Pi})
\]
on the level of stable \(\infty\)-categories, which is compatible
with the inclusions from and projections to \(\Lag\) and \(\Lag\op\) on both sides.

We now address the comparison of Poincaré structures.
For every \(\x \in \D\) the counit of \(i \dashv p\) and unit of \(q \dashv j\) yield a sequence
\begin{equation}
\label{equation:canonical-sequence-lagrangian}%
ip(\x) \to \x \to jq(\x)
\end{equation}
whose composite admits an essentially unique null-homotopy, since \(jq(\x) \in \ker(p)\) and hence \(\map_{\C}(ip(\x),jq(\x))= \map_{\Lag}(p(\x),pjq(\x))= 0\). This null-homotopy exhibits~\eqrefone{equation:canonical-sequence-lagrangian} as exact. Indeed, this sequence maps to an exact sequence by both \(p\) and \(q\) and these two functors are jointly conservative since \(\ker(p) \cap \ker(q) \simeq \im(j) \cap \ker(q) =0\).
Since \(\QFD(ip(\x)) =0\) the shifted exact sequence \(\Om jq(\x) \to ip(\x) \to \x\) determines an exact square
\[
\begin{tikzcd}
\QFD(\x) \ar[r] \ar[d] & \Om\QFD(\Om(jq(\x)) \ar[d] \\
\Om\Bil_{\QFD}(ip(\x),\Om jq(\x)) \ar[r] & \Om\Bil_{\QFD}(\Om jq(\x),\Om jq(\x))
\end{tikzcd}
\]
which, having set \(\Pi(\z) = \Om\QFD(\Om j(\z))\), we can write as
\begin{equation}
\label{equation:defining-square-generic}%
\begin{tikzcd}
\QFD(\x) \ar[r] \ar[d] & \Pi(q(\x)) \ar[d] \\
\Bil_{\QFD}(ip(\x),jq(\x)) \ar[r] & \Bil_{\Pi}(q(\x),q(\x)) \ .
\end{tikzcd}
\end{equation}
Comparing the exact square~\eqrefone{equation:defining-square-generic} with the exact square~\eqrefone{equation:defining-square-swap} we then conclude that the equivalence \(\D \simeq \Pairings(\Lag\op,\Bil_{\Pi})\) refines to an equivalence of Poincaré \(\infty\)-categories \((\D,\QFD) \simeq \Pairings(\Lag\op,\Pi)\). By Corollary~\refone{corolary:left-kan} we then get that the natural transformation \(\QFD \Rightarrow q^*\Pi\) furnished by the top row of the square~\eqrefone{equation:defining-square-generic} exhibits \(\Pi\) as the left Kan extension of \(\QFD\) along \(q\op\colon \D\op \to \Lag\). We hence an equivalence
\[
(\D,\QFD) \simeq \Pairings(\Lag,q_!\QFD) = \Pairings(\Lag\op,\Pi)
\]
compatible with the projections to (and hence also the embedding of) the \(\infty\)-categories \(\Lag\) and \(\Lag\op\) on both sides.
\end{proof}

Given a Poincaré \(\infty\)-category \((\D,\QFD)\) admitting a Lagrangian
\(\begin{tikzcd}
\Lag \ar[r] &
\D \ar[l,phantom,shift right=.8ex,"{\resizebox{1ex}{!}{\rotatebox{180}{$\perp$}}}"] \ar[l,bend right=25,"p"']
\end{tikzcd}\)
with associated projection \(q=\nolinebreak p\Dual\colon \D \to \Lag\op\), Proposition~\refone{proposition:recognize-poincare} yields an equivalence of Poincaré \(\infty\)-categories
\begin{equation}
\label{equation:D-is-pairings}%
(\D,\QFD)\simeq \Pairings(\Lag\op,q_!\QFD) .
\end{equation}
In particular, the association \((\C,\QF) \mapsto \Pairings(\C,\QF)\) takes values in metabolic Poincaré \(\infty\)-categories, and every metabolic Poincaré \(\infty\)-category is obtained in this manner.

In what follows, it will be useful to observe that the hermitian structure \(\Pi := q_!\QFD\) on \(\Lag\op\) can also be recovered using the inclusion \(j\colon \Lag\op \to \D\) right adjoint to \(q\). To avoid a potential confusion we emphasise that \(\Pi\) does not coincide with the restriction of \(\QFD\) along \(j\). Instead, let \(\Pi_{\pair}(\x,\y,\beta) = \Pi(\x) \times_{\Bil(\x,\x)}\map(\x,\y)\) be the quadratic functor on \(\Pairings(\Lag\op,\Pi)\) as in Construction~\refone{construction:hermitian-pairings}.
Under the equivalence~\eqrefone{equation:D-is-pairings} the functor \(j\colon \Lag\op \hrar \D\) sends \(\x\) to \((\x,0,0) \in \Pairings(\Lag\op,\Bil_{\Pi})\) and we have
\begin{equation}
\label{equation:restriction-swap}%
\QFD(j(\x)) = \Pi_{\pair}(\x,0,0) = \Pi(\x) \times_{\Bil_{\Pi}(\x,\x)}\map(\x,0) = \fib[\Pi(\x) \to \Bil_{\Pi}(\x,\x)] \simeq \Sig\Pi(\Sig\x).
\end{equation}
where the last equivalence is issued from Lemma~\refone{lemma:goodwillie} and Example~\refone{example:usual-sequence}. It will consequently be convenient to introduce the following terminology:

\begin{definition}
\label{definition:sigma-omega}%
Let \(\C\) be a stable \(\infty\)-category and \(\QF\) a quadratic functor on \(\C\). We will denote by \(\QF^{[\sig]}(\x) := \Om\QF(\Om\x)\) and
\(\QF^{[-\sig]}(\x) := \Sig\QF(\Sig\x)\) the quadratic functors obtained by pre- and post-composing with \(\Om\) and \(\Sig\), respectively. We note that the operations \(\QF \mapsto \QF^{[\sig]}\) and \(\QF \mapsto \QF^{[-\sig]}\) are inverse to each other, and in particular adjoint (in both directions). By Lemma~\refone{lemma:goodwillie} we have natural equivalences
\[
\QF^{[-\sig]}(\x) \simeq \fib[\QF(\x) \to \Bil_{\QF}(\x,\x)] \quad\text{and}\quad \QF(\x) \simeq \fib[\QF^{[\sig]}(\x) \to \Bil_{\QF^{[\sig]}}(\x,\x)] ,
\]
yielding in particular a natural transformation \(\QF \Rightarrow \QF^{[\sig]}\) and an adjoint transformation \(\QF^{[-\sig]} \Rightarrow \QF\).
\end{definition}

\begin{remark}
\label{remark:shifted-pairing}%
For a hermitian \(\infty\)-category \((\C,\QF)\), the underlying symmetric bilinear form of \(\QF^{[\sig]}\) is given by \((\x,\y) \mapsto \Om\Bil(\Om\x,\Om\y) \simeq \Sig^{\sig}\Bil(\x,\y)\), where \(\Sig^{\sig}\) is the operation of tensoring by the sign representation sphere, see discussion in \S\refone{subsection:herm-shifts}. In other words, \(\Bil_{\QF^{[\sig]}}\) has as underlying bilinear form \(\Sig\Bil_{\QF}\), but the symmetric structure is twisted by a sign, see Remark~\refone{remark:twistsign}. Similarly, \(\Bil_{\QF^{[-\sig]}} = \Sig^{-\sig}\Bil_{\QF} = \Om^{\sig}\Bil_{\QF}\) has underlying bilinear form \(\Om\Bil\), but the symmetric structure is twisted by a sign.
\end{remark}

The identification~\eqrefone{equation:restriction-swap} can be succinctly stated as \(j^*\QFD \simeq \Pi^{[-\sig]} = q_!\QFD^{[-\sig]}\), or equivalently, \(q_!\QFD \simeq j^*\QFD^{[\sig]}\).
We may summarise this discussion by extending~\eqrefone{equation:D-is-pairings} to
\begin{equation}
\label{equation:D-is-pairings-2}%
(\D,\QFD)\simeq \Pairings(\Lag\op,q_!\QFD) \simeq \Pairings(\Lag\op,j^*\QFD^{[\sig]}) .
\end{equation}
We now use the pairing construction in order to form left and right adjoints to the forgetful functor \(\Catp \to \Cath\).

\begin{proposition}
\label{proposition:first-adj-hermitian}%
For every hermitian \(\infty\)-category \((\C,\QF)\) and Poincaré \(\infty\)-category \((\E,\QFE)\), the map
\begin{equation}
\label{equation:first-adj-equivalence}%
\Map_{\Catp}((\E,\QFE),\Pairings(\C,\QF)) \to \Map_{\Cath}((\E,\QFE),(\C,\QF))
\end{equation}
induced by post-composition with the hermitian functor \((q,\eta)\colon \Pairings(\C,\QF) \to (\C,\QF)\) of~\eqrefone{equation:counit-hermitian}, is an equivalence of spaces.
In particular, the association \((\C,\QF) \mapsto \Pairings(\C,\QF)\) assembles to form a functor \(\Cath \to \Catp\) which is right adjoint to \(\Catp \to \Cath\).
\end{proposition}

\begin{remark}
Though one can show directly that the association \((\C,\QF) \mapsto \Pairings(\C,\QF)\) organises into a functor \(\Cath \to \Catp\),
Proposition~\refone{proposition:first-adj-hermitian}
is formulated in a way that does not require knowing this in advance, and on the other hand implies this functoriality via general principles of adjunctions; indeed, knowing that the comma category \(\Catp \times_{\Cath} (\Cath)_{/(\C,\QF)}\)
has a final object
is enough to imply the existence of the desired adjoints, which then must coincide with the given formula on objects.
\end{remark}

\begin{remark}
\label{remark:unit-fully-faithful}%
In the situation of Proposition~\refone{proposition:first-adj-hermitian}, if \((\C,\QF)\) is also Poincaré then we have a natural equivalence \(\Pairings(\C,\QF) \simeq \Ar(\C,\QF)\) (see Example~\refone{examples:typical-pairings}\refoneitem{item:pairing-arr-met}) under which the hermitian functor of~\eqrefone{equation:counit-hermitian} becomes the domain projection \(\Ar(\C,\QF) \to (\C,\QF)\). The unit of the adjunction furnished by Proposition~\refone{proposition:first-adj-hermitian} then corresponds to the essentially unique Poincaré functor \((\C,\QF) \to \Ar(\C,\QF)\) for which the composite with the domain projection is the identity on \((\C,\QF)\). In particular, the unit must coincide with the fully-faithful inclusion
\[
(\C,\QF) \to \Ar(\C,\QF) \quad\quad \x \mapsto [\id\colon \x \to \x]
\]
endowed with the natural equivalence \(\QF(\x) \simeq \QF_{\arr}([\id\colon \x \to \x])\).
\end{remark}

\begin{remark}
\label{remark:pairings-monoidal}%
By~\cite[Corollary 7.3.2.7]{HA} the functor \(\Pairings(-)\) inherits a lax symmetric monoidal structure by virtue of being right adjoint to the symmetric monoidal functor \(\iota\colon \Catp \to \Cath\) (see Theorem~\refone{theorem:tensorpoincare}). In particular, the adjunction \(\iota \dashv \Pairings(-)\) of Proposition~\refone{proposition:first-adj-hermitian} is a symmetric monoidal adjunction. It then follows from the equivalence of Example~\refone{examples:typical-pairings}(ii) that the functor \(\Ar(-)\colon \Catp \to \Catp\) is lax symmetric monoidal as well.
\end{remark}

\begin{remark}
\label{remark:first-adj-bilinear}%
In the situation of Proposition~\refone{proposition:first-adj-hermitian}, if \((\E,\QFE) = (\E,\QF^{\qdr}_{\Bil})\) is the quadratic Poincaré \(\infty\)-category associated to a symmetric bilinear form \(\Bil \in \Funs(\E)\), then by Propositions~\refone{proposition:adjunction-for-hermitian-categories} and~\refone{proposition:adjunction-for-poincare-categories}
the arrow~\eqrefone{equation:first-adj-equivalence} identifies with the arrow
\[
\Map_{\Catps}((\E,\Bil),(\Pairings(\C,\Bil_{\QF}),\Bil_{\pair})) \to \Map_{\Catsb}((\E,\Bil),(\C,\Bil_{\QF})),
\]
and so we may conclude that the association \((\C,\Bil) \mapsto (\Pairings(\C,\Bil),\Bil_{\pair})\) assembles to form a right adjoint to the inclusion \(\Catps \to \Catsb\). Identifying \(\Catps\) with \((\Catx)^{\hC}\) via Corollary~\refone{corollary:action-perfect-bilinear} and using Lemma~\refone{lemma:pair-structure-perfect} we may also reformulate this as saying that the association \((\C,\Bil) \mapsto (\Pairings(\C,\Bil),\Dual_{\pair})\) gives a right adjoint to the functor \((\Catx)^{\hC} \to \Catsb\) sending \((\C,\Dual)\) to \((\C,\Bil_{\Dual})\). This last conclusion could also be obtained from the opposite direction by showing that the association \((\A,\cB,\Bil) \mapsto \Pairings(\A,\cB,\Bil)\) gives a \(\Ct\)-equivariant right adjoint to the functor \(\Catx \to \Catb\) sending \(\C\) to \((\C,\C,m_{\C})\), and hence induces a right adjoint on the level of \(\Ct\)-fixed objects on both sides. In fact, \(\Pairings(-,-,-)\) being right adjoint to \(\C \mapsto (\C,\C,m_{\C})\) is a statement that holds also in the non-stable setting and can be proven using the setting of bifibrations as described in \S\refone{subsection:bifibrations}. Alternatively, an argument in the stable setting can be mounted along the lines of the proof of Proposition~\refone{proposition:first-adj-hermitian} below, using Remark~\refone{remark:recognition-bilinear} in place of Proposition~\refone{proposition:recognize-poincare}.
\end{remark}

\begin{proof}[Proof of Proposition~\refone{proposition:first-adj-hermitian}]
Fix a hermitian \(\infty\)-category \((\C,\QF)\) and a Poincaré \(\infty\)-category \((\E,\QFE)\), and let
\[
(\D,\QFD) := \Funx((\E,\QFE),\Pairings(\C,\QF)) = (\Funx(\E,\Pairings(\C,\Bil_{\QF})),\nat^{\QF_{\pair}}_{\QFE})
\]
be the corresponding internal hom Poincaré \(\infty\)-category constructed in \S\refone{subsection:internal}. We will use Proposition~\refone{proposition:recognize-poincare} in order to identify \((\D,\QFD)\)
with the pairing Poincaré \(\infty\)-category
associated to the internal hom hermitian \(\infty\)-category \(\Funx((\E,\QFE),(\C,\QF)) := (\Funx(\E,\C),\nat_{\QFE}^{\QF})\), thus reducing Proposition~\refone{proposition:first-adj-hermitian} to the algebraic Thom isomorphism of Proposition~\refone{proposition:algebraic-thom}. Indeed, define
\[
\Lag := \Funx(\E,\C\op)
\]
and let \(i_*\colon \Lag \to \D\) stand for post-composition with \(i\colon \C\op \hrar \Pairings(\C,\Bil_{\QF})\). Since \(i\) is fully-faithful and admits a right adjoint \(p \colon \Pairings(\C,\Bil_{\QF}) \to \C\op\) we have that \(i_*\) is fully-faithful and admits a right adjoint \(p_*\colon \D \to \Lag\) obtained by post-composing with \(p\). In addition, the restriction of \(\QFD = \nat_{\QFE}^{\QF_{\pair}}\) to \(\Lag\) vanishes because if \(f\colon \E \to \Pairings(\C,\Bil_{\QF})\) is an exact functor which factors through \(\C\op\) then \(f^*\QF_{\pair} = 0\) (since \(\QF_{\pair}\) vanishes on \(\C\op\)) and so \(\nat_{\QFE}^{\QF_{\pair}}(f) = \nat(\QFE,f^*\QF_{\pair}) \simeq  0\). Finally, the orthogonal complement \(\Lag^{\perp} \subseteq \D\) consists of those exact functors \(f\colon \E \to \Pairings(\C,\Bil_{\QF})\) such that \(p_*\Dual_{\D}(f) = 0\), that is, such that \(p\Dual_{\pair}f\Dual_{\E}(\x) = 0\) for every \(\x \in \E\). This is just equivalent to saying that \(f\) takes values in the orthogonal complement \((i\C\op)^{\perp}\), which coincides with \(i\C\op\) itself since \(i\C\op\) is a Lagrangian. We may then conclude that \(\Lag\) is a Lagrangian in \(\D\).

Now, if we identify \(\Lag\op = \Funx(\E,\C\op)\op = \Funx(\E\op,\C)\) with \(\Funx(\E,\C)\) via pre-composition with the duality of \(\E\), then the inclusion \(\Lag\op \to \D \) sending \(f\) to \(\Dual_{\D}(i_*f)\) identifies with \(f \mapsto j_*f\), where \(j_*\) denotes post-composition with \(j\colon \C \hrar \Pairings(\C,\Bil_{\QF})\). The left adjoint of \(j_*\) is then given by post-composition with \(q\colon \Pairings(\C,\Bil_{\QF}) \to \C\), which we denote by \(q_*\). Let \(\Pi = (q_*)_!\QFD\in \Funq(\Lag\op)\) be the quadratic functor obtained by left Kan extending \(\QFD\) along \(q_*\colon \D \to \Lag\op\). As in~\eqrefone{equation:D-is-pairings-2} we may also identify \(\Pi\) with the quadratic functor \(\Pi(g) = \QFD^{[\sig]}(jg)\). We may then compute
\[
\Pi(g) = \QFD^{[\sig]}(jg) \simeq \nat(\QFE,g^*j^*\QF_{\pair}^{[\sig]}) \simeq \nat(\QFE,g^*\QF) ,
\]
and identify the natural map
\[
\QFD(f) \to \Pi(qf)
\]
for \(f \in \D\)
with the map
\[
\eta_*\colon \nat(\QFE,f^*\QF_{\pair}) \to \nat(\QFE,f^*q^*\QF)
\]
obtained post-composition with \(f^*\eta \colon f^*\QF_{\pair} \Rightarrow f^*q^*\QF\).
Invoking Proposition~\refone{proposition:recognize-poincare} we now get an identification
\[
\begin{tikzcd}
\Funx((\E,\QFE),(\C,\QF)) \ar[d,equal] & \Pairings(\Funx((\E,\QFE),(\C,\QF))) \ar[l] \ar[r] \ar[d,"{\simeq}"] & (\Fun(\E,\C\op),0) \ar[d,equal] \\
\Funx((\E,\QFE),(\C,\QF)) & \Funx((\E,\QFE),\Pairings(\C,\QF)) \ar[l,"{(q_*,\eta_*)}"'] \ar[r,"{(p_*,0)}"] & (\Funx(\E,\C\op),0)
\end{tikzcd}
\]
of \((\D,\QFD) := \Funx((\E,\QFE),\Pairings(\C,\QF))\) as the pairings Poincaré category of the hermitian \(\infty\)-category \(\Funx((\E,\QFE),(\C,\QF))\), under which the associated cartesian
fibration
\[
\Pairings(\Funx((\E,\QFE),(\C,\QF))) \to \Funx((\E,\QFE),(\C,\QF))
\]
identifies with post-composition with \((q,\eta)\colon \Pairings(\C,\QF) \to (\C,\QF)\). Proposition~\refone{proposition:second-adj-hermitian} consequently follows from Proposition~\refone{proposition:algebraic-thom}.
\end{proof}

To obtain a left adjoint to the forgetful functor \(\Catp \to \Cath\) we first promote \(j\) to a hermitian functor
\begin{equation}
\label{equation:unit-hermitian}%
(j,\vartheta)\colon (\C,\QF^{[-\sig]}) \to \Pairings(\C,\QF) \quad\quad\quad j(\x) = (\x,0,0) \;,\; \vartheta\colon \QF^{[-\sig]}(\x) \simeq \QF_{\pair}(j(\x)).
\end{equation}
see Definition~\refone{definition:sigma-omega}.

\begin{proposition}
\label{proposition:second-adj-hermitian}%
For every Poincaré \(\infty\)-category \((\E,\QFE)\), the map
\begin{equation}
\label{equation:second-adj-equivalence}%
\Map_{\Catp}(\Pairings(\C,\QF),(\E,\QFE)) \to \Map_{\Cath}((\C,\QF^{[-\sig]}),(\E,\QFE))
\end{equation}
induced by pre-composition with~\eqrefone{equation:unit-hermitian}, is an equivalence of spaces. In particular, substituting \(\QF^{[\sig]}\) instead of \(\QF\) we deduce that the association \((\C,\QF) \mapsto \Pairings(\C,\QF^{[\sig]})\) assembles to form a functor \(\Cath \to \Catp\) which is left adjoint to the forgetful functor \(\Catp \to \Cath\).
\end{proposition}

\begin{remark}
\label{remark:second-counit}%
In the situation of Proposition~\refone{proposition:second-adj-hermitian}, if \((\C,\QF)\) is also Poincaré then by Example~\refone{examples:typical-pairings}\refoneitem{item:pairing-arr-met}
we have a natural equivalence
\[
\Pairings(\C,\QF^{[\sig]}) \simeq \Ar(\C,\QF^{[\sig]}) \simeq \Met(\C,\QF \circ \Om) \simeq \Met(\C,\QF) ,
\]
where the last equivalence covers the exact functor \([\cob \to \x] \mapsto [\Om\cob \to \Om\x]\). Under this equivalence the hermitian functor~\eqrefone{equation:unit-hermitian} becomes the inclusion
\[
\triv\colon (\C,\QF) \to \Met(\C,\QF) \quad\quad \x \mapsto [0 \to \x] .
\]
The counit of the adjunction furnished by Proposition~\refone{proposition:second-adj-hermitian} then corresponds to the essentially unique Poincaré functor \(\Met(\C,\QF) \to (\C,\QF)\) for which the pre-composing with this unit gives the identity on \((\C,\QF)\). In particular, the counit must coincide with the projection
\[
\met\colon \Met(\C,\QF) \to (\C,\QF) \quad\quad [\cob \to \x] \mapsto \x ,
\]
of Lemma~\refone{lemma:maps-with-metabolic-category}.
\end{remark}

\begin{remark}
\label{remark:dlag-is-canonical}%
For a Poincaré \(\infty\)-category \((\C,\QF)\), applying the functor \(\Pairings(-,(-)^{[\sig]})\) to the canonical hermitian functors \((\C,\QF) \to (\C,0)\) and \((\C,0) \to (\C,\QF)\) yields, using Remark~\refone{remark:second-counit}, the Poincaré functors
\[
\dlag\colon\Met(\C,\QF) \to \Hyp(\C)
\]
and
\[
\dilag\colon\Hyp(\C,\QF) \to \Met(\C)
\]
of Construction~\refone{construction:hyp-to-met}, respectively. This gives, in particular, a certain abstract justification for the appearance of these Poincaré functors.
\end{remark}

\begin{remark}
\label{remark:comonad}%
It follows from Proposition~\refone{proposition:second-adj-hermitian} and Remark~\refone{remark:second-counit} that the association \((\C,\QF) \mapsto \Met(\C,\QF)\) carries the structure of a comonad on \(\Catp\), with the Poincaré functors \(\met\colon \Met(\C,\QF) \to (\C,\QF)\) assembling to form the counit of this monad. A formal consequence of this which we record here for later use is that the resulting comultiplication Poincaré functor \(\Met(\C,\QF) \to \Met(\Met(\C,\QF))\) gives a section for either of the two projections
\[
\Met(\Met(\C,\QF)) \to \Met(\C,\QF) ,
\]
the first being the counit evaluated at \(\Met(\C,\QF)\) and the second obtained by applying \(\Met\) to the counit evaluated at \((\C,\QF)\).
\end{remark}

\begin{remark}
\label{remark:second-adj-bilinear}%
In the situation of Proposition~\refone{proposition:second-adj-hermitian}, if \((\E,\QFE) = (\E,\QF^{\sym}_{\Bil})\) is the symmetric Poincaré \(\infty\)-category associated to a symmetric bilinear form \(\Bil \in \Funs(\E)\), then by Propositions~\refone{proposition:adjunction-for-hermitian-categories} and~\refone{proposition:adjunction-for-poincare-categories} and Remark~\refone{remark:shifted-pairing},
the arrow~\eqrefone{equation:second-adj-equivalence} identifies with the arrow
\[
\Map_{\Catps}((\Pairings(\C,\Sig^{\sig}\Bil_{\QF}),\Bil_{\pair}),(\E,\Bil)) \to \Map_{\Catsb}((\C,\Bil_{\QF}),(\E,\Bil)).
\]
As in Remark~\eqrefone{remark:first-adj-bilinear} we may then conclude that the association \((\C,\Bil) \mapsto (\Pairings(\C,\Sig^{\sig}\Bil),\Dual_{\pair})\) assembles to form a left adjoint to the functor \((\Catx)^{\hC} \to \Catsb\) sending \((\C,\Dual)\) to \((\C,\C,m_{\C})\).
This conclusion could also be obtained differently by showing first that the association \((\A,\cB,\Bil) \mapsto \Pairings(\A,\cB,\Sig\Bil)\) gives  \(\Ct\)-equivariant left adjoint to the functor \(\Catx \to \Catb\) sending \(\C\) to \((\C,\C,m_{\C})\) (though the \(\Ct\)-equivariant structure here involves a somewhat subtle sign). The last claim can be proven using an argument similar to that of the proof of Proposition~\refone{proposition:second-adj-hermitian} below, by replacing the recognition principle of Proposition~\refone{proposition:recognize-poincare} by its bilinear version (see Remark~\refone{remark:recognition-bilinear}). We leave the details to the motivated reader.
\end{remark}

\begin{proof}[Proof of Proposition~\refone{proposition:second-adj-hermitian}]
Fix a hermitian \(\infty\)-category \((\C,\QF)\) and a Poincaré \(\infty\)-category \((\E,\QFE)\), and let
\[
(\D,\QFD) := \Funx(\Pairings(\C,\QF),(\E,\QFE)) = (\Funx(\Pairings(\C,\Bil_{\QF}),\E),\nat_{\QF_{\pair}}^{\QFE})
\]
be the corresponding internal hom Poincaré \(\infty\)-category. As in the proof of Proposition~\refone{proposition:first-adj-hermitian} we will use Proposition~\refone{proposition:recognize-poincare} in order to identify \(\D\) with a pairing Poincaré \(\infty\)-category associated to the internal hom hermitian \(\infty\)-category \(\Funx((\C,\QF),(\E,\QFE)) := (\Funx(\C,\E),\nat_{\QF}^{\QFE})\), thus reducing Proposition~\refone{proposition:second-adj-hermitian} to the algebraic Thom isomorphism of Proposition~\refone{proposition:algebraic-thom}. For this, define
\[
\Lag := \Funx(\C\op,\E)
\]
and let \(p^*\colon \Lag \to \D\) stand for pre-composition with \(p\colon \Pairings(\C,\Bil_{\QF}) \to \C\op\). Since \(p\) has a fully-faithful left adjoint \(i\colon \C\op \to \Pairings(\C,\Bil_{\QF})\) we get that \(p^*\) is fully-faithful and admits a right adjoint \(i^*\colon \D \to \Lag\) given by pre-composition with \(i\). In addition, the restriction of \(\QFD = \nat_{\QF_{\pair}}^{\QFE}\) to \(\Lag\) vanishes because if
\(f = g \circ p\colon \Pairings(\C,\Bil_{\QF}) \to \E\) for some \(g\colon \C\op \to \E\) then
\[
\nat(\QF_{\pair},f^*\QFE) = \nat(\QF_{\pair},p^*g^*\QFE) = \nat(p_!\QF_{\pair},g^*\QFE) = \nat(i^*\QF_{\pair},g^*\QFE) = 0
\]
where the identification \(p_!\QF_{\pair} \simeq i^*\QF_{\pair}\) is since \(i\op\) is right adjoint to \(p\op\).
Finally, the orthogonal complement \(\Lag^{\perp} \subseteq \D\) consists of those exact functors \(f\colon \Pairings(\C,\Bil_{\QF}) \to \E\) such that \(i^*\Dual_{\D}(f) = 0\), that is, such that \(\Dual_{\E}f\Dual_{\pair}(i(\x)) = 0\) for every \(\x \in \C\op\). This is just equivalent to saying that \(f\) vanishes on \(\im(j)=\ker(p)\), which is equivalent to saying that \(f\) factors through \(p\).
We may then conclude that \(\Lag\) is a Lagrangian in \(\D\).

Let us now identify \(\Lag\op = \Funx(\C\op,\E)\op = \Funx(\C,\E\op)\) with \(\Funx(\C,\E)\) via post-composition with the duality of \(\E\). Then the inclusion \(\Lag\op \to \D\) sending \(f\) to \(\Dual_{\D}(p^*f)\) identifies with \(f \mapsto q^*f\), where \(q^*\) denotes pre-composition with the cartesian projection \(q\colon \Pairings(\C,\Bil_{\QF}) \to \C\), and the left adjoint of \(q^*\) is given by pre-composition with \(j\colon \C \to \Pairings(\C,\Bil_{\QF})\).
Let \(\Pi := (j^*)_!\Psi \in \Funq(\Lag\op)\) be the quadratic functor obtained by left Kan extension \(\Psi\) along \(j^*\colon \D \to \Lag\op\), so that by~\eqrefone{equation:D-is-pairings-2} we can also write as \(\Pi(g) = \QFD^{[\sig]}(q^*(g))\) for \(g \in \Lag\op= \Funx(\C,\E)\). Using Corollary~\refone{corolary:left-kan} we then compute
\begin{equation*}\begin{split} \Pi(g) &= \QFD^{[\sig]}(gq) \simeq \nat(\QF_{\pair},q^*g^*\QFE^{[\sig]}) \simeq\\
&\simeq\nat(q_!\QF_{\pair},g^*\QFE^{[\sig]}) \simeq \nat(\QF,g^*\QFE^{[\sig]})\simeq \nat(\QF^{[-\sig]},g^*\QFE).\end{split}\end{equation*}
The canonical map
\[
\QFD(f) \to \Pi(fj)
\]
for \(f \in \D\)
then identifies with the map
\[
\vartheta^*\colon \nat(\QF_{\pair},f^*\QFE) \to \nat(\QF^{[-\sig]},j^*f^*\QFE)
\]
obtained by restricting along \(j\) using the equivalence \(\vartheta\colon \QF^{[-\sig]} \simeq j^*\QF_{\pair}\).
Invoking Proposition~\refone{proposition:recognize-poincare} we now get an identification
\[
\begin{tikzcd}
\Funx((\C,\QF^{[-\sig]}),(\E,\QFE)) \ar[d,equal] & \Pairings(\Funx((\C,\QF^{[-\sig]}),(\E,\QFE))) \ar[l] \ar[r] \ar[d,"{\simeq}"] & (\Fun(\C\op,\E),0)\ar[d,equal] \\
\Funx((\C,\QF^{[-\sig]}),(\E,\QFE)) & \Funx(\Pairings(\C,\QF),(\E,\QFE)) \ar[l,"{(j^*,\vartheta^*)}"'] \ar[r,"{(i^*,0)}"] & (\Funx(\C\op,\E),0)
\end{tikzcd}
\]
of \((\D,\QFD) := \Funx(\Pairings(\C,\QF),(\E,\QFE))\) as the pairings Poincaré category of the hermitian \(\infty\)-category \(\Funx((\C,\QF^{[-\sig]}),(\E,\QFE))\), under which the associated cartesian fibration
\[
\Pairings(\Funx((\C,\QF^{[-\sig]}),(\E,\QFE))) \to \Funx((\C,\QF^{[-\sig]}),(\E,\QFE))
\]
identifies with pre-composition with \((j,\vartheta)\colon (\C,\QF) \to \Pairings(\C,\QF)\). Proposition~\refone{proposition:first-adj-hermitian} consequently follows from Proposition~\refone{proposition:algebraic-thom}.
\end{proof}

We take the point of view that the functor \((\C,\QF) \mapsto \Pairings(\C,\QF^{[\sig]})\) is a categorical analogue of the algebraic Thom construction studied in \S\refone{subsection:algebraic-thom}. Recall that the latter takes a hermitian object and returns a metabolic Poincaré object with respect to a shifted Poincaré structure. In the algebraic setting we saw that this association determines an equivalence
\[
\spsforms(\C,\QF) \simeq \Poincdel(\C,\QF\qshift{1})
\]
between hermitian object in \((\C,\QF)\) and Poincaré objects in \((\C,\QF\qshift{1})\) equipped with a prescribed Lagrangian. To make the categorical analogue complete we would like to argue that \(\Pairings(-)\) determines an equivalence between \(\Cath\) and a suitable \(\infty\)-category whose objects are Poincaré \(\infty\)-categories \((\D,\QFD)\) equipped with a Lagrangian \(\Lag \subseteq \D\), and whose maps are Poincaré functors which preserve the given Lagrangians. While we will not make this completely precise, the gist of this claim amounts to the following two facts:
\begin{enumerate}
\item
The essential image of the functor \((\C,\QF) \mapsto \Pairings(\C,\QF^{[\sig]})\) consists of the metabolic Poincaré \(\infty\)-categories. This follows from Proposition~\refone{proposition:recognize-poincare}.
\item
Given two hermitian \(\infty\)-categories \((\C,\QF),(\Ctwo,\QFtwo)\), Poincaré functors
\[
\Pairings(\C,\QF^{[\sig]}) \to \Pairings(\Ctwo,\QFtwo^{[\sig]})
\]
sending the Lagrangian \(\C\op \subseteq \Pairings(\C,\QF^{[\sig]})\) to the Lagrangian \(\Ctwo\op \subseteq \Pairings(\Ctwo,\QFtwo^{[\sig]})\) are in bijection with hermitian functors \((\C,\QF) \to (\Ctwo,\QFtwo)\). Indeed, by Proposition~\refone{proposition:second-adj-hermitian} Poincaré functors \(\Pairings(\C,\QF^{[\sig]}) \to \Pairings(\Ctwo,\QFtwo^{[\sig]})\) correspond to hermitian functors \((\C,\QF) \to \Pairings(\Ctwo,\QFtwo^{[\sig]})\), and such a hermitian functor takes values in the full subcategory \(\Ctwo \subseteq \Pairings(\Ctwo,\QFtwo^{[\sig]})\) if and only if the corresponding Poincaré functor \(\Pairings(\C,\QF^{[\sig]}) \to \Pairings(\Ctwo,\QFtwo^{[\sig]})\) sends \(\C\) to \(\Ctwo\), which is equivalent to sending the Lagrangian \(\C\op \subseteq \Pairings(\C,\QF^{[\sig]})\) to the Lagrangian \(\Ctwo\op \subseteq \Pairings(\Ctwo,\QFtwo^{[\sig]})\) since \(\C\) and \(\C\op\) are two full subcategories of \(\Pairings(\C,\Bil_{\QF})\) which are switched by the duality, and the same holds for \(\Ctwo\) and \(\Ctwo\op\).
\end{enumerate}

\subsection{Genuine semi-additivity and spectral Mackey functors}
\label{subsection:mackey-functors}%

In this section we will explain how various categorical structures appearing in the theory of Poincaré \(\infty\)-categories can be neatly encoded in framework of \(\Ct\)-categories and Mackey functors, as developed by Barwick and collaborators in the setting of \emph{parametrised higher category theory}, see \cite{Barwick-MackeyI}, \cite{barwick2016parametrized}, \cite{Nardin-stability}, \cite{shah2018parametrized}. The material of this section, and in particular the \emph{hyperbolic Mackey functor} constructed in Corollary~\refone{corollary:hyperbolic-mackey} below, will form the basis to the formation of the real \(\K\)-theory spectrum in \papertwo.

To begin, let \(\OCt\) be the orbit category of \(\Ct\), that is, the category of transitive \(\Ct\)-sets and \(\Ct\)-equivariant maps. We note that \(\OCt\) has two objects, \(\Ct/\Ct = *\) and \(\Ct/e = \Ct\), such that \(\ast\) is terminal, \(\Hom_{\OCt}(\Ct,\Ct) = \Ct\) and there are no maps from \(\ast\) to \(\Ct\). A \(\Ct\)-category is by definition a cocartesian fibration
\[
\pi\colon \E \to \OCt\op ,
\]
which by the straightening-unstraightening equivalence is the same data as a functor \(\OCt\op \to \Cat\). A \(\Ct\)-functor between \(\Ct\)-categories is then a functor over \(\OCt\) which preserves cocartesian edges. We note that by the above explicit description of \(\OCt\) we see that it is isomorphic to the categorical cone on the category \(\BC\). As a result, the data of a functor \(\OCt\op \to \Ct\) is equivalent to that of an \(\infty\)-category \(\E_{\ast}\) (the image of the ``cone point'' \(\ast\)), an \(\infty\)-category \(\E_{\Ct}\) with \(\Ct\)-action (the image of \(\Ct\) with the \(\Ct\)-action induced by its automorphisms) and a \(\Ct\)-equivariant map \(\E_{\ast} \to \E_{\Ct}\) where the domain is considered with the trivial \(\Ct\)-action. Since \(\Cat\) admits limits the data of a \(\Ct\)-equivariant map \(\E_{\ast} \to \E_{\Ct}\) can equivalently be encoded via a map \(\E_{\ast} \to \E_{\Ct}^{\hC}\).

\begin{example}
\label{example:borel}%
If \(\C\) is an \(\infty\)-category with a \(\Ct\)-action then we can right Kan extend the functor \(\BC \to \Cat\) encoding this action to a functor \(\OCt\op \to \Cat\), which we can then straighten to obtain a \(\Ct\)-category \(\E \to \OCt\op\) with fibres \(\E_{\Ct} \simeq \C\) and \(\E_{\ast} \simeq \C^{\hC}\), and structure map \(\E_{\ast} \to \E_{\Ct}^{\hC}\) the identity. This construction embeds \(\Fun(\BC,\Cat)\) as a full subcategory of \(\Ct\)-categories.
\end{example}

\begin{example}
\label{example:diagonal}%
In the situation of Example~\refone{example:borel}, if \(\C\) is of the form \(\D \times \D\) with the flip action then \(\C^{\hC} \simeq \D\) and the \(\Ct\)-equivariant functor \(\E_{\ast} \to \E_{\Ct}\) is the diagonal \(\D \to \D \times \D\).
\end{example}

The examples we will be interested in are the following:
\begin{examples}
\label{examples:C2-categories-of-interest}%
\ %
\begin{enumerate}
\item
\label{item:funq}%
For a stable \(\infty\)-category \(\C\) the functor of taking symmetric bilinear parts \(\Bil_{(-)}\colon \Funq(\C) \to \Funs(\C) = \Funb(\C)^{\hC}\) determines a \(\Ct\)-category
\[
\uFunq(\C) \to \OCt\op
\]
whose fibre over \(\ast\) is \(\Funq(\C)\) and whose fibre over \(\Ct\) is \(\Funb(\C)\).
\item
\label{item:cath}%
The functor \(\Cath \to \Catsb = (\Catb)^{\hC}\) sending a hermitian \(\infty\)-category \((\C,\QF)\) to its underlying symmetric category \((\C,\Bil)\) determines a \(\Ct\)-category
\[
\uCath \to \OCt\op
\]
whose fibre over \(\ast\) is \(\Cath\) and whose fibre over \(\Ct\) is \(\Catb\).
\item
\label{item:catp}%
The functor \(\Catp \to \Catps = (\Catx)^{\hC}\) sending a Poincaré \(\infty\)-category \((\C,\QF)\) to its underlying \(\infty\)-category with perfect duality determines a \(\Ct\)-category
\[
\uCatp \to \OCt\op
\]
whose fibre over \(\ast\) is \(\Catp\) and whose fibre over \(\Ct\) is \(\Catx\).
\end{enumerate}
\end{examples}

For a \(\Ct\)-category \(\E \to \OCt\op\) one may consider the \(\Ct\)-variants of the usual notions of limits and colimits,
defined for a given \(\Ct\)-functor \(p\colon\I \to \E\). If \(\I \to \OCt\op\) is equivalent to a projection \(K \times \OCt\op \to \OCt\op\) for some \(K\) then a \(\Ct\)-colimit of \(p\colon \I \to \E\) is given by a cocartesian section \(s\colon\OCt\op \to \E\) together with a natural transformation \(\eta\colon p \Rightarrow s|_{K}\) which exhibits \(s\) as a colimit fibrewise, and dually for limits. For example, a \(\Ct\)-initial object is given by a section \(s\colon \OCt\op \to \C\) which is fibrewise initial. We shall refer to these as \defi{fibrewise} \(\Ct\)-colimits. By~\cite[Proposition~2.11]{Nardin-stability} a \(\Ct\)-category \(\E\) has all fibrewise \(\Ct\)-colimits indexed by \(K \times \OCt\op \to \OCt\op\) if and only if the fibres \(\E_{\Ct}\) and \(\E_{\ast}\) both have \(K\)-indexed colimits and the functor \(\E_{\ast} \to \E_{\Ct}\) preserves \(K\)-indexed colimits. In particular, all the examples in~\refone{examples:C2-categories-of-interest} have all fibrewise \(\Ct\)-limits and \(\Ct\)-colimits, since in all three cases the individual fibres have all limits and colimits and the cocartesian transition functor preserves all limits and colimits.

For an indexing \(\Ct\)-category \(\I \to \OCt\op\) which is not of the form \(K \times \OCt\op\) the notion of a \(\Ct\)-colimit is a bit more involved. The underlying data is still given by that of \(s\) and \(\eta\), but the condition these are required to satisfy is more complicated, and is neither weaker nor stronger than being a fibrewise colimit.
To avoid a technical digression let us avoid giving the general definition, referring the reader to~\cite[Definition 5.2]{shah2018parametrized}. The simplest type of non-fibrewise \(\Ct\)-(co)limits are finite \(\Ct\)-(co)products, and these will be the only type of non-fibrewise \(\Ct\)-(co)limits that we will consider here. These are \(\Ct\)-(co)limits indexed by finite \(\Ct\)-sets, that is, \(\Ct\)-categories \(\I \to \OCt\op\) which are finite direct sums of corepresentable left fibrations.
We may decompose them as a disjoint union of the finite \(\Ct\)-set \([\ast \to \ast]\) (standing for the left fibration over \(\OCt\op\) corepresented by \(\ast\)) and the finite \(\Ct\)-set \([\emptyset \to \Ct]\) (standing for the left fibration corepresented by \(\Ct\)).

We wish to verify that our examples of interest~\refone{examples:C2-categories-of-interest} all have finite \(\Ct\)-products and coproducts. For this we will use a convenient criterion  from~\cite{Nardin-stability}. Before we can state it, we point out the following observation: if \(\E \to \OCt\op\) is a \(\Ct\)-category with associated \(\Ct\)-equivariant functor \(f\colon \E_{\ast} \to \E_{\Ct}\), and \(g\colon \E_{\Ct} \to \E_{\ast}\) is a left or right adjoint to \(f\), then \(g\) inherits a canonical \(\Ct\)-equivariant structure. In fact, the entire adjunction carries a \(\Ct\)-action, so that the unit and counit are \(\Ct\)-equivariant natural transformations. This essentially follows from the uniqueness of adjoints given their existence. Otherwise put, the functor that forgets an adjunction to its left adjoint is fully-faithful and hence any \(\Ct\)-action can be lifted along it. One can also see this as follows. If \(f\) admits a right adjoint then the cocartesian fibration \(\E \to \OCt\op\) is also locally cartesian (since \(\Ct \to \ast\) is the only arrow that is not an isomorphism in \(\OCt\)), and hence a cartesian fibration. This cartesian fibration then encodes the data of a \(\Ct\)-equivariant functor \(\E_{\Ct} \to \E_{\ast}\), which is right adjoint to \(f\). If a left adjoint to \(f\) is considered then the same argument can be made using the dual cartesian fibration \(\hat{\E} \to \OCt\), that is, the cartesian fibration classified by the same functor as \(\E \to \OCt\op\).

The following lemma is just an adaptation of \cite[Proposition 2.11]{Nardin-stability} to the case at hand:
\begin{lemma}
\label{lemma:extracted-from-nardin}%
Let \(\E \to \OCt\op\) be a \(\Ct\)-category with such that the fibres \(\E_{\Ct},\E_{\ast}\) admit finite coproducts and the functor \(\E_{\ast} \to \E_{\Ct}\) preserves finite coproducts. Write \(\sig\colon \E_{\Ct} \to \E_{\Ct}\) for the action of the generator of \(\Ct\). Then the following are equivalent:
\begin{enumerate}
\item
\label{item:E-finite-Ctwo-copropd}%
\(\E\) admits all finite \(\Ct\)-coproducts.
\item
\label{item:E-Ctwo-colim}%
\(\E\) admits \(\Ct\)-colimits for \(\Ct\)-diagrams indexed by the corepresentable \(\Ct\)-set \([\emptyset \to \Ct]\).
\item
\label{item:left-adjoint}%
The functor \(f\colon \E_{\ast} \to \E_{\Ct}\) admits a left adjoint \(g\colon \E_{\Ct} \to \E_{\ast}\) such that for \(\x \in \E_{\Ct}\) the map
\[
\x \coprod \sig(\x) \to fg(\x)
\]
adjoint to the fold map \(g(\x \coprod \sig(\x)) \simeq g(\x) \coprod g(\sigma(\x)) \simeq g(\x) \coprod g(\x) \to g(\x)\) is an equivalence.
\end{enumerate}
\end{lemma}
\begin{proof}
The equivalence of \refoneitem{item:E-finite-Ctwo-copropd} and \refoneitem{item:E-Ctwo-colim} follows from the fact that under the assumptions of the lemma \(\E \to \OCt\op\) has fibrewise coproducts (that is, coproducts for diagrams indexed by finite direct sums of \([\ast \to \ast]\)) by~\cite[Proposition 2.11]{Nardin-stability} and so the existence of all \(\Ct\)-coproducts reduces to the case of the corepresentable ones, of which \([\ast \to \ast]\) is trivial.
The equivalence of \refoneitem{item:E-finite-Ctwo-copropd} and \refoneitem{item:left-adjoint} follows from~\cite[Proposition 2.11]{Nardin-stability} since \refoneitem{item:left-adjoint} is simply a reformulation of the Beck-Chevalley criterion given there in the case of the unique non-invertible edge \(\Ct \to \ast\) of \(\OCt\).
\end{proof}

\begin{remark}
Lemma~\refone{lemma:extracted-from-nardin} has a dual version which is proven exactly the same way. It says that \(\E\) has a all finite \(\Ct\)-products if and only if it has \(\Ct\)-limits for diagrams indexed by \([\emptyset \to \Ct]\), and that the latter is equivalent to \(f\) having a right adjoint \(g\colon \E_{\Ct} \to \E_{\ast}\) such that for \(\x \in \E_{\Ct}\) the map
\[
fg(\x) \to \x \times \sig(\x)
\]
adjoint to the diagonal \(g(\x) \to \g(\x) \times g(\x) \simeq g(\x) \times g(\sig(\x)) \simeq g(\x \times \sigma(\x))\) is an equivalence.
\end{remark}

\begin{remark}
\label{remark:explicit-coproduct}%
Diagrams in \(\E\) indexed by \([\emptyset \to \Ct]\) are determined by the data of an object \(\x \in \E_{\Ct}\). When the equivalent conditions of Lemma~\refone{lemma:extracted-from-nardin} hold then the \(\Ct\)-coproduct of such a diagram is given by the cocartesian section \(s\colon \OCt\op \to \E\) whose value at \(\ast\) is \(g(x)\) and whose value at \(\Ct\) is \(x \coprod \sigma(\x)\). The existence of such a cocartesian section is insured by \refoneitem{item:left-adjoint} above. Similarly, the \(\Ct\)-product of such a diagram, when exists, is given by a cocartesian section \(s\colon \OCt\op \to \E\) whose value at \(\ast\) is \(g(x)\) and whose value at \(\Ct\) is \(x \times \sigma(\x)\).
\end{remark}

\begin{example}
\label{example:diagonal-products}%
In the situation of Example~\refone{example:diagonal}, the resulting \(\Ct\)-category \(\E \to \OCt\op\) has finite \(\Ct\)-products (resp.\ \(\Ct\)-coproducts) if and only if \(\D\) has finite products (resp.\ coproducts).
\end{example}

\begin{proposition}
\label{proposition:examples-have-C2-coproducts}%
The \(\Ct\)-categories of Examples~\refone{examples:C2-categories-of-interest} all have finite \(\Ct\)-products and coproducts.
\end{proposition}
\begin{proof}
We will use the criterion of Lemma~\refone{lemma:extracted-from-nardin}. For this we first verify that the fibres over \(\Ct\) and \(\ast\) have finite (co)products and that the transition functors preserve finite (co)products. In Example~\refoneitem{item:funq}
the fibres are stable and the transition functors are exact.
For Example~\refoneitem{item:cath}
the existence of limits and colimits is established in Proposition~\refone{proposition:cath-has-limits-and-colimits} and Remark~\refone{remark:catb-has-all-limits-and-colimits}, and the preservation of limits and colimits by the transition functors is a consequence of Proposition~\refone{proposition:generalized-hyp}.
For Example~\refoneitem{item:catp} the same is established in Proposition~\refone{proposition:catx-has-limits-and-colimits}, Proposition~\refone{proposition:Catp-cocomplete} and Corollary~\refone{corollary:hyp-is-adjoint}. We now verify that all three examples satisfy Criterion~\refoneitem{item:left-adjoint} of Lemma~\refone{lemma:extracted-from-nardin}. In the case of Example~\refoneitem{item:funq} it follows from Lemma~\refone{lemma:universal-crs} and Remark~\refone{remark:adj-bilinear-part-diag} that the functor \(\Bil \mapsto \Bil^{\Del}\) gives both a left and a right adjoint to the bilinear part functor \(\QF \mapsto \Bil_{\QF}\).
We now need to verify that for \(\Bil \in \Funb(\C)\) the maps
\[
\Bil(\x,\y) \oplus \Bil(\y,\x) \to \fib[\Bil(\x\oplus\y) \to \Bil(\x,\x)\oplus \Bil(\y,\y)]
\]
and
\[
\cof[\Bil(\x,\x) \oplus \Bil(\y,\y) \to \Bil(\x\oplus \y)] \to \Bil(\x,\y) \oplus \Bil(\y,\x)
\]
are equivalences. Indeed, this follows directly from the bilinearity of \(\Bil\). For Example~\refoneitem{item:cath} we have by Proposition~\refone{proposition:generalized-hyp} that the association \((\A,\cB,\Bil) \mapsto (\A \times \cB\op,\Bil)\) gives both a left and a right adjoint to the functor \(\Cath \to \Catb\). Criterion~\refoneitem{item:left-adjoint} can then be deduced from its validity for Example~\refoneitem{item:funq} and for Example~\refone{example:diagonal-products} with \(\D = \Catx\). Finally, the case of Example~\refoneitem{item:catp} follows from that of~\refoneitem{item:cath} since the former maps to the latter via a \(\Ct\)-functor which is a fibrewise a replete subcategory inclusion, so the two-sided adjoint of \(\Cath \to \Catb\) restricts to give two sided adjoints for \(\Catp \to \Catx\) by Proposition~\refone{proposition:generalized-hyp}.
\end{proof}

It will be important for us in subsequent instalments of this project to know that the \(\Ct\)-categories of Examples~\refone{examples:C2-categories-of-interest} don't just admit finite \(\Ct\)-products and coproducts but that they are furthermore \defi{\(\Ct\)-semiadditive}, see~\cite[Definition 5.3]{Nardin-stability}. To explain what this means
let \(\E \to \OCt\) be a \(\Ct\)-category which admits finite \(\Ct\)-products and coproducts, so that by Lemma~\refone{lemma:extracted-from-nardin} \(f\) admits both a left adjoint \(g\colon \E_{\Ct} \to \E_{\ast}\) and a right adjoint \(h\colon \E_{\Ct} \to \E_{\ast}\), and these satisfy \(fg(\x) \simeq \x \coprod \sigma(\x)\) and \(fh(\x) \simeq \x \times \sigma(\x)\). Suppose that the fibres \(\E_{\Ct}\) and \(\E_{\ast}\) are both semiadditive, so that we may identify products and coproducts and write them as direct sums
\[
fg(\x) \simeq \x \oplus \sig(\x) \simeq fh(\x) .
\]
Then we have natural candidate for a comparison map
\begin{equation}
\label{equation:comparison-left-right}%
g(\x) \to h(\x)
\end{equation}
which is adjoint to the map
\[
\x \to fh(\x) \simeq \x \oplus \sig(\x),
\]
corresponding to the inclusion of the component \(\x\).

The following definition is an adaptation of~\cite[Definition 5.3]{Nardin-stability} to the particular case where the base is \(\OCt\), making use of the fact that \(\OCt\) has a unique non-invertible arrow, given by \(\Ct \to \ast\).

\begin{definition}
\label{definition:Ct-semiadditive}%
Let \(\E \to \OCt\op\) be a \(\Ct\)-category which admits finite \(\Ct\)-products and finite \(\Ct\)-coproducts. Then \(\E\) is called \defi{\(\Ct\)-semiadditive} if the following holds:
\begin{enumerate}
\item
The fibres \(\E_{\Ct}\) and \(\E_{\ast}\) are semiadditive.
\item
The comparison map~\eqrefone{equation:comparison-left-right} between the left and right adjoints of \(f\) is an equivalence.
\end{enumerate}
\end{definition}

\begin{example}
\label{example:diagonal-products-semiadditive}%
In the situation of Example~\refone{example:diagonal}, the resulting \(\Ct\)-category \(\E \to \OCt\op\) is semiadditive if and only if \(\D\) is semiadditive.
\end{example}

\begin{proposition}
\label{proposition:examples-are-C2-semi-additive}%
The \(\Ct\)-categories of Examples~\refone{examples:C2-categories-of-interest} are all \(\Ct\)-semiadditive.
\end{proposition}
\begin{proof}
We first verify that in all the examples in~\refone{examples:C2-categories-of-interest} the fibres are semi-additive. For Example~\refoneitem{item:funq} the fibres are stable and in particular semiadditive. For Examples~\refoneitem{item:cath} and~\refoneitem{item:catp} this was established in Proposition~\refone{proposition:catp-pre-add} and Remark~\refone{remark:catb-pre-add}.

We now establish the second condition of Definition~\refone{definition:Ct-semiadditive}. Arguing as in the proof of Proposition~\refone{proposition:examples-have-C2-coproducts} using Example~\refone{example:diagonal-products-semiadditive} in place of Example~\refone{example:diagonal-products} we see that it will suffice to establish the second condition for Example~\refoneitem{item:funq}. Now by Remark~\refone{remark:adj-bilinear-part-diag} the diagonal restriction functor \(\Del^*\colon \Funb(\C) \to \Funq(\C)\) is both left and right adjoint to the cross effect functor \(\Bil_{(-)}\colon \Funq(\C) \to \Funb(\C)\) and we already saw in the proof of Proposition~\refone{proposition:examples-have-C2-coproducts} that the composite \(\Funb(\C) \to \Funq(\C) \to \Funb(\C)\) is naturally equivalent to the functor \(\Bil \mapsto \Bil \oplus \Bil_{\swap}\). Unwinding the definitions, to establish the second condition it will suffice to show that unit of the adjunction \(\Del^* \dashv \Bil_{(-)}\) is given by the component inclusion
\[
\Bil(\x,\y) \to \Bil(\x,\y) \oplus \Bil(\y,\x),
\]
and the counit of the adjunction \(\Bil_{(-)} \dashv \Del^*\) is given by the component projection
\[
\Bil(\x,\y) \oplus \Bil(\y,\x) \to \Bil(\x,\y) .
\]
Indeed, this is established in Remark~\refone{remark:adj-bilinear-part-diag}.
\end{proof}

\begin{remark}
\label{remark:explicit-biproduct}%
If \(\E \to \OCt\op\) is a \(\Ct\)-semiadditive \(\Ct\)-category then the \(\Ct\)-product and \(\Ct\)-coproduct of a \([\emptyset \to \Ct]\)-indexed \(\Ct\)-diagram in \(\E\) corresponding to an object \(\x \in \E_{\Ct}\) are both given by the cocartesian section \(s\colon \OCt\op \to \E\) whose value at \(\ast\) is \(g(x)\) and whose value at \(\Ct\) is \(x \oplus \sigma(\x)\), cf.\ Remark~\refone{remark:explicit-coproduct}.
\end{remark}

\begin{remark}
If \(\E \to \OCt\op\) is a \(\Ct\)-semiadditive \(\Ct\)-category then the functor \(f\colon \E_{\ast} \to \E_{\Ct}\) admits a two sided adjoint \(g\). The \(\Ct\)-equivariant structure induces a \(\Ct\)-equivariant structure on \(g\) in a priori two different ways: one by the uniqueness of \(g\) as a left adjoint of \(f\) and once by its uniqueness as a right adjoint. The comparison map~\eqrefone{equation:comparison-left-right} is however a natural transformation of \(\Ct\)-equivariant functors (since the component inclusion \(\x \to \x \oplus \sig(\x)\) is such), and so it identifies the left and right adjoints of \(f\) also as \(\Ct\)-equivariant functors. Similarly, the induced functor \(\ovl{g}\colon \E_{\Ct}^{\hC} \to \Fun(\BC,\E_{\ast})\) is a two sided adjoint to the induced functor \(\ovl{f}\colon \Fun(\BC,\E_{\ast}) \to \E_{\Ct}^{\hC}\).
\end{remark}

\begin{remark}
\label{remark:hyp-equivariant}%
Specializing to the case of the \(\Ct\)-category \(\uCatp \to \OCt\op\) we now get that the functor \(\Hyp\) inherits a \(\Ct\)-equivariant structure making it a two-sided \(\Ct\)-equivariant adjoint to \(\rU\colon \Catp \to \Catx\), and similarly the induced functor \(\Hyp^{\hC}\colon (\Catx)^{\hC} \to \Fun(\BC,\Catp)\) is a two-sided adjoint to \(\rU^{\hC}\colon \Fun(\BC,\Catp) \to (\Catx)^{\hC}\).
The composite
\[
\ol{\Hyp} \colon \Catp \to (\Catx)^{\hC} \xrightarrow{\Hyp^{\hC}} \Fun(\BC,\Catp)
\]
then determines a \(\Ct\)-action on \(\Hyp(\rU(\C,\QF)) = \Hyp(\C)\) for a Poincaré \(\infty\)-category. For a fixed \((\C,\QF)\), this is the \(\Ct\)-action of Construction~\refone{construction:equivariant-hyp}, but now promoted to be natural in \(\C\). Similarly, the \(\Ct\)-equivariance of the maps
\[
\Hyp(\C) \xrightarrow{\hyp} (\C,\QF) \xrightarrow{\fgt} \Hyp(\C)
\]
constructed in Lemma~\refone{lemma:fgt-hyp-equivariant} is now exhibited as the components of two natural transformations of \(\Ct\)-equivariant functors in \((\C,\QF)\).
\end{remark}

\begin{remark}
\label{remark:hyp-of-hyp}%
In the situation of Remark~\refone{remark:hyp-equivariant}, the functors \(\rU\) and \(\rU^{\hC}\) participate in a commutative square of forgetful functors
\[
\begin{tikzcd}
\Fun(\BC,\Catp) \ar[r] \ar[d,"{\rU^{\hC}}"'] & \Catp  \ar[d,"{\rU}"] \\
(\Catx)^{\hC} \ar[r] & \Catx  .
\end{tikzcd}
\]
Passing to left adjoints, we obtain a commutative square
\begin{equation}
\label{equation:square-of-left-adjoints}%
\begin{tikzcd}
\Catx \ar[r] \ar[d,"{\Hyp}"'] & (\Catx)^{\hC} \ar[d,"{\Hyp^{\hC}}"] \\
\Catp \ar[r] & \Fun(\BC,\Catp) .
\end{tikzcd}
\end{equation}
At the same time, since \(\Catx\) is semi-additive the top horizontal functor in~\eqrefone{equation:square-of-left-adjoints} is given by the symmetrisation \(\C \mapsto \C \times \C\op\) associated to the \(\mop\)-action on \(\Catx\). But \(\uCatp\) having finite \(\Ct\)-coproducts (Proposition~\refone{proposition:examples-have-C2-coproducts}) means that this symmetrisation is identified with \(\rU \circ \Hyp\) as a functor \(\Catx \to (\Catx)^{\Ct}\). Since \(\Catp\) is also semi-additive (Proposition~\refone{proposition:catp-pre-add}) it then follows that from the commutativity of the above square that the functor
\[
\Catp \to \Fun(\BC,\Catp) \quad\quad (\C,\QF) \mapsto \ol{\Hyp}(\Hyp(\C))
\]
is naturally equivalent to the functor sending \((\C,\QF)\) to \(\Hyp(\C) \times \Hyp(\C)\), equipped with the flip \(\Ct\)-action.
\end{remark}

We now wish to use the \(\Ct\)-semiadditivity of the Examples in~\refone{examples:C2-categories-of-interest} in order to extract extra structures in terms of Mackey functors. For this, let \(\Span(\Ct)\) be the span \(\infty\)-category of finite \(\Ct\)-sets, as defined in \cite[Df.~3.6]{Barwick-MackeyI}.
A \defi{Mackey object} in an additive \(\infty\)-category \(\A\) is by definition a product preserving functor from \(\Span(\Ct) \rightarrow \A\). If \(\A\) is taken to be \(\Spa\), then~\cite[Theorem A.4]{Nardin-stability} %
shows that the arising \(\infty\)-category underlines the various (Quillen equivalent) model categories classically used for the definition of genuine \(\Ct\)-spectra. Taking spectral Mackey functors as the definition of the latter we set
\[
\Spagc:=\Fun^\times(\Span(\Ct),\Spa).
\]
Evaluation at the finite \(\Ct\)-sets \(\Ct\) then defines the functor \(u \colon \Spagc \rightarrow \Spa^{\hC}\), by retaining the action of the span \(\Ct \xleftarrow{\id} \Ct \xrightarrow{+1} \Ct\). Evaluation at the one-point \(\Ct\)-set defines the genuine fixed points \(-^\gCt \colon \Spagc \rightarrow \Spa\). The datum of a genuine \(\Ct\)-spectrum thus is equivalent to the datum of the pair of spectra \((E^{g\Ct},E)\), together with a \(\Ct\)-action on \(E\) and restriction and transfer maps
\[
\mathrm{res}:E^{\gCt}\to E\qquad \mathrm{tr}:E\to E^{\gCt}
\]
coming from the spans
\begin{equation}
\label{equation:Ct-transfer-span}%
* \leftarrow \Ct \xrightarrow{\id} \Ct \quad \text{and} \quad \Ct \xleftarrow{\id} \Ct \rightarrow *
\end{equation}
with a host of compatibility data, and similarly for other target categories.

\begin{proposition}
\label{proposition:canonical-mackey}%
Let \(\E \to \OCt\op\) be a \(\Ct\)-semiadditive \(\infty\)-category with transition functor \(f\colon \E_\ast \to \E_{\Ct}\) and two-sided adjoint \(g\colon \E_{\Ct} \to \E_{\ast}\).
Then the identity functor \(\E_{\ast} \to \E_{\ast}\) canonically lifts to functor
\[
\begin{tikzcd}
 & \Fun^{\times}(\Span(\Ct),\E_{\ast}) \ar[d] \\
\E_{\ast} \ar[ur,dashed] \ar[r,equal] & \E_{\ast}
\end{tikzcd}
\]
where the vertical arrow is given by evaluation at \(\ast\). In addition the composite
\[
\E_* \to \Fun^{\times}(\Span(\Ct),\E_{\ast}) \xrightarrow{\ev_{\Ct}} \E_*,
\]
where \(\ev_{\Ct}\) denotes evaluation at \(\Ct \in \Span(\Ct)\), is  naturally equivalent to the functor \(\x \mapsto gf(\x)\).
\end{proposition}
\begin{proof}
Let \(\uAeff(\Ct) \to \OCt\op\) be the \(\Ct\)-Burnside \(\infty\)-category of \cite[Df.~4.12]{Nardin-stability}. The objects of \(\uAeff(\Ct)\) are given by arrows \(U \to V\) where \(U\) is a finite \(\Ct\)-set and \(V \in \OCt\) is a \(\Ct\)-orbit, and morphisms in \(\uAeff(\Ct)\) from \([U \to V]\) to \([U' \to V']\) are given by diagrams of the form
\[
\begin{tikzcd}
U \ar[d] & U''\ar[d] \ar[l] \ar[r] & U' \ar[d] \\
V & V' \ar[l]\ar[r,equal] & V'.
\end{tikzcd}
\]
The functor \(\uAeff(\Ct) \to \OCt\op\) is then given by \([U \to V] \mapsto V\), and is a cocartesian fibration whose fibre over \(V \in \OCt\) is the span \(\infty\)-category of finite \(\Ct\)-sets over \(V\).
Combining \cite[Pr.~5.11]{Nardin-stability} and \cite[Th.~6.5]{Nardin-stability} we have the evaluation at the object \([\ast\to \ast] \in \uAeff(\Ct)\) yields an equivalence
\[
\Fun_{\Ct}^\times(\uAeff(\Ct),\E) \st{\simeq}{\lrar} \E_{\ast} \,.
\]
Now the action of every \(\Ct\)-product preserving \(\Ct\)-functor \(\uAeff(\Ct) \to \E\) on fibres over \(\ast\) is again a product-preserving functor from \(\Span(\Ct)\) to \(\E_{\ast}\). Base change along \(\{\ast\} \subseteq \OCt\op\) then determines a functor
\[
\mathcal{R}\colon \E_{\ast}\simeq \Fun_{\Ct}^\times(\uAeff(\Ct),\E)\to \Fun^{\times}(\Span(\Ct),\E_{\ast}).
\]
equipped with a natural equivalence
\[
\mathcal{R}(\x) (\ast) \simeq \x\,,
\]
by construction.
Furthermore, for \(\x \in \E_{\ast}\) the Mackey functor \(\mathcal{R}(\x)\colon \Span(\Ct) \to \E_{\ast}\) is obtained in particular by restricting a \(\Ct\)-functor \(\uline{\mathcal{R}(\x)}\colon \uAeff(\Ct) \to \E\), i.e., a functor over \(\OCt\op\) which preserves cocartesian edges. Since \(\uline{\mathcal{R}(\x)}\) sends \([\ast \to \ast]\) to \(x \in \E_\ast \subseteq \E\) by construction it must
send the object \([\Ct \to \Ct] \in \uAeff(\Ct)\)
to \(f(\x) \in \E_{\Ct} \subseteq \E\).
Since \(\uline{\mathcal{R}(\x)}\) furthermore preserves \(\Ct\)-biproducts it must therefore send \([\Ct \to \ast] \in \uAeff(\Ct)\) to \(gf(\x)\), see Remark~\refone{remark:explicit-biproduct}.
\end{proof}

Applying Proposition~\refone{proposition:canonical-mackey} in the case of the \(\Ct\)-category \(\uFunq(\C) \to \OCt\) of Examples~\refone{examples:C2-categories-of-interest}\refoneitem{item:funq} we obtain:

\begin{corollary}
\label{corollary:quadratic-genuine}%
The inclusion \(\Funq(\C) \subseteq \Fun(\C\op,\Spa)\) admits a canonical lift to a functor
\[
\Funq(\C) \to \Fun(\C\op,\Spagc),
\]
In particular, every quadratic functor \(\QF\colon \C \to \Spa\) lifts canonically to a functor \(\wtl{\QF}\colon \C\op \to \Spagc\) valued in genuine \(\Ct\)-spectra, such that \(\wtl{\QF}(\x)\) has underlying $\Ct$-spectrum \(\Bil_{\QF}(\x,\x)\), geometric fixed points $\Lambda_\QF(\x)$ and genuine fixed points \(\QF(\x)\).
\end{corollary}
\begin{proof}
Only the claim about the geometric fixed point is not immediate from the application of Proposition~\refone{proposition:canonical-mackey} to Example~\refone{examples:C2-categories-of-interest}\refoneitem{item:funq}, but it follows from the defining cofibre sequence
\[\Bil_\QF(\x,\x)_\hC \longrightarrow \QF(\x) \longrightarrow \Lambda_\QF(\x),\]
which by construction matches the cofibre sequence
\[\wtl{\QF}(\x)_\hC \longrightarrow \wtl{\QF}(\x)^\gC \longrightarrow \wtl{\QF}(\x)^\geofix.\]
\end{proof}

\begin{corollary}[The hyperbolic Mackey functor]
\label{corollary:hyperbolic-mackey}%
The construction of hyperbolic categories canonically refines to a functor
\[
\gHyp\colon\Catp\lrar \Fun^{\times}(\Span(\Ct),\Catp)
\]
together with natural equivalences of Poincaré \(\infty\)-categories
\[
[\gHyp(\C,\QF)] (*) \simeq (\C,\QF)\,,
\]
    and a natural \(\Ct\)-equivariant equivalence of Poincaré \(\infty\)-categories
\[
[\gHyp(\C,\QF)](\Ct) \simeq \Hyp\C\,.
\]
In addition, the resulting \(\Ct\)-equivariant functors
\[
\begin{tikzcd}
[row sep=1ex]
\Hyp(\C) \ar[r] \ar[d,phantom,"{\rotatebox{270}{$\simeq$}}"] & (\C,\QF) \ar[d,phantom,"{\rotatebox{270}{$\simeq$}}"] \ar[r] & \Hyp(\C) \ar[d,phantom,"{\rotatebox{270}{$\simeq$}}"] \\
\left[\gHyp(\C,\QF)\right](\Ct) \ar[r] & \left[\gHyp(\C,\QF)\right](\ast) \ar[r] & \left[\gHyp(\C,\QF)\right](\Ct)
\end{tikzcd}
\]
associated to the spans of \eqrefone{equation:Ct-transfer-span} are given by the functors \(\hyp\) and \(\fgt\) of~\eqrefone{equation:hyp-forget} which are the unit and counit of the two-sided adjunctions between \(\Catp\) and \(\Catx\) of Corollary~\refone{corollary:hyp-is-adjoint}.
\end{corollary}

Using \refone{corollary:quadratic-genuine} we can now also extend the periodicity considerations for Poincaré structures on module categories of Section \refone{subsection:herm-shifts} to the general case. Specifically, we have:

\begin{proposition}
\label{prop:shiftofqingeneral}%
For a hermitian structure $\QF$ on a stable category $\C$, there is a natural equivalence
\[\wtl{\QF}\qshift{n+m} \circ (\Sigma^n)\op \simeq \SS^{m-n\sigma} \otimes \wtl{\QF}\]
of functors $\C\op \rightarrow \Spagc$, and so in particular $\Omega^n \colon \C \rightarrow \C$ refines to an equivalence
\[(\C,\QF\qshift{2n}) \simeq (\C,(\SS^{n-n\sigma}\otimes \wtl{\QF})^\gC).\]
\end{proposition}

\begin{proof}
The proof is essentially identical to that of \refone{proposition:general-equivalence-of-poincare-infty-categories}. It clearly suffices to treat the case \(m=-n\), as both sides display the same behaviour under shifting, and since pre-composition with \(\Sig\) can be inverted and iterated, one may reduce to proving for \(n=1=-m\). Let \(\alp\colon \Lam_{\QF} \Rightarrow (\Bil_{\QF} \circ \Del)^{\tC}\) be the reference map of \(\QF\). Arguing as in the proof of Proposition~\refone{proposition:general-equivalence-of-poincare-infty-categories} one obtains a natural commutative diagram
\[
\begin{tikzcd}
\QF \circ \Sigma\op \ar[d]\ar[r] &
\Lambda_\QF \circ \Sigma\op \ar[d]\ar[r,equal] &
\Lambda_\QF \circ \Sigma\op \ar[d]\ar[r,"\sim"] &
\Om \circ \Lambda_\QF  \ar[d, "\Om\alp"] \\
(\Bil_\QF \circ \Delta  \circ \Sigma\op)^\hC \ar[d,"\sim"]\ar[r] &
(\Bil_\QF  \circ \Delta \circ \Sigma\op)^\tC\ar[d,"\sim"]\ar[r,"\sim"] &
(\Om \circ \Bil_\QF \circ\Delta)^\tC \ar[r,"\sim"]\ar[d,equal] &
\Om \circ (\Bil_\QF \circ\Delta)^\tC \ar[d,equal] \\
(\Om^{\rho} \circ \Bil_\QF \circ \Delta)^\hC \ar[r] &
(\Om^{\rho} \circ \Bil_\QF \circ \Delta)^\tC \ar[r,"\sim"] &
(\Om \circ \Bil_\QF \circ\Delta)^\tC \ar[r,"\sim"] &
\Om \circ (\Bil_\QF \circ\Delta)^\tC
\end{tikzcd}
\]
which may be viewed as an equivalence between the top left square and the external rectangle. As such, it identifies the classifying square of \(\QF \circ \Sig\op\) with one having bilinear term \(\Om^{\rho}\circ \Bil_{\QF}\), linear term \(\Om \circ \Lam_{\QF}\), and reference map \(\Om\alp\).
But these are precisely the underlying spectrum, geometric fixed point, and glueing map of \(\Sigma^{-\rho} \wtl{\QF}\). We immediately obtain that \(\QF\circ \Sigma\op\) agrees with the genuine fixed points of \(\Sigma^{-\rho} \wtl{\QF}$, and thus the second part (which suffices for all our applications) follows. But it is a folklore result that a genuine $\Ct$-spectrum is in fact naturally determined by the data we have computed above, namely, the diagram
\[\begin{tikzcd}
\Spa^\gC \ar[rr,"(-)^\geofix \Rightarrow (-)^\tC"] \ar[d,"u"] && \Ar(\Spa) \ar[d,"t"] \\
\Spa^\hC \ar[rr,"(-)^\tC"] && \Spa,
\end{tikzcd}
\]
where \(u\) takes the underlying \(\Ct\)-spectrum and \(t\) extracts the target of an arrow is cartesian; a proof of this fact can also be found in \reftwo{remark:descriptionofgenc2}. This implies the full statement of the proposition. %
\end{proof}

\subsection{Multiplicativity of Grothendieck-Witt and \(\L\)-groups}
\label{subsection:GW-L-multiplicative}%

In this section we will prove that the invariants \(\L_0(-)\) and \(\GW_0(-)\) defined in \S\refone{subsection:metabolic-and-L} and \S\refone{subsection:GW-group} are lax symmetric monoidal functors. In addition, the higher \(\L\)-groups organise into a symmetric monoidal functor to the category of graded abelian groups (with its Koszul symmetric monoidal structure). As a result, they carry a graded-commutative algebra structure when applied to symmetric monoidal Poincaré \(\infty\)-categories, such as those described in \S\refone{subsection:examples-monoidal}.

To begin, we first note that by Corollary~\refone{corollary:poinclax} the functor \(\pi_0\Poinc\colon \Catp \to \Set\) admits a canonical lax symmetric monoidal structure. This lax symmetric monoidal structure can be made quite explicit. Indeed, consider the assignment
\[
\pi_0 \Poinc(\C,\QF) \times \pi_0 \Poinc(\Ctwo,\QFtwo) \to \pi_0 \Poinc(\C \otimes \Ctwo,\QF \otimes \QFtwo) \ .
\]
For a pair of Poincaré objects \((x,q)\) and \((x',q')\) we get the Poincaré object \((x \otimes x', q \otimes q')\) in \((\C \otimes \Ctwo,\QF \otimes \QFtwo)\), where \(x \otimes x' \in \C \otimes \Ctwo\) is the image of \((x,x') \in \C \times \Ctwo\) under the universal bilinear functor \(\beta\colon \C \times \Ctwo \to \C \otimes \Ctwo\), and \(q \otimes q'\) denotes the map \(\SS \to [\QF \otimes \QFtwo](x \otimes x')\) obtained as the composite
\[
\SS = \SS \otimes \SS \xrightarrow{q \otimes q'} \QF(x) \otimes \QFtwo(x') = [\QF \boxtimes \QFtwo](x, x') \to
\App_2\beta_![\QF \boxtimes \QFtwo](x \otimes x') = [\QF \otimes \QFtwo](x \otimes x') \ .
\]
This object is Poincaré because its underlying bilinear form is given in light of Proposition~\refone{proposition:tensor} by the combination of the underlying bilinear forms of \(q\) and \(q'\).

We now wish to upgrade the above lax symmetric monoidal structure to the level of \(\Einf\)-spaces. For this, first note that since \(\Catp\) and \(\Cath\) are semi-additive (Proposition~\refone{proposition:catp-pre-add}) the corepresentable functors \(\Poinc\) and \(\spsforms\) canonically refine to functors with values in \(\Einf\)-spaces. Recall (see, e.g.,~\cite[Proposition 5.6]{Nik_Yoneda}) that the \(\infty\)-category \(\Mon_{\Einf}\) of \(\Einf\)-monoids carries a canonical symmetric monoidal structure such that the free-forgetful adjunction
\[
\rF\colon \Sps \adj {\Mon_{\Einf}}{\cocolon \rU}
\]
becomes symmetric monoidal (that is, its left adjoint is symmetric monoidal from which the right adjoint inherits a lax symmetric monoidal structure). We now claim that the \(\Einf\)-refinement \(\wtl{\Poinc}\colon \Catp \to \Mon_{\Einf}\) and \(\wtl{\spsforms}\colon \Cath \to \Mon_{\Einf}\)
also carry lax symmetric monoidal structures. This is in fact a completely formal consequence of the fact that the monoidal structure on \(\Catp\) and \(\Cath\) preserves direct sums:

\begin{lemma}
\label{lemma:general-mon-monoidal}%
Let \(\E\) be a small semi-additive \(\infty\)-category equipped with a symmetric monoidal structure \(\otimes\) which preserves direct sums in each variable. Then the lax symmetric monoidal structure of \(\Map_{\E}(1_{\E},-)\) canonically lifts to its
\(\Einf\)-refinement \(\wtl{\Map}_{\E}(1_{\E},-)\colon \E \to \Mon_{\Einf}\).
\end{lemma}

\begin{corollary}
The lax symmetric monoidal structure of \(\Poinc\) and \(\spsforms\) canonically lifts to their \(\Einf\)-refinements \(\wtl{\Poinc}\colon \Catp \to \Mon_{\Einf}\) and \(\wtl{\spsforms}\colon \Cath \to \Mon_{\Einf}\).
\end{corollary}
\begin{proof}
By possibly enlarging the universe we may assume that \(\Catp\) and \(\Cath\) are small. The claim then follows from Lemma~\refone{lemma:general-mon-monoidal}.
\end{proof}

\begin{proof}[Proof of Lemma~\refone{lemma:general-mon-monoidal}]
The full subcategory \(\Fun^{\times}(\E,\Mon_{\Einf}) \subseteq \Fun(\E,\Mon_{\Einf})\) spanned by the product-preserving functors is an accessible localization of \(\Fun(\E,\Mon_{\Einf})\) with a left adjoint which we will denote by \(L\colon \Fun(\E,\Mon_{\Einf}) \to \Fun^{\times}(\E,\Mon_{\Einf})\). By Lemma~\refone{lemma:day-localisation} this localization is compatible with Day convolution and extends to a symmetric monoidal localization
\[
L^{\otimes}\colon \Fun(\E,\Mon_{\Einf})^{\otimes} \to \Fun^{\times}(\E,\Mon_{\Einf})^{\otimes},
\]
where the codomain is endowed with the structure inherited from being a full suboperad of \(\Fun(\E,\Mon_{\Einf})^{\otimes}\). Now the identification \(\Map_{\E}(1_{\E},-) \simeq \rU\wtl{\Map}_{\E}(1_{\E},-)\) transposes to give a natural transformation of the form \(\rF\circ \Map_{\E}(1_{\E},-)\Rightarrow \wtl{\Map}_{\E}(1_{\E},-)\). Since \(\wtl{\Map}_{\E}(1_{\E},-)\) is product preserving this natural transformation induces a natural transformation
\[
L(\rF\circ \Map_{\E}(1_{\E},-))\Rightarrow \wtl{\Map}_{\E}(1_{\E},-).
\]
We claim that this last transformation is an equivalence. Note that this implies the desired claim via the symmetric monoidal structures of \(L\) and \(\rF\). Now, to prove the claim, it will suffice to show that for every product-preserving functor \(\G\colon \E \to \Mon_{\Einf}\) the induced map
\[
\Nat(\wtl{\Map}_{\E}(1_{\E},-),\G) \to \Nat(\rF\circ \Map_{\E}(1_{\E},-),\G)
\]
is an equivalence of spaces. Indeed, by adjunction we may also identify this map with the map
\[
\Nat(\wtl{\Map}_{\E}(1_{\E},-),\G) \to \Nat(\rU\circ\wtl{\Map}_{\E}(1_{\E},-),\rU\circ\G)
\]
induced by \(\rU \circ (-)\). This last map is an equivalence since the forgetful functor \(\rU \circ (-)\colon \Fun^{\times}(\E,\Mon_{\Einf}) \to \Fun^{\times}(\E,\Sps)\) is an equivalence on product-preserving functors.
\end{proof}

We now come to the main result of this subsection:

\begin{proposition}
\label{proposition:GW-L-lax}%
The functors \(\L_0, \GW_0\colon \Catp \to \CMon\) admit unique lax symmetric monoidal structures such that the transformations
\[
\pi_0\Poinc \Rightarrow \GW_0 \Rightarrow \L_0
\]
are symmetric monoidal, where \(\CMon\) stands for the (ordinary) symmetric monoidal category of commutative monoids.
\end{proposition}

\begin{remark}
\label{remark:also-abelian}%
The full subcategory \(\Ab \subseteq \CMon\) spanned by abelian groups is a full suboperad and a symmetric monoidal localisation of \(\CMon\). Since \(\GW_0\) and \(\L_0\) take values in \(\Ab\) the lax monoidal structures on \(\GW_0,\L_0\) and the map \(\GW_0 \Rightarrow \L_0\) equally applies if we consider \(\GW_0\) and \(\L_0\) as functors to \(\Ab\). The reason for working with the larger category of commutative monoids is to be able to make arguments pertaining to the natural transformation from \(\pi_0\Poinc\).
\end{remark}

The proof of Proposition~\refone{proposition:GW-L-lax} will require knowing certain multiplicative properties of the adjunctions \(\Cath \adj \Catp\) and \(\Catx \adj \Catp\), which we now verify.

\begin{lemma}
\label{lemma:pairing-projection}%
The adjunction \(\Cath \adj \Catp\) of Proposition~\refone{proposition:second-adj-hermitian},
in which the right adjoint \(\Catp \to \Cath\) is symmetric monoidal, satisfies the projection formula: for \((\C,\QF) \in \Catp\) and \((\Ctwo,\QFtwo) \in \Cath\)
the Poincaré functor
\begin{equation}
\label{equation:projection-formula}%
\Pairings(\C \otimes \Ctwo,(\QF \otimes \QFtwo)^{[\sig]}) \to (\C,\QF) \otimes \Pairings(\Ctwo,\QFtwo{^{[\sig]}})
\end{equation}
associated under this adjunction to the hermitian functor
\[
(\C \otimes \Ctwo,\QF \otimes \QFtwo) = (\C,\QF) \otimes (\Ctwo,\QFtwo) \to (\C,\QF) \otimes \Pairings(\Ctwo,\QFtwo{^{[\sig]}})
\]
induced by the unit hermitian functor \((\Ctwo,\QFtwo) \to \Pairings(\Ctwo,\QFtwo{^{[\sig]}})\), is an equivalence of Poincaré \(\infty\)-categories.
\end{lemma}
\begin{proof}
This is a formal consequence of (and actually equivalent to) the fact that the right adjoint \(\Catp \to \Cath\) is \emph{closed} symmetric monoidal, that is, preserves internal mapping objects, see Remark~\refone{remark:inclusion-closed}. Explicitly,
it will suffice to show that for every Poincaré \(\infty\)-category \((\E,\QFE)\) the restriction map
\[
\Map_{\Catp}((\C,\QF) \otimes \Pairings(\Ctwo,\QFtwo{^{[\sig]}}),(\E,\QFE)) \to \Map_{\Cath}((\C,\QF) \otimes (\Ctwo,\QFtwo),(\E,\QFE))
\]
is an equivalence of spaces. Indeed, since the Poincaré \(\infty\)-category \(\Funx((\C,\QF),(\E,\QFE))\) serves as an internal mapping object in both \(\Cath\) and \(\Catp\) we may identify this map with the map
\[
\Map_{\Catp}(\Pairings(\Ctwo,\QFtwo{^{[\sig]}}),\Funx((\C,\QF),(\E,\QFE))) \to \Map_{\Cath}((\Ctwo,\QFtwo), \Funx((\C,\QF),(\E,\QFE))),
\]
which is an equivalence by adjunction.
\end{proof}

Taking \(\QFtwo = 0\) in Lemma~\refone{lemma:pairing-projection} and using Remark~\refone{remark:units-and-counits} we immediately find:

\begin{corollary}
\label{corollary:hyp-projection}%
The adjunction \(\Hyp \dashv \rU\)
satisfies the projection formula: for \((\C,\QF) \in \Catp\) and \(\Ctwo \in \Cath\) the Poincaré functor
\[
\Hyp(\C \otimes \Ctwo) \to (\C,\QF) \otimes \Hyp(\Ctwo)
\]
induced by the component inclusion
\[
\C \otimes \Ctwo \to \C \otimes \rU\Hyp(\Ctwo) = \C \otimes [\Ctwo \oplus \Ctwo{\op}] = [\C \otimes \Ctwo] \oplus [\C \otimes \Ctwo{\op}]
\]
is an equivalence of Poincaré \(\infty\)-categories.
\end{corollary}

\begin{remark}
One can also deduce Corollary~\refone{corollary:hyp-projection} from the fact that the forgetful functor \(\rU\colon \Catp \to \Catx\) is closed symmetric monoidal, being the composition of the inclusion \(\Catp \to \Cath\) and the closed symmetric monoidal projection \(\Cath \to \Catx\), see Remark~\refone{remark:forgetful-closed} and the final part of Remark~\refone{remark:inclusion-closed}. In particular, the adjunction \(\Catx \adj \Cath\) also satisfies the projection formula, as is visible from the equivalence \((\C,\QF) \otimes (\Ctwo,0) \simeq (\C \otimes \Ctwo,0)\).
\end{remark}

\begin{remark}
\label{remark:projection-compatible}%
In the situation of Lemma~\refone{lemma:pairing-projection}, when \((\C,\QF)\) is also Poincaré then by the triangle identities the projection formula equivalence~\eqrefone{equation:projection-formula}  fits into a commutative triangle
\[
\begin{tikzcd}
\Pairings(\C \otimes \Ctwo,(\QF \otimes \QFtwo)^{[\sig]}) \ar[dr] \ar[rr,"{\simeq}"] &&  (\C,\QF) \otimes \Pairings(\Ctwo,\QFtwo{^{[\sig]}}) \ar[dl] \\
& (\C,\QF) \otimes (\Ctwo,\QFtwo)& \ ,
\end{tikzcd}
\]
where the diagonal arrows are obtained from the counit of the adjunction \(\Cath \adj \Catp\) at \((\C,\QF) \otimes (\Ctwo,\QFtwo)\) and \((\Ctwo,\QFtwo)\), respectively.
Using Example~\eqrefone{examples:typical-pairings}\refoneitem{item:pairing-arr-met} and Remark~\refone{remark:second-counit} we may also write this commutative triangle as
\begin{equation}
\label{equation:projection-compatible-met}%
\begin{tikzcd}
\Met((\C,\QF) \otimes (\Ctwo, \QFtwo)) \ar[dr,"{\met}"'] \ar[rr,"{\simeq}"] && (\C,\QF) \otimes \Met(\Ctwo,\QFtwo) \ar[dl,"{(\C,\QF) \otimes \met}"] \\
& (\C,\QF) \otimes (\Ctwo,\QFtwo) & \ .
\end{tikzcd}
\end{equation}
Applying the same argument for the projection formula of Corollary~\refone{corollary:hyp-projection} we similarly have the commutative triangle
\begin{equation}
\label{equation:projection-compatible-hyp}%
\begin{tikzcd}
\Hyp(\C \otimes \Ctwo) \ar[dr,"{\hyp}"'] \ar[rr,"{\simeq}"] && (\C,\QF) \otimes \Hyp(\Ctwo) \ar[dl,"{(\C,\QF) \otimes \hyp}"] \\
& (\C,\QF) \otimes (\Ctwo,\QFtwo) & \ .
\end{tikzcd}
\end{equation}
\end{remark}

\begin{remark}
\label{remark:hyp-met-modules}%
The commutative triangle~\eqrefone{equation:projection-compatible-met} taken with \((\Ctwo,\QFtwo) = (\Spaf,\QF^{\uni})\) yields an equivalence \(\Met(\C,\QF) \simeq (\C,\QF) \otimes \Met(\Spaf,\QF^{\uni})\) of Poincaré \(\infty\)-categories over \((\C,\QF)\). It then follows that when \((\C,\QF)\) is a symmetric monoidal Poincaré \(\infty\)-category the Poincaré \(\infty\)-category \(\Met(\C,\QF)\) acquires the structure of a module object over \((\C,\QF)\) (specifically, the free \((\C,\QF)\)-module generated from \(\Met(\Spaf,\QF^{\uni})\)) such that the functor
\[
\met\colon \Met(\C,\QF) \to (\C,\QF)
\]
is a map of (free) \((\C,\QF)\)-modules. Similarly, the commutative triangle~\eqrefone{equation:projection-compatible-hyp} taken with \((\Ctwo,\QFtwo) = (\Spaf,\QF^{\uni})\) yields an equivalence \(\Hyp(\C) \simeq (\C,\QF) \otimes \Hyp(\Spaf)\) over \((\C,\QF)\), and so when \((\C,\QF)\) is symmetric monoidal Poincaré we get that \(\Hyp(\C)\) acquires the structure of a module object over \((\C,\QF)\) (freely generated by \(\Hyp(\Spaf)\)) and
\[
\hyp\colon \Hyp(\C) \to (\C,\QF)
\]
is a map of free \((\C,\QF)\)-modules.
\end{remark}

\begin{proof}[Proof of Proposition~\refone{proposition:GW-L-lax}]
The uniqueness is clear since the maps \(\pi_0\Poinc(\C,\QF) \to \GW_0(\C,\QF) \to \L_0(\C,\QF)\) are surjective (and surjectivity is stable under tensor products in \(\CMon\)).
It will hence suffice to show that for every pair of Poincaré \(\infty\)-categories \((\C,\QF),(\Ctwo,\QFtwo)\) the dotted arrows
\begin{equation}
\label{equation:well-defined}%
\begin{tikzcd}
\pi_0 \Poinc(\C,\QF) \otimes \pi_0\Poinc(\Ctwo,\QFtwo) \ar[r] \ar[d] & \GW_0(\C,\QF) \otimes \GW_0(\Ctwo,\QFtwo) \ar[d,dashed] \ar[r] & \L_0(\C,\QF) \otimes \L_0(\Ctwo,\QFtwo) \ar[d,dashed] \\
\pi_0\Poinc((\C,\QF) \otimes (\Ctwo,\QFtwo)) \ar[r] & \GW_0((\C,\QF) \otimes (\Ctwo,\QFtwo)) \ar[r] & \L_0((\C,\QF) \otimes (\Ctwo,\QFtwo))
\end{tikzcd}
\end{equation}
exists to make the diagram commute, where the tensor products in the top row is that of commutative monoids.

Let us first treat the case of the functor \(\L_0\). For this it will suffice to show that for a pair of Poincaré objects \((\x, \qone)\) and \((\xtwo,\qtwo)\) such that \((\xtwo,\qtwo)\) is metabolic the associated Poincaré object \((\x,\qone) \otimes (\xtwo,\qtwo)\) is metabolic in \((\C,\QF) \otimes (\Ctwo,\QFtwo)\)
(the rest follows by symmetry).
But this follows directly from the commutativity of the diagram
\begin{equation}
\label{equation:well-defined-L}%
\begin{tikzcd}
[column sep=3ex]
\pi_0\Poinc\Met((\C,\QF) \otimes (\Ctwo, \QFtwo)) \ar[dr,"{[\met]}"'] \ar[r,"{\simeq}"] & \pi_0\Poinc((\C,\QF) \otimes \Met(\Ctwo,\QFtwo)) \ar[d] & \pi_0\Poinc(\C,\QF) \otimes \pi_0\Poinc\Met(\Ctwo,\QFtwo) \ar[l] \ar[d,"{\id \otimes [\met]}"] \\
 & \pi_0\Poinc((\C,\QF) \otimes (\Ctwo,\QFtwo)) & \pi_0\Poinc(\C,\QF) \otimes \pi_0\Poinc(\Ctwo,\QFtwo) \ar[l]
\end{tikzcd}
\end{equation}
given by the commutative triangle~\eqrefone{equation:projection-compatible-met}.

We now turn to \(\GW_0\). For the well-definedness of the middle dotted arrow in~\eqrefone{equation:well-defined}
we have to show that the relation \([\hyp(w)] \sim [x,q]\) for a Lagrangian \(w \to x\) is preserved under tensoring with some Poincaré object \((x',q')\). Given the commutativity of the diagram~\eqrefone{equation:well-defined-L} above, it will suffice to show that the diagram
\[
\begin{tikzcd}
\pi_0\Poinc\Met((\C,\QF) \otimes (\Ctwo, \QFtwo)) \ar[d,"{\dlag_*}"] \ar[r,"{\simeq}"] & \pi_0\Poinc((\C,\QF) \otimes \Met(\Ctwo,\QFtwo)) \ar[d,"{(\id \otimes \dlag)_*}"] & \pi_0\Poinc(\C,\QF) \otimes \pi_0\Poinc\Met(\Ctwo,\QFtwo) \ar[l] \ar[d,"{\id \otimes [\dlag]}"] \\
\pi_0\Poinc\Hyp(\C \otimes \Ctwo)\ar[dr,"{[\hyp]}"'] \ar[r,"{\simeq}"] & \pi_0\Poinc((\C,\QF) \otimes \Hyp(\Ctwo)) \ar[d] & \pi_0\Poinc(\C,\QF) \otimes \pi_0\Poinc\Hyp(\Ctwo,\QFtwo) \ar[l] \ar[d,"{\id \otimes [\hyp]}"] \\
 & \pi_0\Poinc((\C,\QF) \otimes (\Ctwo,\QFtwo)) & \pi_0\Poinc(\C,\QF) \otimes \pi_0\Poinc(\Ctwo,\QFtwo) \ar[l]
\end{tikzcd}
\]
is commutative as well, where the vertical arrows in the top rows are induced by the Poincaré functor \(\dlag\colon \Met(-)\to \Hyp(-)\) of Construction~\refone{construction:hyp-to-met}. Here, the squares on the right hand side are commutative as they are given by the lax monoidal structure on \(\pi_0 \Poinc\), and the bottom left triangle is induced by the commutative triangle~\eqrefone{equation:projection-compatible-hyp}. It will suffice to show that the top left square commutes. Indeed, this square is obtained by applying \(\pi_0\Poinc\) to the square of Poincaré \(\infty\)-categories
\[
\begin{tikzcd}
[column sep=7ex]
\Met((\C,\QF) \otimes (\Ctwo, \QFtwo)) \ar[d,"{\dlag}"'] \ar[r,"{\simeq}"] & (\C,\QF) \otimes \Met(\Ctwo,\QFtwo) \ar[d,"{(\C,\QF) \otimes \dlag}"] \\
\Hyp(\C \otimes \Ctwo) \ar[r,"{\simeq}"] & (\C,\QF) \otimes \Hyp(\Ctwo)  \ ,
\end{tikzcd}
\]
which commutes since it is obtained by evaluating the natural transformation~\eqrefone{equation:projection-formula} of the projection formula at the arrow \((\C,\QF) \to (\C,0)\) in \(\Cath\), see Remark~\refone{remark:dlag-is-canonical}.
\end{proof}

\begin{corollary}
\label{corollary:L-graded-lax}%
The functor \(\L_*\colon \Catp \to \mathrm{gr}\Ab\) admits a lax symmetric monoidal structure, where \(\mathrm{gr}\Ab\) denotes the category of \(\mathbb{Z}\)-graded abelian groups with the symmetric monoidal structure using the Koszul sign rules.
\end{corollary}
\begin{proof}
By definition we have that \(\L_n(\C, \QF) = \L_0(\C, \QF^{[-n]}) = \L_0(\C, \QF \otimes \mathbb{S}^{-n})\). Thus the claim follows by combining Proposition~\refone{proposition:GW-L-lax}, Remark~\refone{remark:minus-operation}, and the fact that \(\mathbb{S}^{-\bullet}\) is a graded commutative algebra in the homotopy category of spectra. Concretely, the structure maps are simply given by
\begin{align*}
\L_n(\C,\QF) &\otimes \L_m(\Ctwo,\QFtwo) = \L_0\left(\C, \QF^{[-n]}\right) \otimes \L_0\left(\Ctwo, \QFtwo^{[-m]}\right)\\
&\to \L_0\left(\C \otimes \Ctwo, \QF^{[-n]} \otimes\QFtwo^{[-m]}\right) = \L_{0}\left(\C \otimes\C, (\QF \otimes \QFtwo)^{[-n-m]}\right) = \L_{n+m}(\C \otimes \C, \QF \otimes \QFtwo)
\end{align*}
and the fact that it is symmetric follows as explained above.
\end{proof}

\begin{corollary}
\label{corollary:GW-L-rings}%
If \((\C,\QF)\) is a monoidal \(\infty\)-category then \(\GW_0(\C,\QF)\) acquires a ring structure and \(\L_\bullet(\C,\QF)\) a graded ring structure such that the natural map
\[
\GW_0(\C,\QF) \to \L_0(\C,\QF)
\]
is a ring homomorphism. If the monoidal structure is symmetric then \(\GW_0(\C,\QF)\) is commutative and \(\L_\bullet(\C,\QF)\) is graded commutative.
\end{corollary}

\begin{example}
\label{example:L-GW-rings}%
Let \(A \in \Alg_{\Einf}^{\hC}\) be a commutative ring spectrum equipped with a \(\Ct\)-action. Then
the Grothendieck-Witt groups \(\GW_0(\Modp{A},\QF^{\s}_A)\), \(\GW_0(\Modp{A},\QF^{\geq 0}_A)\) and \(\GW_0(\Modp{A},\QF^{\tate}_A)\) of the symmetric monoidal Poincaré \(\infty\)-categories of Examples~\refone{examples:monoidal} carry natural commutative ring structures, and similarly the corresponding graded \(\L\)-groups \(\L_\bullet(\Modp{A},\QF^{\s}_A)\), \(\L_\bullet(\Modp{A},\QF^{\geq 0}_A)\) and \(\L_\bullet(\Modp{A},\QF^{\tate}_A)\) carry canonical graded-commutative ring structures.
\end{example}

Combining Proposition~\refone{proposition:GW-L-lax} with Example~\refone{example:GW-of-hyp}, Example~\refone{example:hyp-monoidal} and Remark~\refone{remark:hyp-met-modules} we also get:

\begin{corollary}
\label{corollary:hyp-fgt-GW}%
If \((\C,\QF)\) is a symmetric monoidal \(\infty\)-category then
\[
[\fgt]\colon \GW(\C,\QF) \to \K_0(\C)
\]
is a map of rings and
\[
[\hyp]\colon \K_0(\C) \to \GW_0(\C,\QF)
\]
is a map of \(\GW_0(\C)\)-modules.
\end{corollary}

\begin{remark}
In the situation of Corollary~\refone{corollary:hyp-fgt-GW}, the \(\GW_0(\C,\QF)\)-module structure on \(\K_0(\C)\) could be considered ambiguous: on the one hand we have the module structure determined by the ring structure of \(\K_0(\C)\) via the ring map \(\fgt\colon \GW_0(\C,\QF) \to \K_0(\C)\), and on the other we have the module structure induced by the \((\C,\QF)\)-module structure on \(\Hyp(\C)\) of Remark~\refone{remark:hyp-met-modules} via the identification \(\K_0(\C) \cong \GW_0(\Hyp(\C))\).
These two modules structures however coincide. Indeed,
unwinding the definitions we see that the former structure is induced via the lax monoidal structure of \(\GW_0\) and the symmetric monoidal structure of \((\C,\QF)\) by the Poincaré functor
\[
(\C,\QF) \otimes \Hyp(\C) \to \Hyp(\C \otimes \C),
\]
corresponding via the adjunction \(\rU \dashv \Hyp\) to the exact functor \(\C \otimes [\C \oplus \C\op] \to \C \otimes \C\) induced by the projection \(\C \oplus \C\op \to \C\), while the latter module structure is induced in the same manner by the inverse of the Poincaré equivalence
\[
\Hyp(\C \otimes \C) \xrightarrow{\simeq} (\C,\QF) \otimes \Hyp(\C)
\]
of Corollary~\refone{corollary:hyp-projection}, corresponding via the adjunction \(\Hyp \dashv \rU\) to the exact functor \(\C \otimes \C \to \C \otimes [\C \oplus \C\op]\) induced by the inclusion \(\C \to \C \oplus \C\op\). It will hence suffice to verify that these Poincaré functors determines inverse equivalences between \(\Hyp(\C \otimes \C)\) and \((\C,\QF) \otimes \Hyp(\C)\). However, since Poincaré functors from (or to) \(\Hyp\) are determines by their underlying exact functors, it suffices to check that these underlying exact functors determine inverse equivalences between \(\C \otimes [\C \oplus \C\op]\) and \([\C \otimes \C] \oplus [\C\op \otimes \C\op]\). The latter is however a formal consequence of the fact that the monoidal structure on \(\Catx\) preserves direct sums in each variable.
\end{remark}

Invoking Examples~\refone{example:quadratic-module} and Example~\refone{example:module-truncated} we also have the following two corollaries:

\begin{corollary}
\label{corollary:quadratic-L-GW-module}%
Let \((\C,\QF)\) be a symmetric monoidal hermitian \(\infty\)-category with underlying bilinear part \(\Bil = \Bil_{\QF}\). Then the quadratic Grothendieck-Witt group \(\GW_0(\C,\QF^{\qdr}_{\Bil})\) is canonically a module over the ring \(\GW_0(\C,\QF)\) and the map
\[
\GW_0(\C,\QF^{\qdr}) \to \GW_0(\C,\QF)
\]
is a map of \(\GW_0(\C,\QF)\)-modules. Similarly, the quadratic \(\L\)-groups \(\L_\bullet(\C,\QF^{\qdr}_{\Bil})\) form a graded module over the graded ring \(\L_\bullet(\C,\QF)\) and the map
\[
\L_\bullet(\C,\QF^{\qdr}) \to \L_\bullet(\C,\QF)
\]
is a map of graded \(\L_\bullet(\C,\QF)\)-modules.
\end{corollary}

\begin{corollary}
\label{corollary:truncated-L-GW-truncated}%
Let \(A\) be a connective commutative ring spectrum equipped with an involution. Then for every \(m \in \ZZ\) the Grothendieck-Witt group \(\GW_0(\Modp{A},\QF^{\geq m}_A)\) of the truncated Poincaré structure of Example~\refone{example:truncation} is canonically a module over the \(\GW_0\)-ring \(\GW_0(\Modp{A},\QF^{\geq 0}_A)\) of the Poincaré structure of Example~\refone{examples:monoidal}\refoneitem{item:truncated}.
\end{corollary}